\documentclass[11pt,openany]{book}
\usepackage{amssymb}
\usepackage{latexsym, amsmath, amscd,amsthm}
\usepackage[fullpage]{uiucthesis}
\usepackage{diagrams}
\usepackage{ifthen}

\theoremstyle{plain}
\newtheorem{main}{Main~Theorem}
\newtheorem{theorem}[equation]{Theorem}
\newtheorem{proposition}[equation]{Proposition}
\newtheorem{claim}[equation]{Claim}
\newtheorem{corollary}[equation]{Corollary}

\newtheorem{conj}[equation]{Conjecture}
\newtheorem{cor}[equation]{Corollary}
\newtheorem{hyp}[equation]{Hypothesis}
\newtheorem{ia}[equation]{Inductive Assumption}   
\newtheorem{lemma}[equation]{Lemma}
\newtheorem{step}{Step}
\newtheorem{property}[equation]{Property}

\theoremstyle{definition}
\newtheorem{defn}[equation]{Definition}

\newtheorem{quest}[equation]{Question}
\newtheorem{remark}[equation]{Remark}

\newcommand{\hap}{\widehat{P}}
\newcommand{\hapn}{\widehat{P^{\nu}}}
\newcommand{\pw}{\widetilde{P}}
\newcommand{\qw}{\widehat{Q}}
\newcommand{\hb}{\hat{\beta}}
\newcommand{\ha}{\hat{\alpha}}
\newcommand{\qwn}{\widehat{Q^{\nu}}}
\newcommand{\hbn}{\hat{\beta}^{\nu}}
\newcommand{\cn}{\chi^{\nu}}

\newcommand{\pn}{P^{\nu}}
\newcommand{\qn}{Q^{\nu}}
\newcommand{\map}{\mathcal{P}}
\newcommand{\maq}{\mathcal{Q}}

\newcommand{\an}{\alpha^{\nu}}
\newcommand{\ans}{\alpha^{ \nu,*}}
\newcommand{\bns}{\beta^{ \nu,*}}
\newcommand{\bn}{\beta^{\nu}}
\newcommand{\pns}{P ^{ \nu,*}}
\newcommand{\qns}{Q ^{ \nu,*}}

\newcommand{\tria}{\textit{triangular}}

\newcommand{\ZZ}{\mathbb{Z}}
\newcommand{\CC}{\mathbb{C}}
\newcommand{\Syl}{\operatorname{Syl}}
\newcommand{\Hall}{\operatorname{Hall}}
\newcommand{\tin}{\text{\textrm{\textup{ in }}}}
\newcommand{\Gap}[1]{G_{\infty,#1}}
\newcommand{\Aut}{\operatorname{Aut}}

\newcommand{\DCC}{\text{DCC}}

\newcommand{\slgs}[1]{G^*_{#1, \lambda}}
\newcommand{\slg}[1]{G_{#1, \lambda}}
\newcommand{\slqs}[1]{Q^*_{#1, \lambda}}
\newcommand{\slps}[1]{P^*_{#1, \lambda}}
\newcommand{\sla}[1]{\alpha_{#1, \lambda}}
\newcommand{\slb}[1]{\beta_{#1, \lambda}}
\newcommand{\slq}[1]{Q_{#1, \lambda}}
\newcommand{\slqw}{\widehat{Q}_{\lambda}}
\newcommand{\slp}[1]{P_{#1, \lambda}}
\newcommand{\slc}[1]{\chi_{#1, \lambda}}
\newcommand{\subl}[2]{#1_{#2, \lambda}}

\newcommand{\smgs}[1]{G^*_{#1, \mu}}
\newcommand{\smg}[1]{G_{#1, \mu}}
\newcommand{\smqs}[1]{Q^*_{#1, \mu}}
\newcommand{\smps}[1]{P^*_{#1, \mu}}
\newcommand{\sma}[1]{\alpha_{#1, \mu}}
\newcommand{\smb}[1]{\beta_{#1, \mu}}
\newcommand{\smq}[1]{Q_{#1, \mu}}
\newcommand{\smp}[1]{P_{#1, \mu}}
\newcommand{\smc}[1]{\chi_{#1, \mu}}
\newcommand{\smqw}{\widehat{Q}_{\mu}}

\newcommand{\skgs}[1]{G^*_{#1, K}}
\newcommand{\skg}[1]{G_{#1, K}}

\newcommand{\ska}[1]{\alpha_{#1, K}}
\newcommand{\skb}[1]{\beta_{#1, K}}
\newcommand{\skq}[1]{Q_{#1, K}}
\newcommand{\skp}[1]{P_{#1, K}}
\newcommand{\skc}[1]{\chi_{#1, K}}
\newcommand{\skqw}{\widehat{Q}_{K}}

\newcommand{\pps}[1]{(P'_{#1})^*}
\newcommand{\aas}[1]{(\alpha_{#1}')^*}

\newcommand{\mbr}{\mathbb{R}}
\newcommand{\mbt}{\mathbb{T}}

\newcommand{\mbp}{\mathbb{P}}
\newcommand{\mbg}{\mathbb{G}}
\newcommand{\mba}{\mathbb{A}}
\newcommand{\mbn}{\mathbb{N}}

\newcommand{\T}{\Theta}

\newcommand{\VP}{\boldsymbol{ \Phi}}
\newcommand{\VPS}{\boldsymbol{ \Psi}}
\newcommand{\ba}{\boldsymbol{ \alpha}}

\newcommand{\be}{\boldsymbol{ \eta}}

\newcommand{\ma}{\mathbf{A}}
\newcommand{\mb}{\mathbf{B}}
\newcommand{\mai}[1]{\mathbf{A}^{(#1)}}
\newcommand{\mbi}[1]{\mathbf{B}^{(#1)}}

\newcommand{\mbni}[1]{\mathbf{B}^{(#1), \nu}}

\newcommand{\mfg}{\mathcal{G}}
\newcommand{\mfa}{\mathcal{A}}
\newcommand{\mfn}{\mathcal{N}}

\newcommand{\mbti}[1]{\mathbb{T}^{(#1)}}
\newcommand{\mbui}[1]{\mathbb{U}^{(#1)}}
\newcommand{\mbsi}[1]{\mathbb{S}^{(#1)}}
\newcommand{\mbri}[1]{\mathbb{R}^{(#1)}}
\newcommand{\mbgi}[1]{\mathbb{G}^{(#1)}}
\newcommand{\mbqi}[1]{\mathbb{Q}^{(#1)}}
\newcommand{\mbqsi}[1]{\mathbb{Q}^{(#1), *}}
\newcommand{\mbpsi}[1]{\mathbb{P}^{(#1), *}}
\newcommand{\mbpi}[1]{\mathbb{P}^{(#1)}}
\newcommand{\mbii}[1]{\mathbb{I}^{(#1)}}
\newcommand{\mbji}[1]{\mathbb{J}^{(#1)}}
\newcommand{\mbqwi}[1]{\mathbb{\widehat{Q}}^{(#1)}}
\newcommand{\mbhapi}[1]{\mathbb{\widehat{P}}^{(#1)}}
\newcommand{\bai}[1]{\boldsymbol{\alpha}^{(#1)}}
\newcommand{\basi}[1]{\boldsymbol{\alpha}^{(#1), *}}
\newcommand{\bbi}[1]{\boldsymbol{\beta}^{(#1)}}
\newcommand{\bbsi}[1]{\boldsymbol{\beta}^{(#1),*}}
\newcommand{\bzi}[1]{\boldsymbol{\zeta}^{(#1)}}
\newcommand{\bei}[1]{\boldsymbol{\eta}^{(#1)}}

\newcommand{\Ti}[1]{\Theta^{(#1)}}

\newcommand{\mbgni}[1]{\mathbb{G}^{(#1), \nu}}
\newcommand{\mbqni}[1]{\mathbb{Q}^{(#1), \nu}}

\newcommand{\mbpsni}[1]{\mathbb{P}^{(#1), \nu, *}}
\newcommand{\mbpni}[1]{\mathbb{P}^{(#1), \nu}}
\newcommand{\mbqwni}[1]{\mathbb{\widehat{Q}}^{(#1), \nu}}

\newcommand{\mbsni}[1]{\mathbb{S}^{(#1), \nu}}
\newcommand{\mbrni}[1]{\mathbb{R}^{(#1), \nu}}
\newcommand{\mbini}[1]{\mathbb{I}^{(#1), \nu}}

\newcommand{\mbtni}[1]{\mathbb{T}^{(#1), \nu}}

\newcommand{\mbtn}{\mathbb{T}^{\nu}}
\newcommand{\mbrn}{\mathbb{R}^{\nu}}
\newcommand{\mbpn}{\mathbb{P}^{\nu}}
\newcommand{\ben}{\mathbb{\eta}^{\nu}}
\newcommand{\ban}{\mathbb{\alpha}^{\nu}}

\newcommand{\bani}[1]{\boldsymbol{\alpha}^{(#1), \nu}}
\newcommand{\basni}[1]{\boldsymbol{\alpha}^{(#1),  \nu, *}}
\newcommand{\bbni}[1]{\boldsymbol{\beta}^{(#1), \nu}}
\newcommand{\Tni}[1]{\Theta^{(#1), \nu}}
\newcommand{\bzni}[1]{\boldsymbol{\zeta}^{(#1), \nu}}
\newcommand{\beni}[1]{\boldsymbol{\eta}^{(#1), \nu}}

\newarrow{corr}{<}{-}{-}{-}{>}
\newarrow{nor}{}{-}{-}{-}{>}
\newarrow{dash}{}{dash}{}{dash}{}
\newarrow{DashTo}{}{dash}{}{dash}{>}
\newarrow{coto}{|}{-}{-}{-}{>}

\DeclareMathOperator{\Irr}{Irr}
\DeclareMathOperator{\Lin}{Lin}

\DeclareMathOperator{\Hom}{Hom}
\DeclareMathOperator{\Ker}{Ker}

 \provideboolean{showlabels}
 \setboolean{showlabels}{false}
\newcommand{\mylabel}[1]{\ifthenelse{\boolean{showlabels}}{{\tt{[{#1}]}}\label
{#1}}{\label{#1}}}

\begin{document}

\title{ NORMAL SUBGROUPS OF ODD ORDER MONOMIAL  $p^aq^b$-GROUPS  }
     
\author{Maria I. Loukaki}

\degreeyear{2001}

%\department{Mathematics}
%\schools{B.A., University of Crete, 1991\\
 %        M.S., University of Crete, 1993}
          %M.S., University of Illinois, 2000}
 %\phdthesis
 %\degreeyear{1999}
 \maketitle

 \frontmatter

\chapter*{Abstract}
A finite group $G$ is called monomial if every irreducible character 
of $G$ is induced from a linear character  of some subgroup of $G$. 
One of the main questions regarding monomial groups is  whether  or not 
a normal subgroup $N$ of a monomial group $G$   is itself monomial. 
In the case that $G$ is a group of even order,  it has been proved (Dade, 
van der Waall) that $N$ need not be monomial. 
Here we show that, if $G$ is a monomial  group 
of order $p^aq^b$, where  $p$ and 
$q$ are distinct odd primes, 
then any normal subgroup $N$  of $G$ is also monomial.

%% Create a dedication with no heading, centered vertically
%% on the page.
 \newpage
\leavevmode\vfill
 \begin{center} 
to my parents Anna and Giannis, 
\end{center}

\begin{center} 
and 
\end{center}

 \begin{center} 
to my teacher, Everett  C. Dade, 

without whom    none of this would have been done,
 
and everything  would have been written  faster.
\end{center}
 \vfill

\chapter*{Acknowledgments}
First and foremost I would like to thank my advisor  Prof.  Everett C. Dade
not only for his endless patience, his continued support and encouragement,    
 his creativity and his humor, but  also for all those afternoon meetings in his 
office, where I saw how mathematics can become pure art.  
Thank you. This would have been only a dream without your help.

I would also like to thank the members of my committee that went through the 
trouble of reading this thesis, and Prof. Marty Isaacs for his sincere
 interest in  this work. I also thank my first math teacher Manolis Mpelivanis, 
Prof. George Akrivis  for his faith  in me,  and Prof. Giannis Antoniadis 
for introducing me to Representation Theory.

Special thanks go to 
my parents Anna and Giannis and my sister Marianthi 
for always being there, and to  Michalis for the verses  at Daily Grind and 
 the  ``sweet and  sour''  evenings.

I'm grateful  to  all my friends that 
made my stay in Urbana memorable,  especially Sasa ($\gamma \iota \alpha$ $o \lambda \alpha$),
George  ($\gamma \iota \alpha$  $\kappa \alpha\tau \iota$  ``$\epsilon \nu \delta
\iota \alpha \phi \epsilon \rho o \nu$ $\kappa \alpha \iota$ 
$\delta \iota \alpha \sigma 
\kappa \epsilon \delta \alpha \sigma \tau \iota \kappa o $''), 
  Michalino ($\gamma \iota \alpha$ $\tau \iota 
s$ $\kappa o \upsilon \beta  \epsilon \nu \tau \epsilon s$  $\sigma 
\tau o $ $\kappa \alpha \phi \epsilon \nu \epsilon \iota o$),  Daniella (for Mr. Maxx), 
Leonida ($\gamma \iota \alpha$ $\tau \alpha  $  $\mu \alpha \theta \eta \mu
\alpha \tau \alpha $ $\chi o \rho o \upsilon$), 
Tina ($\gamma \iota \alpha$ $\tau o$ $\tau \epsilon \nu \nu \iota s$), 
 Niko ($\gamma \iota \alpha$ $\tau \alpha$  
$\beta \rho \alpha \delta \iota \alpha$ $\mu \epsilon $ 
 $\pi \alpha \sigma \tau \iota \tau \sigma \iota o $),
 Edith (for the hours in the library), Nader (for the ``Jazz''), 
 and   Anne  who  kept her   confidence in me,   long after I had lost mine.

Finally I wish to acknowledge the University of  Illinois for the financial support 
 (Fellowships),  and the 
  National Science Foundation that partially  supported  the writing of this thesis 
by grants DMS 96-00106 and DMS 99-70030. 
Of course, the correctness of the results in this thesis is the responsibility 
of the author and not of the National Science Foundation.

%% The thesis format requires these lists to come in the following
%% order:
%%
%% Table of Contents
%% List of Tables
%% List of Figures or List of Illustrations
%% List of Symbols and/or Abbreviations

\tableofcontents
%\listoftables
%\listoffigures

%% Create a List of Abbreviations.

%\chapter*{List of Abbreviations}
%\begin{description}
%\item[CA] Caffeine Addict.
%\item[CD] Coffee Drinker.
%\end{description}

\mainmatter 

 \chapter{ Introduction}
\mylabel{intro}

\section{The problem}
Around 1930,  Taketa (see  \cite{ta})  introduced the notion of a
 monomial finite  group, i.e., a group 
for which every irreducible character is induced from a linear character 
of some subgroup. He also proved that any monomial group is solvable.
As  any supersolvable  group is monomial, Taketa's theorem 
places the class of  monomial groups  between the classes  of solvable 
and of supersolvable  groups. 
 
In 1967, Dornhoff (in \cite{do}) proved that every normal Hall subgroup    
of a monomial group is monomial. Furthermore, he 
 asked whether or not every normal subgroup of a monomial
 group is monomial, a question that arises  very naturally
  after his result on 
normal Hall subgroups.   A negative  answer to Dornhoff's 
question was found   independently by both  Dade \cite{da1} and van der Waall
\cite{wa} in 1973.
But  their counterexamples have even order. Furthermore,
 the prime 2 plays an important role in their construction, so that 
the examples can't be modified to give an answer to Dornhoff's question 
in the case of an odd order group. Thus Dornhoff's question remains open 
in the case of odd order monomial groups.
In the 1980's  Dade and Isaacs,  among others, tried to solve the
 remaining part of the  problem. 
They produced   many
beautiful theorems  suggesting that   the following conjecture might be true:
\begin{conj}\mylabel{conj1}
Let  $G$ be a monomial group of odd order.  Assume 
that $N$ is a normal subgroup of $G$. Then $N$ is also monomial.
\end{conj}
Among their results of that period,  the following are 
 the most useful  for this thesis.
\begin{theorem}[Isaacs] {\rm[Problem 6.11 in \cite{is}] }\mylabel{is1} 
Assume that  $G$ is 
 a monomial group and that  $A$ is any normal
 subgroup of $G$. Let $\lambda$ be a  linear character 
of $A$. We write   $G(\lambda)$ for  the stabilizer of $\lambda$ in $G$.
Assume that    $\chi$ is  any irreducible character of 
$G$ that lies over $\lambda$. Let   $\chi_{\lambda}$ be  its
 Clifford correspondent, i.e., the unique irreducible character
of $G(\lambda)$ that lies above $\lambda$ and induces $\chi$ in $G$.
 Then $\chi_{\lambda}$ is monomial.
\end{theorem}

%\begin{theorem}
%[Isaacs]   Elementary stabilizer limits....
%\end{theorem}

\begin{theorem}[Dade]{\rm [Theorem 3.2 in \cite{da}]} \mylabel{theo1}
 Let $G$ be a $p$-solvable group for some odd prime $p$. Assume that 
$U$ is a $\ZZ_p(G)$-module that affords a nondegenerate alternating 
 $G$-invariant  $\ZZ_p$-bilinear form. Let $H$ be a subgroup of $G$ 
that has $p$-power index in $G$. Then $U$ is $\ZZ_p(G)$-hyperbolic
if and only if  $U$ is $\ZZ_p(H)$-hyperbolic.
\end{theorem}
(For the definition of  hyperbolic modules, see Section \ref{nota}.)

Based on this theorem Dade was able to prove 
\begin{theorem}[Dade]{\rm [Theorem (0) in \cite{da}]}\mylabel{theo2}
 Let $G$ be a $p$-solvable group for some odd prime $p$,  
and let $N$ be a 
normal subgroup of $G$. Assume that  $\chi\in \Irr(G)$ is a
 monomial character of $G$ that has $p$-power degree. Then every 
irreducible constituent of the restriction $\chi_N$ of $\chi$ to $N$
is monomial.
\end{theorem}

Dade's Theorem (Theorem \ref{theo2}) is a very powerful result as it
 manages to handle monomial characters individually. That is, he doesn't 
 assume that $G$ is monomial or that there are any monomial
 characters in $\Irr(G)$ except the specific prime-power degree character 
that is analyzed.
Unfortunately, if we have characters with degrees that are divisible  by more
 than one prime, then we can't hope to prove something similar to Theorem 
\ref{theo2}, as Dade has given counterexamples (in \cite{da})
where his theorem fails when two odd primes divide $\chi(1)$.

But if we have enough monomial characters so that we  have the freedom
to replace, in some sense, the ``bad''ones with ``good'' ones, 
we can do more. We can actually prove that a big enough
section inside any  normal subgroup $N$ of a monomial  $p^a q^b$-group $G$ 
is nilpotent, provided that $p$ and $q$  are odd primes.
What we show   is 
\begin{main}\mylabel{main}
Let $G$ be a finite $p^a q^b$-monomial  group,
 for some odd primes $p$ and $q$.
Assume that $N$ is a normal subgroup of $G$ and that 
 $\chi$ is an irreducible character of $N$. 
Then there exists a faithful linear limit   $\chi^*$  of
 $\chi$, such that  the domain $Dom(\chi^*)$ of $\chi^*$ 
 is a nilpotent group.
\end{main}
(For the definition of  ``faithful linear limits'',  see Section \ref{lim1}.)

As corollaries  of the  Main Theorem \ref{main}  we have 
\begin{main}\mylabel{main2}
Let  $G$ be  an odd order monomial $p^a q^b$-group. If $N$ is a normal 
subgroup of $G$, then $N$ is monomial.
\end{main}

and 

\begin{main}\mylabel{main3}
Let $G$ be an odd order monomial $p^a q^b$-group and let $\chi \in \Irr(G)$.
Then there exists a faithful linear limit  $\chi^*$  of $\chi$ such that 
$\chi^*(1)=1$, i.e., $\chi^*$  is a linear character. 
\end{main}

What is the desired  property a ``good'' character   has?  We can control its 
degree. If this degree has the right properties 
then nilpotent subgroups appear,  as the following result shows:
\begin{theorem}\mylabel{dd}
Assume that $G$ is a $p, q$-group, where $p$ and $q$ are distinct odd primes, 
 and that $N, M$ are  normal subgroups of $G$.
Let $M=P \times S$ and  $N = P \rtimes Q$, where $P$ is a $p$-group, 
and $S,  Q $ are $q$-groups with $S \leq Q$.
  Assume that the center  $Z(P)$ of $P$ is  maximal among the  abelian 
$G$-invariant subgroups of $P$.
 Let $\chi, \alpha,\beta $ and $  \zeta$  be irreducible characters of
 $G, P, S$ and 
$Z(P)$ respectively that satisfy 
\begin{align*}
\chi \in \Irr(G | \alpha \times \beta)
 &\text{ and } \alpha \in \Irr(P| \zeta), \\ 
 \zeta \text{ is a faithful } &G\text{-invariant character of  }Z(P), \\
G(\beta ) &= G,  \\
\chi \text{ is a monomial } &\text{character of $G$   with } \chi(1)_{q} = \beta(1),
\end{align*}
where $\chi(1)_{q}$ denotes the $q$-part of the integer $\chi(1)$.
Then $Q$  centralizes $P$.
\end{theorem}
The above theorem, that appears as Theorem \ref{m.t2} in this thesis, 
is heavily based on Theorem \ref{theo1}, and the work done  in Dade's paper 
\cite{da}.

The way we use  in this thesis to  approach  these ``good'' characters 
is by constructing a character of known degree in a subgroup of $G$ 
that extends to its own stabilizer in $G$.  So the key tool for the proof of 
Theorem \ref{main} is the following  result (that appears as 
   Theorem \ref{cc:co2} in   Chapter \ref{cc}):
\begin{theorem}\mylabel{key step}
Let $P$ be a $p$-subgroup, for some prime $p$, of a finite odd order group $G$.
Let $Q_1,  Q$ be   $q$-subgroups of  $G$, for some prime 
$q\ne p$, with $Q_1 \leq Q$.
Assume that $P$  normalizes $Q_1$, 
while $Q$ normalizes the product  $P \cdot Q_1$.
Assume further that  $\beta_1$ is an irreducible character of $Q_1$.
 Then there exist  irreducible characters $\bn_1$ of $Q_1$ 
and $\bn$ of $Q(\bn_1)$ such that
\begin{align*}
P(\beta_1 ) &= P(\bn_1), \\
Q(\beta_1) &\leq Q(\bn_1) \, \text{and}\\
\bn |_{Q_1} &= \bn_1.
\end{align*}
Therefore $\bn$ is an  extension of $\bn_1$  to $Q(\bn_1)$.
\end{theorem}
A possible generalization of Theorem \ref{key step} to  the case where $Q_1$
and $Q$ are arbitrary $p'$-subgroups of $G$ 
would give  a generalization of the Main Theorem \ref{main}.
In this way the original problem (Conjecture \ref{conj1}) can be  transformed
to the following question:
\begin{quest}
Does Theorem \ref{key step} hold
if its second sentence is replaced by 
``Let $Q_1, Q$ be $p'$-subgroups of $G$''?  
\end{quest}

 {\bf Note: }  
Shortly after this thesis was submitted  E. C. Dade
 found a counterexample 
to the above question. His example is not a counterexample 
for  Conjecture \ref{conj1}. Nevertheless, it suggests strongly  that 
 a generalization of Theorem \ref{key step} requires some new ideas.

\newpage
\section{Notation}\mylabel{nota}
Let  $G$ be a finite group  with  
 $M ,  N, K,  K_1, \dots, K_t$ subgroups of $G$, 
for some $t \geq 1$, 
and $g,h$ elements of $G$.  
 Assume further that   $M$  normalizes $N$.
Let $\chi$ be  a (complex) irreducible
 character of $M$, (we assume that all the characters we use
in this text  are over the field $\mathbb{C}$ of the complex numbers) 
and $\mathfrak{X}$ an  irreducible $\mathbb{C}$-representation 
of $G$ that affords $\chi$. Assume further that $\phi_1, \dots, \phi_t $ are 
 irreducible characters of $K_1, \dots, K_t$, respectively.
Let $p$ be a prime number,  and $\pi$ a set of primes.

The  list  that follows describes the notation we will be using 
for the rest of this thesis.
\begin{align*}
&\ZZ_p:          &\quad  &\text{the field of $p$ elements $\ZZ/p\ZZ$}\\
&Z(G):           &\quad  &\text{the center of $G$}\\
&\Syl_p(G):      &\quad  &\text{the set of all Sylow $p$-subgroups of $G$}\\
&\Hall_{\pi}(G): &\quad  &\text{the set of all Hall $\pi$-subgroups of $G$}\\
&O_p(G):        &\quad  &\text{the largest normal $p$-subgroup of $G$}\\
&\Phi(G):        &\quad  &\text{the Frattini subgroup of $G$ }\\
&g^h:            &\quad  &\text{the $h$-conjugate $h^{-1}g h$ of $g$}\\
&K^h:            &\quad  &\text{the $h$-conjugate group $h^{-1} K h$ of $K$}\\
&N \leq G:  &\quad &\text{$N$ is a subgroup of $G$ }\\
&N \unlhd G:  &\quad &\text{$N$ is a normal subgroup of $G$ }\\
&N(M \tin G)=N_G(M):     &\quad  &\text{the normalizer of $M$  in $G$}\\
&N(K \tin M )=N_M(K):     &\quad  &\text{the normalizer of $K$ in $M$ }\\
&N(K_1, \dots, K_t  \tin M):        & &\text{the normalizer of  all $K_i$, for $i=1, \dots, t$
 in $M$, i.e., }\\ 
&                             &&\text{ the intersection           $\bigcap_{i=1}^t N(K_i  \tin M) $  }\\
&C(M \tin G)=C_G(M):     &\quad  &\text{the centralizer of $M$ in $G$ }\\
&C(K \tin M)=C_M(K):     &\quad  &\text{the centralizer of $K$ in $M$}\\
%\end{align*}
%\begin{align*}
&C(K_1, \dots, K_t  \tin M):        & &\text{the centralizer of  all $K_i$, for $i=1, \dots, t$
 in $M$, i.e., }\\ 
&                             &&\text{ the intersection           $\bigcap_{i=1}^t C(K_i  \tin M) $  }\\
&M\ltimes N:     &\quad  &\text{the semidirect product of $M$ and  $N$
                                  when $M$ acts on $N$}\\
&K\rtimes N:     &\quad  &\text{the semidirect product of $K$ and $N$
                                   when $N$ acts on $K$}\\ 
&[M,N]:          &\quad  &\text{the commutator subgroup of $M,N$}\\
&\Irr(G):        &\quad  &\text{the set of all complex  irreducible 
                              characters of $G$ }\\
&\Lin(G):        &\quad  &\text{the set of all linear complex characters
                                    of $G$ }\\
&\Ker(\chi):     &\quad  &\text{the Kernel of $\chi$, i.e.,
                              $\Ker(\chi)=\{m \in M|\chi(m) = \chi(1)\}$}\\
&\chi^g:         &\quad  &\text{the $g$-conjugate of $\chi$, i.e.,
                                 $\chi^g$ is a character of $M^g$ defined as}\\
&                &       &\text{ $\chi^g(m^g):= \chi(m)$ for all $m \in M$}\\
&\chi^G: & &\text{the induced character on $G$} \\
\end{align*}
\begin{align*}
&[\chi, \psi]: & &\text{the inner product of $\chi, \psi \in \Irr(M)$, i.e. } \\
      &                    &     &\text{ $[\chi, \psi]= (1 /  |G|) \sum_{g \in G}\chi(g)\psi(g^{-1})$ }\\
&\Irr^M(N):      &\quad  &\text{the set of all $M$-invariant irreducible 
                            characters of $N$}\\
&\Irr^M_N(G): &\quad &\text{the set of all $\chi \in \Irr(G)$ such that $\chi$ lies above }\\
 &                    &          &\text{   at least one character $\theta \in \Irr^M(N)$ }\\
&\Irr(G|\chi):   &\quad  &\text{the set of all irreducible characters of $G$
                           that lie above $\chi$ }\\
&G(\chi):        &\quad  &\text{the stabilizer of $\chi$ in $G$, i.e., 
                         the   set of all  elements $g$ of $G$}\\ 
&                &        &\text{         that satisfy $\chi^g = \chi$}\\
&G(\phi_1, \dots,  \phi_t): &  & \text{the stabilizer of $\phi_i$, for $i=1, \dots, t$, in  
$G$, i.e.,                            the set $\bigcap_{i=1}^tG(\phi_i)$ }\\
&K(\chi):        &\quad  &\text{the stabilizer of $\chi$ in $K$}\\
&K(\phi_1, \dots,  \phi_t): &  & \text{the stabilizer of $\phi_i$, 
for $i=1, \dots, t$, in  $K$, i.e., 
                           the set $\bigcap_{i=1}^tK(\phi_i)$ }\\
&\det(\mathfrak{X}(g)):  &\quad  &\text{the determinant of the matrix 
                                 $\mathfrak{X}(g)$ for some $g\in G$}\\
&\det(\chi):      &\quad  &\text{the linear character of $G$ defined as}\\ 
&                &       &\text{   $(\det(\chi))(g) = \det(\mathfrak{X}(g))$
                                  for all  $g \in G$}\\
&o(\chi):        &\quad  &\text{the determinantal order of $\chi$, i.e., 
                               the order of the linear character}\\
&                &       &\text{ $\det(\chi)$  as an element of the group $\Lin(G)$}\\ 
&[r]:   &       &\text{  the integral part of a real number $r$, i.e., the
 largest integer $t$     such that $t \leq r$.  }               
\end{align*}

Assume that $G_0 \unlhd G_1 \unlhd \dots \unlhd G_n$  is  a series of normal subgroups
of $G$, for some integer $n \geq 0$.
  Assume further that $\chi_i$ is an  irreducible character of $G_i$, for 
all $i=0, 1, \dots, n$,  such that  $\chi_i \in \Irr(G_i)$ lies above $\chi_{i-1} \in \Irr(G_{i-1})$, 
whenever $i =1, \dots, n$.  Then we call  the set  $\{ \chi_i \}_{i=0}^n$ 
{ \em a character tower }  for the series  $\{ G_i \}_{i=0}^n$.

If $T$ is a finite--dimensional $\ZZ_q(G)$-module
that affords a $G$-invariant  symplectic form $<\cdot , \cdot>$, then we will
use the terminology  introduced in \cite{da}.
So,  if $S$ is a $\ZZ_q(G)$-submodule of $T$, then 
$S^{\perp}$ is the { \em perpendicular  }  subspace 
of $S$, i.e.,  $S^{\perp}:= \{t \in T |<S , t> = 0\}$. Furthermore, $S$ is { \em isotropic }  
if $S \leq S^{\perp}$.
$T$ is  { \em  anisotropic } if it contains no non--trivial isotropic 
$\ZZ_q(G)$-submodules, and is { \em hyperbolic } if it contains some
 self--perpendicular $\ZZ_q(G)$-submodule $S$, i.e., $S$ is a
 $\ZZ_q(G)$-submodule satisfying  $S = S^{\perp}$.

If $U_i$, for $i=1,\dots,n$, are  $F(G)$-submodules of
an $F(G)$-module $S$, where $F$ is a field and $n$ a positive  integer,
we write
\begin{align*}
&U_1 \dotplus U_2:   &\quad   &\text{for the internal direct
                                           sum of $U_1$ and  $U_2$, and }\\
&\sum_{1\leq i \leq n}^{\cdot}U_i: &\quad  &\text{for the internal direct sum of the 
                                   $U_i$,      for all $i=1,\dots,n$.}\\
\end{align*}

\newpage
\section{The general ideas in  the proof}
As this thesis turned out to be much longer and  
complicated than expected, the author would like to apologize to the reader 
for all the  mysterious groups that appear suddenly,  and  for no apparent reason, 
in the chapters that follow. We will attempt in this section  to give the main ideas
of the proof,   trying  to avoid, as much as possible, the technical parts.

Assume that $G$ is a finite monomial group. 
Assume further that $G$ has order $p^aq^b$ for 
some distinct odd primes   $p$ and $q$. Let  $N$ be a normal subgroup of $G$.
Fix an irreducible character  $\psi$  of $N$. 
If $A$ is any  normal  subgroup of $G$ contained in $N$, 
and $\lambda$ is  a linear character of $A$ lying under $\psi$,
then Isaacs observation (Theorem \ref{is1}), implies that we may pass from  $G$ to
the stabilizer $G' = G(\lambda)$  of $\lambda$ in $G$ without losing the
 monomiality of those irreducible characters of $G'$ that lie above $\lambda$.
This way we reduce the order of the group $G$, possibly 
loosing  some of the monomial
 characters,   but still keeping track of those 
that lie above $\psi$. (This  reduction procedure, that is  described in
 Chapter \ref{lim}, produces a ``linear limit''.)
What we prove by induction is that, by applying this
 procedure many times,  and choosing the $A$ and $\lambda$ carefully,
(see Section \ref{lmt}  for an explanation of ``careful''),  we can reduce $G, N$ and 
$\psi$ to $G' , N'$ and $\psi '$, respectively that satisfy 
\begin{align}\mylabel{qq}
G' \leq G, \text{ and } N '&= N \cap G'   \leq N,   \text{ while  } \psi' \in \Irr(N'), \notag \\
N'/ \Ker(\psi') &\text{ is nilpotent,  } \notag \\
(\psi')^N &= \psi, \\ 
\text{ all $\chi'  \in \Irr(G' | \psi') $ } &\text{    are monomial. }  \notag 
\end{align}
This is equivalent to  our Main Theorem \ref{main}. 
It is proved using induction on the order of $N$.
Note that \eqref{qq} easily implies that the character $\psi$ is monomial, 
since $\psi$  is induced from $\psi'$, and  the 
unique character $\psi'/ \Ker(\psi')$  of the factor group $N'/ \Ker(\psi')$
that inflates to  $\psi'$ is monomial 
as an irreducible  character of a nilpotent group.

First notice that if $N$ is a nilpotent group, \eqref{qq} holds trivially with 
$N$ in the place of $N'$. (We don't even need to apply Clifford's  theorem
to  any normal subgroup $A$ and any linear character $\lambda \in \Irr(A)$.)
Now suppose that $N$ is not nilpotent.
Because $G$ is solvable, there exists a normal subgroup $L$ of $G$ such 
 that $ L  < N$,  while 
$N/ L$ is either a $p$-  or a $q$-group. 
We may assume that $N/ L$ is a $q$-group, and that  $p$ divides $|L|$. 
Let $\mu \in \Irr(L)$ be an irreducible  character 
of $L$ lying under $\psi \in \Irr(N)$. We apply the reductions described above 
with $L$ and $\mu$ in the place of $N$ and $\psi$, i.e., 
we reduce $G$ using   linear characters of normal 
subgroups of $G$ that are contained inside
 $L$ (again choosing these linear characters carefully).
Every time $G$ gets reduced $N$ and $L$ are    also reduced.
Furthermore, the  fixed character $\psi$ of $N$ also gets reduced, 
at every step,  to a Clifford correspondent. So does  the character $\mu \in \Irr(L)$.
According to the inductive hypothesis applied to  $L$, at some point  these 
reductions lead us to  groups   $G''$ and $ L''$,  and to an  
 irreducible character $\mu'' \in \Irr(L'')$ 
that induces $\mu$ while $L'' / \Ker(\mu'')$ is nilpotent. 
Of course 
if  we reduce  the group $G'' $  more,   using normal subgroups $A$ of 
 $G''$ contained in $L''$  and their linear characters $\lambda \in \Lin(A)$, the
 above properties remain valid.
So we continue reducing  until there is no 
normal subgroup $A \leq L'' $ and linear 
character $\lambda \in \Irr(A)$ lying under $\mu''$ with $G''(\lambda) \leq G''$.
 At the end of this procedure we reach 
groups  $L'$ and $ G' $,  and  a character $\mu'\in \Irr(L')$,  that satisfy the 
equivalent of \eqref{qq} for $L$.
 In addition to those,  the group $N$ gets reduced to $N' = G' \cap N$,
 and its irreducible character  $\psi$ is reduced to $\psi' \in \Irr(N')$.
Therefore we  have 
\begin{align}\mylabel{qqq} 
G' \leq G, \text{ and } L '= L  \cap G'   \leq L,   &\text{ while  } \mu' \in \Irr(L'), \notag \\
N' = N \cap G'  \leq N &\text{ and }  \psi' \in \Irr(N' | \mu'), \notag \\
L'/ \Ker(\mu') &\text{ is nilpotent, } \\
(\mu')^L &= \mu, \notag \\ 
\text{ all $\chi'  \in \Irr(G' | \mu') $ } &\text{    are monomial, }  \notag \\
\text{ for any $A \unlhd G'$ contained in $L'$,  and any   $\lambda \in \Lin(A)$ }
&\text{  that lies under $\mu'$, we get $G'(\lambda) = G'$. } \notag 
\end{align}
 It is easy to see (after all,  these reductions are just
repeated applications of Clifford theory) that $\psi'$  induces $\psi $,  and lies above 
$\mu'$. 

The fact $L'= L/ \Ker(\mu')$ is nilpotent implies that the factor groups 
$L/ \Ker(\mu')^g$ are all nilpotent, whenever $g \in G'$. 
Let  $K' = \bigcap_{g \in G'}\Ker(\mu')^g$ and $L_1 = L'/ K'$.  
Then $L_1$ is also nilpotent.
Therefore
 $L_1$  splits as the direct product 
$L_1 = P_1 \times Q_1$ of  its $p$- and $q$-Sylow subgroups  $P_1$ 
and $Q_1$, 
respectively. If $\mu_1$ is the unique irreducible character of the factor group 
 $L_1$ that inflates to $\mu' \in \Irr(L')$, then  $\mu_1$ also decomposes as 
$\mu_1 = \alpha \times \beta$, where $\alpha \in \Irr(P_1)$ and $\beta \in \Irr(Q_1)$.
At this point we can say more for $P_1$ and $Q_1$. 
If $A_1$  is any normal  subgroup of $G_1= G'/K'$ contained in $L_1$, 
  then $ A _1 = A'/ K'$, 
where $A'$ is a  normal  subgroup of $G'$. 
If, in addition,  $\lambda_1 \in \Lin(A_1)$ lies under $\mu_1$,  
then $\lambda_1$ 
inflates to a unique character $\lambda' \in \Lin(A')$ that lies under $\mu'$.
As $L'$ was as  reduced as possible,  
the character $\lambda'$ is $G'$-invariant. It also lies under $\mu'$. 
Hence $\Ker(\lambda') = \Ker(\mu'_{A'})$. So 
$$
\Ker(\lambda') = 
\bigcap_{g \in G'}\Ker(\lambda')^g  \leq \bigcap_{g \in G'}\Ker(\mu')^g=K'.
$$
 Therefore
 $\lambda_1$ is a faithful linear character of $A_1$,  and is $G_1$-invariant. 
Thus every characteristic abelian  subgroup 
of $P_1$  is cyclic and is contained in the center of $P_1$.
 Similarly for $Q_1$.
In particular the center $Z(P_1)$ of $P_1$ is cyclic, and  affords a 
faithful $G_1$-invariant  linear  
character that lies under $\alpha$. Similarly the center $Z(Q_1)$ 
is cyclic,  and affords a faithful $G_1$-invariant 
linear character lying under $\beta$. 
Furthermore, 
 $P_1$ (and similarly $Q_1$)  is of a very specific type.
 Either  it is  cyclic  or  it is   the central product of 
an extra special  $p$-group 
of exponent $p$ with  the  cyclic  $p$-group $Z(P_1)$.
Similarly for the $q$-group $Q_1$.
Then  clearly  the characters $\alpha$ and $\beta$ are 
$G_1$-invariant  and faithful.  Hence $\mu_1 = \alpha \times \beta$ is 
$G_1$ invariant. 
We conclude that the character $\mu'$ is $G'$-invariant.
 So $ \Ker(\mu')= K'$  is a normal 
subgroup of $ G'$. Thus $G_1 = G' / \Ker(\mu')$.
   Furthermore,  because  $\psi' \in \Irr(N')$  lies 
above $\mu'$, the group $K'$ 
is  contained in  the kernel 
of $\psi' \in \Irr(N')$. Therefore  there exists a unique irreducible character 
$\psi_1$ of the factor group $N_1 = N'/K'$ 
that inflates  to  $\psi' \in \Irr(N')$.   Note that 
$Z(P_1)$ is maximal among the abelian normal 
subgroups of $G_1$ contained in $P_1$. This makes the factor group 
 $P_1 / Z(P_1)$ naturally  a symplectic space
 (see \eqref{m.d1} for the definition of the symplectic form). It is  actually a 
symplectic 
$\ZZ_p(G_1/ L_1)$-module. 
A  picture of the situation is 
\begin{align}\mylabel{.0a}
&G_1                    &\qquad         & \notag \\
& \big|                         &       & \notag  \\
&N_1                    &       &\psi_1\\
&\big|                  &       &\big|  \notag \\
L_1= P_1 &\times Q_1    &   \mu_1 &= \alpha \times \beta \notag  
\end{align}
Note that every character of $G_1$ that lies above $\mu_1$ is still monomial. 
Furthermore the factor group $N_1 / L_1$ is a $q$-group.

We could continue the reductions with normal subgroups of $G'$ that lie 
inside $N'$. But it turns out that it is not needed. The group $N'$ and the character 
$\psi'$ already achieved when doing the reductions for $L$, 
satisfy \eqref{qq}.  Indeed  $K'= \Ker(\mu' ) \leq \Ker(\psi')$, while 
any character of $G'$ lying above 
$\psi'$ is monomial, because $\psi'$  lies above $\mu'$.
 Therefore if we could prove that $N_1 = N' / K'= N'/ \Ker(\mu')$ is nilpotent, 
then we would have that $N ' / \Ker(\psi')$ is also nilpotent. This would
 prove the 
inductive step for  \eqref{qq}.  Actually,  that is what  we prove.

Let  $Q$ be  a $q$-Sylow subgroup of $N_1$. 
Then $N_1$ is the semidirect product
$N_1 = P_1 \rtimes Q$. We prove that $Q$ centralizes $P_1$, using 
Theorem \ref{dd}.
Observe that the situation for the groups $G_1, N_1$  and $L_1$ 
looks similar to that described in Theorem \ref{dd}.
Indeed,   the center  $Z(P_1)$ of $P_1$ 
is  maximal among the $G_1$-invariant abelian subgroups of 
$P_1$, while the  characters $\alpha$ and $\beta$ are  $G_1$-invariant.
There is one ingredient missing from the hypothesis 
of the above theorem.  That is,   a monomial character
$\chi_1$  of $G_1$ lying above $\mu_1$ and satisfying $\chi_1(1)_q = \beta(1)$.
We do know that every character of $G$ above $\mu_1$ is monomial, 
but there is no reason for  one of those monomial characters to have 
the desired degree. And that is the obstacle.

If  the irreducible character $\beta$ of $Q_1$  extends to $G_1$,    then
the product  of this extended character with any $p$-special character of $G_1$  
lying above $\alpha$, is  an irreducible  character of $G_1$  that lies above $\mu_1$ 
and whose degree has the $q$-part equal to $\beta(1)$.
Of course there is no reason for the character $\beta \in \Irr(Q_1)$ to extend.

 A way to resolve this problem is to replace  the character $\beta$ 
with a new character $\bn$ that extends. Actually Theorem \ref{key step} 
offers a way to replace characters.  But now another problem appears. 
If we replace $\beta$ in $Q_1$ with a new character $\bn$ of the same group
 $Q_1$, then we can get an  irreducible  character of $G_1$  
lying above $\alpha \times \bn$ with the correct degree, but which may not 
(and most probably is not)  monomial. 
After all the only characters   we know to be monomial  in 
$G_1$ are those lying above $\mu_1$. 
Since the only source of monomial characters is back in $G$,  
we must change  the original character $\mu \in \Irr(L)$.

Suppose  we could find a $\mu^{\nu} \in \Irr(L)$ so that,   
when we reduce  $G$ as above, using Isaacs observation as much as
 possible for subgroups of $L$, 
we end up with a system 
\begin{align}
&G_1^{\nu}                      &\qquad         & \notag \\
& \big|                         &       &  \notag \\
&N_1^{\nu}                      &       &\psi_1^{\nu}\\
&\big|                  &       &\big|  \notag \\
L_1^{\nu}= P_1^{\nu} &\times Q_1^{\nu}  &   \mu_1^{\nu} &= 
\alpha_1^{\nu} \times \beta_1^{\nu} \notag 
\end{align}
having same  properties as \eqref{.0a}, plus
\begin{itemize}
\item{}  the character $\bn$ extends to a $q$-Sylow subgroup $\maq^{\nu}$ 
of $G_1^{\nu}$, and thus $\mu_1^{\nu}$ extends  to \
$ \maq^{\nu} \cdot L_1^{\nu}$, \\
\item{} the symplectic space $P_1/ Z(P_1)$ is isomorphic 
to $P_1^{\nu} / Z(P_1^{\nu})$, and \\
\item{} this isomorphism carries
the commutator $[ P_1/Z(P_1) , Q]$   into 
$[ P_1^{\nu}/Z(P_1^{\nu}) , Q^{\nu}]$, where $Q^{\nu}$ is a 
$q$-Sylow subgroup of $N_1^{\nu}$. 
\end{itemize}
 The fact that $\bn$ extends implies, as we saw, that 
the $q$-Sylow subgroup $Q^{\nu}$ of $N_1^{\nu}$  centralizes $P_1^{\nu}$. 
Hence $[ P_1^{\nu}/Z(P_1^{\nu}) , Q^{\nu}]= 1$. 
We conclude that the group $[ P_1/Z(P_1) , Q]$,   which is isomorphic to 
a subgroup of $[ P_1^{\nu}/Z(P_1^{\nu}) , Q^{\nu}]$,   is also trivial.
Hence $Q$ centralizes $P_1$. So the inductive step would be complete, 
provided that we could find such a character $\mu^{\nu}$.

The ``miracle''  is that such a character $\mu^{\nu}$ does exist.
To find it,  we had to introduce the notion of triangular 
sets and go through all the
 complicated machinery described in Chapters  \ref{pq},  \ref{pq:sec5} 
and \ref{n}. The reason is that 
  we have a replacement theorem (Theorem \ref{key step}) 
that works under special hypotheses. The character $\mu$ doesn't satisfy these
hypotheses. 
The character $\beta$ does satisfy  them,  but, for the reasons explained above, 
if we change $\beta$   we don't get  monomial characters. 
So instead of $\mu$ we change a character $\mu^*$ 
 corresponding  to $\mu$  in a certain  subgroup $L^*$  of $L$. 
Both  $\mu^*$ and $L^*$  satisfy the hypotheses of Theorem \ref{key step}, 
so that we can replace $\mu^*$ with another  $\mu^{*, \nu} \in \Irr(L^*)$ 
that extends to a $q$-Sylow subgroup of $G(\mu^*) \cdot L^*$. 
In addition, $\mu^*$ and $L^*$  are picked in such a way that we can get back 
from  the new character $\mu^{*, \nu}$  of $L^*$ to a corresponding new character 
$\mu^{\nu}$ of $L$.

In Chapters \ref{elm} and  \ref{lim}  we show that this new character retains 
its extendibility properties throughout   our reductions.

In Chapter  \ref{m} we put all the pieces together to prove the Main Theorem 
\ref{main}.

%%% Local Variables:            
%%% TeX-master: "lemma1" go back  fro 

%%% End: 

\chapter{ Preliminaries}

\section{Glauberman  correspondence}
Let $A$ and $G$ be two finite groups with orders that are relatively prime.
In the case that $A$ is a solvable group, Glauberman  \cite{gl}
constructed a    `natural' bijection between the set 
$\Irr^A(G)$ of  $A$-invariant irreducible characters of $G$,
 and the set $\Irr(C(A \tin G))$  of irreducible  characters of the 
fixed points $C(A \tin G)$ of $A$ in $G$.
Two well known facts (see Chapter 13 (page 299) in \cite{is})
about the Glauberman correspondence are 
\begin{lemma}\mylabel{cc:l2}
Suppose that $A,G$ are finite groups such that $(|A|, |G|) = 1$ and that $A$
acts on $G$. 
Let a group $S$ act on the semidirect product  $AG$,  leaving
both $A$ and $G$ invariant. If $\chi$ is an $A$-invariant irreducible
character of $G$, 
and $\chi^*$ is its $A$-Glauberman correspondent in 
$\Irr(C_G(A))$, then for any element $s\in S$, the Glauberman correspondent
$(\chi^s)^*$ of $\chi^s$ equals $(\chi^*)^s$.   
\end{lemma}

Lemma \ref{cc:l2} obviously implies 
\begin{corollary}\label{cc:co1}
In the situation of Lemma \ref{cc:l2}, let $T$ be a subgroup of $S$.
Then  $\chi$ is fixed by $T$ if and only if  its Glauberman 
 correspondent 
$\chi^*$ is also fixed by $T$.
\end{corollary}

\section{Groups}
\begin{proposition}\mylabel{prel:p1}
Assume $Q, P$ are two finite groups with coprime orders,  and  that $Q$ acts 
as  automorphisms of $P$. Assume further that  
$P$ is  the product $P = P_1 \cdot P_2$, of its normal subgroup $P_1$ and 
some subgroup $P_2$ of $P$. If both $P_1$ and $P_2$ are $Q$-invariant, then 
$$N(Q \tin P) = C(Q \tin P) = C(Q \tin P_1) \cdot C(Q \tin P_2).$$ 
\end{proposition}

\begin{proof}
As $Q$ and $P$ have coprime orders, and $Q$ acts on $P$, we obviously have that
$$N(Q\tin P) = C(Q \tin P) \geq C(Q \tin P_1) \cdot C(Q \tin P_2).$$   
Hence it remains to show that $ C(Q \tin P) \leq C(Q \tin P_1) \cdot C(Q
 \tin P_2)$.

As $P_1 \unlhd P$, we have $C(Q \tin P)/C(Q \tin P_1) \cong C(Q \tin P/P_1)$, 
by Glauberman's Lemma, 13.8 in \cite{is}. 
The natural isomorphism of $P/ P_1= (P_1 P_2)/P_1$ onto
 $P_2/ (P_1 \cap P_2)$
 preserves the action of $Q$. So it sends $C(Q \tin P/P_1)=
 C(Q \tin P)\cdot P_1/ P_1$ onto $C(Q \tin P_2/(P_1 \cap P_2)) = 
C(Q \tin P_2)\cdot (P_1 \cap P_2)/ (P_1 \cap P_2)$. So  
$C(Q \tin P_2)$ covers $C(Q \tin P) $ modulo  $C(Q \tin P) \cap P_1= 
C(Q \tin P_1)$.
 We conclude that
 $C(Q \tin P) = C(Q \tin P_2) \cdot C(Q \tin P_1)$.  
Hence the  proposition holds.
\end{proof}

As an easy consequence of Proposition \ref{prel:p1} we have
\begin{corollary} \mylabel{prel:co1}
Assume $Q, P$ are two finite groups with coprime orders,  and that $Q$ acts as 
 automorphisms of $P$. For every $i=1, \dots, n$, let  $P_i$ be a 
 $Q$-invariant sugroup of $P$. Assume further that 
 $P_j$ normalizes $P_i$, whenever $1\leq i \leq j \leq n$, while 
$P$ is the product $P =
P_1 \cdot P_2  \dots P_n$.  Then 
$$N(Q \tin P) = C(Q \tin P) = C(Q \tin P_1) \cdot C(Q \tin P_2) \cdot \dots
\cdot C(Q \tin P_n).$$ 
\end{corollary}

The following theorem is a multiple application  of Clifford's Theorem:
\begin{theorem}\mylabel{prel:t1}
Let $G$ be a finite group. Assume that 
$1 \unlhd G_1 \unlhd \dots \unlhd G_{m} \unlhd  G$ is a 
series of normal subgroups  of $G$.
 Assume further that we have fixed a character tower $\{ \chi_i\}_{i=1}^m $ 
for the above series, i.e., 
$\chi_i \in \Irr(G_i)$ for $i=1, \dots, m$, such that  $\chi_i$ 
 lies above $\chi_{i-1} $,  whenever $i=2, \dots, m$. 
 For any $i=1, \dots ,m$ we  write
$G_i ^{m} = G_i(\chi  _1, \chi  _2 , \dots ,\chi  _{m})$
for the stabilizer of $\chi  _1, \chi  _2 , \dots ,\chi  _{m}$ 
in $G_i$,  and  $G^{m} = G(\chi_1,\dots,\chi_{m})$
for the corresponding stabilizer in $G$.
Then $ G_1^{m} = G_1 $ and $G_k^{m}
 =G_i^{m}\cap G_k  = G^{m} \cap G_k \unlhd G_i^{m}$,  whenever
 $1\leq k\leq i\leq m$.
Furthermore, there exist  unique characters $\chi _{i}^{m}$, for 
$i=1, \dots ,m$, such that 
\begin{multline}\mylabel{pr01}
\text{
 $\chi_1^{m}= \chi_1$, while 
   $\chi_{i}^{m}  \in \Irr(G_{i}^m)$  lies over $\chi_{1}^{m}, \dots 
,\chi_{i-1}^{m}$ }\\
\text{ and induces  $\chi  _i \in \Irr(G_i)$, for all $i=1, \dots, m$. }
\end{multline}
Furthermore, these characters satisfy
\begin{itemize} 
\item[(1)]$G^{m} = G( \chi  _{1}^{m}, \chi  _{2}^{m},  \dots ,
\chi _{m}^{m})$,  and 
\item[(2)] If   $G_m \leq H \leq G$ and  $H^m= G^m \cap H$, then 
 for any $\phi \in \Irr(H |\chi_m)$ and 
$\chi \in \Irr(G| \phi)$, 
there exist unique   characters $ \phi^m \in \Irr(H^m|\chi^m_m )$ and 
  $\chi^{m} \in \Irr(G^{m} |\phi^m)$ that induce $\phi^m$ and  $\chi$, respectively.
Conversely, if $ \phi^m \in \Irr(H^m|\chi^m_m )$ and 
  $\chi^{m} \in \Irr(G^{m} |\phi^m)$, 
 then $(\phi^m)^H  \in \Irr(H| \chi_m)$  while $(\chi^{m} )^G\in \Irr(G |(\phi^m)^H)$.
 \end{itemize}
\end{theorem}

 \begin{proof}
 Since $G_i \unlhd G_j$ and $\chi_j \in \Irr(G_j |\chi_i)$,
 whenever $1\leq i\leq j \leq m $,
we obviously have
\begin{equation}\mylabel{pre1}
G_i^{m} = G_i(\chi  _1, \chi  _2 , \dots ,\chi  _{m}) = 
 G_i(\chi  _1, \chi  _2 , \dots ,\chi  _{i-1})=G_i(\chi_1, \dots, \chi_i), 
\end{equation}
for any $i=1, 2, \dots, m$.
Thus $G_k^{m}  = G_i^{m} \cap G_k = G^{m}\cap G_k$, 
for all $k$ with $1\leq k \leq i$.
It is also clear that  
$G_1^{m}  = G_1$.

To prove the rest of the theorem we will use induction on $m$.
First assume that $m=1$. Because $G_1^{1}= G_1(\chi_1) = G_1$,  
the only possible choice of $\chi_1^{1}$ satisfying \eqref{pr01} for $m=i = 1$ is 
\begin{equation}\mylabel{prel:2c}
\chi_{1}^{1} = \chi_1 \in \Irr(G_1^{1}).
\end{equation}
Observe that,  in this case, $G(\chi_1^{1}) = G(\chi_1) = G^{1}$.
Furthermore, if $G_1 \leq H \leq G$,  then  
Clifford's  Theorem provides a bijection between 
the irreducible characters $\phi $ of $H$ lying above $\chi_1$,  and 
the irreducible characters $\phi^1$ of $H^1=H(\chi_1)$ that lie above $\chi_1$. 
Any two such characters $\phi$ and $\phi^1$ correspond   
 if and only if  $\phi^1$ induces $\phi$.
Note also that,  since Clifford's  theory 
   respects multiplicities, any character $\chi \in \Irr(G)$ that lies above $\phi$ 
 corresponds to some $\chi^1 \in \Irr(G^1)$ that lies above 
$\phi^1 \in \Irr(H^1)$, and vice versa.  So (1) and (2) hold for $m=1$.

Now assume that the theorem holds for all $m$ with $m < n$, 
and some integer $n \geq 2$. We will prove it 
is also true when $m=n$. So assume that 
 $1 \unlhd G_1 \unlhd \dots \unlhd G_{n-1}  \unlhd G_n \unlhd   G$ is a
 normal series of $G$, while $\{ \chi_i \}_{i=1}^n$ is a character tower for that series.
 The inductive hypothesis,  applied to the normal  series 
 $1 \unlhd G_1 \unlhd \dots \unlhd G_{n-1}   \unlhd   G$ of $G$ and its character 
tower  $\{ \chi_i \}_{i=1}^{n-1}$, implies the existence of unique 
irreducible characters 
$\chi_i^{n-1}$ of $ G_i^{n-1} = G(\chi_1, \dots, \chi_{n-1})$ that satisfy
the conclusions   of the theorem, for $m=n-1$.
  Let $G_n^{n-1} = G_n(\chi_1, \dots, \chi_{n-1}) 
=G^{n-1} \cap G_n$. Then we  
 clearly have  
\begin{equation}\mylabel{pre5}
G_i^n = G_i(\chi_1, \dots, \chi_n) = G_i(\chi_1, \dots, \chi_{n-1}) =G_i^{n-1},  
\end{equation}
for all $i=1, \dots, n$. 

The character $\chi_n \in \Irr(G_n)$  lies above $\chi_{n-1}$, while the series 
 $1 \unlhd G_1 \unlhd \dots \unlhd G_{n-1}  \unlhd G_n$ is a normal 
series of $G_n$. 
Hence (2) of the theorem  for $m=n-1$ 
implies that there exists a unique irreducible character 
$\chi_n^{n-1}  \in \Irr(G_n^{n-1})$, that
 lies above $\chi_{n-1}^{n-1}$ and induces $\chi_n \in \Irr(G_n)$. 
We set 
\begin{equation}\mylabel{pre4}
\chi_i^n:= \chi_i^{n-1} \in \Irr(G_i^{n-1}) = \Irr(G_i^{n}),  
\end{equation}
for all $i=1, \dots, n$.
Clearly the characters $\chi_i^n$ satisfy \eqref{pr01},  for
$m=n$ and  all $i=1, 2, \dots, n$. 
(Note that $\chi_1^n = \chi_1$.) 
Furthermore, these characters are unique among those that satisfy \eqref{pr01}
for $m=n$. 
Indeed,  assume that 
 $\{ \psi_i^n \} _{i=1}^n$ is a character tower for $\{ G_i^n\}_{i=1}^n$, 
so that $\psi_i^n \in \Irr(G_i^n)$ induces $\chi_i \in \Irr(G_i)$,
 for all $i=1, \dots, n$.
As $G_i ^n = G_i^{n-1}$ for all such $i$, 
the uniqueness  of the characters $\{ \chi_i^{n-1}\}_{i=1}^{n-1}$ achieved from 
the inductive argument, 
implies that $\psi_i^n = \chi_i^{n-1}= \chi_i^n$, 
whenever $i=1, \dots, n-1$. In addition, the character $\chi_n^n$  was picked 
as the  unique character of $G_i^n$ that induces $\chi_n$ and lies above 
$\chi_{n-1}^{n-1}$. Thus $\chi_n^n  = \psi_n^n$.
This proves that there exist unique characters $\{\chi_{i}^n\}_{i=1}^n$ 
that satisfy \eqref{pr01}, for $m=n$.

To prove that (1) also holds for the characters $\{ \chi_i^n\}_{i=1}^n$,  first 
observe that  $G(\chi_1, \dots, \chi_{n-1})= 
 G^{n-1} = G(\chi_1^{n-1}, \dots, \chi_{n-1}^{n-1}) $, by (1) for $m=n-1$.
Hence 
\begin{equation}\mylabel{pre6}
G^{n} = G(\chi_1, \dots, \chi_{n-1}, \chi_n) = 
 G(\chi_1^{n-1}, \dots, \chi_{n-1}^{n-1}, \chi_n).
\end{equation}
Furthermore,   $G(\chi_1^{n-1}, \dots, \chi_{n-1}^{n-1})$
normalizes  both $G_{n-1}^{n-1}$  and $G_n^{n-1}$, 
and fixes the character  
$\chi_{n-1}^{n-1}$.  As $\chi_n^{n-1}$ is the unique character 
of $G_n^{n-1}$ 
that lies above $\chi_{n-1}^{n-1}$ and induces $\chi_n$, we conclude that 
\linebreak 
$G(\chi_1^{n-1}, \dots, \chi_{n-1}^{n-1})( \chi_n)= 
G(\chi_1^{n-1}, \dots, \chi_{n-1}^{n-1})( \chi_n^{n-1})$. This, along with 
\eqref{pre6} and \eqref{pre4},  implies 
\begin{equation}\mylabel{pre7}
G^n=  G(\chi_1^{n-1}, \dots, \chi_{n-1}^{n-1})( \chi_n^{n-1})= 
 G(\chi_1^{n}, \dots, \chi_{n-1}^{n}, \chi_n^n)= G^{n-1}(\chi_n^n).
\end{equation}
So (1) holds for the inductive step.

Assume  now that $H$ is any subgroup of $G$ with $G_n \leq H \leq G$, while 
$ \phi, \chi$ are irreducible characters of $H$ and $G$, respectively,
so that $\chi$ lies above $\phi$ and $\phi$  above $\chi_n$.      
Then they both  lie above $\chi_{n-1}$ as well. 
Hence the inductive hypothesis implies that there exist 
 unique characters $\phi^{n-1} \in \Irr(H^{n-1} | \chi_{n-1}^{n-1} )$ 
and  $\chi^{n-1}  \in \Irr(G^{n-1} | \phi^{n-1})$
that induce $\phi \in \Irr(H)$ and $ \chi \in \Irr(G)$, respectively.
Observe also that the inductive hypothesis for (2) guarantees that 
the characters $\chi^{n-1}$ and $\phi^{n-1}$ lie above 
$\chi_n^{n-1}  \in \Irr(G_n^{n-1})$, since both $\chi$ and $\phi$ lie above $\chi_n$,
and $\chi_n^{n-1}$ is the unique character of $G_n^{n-1}$   
that induces $\chi_n$ and lies above $\chi_{n-1}^{n-1}$.
Because $G_n^{n-1} $ is a normal subgroup of both 
 $H^{n-1} $ and $ G^{n-1}$, 
  Clifford's Theorem implies the existence of  unique irreducible 
characters  $\chi^n \in \Irr(G^{n-1}(\chi_n^{n-1}))$ 
and $\phi^n \in \Irr(H^{n-1}(\chi_n^{n-1})$ that 
that lie above $\chi_n^{n-1}$ 
and induce $\chi^{n-1}$ and $\phi^{n-1}$,   respectively.
Hence $\chi^n$ induces $\chi$,  and $\phi^n$ induces $\phi$.  
Furthermore, Clifford's Theorem implies that $\chi^{n}$ lies above $\phi^n$,
because $\chi^{n-1}$ lies above $\phi^{n-1}$. 
In addition, $G^{n-1}(\chi_n^{n-1})= G^{n-1}(\chi_n^n)= G^n$, by \eqref{pre7} 
and the fact that $\chi_n^n:= \chi_n^{n-1}$. So $H^n = H^{n-1}(\chi_n^{n-1})$.
We conclude that $\chi^n$ and $\phi^n$ satisfy (2) for the inductive step.
Furthermore, they are unique with that property. 
Indeed, assume that  $\eta \in \Irr(H^n| \chi_n^n)$ and   $\psi \in \Irr(G^n | \eta)$ induce 
 $\phi \in \Irr(H)$ and $\chi \in \Irr(G)$,  respectively. 
Then $\eta^{H^{n-1}}$ and $\psi^{G^{n-1}}$  are 
irreducible characters of  $H^{n-1}$ and $G^{n-1}$,
(since  $H^n \leq  H^{n-1} \leq H$ and  $G^n \leq  G^{n-1} \leq G$),  
that induce  $\phi$ and $\chi$, repsectively.
Also
 $\psi^{G^{n-1}}$ lies above $\eta^{H^{n-1}}$,   and $\eta^{H^{n-1}}$
  above $ \chi_{n-1}^{n-1} $, 
because $\chi_n^n = \chi_{n}^{n-1}$ lies above this last character, 
and $\eta$ lies above $\chi_n^n$.  Hence the inductive hypothesis forces 
 $\eta^{H^{n-1}} = \phi^{n-1}$ and 
$\psi^{G^{n-1}} = \chi^{n-1}$. 
Therefore, $\psi \in \Irr(G^{n})= 
\Irr(G^{n-1}(\chi_n^n))$  induces $\chi^{n-1}$  and lies above $\chi_n^n= 
\chi_n^{n-1}$, 
while $\eta \in \Irr(H^n)$  induces $\phi^{n-1}$ and lies above $\chi_n^{n-1}$.
 In conclusion, $\eta, \psi $ are  the $\chi_n^{n-1}$-Clifford correspondents
of  $\phi^{n-1}$ and  $\chi^{n-1}$, respectively. So $\eta = \phi^n$ and 
$\psi = \chi^{n}$.  Thus $\phi^n$ and $\psi^n$ are unique.

Conversely, if $\phi^n  \in \Irr(H^n  | \chi_n^n)$ then  Clifford's Theorem implies that 
$(\phi^n)^{H^{n-1}}$ is  an irreduible character of $H^{n-1}$. 
since $H^n = H^{n-1}(\chi_n^n)$. Furthermore, $(\phi^n)^{H^{n-1}}$  lies above
$\chi_{n-1}^{n-1}$, as $\chi_n^n = \chi_n ^{n-1}$ lies above $\chi_{n-1}^{n-1}$.
This,  along with the inductive argument, implies  that 
$ \phi^n$  induces an irreducible character of $H$ lying above 
$\chi_n$.  Similarly we can work with any character $\chi^n \in \Irr(G^n | \phi^n)$
to show that it induces irreducibly 
to a character of $G$ lying above $(\phi^n)^H$.
Hence (2)  for the inductive step  is proved. 
Therefore  Theorem \ref{prel:t1} holds.
\end{proof}

\section{Character extensions}
The following is a well known and useful result: 
\begin{theorem}\mylabel{extension}
Assume that  $N$ is a normal subgroup of $ G$  such that $(|N|, |G/N|)=1$. Let 
$\chi$ be a $G$-invariant irreducible character of $N$. 
Then there exists a unique extension $\chi^e$ of $\chi$ to $G$ such that 
$(o(\chi^e), |G/N|)=1$. This is called the canonical extension of $\chi$, 
and has the additional property that  it 
is the unique extension of $\chi$ to $G$ such that 
$o(\chi)= o(\chi^e)$.
\end{theorem}
\begin{proof}
The proof follows easily from Lemma 6.24 and Theorem 11.32 of  \cite{is}, as both 
$o(\chi) $ and $\chi(1)$ divide $|N|$,  and thus  are coprime to $|G:N|$.
\end{proof} 

In the situation of the preceding theorem we can describe all the irreducible 
constituents of $\chi^G$,  as the next result shows.
\begin{theorem}[Gallagher]\mylabel{Galla}
Assume that $N$ is a normal subgroup of $G$. 
 Let $\chi \in \Irr(N)$ be a $G$-invariant character of $N$ that extends 
to an irreducible
 character $\psi$ of $G$. 
 Then there is a one--to--one correspondence
 between the irreducible characters $\gamma$ of $G$ lying above $\chi$
and the irreducible characters $\lambda$ of $G/N$. 
Two such $\gamma$ and $\lambda$ correspond  if and only if 
$\gamma = \lambda \cdot \psi$. The latter product is defined as  
$$
(\lambda \cdot \psi)(g) = \lambda(gN) \cdot \psi(g),
$$
for any $g \in G$.
\end{theorem}
\begin{proof}
The proof follows immediately from Theorem \ref{extension} and 
Corollary 6.17 in \cite{is}.
\end{proof}

Note that,  under the addition hypothesis $(|N|, |G/N|)=1$, 
 we could have used the unique canonical 
extension $\chi^e$ in the place of $\psi$ in 
 Theorem \ref{Galla}.

\begin{lemma}\mylabel{lemma0} 
Assume that  $G$ is  a finite group, $H$ is a normal subgroup of $G$,  and $N$ is a 
normal subgroup of $H$. Assume further that $\theta$ is an irreducible 
character of $N$, and that its stabilizer $H(\theta)$ in $H$ is the 
 product $H(\theta)= N \cdot A$ of $N$ with a subgroup $A$ of $H$, such 
that $(|A|,|N|) = 1$.

\noindent
Then $\theta$ has a canonical extension $\theta^e$ to $H(\theta)$. Furthermore,
 any irreducible character $\Psi_{\gamma}$ of $H(\theta)$ lying above
 $\theta$ is of the form   $\Psi_{\gamma} = \gamma \cdot \theta^e $  where  
 $\gamma  \in \Irr(A)$ and the product is defined as 
$(\gamma \cdot \theta^e ) (s\cdot t) = \gamma(t) \cdot \theta^e(st) $
 for any $s \in N$  and any $t \in A$.
As for  the stabilizer $G(\gamma , \theta) $
 of $ \gamma$ and  $\theta$ in $G$,  we have:
$$
G(\gamma, \theta) =N(A \tin G(\theta,\Psi_{\gamma}))
 =N(A \tin G(\theta,\Psi_{\gamma}^{H})),
$$
where $\Psi_{\gamma}^H$ is the irreducible   character of $H$ 
induced by  $\Psi_{\gamma}$.
\end{lemma}

\begin{proof}
Theorem \ref{extension} above implies that $\theta$ has a unique
 canonical extension  to  $H(\theta)$. This, along with 
 Gallagher's Theorem  \ref{Galla}, implies all but the last statement 
 of the lemma.

As $G(\gamma, \theta)$ fixes $\gamma$ and $\theta$, it normalizes $A$ and $N$ 
respectively. Hence it  normalizes the product $H(\theta) = N A$.
The  canonical extension, $\theta^e$,  of $\theta$ to $H(\theta)$  is fixed by 
 $G(\gamma, \theta)$, as the latter fixes $\theta$ and normalizes $H(\theta)$.
Therefore  $G(\gamma, \theta)$ fixes the product 
$\gamma \cdot \theta^e = \Psi_{\gamma}$. So  $G(\gamma, \theta)$ is a
 subgroup of $N( A \tin G(\theta, \Psi_{\gamma}))$. 
 Because   $H$ is a normal subgroup of $G$, any subgroup of $G$  
that fixes $\Psi_{\gamma}$ also fixes the induced character 
$\Psi_{\gamma}^H$. Hence $N(A \tin G(\theta,\Psi_{\gamma}))$ is a subgroup of 
$N(A \tin G(\theta,\Psi_{\gamma}^{H}))$. 
Therefore
$$
G(\gamma, \theta) \leq N(A \tin G(\theta,\Psi_{\gamma})) \leq 
 N(A \tin G(\theta,\Psi_{\gamma}^H)). 
$$
For the other inclusions we note that any element $g \in N(A \tin G(\theta, 
\Psi_{\gamma}^H))$  fixes $\theta$
and $\Psi_{\gamma}^H$,  normalizes $H(\theta)$, and fixes the unique 
Clifford correspondent  $\Psi_{\gamma} \in  \Irr(H(\theta)| \theta)$ of 
$\Psi_{\gamma}^H$.  Furthermore, $g$ 
 fixes   $\theta^e \in \Irr(H(\theta))$, because it fixes $\theta$ and normalizes 
$H(\theta)$.  As $g$ 
 normalizes $A$, and fixes the product character $\Psi_{\gamma}=
\gamma \cdot \theta^e$,  it  also fixes $\gamma$.  Hence $g \in G(\gamma, \theta)$.
So  $N(A \tin G(\theta,\Psi_{\gamma}^H))
\leq  G(\gamma, \theta)$,   and the proof of 
the lemma is complete.
\end{proof}

\begin{proposition}\mylabel{prel:100B}
Let $G$ be  finite group of odd order such that $G = N\cdot K$,  where 
$N$ is a normal subgroup of $G$ and $(|G/N|, |N|) = 1$. Let  
$H = N\cap K$ and let $\theta$  be any irreducible  $K$-invariant character 
of $H$ that induces  a $G$-invariant irreducible  character $\theta^N$ of $N$.
Then $\theta^N$ has a unique canonical extension, $(\theta^N)^e$,  to  $G$
 such that  
$(|G/N|, o((\theta^N)^e)) = 1$, while $\theta$ has a unique
 canonical extension,
$\theta^e$, to  $K$ such that $(|K/H|, o(\theta^e) ) = 1$.
Furthermore, $\theta^e$ induces 
$$(\theta^e)^G = (\theta^N)^e.$$
\end{proposition}

\begin{proof}
Let $\pi$ be the set of primes that divide $|N|$.  Then 
$|K/H|= |G/N|$ is a $\pi'$-number,   and thus is coprime to $|H|$. 
As $\theta \in \Irr(H)$ is $K$-invariant, Theorem \ref{extension}
implies that $\theta$ has a unique extension  $\theta^e$ to $K$, with 
\begin{equation}\mylabel{pr.e0}
o(\theta)= o(\theta^e).
\end{equation}
 
According to Corollary 4.3 in \cite{isa3},  induction defines a bijection \,$\Irr(K |\theta) \to \Irr(G |\theta^N)$. 
Therefore,  
\begin{equation}\mylabel{pr.e1}
\chi:=( \theta^e)^G  \in \Irr(G | \theta^N).
\end{equation}
But $\theta^N$ extends to $G$, as it is $G$-invariant and $(|N|, |G/N|)=1$. 
Let $\Psi = (\theta^N)^e \in \Irr(G)$ be the unique extension of $\theta ^N$ 
such that  $o(\Psi)= o(\theta^N)$ 
is a $\pi$-number (see Theorem \ref{extension}). Since $\chi$ lies above $\theta^N$,  
 Theorem \ref{Galla} implies that 
$$\chi= \mu \cdot \Psi,$$
for some $\mu \in \Irr(G/N)$. We compute the degree $\deg(\chi)$ in two ways.
First 
$$
\deg(\chi)= \deg(\mu) \cdot \deg(\Psi)= \deg(\mu) \cdot \deg(\theta^N)= \deg(\mu)  \cdot |N:H| \cdot \deg(\theta). 
$$
As $\chi= (\theta^e)^G$ we also have that 
$$
\deg(\chi) = |G:K| \cdot \deg(\theta^e)= |G:K|\cdot  \deg(\theta) = |N:H| \cdot \deg(\theta).
$$
We conclude that $\deg(\mu)=1$. Thus $\mu \in \Lin(G/N)$.
Therefore
\begin{equation}\mylabel{pr.e2}
\det(\chi) = \mu^{\Psi(1)} \det(\Psi).
\end{equation}
  
We can now compute $o(\chi) $ in two ways.  First,
$o(\Psi) = o(\theta^N) $  and  $\Psi(1) = \theta^N(1)$ are  $\pi$-numbers. 
Since $\mu \in \Irr(G/N)$,  we get that $o(\mu)$ is a $\pi'$-number. 
Therefore,  \eqref{pr.e2} implies that  the $\pi'$-number $o(\mu)$ divides $o(\chi)$.

On the other hand, \eqref{pr.e1} and  Lemma 2.2 in \cite{isa4} imply that 
$$
o(\chi) = o((\theta^e)^G)  \text{ divides }  2\cdot  o(\theta^e).
$$
As $G$ has odd order, we get that 
$o(\chi)$ divides $o(\theta^e) $. In view of \eqref{pr.e0},  we have $o(\theta^e) = 
o(\theta)$, while $ o(\theta)  \bigm |   |H|$. We conclude that 
$ o(\chi) $ is a $\pi$-number.

Hence  the only way  the $\pi'$-number $o(\mu)$ can divide $o(\chi)$, 
 is if $o(\mu) = 1$.  So $\mu=1$,   and 
$$
(\theta^e)^G = \chi = \Psi  = (\theta^N)^e,
$$
as desired.
\end{proof}

\begin{lemma}\mylabel{100c}
Let $G$ be any finite group, and $H$ be  any subgroup of $G$. If $\theta \in \Irr(G)$
and $\Psi \in \Irr(H)$ then:
$$ (\theta_{H} \cdot \Psi )^{G} = \theta \cdot (\Psi^G).$$ 
\end{lemma}

\begin{proof}
See Exercise 5.3 in \cite{is}.
\end{proof}
As a corollary  of Lemma \ref{100c} we can prove 
\begin{corollary}\mylabel{100d}
Let $G= H \ltimes M$ be a finite group,   and  $S$ be an $H$-invariant 
subgroup of $M$.
Let $\alpha, \theta$ be irreducible characters of $H$ and $H\ltimes S$   
 respectively. 
Then $\alpha \cdot \theta$ is a character of $HS$ defined as 
$(\alpha \cdot \theta) (x\cdot y) = \alpha(x) \cdot \theta(x\cdot y )$, whenever $x \in H$ and $y \in S$.
 Furthermore,  
$$
(\alpha\cdot \theta)^{G} = \alpha \cdot \theta^G,
$$
where $ (\alpha \cdot \theta^G)(x \cdot y) = \alpha(x)\cdot \theta^G(x \cdot y)$, 
whenever $x \in H$ and $y \in M$. 
\end{corollary}
\begin{proof}
Using the isomorphism $H \cong HS /S$,  we regard  $\alpha$ as a  character of $HS$, defined 
 as $\alpha(x\cdot y) = \alpha(x)$, for all $x \in H$ and $y\in S$.
It is obvious that the product $\alpha \cdot \theta$ is a character of $HS$. 

Furthermore, as $H \cong HM / M = G /M$, we can define an 
irreducible character $\alpha' \in \Irr(G)$ 
as 
\begin{equation}\mylabel{100cd}
\alpha'(x\cdot m)= \alpha(x),
\end{equation}
 for all $x\in H$ and $m \in M$.  
Thus  the restriction $\alpha'|_{HS}$ of $\alpha'$ to $HS$ is 
$\alpha \in \Irr(HS)$, i.e., 
$$\alpha'|_{HS} = \alpha.
$$
Therefore
\begin{align*}
(\alpha\cdot \theta)^{G} &=(\alpha'|_{HS} \cdot \theta)^G &\\
&=\alpha'\cdot \theta^G  &\text{by Lemma \ref{100c}} \\
&= \alpha \cdot \theta^G  &\text{ by \eqref{100cd}.} 
\end{align*}
This completes the proof of the corollary.
\end{proof}

%%% Local Variables: 
%%% mode: latex
%%% TeX-master: "intr"
%%% End: 

\chapter{ A Key Theorem }

One of the main ideas of this thesis is the way
an irreducible character of a finite group $G$ may correspond to 
an irreducible character of a subgroup of $G$.
The  most common such example is the Clifford correspondence. 
Other interesting and fruitful examples are the 
 Glauberman correspondence and the Isaacs correspondence,  
that coincide when applied to groups we are dealing with in
 this thesis, groups  of odd order. According to these 
 correspondences, whenever an odd group $A$ acts on an odd group $G$, 
with $(|A|, |G|)= 1$, there is  `natural' bijection between the set 
$\Irr^A(G)$ of  $A$-invariant irreducible characters of $G$,
 and the set $\Irr(C(A \tin G))$  of irreducible  characters of the 
fixed points $C(A \tin G)$ of $A$ in $G$.

The Glauberman-Isaacs correspondence can be easily generalized
to involve a normal series of subgroups of $G$ and not only $G$. That is,  if 
$G_1 \unlhd  G_2 \unlhd  \dots \unlhd  G_n = G$ is a normal series of $G$,  and $A$ 
acts coprimely on  $G_i$, for all $i=1, \dots, n$, (both $A$ and $G$ are 
assumed  of odd order), 
then we  still have a bijection between the set of  towers $\{ \chi_i \}_{i=1}^n$
of $A$-invariant 
 characters for the series of $G_i$,   and the set of towers 
of irreducible characters of the normal series 
$C(A \tin G_1) \unlhd  C(A \tin G_2) \unlhd  \dots \unlhd  C(A \tin G_n)$ of 
$C(A \tin G)$. 

However,  whether we use the Glauberman-Isaacs correspondence on
a single group  $G$ or on  a  normal series of $G$,  
 the condition that wants   the acting group $A$ to 
normalize every group $G_i$ involved in  the series can't be avoided. 
The solution to this problem is given by E.C.Dade
in  \cite{da5}. Here we only  state the results needed from that paper.
The easy case, where the series is replaced by a single group, is 
done in Theorem \ref{dade:t1}, while the general case is described 
in Theorem \ref{dade:t2}.  

Using the notation introduced in Section \ref{nota}, 
 we write $\Irr^A(G)$ for the $A$-invariant
 characters of $G$,  whenever $A$ acts on $G$. 
In addition,  if $A$ acts on a subgroup $B$ of $G$, we write $\Irr^A_B(G)$ for the irreducible characters of $G$ 
that lie above at least one $A$-invariant character of $B$. 
Furthermore, if $M$ is a subgroup of $G$ and $\mu \in \Irr(M)$ we denote by  
$\Irr(G|\mu)$ the set of irreducible characters of $G$ that lie above $\mu$, 
and by $\Irr^A_B(G|\mu)$ the intersection
$$\Irr^A_B(G|\mu):= \Irr^A_B(G) \cap \Irr(G|\mu).$$
If $\chi \in \Irr(G| \mu)$,  we 
 write $m(\mu \tin \chi)$ for the multiplicity that  $\mu$ appears 
as a constituent of the restriction 
$\chi|_{M}$ of $\chi$ in $m$, i.e., 
 $m(\mu \tin \chi) = [ \chi|_{M}, \mu]$. 

\begin{theorem}\mylabel{dade:t1}
Assume that $G$ is a finite group of odd order, 
 and that  $B$ is a normal subgroup of
$G$. Let $A ,H$ be subgroups of $G$  such that  $(|A|, |B|) = 1$,  while 
$B$ is contained in $H$.
Assume further that the subgroup 
$AB = A\ltimes B$ is normal in $G$.
Then there  is a one to one 
 correspondence   
$$
\psi \in \Irr^A_B(H) \leftrightarrow  \psi_{(A)} \in \Irr(N(A \tin H)), 
$$
between the set $\Irr^A_B(H)$ of all  characters $\psi \in \Irr(H)$ 
such that $\psi$ lies 
above at least one irreducible  $A$-invariant character of $B$, and the set 
$\Irr(N(A \tin H))$ of all irreducible characters of the normalizer
$N(A\tin H)$ of  $A$ in $H$. We call the correspondence 
$\psi \leftrightarrow \psi_{(A)}$ an $A$-correspondence,  
and  say that the characters
   $\psi$ and  $\psi_{(A)}$ are 
$A$-correspondents of each other.

If  $H = B$, then the above  $A$-correspondence coincides with the 
Glauberman-Isaacs correspondence between 
$\Irr^A(B)$ and $\Irr(C(A \tin B))$.

Furthermore, for any subgroup $K$ of $G$ that normalizes both 
$A$ and $H$, the stabilizer, $K(\psi)$, of any $\psi \in \Irr^A_B(H)$ in $K$
equals the corresponding stabilizer, $K(\psi_{(A)})$, of $\psi_{(A)}$ in $K$.   
\end{theorem}
\begin{proof}
The $A$-correspondence is done in  Theorem 17.4 in \cite{da5}, with 
$A, B, H$ here in the place of $K, L , H$ there. (Observe that, 
 in that theorem,  the $A$-correspondence is 
described in a more general setting.)
 
Theorem 17.29 in \cite{da5} implies that  $A$-correspondence 
and Glauberman correspondence coincide when $H= B$.

The last part of the theorem follows easily  from Proposition 17.10 
in \cite{da5}.
\end{proof}

Before we continue to the general case, we note that the group $H$ 
in Theorem \ref{dade:t1} doesn't need to be normal subgroup of $G$. 
Furthermore, the $A$-correspondence described in the above theorem
  depends only on $A$ and $H$ and not on the choice of $B$ (see  
 Proposition 17.13 in \cite{da5}).

Theorem  17.4  in \cite{da5} not only  provides the character 
 correspondence  we describe in Theorem \ref{dade:t1}, 
 but also  gives a specific algorithm 
we can use to obtain   this correspondence.  It tells us that 
\begin{theorem}\mylabel{dade.14}
Assume that $A, B, H$ and $G$ satisfy the hypotheses 
of  Theorem \ref{dade:t1}. Assume  further that $\psi$ is an irreducible 
character of $H$ that lies above at least one irreducible $A$-invariant 
 character of $B$, i.e., $\psi  \in \Irr^A_B(H)$.  Then there exists 
some sequence $M_0, M_1, \dots, M_n$ of subgroups  of $G$ satisfying 
\begin{subequations}\mylabel{dd.1}
\begin{gather}
n\geq 1, \quad  M_0=B \geq M_1 \geq \dots \geq M_n = 1, \\
M_i \unlhd H, \quad  [M_i , A] \leq M_i,  \text{ and } \\
M_{i-1}/M_i \text{ is abelian, }  
\end{gather}
\end{subequations}
for all $i=1, 2, \dots, n$. Any such sequence of subgroups determines a
 unique sequence of characters $\theta_0, \theta_1, \dots, \theta_n$ 
such that  
\begin{subequations}\mylabel{dd.2}
\begin{equation}
\theta_{0} = \psi \in \Irr^A_B(H) = \Irr^A_{M_0}(N(AM_0 \tin H))
\end{equation}
  and  
\begin{multline}
\text{ If  $i=1, 2, \dots, n$,  then $\theta_i$ is the unique character in 
$\Irr^A_{M_i}(N(A M_i \tin H))$ } \\
\text{  such that  $m(\theta_i \tin \theta_{i-1})$  is odd.} 
\end{multline}
\end{subequations}
The character $\theta_n \in \Irr^A_{1}(N(A  \tin H)) = \Irr(N(A \tin H))$
 is independent of the choice of the 
sequence $M_0, M_1, \dots, M_n$ satisfying \eqref{dd.1}. Furthermore, 
$\theta_n$ is the  $A$-correspondent  $\psi_{(A)} \in \Irr(N(A \tin H))$ 
of $\psi$, as used in Theorem \ref{dade:t1}.
\end{theorem}

It is clear from the above construction that 
the $A$-correspondence is preserved under   epimorphic images. So we  have 
\begin{proposition}\mylabel{daco}
Assume  that  $A, B,  H$  and $G$ satisfy the hypotheses of Theorem 
\ref{dade:t1}. Let 
$\rho$ be an epimorphism of  $G$ onto some group $G'$.
If $A', B'$ and $H'$ are the images under $\rho$ 
of $A, B$ and $H$, respectively,  then 
$A', B', H'$ and $G'$ satisfy the hypotheses of Theorem \ref{dade:t1}. 
Furthermore, $\rho$ maps $N(A \tin H) $ onto 
$N(A' \tin H')$.
Assume further that 
$\psi'  \in \Irr(H') $ and   $\psi = \psi' \circ \rho_H \in \Irr(H)$. Then 
 $\psi' \in \Irr^{A'}_{B'}(H')$ if and only if $\psi \in \Irr^A_B(H)$. 
In that case  $\psi_{(A)} = \psi'_{(A')} \circ  \rho_{N(A \tin H)} \in \Irr(N(A \tin H))$.
\end{proposition}

\begin{proof}
Clearly the groups $G', A', B'$ and $H'$ satisfy the hypotheses in 
Theorem \ref{dade:t1}.

Let  $X$ be  any subgroup of $G$  
and $X'= \rho(X)$ be its image under $\rho$ in $G'$.
 Then $\rho$ restricts to an epimorphism $\rho_X$ 
of $X$ onto $X'$. This induces an injection $\phi' \to \phi' \circ \rho_X$ 
of $\Irr(X')$ into  $\Irr(X)$. 
For any  $\psi'  \in \Irr(H')$ the above injection determines the  irreducible 
 character $\psi \in \Irr(H)$ 
 satisfying $\psi = \psi' \circ \rho_H$. 
 We remark here
that any $ \psi \in \Irr(H)$ with $\Ker(\rho_H) \leq \Ker(\psi)$ 
corresponds,  under the above injection, to a character $\psi' \in \Irr(H')$, 
satisfying 
$\psi = \psi' \circ \rho_H$. 

Evidently if  $\psi'$ lies over a character $\phi' \in \Irr(B')$ then  
$\psi$ lies over the corresponding character 
$\phi = \phi' \circ \rho_B \in \Irr(B)$.  In addition, if $A'= \rho(A)$ 
fixes $\phi'$, then $A$ fixes $\phi$. Hence if 
 $\psi' \in \Irr^{A'}_{B'}(H')$  then  $\psi \in \Irr^A_B(H)$. 
Conversely assume that 
 $\psi \in  \Irr^A_B(H)$ satisfies  $\Ker(\rho_H ) \leq \Ker(\psi)$. Then 
$\psi$ has a corresponding character 
 $\psi' \in \Irr(H')$ satisfying $\psi = \psi' \circ \rho_H$. 
Because $\psi \in \Irr^A_B(H)$, 
 there exists an irreducible $A$-invariant character $\phi \in \Irr^A(B)$
of $B$ that lies under $\psi$.
Then $\Ker(\phi) \geq \Ker(\psi) \cap B \geq \Ker{\rho}  \cap B$..
 Hence $\phi$ has a corresponding character 
$\phi' \in \Irr(B')$  satisfying $\phi = \phi' \circ \rho_B$. 
Furthermore, 
$\psi'$  lies above $\phi'$, because $\psi$ lies above $\phi$.
In addition,  $A'$ fixes $\phi'$ because 
$A$ fixes $\phi$. Hence $\psi'  \in \Irr^{A'}_{B'}(H')$.
In conclusion, 
$\psi' \in \Irr^{A'} _{B'} (H)$  if and only if $\psi \in \Irr^A_B(H)$. 

If $\psi'  \in \Irr^{A'}_{B'}(H')$, then Theorem \ref{dade:t1} applies.  So $\psi'$  has 
 an $A'$-correspondent irreducible character 
$\psi'_{(A')}$ in  $N(A' \tin H')$.
To complete the proof of the proposition it suffices to show that 
$\rho(N(A \tin H)) = N(A' \tin H')$   while  $\psi_{(A)} = \psi'_{(A')} \circ
 \rho_{N(A \tin H)}$.

Let $M$  be any subgroup of $B$ such that $ M \unlhd H$ and $[M, A] \leq M$.
Then $M' = \rho(M)$ is a subgroup of $B'$ such that $M' \unlhd H'$ and 
$[M', A'] \leq M'$. We claim that 
\begin{equation}\mylabel{dae}
\rho(N(AM \tin H))  = N(A'M' \tin H').
\end{equation}
Clearly $\rho(N(AM \tin H))  \leq  N(A'M' \tin H')$. Thus to prove  \eqref{dae} 
it is enough to show that any $t' \in N(A'M' \tin H')$ is in the image of 
$N(AM \tin H)$.  Since $t' \in H'= \rho(H)$ 
 there exists a $t \in H$ with $\rho(t)= t'$.
If $K = \Ker(\rho)$, then $t$ normalizes $AMK$, 
since $t'$ normalizes $A'M' = \rho(AM)$. 
In addition, $t$ normalizes $AB \unlhd G$. 
Hence $t$ normalizes the intersection 
$AB \cap AMK= A(B \cap AMK)$.
 The group $A$ normalizes $B \cap AMK$, since it normalizes 
$B, M$ and  $K$. In addition, $(|A|, |B \cap AMK|)= 1$, because 
$(|A|, |B|)=1$.  Hence $AB \cap AMK= 
A (B \cap AMK) = A \ltimes (B \cap AMK)$, and the
 $B \cap AMK$-conjugates of $A$ are the only subgroups of order $|A|$ 
in $A (B \cap AMK)$. It follows that there is some element $s \in B \cap AMK$ 
such that $A ^{ts}=A$. But $ts \in H$  normalizes $M \unlhd H$. Thus 
$(AM)^{ts}= AM$. So $ts \in N(AM \tin H)$ has image $\rho(ts)= t' \rho(s) \in 
N(A'M' \tin  H')$.  In addition, the image $\rho(s)$ of $ s \in AMK$ is 
 an element of 
$\rho(AMK) =A' M'= \rho(AM)$ and thus lies in $\rho(N(AM \tin H))$. 
We conclude that  $t' = \rho(ts) \rho(s)^{-1}$ lies in $\rho(N(AM \tin H))$. 
Therefore \eqref{dae} holds.

Observe that \eqref{dae}  
for $M=1$, implies that $\rho(N(A \tin H)) = N(A' \tin H')$. 
Furthermore,  if $M_0, M_1, \dots, M_n$ are any subgroups of $G$ satisfying 
\eqref{dd.1}, their images $M_i' = \rho(M_i)$, for $i=0, 1, \dots, n$, satisfy the
 equivalent of \eqref{dd.1} for  $G', H', B'$ and $A'$.  According to Theorem 
\ref{dade.14}, the character $\psi' \in \Irr^{A'}_{B'}(H')$  determines characters 
$\theta_0' , \dots, \theta_n'$ such that 
$\theta_0' = \psi'$  and $\theta_i'$ , for any $i=1, \dots, n$,  is the unique character 
in $\Irr^{A'}_{M_i'}(N(A'M_i' \tin H'))$ such that $m(\theta_i' \tin \theta_{i-1}')$ is odd.
Since $\rho$ sends $N(AM_i \tin H)$ onto $N(A' M_i' \tin H')$, by \eqref{dae}, 
it follows that $\theta_i = \theta_i' \circ \rho_{N(AM_i \tin H)}$ lies in 
$\Irr^{A}_{M_i}(N(AM_i \tin H))$, for each $i=0, 1, \dots, n$. Furthermore, 
$\theta_0 = \psi' \circ \rho_H  = \psi$  and $m(\theta_i \tin \theta_{i-1})=
 m(\theta_i' \tin \theta_{i-1}')$ is odd, 
for each $i=1, \dots, n$. Theorem \ref{dade.14}  then implies 
$$
\psi_{(A)} = \theta_n = \theta_n' \circ \rho_{N(AM_n \tin H)} = 
\psi'_{(A')} \circ \rho_{N(A \tin H)}.
$$
Hence Proposition \ref{daco} holds.
\end{proof}

Proposition \eqref{daco} easily implies 
\begin{corollary}\mylabel{daco2}  
Assume  that  $A, B,  H$  and $G$ satisfy the hypotheses of Theorem 
\ref{dade:t1}. Let 
$\rho$ be an isomorphism of  $G$ onto some group 
 $G'$. If $A', B'$ and $H'$ are the 
isomorphic images, under $\rho$, of $A, B$ and $H$ respectively,
then $ A', B', H'$ and $G'$ satisfy the hypotheses of Theorem \ref{dade:t1}. 
Furthermore, $N(A' \tin H')$ is the isomorphic image under $\rho$ of 
$N(A \tin H) $. In addition, for any $\psi \in \Irr(H)$ there exists a 
$\psi'  \in \Irr(H') $ such that  $\psi = \psi' \circ \rho_H$. 
Then 
 $\psi' \in \Irr^{A'}_{B'}(H')$ if and only if $\psi \in \Irr^A_B(H)$. 
In that case  $\psi_{(A)} = \psi'_{(A')} \circ  \rho_{N(A \tin H)}$.
\end{corollary}

Proposition 17.12 in \cite{da5} implies
\begin{proposition}\mylabel{dade:l1}
If the group $A$ in Theorem \ref{dade:t1} centralizes $B$, 
then  $\Irr^A_B(H) = \Irr(H)$
\linebreak
$ = \Irr(N(A \tin H))$ and the 
 $A$-correspondence is the identity map of these  equal sets onto themselves.
\end{proposition}

From  Proposition 17.14 in \cite{da5} we obtain 
\begin{proposition}\mylabel{dade:p1.5}
Let $A, B, H, G$ be as in Theorem \ref{dade:t1}. Let $A'$ be a subgroup 
of $A$ such that $N(A \tin H) = N(A' \tin H)$,  and thus $N(A \tin B) = N(A'
\tin B)$. Then $\Irr^A_B(H)= \Irr^{A'}_B(H)$ and 
$$
 \psi_{(A)} = \psi_{(A')},
$$
for any $\psi \in \Irr^A_B(H) = \Irr^{A'}_B(H)$.
\end{proposition}

In the special case where the $A$-correspondence is the Glauberman
correspondence
(that is the case $H= B$), Proposition \ref{dade:p1.5} translates to 
\begin{corollary}\mylabel{coo.d1}
Assume that   $A$ acts coprimely on $B$, where $A$ and $B$ are both finite 
groups of odd order.  If $A'$ is a subgroup of $A$ satisfying 
$C(A \tin B) = C(A' \tin B)$  then the $A$-Glauberman 
and the $A'$-Glauberman correspondences coincide. 
\end{corollary}

In the special case that $H= AB$, Theorem 17.36  in \cite{da5} describes 
clearly the $A$-correspondence of Theorem \ref{dade:t1}. So we get 
\begin{theorem} \mylabel{dadet3}
Assume that  $A, B, G$ satisfy the hypotheses  of Theorem \ref{dade:t1}, 
with the additional condition  that $H= AB$.
Assume further that $\chi \in \Irr(H)$ is of the form 
$\chi = \alpha \cdot \beta^e$, where $\alpha \in \Irr(A)$ and 
$\beta^e$ is the canonical extension to $H$  of an irreducible 
$A$-invariant character $\beta \in \Irr^A(B)$.
Then $\chi \in \Irr^A_B(H)$. 
Furthermore,  $N(A \tin H) = A \times C(A \tin B)$, 
where $C(A \tin B)$ is the centralizer of $A $ in $B$. 
In addition, the $A$-correspondent $\chi_{(A)} \in \Irr(N(A \tin H))$   of $\chi$, 
is of the form 
$$
\chi_{(A)}= \alpha \times \gamma, 
$$
where $\gamma \in \Irr(C(A \tin B))$ is the  $A$-Glauberman correspondent of 
$\beta \in \Irr^A(B)$.               
\end{theorem}
\begin{proof}
See Theorem 17.36  in \cite{da5}.
\end{proof}

The next proposition shows that the $A$-correspondence is compatible with 
Clifford correspondence.

\begin{proposition} \mylabel{dade:p1.6}
Let $A, B, H$ and $G$ 
 be as in Theorem \ref{dade:t1} and let  $M$ be an $A$-invariant
normal  subgroup of $G$ contained in $B$.
 Assume further that  $\mu$ is an $A$-invariant irreducible character 
of  $M$, and let $G(\mu), H(\mu)$ and $B(\mu)$  be the stabilizers of $\mu$ in 
$G,H$ and $B$ respectively.  Let $\mu_{(A)} \in \Irr(N(A \tin M))$  be the $A$-correspondent of $\mu$, 
(note this is the $A$-Glauberman correspondent of $\mu$).
If $\psi \in \Irr^A_B(H )$  lies above $\mu$, that is, 
 $\psi \in \Irr^A_B(H | \mu)$, then the $\mu$-Clifford
 correspondent $\psi_{\mu} \in \Irr(H(\mu))$
 of $\psi$ lies in $\Irr^A_{B(\mu)}(H(\mu) | \mu)$.
Furthermore, the $A$-correspondent
 $\psi_{\mu, (A)} \in \Irr(N(A \tin H(\mu)))$ of $\psi_{\mu}$ is the 
$\mu_{(A)}$-Clifford 
correspondent of the $A$-correspondent
 $\psi_{(A)} \in \Irr(N(A \tin H))$ of $\psi$. 
\end{proposition}

\begin{proof}
See Propositions 17.19 17.20, 17.22, 17.23 and Theorem 17.24 in \cite{da5}.
\end{proof}

We conclude with a generalization of Theorem \ref{dade:t1}.
\begin{theorem}\mylabel{dade:t2}
Let $A, B, G$ be as in Theorem \ref{dade:t1}.
Assume further that $G_0=B \unlhd G_1 \unlhd 
\dots \unlhd G_{n-1} \unlhd  G_n  \unlhd G$ is a series of normal subgroups of 
$G$.
Then the groups $N(A \tin G_0)= 
N(A \tin B) \unlhd  N(A \tin G_1) \unlhd  \dots \unlhd  N(A \tin G_{n}) \unlhd  N(A \tin G)$, 
form a  series of normal subgroups of $N(A \tin G)$.
Let  $\psi_0, \psi_1, \dots, \psi_n$ 
be a tower of irreducible characters for the chain 
$G_0 \unlhd G_1 \unlhd \dots \unlhd G_{n-1} \unlhd  G_n$, while 
$\psi_0 \in \Irr^A(B)$. Then 
 $\psi_i \in \Irr^A_B(G_i)$, for all $i=1,2, \dots, n$. Let  
$\psi_{i, (A)} \in \Irr(N(A \tin G_i))$ be the  
$A$-correspondent of $\psi_i \in \Irr^A_B(G_i)$, for all $i=0, 1, \dots, n$.  
Then the $\psi_{i, (A)}$ form a character tower
for the chain $N(A \tin G_0)= 
N(A \tin B) \unlhd  N(A \tin G_1) \unlhd  \dots \unlhd  N(A \tin G_{n}) $.
This way we get a bijection between  character 
towers $\{ \psi_i \}_{i=0}^n$
 for the series  $\{ G_i \}_{i=0}^n$  with  $\psi_0 \in \Irr^A(B)$, 
and character towers
$\{  \psi_{i, (A)} \}_{i=0}^n$ 
for the series $\{ N(A \tin G_i ) \}_{i=0}^n$.
In addition,  this correspondence respects restrictions and inductions, i.e., 
\begin{itemize}
\item[(a)] $\psi_i^{G_{i+1}} = \psi_{i+1}$ if and only if
$\psi_{i, (A)}^{N(A \tin G_{i+1})} = \psi_{i+1, (A)}$,
 while
\item[(b)]  $\psi_{i+1}\big|_{G_{i}} = \psi_{i}$ if and only if
$\psi_{i+1, (A)}\big|_{N(A \tin G_{i})}= \psi_{i, (A)}$,
 \end{itemize}
for any  $i=1, \dots, n-1$.

 Furthermore, for any subgroup $K$ of $G$ that normalizes $A$, 
 the stabilizer, $K(\psi_i)$, of $\psi_i$ 
in $K$ equals the corresponding stabilizer, $K(\psi_{i, (A)})$, of
 $\psi_{i, (A)}$ in $K$,  for all $i= 0,1, \dots ,n$.
\end{theorem}

\begin{proof}
See Theorem 17.15 and Propositions 17.16 and 17.17 in \cite{da5}.
\end{proof}

%%% Local Variables: 
%%% mode: latex
%%% TeX-master: t
%%% End: 

\chapter{ Changing Characters}\mylabel{cc}

\section{A ``petite'' change}
Assume that  $G, Q_1 ,Q_2$ and $P$ are finite groups that satisfy 
\begin{hyp}\mylabel{cc:hyp1}
$Q_1, Q_2$ are $q$-subgroups of   $G$,
 for some odd prime $q$, with   $Q_1 \unlhd Q_2$. Furthermore,
$P$ is a   $p$-group, for some prime $p$ with  $2 \ne p\ne q$, that   normalizes $Q_1$,
while $P\cdot Q_1$ is normalized
 by $Q_2$.
\end{hyp}
 Let $\beta_1$ be an irreducible character of $Q_1$.     
Our main goal in this section is to show how we can change the  character,
$\beta_1 \in 
\Irr(Q_1)$, to a new one, $\bn_1 \in \Irr(Q_1)$, so that  
\begin{itemize}
\item[(1)]  $P(\beta_1) = P(\bn_1)$ and 
$Q_2(\beta_1) \leq Q_2(\bn_1)$, while 
\item[(2)] $\bn_1$ can be extended to $Q_2(\bn_1)$.
\end{itemize}

The following diagram describes that  situation:  
\begin{equation}\mylabel{cc:d1}
\begin{diagram}[small]
   &      &  Q_2    &\qquad   &      &       &    \\
   &\ldTo &\vLine   &\qquad   &      &       &          \\ 
P  &      &Q_2(\bn_1)&\qquad   &      &       &(\bn_1)^e  \\ 
   &\rdTo &\dTo     &\qquad   &      &       &\vLine          \\
   &      &Q_1      &\qquad   &\beta_1 &\rcoto &\bn_1       \\
\end{diagram}
\end{equation}
Most of the work towards that direction is done in
\begin{lemma}\mylabel{cc:l0.5}
Let $Q$ be a $q$-group acting on a $p$-group $P$, with $p \ne q$ odd primes.
Let $T$ be a finite-dimensional right 
 $\ZZ _q(Q\ltimes P)$-module such that the action of $P$ on
$T$ is faithful.
Then there exists an element $\tau \in T$
  such that  its stabilizer $(QP)(\tau)$ in $Q\ltimes P$ equals $Q$.
\end{lemma}

\begin{proof}
We will prove a series of claims under the  
\begin{ia}\mylabel{cc:1}
  $Q ,P, T$  are  chosen among all the triplets satisfying the hypothesis,
but not the conclusion, of Lemma \ref{cc:l0.5}, so as to minimize first the
 order $|QP|$  of the semidirect product $Q\ltimes P$, and then the
 $\ZZ_q$-dimension 
$\dim_{\ZZ_q}T$ of $T$.
\end{ia}
These claims will lead to a contradiction, thus proving the lemma.
\begin{claim}\mylabel{cc:c1}
$T$ is an indecomposable $\ZZ_q(QP)$-module.
\end{claim}

\begin{proof}
Suppose not. Let $T=T_1\, \dot{+} \, T_2$ be a direct decomposition of $T$, 
where $T_1 ,T_2$ are nontrivial $\ZZ_q (QP)$-submodules of $T.$
For $i=1,2$ let $K_i$ be the kernel of the action of $P$ on $T_i.$
Hence $T_i$ is a $\ZZ_q(Q\ltimes P/K_i)$-module such that $P/K_i$ acts 
faithfully on it. As $\dim _{\ZZ_q} (T_i)$ is strictly 
smaller than  $\dim _{\ZZ_q}T$, the minimality in Inductive Assumption
\ref{cc:1}
provides  an  element  $\tau_i \in T_i$ 
such that $(Q \ltimes {P/ {K_i}})(\tau _i)=Q.$
If we take as $\tau$ the sum, $\tau =\tau_1 +\tau_2$,
then $\tau$ is an element of $T$ fixed by $Q$, as $Q$ fixes 
each one of the $\tau_i$ for $i=1,2.$
Furthermore for the stabilizer of $\tau$ in $P$ we have
$$
P(\tau)=\cap_{i=1}^2 P(\tau_i) = \cap_{i=1}^2 K_i.
$$
Since $P$ acts faithfully on $T$ the last intersection is trivial. Therefore
$(QP)(\tau) =Q$, which contradicts Inductive Assumption \ref{cc:1}. Hence
$T$ is an indecomposable $\ZZ_q(QP)$-module.
\end{proof}

\begin{claim}\mylabel{cc:c2}
The restriction $T_P$ of $T$ to  $P$ is a multiple of an irreducible 
$Q$-invariant 
$\ZZ_q(P)$-module.
\end{claim}

\begin{proof}
Claim \ref{cc:c1} and Clifford's Theorem  imply that
$T_p$ can be written as a direct sum of its $\ZZ_q(P)$-homogeneous
components, i.e.,
$$
T_P  = U_1 \dotplus  U_2 \dotplus \dots \dotplus  U_r,
$$ 
where $U = U_1 \cong m V = m V^{\sigma_1} ,U_2 \cong m V^{\sigma_2}, \dots , U_r \cong
mV^{\sigma_r}$. Here 
 $V= V^{\sigma_1}, \dots ,V^{\sigma_r}$, are the distinct conjugates of a
simple 
$\ZZ_q(P)$-submodule, $V$, of $T_p$, and ${1=\sigma_1, \dots, 
\sigma_r}$ are  representatives for the cosets in $Q\cdot P$ of the
 stabilizer,$ (QP)_V$,  of the isomorphism class of $V$ in $QP$.
We may pick ${\sigma_1 , \dots ,\sigma_r}$ to be representatives
of the cosets in $Q$ of  the stabilizer, $Q_V$, of that isomorphism class in
 $Q$.
Note that $Q_V = Q_U$ as $U \cong mV$, where $Q_U$ is the stabilizer in 
$Q$  of $U$ under multiplication in $T$.
If $T_P$ is not homogeneous, then $r > 1$ and $Q_U = Q_V < Q.$
For $i=1, \dots , r$ let $\widehat{K_i}$  be the kernel of the action of $P$
on $U_i$.
 Then  for every $i=1, \dots , r$ the stabilizer $Q_{U_i}$ of $U_i$ in $Q$
equals the $\sigma_i$-conjugate, $Q_U^{\sigma_i}$ ,  of $Q_U = Q_V$.  
For the corresponding kernels we similarly have
 $\widehat{K_i} = {\widehat{K_1}} ^{\sigma_i}$. 

As  $U$ is a faithful $\ZZ_q(P/{\widehat{K_1}})$-module
and $Q_U \lneq Q$,  the minimality in  Inductive Assumption \ref{cc:1}
implies that there exists an element  $\mu \in  U$ such that 
$$
(Q_U \ltimes ( P/{\widehat{K_1}}))(\mu) = Q_U.
$$
For every $i=1, \dots ,r$ we can define an element
$\mu_i = \mu  \sigma_i$ of
$U_i$.
 Then $Q_{U_i}$ fixes $\mu_i$ as $Q_U$ fixes $\mu$.
Furthermore if $x$ is any element of $P$ fixing $\mu_i$ then 
$x^{\sigma_i^{-1}}$ is an element of $P$ fixing $\mu$. Therefore 
$x^{\sigma_i^{-1}} \in \widehat{K_1}$, which implies that
 $x\in \widehat{K_i}$. Thus 
$$
(Q_{U_i}\ltimes (P/{\widehat{K_i}}))(\mu_i) = Q_{U_i}
$$
for every $i=1, \dots , r$.

Let $\tau$ be the sum of the  $\mu_i$ for $i=1, \dots , r.$
Then $\tau $ is an element  of $T$ fixed by $Q$, since multiplication by
any element in $Q$ permutes the $U_i$  and the  $\mu_i$ 
among themselves.
The stabilizer $P(\tau)$ of $\tau$ in $P$ equals the intersection of
the stabilizers of $\mu_i$ in $P$ for $i=1, \dots ,r$.
Since $(Q_{U_i}\ltimes (P/{\widehat{K_i}}))(\mu_i) = Q_{U_i} $ for every such $i$, the
latter equals the intersection of $\widehat{K_i}$ for $i=1, \dots ,r.$
The faithful action of $P$ on $T$ implies that 
$$
P(\tau) = \cap_{i=1}^r \widehat{K_i} = {1}.
$$
Hence $T$ has an element $\tau$ with $(QP)(\tau) = Q$,
contradicting Inductive Assumption  \ref{cc:1}.
This contradiction proves Claim  \ref{cc:c2}.
\end{proof}

\begin{claim}\mylabel{cc:c3}
There are  no $Q$-invariant subgroup, $H < P$, and $\ZZ_q(QH)$-submodule, $S$,
of 
$T_{QH}$ such that $T$ is the $\ZZ_q(QP)$-module
 $S^{QP}$ induced from $S$, i.e.,
$$
T = \sum_{1\leq i\leq n}^{\cdot}S\sigma_i,
$$
where the  $\sigma_i$ are representatives for the cosets $H \sigma_i$  of $H$ in $P$.
\end{claim}

\begin{proof}
Suppose Claim \ref{cc:c3} is false. We choose $H$ to have maximal order among
all those $Q$-invariant subgroups of $P$ that contradict Claim \ref{cc:c3}.
Hence $T_{QH}$ has a 
$\ZZ_q(QH)$-submodule, $S$, such that 
$S^{QP} = T$.
If $H$ is not normal in $P$ then its normalizer, $N_P(H)$, 
in $P$ satisfies $H \lhd N_P(H) < P$. Since $H$ is $Q$-invariant, $N_P(H)$ 
is also $Q$-invariant. Hence $S^{Q N_P(H)}$ is a $\ZZ_q(QN_P(H))$-submodule of
 $T_{QN_P(H)}$. Furthermore $S^{QN_P(H)}$ induces $T$.
Thus $N_P(H)$ is among the $Q$-invariant subgroups of $P$ that contradict
Claim \ref{cc:c3} with $|N_P(H)| > |H|$. So the maximality of $|H|$ implies that
$H$ is normal in $P$. 

Let $1=r_1 , \dots , r_k$ be coset representatives of
$H$ in $P$, and let $\bar{r}_m$ denote the
 image of $r_m$ in $P/H$ for $m=1,\dots,k$.
 Then $\bar{1} = \bar{r_1}, \bar{r_2},\dots,\bar{r_k}$
are the distinct elements of $P/H$.
As $Q$ acts on $P/H$, it has to divide  the 
$\bar{r}_m$, for $ m=1, \dots , k$, into  orbits, $\overline{R_1}, \overline{R_2}, 
\dots,\overline{R_l}$, for some
$l\in \{1, \dots ,k\}$. We may choose $\overline{R_1}$ to be 
equal to  $\{{\bar{r}_1}\} = \{{1}\} $. For every $i = 2, \dots ,l$, we
pick some element $\bar{r}_{i,1} \in \overline{R_i}$. Then 
 $\overline{R_i} =\{ \bar{r}_{i,1}^{q_j}\}_{j=1}^{j=k_i} $  where  
$k_i = |\overline{R_i}|$ and 
$q_j $ runs over a set $Q_j$ of coset  representatives of the stabilizer,
$C_Q(\bar{r}_{i,1})$, in $Q$.  
For every $i=2, \dots , l$ the stabilizer $C_Q(\bar{r}_{i,1})$ acts
by conjugation on $H$ and on $r_{i,1}H$, where $r_{i,1} \in P$ has
image $\bar{r}_{i,1} \in P/H$.
Furthermore, $H$ acts transitively by right multiplication on 
$r_{i,1}H$ and $(x h)^c = x^c h^c$  for all
 $x\in r_{i,1}H , h\in H , c \in C_Q(\bar{r}_{i,1}) $.
Hence Glauberman's Lemma (13.8 in \cite{is}) provides an element 
$t_{i,1} \in r_{i,1}H$ that is fixed by $C_Q (\bar{r}_{i,1})$.
So $C_Q(t_{i,1}) \geq  C_Q(\bar{r}_{i,1})$.
Furthermore, the opposite inclusion, $C_Q(t_{i,1}) \leq 
C_Q(\bar{r}_{i,1})$, 
 also holds as $\bar{r}_{i,1} = t_{i,1}H$. Hence, 
$$
C_Q(t_{i,1}) = C_Q(\bar{r}_{i,1}).
$$
In this way we can pick a $t_{i,1} \in r_{i,1}H$, for every ${i= 1, \dots , l}$, 
 such that $C_Q (t_{i,1}) = C_Q (\bar{r}_{i,1})$.
We can even assume that  $t_{1,1} = 1$. 
Let $t_{i,j}$ denote the $q_j$-conjugate, $ t_{i,1}^{q_j}$, 
of $t_{i,1}$ for every $j=1, \dots , k_i$.
Hence the set of all $t_{i,j}$,  for  $ i=1, \dots l $ and for 
 $ j=1, \dots , k_i$,  is a complete 
set of coset representatives of $H$ in $P$.
Furthermore the  $Q$-orbit  $\overline{R_i}$  corresponds
  to a $Q$-orbit  $R_i = \{t_{i,1}, \dots , t_{i_{k,i}} \}$, 
for every $i=1, \dots , l$.

Let $K_S$ be the kernel of the action of $H$ on $S$.
As $|H/K_S| < |P|$, the minimality in  Inductive Assumption \ref{cc:1}
implies that there exists  $\mu  \in S$ such that its stabilizer, 
$(Q\ltimes H/K_S)(\mu)$, in $Q\ltimes H/K_S$ equals $Q$, or equivalently 
$(QH  )(\mu) =QK_S$. 
We note here that $K_S < H$. Indeed,  if $H$ acts trivially on $S$, then 
$T$ is  induced from a trivial module and thus  contains both trivial and
non--trivial irreducible 
$P$-submodules, contradicting Claim \ref{cc:c1}. We also
have that 
 $\mu \ne 0$ since $Q= (Q\ltimes H/K_S)(\mu) < QH/K_S$.
We denote by $\mu t_{i,j}$ the $t_{i,j}$-translation  
of $\mu$, for every $i= 1, \ldots, l $ and for every 
$j=1, \dots , k_i$. Then $\mu t_{i,j}$ is an element of 
 $St_{i,j}$ such that 
$$
(Q^{t_{i,j}}H)(\mu t_{i,j}) = Q^{t_{i,j}}K_S^{t_{i,j}}
$$
Since  $S^{QP} = T$  we get that
\begin{subequations}\mylabel{cc:eq1}
\begin{equation}\mylabel{cc:eq1a}
 T = S^{QP} =  \sum_{1\leq i\leq l}^{\cdot} \sum_{1\leq j \leq k_i}^{\cdot}
St_{i,j}=S \dotplus \sum_{2\leq i\leq l}^{\cdot} \sum_{1\leq j \leq k_i}^{\cdot}
St_{i,j}.
\end{equation}

Let $\tau$ be the element  of $T$ defined by 
\begin{equation}\mylabel{cc:eq2}
\tau = -\mu + \sum_{i=2}^{l} \sum_{j=1}^{k_i} \mu{t_{i,j}}=
 -\mu + \sum_{i=2}^l  \sum_{t_{i,j} \in R_i} \mu{t_{i,j} }.
\end{equation}
\end{subequations}
We claim that $\tau $ satisfies the condition in Lemma \ref{cc:l0.5}, i.e., that
$(QP)(\tau) = Q$. This will  contradict Inductive Assumption \ref{cc:1},
and thus prove Claim \ref{cc:c3}.
Indeed, $R_i = \{t_{i,1}, \dots ,t_{i,k_i} \}$ is a $Q$-orbit for every
$i=2, \dots , l$.
Also $\mu$ and  $-\mu$ are $Q$-invariant as 
$(QH)(-\mu) = (QH)(\mu) = QK_S$.
Hence $\sum_{t_{i,j} \in R_i} \mu t_{i,j}$  is   $Q$-invariant.
Thus $\tau$ is  a $Q$-invariant element of $T$. 

If $x\in  H(\tau)$ then, since $H\lhd P$, we get
that $(\mu {t_{i,j}})x =\mu x^{(t_{i,j})^{-1}}  t_{i,j}$ 
is an element of $St_{i,j}$,  for all $i=2, \dots, l$ and $j= 1,\dots,k_i$,
while $(-\mu)x$ is an 
element  of $S$.
Since $\tau x = \tau$, it follows from \eqref{cc:eq1}  that   $(-\mu)x =
-\mu$ and 
$(\mu {t_{i,j}})x =\mu {t_{i,j}}$ for every $i= 2, \dots , l$
and for every $j=1, \dots , k_i$.
Hence $x$ is an element of:
$$
H(\mu) \cap   \bigcap_{i=2}^l \bigcap_{j=1}^{k_i}(P(\mu {t_{i,j}}) \cap H) = 
\bigcap_{i=1}^l\bigcap_{j=1}^{k_i} H(\mu {t_{i,j}}) = \bigcap_{i=1}^l
\bigcap_{j=1}^{k_i} K_S^{t_{i, j}}.
$$
As $H$ acts faithfully on $T$, we get that $\bigcap_{i=1}^l
\bigcap_{j=1}^{k_i} K_S^{t_{i, j}} = {1}$.
Hence $H(\tau)  = {1}$.

Now let  $x \in P \backslash H$. We claim that $\tau x \ne \tau$.
Indeed any  $x \in P $  permutes the $St_{i, j}$ among themselves. 
If $x$ fixes $\tau$, then  it also  permutes 
among themselves the  summands  $-\mu$ and $\mu t_{i,j}$, for $i \ne 1$, of 
$\tau$. Since $Sx \ne S$ we have 
 $(-\mu) x = \mu t_{i,j}$ for some $i=2,\dots,l$ and 
some $j = 1,\dots,k_i$. But as
 $x \in P \backslash H$ we have that  $x = h t$
 for some coset representative 
$t = t_{i_0,j_0}$ of $H$ in $P$ with $i_0 = 2, \dots ,l$ 
and some element  $h \in H$.
 Hence $\mu t_{i,j}= (-\mu)x = (-\mu)h t \in S t$, which implies that 
$t_{i,j} = t$ and  $(-\mu) h = \mu$.
This last equation  leads to a contradiction  as $h$ has odd order 
($|P|$ is odd)
and $\mu \ne -\mu$ ( as $S \leq T$ has odd order).  
Therefore  $\tau^x \ne \tau $ whenever $x \in P\backslash H$.
Hence $P(\tau) = H(\tau) =1$ and $(QP)(\tau) =Q$, 
contradicting Inductive Assumption \ref{cc:1}.
This contradiction proves Claim \ref{cc:c3}.
\end{proof}

\begin{claim}\mylabel{cc:c4}
Every normal abelian subgroup $A$ of $QP$ contained in $P$ is cyclic.
\end{claim}

\begin{proof}
Let $A$ be a normal subgroup of $QP$ contained in $P$, and let $T_A$ be the
restriction of
 $T_P$ to $A$. According to  Claim \ref{cc:c2}, 
 and Clifford's Theorem (11.1 in CR), we have that 
$T_A$ can be written as a direct sum of its $\ZZ_q(A)$-homogeneous
components, i.e.,
$$
T_A = W_1 \dotplus W_2 \dotplus \dots \dotplus W_s. 
$$
Furthermore, $P$ acts transitively on  the $W_i$ for all $i=1,\dots,s$,
while $Q$ permute the $W_i$ among themselves (as $T$ is a $\ZZ_q(QP)$-module).
Hence Glauberman's lemma implies that $Q$ fixes some
$\ZZ_q(A)$-homogeneous component, $W$,  of $T_A$.
Thus $W$ is a $\ZZ_q(QA)$-submodule of $T_{QA}$. Even more, Clifford's 
Theorem implies that $W^{QP} = T$. This, along with Claim \ref{cc:c3},
implies that $W =  T$.
Hence $T_A \cong e V$  where
 $V$ is an irreducible $P$-invariant 
$\ZZ_q(A)$-submodule of $T$. As $P$ acts faithfully  on $T$, the 
$\ZZ_q(A)$-module $V$ is also faithful. If $A$ is abelian, the existence
of a faithful irreducible
 $\ZZ_q(A)$-module implies that $A$ is cyclic.
Therefore,  the claim is proved.
\end{proof}

The $q$-group $Q$ acts on the $q$-group 
$T$, fixing the trivial element  $0$ of $T$. Hence the group 
$Q$ fixes at least  $q$ elements of $T$.
So $Q$ fixes some $\tau$ with
\begin{equation}\mylabel{cc:eq3}
 \tau \in T \text{ and }  \tau \ne 0.
\end{equation}
Hence, to complete the proof of Lemma \ref{cc:l0.5}, by contradicting 
the Inductive Assumption \ref{cc:1},
  is enough to show that $P(\tau) = 1$

By Claim \ref{cc:c4} every characteristic abelian subgroup of $P$ is
cyclic. Since $p$ is odd, Theorem 4.9 in section 5.4 of \cite{go} implies
that either $P$ is cyclic or $P$
 is the central product $E \odot C$,  of an extra--special $p$-group $E$
of exponent $p$
 and a cyclic
group $C$. 
 
If $P$ is cyclic then $Z(P) = P$. 
According to Claim \ref{cc:c2}, the $\ZZ_q(P)$-module $T_P$ is a multiple of
an irreducible $Q$-invariant $\ZZ_q(P)$-module $V$, i.e.,
$T_P = m V$. Hence $Z(P)$ acts fix point free on $T$ as it acts fix 
point free on $V$ (or else $V$ wouldn't be simple).
This implies that no element of $ P = Z(P) - \{ 1\}$ could fix $\tau$. Hence 
$P(\tau) = 1$. So $(QP)(\tau) = Q$, 
contradicting Inductive Assumption \ref{cc:1}.
Therefore,  $P$ can't be cyclic. 

Hence, $P = E\odot C$, where
 $E = \Omega _1(P)$ is an extra special of exponent $p$  
 and $C = Z(P)$ is cyclic. So,
\begin{equation}\mylabel{cc:e2}
P = E \odot C = \Omega _1(P) \odot Z(P).
\end{equation}
Therefore the factor group $\overline{P} = P/ Z(P)$ is an elementary 
abelian $p$-group. Furthermore it affords a bilinear form
 $c: \overline{P} \times \overline{P} \to Z(E)$ defined,  for every 
$\bar{x}, \bar{y} \in \overline{P}$, as $c(\bar{x} , \bar{y}) = [x,y]$,
 where
$x, y$ are any elements of $P$ having images $\bar{x}, \bar{y}$ respectively, in 
$\overline{P}$.
With respect to that form $\overline{P}$ is a symplectic $\ZZ _p (Q)$-module.

\begin{claim}\mylabel{cc:c5}
The symplectic $\ZZ _p(Q)$-module $\overline{P}$ is anisotropic.
 \end{claim}

\begin{proof} 
Assume not. Then there is an isotropic non--zero $\ZZ _p(Q)$-submodule
$\bar{A}$ of $\overline{P}$. Hence $c(\bar{a} , \bar{b}) = 0$ for every
$\bar{a}, \bar{b} \in
 \bar{A}$, as $\bar{A} \subset \bar{A}^{\perp}$.
Therefore, by the definition of the symplectic
form $c$, we get that the inverse image $A$ of $\bar{A}$ in $P$ is
an abelian subgroup of $P$ containing $Z(P)$.
Since $\bar{A}$ is a $\ZZ _p(Q)$-submodule of $\overline{P}$, the abelian
group $A$ is a normal
 subgroup of $QP$ contained in $P$. 
Hence by Claim \ref{cc:c4} , $A$ is cyclic and properly
 contains the $Z(P)$.
Therefore there exists an element $a  \in A\smallsetminus Z(P)$
such that $a^p$ is a generator of $Z(P)$. On the other hand 
according to  \eqref{cc:e2} $a = \omega \cdot c $
where $\omega \in \Omega_1(P)$ and $c \in C=Z(P)$.
Hence $a^p = \omega^p \cdot c^p = c^p$. Since $a^p$ is a generator 
of the cyclic 
 $p$-group $Z(P)$ and $c \in Z(P)$, this last equation leads  to a
contradiction.
This proves the claim.
\end{proof}

Now we can complete the proof of Lemma \ref{cc:l0.5}. 
If $(QP)(\tau ) \ne Q$ then there exists a $Q$-invariant subgroup 
$D = P(\tau)\ne 1$
 of $P$ such that $(QP)(\tau) = QD$. 
Hence the center $Z(D)$ of $D$ is a non--trivial $Q$-invariant abelian subgroup of $P$.
Therefore its
 image $\overline{Z(D)} = Z(D)Z(P)/Z(P)$ in
 $\overline{P}$ is an isotropic $\ZZ _p(Q)$-submodule of $\overline{P}$. 
Since $\overline{P}$ is anisotropic, $\overline{Z(D)} = \bar{1}$, i.e., 
$Z(D)$ is contained in $Z(P)$.

As we saw in the first case, $Z(P)$ acts fix point free on $T$.
This implies that no element of $Z(P) - \{1\}$ could fix $\tau$.
Hence $Z(D) = 1$, contradicting the fact that $Z(P) \ne 1$.
 So $(QP)(\tau) = Q$, contradicting Inductive Assumption  \ref{cc:1}.
This final  contradiction completes the proof of Lemma \ref{cc:l0.5}. 
\end{proof}

In terms of characters, Lemma \ref{cc:l0.5} implies
\begin{corollary}\mylabel{cc:c0.5}
Let $Q$ be a $q$-group acting on a $p$-group $P$ with $p\ne q$ odd primes. 
Suppose that the semi--direct product $Q\ltimes P$ acts on a $q$-group $Q_1$
such that the action of $P$ on $Q_1$ is faithful. Then there exists a linear
character $\lambda$ of $Q_1$ 
whose kernel $\Ker(\lambda)$ contains the Frattini subgroup 
$\Phi(Q_1)$ and whose  stabilizer $(QP)(\lambda)$ in $Q\ltimes P$
is $Q$.
\end{corollary}

\begin{proof}
Let $T$ be the factor group $T:= Q_1/ \Phi(Q_1)$.
 Then $T$ is a $\ZZ_q(QP)$-module.
 We write  $T^*$ for  its  dual $\ZZ_q(QP)$-module, i.e.,
$T^*= \Hom_{\ZZ_q}(T,\ZZ_q)$. Then $P$ acts faithfully on both,
$T$ and $T^*$. Furthermore, according to  Lemma \ref{cc:l0.5} 
there is an element $\tau  \in T^*$ whose stabilizer in $QP$ equals $Q$.
Since the linear characters of $T$ can be considered as the elements of $T^*$
 composed with some faithful  linear character of $\ZZ_q$, we conclude that 
there is a linear character $\lambda^* \in \Lin(T)$ whose stabilizer in 
$QP$ is $Q$. 
Let $\lambda$ be the linear character of $Q_1$ to which 
$\lambda^*$ inflates. Then $\Phi(Q_1) \leq \Ker(\lambda)$.
Furthermore, $(QP)(\lambda) = (QP)(\lambda^*) = Q$, 
and the corollary follows.
\end{proof}

The following  straightforward lemma is necessary for the rest of the chapter,
and gives a stronger result than Corollary \ref{cc:c0.5}.

\begin{lemma}\mylabel{lemmaA}
Let $P$ be a $p$-subgroup of a finite group $G$ and let  $Q_1 \unlhd  Q_2$
be $q$-subgroups 
of $G$, for some distinct odd primes $p$ and $q$.
If $P$ normalizes $Q_1$, and $Q_2$ normalizes their product $Q_1P$, then 
$Q_2 P$ is also a subgroup of $G$ with  $Q_2 \in \Syl_q(Q_2 P) , 
P \in \Syl_p(Q_2 P)$  and 
$Q_1 P \unlhd Q_2 P$. Furthermore,  $Q_2$ is the product 
$Q_2 = [Q_1 , P] N(P \tin Q_2)$, where $[Q_1,P] \unlhd Q_2P$ and 
 $[Q_1, P] \cap N(P \tin Q_2) = C( P
\tin [Q_1, P]) \leq 
\Phi([Q_1, P]).$     
\end{lemma}
\begin{proof}
That the product,  $Q_2 P = Q_2 (Q_1 P)$, is a subgroup of $G$ is clear as  
$Q_2$ normalizes the semidirect product $Q_1 \rtimes P$. 
That same  product $Q_1 P$ is a normal subgroup of $Q_2 P =
Q_2(Q_1 P)$.
We obviously have that  $Q_2 \in \Syl_q(Q_2 P)$ and  $ P
\in \Syl_p(Q_2 P)$.

By Frattini's argument for the Sylow $p$-subgroup  $P$ of  $Q_1P\unlhd Q_2P$
we get
\begin{subequations}
\begin{equation}\mylabel{A1} 
Q_2 P = Q_1 P N( P \tin Q_2 P ).
\end{equation}
The normalizer, $N( P \tin Q_2 P)$, of $P$  in $Q_2 P$  contains $P$. So 
it is equal to $P N(P\tin Q_2)$. 
 Hence  \eqref{A1} can be written as  $Q_2 P = Q_1 N(P \tin Q_2 ) P$.
Since $Q_1 N(P \tin Q_2) \leq Q_2$  and $Q_2 \cap P = 1$,  we 
conclude  that 
\begin{equation}\mylabel{A2}
Q_2 = Q_1 N(P \tin Q_2).
\end{equation}
\end{subequations}
Because $(|Q_1|, |P|) = 1$, and $P$ acts on $Q_1$, we can write $Q_1$ as the
product  $Q_1 =
[Q_1 , P] N( P \tin Q_1)$.
The commutator subgroup $[Q_1, P]$ is a characteristic subgroup of $Q_1 P$
and thus is also a
normal subgroup of $Q_2$,  as  $Q_2$ normalizes $Q_1 P$.
Therefore, \eqref{A2} implies
$$Q_2 = [ Q_1, P] N(P \tin  Q_2).$$
That $ [Q_1, P] \cap N(P \tin Q_2) = C( P \tin [Q_1, P])$ is obvious as
$(|Q_1|, |P|) = 1$.
 Also
the factor group $K:= [Q_1,P]/\Phi([Q_1,P])$ is abelian and thus 
$K=[K,P] \times C(P \tin K)$. As $[Q_1,P,P] = [Q_1,P]$
(by Theorem 3.6 in section  3.5 in \cite{go}),
 we get that $K = [K,P]$ and 
$C(P \tin K)= 1$. This  implies that 
 $C(P \tin [Q_1,P]) \leq \Phi([Q_1,P])$.
\end{proof}

As an  easy consequence of Corollary  \ref{cc:c0.5} and
 Lemma \ref{lemmaA} we have:
\begin{proposition}\mylabel{cc:l1}
Let $Q$ be a $q$-group acting on a $p$-group $P$ with $p\ne q$ odd primes. 
Suppose that the semi--direct product $Q\ltimes P$ acts on a $q$-group $Q_1$
such that the action of $P$ on $Q_1$ is faithful. Then there exists a linear
character $\lambda$ of $Q_1$ 
such that  $C(P \tin Q_1 )  \leq \Ker(\lambda)$ and 
 $(QP)(\lambda) = Q$. 
\end{proposition}

\begin{proof}
As $P$ acts on $Q_1$ we can write $Q_1$ as the product $Q_1 =
[Q_1,P] \cdot C(P \tin Q_1)$. It is clear that 
the product $Q C(P \tin Q_1)$ forms a group. Furthermore, 
$Q C(P \tin Q_1)$ normalizes $P$ and 
the semidirect product $(QC(P \tin Q_1))\ltimes P$ 
acts on $[Q_1,P]$, while the action of $P$ on $[Q_1,P]$ is faithful. 
Then according to Corollary \ref{cc:c0.5} there exists a linear character
$\lambda_1$ of $[Q_1,P]$ such that 
$(QC(P \tin Q_1)P)(\lambda_1) =  Q C(P \tin Q_1)$, while 
 $\Phi([Q_1,P]) \leq \Ker(\lambda_1)$.

As we have seen in Lemma \ref{lemmaA}
$$
[Q_1,P] \cap C(P \tin Q_1) = C(P \tin [Q_1,P]) \leq \Phi([Q_1,P]).
$$
Since $\lambda_1$ is a  linear character of $[Q_1,P]$ 
that is trivial on $\Phi([Q_1,P])$ and $C(P \tin Q_1)$-invariant,
 the above inclusion implies that
 $\lambda_1$ has a unique extension to a linear character
$\lambda$ of $Q_1$ trivial on $C(P \tin Q_1)$.
Furthermore,
$(QP)(\lambda) = (QP)(\lambda_1) = Q$,
 and  the proposition follows.
\end{proof}

We are now ready to show how our first change works:

\begin{theorem}\mylabel{cc:co2}
Let $Q_1 \unlhd Q_2 = Q \leq G$ and $P$ satisfy Hypothesis \ref{cc:hyp1}.
Assume further that  $\beta_1$ is an irreducible character of $Q_1$.
 Then there exist  irreducible characters $\bn_1$ of $Q_1$ 
and $\bn$ of $Q_2(\bn_1)$ such that
\begin{align*}
P(\beta_1 ) &= P(\bn_1), \\
Q(\beta_1) &\leq Q(\bn_1) \, \text{and}\\
\bn |_{Q_1} &= \bn_1.
\end{align*}
Therefore $\bn$ is an  extension of $\bn_1$  to $Q(\bn_1)$.
\end{theorem}

\begin{proof}
Let $P(\beta_1)$ be the stabilizer of $\beta_1$ in $P$ and $P_1$ be the
normalizer of 
$P(\beta_1)$ in $P$.  Let $\overline{P_1}$ denote the factor group $P_1 /
P(\beta_1)$.
We write  $C_1$ for  the centralizer, $C_1 = C(P(\beta_1) \tin Q_1)$,  of $P(\beta_1)$
in $Q_1$. Then it is clear that $\overline{P_1}$ acts 
 on $C_1$.

The Glauberman--Isaacs correspondence (Theorem 13.1 in \cite{is}), applied to 
the groups $P(\beta_1)$ and $Q_1$, provides an irreducible character 
$\theta$ of $C_1$ corresponding to the irreducible character $\beta_1$ 
of $Q_1$.
As $P_1$ normalizes both  $P(\beta_1)$ and $Q_1$  
 we get that  $(P_1)(\theta) = (P_1)(\beta_1) = P(\beta_1)$.
If $(P_1)(\theta)  < P(\theta)$ then  $N(P(\beta_1) \tin
P(\theta)) =
 (P_1)(\theta) >(P_1)(\theta) = P(\beta_1)$.
Therefore 
$$
P(\theta) = (P_1)(\theta) = P(\beta_1).
$$ 
Since $P(\beta_1)$ centralizes $C_1 =C(P(\beta_1) \tin Q_1)$, we have 
$P(\beta_1) \leq C(C_1\tin P_1)\leq (P_1)(\theta) = P(\beta_1)$.
Hence $C(C_1\tin P_1) = P(\beta_1)$ and $\overline{P_1}$ acts faithfully 
on $C_1$. 

Let $C_2 :=N(P(\beta_1) \tin Q)$ be the normalizer of $P(\beta_1)$ in $Q$.
Then $C_1$ is a normal subgroup of $C_2$ as $Q_1 \unlhd Q_2$.
 Furthermore, $C_2$  normalizes    $N(P(\beta_1) \tin PQ_1)$ 
 as $Q_2$ normalizes  the product $PQ_1$. 
 As $P_1 C_1 = N(P(\beta_1) \tin PQ_1)$ 
 we conclude that $C_2$ normalizes the product 
$P_1 C_1$. Hence Frattini's argument implies that 
$C_2 = N(P_1 \tin C_2) C_1$.
Let $C_2' := N(P_1 \tin C_2)$.  Then $C_2'$ normalizes $\overline{P_1}$ 
and the semidirect product $C_2' \ltimes \overline{P_1}$ acts 
on $C_1$. Furthermore,
 the action of  $\overline{P_1}$
on $C_1$ is faithful.
By Proposition \ref{cc:l1}, there exists a linear character 
$\lambda \in \Lin(C_1)$  such that $C(\overline{P_1} \tin C_1)
 \leq \Ker(\lambda)$ and  $(C_2' \overline{P_1})(\lambda) = C_2'$.

The last equation implies that $P_1(\lambda) = P(\beta_1)$.
 Thus $P(\beta_1) = P_1(\lambda) \leq P(\lambda)$.
We actually have that 
 $$
P(\lambda) =P(\beta_1).
$$
Indeed,  if $P(\beta_1) < P(\lambda)$, 
 then $P(\beta_1)$ would be a proper subgroup of  
$N(P(\beta_1) \tin P(\lambda))$.
Thus $P(\beta_1) < N(P(\beta_1) \tin P(\lambda)) = N(P(\beta_1) \tin P)
(\lambda) = P_1(\lambda) = P(\beta_1)$.

Since $C_2 = C_2' C_1$ and $C_2'$ fixes $\lambda_1$, we conclude that
 $C_2$ also fixes $\lambda_1$. Furthermore,
  for the intersection $C_1 \cap C_2'$ we get 
$$
C_1 \cap C_2' = C_1 \cap N(P_1 \tin C_2) = 
C(P_1 \tin C_1) \leq C(\overline{P_1} \tin C_1) \leq \Ker(\lambda).
$$
Therefore $\lambda$ can be extended to $C_2$.
Furthermore,
according to  Theorem 6.26 in \cite{is}  and the 
fact that  $C_2 P(\beta_1)$ fixes $\lambda$, we get that  $\lambda$ can be
extended to  $C_2P(\beta_1)$.   

 Let $\bn_1 \in \Irr(Q_1)$ be the Glauberman--Issacs 
$P(\beta_1)$-correspondent  to $\lambda$.
Then  as $C_2 P_1$ normalizes both $P(\beta_1)$ and $Q_1$, Corollary \ref{cc:co1} 
implies  that 
$$
(C_2P_1)(\bn_1) = (C_2P_1)(\lambda) = C_2P(\beta_1). 
$$
Hence $P(\bn_1) \geq P_1(\bn_1) = P(\lambda) = P(\beta_1)$.
If $P(\bn_1) >  P(\beta_1)$ then  
$$
P(\beta_1)  < N(P(\beta_1) \tin P(\bn_1)) =
P_1(\bn_1)= P(\beta_1).
$$
Thus    $P(\bn_1) = P(\lambda) = P(\beta_1)$ and 
$$
(C_2(PQ_1))(\bn_1) = C_2 P(\beta_1)Q_1.
$$
Since $C_2$ fixes $\bn_1$ and normalizes $P(\beta_1)$ we have
$C_2 \leq N(P(\beta_1) \tin Q(\bn_1)) \leq N(P(\beta_1) \tin Q) = C_2$.
Hence $C_2 = N(P(\beta_1) \tin Q(\bn_1))$.
Furthermore, $P(\beta_1) Q_1 = (P_1 Q_1) (\bn_1)$ as $P(\bn_1) = P(\beta_1)$.
Hence the group $ P(\beta_1) Q_1 = (PQ_1)(\bn_1)$ is a normal subgroup
 of $P(\beta_1) Q(\bn_1)$ as
$Q$ normalizes the product $P_1 Q_1$. So we can apply the 
 Main Theorem in \cite{lew} to the groups $P(\beta_1) Q(\bn_1)$,
 $ P(\beta_1) Q_1$ and $Q_1$.
We conclude that $\bn_1$ extends  to $P(\beta_1) Q(\bn_1)$ as 
its $P(\beta_1)$-Glauberman correspondent
$\lambda$ can be extended to $ P(\beta_1) C_2 =
P(\beta_1) N(P(\beta_1) \tin Q(\bn_1))$.
We write $\bn$ for the extension of $\bn_1$ to $Q(\bn_1)$.

To complete the proof of the theorem it remains to show that 
$Q(\beta_1) \leq Q(\bn_1)$.
The group $(PQ_1) (\beta_1) = P(\beta_1) Q_1$ is a normal subgroup of 
$P(\beta_1) Q(\beta_1)$, as $Q$ normalizes $P_1 Q_1$.
Hence Frattini's argument implies that 
$$Q(\beta_1) = Q_1 N(P(\beta_1) \tin Q(\beta_1)).
$$
Therefore  $Q(\beta_1) \leq Q_1 N(P(\beta_1) \tin Q)= Q_1 C_2$.
But we have already seen that $ P(\beta_1) Q_1$ is a normal subgroup of
 $P(\beta_1) Q(\bn_1)$.
 Hence the  Frattini argument implies that 
$$
Q(\bn_1) = Q_1 N(P(\beta_1) \tin Q(\bn_1))
= Q_1 C_2.
$$ 
Thus $Q(\beta_1) \leq Q(\bn_1)$ and the theorem follows.
\end{proof}

\section{A ``multiple'' change}
In the previous section we saw how we can make a change of a character
whenever we have only
two $q$-groups, $Q_1$ and $Q_2$
 involved. The natural question that follows from that restricted case
is whether or not 
we can prove a similar theorem when a chain of $q$-groups, $Q_1 \unlhd Q_2
\unlhd \dots 
\unlhd Q_{n} \unlhd Q_{n+1}$, is involved.
What we will show is that, eventhought we can not do a character  replacement 
as in the two group case,  we can still find   
new linear characters having enough of the desired properties.

So assume that, along with the  series $Q_1 \unlhd Q_2 \unlhd \dots \unlhd
Q_{n} \unlhd Q_{n+1}=Q$ of normal subgroups of $Q_{n+1}= Q$, 
where  $Q_{n+1}$ is a $q$-subgroup of 
a finite group $G$,  we also have $p$-subgroups
 $P_1, P_2, \dots  ,P_{n-1}, P_n$ of $G$,
 such that $P_i$ normalizes  the groups $P_j$
and $Q_j$ 
whenever $1\leq j \leq i \leq n$, while $Q_i$ normalizes the semidirect
 product $P_j \ltimes Q_j$ whenever $1\leq j < i\leq n+1$.
Assume further that $K_i = C(Q_i \tin P_i) $ for all $i=1, \dots, n$. 
Then, as we will see by the end of the section, we can find 
 a chain  of  linear  characters,
 $\bn_1, \bn_2, \dots ,\bn_n$, of $Q_1, Q_2, \dots ,Q_n$  respectively, such
that 
\begin{itemize}
\item[(1)] $P_i(\bn_i) =K_i$ and $Q_{n+1}(\bn_i)= Q_{n+1}$  while
\item[(2)] $\bn_{n}$ can be extended to $Q_{n+1}(\bn_{n})= Q$.
\end{itemize}
The  equivalent  of the diagram  \ref{cc:d1}  in this case is:
\begin{equation}\mylabel{cc:d2}
\begin{diagram}[small]           
        &&Q=Q_{n+1}= Q(\bn_n)&\qquad   &             &       &(\bn_n)^e   & \\
        &\ldTo     &\dTo     &\qquad   &             &       &          \\
P_{n}   &          &         &\qquad   &             &       &           \\
        &\rdTo     &         &\qquad   &             &       &\vLine    \\ 
        &          & Q_n     &\qquad   &             &       &\bn_n     \\
        &\ldTo     &         &\qquad   &             &       &          \\
P_{n-1} &          &         &\qquad   &             &       &\vLine    \\
        &\rdTo     &\dTo     &\qquad   &             &       &          \\
        &          &Q_{n-1}  &\qquad   &             &       &\bn_{n-1} \\
        &\ldDashTo &         &\qquad   &             &       &\ddash    \\
P_2     &          &\dDashTo &\qquad   &             &       &          \\
        &\rdTo     &         &\qquad   &             &       &          \\ 
        &          &  Q_2    &\qquad   &             &      &\bn_2     \\
        &\ldTo     &         &\qquad   &             &       &          \\
P_1     &          &         &\qquad   &             &       &\vLine    \\
        &\rdTo     &\dTo     &\qquad   &             &       &          \\
        &          &Q_1      &\qquad   &             &       &\bn_1      \\
\end{diagram}
\end{equation}
We will prove in this section that such a generalization is possible when 
the primes $p,q$ are odd.

\noindent
We first need some  lemmas:
\begin{lemma}\mylabel{lemmaB}
Assume $Q_1 , Q_2$ are two $q$-subgroups of a finite group $G$, for some odd
prime $q$.
Assume further that $Q_1 \unlhd Q_2$ and that $A$ is a normal subgroup of
$Q_1$ normalized 
by $Q_2$ as well.
Let $P$ be a $p$-subgroup of $G$, where $2\ne p\ne q$, such that 
\begin{itemize}
\item[(a)] $P$ normalizes $A$ and $Q_1$, and
\item[(b)]$Q_2$ normalizes   $ Q_1 P$  
\end{itemize} 
If $\pw$ denotes the centralizer, $C(Q_1/A  \tin P)$, of $Q_1/ A$ in $P$, then 
$Q_2$ normalizes the semidirect product $\pw \ltimes A$. Thus $Q_2 \pw$ 
is a group with $Q_1 \pw$ and $A$ as normal subgroups.
The factor group $ Q_2 \pw /A =Q_2(\pw \ltimes A)/A$  equals the semidirect product
$(Q_2/A) \,\ltimes \,((\pw A)/A)$, of its $q$-subgroup $Q_2/A$ and its normal
 $p$-subgroup
$(\pw A)/A \cong \pw$.
\end{lemma} 

\begin{proof}  
Since $Q_2$ normalizes $Q_1 P$, the product  $Q_2 P = Q_2 \cdot (Q_1 P)$ is a group.
Note that $A \unlhd Q_2 P$.
We will use bars to denote the image in $\overline{PQ}_2 = (PQ_2)/A$ of any
subgroup of 
$PQ_2$ containing $A$. 
Then $\pw = C(\overline{Q}_1 \tin P)$, is a normal subgroup of $P$, as the latter
normalizes both
$A$ and $Q_1$.
Furthermore, $\overline{PQ}_1 = \overline{P} \ltimes \overline{Q}_1$ and $\overline{\pw}$ is  the 
centralizer in $\overline{P}$ of $\overline{Q}_1$.  It follows that $\overline{\pw} = (\pw
A)/A$  is
 the maximal normal $p$-subgroup of $\overline{PQ}_1$. Therefore, $\overline{Q}_2$
normalizes 
$\overline{\pw}$ as it normalizes $\overline{PQ}_1$. Hence $Q_2$ normalizes the inverse
image $\pw A$
of $\overline{\pw}$ in $PQ_2$. So $\pw Q_2 = (\pw A)Q_2$ is a group  and thus a
subgroup of $PQ_2$.
Furthermore, as $Q_1 \unlhd Q_2$ and $Q_2$ normalizes $\pw A$ we get that 
$Q_2$ normalizes their product $Q_1 (\pw A) = Q_1 \pw$. Therefore $Q_1 \pw \unlhd Q_2 \pw$.

As  $Q_2$ normalizes $\pw A$ we get  that $Q_2/A$ normalizes $\pw A/A$.
Therefore $(Q_2/A) \ltimes  (\pw A/A)$ is a group. Clearly 
$Q_2(\pw A)/A = (Q_2/A) \ltimes (\pw A/A)$, 
and the lemma follows.
\end{proof}

\begin{lemma}\mylabel{lemmaC}
Assume that $Q_1,Q_2, G$ and $P$ satisfy Hypothesis \ref{cc:hyp1}.
Let  $K := C(Q_1 \tin P)$ be the kernel of the $P$-action  on $Q_1$.
Then $Q_2 = [Q_1, P] \cdot N(P \tin Q_2)$
 and $Q_1 = [Q_1,P]\cdot C(P \tin Q_1)$.
 Furthermore, there exist linear characters
$\lambda_1 \in \Lin(Q_1)$ and 
$\lambda_2 \in \Lin(Q_2)$ that satisfy the following three conditions: 
\begin{itemize}
\item[(1)]   $\lambda_2 |_{Q_1} = \lambda_1$.
\item[(2)] $(Q_2\cdot P)(\lambda_1) =Q_2\cdot K$.
\item[(3)] $C(P \tin Q_1) \leq \Ker(\lambda_1)$ and $N(P \tin Q_2) 
\leq \Ker(\lambda_2)$.
\end{itemize}
Therefore, $\lambda_1(s\cdot u ) = \lambda_1(s)= \lambda_1|_{[Q_1,P]}(s)$ and 
 $\lambda_2(s \cdot t) = \lambda_1(s)$, for all 
$s \in [Q_1,P], u \in 
C(P \tin Q_1)$ and $t \in N(P \tin Q_2)$.
\end{lemma}

\begin{proof}
Since $Q_2$ normalizes the semidirect product $Q_1 \rtimes P$ and $Q_1 \unlhd Q_2$, 
Lemma \ref{lemmaA} implies 
\begin{equation}\mylabel{cc:eC1}
Q_2 = Q_1 N(P\tin Q_2) =  [Q_1, P] \cdot N(P \tin Q_2),
\end{equation}
where  $[Q_1,P]\unlhd Q_2 P$,(and thus $[Q_1, P] \unlhd Q_2$). Also 
$N(P\tin Q_2 ) \cap Q_1 = C(P \tin Q_1)$.
Furthermore, as the $p$-group $P$  normalizes the $q$-group $Q_1$, we have
\begin{equation}\mylabel{cc:eC2}
Q_1 = [Q_1,P] \cdot C(P \tin Q_1).
\end{equation}

Let $C_2 = N(P \tin Q_2)$.
Then $C_2$  normalizes $P$, while their semidirect product
$C_2 \ltimes P$ normalizes $Q_1$.

\emph{Case 1}: Assume that $K = {1}$.

Since $K=1$, the $p$-group $P$ acts faithfully on $Q_1$.
Therefore, in view of Proposition \ref{cc:l1},
 there exist  a linear character
$\lambda_1 \in \Lin(Q_1)$ suct that 
$(C_2 P) (\lambda_1) = C_2$, and $C(P \tin Q_1) \leq \Ker(\lambda_1)$.

Since $Q_2= Q_1 C_2$ (according to \eqref{cc:eC1}) and 
$C_2$ fixes $\lambda_1$, we conclude that $\lambda_1$ is $Q_2$ invariant. 
Furthermore, the fact that  $Q_1 \cap C_2 = C(P\tin Q_1)
 \leq \Ker(\lambda_1)$,
implies that  $\lambda_1$ extends canonicaly
 to a linear character $\lambda_2$ of $Q_2$ such that 
$N(P \tin Q_2)\leq  \Ker(\lambda_2)$.
This along with \eqref{cc:eC1} and \eqref{cc:eC2} imply that 
$$
 \lambda_1(s\cdot u ) = \lambda_1(s) =\lambda_1|_{[Q_1,P]}(s),
$$
while
$$ \lambda_2(s \cdot t) = \lambda_2(s) = \lambda_1(s),
$$
for all $s \in [Q_1,P],  u \in C(P \tin Q_1)$ and $t \in N(P \tin Q_2)$.

This completes the proof of Lemma \ref{lemmaC} when $K = C( Q_1 \tin P) = 1$.

\emph{Case 2}: Assume that $1 < K \leq P$.

\noindent 
In this case we work with the group $P' = P/K$ in the place of $P$.
Note that in  view of \eqref{cc:eC2}, we have
$$K = C(Q_1 \tin P) = C([Q_1,P] \tin P).$$
As $K= O_p(PQ_1)$ and $PQ_1 \unlhd PQ_2$ we have that $K$ is a normal
 subgroup  of $PQ_2$.
Hence $K$ is a $Q_2$-invariant subgroup of $P$.
Therefore
the hypothesis of Lemma \ref{lemmaC} are satisfied for $P'$ in
 the place of  $P$,
  as $P'$ normalizes $Q_1$ while $Q_2$ normalizes their product
$Q_1\rtimes P'$.

\noindent
Hence the previous case provides linear  characters
$\lambda_1 \in \Lin(Q_1)$ and $\lambda_2 \in \Lin(Q_2)$ such that
$\lambda_2$ is an extension of $\lambda_1$ to $Q_2$. Furthermore,
\begin{subequations}\mylabel{cc:eC5}
\begin{equation}\mylabel{cc:eC5a}
(Q_2\cdot P')(\lambda_1) = Q_2
\end{equation}
 while 
\begin{equation}\mylabel{cc:eC5b}
C(P' \tin Q_1)\leq \Ker(\lambda_1) \text{ and }
 N(P' \tin Q_2) \leq \Ker(\lambda_2).
\end{equation}
\end{subequations}
Even more, 
 $\lambda_1(s\cdot u ) = \lambda_1(s)$ and 
 $\lambda_2(s \cdot t) = \lambda_1(s)$, for all $s \in [Q_1,P],  u \in 
C(P' \tin Q_1)$ and $t \in N(P' \tin Q_2)$.

We observe that  \eqref{cc:eC5a} implies
$$
(Q_2 \cdot P)(\lambda_1) = Q_2 \cdot K,
$$
 as  $K$ centralizes $Q_1$. 
Furthermore, we note that $[Q_1,P'] = [Q_1,P]$  while
$N(P \tin Q_2) \leq N(P' \tin Q_2)$ and  $C(P \tin Q_1) \leq C(P' \tin Q_1)$.
Hence in view of \eqref{cc:eC5b} we have
 $N(P \tin Q_2) \leq \Ker(\lambda_2)$ and $C(P \tin Q_1) \leq 
\Ker(\lambda_1)$. 
Therefore the lemma follows.
\end{proof}

\begin{theorem}\mylabel{cc:p1}
Assume $G$ is a finite group of order $p^a q^b$  for distinct odd primes
$p$ and $q$, and 
non--negative integers $a$ and $b$. Let  $Q_1 \unlhd Q_2 \unlhd \dots \unlhd
Q_{n} \unlhd Q_{n+1} = Q$ be a series of normal subgroups of
a $q$-subgroup  $Q \leq G$, and let 
$P_1, P_2, \dots  ,P_{n-1}, P_n$ be $p$-subgroups of $G$  such that  the
following hold:
\begin{itemize}
\item[(1)] $P_i$ normalizes  the groups $P_j$ and $Q_j$ 
whenever $1\leq j \leq i \leq n$, while
\item[(2)] $Q_i$ normalizes the semidirect product $P_j \ltimes Q_j$
whenever $1 \leq j < i\leq n+1$.
\end{itemize}
Let  $K_i$  denote the kernel of the $P_i$-action on $Q_i$, i.e., $K_i =
C(Q_i \tin P_i)$ for every 
$i= 1, \dots n$. 
Then there exist linear characters $\beta_i$ of $Q_i$, for all $i = 1, \dots
,n+1$, such that: 
\begin{itemize}  
\item[(a)] the restriction $\beta_i|_{Q_j}$ of $\beta_i$ to $Q_j$ equals
$\beta_{j}$ if $1\leq j \leq i \leq n+1$, and 
\item[(b)] the stabilizer $(Q P_i)(\beta_i)$ of $\beta_i$ in $Q P_i$ 
equals $QK_i$ if $1\leq i \leq n$.
\end{itemize}
Thus $\beta_{n+1}$ is an extension to $Q = Q_{n+1}$ of $\beta_1, \beta_2, 
\dots  ,\beta_n$.
\end{theorem}

\begin{proof}
For the proof we will use induction on $n$.
The case $n =1$ is done in Lemma \ref{lemmaC}.

We assume   that the proposition holds for all $n$ with $1\leq n <k$ 
and  some  $k\geq  2$. We will prove it also holds for $n =k$.
Since $Q_i$ normalizes $P_1\ltimes Q_1$, for all $i=1,\dots,k+1$,
 while $Q_1 \unlhd Q_i$, Lemma \ref{lemmaA}
implies that $Q_i$ is the  product 
\begin{subequations}\mylabel{cc:p4}
\begin{equation}\mylabel{cc:p4a}
Q_i = Q_1\cdot N(P_1 \tin Q_i),
\end{equation}
of its normal subgroup $Q_1$ with $N(P_1 \tin Q_i)$, where
\begin{equation}\mylabel{cc:p4b}
Q_1\cap N(P_1 \tin Q_i) \leq C(P_1 \tin Q_1).
\end{equation}
\end{subequations}

As $P_i$ normalizes the groups $Q_i, Q_1$,
 it also normalizes the factor group  $\overline{Q}_{i,1}:= Q_i/Q_1$, 
whenever $1\leq i \leq k$.

\noindent
We define:
\begin{equation}\mylabel{cc:p6}
P_{i,1}=C(\overline{Q}_{i,1} \tin P_i).
\end{equation}
Note that $P_{1,1}= C(Q_1/Q_1 \tin P_1)= P_1$.

If we  apply Lemma \ref{lemmaB} to the groups $Q_i = Q_1\cdot 
N(P_1 \tin Q_i), Q_1, Q$ and $P_i$, for some $i=1,\dots,k$,
 in the place of $Q_1, A, Q_2$ and $P$ 
respectively, we conclude that 
\begin{equation}\mylabel{cc:p9a}
Q \text{ normalizes the semidirect product }
 P_{i,1}\ltimes Q_1. 
\end{equation}

\noindent
Furthermore,  the group $P_{i,1}$  normalizes $Q_1$, 
as $P_i$ does. Since   $P_i$ normalizes both $P_j$ and 
 $Q_j$ for all $j=1,\dots,i$  we have that $P_{i,1}$   normalizes
both $P_j$ and $Q_j$ as well. Therefore $P_{i,1}$ normalizes both
 the factor group 
$\overline{Q}_{j,1}$ and the centralizer $C(\overline{Q}_{j,1} \tin P_j)=  P_{j,1}$,  for all such $j$.
Hence  the product
\begin{equation}\mylabel{cc:p7}
\map_1 = P_{1,1}\cdot P_{2,1}\cdots P_{n,1}= P_1\cdot P_{2,1}\cdots P_{n,1}
\end{equation} 
is a $p$-subgroup of $G$ that normalizes $Q_1$. Thus $\map_1 \ltimes Q_1$ 
is a group.
In view of \eqref{cc:p9a}, the group $\map_1 \ltimes Q_1$ 
is normalized by $Q = Q_{k+1}$. 
Let $C_1:=C(Q_1 \tin \map_1)$ be the centralizer of $Q_1$ in 
$\map_1$. Then Lemma \ref{lemmaC} implies that there exists a linear character 
$\mu_1 \in \Lin(Q_1)$ that can be extended to a linear character 
 $\mu_1^e  \in \Lin(Q)$, with the following properties:
\begin{subequations}\mylabel{cc:p14}
\begin{equation}\mylabel{cc:p14a}
Q(\mu_1) = Q,
\end{equation}
\begin{equation}\mylabel{cc:p14b}
\map_1(\mu_1) = C_1 = C(Q_1 \tin \map_1)
\end{equation}
and
\begin{equation}\mylabel{cc:p14c}
C(P_1 \tin  Q_1)\leq \Ker(\mu_1).
\end{equation}
\end{subequations}
 Furthermore, for the extension character $\mu_1^e$, we have:
\begin{subequations}\mylabel{cc:p16}
\begin{equation}
\mu_1^e|_{Q_1} = \mu_1,
\end{equation}
while
\begin{equation}\mylabel{cc:p16a}
\mu_1^e(s\cdot t) = \mu_1(s) 
\end{equation}
for all $s \in Q_1$ and $t \in N(P_1 \tin Q_{k+1})$.
Clearly for all $i=1,\dots,k$ we have that $\mu_1^e|_{Q_i}\in \Lin(Q_i)$.
Furthermore \eqref{cc:p16a}, along with \eqref{cc:p4a}, implies 
\begin{equation}\mylabel{cc:p16b}
\mu_1^e|_{Q_i}(s\cdot t_i) = \mu_1(s),
\end{equation}
for all $s \in Q_1$ and $t_i \in N(P_1 \tin Q_i)$.
\end{subequations}

We will use our inductive argument on the groups
$$Q_2/Q_1 \unlhd Q_3/Q_1\unlhd \dots \unlhd Q_{k+1}/Q_1
 = Q/Q_1.$$ 
Note that the above groups form a series of normal subgroups of
the $q$-group  $Q/Q_1$,
as $Q_2 \unlhd Q_3 \unlhd \dots \unlhd Q = Q_{k+1}$ is a normal series of 
$Q$. Furthermore the group $P_i$ normalizes  $Q_j/Q_1$,
whenever $1\leq j \leq i\leq k$,
 as $P_i$ normalizes both  $P_1$ and $Q_j$.
Thus $P_i \ltimes Q_i/Q_1$ is a group.
 Also $Q_i/Q_1$  normalizes the semidirect product 
 $P_j \ltimes (Q_j/Q_1)$, whenever $1\leq j \leq i \leq k+1$,
 as $Q_i$ normalizes the semidirect product
$P_j\ltimes Q_j$.
Hence by induction, there exist linear characters  $\lambda_2^* \in \Irr(Q_2/Q_1),
\dots,\lambda_{k}^* \in \Irr(Q_k/Q_1)$, and $\lambda_{k+1}^* \tin 
\Irr(Q/Q_1)$ such that 
\begin{equation}\mylabel{cc:p11}
\lambda_{k+1}^*|_{Q_i/Q_1} = \lambda_i^*
\end{equation}
and
\begin{equation}\mylabel{cc:p17}
(Q/Q_1 \cdot P_i )(\lambda_i^*) = Q/Q_1 \cdot  C(Q_i/Q_1 \tin P_i)
= Q/Q_1 \cdot P_{i,1},
\end{equation}
 for all $i=2,\dots,k$.

\noindent
Let $\lambda_i \in \Lin(Q_i)$ be the linear character of $Q_i$ inflated from
$\lambda_i^* \in \Lin(Q_i/ Q_1)$.
Then \eqref{cc:p11} and \eqref{cc:p17} imply:
\begin{subequations}\mylabel{cc:p18}
\begin{equation}\mylabel{cc:p18a}
\lambda_{k+1}|_{Q_i} = \lambda_i,
\end{equation}
\begin{equation}\mylabel{cc:p18b}
Q_1 \leq \Ker(\lambda_i),
\end{equation}
and
\begin{equation}\mylabel{cc:p18c}
Q(\lambda_i) = Q,
\end{equation}
for all $i=2,\dots,k+1$. Furthermore,
\begin{equation}\mylabel{cc:p18d}
 P_i(\lambda_i)=P_{i,1}
\end{equation}
for all $i=2,\dots,k$.
\end{subequations}

As $\lambda_{k+1}$ is a linear character of $Q_{k+1}$ and $Q_1 \leq Q_{k+1}$,
the restriction $\lambda_1:= \lambda_{k+1}|_{Q_1}$, is a linear character of $Q_1$.
Furthermore, \eqref{cc:p18b} implies
 \begin{equation}\mylabel{cc:p19}
\lambda_1 = 1_{Q_1} \text{ and thus \, } Q(\lambda_1) = Q.
\end{equation}

Since $Q_i = Q_1 \cdot N(P_1 \tin Q_i)$ for every $i =1,\dots,k+1$, 
(see \eqref{cc:p4a}),  equations \eqref{cc:p18a} and 
 \eqref{cc:p19} imply
\begin{equation}\mylabel{cc:p20}
\lambda_i(s\cdot t_i) = \lambda_i(t_i) = \lambda_{k+1}(t_i),
\end{equation}
for all $s \in Q_1$ and $t_i \in N(P_1 \tin Q_i)$.

Using the equation \eqref{cc:p4a}, we define for all $i=1,\dots,k+1$, 
\begin{equation}\mylabel{cc:p21}
\beta_i(s\cdot t_i) := \mu_1(s) \cdot \lambda_i(t_i),
\end{equation}
whenever  $s\in Q_1$ and  $t_i \in N(P_1 \tin Q_i)\leq N(P_1 \tin Q_{k+1})$. 
According to \eqref{cc:p16b} and \eqref{cc:p20}, we can rewrite the characters $\beta_i$  as
$$\beta_i = \mu_1^e|_{Q_i} \cdot \lambda_i = 
\mu_1^e|_{Q_i} \cdot \lambda_{k+1}|Q_i = 
(\mu_1^e \cdot \lambda_{k+1})|_{Q_i}.$$
Hence  
\begin{equation}\mylabel{cc:p22}
\beta_i|_{Q_j} = \beta_j,
\end{equation}
whenever $1\leq j \leq i \leq k+1$.
Therefore, $\beta_{k+1}$ is an extension of $\beta_1$ to $Q$.

As $Q$ fixes $\mu_1$ by \eqref{cc:p14a}, and fixes $\lambda_i$ by \eqref{cc:p18c},
it also  fixes  $\beta_i$, in view of  \eqref{cc:p21}.
Furthermore, \eqref{cc:p21} implies that
 $P_i(\beta_i) = P_i(\mu_1) \cap P_i(\lambda_i)$  for all $i=1,\dots,k$.
In view of \eqref{cc:p14b}  and \eqref{cc:p18d} we conclude that 
$$
P_i(\beta_i)= C(Q_1 \tin P_i) \cap P_{i,1}.
$$
But $P_{i,1}= C(Q_i/Q_1 \tin P_i)$ by  \eqref{cc:p6}.
Hence  
$$
P_i(\beta_i)= C(Q_1 \tin P_i) \cap C(Q_i/Q_1 \tin P_i) =
C(Q_i \tin P_i) = K_i.
$$
This completes the proof of the inductive step for $n=k$. Thus 
Proposition \ref{cc:p1} follows.
\end{proof}

%%% Local Variables: 
%%% mode: latex
%%% TeX-master: "lemma1"
%%% End: 

\chapter{ Triangular Sets }\mylabel{pq}  

\noindent     
\section{The correspondence }\mylabel{pq:sec1}
Assume we have the following situation:
\begin{hyp}\mylabel{hyp1}
Let $G$ be an odd order group and $\pi$ any set of  primes.
  Let 
\begin{equation}\mylabel{*}
1=G_0\unlhd G_1 \unlhd \dots \unlhd G_{m} \unlhd G
\end{equation}
 be a series of normal subgroups  of $G$, for some arbitrary
integer $m >0$,  such that  $G_i / G_{i-1}$, for $i=1,2,\dots,m$,
  is a $\pi'$-group when $i$ is odd and  a $\pi$-group when $i$ 
 is even. 
\end{hyp}

Recall  the definition given in Chapter \ref{intro} 
\begin{defn}
Let  $\chi  _i $ be an irreducible character of $G_i$,  for all $i=0,1,\dots,m$,
 such that $\chi _i$  lies over  $ \chi_k$
 for all $k=0,1, \dots  ,i$.  
Any such collection of irreducible  characters  
$\{\chi_i \}_{i=0}^m$ 
is said to be a { \em character tower }  for the series 
 $\{ G_i \}_{i=0}^m$. 
\end{defn}

Suppose further that there exist
 $\pi$- and $\pi'$-groups $P_{2r}$ and $Q_{2i-1}$ respectively,
 along with irreducible characters $\alpha_{2r}$ and $\beta_{2i-1}$, 
such that 
\begin{equation}\mylabel{x}
\begin{gathered}
P_0=1 \text{ and } \alpha_0=1,\\
Q_1=G_1 \text{ and } \beta_1=\chi_1, \\
P_{2r} \in \Hall_{\pi}(G_{2r}(\alpha_2,\dots,\alpha_{2r-2},\beta_1,
\dots,\beta_{2r-1})), \\
\alpha_{2r} \in \Irr(P_{2r}) \text{  lies 
above the $Q_{2r-1}$-Glauberman correspondent
of  $\alpha_{2r-2}$}, \\
Q_{2i-1} \in \Hall_{\pi'}(G_{2i-1}(\alpha_2,\dots,\alpha_{2i-2},
\beta_1,\dots,\beta_{2i-3})), \\
\beta_{2i-1} \in \Irr(Q_{2i-1}) \text{ lies above the $P_{2i-2}$-Glauberman
 correspondent  of $\beta_{2i-3}$  }
\end{gathered}
\end{equation}
for all odd  $2i-1$ and even $2r$  with $1< 2i-1\leq m$ 
and $1<2r \leq m$.

Depending on the parity of $m$ 
 the collection of  groups and characters appearing in \eqref{x}
 consists of 
\begin{align*}
\{P_0,P_2, \dots, P_{m-1}, Q_1,Q_3,\dots ,Q_{m}| \alpha_0, \dots,\alpha_{m-1},\beta_1,\dots,\beta_{m}\} \text{ if $m$ is odd, and  } \\ 
\{P_0,P_2,\dots,P_m,Q_1,Q_3,\dots,Q_{m-1}|\alpha_0,\dots,\alpha_m,\beta_1,\dots,\beta_{m-1}\} \text{ if $m$ is even. }
\end{align*}
By convention we will write, for both cases, the collection of groups and characters 
 as
\linebreak 
 $\{P_{2r},Q_{2i-1}|\alpha_{2r},\beta_{2i-1}\}$
where $1\leq 2i-1  \leq m$ and $0\leq 2r \leq m$. 
We also write $\{Q_{2i-1}|\beta_{2i-1}\}$
for the $\pi'$-subset of the above collection and similarly, 
$\{P_{2r} |\alpha_{2r}\}$ for the $\pi$-subset.
\begin{defn}\mylabel{trian}
Any  set of groups and characters 
 that satisfies \eqref{x} will be called  a { \em   triangular set }
  for   the normal series \eqref{*}.
\end{defn}

At this point it is not clear at all that such a collection of 
groups and characters exists. Even worse, it is not at all obvious  that 
\eqref{x} is well defined. For all we know,  the group $C(P_{2i}
\tin Q_{2i-1})$, which is the support of 
 the $P_{2i}$-Glauberman correspondent 
of $\beta_{2i-1}$, need not be a subgroup of $Q_{2i+1}$. Thus, 
 to ask for the  character $\beta_{2i+1} \in \Irr(Q_{2i+1})$ to lie 
above a character of $C(P_{2i} \tin Q_{2i-1})$ seems out of place.
(Of course the same problem appears for the $\pi$-groups
 $C(Q_{2r-1} \tin P_{2r-2})$
and the character $\alpha_{2r}$). 
But in fact these collections of groups and characters do exist, as we will 
see in Section \ref{pq:sectrian}.  
Furthermore, we  will prove in the rest of this chapter  not only that 
triangular sets exist but also that they  correspond uniquely, up to conjugation, 
to character towers.  In particular we will prove

\begin{theorem}\mylabel{cor:t}
Assume that Hypothesis \ref{hyp1} holds.
Then there is a one--to--one correspondence between
$G$-conjugacy classes of 
character towers  of \eqref{*} and $G$-conjugacy classes of triangular sets 
for \eqref{*}.
\end{theorem}

\section{Triangular-sets: existence and  properties }
\mylabel{pq:sectrian}

Assume that  a finite group $G$ and a normal series 
\eqref{*} are given so that Hypothesis \ref{hyp1} is satisfied.
Recall  (see Chapter \ref{intro})  that for any real number $x$, we denote by 
  $[x]$  the greatest integer $n$ such that $ n \leq x$. 
If we write
\begin{subequations}\mylabel{kl:def}
\begin{align}
l&=[(m+1)/2],text{ and }  \\
k&=[m/2], 
\end{align}
\end{subequations}
then  $2l-1$ is the greatest odd integer  in the set $\{1,\dots,m\}$ while 
$2k$ is the greatest even integer in the same set.
Furthermore,
\begin{equation}\mylabel{kl}
k\leq l \leq k+1,
\end{equation}
where for $m$ even we get  $k=l$,  while  for $m$ odd we have  $k= l-1$.
Then it is easy to construct, in a recursive way, 
 a collection of groups and characters
$Q_{2i-1}, P_{2r} ,\beta_{2i-1}$ and $\alpha_{2r}$ so
 that the following holds: 
\begin{subequations}\mylabel{xx}
\begin{gather}
P_0=1 \text{ and } \alpha_0=1,\\
Q_1=G_1 \text{ and } \beta_1=\chi_1, \\
P_{2r} \in \Hall_{\pi}(G_{2r}(\alpha_2,\dots,\alpha_{2r-2},\beta_1,
\dots,\beta_{2r-1})),\mylabel{pq13a} \\
\alpha_{2r} \in \Irr(P_{2r}), \mylabel{pq13aa} \\
Q_{2i-1} \in \Hall_{\pi'}(G_{2i-1}(\alpha_2,\dots,\alpha_{2i-2},
\beta_1,\dots,\beta_{2i-3})), \mylabel{pq14a} \\
\beta_{2i-1} \in \Irr(Q_{2i-1}), \mylabel{pq14aa}
\end{gather}
\end{subequations}
whenever $2\leq i \leq l$ and $1\leq r \leq k$.

Notice that \eqref{xx}  is   a part of  \eqref{x}.
So to prove existence of triangular sets we  need to show that 
the characters 
$\beta_{2i-1}$ and $\alpha_{2r}$ that appear in \eqref{xx}
 can be chosen to satisfy  the additional conditions required in \eqref{x}.
Before we prove this,
let us see what conclusions we can draw from \eqref{xx}.
The first obvious remark is that, according to 
\eqref{pq14a} and \eqref{pq13a}, we have:
 \begin{subequations}\mylabel{pq14'}
\begin{equation}\mylabel{pq14'a}
 Q_{2i-1} \text{ normalizes the groups } P_2, \dots ,P_{2i-2}, Q_1, \dots 
,Q_{2i-3},
\end{equation}
while 
\begin{equation}\mylabel{pq14'b}
P_{2r} \text{ normalizes the groups } P_0, P_2, \dots ,P_{2r-2}, Q_1, \dots,
 Q_{2r-1}, 
\end{equation}
whenever $2\leq i \leq l$ and  $1\leq r \leq k$.
\end{subequations}

Since $P_{2r}$ normalizes $Q_{2r-1}$, the group 
\begin{equation}\mylabel{anno1}
Q_{2r-1, 2r}:=C(P_{2r} \tin Q_{2r-1} ) = N(P_{2r} \tin Q_{2r-1})
\end{equation}
is defined whenever  $1\leq r \leq k$ (note that  the group 
$Q_{2r-1}$  is defined for all such $r$ as $k\leq l$).
Furthermore,  \eqref{xx} implies two lemata that  lead to the existence 
of triangular sets. We start with  
\begin{lemma}\mylabel{yq}
For every $i$ with $1\leq  i \leq l-1$ we have 
$$
Q_{2i-1, 2i} = Q_{2i+1} \cap G_{2i-1}.
$$ 
Hence $Q_{2i-1,2i}$ is a normal subgroup of $Q_{2i+1}$.
\end{lemma}
Note that this lemma gives  us no information about  the group
$Q_{2k-1,2k}$ in the case of an even $m= 2l= 2k$, as the group 
$Q_{2k+1}$ is not  defined in that case.

\begin{proof}
The group $G_{2i-1}$ is a normal subgroup of $G_{2i+1}$ whenever 
$1\leq i <l$. Hence the definition  \eqref{pq14a} of $Q_{2i+1}$ 
implies that 
\begin{equation}\mylabel{yq:e}
Q_{2i+1} \cap G_{2i-1} \in \Hall_{\pi'}(G_{2i-1}(\alpha_2, \dots,
\alpha_{2i}, \beta_{1}, \dots,\beta_{2i-1})).
\end{equation}
In particular, whenever  $1< i <l$ the character $\beta_{2i-3}$ is defined, and
we have 
$$
Q_{2i+1}\cap G_{2i-1} \leq  G_{2i-1}(\alpha_2, \dots,
\alpha_{2i}, \beta_{1}, \dots,\beta_{2i-1})\leq G_{2i-1}(\alpha_2, \dots,
\alpha_{2i-2}, \beta_{1}, \dots,\beta_{2i-3}).
$$
Furthermore, $Q_{2i+1}$  normalizes $Q_{2i-1}$ (according to  \eqref{pq14'a}).
Also  $Q_{2i-1}$ is a $\pi'$-Hall subgroup  of the group 
$G_{2i-1}(\alpha_2, \dots,\alpha_{2i-2}, \beta_{1},
\dots,\beta_{2i-3})$.
So the intersection $Q_{2i+1} \cap G_{2i-1}$ is a $\pi'$-subgroup of $G_{2i-1}(\alpha_2, \dots,
\alpha_{2i-2}, \beta_{1},\dots,\beta_{2i-3})$,  and  normalizes the $\pi'$-Hall subgroup 
$Q_{2i-1}$ of that group. Therefore
$Q_{2i+1} \cap G_{2i-1} \leq Q_{2i-1}$ whenever $1<i <l$.
But the last inclusion  is still valid when $i = 1$, as $Q_1 =G_1$ and therefore
$Q_3 \cap G_1 \leq Q_1$. Hence, $Q_{2i+1} \cap G_{2i-1} \leq Q_{2i-1}$
 whenever $1\leq  i  < l$.
As $Q_{2i+1}$ normalizes the group $P_{2i}$ by \eqref{pq14'a},
we conclude that $Q_{2i+1}\cap G_{2i-1}$ is a subgroup of 
$ N(P_{2i} \tin Q_{2i-1})=Q_{2i-1,2i}$, i.e., that  
$$
Q_{2i+1} \cap G_{2i-1} \leq Q_{2i-1,2i}
$$
whenever $1\leq i<l$.

To prove the opposite inclusion we remark that, as 
$Q_{2i-1,2i}$ centralizes $P_{2i}$,
 it fixes the character 
$\alpha_{2i}$. It also fixes the character $\beta_{2i-1} \in \Irr(Q_{2i-1})$, as it is 
a subgroup of $Q_{2i-1}$. This, along with the definition of $Q_{2i-1}$ (see \eqref{pq14a}),
implies that 
$$
Q_{2i-1,2i} \leq G_{2i-1}(\alpha_2,\dots,\alpha_{2i-2},\alpha_{2i},\beta_1,
\dots,\beta_{2i-3},\beta_{2i-1}).
$$
But  the group
$Q_{2i+1}\cap G_{2i-1}$ is a $\pi'$-Hall subgroup of 
$G_{2i-1}(\alpha_2,\dots,\alpha_{2i},\beta_{1},\dots,\beta_{2i-1})$, 
by  \eqref{yq:e}.
 Hence the fact that 
$Q_{2i+1}\cap G_{2i-1}$ is contained in the $\pi'$-group $Q_{2i-1,2i}$ 
implies that 
$Q_{2i+1} \cap G_{2i-1} = Q_{2i-1,2i}$.
As $G_{2i-1}$ is a normal subgroup of $G_{2i+1}$ we conclude that 
$Q_{2i+1} \cap G_{2i-1} \unlhd Q_{2i+1}$,
and the lemma follows.
\end{proof}

Note that Lemma \ref{yq} resolves the problem discussed in Section \ref{pq:sec1}, at 
least for the $\pi'$-groups.
Indeed, the character $\beta_{2r-1}$ is fixed by the $\pi$-group $P_{2r}$, by \eqref{pq13a}.
Thus we can define $\beta_{2r-1,2r} \in \Irr(Q_{2r-1,2r})$ to be 
the $P_{2r}$-Glauberman correspondent 
of $\beta_{2r-1} \in \Irr(Q_{2r-1})$ whenever $1\leq r \leq k$.
Hence, in view of Lemma \ref{yq},  and starting with $\beta_1 = \chi_1$,
it makes sense to pick the character $\beta_{2i-1}$ so that it lies above 
$\beta_{2i-3,2i-2}$ whenever $2\leq i \leq l$.
 
Similarly we can work with the $\pi$-groups. 
So we can define  the group 
\begin{equation}\mylabel{anno2}
P_{2i,2i+1}= C(Q_{2i+1} \tin P_{2i}) = N(Q_{2i+1} \tin P_{2i}),
\end{equation}
whenever $1\leq i <l $. (Note that $Q_{2i+1}$ normalizes $P_{2i}$, by \eqref{pq14'a}). 
Furthermore, in a symmetric way to that we used for the $\pi'$-groups we 
 can prove
\begin{lemma}\mylabel{yp}
For every $r$ with $1\leq r  \leq k-1$ 
$$
P_{2r, 2r+1} = P_{2r+2} \cap G_{2r}.
$$ 
Hence $P_{2r,2r+1}$ is a normal subgroup of $P_{2r+2}$.
\end{lemma}
Note that, as in the case of the $\pi'$-groups, we get no information about
the groups $P_{2l-2,2l-1}= P_{2k,2k+1}$ that appear in the case where 
$m= 2l-1= 2k+1$ is odd.

As the character $\alpha_{2i}$ is fixed by $Q_{2i+1}$ (see \eqref{pq14a}),
 we can define
$\alpha_{2i,2i+1} \in \Irr(P_{2i,2i+1})$ to be the $Q_{2i+1}$-Glauberman 
correspondent of $\alpha_{2i} \in \Irr(P_{2i})$ whenever $1\leq i <l$. 
Hence, in view of Lemma \ref{yp},  and  starting with $\alpha_0 =1$, we can pick 
the character $\alpha_{2r} \in \Irr(P_{2r})$ so that it lies above 
$\alpha_{2r-2,2r-1}$ whenever $1\leq  r\leq k$. (Observe that, 
 since $\alpha_0 =1$,
the only requirement for the character $\alpha_2$ is to be an irreducible character 
of $P_{2}$).

This completes the proof of the existence of triangular sets, as the groups 
and characters we just constructed  satisfy \eqref{x}. Indeed, 
we have  proved 
\begin{proposition}\mylabel{exist}
Assume that Hypothesis \eqref{hyp1} holds for the group $G$.
 Then there exists  a triangular set 
  $\{P_{2r}, Q_{2i+1}|\alpha_{2r},\beta_{2i+1}\}$
for the normal series \eqref{*},   
 so as  to satisfy
the following conditions, whenever $1\leq r \leq k$ and $2\leq i \leq l$:
\begin{subequations}\mylabel{xxx}
\begin{gather}
P_0=1 \text{ and } \alpha_0=1, \mylabel{x1}\\
Q_1=G_1 \text{ and } \beta_1:=\chi_1 \in \Irr(Q_1) = \Irr(G_1),\mylabel{x2} \\
P_{2r} \in \Hall_{\pi}(G_{2r}(\alpha_2,\dots,\alpha_{2r-2},\beta_1,
\dots,\beta_{2r-1})),\mylabel{x3} \\
\alpha_{2r} \in \Irr(P_{2r}|\alpha_{2r-2,2r-1}), \mylabel{x4}\\
Q_{2i-1} \in \Hall_{\pi'}(G_{2i-1}(\alpha_2,\dots,\alpha_{2i-2},
\beta_1,\dots,\beta_{2i-3})),\mylabel{x5} \\
\beta_{2i-1} \in \Irr(Q_{2i-1}|\beta_{2i-3,2i-2}), \mylabel{x6}
\end{gather}
\end{subequations}
where $\alpha_{2r-2,2r-1}$ is the $Q_{2r-1}$-Glauberman correspondent of 
$\alpha_{2r-2}$ and similarly $\beta_{2i-3,2i-2}$ is the $P_{2i-2}$-Glauberman 
correspondent of $\beta_{2i-3}$.
\end{proposition}

From now on,  and until the end  of this section,   
we assume that the set  $\{P_{2r}, Q_{2i+1}|\alpha_{2r},\beta_{2i+1}\}$
is a triangular set for \eqref{*}, and therefore satisfies \eqref{xxx}.

An attempt to give a diagram that describes the relations 
in \eqref{xxx} produces the following ``Double Staircases''
 of groups and characters:
\begin{subequations}\mylabel{AA1}
\begin{equation}\mylabel{Aa}
\begin{diagram}[small]
(Q_{2l-1,2l}&\rLine)     &Q_{2l-1}       &         &            &         &       &   &    & & \\
&&\vLine&         &            &         &       &         &    &  & \\
&&Q_{2l-3,2l-2} &\rLine  &Q_{2l-3}    &         &       &   & &    &  \\
 &&    &\ddots  & &\ddots  & &  &  &    &       \\
&  && &Q_{5,6} &\rLine &Q_5        &        &        &       &        \\
&   &&&       &       &\vLine       &       &        &       &       \\
&   &  &&     &      &Q_{3,4}      &\rLine  &Q_3     &       &        \\
&   &    &&   &      &             &        &\vLine  &       &         \\
&   &      && &      &             &        &Q_{1,2} &\rLine &Q_{1}      
\end{diagram}
\end{equation}
and
\begin{equation}\mylabel{Ab}
\begin{diagram}[small]
(\beta_{2l-1,2l}&\rcorr)  &\beta_{2l-1}       &         &            &    &  &   &  &  & \\
&&\vLine&         &            &         &       &         &      &       &   \\
&&\beta_{2l-3,2l-2}&\rcorr_{P_{2l-2}} &\beta_{2l-3}    &   &    &  &  &   &  \\
  && &\ddots & &\ddots &  &      &       &     &        \\
   &&&  &\beta_{5,6}  &\rcorr_{P_{6}}  &\beta_5            &         &  && \\
   & && &     &       &\vLine         &         &  &&   \\
   & & && &  &\beta_{3,4}      &\rcorr_{P_4} &\beta_3    &    &   \\
   & &  &&&       & &       &\vLine  & &   \\
   & &  & && &    & &\beta_{1,2} &\rcorr_{P_2} &\beta_1  
\end{diagram}
\end{equation}
\end{subequations}
and similarly for the $\pi$-groups and their characters 
\begin{subequations}\mylabel{B1}
\begin{equation}\mylabel{Ba}
\begin{diagram}[small]
(P_{2k,2k+1}&\rLine ) &P_{2k}       &         &            &         &       &         &       &&   \\
&&\vLine&         &            &         &       &         &      & &    \\
&&P_{2k-2,2k-1}&\rLine &P_{2k-2}    &         &       &         &   & &\\
  &&      &\ddots         &     &\ddots         &       &         &  & &  \\
   &&&  &P_{6,7} &\rLine &P_6            &         &      &      & \\
  &   &&      &       &       &\vLine         &         &       &      & \\
   & & &&&  &P_{4,5}      &\rLine  &P_4   &      &    \\
    &     &&    &       &         &   &       &\vLine  &      &   \\
    & & & & &&     & &P_{2,3}  &\rLine &P_2  
\end{diagram}
\end{equation}
and
\begin{equation}\mylabel{Bb}
\begin{diagram}[small]
(\alpha_{2k,2k+1}&\rcorr )  &\alpha_{2k}       &         &            &         &       &         &             &    &  \\
&&\vLine&         &            &         &       &         &             &         &  \\
&&\alpha_{2k-2,2k-1}&\rcorr_{Q_{2k-1}} &\alpha_{2k-2}    &         &       &    &  & &    \\
  &&    &\ddots   &     &\ddots         &       &         &             &         &   \\
    &&& &\alpha_{6,7}  &\rcorr_{Q_{7}}  &\alpha_6            &         &     &      & \\
  &   &&      &       &       &\vLine         &         &      &      & \\
 & &  & &&&\alpha_{4,5}      &\rcorr_{Q_5} &\alpha_4    &      &    \\
 &        && &       &         &  &       &\vLine  &      &   \\
   & & &  & && & &\alpha_{2,3} &\rcorr_{Q_3} &\alpha_2 
\end{diagram}
\end{equation}
\end{subequations}
where the groups and characters in parentheses are those extra groups and characters 
that appear in the case 
of an even $m$ (for the  $\pi'$-group  $Q_{2l-1,2l}$ 
and its character $\beta_{2l-1,2l}$) or  an odd $m$ (for the $\pi$-group
$P_{2k, 2k+1}$ and its character $\alpha_{2k, 2k+1}$) respectively. 
Observe,  that every group appearing in \eqref{Aa} or  \eqref{Ba} 
is contained in all other   groups that lie above or to its right.
Furthermore, any character appearing in \eqref{Ab} or \eqref{Bb} 
is a Glauberman correspondent of the character that lies on its right.

We can actually expand these staircases  into the following 
``Double Triangles''  \eqref{pq12} and \eqref{pq11} 
of groups and characters (which is the reason behind the name 
\tria).
\begin{subequations} \mylabel{pq12}
\begin{equation} \mylabel{pq12a}
\begin{diagram}[small]
(Q_{2l-1,2l}&\rLine) &Q_{2l-1}       &         &            &         & &   &   &  &   &       &  \\
  & &\vLine&         &            &         &       &         &             &  &    &      & \\
(Q_{2l-3,2l}&\rLine)  &Q_{2l-3,2l-2}&\rLine &Q_{2l-3}    &         &       &         &             &    &    &      & \\
& &\vLine&         &\vLine&         &       &         &             &   &    &      &  \\
& &\vdots       &         &\vdots      &         &       &         &             &         &   
     &      &  \\
&&\vLine&        &\vLine&         &       &         &             &         &      
  &      & \\
(Q_{5,2l} &\rLine)  &Q_{5,2l-2}   &\rLine &Q_{5,2l-4}  &\rLine &\dots &\rLine &Q_5            &         &    &      & \\
  & &\vLine&         &\vLine&         &       &       &\vLine        &         &   
     &      & \\
(Q_{3,2l}&\rLine) &Q_{3,2l-2}   &\rLine &Q_{3,2l-4}  &\rLine  &\dots &\rLine &Q_{3,4}      &\rLine &Q_3    &      &    \\
  & &\vLine&         &\vLine&         &       &         &\vLine&       
   &\vLine &      &   \\
(Q_{1,2l}&\rLine) &Q_{1,2l-2}   &\rLine &Q_{1,2l-4}  &\rLine &\dots &\rLine &Q_{1,4}      &\rLine &Q_{1,2}
 &\rLine &Q_1  
\end{diagram}
\end{equation}
and 
\begin{equation} \mylabel{pq12b}
\begin{diagram}[small]
(\beta_{2l-1,2l}&\rcorr)  &\beta_{2l-1}       &         &            &         &       &         &             &         &   
     &       &  \\
& &\vLine&         &            &         &       &         &             &         &   
     &      & \\
(\beta_{2l-3,2l}&\rcorr)  &\beta_{2l-3,2l-2}&\rcorr_{P_{2l-2}} &\beta_{2l-3}    &         &       &         &             &         &    
    &      & \\
& &\uLine&         &\uLine &         &       &         &             &         &   
     &      &  \\
 & &\vdots       &         &\vdots      &         &       &         &             &         &   
     &      &  \\
& &\uLine &        &\uLine &         &       &         &             &         &      
  &      & \\
(\beta_{5,2l}&\rcorr) &\beta_{5,2l-2}   &\rcorr_{P_{2l-2}} &\beta_{5,2l-4}  &\rcorr_{P_{2l-4}}
 &\dots &\rcorr_{P_{6}}  &\beta_5            &         &    
    &      & \\
  & &\vLine   &         &\vLine &         &       &       &\vLine         &         &   
     &      & \\
(\beta_{3,2l}&\rcorr) &\beta_{3,2l-2}   &\rcorr_{P_{2l-2}}  &\beta_{3,2l-2}  &\rcorr_{P_{2l-4}} 
&\dots &\rcorr_{P_6} &\beta_{3,4}      
&\rcorr_{P_4} &\beta_3    &      &    \\
   & &\vLine &         &\vLine &         &       &         &\vLine  &       
   &\vLine  &      &   \\
(\beta_{1,2l}&\rcorr) &\beta_{1,2l-2}   &\rcorr_{P_{2l-2}} &\beta_{1,2l-4}  &\rcorr_{P_{2l-4}} 
&\dots &\rcorr_{P_6} &\beta_{1,4}  
    &\rcorr_{P_4} &\beta_{1,2} &\rcorr_{P_2} &\beta_1  
\end{diagram}
\end{equation}
\end{subequations}
for the $\pi'$-groups and their characters, and  similarly 
  for the $\pi$-groups:
\begin{subequations} \mylabel{pq11}
\begin{equation} \mylabel{pq11a}
\begin{diagram}[small]
(P_{2k,2k+1}&\rLine) &P_{2k}       &         &            &         &       &         &             &         &   
     &       &  \\
& &\vLine&         &            &         &       &         &             &         &   
     &      & \\
(P_{2k-2,2k+1}&\rLine) &P_{2k-2,2k-1}&\rLine &P_{2k-2}    &         &       &         &             &         &    
    &      & \\
& &\vLine&         &\vLine&         &       &         &             &         &   
     &      &  \\
 & &\vdots       &         &\vdots      &         &       &         &             &         &   
     &      &  \\
& &\vLine&        &\vLine&         &       &         &             &         &      
  &      & \\
(P_{6,2k+1}&\rLine) &P_{6,2k-1}   &\rLine &P_{6,2k-3}  &\rLine  &\dots &\rLine&P_6            &         &    
    &      & \\
&   &\vLine&         &\vLine&         &       &       &\vLine        &         &   
     &      & \\
(P_{4,2k+1}&\rLine) &P_{4,2k-1}   &\rLine &P_{4,2k-3}  &\rLine&\dots &\rLine&P_{4,5}      &\rLine &P_4   
  &      &    \\
   & &\vLine&         &\vLine&         &       &         &\vLine&       
   &\vLine &      &   \\
(P_{2,2k+1}&\rLine) &P_{2,2k-1}   &\rLine &P_{2,2k-3}  &\rLine&\dots &\rLine &P_{2,5}      &\rLine &P_{2,3}
 &\rLine &P_2  
\end{diagram}
\end{equation}
and 
\begin{equation} \mylabel{pq11b}
\begin{diagram}[small]
(\alpha_{2k,2k+1}&\rcorr) &\alpha_{2k}       &         &            &         &       &         &             &         &   
     &       &  \\
& &\vLine&         &            &         &       &         &             &         &   
     &      & \\
(\alpha_{2k-2,2k+1}&\rcorr) &\alpha_{2k-2,2k-1}&\rcorr_{Q_{2k-1}} &\alpha_{2k-2}    &         &       &         &             &         &    
    &      & \\
& &\uLine&         &\uLine &         &       &         &             &         &   
     &      &  \\
 & &\vdots       &         &\vdots      &         &       &         &             &         &   
     &      &  \\
& &\uLine &        &\uLine &         &       &         &             &         &      
  &      & \\
(\alpha_{6,2k+1}&\rcorr) &\alpha_{6,2k-1}   &\rcorr_{Q_{2k-1}} &\alpha_{6,2k-3}  &\rcorr_{Q_{2k-3}} 
&\dots &\rcorr_{Q_{7}}  &\alpha_6            &         &    
    &      & \\
  & &\vLine   &         &\vLine &         &       &       &\vLine         &         &   
     &      & \\
(\alpha_{4,2k+1}&\rcorr) &\alpha_{4,2k-1}   &\rcorr_{Q_{2k-1}}  &\alpha_{4,2k-3}  &\rcorr_{Q_{2k-3}} 
&\dots &\rcorr_{Q_7} &\alpha_{4,5}      
&\rcorr_{Q_5} &\alpha_4    &      &    \\
   & &\vLine &         &\vLine &         &       &         &\vLine  &       
   &\vLine  &      &   \\
(\alpha_{2,2k+1}&\rcorr) &\alpha_{2,2k-1}   &\rcorr_{Q_{2k-1}} &\alpha_{2,2k-3}  &\rcorr_{Q_{2k-3}} 
&\dots &\rcorr_{Q_7} &\alpha_{2,5}  
    &\rcorr_{Q_5} &\alpha_{2,3} &\rcorr_{Q_3} &\alpha_2  
\end{diagram}
\end{equation}
\end{subequations}
where,  as before,  $l=[(m+1)/2]$ and $k= [m/2]$.  Furthermore,
 the $\pi'$-groups $Q_{2i-1, 2l}$  and characters $\beta_{2i-1,2l}$ in 
parentheses exist only when $m$ is even,  and the $\pi$-groups 
$P_{2r, 2k+1}$ and characters $\alpha_{2r, 2k+1}$  only when 
$m$ is odd.

Before we give the long list 
of the groups, the characters and their properties
that are involved in the above diagrams, we remark again that 
the groups and characters  that appear in the first two  diagonals of the above 
``Double Triangle'' diagrams form the ``Double Staircase'' diagrams.
 The rest of the groups that appear in the above  diagrams 
are defined as
\begin{subequations}\mylabel{pq1413def}
\begin{equation}\mylabel{pq1413defa}
Q_{2i-1,2j}:= N(P_{2i},P_{2i+2},\dots,P_{2j} \tin Q_{2i-1}),
\end{equation} 
and 
\begin{equation}\mylabel{pq1413defb}
P_{2r,2s+1}:=N(Q_{2r+1},Q_{2r+3},\dots,Q_{2s+1} \tin P_{2r}),
\end{equation}
\end{subequations}
for all $i,j, r,s$ with $1\leq i \leq l$, $i\leq j \leq k$, $1\leq r \leq k $ and 
 $r\leq s \leq l-1$.
Note that, in view of \eqref{pq14'a} and \eqref{pq14'b}, the products
$P_{2i}\cdot P_{2i+2} \cdots P_{2j}$ and $Q_{2r+1}\cdot Q_{2r+3}
 \cdots Q_{2s+1}$
form groups,  for all $i,j,r,s$ as above. This, along with 
\eqref{pq1413def}, implies that  
\begin{subequations}\mylabel{pq1413e}
\begin{multline}\mylabel{pq14e}
Q_{2i-1,2j}= N(P_{2i},P_{2i+2},\dots,P_{2j} \tin Q_{2i-1})\\=
 C(P_{2i},P_{2i+2},\dots,P_{2j} \tin Q_{2i-1})= 
 C(P_{2i}\cdot P_{2i+2}\cdots P_{2j} \tin Q_{2i-1}),
\end{multline}
and
\begin{multline}\mylabel{pq13e}
P_{2r,2s+1}=N(Q_{2r+1},Q_{2r+3},\dots,Q_{2s+1} \tin P_{2r})\\
= C(Q_{2r+1},Q_{2r+3},\dots,Q_{2s+1} \tin P_{2r})=
C(Q_{2r+1}\cdot Q_{2r+3}\cdots  Q_{2s+1} \tin P_{2r}),
\end{multline}
whenever $1\leq i\leq l$, $i \leq j \leq k$, $1\leq r \leq k$ and 
 $r\leq s \leq l-1$. 
\end{subequations}
Furthermore, the way the groups $Q_{2i-1,2j}$ and $P_{2r,2s+1}$
are defined, along with \eqref{pq14'}, implies that 
\begin{subequations}\mylabel{pq15'}
\begin{equation}\mylabel{pq15'a}
Q_{2i-1} \text{ normalizes the groups }  Q_{2t-1,2j} \text{ and } P_{2t,2j+1},
\end{equation}
whenever $1\leq t\leq j \leq i-1 \leq l-1$. Similarly, 
\begin{equation}\mylabel{pq15'b}
P_{2r} \text{ normalizes the groups }   Q_{2t+1,2j} \text{ and } P_{2t,2j-1}, 
\end{equation}
whenever $1\leq t <j\leq r \leq k$.
\end{subequations}

Looking at the diagrams \eqref{pq12a} and \eqref{pq11a}, we seee that 
what \eqref{pq14'} and \eqref{pq15'} say is that any group on
 the main diagonal 
of these diagrams, that is $Q_{2i-1}$ or $P_{2r}$,
 normalizes all the other groups that lie below  or to its  right.
They  also say that $Q_{2i-1}$ normalizes all the groups in \eqref{pq11a}
which are below or to the right of $P_{2i-2, 2i-1}$, while $P_{2r}$ 
normalizes 
all the groups in \eqref{pq12a} which are below or to the right of 
$Q_{2r-1,2r}$.

Furthermore,   in the case that $j >i$
and $s > r$ (with $i,j,s,r$ as in \eqref{pq1413e}),
the groups $Q_{2i-1,2j-2}$ and $P_{2r,2s-1}$  satisfy
the  equations   \eqref{pq14e} and \eqref{pq13e}, respectively. 
Hence 
$$Q_{2i-1,2j} = N(P_{2j} \tin N(P_{2i},P_{2i+2},\dots,P_{2j-2} \tin Q_{2i-1}))
=N(P_{2j} \tin Q_{2i-1,2j-2}).$$
But $P_{2j}$ normalizes  $Q_{2i-1,2j-2}$ by \eqref{pq15'b}.
Therefore
\begin{subequations}\mylabel{pq1413f}
\begin{equation}\mylabel{pq14f}
Q_{2i-1,2j}=N(P_{2j} \tin Q_{2i-1,2j-2}) = C(P_{2j} \tin Q_{2i-1,2j-2}),
\end{equation}
and similarly for the $\pi$-groups
\begin{equation}\mylabel{pq13f}
P_{2r,2s+1} = N(Q_{2s+1} \tin P_{2r,2s-1}) = C(Q_{2s+1} \tin P_{2r,2s-1}),
\end{equation}
whenever $1\leq i \leq l$, \, $i< j \leq k$, $1\leq r \leq k$, and  $r < s \leq l-1$. 
\end{subequations}

According to \eqref{x5} and \eqref{x3},
the groups $Q_{2i-1}$ and $P_{2r}$ were chosen to be $\pi'$-Hall and
$\pi$-Hall subgroups
of specific ``stabilizer''-subgroups of $G_{2i-1}$ and $G_{2r}$, respectively.  
A similar characterization  for the groups
 $Q_{2i-1,2j}$ and $P_{2r,2s+1}$ is 
described and proved in
\begin{proposition}\mylabel{pqremark2}
For every $i,j$ with $1\leq i\leq j \leq l-1$ the following holds
\begin{subequations}\mylabel{pqarca}
\begin{gather}
Q_{2i-1,2j} = Q_{2j+1} \cap G_{2i-1} 
\text{ and therefore } \mylabel{pqarca1}\\
Q_{2i-1,2j} \in 
\Hall_{\pi'}(G_{2i-1}(\alpha_2,\dots,\alpha_{2j},\beta_1,\dots,\beta_{2j-1})) \mylabel{pqarca2} \\
=\Hall_{\pi'}(G_{2i-1}(\alpha_2,\dots,\alpha_{2j},\beta_1,\dots,
\beta_{2j-1},\beta_{2j+1})).\notag
\end{gather}
\end{subequations}
Similarly,  for  all $r,s$ with $1\leq r\leq s \leq k-1$ we have 
\begin{subequations}\mylabel{pqarcb}
\begin{gather}
P_{2r,2s+1} = P_{2s+2} \cap G_{2r} \text{ and therefore } \mylabel{pqarcb1}\\
 P_{2r,2s+1} \in 
\Hall_{\pi}(G_{2r}(\alpha_2,\dots,\alpha_{2s},\beta_1,\dots,\beta_{2s+1})) \mylabel{pqarcb2}\\
=\Hall_{\pi}(G_{2r}(\alpha_2,\dots,\alpha_{2s},\alpha_{2s+2}, 
\beta_1,\dots,\beta_{2s+1})). \notag 
\end{gather}
\end{subequations}
\end{proposition}
Note that the extra groups $Q_{2i-1,2l}$ (when $m=2l=2k$), and $P_{2r,2k+1}$ (when
 $m=2l-1=2k+1$) are not covered in Proposition \ref{pqremark2}.

\begin{proof}
The definition of $Q_{2j+1}$ in  \eqref{pq14a}, along  with the fact that $G_{2i-1}$ 
is a normal subgroup of $G_{2j+1}$ whenever $1\leq i \leq j \leq l-1$, implies that 
\begin{equation}\mylabel{pqarcae1}
Q_{2j+1} \cap G_{2i-1} \in \Hall_{\pi'}(G_{2i-1}(\alpha_2, \dots,\alpha_{2j}, \beta_{1}, 
\dots,\beta_{2j-1})).
\end{equation}
But $Q_{2j+1}$ also fixes $\beta_{2j+1}$. Hence 
\begin{equation}\mylabel{pqarcae1'}
Q_{2j+1} \cap G_{2i-1} \in \Hall_{\pi'}(G_{2i-1}(\alpha_2, \dots,\alpha_{2j}, \beta_{1}, 
\dots,\beta_{2j-1},\beta_{2j+1})).
\end{equation}

In particular, $Q_{2j+1}\cap G_{2i-1}$ is a $\pi'$-subgroup of 
$G_{2i-1}(\alpha_2,\dots,\alpha_{2i-2},\beta_1,\dots,\beta_{2i-3})$.
Furthermore, $Q_{2j+1}$ fixes $\beta_{2i-1}$, and  so normalizes $Q_{2i-1}$
But the latter is  a $\pi'$-Hall subgroup 
 of  the group $G_{2i-1}(\alpha_2, \dots,\alpha_{2i-2}, \beta_{1},
\dots,\beta_{2i-3})$.
Therefore $Q_{2j+1} \cap G_{2i-1} \leq Q_{2i-1}$.
As $Q_{2j+1}$ normalizes the groups $P_2,\dots,P_{2j}$, by \eqref{pq14'a},
we conclude that $Q_{2j+1}\cap G_{2i-1}$ is a subgroup of 
$Q_{2i-1,2j}= N(P_{2i},\dots,P_{2j} \tin Q_{2i-1})$.
 Hence for all $i,j$ with 
$1\leq i \leq j \leq l-1$ we have 
\begin{equation}\mylabel{pqarcae2}
Q_{2j+1}\cap G_{2i-1} \leq Q_{2i-1,2j}.
\end{equation}

To prove the opposite inclusion, and complete the proof of \eqref{pqarca},
we will use induction on $j$.
According to Lemma \ref{yq}  we have 
$Q_{2i-1,2i} = Q_{2i+1}\cap G_{2i-1}$,  for all $i$ with $1\leq i \leq l-1$. 
Hence the proposition holds  in the case that $i=j$.

Suppose that, for some fixed $r=i+1,\dots,l-1$ and  for all $j$ with 
$1\leq i \leq j < r$, we have 
 $Q_{2i-1,2j} \leq Q_{2j+1} \cap G_{2i-1}$ (and thus equality as the other inclusion is
proved). 
Then according to \eqref{pq14f},  we have $Q_{2i-1,2r} = C(P_{2r} \tin Q_{2i-1,2r-2})$.
By our supposition $Q_{2i-1,2r-2} $ is a subgroup of $Q_{2r-1}$.
Hence $Q_{2i-1,2r} \leq C(P_{2r} \tin Q_{2r-1})$. But $C(P_{2r} \tin Q_{2r-1})= Q_{2r-1,2r}$,
by \eqref{pq14e}. Therefore $Q_{2i-1,2r} \leq Q_{2r-1,2r}$.
Furthermore, $Q_{2r-1,2r} \leq Q_{2r+1}$, by Lemma \ref{yq}.
Hence $Q_{2i-1,2r} \leq Q_{2r+1}$.
This proves the inductive argument in the case that $j=r$.
Hence $Q_{2i-1,2j} \leq Q_{2j+1} \cap G_{2i-1}$ whenever 
$1\le i \le j \le l-1$. This, along with 
\eqref{pqarcae1}, \eqref{pqarcae1'} and \eqref{pqarcae2}, completes the proof of \eqref{pqarca}.

The proof for \eqref{pqarcb} is similar, so we omit it. 
As a final remark, we observe that the only tools we used for the  proof of Proposition \ref{pqremark2} are the definitions of the groups $Q_{2i-1}, P_{2r}$, in  \eqref{x5} and \eqref{x3},  and the definitions of the groups $Q_{2i-1,2j}$ and $P_{2r,2s+1}$
in  \eqref{pq1413def}.
\end{proof}

Proposition \ref{pqremark2} implies
\begin{cor}
For all $i,j,r,s$ with  $1\leq i \leq j \leq l-1$ and $1\leq r \leq s \leq k-1$
we have 
$$
Q_{2i-1} \cap Q_{2j+1} = Q_{2i-1,2j} \text{ and }
P_{2r} \cap P_{2s+2} = P_{2r, 2s+1}.
$$
Therefore
\begin{multline}\mylabel{pq14b}
Q_{2i-1,2j}= N(P_{2j} \tin Q_{2i-1,2j-2})=
 C(P_{2j} \tin Q_{2i-1,2j-2})\\
=N(P_{2i},\dots,P_{2j} \tin Q_{2i-1})=
 C(P_{2i},\dots,P_{2j} \tin Q_{2i-1})\\
=C(P_{2i}\cdot P_{2i+2}\cdots P_{2j} \tin Q_{2i-1}) = Q_{2i-1} \cap Q_{2j+1},
\end{multline}
and 
\begin{multline} \mylabel{pq13b}
   P_{2r, 2s+1} = N(Q_{2s+1} \tin P_{2r,2s-1}) = 
 C(Q_{2s+1} \tin P_{2r,2s-1}) \\
  = N(Q_{2r+1}, \dots ,Q_{2s+1} \tin P_{2r}) = 
     C(Q_{2r+1}, \dots ,Q_{2s+1} \tin P_{2r} )\\
   = C(Q_{2r+1} \cdot Q_{2r+3} \cdots Q_{2s+1} \tin P_{2r})
                = P_{2r} \cap P_{2s+2},
 \end{multline}
where, by convention,  we write $Q_{2i-1,2i-2}= Q_{2i-1}$ and $P_{2r,2r-1} = P_{2r}$.
Furthermore,  for any $t,t'$ with $i\leq t \leq j$ and 
$r\leq  t' \leq s $, 
where $i,j,r,s $ are as above,  we have
\begin{equation}\mylabel{pq:norm}
Q_{2i-1,2j} \unlhd Q_{2t-1,2j} \unlhd Q_{2j+1}\text{ and } 
P_{2r,2s+1} \unlhd P_{2t',2s+1}\unlhd P_{2s+2}. 
\end{equation}
Similarly for the extra groups $Q_{2i-1,2l}$ and $P_{2i,2k+1}$ 
we have 
\begin{equation}\mylabel{pq:norm'}
\begin{aligned}
Q_{2i-1,2l}&\unlhd Q_{2t-1,2l} &\text{ when $m=2k$ and thus $k=l$, }\\
P_{2r,2k+1} &\unlhd P_{2t',2k+1} &\text{ when $m=2l-1$ and thus $k=l-1$},
\end{aligned}
\end{equation}
whenever $1\leq i \leq t \leq l$ and $1\leq r\leq t'\leq k$.
\end{cor}

\begin{proof}
The first part follows easily from Proposition \ref{pqremark2} and the two sets of 
inclusions $Q_{2i-1,2j} \leq Q_{2i-1} \leq G_{2i-1}$ and $P_{2r,2s+1} \leq P_{2r} \leq G_{2r}$.
The multiple equations \eqref{pq14b} and \eqref{pq13b}
are a collection of \eqref{pq14e}, \eqref{pq14f}, \eqref{pq13e} and
\eqref{pq13f}.
Also \eqref{pq:norm} follows directly  from 
Proposition \ref{pqremark2},  since  $G_{2i-1} \unlhd  G_{2t-1}$
and $G_{2r} \unlhd G_{2t'}$ whenever  $1\leq i \leq t \leq l-1$ and 
$1\leq  r \leq t' \leq k-1$.
It remains to show that \eqref{pq:norm'} also holds for the extra groups 
(whenever these exist)
$Q_{2i-1,2l}$ and $P_{2r,2k+1}$.
Indeed, in the case that $m=2k$ is even (and so $k=l$) the groups $Q_{2i-1,2l}$
 are well defined (see \eqref{pq1413defa}) for all $i=1,\dots,l$. Furthermore,
\eqref{pq14f} implies that $Q_{2i-1,2l}= N(P_{2l} \tin Q_{2i-1,2l-2})$ for
 all $i=1,\dots,l-1$.
Since $Q_{2i-1,2l-2} \unlhd Q_{2t-1,2l-2} \unlhd Q_{2l-1}$
whenever  $1\leq i\leq t \leq l-1$, we easily have 
that 
$$
Q_{2i-1,2l}=N(P_{2l} \tin Q_{2i-1,2l-2})\unlhd N(P_{2l} \tin Q_{2t-1,2l-2})
=Q_{2t-1,2l} \unlhd N(P_{2l} \tin Q_{2l-1}) =Q_{2l-1,2l},
$$ 
for all such $i$ and $t$. This proves \eqref{pq:norm'} for the $\pi'$-groups.
The proof for the $\pi$-groups (that occurs when $m=2l-1=2k+1$)
 is similar. So we omit it.
\end{proof}

The following proposition covers the extra groups that Proposition 
\ref{pqremark2}
left out. 
\begin{proposition}\mylabel{remacor}
For every $i,r$ with $1\leq i \leq l$ and $1\leq r \leq k$ we have 
\begin{equation}\mylabel{arcab}
\begin{aligned}
Q_{2i-1,2l} &\in 
\Hall_{\pi'}(G_{2i-1}(\alpha_2,\dots,\alpha_{2l},\beta_1,\dots,\beta_{2l-1}))
&\text{ if  $m=2k=2l$,  }\\
P_{2r,2k+1} &\in 
\Hall_{\pi}(G_{2r}(\alpha_2,\dots,\alpha_{2k},\beta_1,\dots,\beta_{2k+1}))
&\text{ if   $m=2l-1=2k+1$ }.\\
\end{aligned}
\end{equation}
\end{proposition}

\begin{proof}
Assume that $m=2k=2l$ is even. Then   for all $i=1,\dots,l$
the groups $Q_{2i-1,2l}= C(P_{2i},\dots,P_{2l} \tin Q_{2i-1})$ are well
 defined (see \eqref{pq1413defa}).
By \eqref{x5} we have that 
 $Q_{2l-1}$ is a $\pi'$-Hall
subgroup of $G_{2l-1}(\alpha_2,\dots,\alpha_{2l-2},\beta_1,
\dots,\beta_{2l-3})$. Thus 
  $Q_{2l-1}$ is also a $\pi'$-Hall subgroup of 
$G_{2l-1}(\alpha_2,\dots,\alpha_{2l-2},\beta_1,
\dots,\beta_{2l-1})$,  as $Q_{2l-1}$  fixes $\beta_{2l-1}$.
Furthermore, according to \eqref{pqarcb} for $r=s=k-1$,  we get  
$P_{2l-2,2l-1}=P_{2k-2,2k-1} \in 
\Hall_{\pi}(G_{2l-2}(\alpha_2,\dots,\alpha_{2l-2},\beta_1,
\dots,\beta_{2l-1}))$.  Thus 
$P_{2l-2,2l-1}$ is also a $\pi$-Hall subgroup of $G_{2l-1}(\alpha_2,\dots,
\alpha_{2l-2},\beta_1,\dots,\beta_{2l-1})$, since 
 $G_{2l-1}/G_{2l-2}$ is a $\pi'$-group.
Since   $P_{2l-2,2l-1} =C(Q_{2l-1} \tin P_{2l-2})$,  we have 
$$
G_{2l-1}(\alpha_2,\dots,\alpha_{2l-2},\beta_1,\dots,\beta_{2l-1})=Q_{2l-1} \times 
P_{2l-2,2l-1}.
$$
This implies that 
$Q_{2l-1}(\alpha_{2l})$ is a  $\pi'$-Hall subgroup of $G_{2l-1}(\alpha_2,\dots,
\alpha_{2l-2},\alpha_{2l},\beta_1,\dots,\beta_{2l-1})$.
Furthermore, 
$Q_{2l-1}(\alpha_{2l})\leq N(P_{2l} \tin Q_{2l-1}) = C(P_{2l} \tin Q_{2l-1}) \leq Q_{2l-1}(\alpha_{2l})$.
Hence $Q_{2l-1,2l} = N(P_{2l} \tin Q_{2l-1})$ is a $\pi'$-Hall subgroup 
of $G_{2l-1}(\alpha_2,\dots,\alpha_{2l},\beta_1,\dots,\beta_{2l-1})$. Thus 
\eqref{arcab} holds for $i=l$.
Also for any $i=1,\dots ,l-1$ we have 
\begin{align}\mylabel{extra1}
Q_{2i-1,2l} &=N(P_{2l} \tin Q_{2i-1,2l-2})  &\quad  (&\text{ by 
\eqref{pq14f}})\notag \\
                &=N(P_{2l} \tin Q_{2l-1} \cap G_{2i-1}) &\quad (&\text{ by 
\eqref{pqarca1}})\notag \\
                &=N(P_{2l} \tin Q_{2l-1}) \cap G_{2i-1} &\quad  &\notag \\
                &=Q_{2l-1,2l} \cap G_{2i-1}.    
\end{align}
This,  along with the facts that $G_{2i-1}\unlhd G_{2l-1}$ and $Q_{2l-1,2l}\in \Hall _{\pi'}
(G_{2l-1}(\alpha_2,\dots,\alpha_{2l},\beta_1,\dots,\beta_{2l-1}))$, implies 
that 
$$
Q_{2i-1,2l}\in \Hall _{\pi'}
(G_{2i-1}(\alpha_2,\dots,\alpha_{2l},\beta_1,\dots,\beta_{2l-1}))
$$
whenever $1\leq i \leq l$. 
Hence \eqref{arcab} holds  in the case $m=2k=2l$.

Similarly we can work with the $\pi$-groups in the case of an odd $m=2l-1$.
\end{proof}

As a straight forward consequence of \eqref{extra1} and \eqref{pqarca1}
we have    
\begin{remark}\mylabel{extra}
For every $i=1,\dots,l$ and every $j,s$ with $1\le i \le j \le s \le k$ the
following holds:
$$ Q_{2i-1,2s}= Q_{2j-1,2s} \cap G_{2i-1}.
 $$
\end{remark}

Regarding the possible products of the  groups $Q_{2i-1,2j}$ and $P_{2r,2s+1}$
we have  
\begin{proposition}\mylabel{prodd}
For every $i,r$ with $i=2,\dots,l$ and $r=1,\dots,k$,  we have 
\begin{subequations}
\begin{gather}
 G_{2r}(\alpha_2, \dots  ,\alpha_{2r-2}, \beta_1,  \dots ,  \beta_{2r-1})
= P_{2r}\ltimes Q_{2r-1}, \mylabel{pq15a}\\ 
G_{2r}(\alpha_2, \dots  ,\alpha_{2r}, \beta_1,  \dots ,  \beta_{2r-1})
= P_{2r}\times Q_{2r-1,2r}, \mylabel{pq-15a}\\ 
G_{2i-1}(\alpha_2, \dots  ,\alpha_{2i-2}, \beta_1, \dots ,\beta_{2i-3})
=P_{2i-2} \rtimes Q_{2i-1}, \mylabel{pq16a}\\
G_{2i-1}(\alpha_2, \dots  ,\alpha_{2i-2}, \beta_1, \dots ,\beta_{2i-1})
=P_{2i-2,2i-1} \times Q_{2i-1}.\mylabel{pq-16a}
\end{gather}
\end{subequations}
Furthermore 
\begin{subequations}
\begin{gather}
 G_{2r}(\alpha_2, \dots  ,\alpha_{2s-2}, \beta_1,  \dots ,  \beta_{2s-1})
= P_{2r,2s-1} \times Q_{2r-1,2s-2}, \mylabel{prod1}\\
 G_{2r}(\alpha_2, \dots  ,\alpha_{2s}, \beta_1,  \dots ,  \beta_{2s-1})
= P_{2r,2s-1} \times Q_{2r-1,2s}, \mylabel{prod2}\\
G_{2i-1}(\alpha_2, \dots  ,\alpha_{2j-2}, \beta_1, \dots ,\beta_{2j-1})=
P_{2i-2,2j-1} \times Q_{2i-1,2j-2}, \mylabel{prod3}\\
G_{2t-1}(\alpha_2, \dots  ,\alpha_{2v}, \beta_1, \dots ,\beta_{2v-1})=
P_{2t-2,2v-1} \times Q_{2t-1,2v}, \mylabel{prod4}
\end{gather}
\end{subequations}
whenever $1\leq  i <j \leq l$, \, $1\leq r < s\leq k$ and $1\leq t \leq v \leq k$.
\end{proposition}

Note that, according to \eqref{kl},  we have $k\leq l\leq k+1$.  Hence all the 
above groups are well defined.

\begin{proof}
Clearly \eqref{x1} and \eqref{x2}, along with the 
fact that $G_2 / G_1$ is a $\pi$-group,  imply 
\begin{equation}\mylabel{ase}
Q_1 = G_1  \in \Hall_{\pi'}(G_1) \cap \Hall_{\pi'}(G_2(\beta_1)) \cap 
\Hall_{\pi'}(G_2(\alpha_0 , \beta_1)).
\end{equation}
In addition,  for all $i=2,\dots,k$  the factor group 
  $G_{2i}/G_{2i-1}$ is a $\pi$-group. Furthermore,  in view of  \eqref{x5} 
the group 
$Q_{2i-1}$ is a $\pi'$-Hall subgroup of $G_{2i-1}(\alpha_2,\dots,
\alpha_{2i-2},\beta_1,
\dots,\beta_{2i-3})$. Hence 
$Q_{2i-1}$ is also a $\pi'$-Hall subgroup of 
$G_{2i}(\alpha_2,\dots,\alpha_{2i-2},\beta_1,
\dots,\beta_{2i-3})$.
As $Q_{2i-1}$ obviously fixes the character  $\beta_{2i-1} \in 
\Irr(Q_{2i-1})$, we conclude that 
\begin{subequations}\mylabel{hallq}
\begin{equation}\mylabel{hallq1}
\begin{aligned}
Q_{2i-1} \in &\Hall_{\pi'}(G_{2i-1}(\alpha_2,\dots,\alpha_{2i-2},\beta_1,
\dots,\beta_{2i-3}))\\
&\cap  \Hall_{\pi'}(G_{2i-1}(\alpha_2,\dots,\alpha_{2i-2},\beta_1,
\dots,\beta_{2i-1}))\\
&\cap  \Hall_{\pi'}(G_{2i}(\alpha_2,\dots,\alpha_{2i-2},\beta_1,
\dots,\beta_{2i-3}))\\
&\cap  \Hall_{\pi'}(G_{2i}(\alpha_2,\dots,\alpha_{2i-2},\beta_1,
\dots,\beta_{2i-1})),
\end{aligned}
\end{equation}
whenever $1 \leq i \leq k$, while 
\begin{equation}\mylabel{hallq2}
\begin{aligned}
Q_{2l-1} \in &\Hall_{\pi'}(G_{2l-1}(\alpha_2,\dots,\alpha_{2l-2},\beta_1,
\dots,\beta_{2l-3}))\\
&\cap  \Hall_{\pi'}(G_{2l-1}(\alpha_2,\dots,\alpha_{2l-2},\beta_1,
\dots,\beta_{2l-1})).
\end{aligned} 
\end{equation}
\end{subequations}
Note that we need to include as a special case the group  $Q_{2l-1}$,
since it is not covered when $m=2l-1$ is odd.

Similarly for the $\pi$-groups we have 
\begin{subequations}\mylabel{hallp}
\begin{equation}\mylabel{hallp1}
\begin{aligned}
P_{2r}\in &\Hall_{\pi}(G_{2r}(\alpha_2,\dots,\alpha_{2r-2},
\beta_1,\dots,\beta_{2r-1}))\\
&\cap \Hall_{\pi}(G_{2r}(\alpha_2,\dots,\alpha_{2r},
\beta_1,\dots,\beta_{2r-1}))\\
&\cap  \Hall_{\pi}(G_{2r+1}(\alpha_2,\dots,
\alpha_{2r-2},\beta_1,\dots,\beta_{2r-1}))\\
&\cap \Hall_{\pi}(G_{2r+1}(\alpha_2,\dots,\alpha_{2r},
\beta_1,\dots,\beta_{2r-1})),
\end{aligned}
\end{equation}
whenever $1\leq r \leq l-1$, while
\begin{equation}\mylabel{hallp2}
\begin{aligned}
P_{2k}\in &\Hall_{\pi}(G_{2k}(\alpha_2,\dots,\alpha_{2k-2},
\beta_1,\dots,\beta_{2k-1}))\\
&\cap \Hall_{\pi}(G_{2r}(\alpha_2,\dots,\alpha_{2k},
\beta_1,\dots,\beta_{2k-1})).
\end{aligned}
\end{equation}
\end{subequations}

Furthermore, \eqref{pqarca} and 
 \eqref{arcab}, along with the fact 
that $G_{2i}/G_{2i-1}$ is a $\pi$-group,   imply that
\begin{subequations}\mylabel{hallqq}
\begin{equation}
\begin{aligned}
 Q_{2i-1,2j-2}\in &\Hall_{\pi'}(G_{2i-1}(\alpha_2,\dots,\alpha_{2j-2},
\beta_1,\dots,\beta_{2j-3}))\\
&\cap  \Hall_{\pi'}(G_{2i-1}(\alpha_2,\dots,
\alpha_{2j-2},\beta_1,\dots,\beta_{2j-1}))\\
&\cap \Hall_{\pi'}(G_{2i}(\alpha_2,\dots,\alpha_{2j-2},
\beta_1,\dots,\beta_{2j-3}))\\
&\cap \Hall_{\pi'}(G_{2i}(\alpha_2,\dots,\alpha_{2j-2},
\beta_1,\dots,\beta_{2j-1})),
\end{aligned}
\end{equation}
whenever $1\leq i <j \leq l$,  while  for all $t$ with $1\leq t \leq l$ we have 
\begin{equation}
\begin{aligned}
Q_{2t-1,2l} \in &\Hall_{\pi'}(G_{2t-1}(\alpha_2,\dots,\alpha_{2l},
\beta_1,\dots,\beta_{2l-1}))\\
&\cap \Hall_{\pi'}(G_{2t}(\alpha_2,\dots,\alpha_{2l},
\beta_1,\dots,\beta_{2l-1})).
\end{aligned}
\end{equation}
\end{subequations}

Similarly,  \eqref{pqarcb}, \eqref{arcab} and the fact that $G_{2r+1} /G_{2r}$ is
a $\pi'$-group imply that 
\begin{subequations}\mylabel{hallpp}
\begin{equation}
\begin{aligned}
P_{2r,2s-1}\in &\Hall_{\pi}(G_{2r}(\alpha_2,\dots,\alpha_{2s-2},
\beta_1,\dots,\beta_{2s-1}))\\
&\cap  \Hall_{\pi}(G_{2r}(\alpha_2,\dots,
\alpha_{2s},\beta_1,\dots,\beta_{2s-1}))\\
&\cap \Hall_{\pi}(G_{2r+1}(\alpha_2,\dots,\alpha_{2s-2},
\beta_1,\dots,\beta_{2s-1}))\\
&\cap \Hall_{\pi}(G_{2r+1}(\alpha_2,\dots,\alpha_{2s},
\beta_1,\dots,\beta_{2s-1})),
\end{aligned}
\end{equation}
whenever   $1\leq r <s \leq k $, while for all $t'$ with   $1\leq t' \leq k$ we have 
\begin{equation}
\begin{aligned}
P_{2t',2k+1}\in &\Hall_{\pi}(G_{2t'}(\alpha_2,\dots,\alpha_{2k},
\beta_1,\dots,\beta_{2k+1}))\\
&\cap \Hall_{\pi}(G_{2t'+1}(\alpha_2,\dots,\alpha_{2k},
\beta_1,\dots,\beta_{2k+1})).
\end{aligned}
\end{equation}
\end{subequations}

Furthermore, $P_{2r}$ normalizes $Q_{2r-1}$, 
while $Q_{2i-1}$ normalizes $P_{2i-2}$.
 Therefore \eqref{ase},  \eqref{hallq} and \eqref{hallp} imply that 
$$
   P_{2r} \ltimes Q_{2r-1} =  G_{2r}(\alpha_2, 
\dots  ,\alpha_{2r-2}, \beta_1,  \dots ,  \beta_{2r-1}), 
$$
and 
$$
  P_{2i-2} \rtimes Q_{2i-1} =  G_{2i-1}(\alpha_2, \dots  ,\alpha_{2i-2}, 
\beta_1, \dots ,       \beta_{2i-3}),
$$
for all $i=2,\dots,l$ and $r=1,\dots,k$.
For the same range of $i$ and $r$ equations \eqref{anno1} and \eqref{anno2}
imply that $Q_{2r-1,2r}$ centralizes $P_{2r}$ while $P_{2i-2,2i-1}$ centralizes 
$Q_{2i-1}$. As these groups are $\pi$-and $\pi'$-Hall subgroups of the correct
groups (see \eqref{hallp}, \eqref{hallq}, \eqref{hallpp} and \eqref{hallqq})
equations \eqref{pq-15a} and \eqref{pq-16a} follow. 

We can work similarly  for the rest of the proposition. 
We only remark here that,  whenever $1\leq r <s \leq k$, \, $1\leq i \leq j \leq  l $ and
 $1\leq t\leq v \leq k$, 
equations  \eqref{pq14e}  and \eqref{pq13e}  
imply that $Q_{2r-1,2s-2}$ and $Q_{2r-1,2s}$ centralize 
$P_{2r}$ (and thus $P_{2r,2s-1}$), while $P_{2i-2, 2j-1}$ and $P_{2t-2,2v-1}$
centralize $Q_{2i-1}$ and $Q_{2t-1}$,  respectively.
This, along with \eqref{hallpp} and  \eqref{hallqq},
implies  the rest of the proposition.
\end{proof}

What about the characters that appear in  the diagrams \eqref{pq12b} and
 \eqref{pq11b}? We  have already seen, in  \eqref{x6} and \eqref{x4}, 
that $\beta_{2i-1}$ and $\alpha_{2r}$ are irreducible 
characters of $Q_{2i-1}$ and $P_{2r}$, respectively. 
Furthermore, according to \eqref{x3}, for every  $i=1,\dots,l$
 the character $\beta_{2i-1}$ is fixed by the
$\pi$-groups $P_{2i},P_{2i+2},\dots,P_{2k}$,  and thus 
is also fixed by their product (note that their product forms a group
 according to \eqref{pq14'b}). Similarly whenever 
$r=1,\dots,k$, using \eqref{x5}, 
we see that the character $\alpha_{2r}$ is fixed by the groups 
$Q_{2r+1}, Q_{2r+3},\dots,Q_{2l-1}$, 
and therefore is also fixed by their product. 
Hence, we can naturally make the 
\begin{defn}\mylabel{pq1413def2}
\begin{itemize}
\item[(1)]  We write  $\beta_{2i-1,2j} \in \Irr(Q_{2i-1,2j})$  for  the  
$P_{2i}\cdot P_{2i+2}\cdots  P_{2j}$-Glauberman correspondent of 
$\beta_{2i-1} \in \Irr(Q_{2i-1})$ and 
\item[(2)] we write $\alpha_{2r,2s+1} \in \Irr( P_{2r,2s+1})$ for  the 
$Q_{2r+1} \cdot Q_{2r+3} \cdots Q_{2s+1}$-Glauberman correspondent
of $\alpha_{2r} \in \Irr(P_{2r})$, 
\end{itemize}
whenever $1\leq i \leq l$, $i\leq j \leq k$,  $1\leq r \leq k $  and 
$r \leq s \leq l-1$.
\end{defn}
We remark that  $C(P_{2i}\cdot P_{2i+2} \cdots P_{2j}
 \tin Q_{2i-1}) = Q_{2i-1,2j}$, by \eqref{pq14e}. Hence the 
$\beta_{2i-1,2j}$ are well defined irreducible characters of 
$Q_{2i-1,2j}$.
Similarly we see that the characters $\alpha_{2r,2s+1}$ are also well defined.

As we did with the corresponding groups, 
starting from the above basic properties we will describe the relations
 these characters satisfy.
Towards that direction we state and prove 
 \begin{proposition}\mylabel{pqpropo}
The following holds:
\begin{multline}\mylabel{pq14c}
\beta_{2i-1,2j} \in \Irr(Q_{2i-1,2j}) \text{ is the $P_{2j}$-Glauberman 
correspondent of } \\
\beta_{2i-1,2j-2} \in \Irr^{P_{2j}}(Q_{2i-1,2j-2}),  \text{ and lies 
above } \beta_{2i-3,2j}, \beta_{2i-5,2j},\dots,\beta_{1,2j},
\end{multline}
whenever $1\leq i\leq l$ and $i \leq j \leq k$. By convention we write 
$\beta_{2i-1,2i-2}:=\beta_{2i-1}$ when $j=i$.

Similarly, 
\begin{multline}\mylabel{pq13c}
\alpha_{2r,2s+1} \in \Irr(P_{2r,2s+1}) \text{ is the $Q_{2s+1}$-Glauberman 
correspondent of } \\
\alpha_{2r,2s-1} \in \Irr^{Q_{2s+1}}(P_{2r,2s-1}),  \text{ and lies over } 
\alpha_{2r-2,2s+1}, \dots,\alpha_{2,2s+1},
\end{multline}
whenever $1\leq r \leq k $ and $r \leq  s \leq l-1$. By convention we write 
 $\alpha_{2r,2r-1}:=\alpha_{2r}$ when $r=s$.

Therefore,
\begin{equation}\mylabel{pq14d}
\beta_{2i-1} \in \Irr(Q_{2i-1}|\beta_{2i-3,2i-2}, \beta_{2i-5,2i-2},
\dots,\beta_{1,2i-2}),
\end{equation}
and
\begin{equation}\mylabel{pq13d}
\alpha_{2r} \in \Irr(P_{2r} | \alpha_{2r-2,2r-1}, \alpha_{2r-4,2r-1},
\dots,\alpha_{2,2r-1}), 
\end{equation}
whenever  $i=1,\dots,l$ and $r=1,\dots,k$.
\end{proposition}

\begin{proof}
In view of  Definition \ref{pq1413def2}, it is easy to see that 
$\beta_{2i-1,2j}$ is the $P_{2j}$-Glauberman correspondent of 
$\beta_{2i-1,2j-2}$ (as the latter is the $P_{2i}\cdot P_{2i+2} \cdots  
P_{2j-2}$-Glauberman correspondent of $\beta_{2i-1}$),  for all $i,j$ with 
 $1\leq i \leq l$ and $i<j \leq k$.
We also remark that the same argument implies that 
$\beta_{2i-1,2j}$ is the $P_{2t}\cdot P_{2t+2}\cdots P_{2j}$-Glauberman
correspondent of $\beta_{2i-1,2t-2}$, for any $t$ with $1\leq i <t<j$. 
Furthermore, the same definition tells  us that $\beta_{2j-1,2j}$ is
 the $P_{2j}$-Glauberman correspondent
of $\beta_{2j-1,2j-2} = \beta_{2j-1}$, for all $j=1,\dots,k$.

Thus to prove \eqref{pq14c} it suffices to show  that $\beta_{2i-1,2j}$
lies over $\beta_{2i-3,2j},\dots,\beta_{1,2j}$,  for all $i,j$ with  $1\leq i 
\leq l$ and $i \leq j \leq k $.
For this we will use induction on $i$.
For $i=1$, it holds vacuously, since  the character $\beta_{2i-3,2j}$ 
doesn't exist.
The first interesting case appears when $i=2$. According to \eqref{x6}
the character $\beta_{3}$ lies above $\beta_{1,2}$. Therefore, for any 
$j=2,\dots,k$, 
the $P_4\cdots P_{2j}$-Glauberman correspondent $\beta_{3,2j}$ of $\beta_3$
lies above the $P_4 \cdots P_{2j}$-Glauberman  correspondent $\beta_{1,2j}$
of $\beta_{1,2}$.

For the inductive step the argument is similar. If 
  $i \geq 3$ and $\beta_{2i-3,2j}$ lies above 
$\beta_{2i-5,2j},\dots,\beta_{1,2j}$ for all $j=i-1,\dots,k$ then 
$\beta_{2i-3,2i-2}$ lies above $\beta_{2i-5,2i-2},\dots,\beta_{1,2i-2}$.
According to \eqref{x6}, the character $\beta_{2i-1}$ was picked 
to lie above $\beta_{2i-3,2i-2}$. Therefore the 
 $P_{2i} \cdot P_{2i+2}  \cdots P_{2j}$-Glauberman correspondent 
$\beta_{2i-1,2j}$ of $\beta_{2i-1}$ lies above the
$P_{2i} \cdot P_{2i+2} \cdots P_{2j}$-Glauberman correspondent
$\beta_{2i-3,2j}$ of $\beta_{2i-3,2i-2}$, for any $j$ with $j=i,\dots,k$.
Hence, $\beta_{2i-1,2j}$ lies above $\beta_{2i-3,2j}, \beta_{2i-5,2j},
\dots,\beta_{1,2j}$  whenever $j=i,\dots,k$. This 
completes the inductive argument on $i$, thus  proving \eqref{pq14c}.

As $\beta_{2i-1}$ lies above $\beta_{2i-3,2i-2}$ (by \eqref{x6}),
\eqref{pq14d} is an immediate consequence of \eqref{pq14c}.

The proof of \eqref{pq13c} and \eqref{pq13d} is similar.
\end{proof}

Looking at the ``character triangles ''  \eqref{pq11b} and \eqref{pq12b},
we can translate Proposition \ref{pqpropo} as follows:

  Every horizontal line in the triangles \eqref{pq11b} and \eqref{pq12b}
 (with the characters in  parenthesis included) 
 is formed by taking a character that is a Glauberman correspondent
 of the previous one. 
Also the vertical lines in these two triangles are formed by 
characters that are lying one above  the other.  
We can say even more:
\begin{proposition} \mylabel{pqremark1'}
For every $i,j,t$ with $1\leq i \leq t \leq l$ and $t\leq j \leq k$, 
the group $Q_{2t-1,2j}$  fixes the character $\beta_{2i-1,2j}$. 
Hence $\beta_{2i-1,2j}$ is the unique character in $\Irr(Q_{2i-1,2j})$ 
lying under  $\beta_{2t-1,2j}  \in \Irr(Q_{2t-1,2j})$.
In addition, for every $i, j$ with $1\leq i \leq j \leq l$, the group 
$Q_{2j-1}$  fixes the character $\beta_{2i-1,2j-2}$.
Hence $\beta_{2i-1,2j-2}$ is the unique character
 in $\Irr(Q_{2i-1,2j-2})$  lying under $\beta_{2j-1} \in \Irr(Q_{2j-1})$. 

Similarly,  for every $r, s, t$ with $1\leq r \leq t \leq k$  and $t \leq s \leq l-1$,  
the group $P_{2t, 2s+1}$  fixes the character $\alpha_{2r, 2s+1}$.
Therefore,  $\alpha_{2r,2s+1}$ is the unique character of $P_{2r,2s+1}$ that 
lies under $\alpha_{2t,2s+1} \in \Irr(P_{2t,2s+1})$.
In addition, for every $r, s$ with $1 \leq r \leq s \leq k$, the 
group $P_{2s}$  fixes the character $\alpha_{2r, 2s-1}$. Hence 
$\alpha_{2r, 2s-1}$  is the unique character in $\Irr(P_{2r, 2s-1})$ 
lying under $\alpha_{2s}  \in \Irr(P_{2s})$. 
\end{proposition}

\begin{proof}
Because of  symmetry it suffices to 
prove the proposition for the $q$-groups $Q_{2t-1,2j}$ and the characters $\beta_{2i-1,2j}$, for fixed $i,t,j$ in the range of  the proposition.

If $1\leq i \leq t \leq l$ and $t \leq j \leq k$, 
then equations  \eqref{pq:norm} and \eqref{pq:norm'} imply   that 
$Q_{2i-1,2j}$ is a normal subgroup of $Q_{2t-1,2j}$.
Equation \eqref{pq:norm}  also implies  that $Q_{2i-1,2j}\unlhd Q_{2j-1}$, whenever
 $1\leq i <  j \leq l$.  
Therefore,  according to Clifford's Theorem, 
it is enough to prove that $Q_{2t-1,2j}(\beta_{2i-1,2j}) = 
Q_{2t-1,2j}$ and $Q_{2j-1}(\beta_{2i-1,2j-2})= Q_{2j-1}$,  in order to 
complete the proof of Proposition  \ref{pqremark1'}.

In view of \eqref{x5}
 the group $Q_{2t-1}$ fixes $\beta_{2i-1}$.
According to \eqref{pq14'a}, the group $Q_{2t-1}$ normalizes the groups 
$P_{2i}, \dots,P_{2t-2}$. Hence its subgroup 
 $Q_{2t-1,2j} = C(P_{2t}, \dots,P_{2j} \tin Q_{2t-1})$
normalizes the groups $P_{2i}, \dots,P_{2t-2}$,  centralizes  
$P_{2t}, \dots, P_{2j}$, and fixes $\beta_{2i-1}$.
 Therefore $Q_{2t-1,2j}$ fixes $\beta_{2i-1,2j}$, which is the $P_{2i}, \dots, 
P_{2j}$-Glauberman correspondent of $\beta_{2i-1}$ by \eqref{pq14c}.
So $Q_{2t-1,2j}(\beta_{2i-1,2j}) = Q_{2t-1,2j}$

Similarly \eqref{x5} and \eqref{pq14'a} imply that $Q_{2j-1}$ 
fixes $\beta_{2i-1}$ and normalizes
 $P_{2i}, \dots, P_{2j-2}$, whenever $1\leq i < j \leq l$. Hence $Q_{2j-1}$ fixes $\beta_{2i-1,2j-2}$.
So $Q_{2j-1}(\beta_{2i-1,2j-2}) = Q_{2j-1}$ and the proposition follows.
\end{proof}

%%% Local Variables: 
%%% mode: latex
%%% TeX-master: "t-t"
%%% End: 

\newpage

\section{From towers to triangles}

We are now ready to prove one direction of the correspondence
 in Theorem \ref{cor:t}. In particular, we will prove that 
for any character tower of \eqref{*} there is a corresponding 
$G(\chi_1, \dots, \chi_m)$-conjugacy class of triangular sets for \eqref{*}. The explicit relation 
between the character towers and the triangular sets   is
described  in Theorem \eqref{tow--tri} below.
Before we give the inductive proof of that theorem, we will demonstrate,
  for clarity, how the correspondence works in  the special  cases  
where $m=1,2,3$.
 
We begin with a lemma that is an easy application of Theorems \ref{dadet3} and 
\ref{dade:t2}.
\begin{lemma} \mylabel{dade'}
Let $G$ be a finite group of odd order,  and $\pi$ be any set of primes.
 Suppose that  $N, K_1 ,K_2, \dots,K_r$ are normal subgroups of 
$G$,  for some $r\geq 1$,  such that  $N\unlhd K_1 \unlhd K_2 \unlhd  \dots
\unlhd  K_r$.
 Assume further that $N= A\ltimes B$,  where $B$ is a normal 
$\pi'$-subgroup of $G$,  and $A$ is  any $\pi$-subgroup of $N$.
Let $\chi \in  \Irr(N)$  be a $\pi$-factorable character of $N$.
Assume that   
$\chi = \alpha \cdot \beta^e$ is the decomposition of $\chi$ to its $\pi$- and 
$\pi'$-special parts respectively,  where  $\beta^e$ is the canonical  extension to 
$N$ of an irreducible $A$-invariant character $\beta \in \Irr^A(B)$.
 Let   $K_i(\chi)$ be  the stabilizer of $\chi$ in $K_i$, for $i=1,\dots,r$, 
 and  $C = C(A \tin B)$ be  the centralizer 
of $A$ in $B$. 

Then there is a one-to-one correspondence
between the  character towers $\{ \chi,\chi_1,\dots,\chi_r \}$
of  the series  $N\unlhd K_1 \unlhd K_2 \unlhd  \dots
\unlhd  K_r$,  starting with  $\chi$,  and the character towers
$\{ \alpha\times  \gamma, \Psi_1,\dots,\Psi_r \}$ 
of the series
$N(A \tin N)= C \times A \unlhd N(A \tin K_1(\chi)) \unlhd
N(A \tin  K_2(\chi))\unlhd   \dots \unlhd N(A \tin  K_r(\chi))$,  starting with  
$\alpha \times \gamma \in \Irr(C \times A)$, where  $\gamma \in \Irr(C)$ 
is the $A$-Glauberman  
correspondent of $\beta \in \Irr^{A}(B)$.
Furthermore, for any  subgroup  $M$ of $N(A \tin G)$ we have 
$$M(\chi,\chi_1,\dots,\chi_r) = M(\alpha\times  \gamma, \Psi_1,\dots,\Psi_r).$$
\end{lemma}

\begin{proof}
Let 
\begin{equation}\mylabel{ll1}
\{ \chi, \chi_1,\dots,\chi_{r}\},
\end{equation}
 be a character tower of the normal series 
$N\unlhd K_1 \unlhd K_2 \unlhd  \dots
\unlhd  K_r$,  starting with  $\chi$. 
According to Clifford's theorem, for every $i=1,\dots,r$ there exists
a unique  irreducible character $\chi_i^* \in \Irr(K_i(\chi))$ that induces $\chi_i \in \Irr(K_i)$ 
and lies above $\chi \in \Irr(N)$.
Furthermore, the characters 
\begin{equation}\mylabel{ll2}
\{\chi, \chi_1^*,\dots,\chi_r^*\},
\end{equation}
 form a tower for the normal series 
$N \unlhd K_1(\chi) \unlhd \dots \unlhd K_r(\chi)$.
Hence \eqref{ll1} corresponds to \eqref{ll2}.              
Clifford's Theorem also implies    that
 this correspondence between \eqref{ll1} and \eqref{ll2} 
is invariant under any subgroup of $G(\chi)$. So,  in particular,   
\begin{equation}\mylabel{ll5}
G(\chi, \chi_1,\dots, \chi_r) = G(\chi,\chi_1^*,\dots,\chi_r^*).
\end{equation}
To complete the proof of the lemma we only need to observe that 
Theorem \ref{dade:t2} can be applied to  the tower  \eqref{ll2} 
and the normal series $N= A \rtimes B 
\unlhd K_1(\chi) \unlhd \dots \unlhd K_r(\chi)$.
Note also that, in view of Theorem \ref{dadet3}, the $A$-correspondent of the irreducible 
 character
 $\chi = \alpha \cdot \beta^e \in \Irr(N)$ is ofthe form 
$\chi_{(A)}= \alpha \times \gamma \in \Irr(A \times C)$, where  $\gamma \in \Irr(C)$
is the $A$-Glauberman correspondent of $\beta \in \Irr^A(B)$.    
Hence  the character tower \eqref{ll2} has a unique $A$-correspondent 
character tower 
\begin{equation}\mylabel{ll3}
\{\alpha\times  \gamma, \Psi_1,\dots,\Psi_r \}
\end{equation}
 of the series 
$$
A \times C    \unlhd N(A \tin K_1(\chi)) \unlhd \dots \unlhd N(A \tin K_r(\chi)).
$$
 This way we have created a correspondence, that is 
the combination of a Clifford and an $A$-correspondence,  between 
the tower \eqref{ll1} and the \eqref{ll3}.

Furthermore, as $\gamma$ is the
$A$-Glauberman correspondent of $\beta$ and $B$ is normal in $G$, 
we have that $M(\beta) = M(\gamma)$,
 for any group $M$ with 
$M \leq N(A \tin G)$.
Also $M(\beta) = M(\beta^e)$ as $N \unlhd G$. 
 Hence we conclude that  
\begin{equation}\mylabel{ll6}
M(\chi)= M(\alpha \cdot \beta^e ) = M(\alpha , \beta)=
 M(\alpha \times \gamma).
\end{equation}

Furthermore, Theorem \ref{dade:t2}  implies that for all $M$ with  
 $M \leq N(A \tin G)$ we have 
\begin{equation}\mylabel{ll4}
M(\chi,\chi_1^*,\dots,\chi_r^*) = M(\alpha\times  \gamma, \Psi_1,\dots,\Psi_r),
\end{equation}
as $M$ normalizes $N, K_1, \dots, K_r$.
Therefore we have 
\begin{align*}
M(\chi, \chi_1,\dots,\chi_r)&=M(\chi,\chi_1^*,\dots,\chi_r^*) &\quad &\text{ by \eqref{ll5}}\\
&=M(\chi)(\chi,\chi_1^*,\dots,\chi_r^*) & &\\
&=M(\chi)(\alpha\times  \gamma, \Psi_1,\dots,\Psi_r) 
& &\text{ by \eqref{ll4}}\\
&=M(\alpha \times \gamma)(\alpha\times  \gamma, \Psi_1,\dots,\Psi_r) &
 &\text{ by \eqref{ll6}}\\
&=M(\alpha\times  \gamma, \Psi_1,\dots,\Psi_r). & &
\end{align*}

This completes the proof of the lemma.
\end{proof}

\begin{defn}\mylabel{def.t-t}
For the rest of this thesis, 
the correspondence between  towers
$$\{ \chi=\alpha \cdot \beta^e,
 \chi_1,\dots,\chi_r\} \leftrightarrow 
\{\alpha\times  \gamma, \Psi_1,\dots,\Psi_r \}
$$
 that  is described in Lemma \ref{dade'}, 
will be called a { \em $cA$-correspondence }
(Clifford-$A$).  We call  the tower $\{\alpha\times  \gamma, \Psi_1,\dots,
\Psi_r \}$ the  { \em $cA$-correspondent } of  $\{ \chi=\alpha \cdot \beta^e,
 \chi_1,\dots,\chi_r\}$.
Similarly, we call $\Psi_i$  the { \em $cA$-correspondent } 
of $\chi_i$, for all $i=1,\dots, r$. 
\end{defn}

We can now look at the cases $m=1,2, 3$.
If $m=1$ then the normal series 
\eqref{*} consists of the groups $1=G_0 \unlhd G_1 \unlhd G$.
So   any character tower  $\{1=\chi_0, \chi_1\}$,
 of this series determines  the   triangular set  
$\{1=P_0, Q_1=G_1|1=\alpha_0, \beta_1 =\chi_1\}$.
Furthermore, assume that 
\begin{equation}\mylabel{n*}
1=G_0\unlhd G_1\unlhd \dots \unlhd G_n \unlhd G
\end{equation}
is a normal series of $G$,  for some $n \geq m=1$,   that extends the series
 $1=G_0 \unlhd G_1 \unlhd G$. Assume further that we have an extension of the character 
tower $\{1=\chi_0, \chi_1\}$ to a character tower   $\{1=\chi_0,\chi_1,\dots, \chi_n\}$
for the series \eqref{n*},  so that Hypothesis \ref{hyp1} holds.
As $\chi_1 = \beta_1$, we have that $G_i(\chi_1) = G_i(\beta_1)$
for all $i=1,\dots,n$. Hence we can define the groups 
\begin{subequations}
\begin{equation}\mylabel{i,1}
G_{i,1} := G_i(\beta_1)=G_i(\chi_1) = N(Q_1 \tin G_i(\chi_1)),
\end{equation}
where the last equality holds as $Q_1 = G_1 \unlhd G$.
By convention, whenever  we have a series as in \eqref{n*}, we will write  
\begin{equation}\mylabel{infy}
G_{\infty}=G.
\end{equation}
\end{subequations}
With this notation, we can also write 
$\Gap{1}$ for the stabilizer $G(\chi_1)=G(\beta_1)$. 
Therefore 
the series 
\begin{equation}\mylabel{n1*}
1=G_{0,1} \unlhd G_{1,1} \unlhd \dots \unlhd G_{n,1} \unlhd \Gap{1}, 
\end{equation}
is a normal series of $\Gap{1}$. Furthermore, for any $i=1,2,\dots,n$, 
Clifford's theorem applied to the groups $G_1 \unlhd G_i$ and the characters 
$\chi_1, \chi_i$ implies the existence of a unique character 
$\chi_{i,1} \in \Irr(G_i(\chi_1))$
that lies above $\chi_1$ and induces $\chi_i$, i.e, 
\begin{equation}\mylabel{i.1}
\chi_{i, 1} \in \Irr(G_i(\chi_1)) \text{ is the $\chi_1$-Clifford correspondent of  } 
\chi_{i} \in \Irr(G_{i} | \chi_1). 
\end{equation}
Note that  $\chi_{1,1} = \chi_1$.
We write $\chi_{0,1} =1$.
Then it is clear that $\chi_{i,1}$ lies above $\chi_{k,1}$ whenever
 $1\leq k \leq i \leq n$.
This way we have created a tower  $\{ 1=\chi_0, \chi_{1,1}= \chi_1,
\chi_{2,1},  \dots, \chi_{n,1}\}$ for the series \eqref{n1*}, 
 fully determined by  the character tower 
 $\{1=\chi_0,\chi_1,\dots, \chi_n\}$. 
Furthermore, Clifford's theorem implies that 
for   any subgroup $M$ of $G = N(Q_1 \tin G)$
we have 
\begin{equation}\mylabel{pqn0}
M(\chi_1,\chi_2,\dots, \chi_k ) = M(\chi_{1,1},\chi_{2,1},\dots, \chi_{k,1}),
\end{equation}
for any $k=1,2,\dots,n$.
Therefore, in the case where $m=1$, in addition to the correspondence
between towers and triangular sets for \eqref{*}, we proved that 
any tower of \eqref{n*} determines a unique tower of \eqref{n1*}.
This is a  property that, as we will see in Theorem \ref{tow--tri},
 carries  over to  every $m$. 
By convention, we  write this first correspondence as a $cQ_1$-correspondence
(even though it is a Clifford correspondence). Table \ref{diagr.1} describes exactly 
the above relations.
\begin{table}[h]
\begin{diagram}[small]
G=G_{\infty} & & & &\Gap{1}:=G(\chi_1) & & &\\
\vLine   &\qquad   &     &   &\vLine    &\qquad  &    &       \\
  G_n    &\qquad   &\chi_n  &   &G_{n,1}:=G_n(\chi_1) 
 &\qquad  &\chi_{n,1} & \\
     \vLine   &\qquad   &\vLine     &   &\vLine    &\qquad  &\vLine    &                \\
     \vdots  &\qquad   &\vdots     &   &\vdots    &\qquad  &\vdots   &                \\
\vLine   &\qquad   &\vLine     &   &\vLine    &\qquad  &\vLine    &                \\
  G_3    &\qquad   &\chi_3  &    &G_{3,1}:=   G_3(\chi_1)   &\qquad 
    &\chi_{3,1} &    \\
     \vLine    &\qquad  &\vLine   &\text{\qquad  $\overleftrightarrow{cQ_1}$ \qquad} &\vLine 
       &\qquad  &\vLine   &\qquad  & \\
   G_2  &\qquad   &\chi_2  &   
   &G_{2,1}:= G_2(\chi_1)   &\qquad  &\chi_{2,1}  & \\
    \vLine    &\qquad  &\vLine     & &\vLine    &\qquad  &\vLine   & \\
    G_1=Q_1  &\qquad   &\chi_1 =\beta_1 &  &G_{1,1}= G_1= Q_1 &\qquad
  &\chi_{1,1}=\chi_1=\beta_1 &  \\
\vLine    &\qquad  &\vLine     & &\vLine    &\qquad  &\vLine   & \\
     G_0=1   &\qquad   &\chi_0=1 & &G_{0,1} = 1 &\qquad &\chi_{0,1} = 1 &\quad 
  \end{diagram}
\caption{The $cQ_1$-correspondence.}
\mylabel{diagr.1}
\end{table}

The first interesting case appears when  $m=2$.
Here the normal series \eqref{*} consists of the groups
$1=G_0\unlhd G_1\unlhd G_2$. Let  
\begin{equation}\mylabel{2*}
\{1=\chi_0,\chi_1,\chi_2\}
\end{equation}
be a character tower of that series.
We  have already seen (from the case $m=1$) that 
the subtower $\{1=\chi_0, \chi_1\}$ of \eqref{2*}, 
determines the triangular set $\{P_0=1, Q_1=G_1 | \alpha_0=1, \beta_1=\chi_1 \}$.
We  expand this set  by picking the $\pi$-group $P_2$ to be 
any $\pi$-Hall subgroup of $G_2(\chi_1)= G_2(\beta_1)$. (Note that to get
 a triangular set for the series $1=G_0\unlhd G_1\unlhd G_2$, it  is enough to expand 
the existing   set, $\{P_0=1, Q_1=G_1 | \alpha_0=1, \beta_1=\chi_1 \}$, by 
a $\pi$-group $P_2$ along with  its irreducible character $\alpha_2$, so that 
\eqref{xxx} will be  valid for this new set.)
Before we see how to chose the desired irreducible character $\alpha_2 \in \Irr(P_2)$, we observe the following:
as $G_2/G_1$ is a $\pi$-group and  $P_2 \in \Hall_{\pi}(G_{2}(\chi_1))$, we   
have  that $P_2$   covers $G_2(\beta_1) = G_2(\chi  _1) $  modulo $G_1$. Hence,
\begin{equation} \mylabel{pq3b}
G_{2,1}= G_2(\chi_1)= G_2(\beta_1) = P_2 \ltimes G_1 = P_2 \ltimes Q_1.
\end{equation}

In view of the work we did in the case $m=1$,  we have that,  
for  every $n$ with $n \geq 2$,  the character tower  
$\{\chi_0=1, \chi_1 = \beta_1, \chi_2, \dots , \chi_n \}$ of the normal series 
\eqref{n*}, extending   the tower \eqref{2*}, 
 has a unique $Q_1$-correspondent  character tower 
$\{ \chi_{0,1}=1, \chi_{1,1} = \beta_1, \chi_{2,1}, \dots , \chi_{n,1} \}$,
of the series \eqref{n1*} (see Table \ref{diagr.1}).

Furthermore, equation \eqref{pq3b} permits us to  apply Lemma \ref{lemma0} to the groups 
$G_{1,1} = G_1=Q_1$, $  G_{2,1}= G_{2}(\chi_1) $,  $ \Gap{1} $ and the character 
$\chi_1=\beta_1$ in the place of the groups 
$N, H, G$ and the character $\theta$, respectively. (Note that in this case 
 $H(\theta) = H $). Thus we conclude that 
$\chi_1 = \beta_1$ has a unique canonical extension 
$\beta_1^e \in \Irr(G_{2,1})$. As $\chi_{2,1} \in \Irr(G_{2,1})$ lies above $\chi_1 = \beta_1$, 
Lemma \ref{lemma0}  also implies 
that   there is a unique character
 $\alpha_2 \in \Irr(P_2)$ such that
\begin{subequations}
\begin{equation}\mylabel{pq3c}
\chi_{2,1}= \alpha_2 \cdot \beta_1^e, 
\end{equation}
while \begin{equation}\mylabel{pq4a}
\Gap{1}(\alpha_2, \beta_1) = N(P_2 \tin \Gap{1}(\beta_1,\chi_{2,1})).
\end{equation}
But  
\begin{equation}\mylabel{pq4a'}
G_{\infty,1}(\beta_1) = G_{\infty,1} = G(\beta_1)=G(\chi_1).
\end{equation}
 Furthermore,
\eqref{pqn0} for $k=2$ implies that 
$$
G_{\infty}(\chi_1,\chi_2) = G_{\infty}(\chi_{1,1},\chi_{2,1})=
G_{\infty}(\chi_{1},\chi_{2,1}).
$$
Therefore, \eqref{pq4a} and \eqref{pq4a'} imply 
\begin{multline}\mylabel{pq3d}
G(\alpha_2, \beta_1) = N(P_2 \tin \Gap{1}(\chi_{2,1}))= 
 N(P_2 \tin \Gap{1}(\chi_{2}))\\
= N(P_2 \tin G_{\infty}(\chi_1,\chi_{2,1}))=
N(P_2 \tin G_{\infty}(\chi_1,\chi_{2}))= N(P_2 \tin G(\chi_1,\chi_{2})).
\end{multline}
Hence, by intersecting both sides with $G_i$,  we get
\begin{equation}\mylabel{pq3e}
G_{i}(\alpha_2, \beta_1) = N(P_2 \tin G_{i,1}(\chi_{2,1}))=N(P_2 \tin G_{i,1}(\chi_{2}))=
N(P_2 \tin G_i(\chi_1,\chi_2)),
\end{equation}
whenever $i=0,1,\dots,n,\infty$.
\end{subequations}
As $P_2$ was picked to be a $\pi$-Hall subgroup of 
$G_2(\beta_1)$, we obviously have that the set $\{ P_0=1, P_2, Q_1 |\alpha_0=1,
 \alpha_2,\beta_1 \}$ is a triangular set for the series $1=G_0 \unlhd G_1 \unlhd G_2$.
   Hence \eqref{prod4} for $t=u=1$ and \eqref{pq-15a} for $r=1$ imply 
\begin{align*}
G_2(\alpha_2, \beta_1)&=P_2 \times Q_{1,2},\\
G_1(\alpha_2, \beta_1) &= 1\times Q_{1,2} = Q_{1,2}.
\end{align*}
 
Furthermore, we have a correspondence similar to the one described in 
Table \ref{diagr.1}. Indeed, in view of \eqref{pq3b}
 and \eqref{pq3c}, the normal series 
$ G_{2,1} \unlhd G_{3,1} \unlhd  \dots \unlhd G_{n,1}\unlhd G_{\infty,1}$ ,
 along with the $\pi$-factorable character $\chi_{2,1}$,
satisfies the hypotheses of Lemma \ref{dade'}. Hence, 
there is a $cP_2$-correspondence between the character towers 
of the above  series  and those of the series  
$G_{2,2} \unlhd G_{3,2} \unlhd \dots \unlhd G_{n,2} \unlhd \Gap{2}$, where
\begin{equation}\mylabel{i,2} 
G_{i,2}:= N(P_2 \tin G_{i,1}(\chi_{2,1}))
\end{equation}
 for  all $i=2,\dots,n,\infty$.
Thus, the tower $\{\chi_{2,1}, \chi_{3,1}, \dots, \chi_{n,1}\}$ has a 
$cP_2$-correspondent tower $\{ \chi_{2,2},\chi_{3,2}, \dots, \chi_{n,2} \}$, 
where $\chi_{i,2} \in \Irr(G_{i,2})$, for all $i=2,\dots,n$.
Furthermore,
for any $M$ with $M \leq N(P_2 \tin \Gap{1})= N(P_2 \tin G(\chi_1))$, we have 
\begin{equation}\mylabel{pqn1}
M(\chi_{2,1}, \dots, \chi_{k,1}) = M(\chi_{2,2},\dots,\chi_{k,2}),
\end{equation}
  whenever $2\leq k \leq n $.
The same lemma describes $G_{2,2}$ as well as   $\chi_{2,2}$. So we get that 
\begin{equation}\mylabel{pqn2}
\begin{aligned}
G_{2,2}&= P_2 \times Q_{1,2},\\ 
\chi_{2,2}&= \alpha_2 \times \beta_{1,2},
\end{aligned}
\end{equation}
where $Q_{1,2}= N(P_2 \tin Q_1) = C(P_2 \tin Q_1)$ (see \eqref{anno1}),
and $\beta_{1,2} \in \Irr(Q_{1,2})$ is the $P_2$-Glauberman correspondent 
of $\beta_1 \in \Irr^{P_{2}}(Q_1)$ (see Definition \ref{pq1413def2}).

We  observe that the earlier definition of $G_{i,2}$ (see \eqref{i,2}), works also for $i=1$, 
as $G_{1,1} \leq G_{2,1}$ fixes $\chi_{2,1}$. So, 
   $G_{1,2} :=N(P_2 \tin G_{1,1}) =Q_{1,2}$ while 
the character $\beta_1= \chi_{1,1} \in \Irr(G_1) =\Irr(G_{1,1})$ 
has as a unique $P_2$-Glauberman correspondent
the character  $\beta_{1,2} \in \Irr(Q_{1,2})$.
This, combined with the former $cP_2$-correspondence, provides
a correspondence (that we also  write as a $cP_2$-correspondence) 
between the character towers of the series \eqref{n1*} 
and those of the series
\begin{equation}\mylabel{n2*}
G_{0,2} =1\unlhd G_{1,2}:=N(P_2 \tin G_{1,1})=Q_{1,2} \unlhd 
G_{2,2}\unlhd \dots \unlhd G_{n,2}\unlhd \Gap{2},
\end{equation}
 described in the  Table \ref{diagr.2}.
We remark here that, for every group $M$ with $M \leq N(P_2 \tin G)$, the 
$P_2$-Glauberman correspondence  between $\Irr^{P_2}(Q_1)$ and $\Irr(Q_{1,2})$, 
(with the character  $\chi_1 = \chi_{1,1} = \beta_1$ in the former set corresponding to 
the character 
 $\chi_{1,2}= \beta_{1,2}$), is $M$-invariant. Hence 
\begin{equation}\mylabel{pqn1'}
M(\chi_1) = M(\chi_{1,2}).
\end{equation}
\begin{table}[h]
\begin{footnotesize}
\begin{diagram}[small]
G=G_{\infty,1} & & & &G_{\infty, 2}:=N(P_2 \tin \Gap{1}(\chi_{2,1})) & & & \\
\vLine   &\qquad   &     &   &\vLine    &\qquad  &    &                \\ 
  G_{n,1}    &\qquad   &\chi_{n,1}  &   &G_{n,2}:=N(P_2 \tin G_{n,1}(\chi_{2,1})) 
 &\qquad  &\chi_{n,2} & \\
     \vLine   &\qquad   &\vLine     &   &\vLine    &\qquad  &\vLine    &                \\
     \vdots  &\qquad   &\vdots     &   &\vdots    &\qquad  &\vdots   &                \\
\vLine   &\qquad   &\vLine     &   &\vLine    &\qquad  &\vLine    &                \\
  G_{3,1}    &\qquad   &\chi_{3,1}  &    &G_{3,2}:= N(P_2 \tin  G_{3,1}(\chi_{2,1}))  
 &\qquad     &\chi_{3,2} &    \\
     \vLine    &\qquad  &\vLine   &\text{\qquad  $\overleftrightarrow{cP_2}$ \qquad} &\vLine 
       &\qquad  &\vLine   &\qquad  & \\
   G_{2,1}=P_2\ltimes Q_1  &\qquad   &\chi_{2,1}=\alpha_2 \cdot \beta_1^e  &   
   &G_{2,2}= P_2 \times Q_{1,2}   &\qquad  &\chi_{2,2}=\alpha_2 \times \beta_{1,2}  & \\
    \vLine    &\qquad  &\vLine     & &\vLine    &\qquad  &\vLine   & \\
    G_{1,1}=G_1=Q_1  &\qquad   &\chi_{1,1} =\chi_1=\beta_1 &  &G_{1,2}:= Q_{1,2} &\qquad  &\beta_{1,2} &\\
 \vLine    &\qquad  &\vLine     & &\vLine    &\qquad  &\vLine   & \\
     G_{0,1}=1   &\qquad   &\chi_{0,1}=1 & &G_{0,2} = 1 &\qquad &\chi_{0,2} = 1 &
  \end{diagram}
\caption{The $cP_2$-correspondence}
\mylabel{diagr.2}
\end{footnotesize}
\end{table}
\newpage

According to \eqref{pq3e}, we have that 
\begin{equation}\mylabel{pqn2a}
G_{i,2} = N(P_2 \tin G_{i,1}(\chi_{2,1}))=N(P_2 \tin G_{i}(\chi_1,\chi_2))= G_i(\alpha_2,\beta_1),
\end{equation}
whenever $i=1,\dots , n ,\infty$.
This, along with Tables  \ref{diagr.1} and \ref{diagr.2}, implies 
a $cQ_1, cP_2$-correspondence described in the diagram that follows:
\begin{table}[h]
\begin{scriptsize}
\begin{diagram}[small]
G=G_\infty  & &  & &\Gap{1}=G(\beta_1) & & & &\Gap{2}=G(\alpha_2,\beta_1) & & &\\
\vLine &\qquad & &\qquad  &\vLine   &\qquad   &     &   &\vLine    &\qquad  & & \\
G_n &\qquad &\chi_n &\qquad   &G_{n,1}=G_n(\beta_1)    &\qquad   &\chi_{n,1}  &   
&G_{n,2}=G_n(\alpha_2,\beta_1) 
 &\qquad  &\chi_{n,2} & \\
\vLine &\qquad &\vLine &\qquad  &\vLine   &\qquad   &\vLine     &   &\vLine    &\qquad  &\vLine    &                \\
\vdots &\qquad  &\vdots &    &\vdots  &\qquad   &\vdots     &   &\vdots    &\qquad  &\vdots   &                \\
\vLine &\qquad &\vLine &\qquad  &\vLine   &\qquad   &\vLine     &   &\vLine    &\qquad  &\vLine    &                \\
G_3 & &\chi_3 &   &G_{3,1}=G_3(\beta_1)    &\qquad   &\chi_{3,1}  &   
 &G_{3,2}= G_3(\alpha_2,\beta_1) 
 &\qquad     &\chi_{3,2} &    \\
     \vLine &\qquad &\vLine &\text{\qquad  $\overleftrightarrow{cQ_1}$ \qquad} 
&\vLine    &\qquad  &\vLine   &\text{\qquad  $\overleftrightarrow{cP_2}$ \qquad} &\vLine 
       &\qquad  &\vLine   &\qquad  & \\
G_2 & &\chi_2 &    &G_{2,1}=P_2\ltimes Q_1  &\qquad   &\chi_{2,1}=\alpha_2 \cdot \beta_1^e  &   
   &G_{2,2}= P_2 \times Q_{1,2}   &\qquad  &\chi_{2,2}=\alpha_2 \times \beta_{1,2}  & \\
    \vLine &\qquad &\vLine &\qquad  &\vLine    &\qquad  &\vLine     & &\vLine    &\qquad  &\vLine   & \\
   G_1=Q_1 & &\chi_1 &  &G_{1,1}=Q_1  &\qquad   &\chi_{1,1} =\chi_1=\beta_1 & 
 &G_{1,2}= Q_{1,2} &\qquad  &\beta_{1,2} &\\
 \vLine &\qquad &\vLine &\qquad  &\vLine    &\qquad  &\vLine    
 & &\vLine    &\qquad  &\vLine   & \\
G_0=1 &  & \chi_0=1 &    &G_{0,1}=1   &\qquad   &\chi_{0,1}=1 & &G_{0,2} = 1 &\qquad &\chi_{0,2} = 1 &
  \end{diagram}
\caption{The $cQ_1,cP_2$-correspondence}
\mylabel{diagr.1-2}
\end{scriptsize}
\end{table}

\newpage
Furthermore,  for any group $M$ with 
$M \leq N(P_2 \tin G)$ we have 
\begin{align*}
M(\chi_1,\chi_2,\dots,\chi_k) &= M(\chi_{1}, \chi_{2,1}, \dots, \chi_{k,1})
 &\quad &\text{by  \eqref{pqn0} }\\
&= M(\chi_1)(\chi_{2,1}, \dots, \chi_{k,1}) & &\\
&=M(\chi_1)(\chi_{2,2}, \dots, \chi_{k,2}) & &\text{ by \eqref{pqn1} }\\
&=M(\chi_{1,2})(\chi_{2,2}, \dots, \chi_{k,2}) & &\text{ by \eqref{pqn1'}}\\ 
&=M(\chi_{1,2},\chi_{2,2}, \dots, \chi_{k,2}). & &
\end{align*}
Hence 
\begin{equation}\mylabel{pqn1-2}
M(\chi_1,\chi_2,\dots,\chi_{k}) = M(\chi_{1,1}, \chi_{2,1}, \dots, \chi_{k,1})
= M(\chi_{1,2}, \chi_{2,2}, \dots,\chi_{k,2}),
\end{equation}
whenever $1\leq k \leq n$.

The case $m=3$ is quite similar to $m=2$, so we will only  describe the main steps.
 We first  pick the  $\pi'$-group, $Q_3$,  
 as any 
\begin{equation}\mylabel{pq3'}
Q_3 \in \Hall_{\pi'}(G_{3,2})= \Hall_{\pi'}(G_3(\alpha_2 , \beta_1)). 
\end{equation}
Therefore we get that 
\begin{subequations} \mylabel{pq5}
\begin{equation} \mylabel{pq5a}
G_{3,2}= G_3(\alpha_2, \beta_1) = Q_3 \ltimes P_2.
\end{equation}
  Even more, 
 \eqref{pq5a} and \eqref{pq3e} for $i=3$ imply:
\begin{equation} \mylabel{pq5b}
G_{3,2}= G_3(\alpha_2, \beta_1) = Q_3 \ltimes P_2 = N(P_2 \tin G_3(\chi_1, \chi_2)).
\end{equation}
\end{subequations} 

To pick the character $\beta_3 \in \Irr(Q_3)$, we follow the same steps  
as we did for the character $\alpha_2$. So we  apply Lemma \ref{lemma0}  
to the groups $P_2$,  $  Q_3 \ltimes P_2=G_{3,2}$,  $ \Gap{2}$ and the 
character $\alpha_2$ in the place of the groups $N, H, G$ and the
character $\theta$ respectively.
Thus we conclude that 
that   there is a unique character
 $\beta_3 \in \Irr(Q_3)$ such that
\begin{subequations}
\begin{equation}\mylabel{pq6c}
\chi_{3,2}= \alpha_2^e \cdot \beta_3, 
\end{equation}
where $\alpha_2^e \in \Irr(G_{3,2})$ is the canonical extension of $\alpha_2 \in \Irr(P_2)$  
to  $Q_3 \ltimes P_2 =G_{3,2}$.
We also have that 
\begin{equation}\mylabel{pq4an}
\Gap{2}( \beta_3) = N(Q_3 \tin \Gap{2}(\chi_{3,2})), 
\end{equation}
and thus, in view of \eqref{pqn2a},
\begin{multline}\mylabel{pq4bn}
G(\beta_1,\alpha_2, \beta_3) =  N(Q_3 \tin \Gap{2}(\chi_{3,2}))\\
=N(Q_3 \tin G(\beta_1, \alpha_2)(\chi_{3,2}))=
N(Q_3 \tin N(P_2 \tin G(\chi_1,\chi_2))(\chi_{3,2})). 
\end{multline}
\end{subequations}
But \eqref{pqn1-2}, applied twice (with  $k=2$ and  $M=N(P_2 \tin G)$ for the first equality,
 and $k=3$  and  $M=N(P_2 \tin G)$ for the second one),  implies that 
$$
N(P_2 \tin G(\chi_1,\chi_2))(\chi_{3,2})
=N(P_2 \tin G(\chi_{1,2},\chi_{2,2}))(\chi_{3,2})
=N(P_2 \tin G(\chi_1,\chi_2))(\chi_{3}).
$$
Hence we conclude that
\begin{subequations} 
\begin{multline}\mylabel{pq9a}
G(\beta_1,\alpha_2, \beta_3) =  N(Q_3 \tin \Gap{2}(\chi_{3,2}))\\
=N(Q_3 \tin N(P_2 \tin G(\chi_1,\chi_2, \chi_{3}))=N(P_2, Q_3 \tin G(\chi_1,\chi_2,\chi_3)), 
\end{multline}
and thus
\begin{equation}\mylabel{pq9b}
G_i(\beta_1,\alpha_2, \beta_3)=N(Q_3 \tin G_{i,2}(\chi_{3,2})) 
= N(P_2, Q_3 \tin G_i(\chi_1,\chi_2,\chi_3)), 
\end{equation}
\end{subequations}
 for all $i=0, 1,\dots, n, \infty$.

Note that  $\{ P_0,P_2,Q_1,Q_3|1,\alpha_2,\beta_1,\beta_3\}$ 
is a triangular set. Indeed, as $\chi_{3,2} = \alpha_2^e \cdot \beta_3$
lies above $\chi_{2,2} = \alpha_2 \times \beta_{1,2}$ (see \eqref{pqn2}),
 we conclude that $\beta_3$ lies above the $P_2$-Glauberman correspondent 
$\beta_{1,2}$ of $\beta_1$. This, along with \eqref{pq3'} and the fact that 
$\{P_0, P_2, Q_1| 1,\alpha_2,\beta_1 \}$ is a triangular set, implies that 
 $\{ P_0,P_2,Q_1,Q_3|1,\alpha_2,\beta_1,\beta_3\}$ satisfies
 \eqref{xxx}, and thus is a triangular set.

To expand Table \ref{diagr.1-2} by one more step (that  will be a
 $cQ_3$-correspondence) we will apply (as we did for the
 $cP_2$-correspondence), Lemma \ref{dade'} to the last normal series 
of $\Gap{2}$ that the above table reaches. Notice that the normal series 
$ G_{3,2} \unlhd G_{4,2} \unlhd  \dots \unlhd G_{n,2}\unlhd G_{\infty,2}$ ,
 along with the $\pi$-factorable character $\chi_{3,2}$,
satisfies the hypotheses of Lemma \ref{dade'}. Hence 
there is a $cQ_3$-correspondence between the character towers 
of the above  series  and those of the series  
$G_{3,3} \unlhd G_{4,3} \unlhd \dots \unlhd G_{n,3} \unlhd \Gap{3}$, where 
\begin{equation}\mylabel{i,3}
G_{i,3}:= N(Q_3 \tin G_{i,2}(\chi_{3,2})), 
\end{equation}
 for  all $i=3,\dots,n,\infty$.
Assume that the tower 
 $\{ \chi_{3,3},\chi_{4,3}, \dots, \chi_{n,3} \}$ is  the $cQ_3$-correspondent
of the tower  $\{\chi_{3,2}, \chi_{4,2}, \dots, \chi_{n,2}\}$, 
where $\chi_{i,3} \in \Irr(G_{i,3})$ for all $i=3,\dots,n$.
Furthermore,
for any $M$ with $M \leq N(Q_3 \tin \Gap{2})= N(Q_3 \tin N(P_2 
\tin G(\chi_1,\chi_2)))$ we have 
\begin{equation}\mylabel{pqn3}
M(\chi_{3,2}, \dots, \chi_{k,2}) = M(\chi_{3,3},\dots,\chi_{k,3}),
\end{equation}
whenever $3\leq k \leq n$.
Furthermore,
\begin{equation}\mylabel{pqn3a}
\begin{aligned}
G_{3,3}&= P_{2,3} \times Q_{3},\\ 
\chi_{3,3}&= \alpha_{2,3} \times \beta_{3},
\end{aligned}
\end{equation}
where $P_{2,3}= N(Q_3  \tin P_2)$, and 
 $\alpha_{2,3} \in \Irr(P_{2,3})$ is the $Q_3$-Glauberman correspondent 
of $\alpha_2 \in \Irr^{Q_3}(P_{2})$.

We expand the definition of the $G_{i,3}$ to all $i=1,\dots,n, \infty$, 
that is,  we write 
$$
G_{i,3} = N(Q_3 \tin G_{i,2}(\chi_{3,2}))
$$
 for all such $i$.
Then 
\begin{equation*}
G_{2,3}= N(Q_3 \tin G_{2,2}(\chi_{3,2}))= N(Q_3 \tin G_{2,2}) 
=N(Q_3 \tin P_2 \times Q_{1,2}) = P_{2,3} \times Q_{1,2},
\end{equation*}
where the last equation  holds,  as according to \eqref{pq14b}, 
we have $Q_{1,2} = Q_3 \cap Q_1$.
Furthermore,  the character $\chi_{2,2}= \alpha_2 \times \beta_{1,2}$
corresponds to  the character $\chi_{2,3}:=
\alpha_{2,3} \times \beta_{1,2} \in  \Irr(N(Q_3 \tin G_{2,2}))$, through the $Q_3$-Glauberman
 correspondent $\alpha_2$ of  $\alpha_{2,3}$.

Also,
\begin{equation*}
 G_{1,3} = N(Q_3 \tin G_{1,2}) = N(Q_3 \tin Q_{1,2}) = Q_{1,2}=G_{1,2},
\end{equation*}
 and thus we take $\chi_{1,3}:= \chi_{1,2} = \beta_{1,2}$.

This, combined with the former $cQ_3$-correspondence,  provides
a correspondence (that we also  write as $cQ_3$-correspondence) 
between the character towers of those of the series \eqref{n2*} 
and the series
\begin{equation}\mylabel{n3*}
G_{0,3} =1\unlhd G_{1,3}= Q_{1,2} \unlhd 
G_{2,3}= P_{2,3} \times Q_{1,2} \unlhd \dots \unlhd G_{n,3}\unlhd \Gap{3}.
\end{equation}
 
We remark here that,  for every group $M$ with $M \leq N(P_2, Q_3  \tin G)$,  the 
$Q_3$-Glauberman correspondence  between $\Irr^{Q_3}(P_2)$ and $\Irr(P_{2,3})$
 is $M$-invariant. In particular we have $M(\alpha_2) = M(\alpha_{2,3})$, and thus
$$
M(\chi_{1,2}, \chi_{2,2}) = M(\chi_{1,3}, \chi_{2,3}).
$$
Therefore, in view of \eqref{pqn1-2},  we get 
\begin{equation}\mylabel{pqn1''}
M(\chi_1,\chi_2 ) = M(\chi_{1,2},\chi_{2,2})= M(\chi_{1,3}, \chi_{2,3}),
\end{equation}
So,  for any $M$ with $M \leq N(P_2 , Q_3 \tin G)$,  we have 
\begin{align*}
M(\chi_1, \chi_2, \dots,\chi_n)&= M(\chi_{1,1}, \chi_{2,2}, \dots,\chi_{n,2}) & &\text{ by 
\eqref{pqn1-2}}\\
&=M(\chi_{1,1}, \chi_{2,2})(\chi_{3,2},\dots,\chi_{n,2}) & & \\
&=M(\chi_1,\chi_2)(\chi_{3,2},\dots,\chi_{n,2}) & &\text{ by \eqref{pqn1''}}\\
&=M(\chi_1,\chi_2)(\chi_{3,3},\dots,\chi_{n,3}) & &\text{ by \eqref{pqn3}, since 
$M(\chi_1, \chi_2) \leq N(Q_3 \tin \Gap{2})$. }\\
&=M(\chi_{1,3},\chi_{2,3})(\chi_{3,3},\dots,\chi_{n,3}) & &\text{ by \eqref{pqn1''}}\\
&=M(\chi_{1,3}, \chi_{2,3},\dots,\chi_{n,3}). & &
\end{align*}
The following table gives a clear  picture of the situation when $m=3$:
\begin{footnotesize}
\begin{diagram}[small]
G=G_{\infty} & & & &G_{\infty, 1}=G(\beta_1) & & & \\
\vLine   &\qquad   &     &   &\vLine    &\qquad  &    &                \\ 
  G_{n}    &\qquad   &\chi_{n}  &   &G_{n,1}=G_n(\beta_1)
 &\qquad  &\chi_{n,1} & \\
     \vLine   &\qquad   &\vLine     &   &\vLine    &\qquad  &\vLine    &                \\
     \vdots  &\qquad   &\vdots     &   &\vdots    &\qquad  &\vdots   &                \\
\vLine   &\qquad   &\vLine     &   &\vLine    &\qquad  &\vLine    &                \\
  G_{3}    &\qquad   &\chi_{3}  &    &G_{3,1}=G_3(\beta_1)  
 &\qquad     &\chi_{3,1} &    \\
     \vLine    &\qquad  &\vLine   &\text{\qquad  $\overleftrightarrow{cQ_1}$ \qquad} &\vLine 
       &\qquad  &\vLine   &\qquad  & \\
   G_{2} &\qquad   &\chi_{2}  &   
   &G_{2,1}= P_2 \ltimes Q_{1}   &\qquad  &\chi_{2,1}=\alpha_2 \cdot \beta^e_{1}  & \\
    \vLine    &\qquad  &\vLine     & &\vLine    &\qquad  &\vLine   & \\
   G_1=Q_1  &\qquad   &\chi_1=\beta_1 &  &G_{1,1}= Q_{1} &\qquad  &\chi_{1,1}=\beta_{1} &\\
 \vLine    &\qquad  &\vLine     & &\vLine    &\qquad  &\vLine   & \\
     G_{0}=1   &\qquad   &\chi_{0}=1 & &G_{0,1} = 1 &\qquad &\chi_{0,1} = 1 &
  \end{diagram}
\end{footnotesize}
\begin{table}[h]
\begin{footnotesize}
\begin{diagram}[small]
 & &\Gap{2}=G(\beta_1,\alpha_2) & & & &G_{\infty, 3}=G(\beta_1,\alpha_2,\beta_3) & & & \\
& &\vLine   &\qquad   &     &   &\vLine    &\qquad  &    &        \\ 
  & &G_{n,2}=G(\beta_1,\alpha_2)    &\qquad   &\chi_{n,2}  &  
       &G_{n,3}=G_n(\beta_1,\alpha_2,\beta_3)
 &\qquad  &\chi_{n,3} & \\
     & &\vLine   &\qquad   &\vLine     &   &\vLine    &\qquad  &\vLine    &         \\
    &  &\vdots  &\qquad   &\vdots     &   &\vdots    &\qquad  &\vdots   &          \\
& &\vLine   &\qquad   &\vLine     &   &\vLine    &\qquad  &\vLine    &                \\
  & &G_{3,2}=P_2\rtimes Q_3   &\qquad   &\chi_{3,2}= \alpha_2^e \cdot \beta_3 &     &G_{3,3}=P_{2,3} \times Q_3  &\qquad     &\chi_{3,3}=\alpha_{2,3} \times \beta_3 &    \\
 \text{\qquad  $\overleftrightarrow{cP_2}$ \qquad}  &\qquad &\vLine    &\qquad  &\vLine   &\text{\qquad  $\overleftrightarrow{cQ_3}$ \qquad} &\vLine &\qquad  &\vLine   &\qquad  & \\
  & &G_{2,2}=P_2 \times Q_{1,2} &\qquad   &\chi_{2,2}=\alpha_2 \times \beta_{1,2}  &   
 &G_{2,3}= P_{2,3} \times Q_{1,2}   &\qquad  &\chi_{2,3}=\alpha_{2,3}\times \beta_{1,2} & \\
    & &\vLine    &\qquad  &\vLine     & &\vLine    &\qquad  &\vLine   & \\
   & &G_{1,2}=Q_{1,2}  &\qquad   &\beta_{1,2} &  &G_{1,3}= Q_{1,2} &\qquad  &\beta_{1,2} &\\
 & &\vLine    &\qquad  &\vLine     & &\vLine    &\qquad  &\vLine   & \\
    &  &G_{0,2}=1   &\qquad   &\chi_{0,2}=1 & &G_{0,3} = 1 &\qquad &\chi_{0,3} = 1 &
  \end{diagram}
\caption{The $cQ_1,cQ_3, cP_2$-correspondence.}
\mylabel{diagr.1-2-3}
\end{footnotesize}
\end{table}

\newpage
We stop here with the individual cases $m=1,2,3$, hoping that it has become clear 
how the mechanism  that produces  triangular sets from character towers  works.
We only remark that the role of the $\pi$- and $\pi'$-groups 
is interchanged at every step. So the $\pi$-groups play the protagonistic role when $m$ is even,  and the $\pi'$-groups when $m$ is odd.
 This role  consists 
of two acts:

\smallskip
1) to pick the group and its character (here Lemma \ref{lemma0} is used),  and 

\smallskip 
2) to create the new $cP$- or $cQ$-correspondence (for this we use  Lemma
 \ref{dade'}).

We are ready to state and prove the inductive step of the above mechanism. 
\begin{theorem}\mylabel{tow--tri}
Assume that Hypothesis \ref{hyp1} holds.
Then every character tower $\{\ \chi_i \}_{i=0}^m$
 of \eqref{*} determines   a $G_m(\chi_0, \chi_1, \dots,\chi_m)$-conjugacy class
of triangular sets
\begin{equation}\mylabel{telo}
\{P_0,\dots,P_{2k}, Q_1,\dots,Q_{2l-1}
|\alpha_0,\dots,\alpha_{2k},\beta_1,\dots,\beta_{2l-1}\}
\end{equation}

for  \eqref{*}, where $k=[m/2]$ and $l=[(m+1)/2]$, 
such that 

\smallskip
1)Any subtower  $\{1=\chi_0, \chi_1, \dots,\chi_s \}$  
of the original character tower, for some
$s=1,\dots,m$, determines
a  $G_s(\chi_0, \chi_1, \dots, \chi_s)$-conjugacy class of triangular sets 
$$
\{P_0,\dots,P_{2[s/2]}, Q_1,\dots,Q_{2[(s+1)/2]-1}
|\alpha_0,\dots,\alpha_{2[s/2]},\beta_1,\dots,\beta_{2[(s+1)/2]-1}\},
$$
that are subsets of \eqref{telo}. 

\smallskip
2)For any  
\begin{equation}\mylabel{*1'}
1=G_0\unlhd G_1 \unlhd \dots \unlhd G_m \unlhd G_{m+1}\unlhd \dots \unlhd G_n
\unlhd G
\end{equation} 
extension of \eqref{*} to  a 
 normal series of $G$  that satisfies 
Hypothesis \ref{hyp1},  and  any
extension of $\{ \chi_i \}_{i=0}^m$ to  a 
 character tower 
$\{\chi_i \}_{i=0}^n$
 of this series,  there is a unique $cP_2,\dots,cP_{2k},cQ_1,\dots,
cQ_{2l-1}$-correspondent  
character tower $\{ 1=\chi_{0,m}, \chi_{1,m}, \dots, \chi_{m,m}, \dots ,
\chi_{n,m} \}$ of  the normal series 
$1=G_{0,m} \unlhd G_{1,m} \unlhd \dots \unlhd  G_{m,m}\unlhd 
\dots \unlhd G_{n,m}\unlhd \Gap{m}$, for all $n$ with 
 $1\leq m \leq n$.
Here
\begin{multline}\mylabel{tow--tri3}
G_{i,m}= G_i(\alpha_2,\dots, \alpha_{2k},
 \beta_1,\dots,\beta_{2l-1})\\
= N(P_0, \dots,P_{2k}, Q_1, \dots,Q_{2l-1} \tin 
G_i(\chi_1, \dots, \chi_{m})),
\end{multline}
where      $i=0,1,\dots,n,\infty$.

\smallskip
3)For every $M$ with $M \leq 
N(P_0, \dots,P_{2k}, Q_1, \dots,Q_{2l-1} \tin G)$
we have 
$$
M(\chi_1,\dots,\chi_{n}) = M(\chi_{1,s},\dots,\chi_{n,s}).
$$

\smallskip
4) For every $i=1,2,\dots,m$ and $s = i+1,\dots,m$, these 
groups and characters
 follow the rules  
\begin{equation}\mylabel{tow--tri1}
\begin{aligned}
G_{i,i-1}&= P_{i} \ltimes Q_{i-1}&\quad  &\text{and} &\quad
\chi_{i,i-1}&= \alpha_{i} \cdot \beta_{i-1}^e,  \\
 G_{i,i}&=P_i \times Q_{i-1,i} &\quad  &\text{and} &\quad 
\chi_{i,i}&=\alpha_i \times \beta_{i-1,i}, \\
 G_{i,s}&= P_{i, 2[(s+1)/2]-1} \times Q_{i-1,2[s/2]} &\quad  &\text{and}  &\quad 
\chi_{i,s} &= \alpha_{i, 2[(s+1)/2]-1} \times \beta_{i-1,2[s/2]}, 
\end{aligned}
\end{equation}
whenever $i$ is even, and  
\begin{equation}\mylabel{tow--tri2}
\begin{aligned}
G_{i,i-1}&= P_{i-1} \rtimes Q_{i} &\quad &\text{and} &\quad 
\chi_{i,i-1}&= \alpha_{i-1}^e \cdot \beta_{i}, \\
 G_{i,i}&=P_{i-1,i} \times Q_{i} &\quad &\text{and} &\quad 
\chi_{i,i}&=\alpha_{i-1,i} \times \beta_{i},  \\
G_{i,s} &= P_{i-1, 2[(s+1)/2]-1} \times Q_{i,2[s/2]} &\quad  &\text{and} &\quad 
\chi_{i,s}&= \alpha_{i-1, 2[(s+1)/2]-1} \times \beta_{i,2[s/2]}, 
\end{aligned}
\end{equation}
when $i$ is odd.
(Here $\beta_{i-1}^e$ is the canonical extension of $\beta_{i-1} \in \Irr(Q_{i-1})$ 
to $G_{i,i-1}$  and similarly, $\alpha_{i-1}^e$ is the canonical  extension of 
$\alpha_{i-1} \in 
\Irr(P_{i-1})$ to $G_{i,i-1}$.) 
\end{theorem}

\begin{proof}
We will use induction on $m$. We have already seen that the theorem holds 
when $m=1$ (also when $m=2$ and $m=3$).

So  assume that the theorem holds for all $m=1,\dots ,t$ and some integer
$t \geq 0$. 
We will prove it also holds  when $m=t+1$. So assume that 
the normal series 
\begin{equation}\mylabel{**}
1=G_0 \unlhd G_1 \unlhd \dots \unlhd G_t \unlhd G_{t+1} \unlhd G
\end{equation}
is fixed. Along with that we fix a character tower 
\begin{equation}\mylabel{**ch}
\{1=\chi_0, \chi_1,\dots,\chi_t,\chi_{t+1}\}
\end{equation}
 for  \eqref{**}. 
 As the triangular sets have different 
form depending on whether $t$ is even or odd,  we split 
the proof in two symmetric 
cases.
\begin{description}
\item[Case 1: $t$ is odd]

The series $1=G_0 \unlhd G_1 \unlhd \dots \unlhd G_t  \unlhd G$ 
is also a normal 
series of $G$,  while the irreducible characters  $\{ \chi_i\}_{i=0}^t$
 form a character tower for this series.
Hence the inductive hypothesis implies the existence of a 
$G_{t}(\chi_1, \dots, \chi_t)$-conjugacy
 class of 
triangular sets that is determined by  the last  character tower.
 Let  $\{P_0,\dots,P_{2k},Q_1,\dots,Q_{2l-1}|
\alpha_0,\dots,\alpha_{2k},\beta_1,\dots,\beta_{2l-1}\}$ 
be a representative of 
this conjugacy class.
As $t$ is odd we have that $l=[(t+1)/2] = (t+1) /2$ while 
$k = [t/2] = (t-1)/2$. So $2l-1 = t$ while $2k = t-1$. Therefore the above 
triangular set has the form $\{P_0,\dots,P_{t-1},Q_1,\dots,Q_{t}|
\alpha_0,\dots,\alpha_{t-1},\beta_1,\dots,\beta_{t}\}$. 

Hence,  to prove that the character tower \eqref{**ch} determines  
 a $G_{t+1}( \chi_1, \dots, \chi_t, \chi_{t+1})$-conjugacy class
of triangular sets that respect subtowers, 
it is enough  to prove 
the existence of a $\pi$-group $P_{t+1}$ unique up to conjugations
 by any element of $G_{t+1}(\chi_1, \dots, \chi_t, \chi_{t+1})$, 
  and an irreducible character 
$\alpha_{t+1} \in \Irr(P_{t+1})$ such that 
the set 
$$\{P_0,\dots,P_{t-1},P_{t+1},Q_1,\dots,Q_{t}|
\alpha_0,\dots,\alpha_{t-1},\alpha_{t+1},\beta_1,\dots,\beta_{t}\}
$$
is a  triangular set depending on the tower \eqref{**ch}. 

Let 
\begin{equation}\mylabel{*1}
1=G_0 \unlhd G_1 \unlhd \dots \unlhd G_t \unlhd G_{t+1}
\unlhd \dots \unlhd  G_n \unlhd G,
\end{equation}
be an extension of \eqref{**} to a  normal series of $G$ so that Hypothesis \ref{hyp1} holds,  for some $n \geq t+1$. 
Assume further that 
\begin{equation}\mylabel{*2}
\{ 1=\chi_0, \chi_1, \dots,\chi_t, \chi_{t+1}, \dots, \chi_n \},
\end{equation} 
is a character tower for \eqref{*1} that extends the character tower \eqref{**ch}.
For any $n$ with $n \geq t$,   our inductive hypothesis implies
  that the 
character tower \eqref{*2} of  \eqref{*1}
has a $cQ_1,cP_2,\dots,cP_{t-1},cQ_t$-correspondent   character tower
\begin{equation}\mylabel{*3}
\{ 1=\chi_{0,t}, \chi_{1,t}, \dots, \chi_{t,t}, \dots ,
\chi_{n,t} \},
\end{equation}
 of   the normal series 
$1=G_{0,t} \unlhd G_{1,t} \unlhd \dots \unlhd  G_{t,t}\unlhd 
\dots \unlhd G_{n,t}\unlhd \Gap{t}$,
 where 
\begin{multline}\mylabel{t-t1}
G_{i,t}= G_i(\alpha_2,\dots, \alpha_{t-1},
 \beta_1,\dots,\beta_{t})\\
= N(P_0, \dots,P_{t-1}, Q_1, \dots,Q_{t} \tin 
G_i(\chi_1, \dots, \chi_{t})).
\end{multline}
For every $M \leq N(P_2,\dots,P_{t-1},Q_1,\dots,Q_{t} \tin G)$ we also  have 
\begin{equation}\mylabel{t-t34}
M(\chi_1,\dots,\chi_n)=M(\chi_{1,t},\dots,\chi_{n,t}).
\end{equation}
Furthermore, \eqref{tow--tri1} and \eqref{tow--tri2}  for $i=t$ imply that 
$G_{t,t}= P_{t-1,t} \times Q_{t}$,  while
$\chi_{t,t} = \alpha_{t-1,t} \times \beta_t$. (Note that $\alpha_{t-1,t} \in 
\Irr(P_{t-1,t})$ is the $Q_t$-Glauberman correspondent of 
$\alpha_{t-1} \in \Irr^{Q_t}(P_{t-1})$).
The following diagram describes the situation. 

\newpage
\begin{footnotesize}
\begin{diagram}[small]
G=G_{\infty} & & & &G_{\infty, 1}=G(\beta_1) & & & \\
\vLine   &\qquad   &     &   &\vLine    &\qquad  &    &                \\ 
  G_{n}    &\qquad   &\chi_{n}  &   &G_{n,1}=G_n(\beta_1)
 &\qquad  &\chi_{n,1} & \\
     \vLine   &\qquad   &\vLine     &   &\vLine    &\qquad  &\vLine    &    \\
     \vdots  &\qquad   &\vdots     &   &\vdots    &\qquad  &\vdots   &                \\
\vLine   &\qquad   &\vLine     &   &\vLine    &\qquad  &\vLine    &                \\
 G_{t+1}    &\qquad   &\chi_{t+1}  &   &G_{t+1,1}=G_{t+1}(\beta_1)
 &\qquad  &\chi_{t+1,1} & \\
  \vLine   &\qquad   &\vLine     &   &\vLine    &\qquad  &\vLine    & \\
      G_{t}    &\qquad   &\chi_{t}  &   &G_{t,1}=G_{t}(\beta_1)
 &\qquad  &\chi_{t,1} & \\ 
 \vLine   &\qquad   &\vLine     &   &\vLine    &\qquad  &\vLine    & \\
      G_{t-1}    &\qquad   &\chi_{t-1}  &   &G_{t-1,1}=G_{t-1}(\beta_1)
 &\qquad  &\chi_{t-1,1} & \\ 
    \vLine   &\qquad   &\vLine     &   &\vLine    &\qquad  &\vLine    & \\
\vdots  &\qquad   &\vdots     &   &\vdots    &\qquad  &\vdots   &                \\
     \vLine    &\qquad  &\vLine   &\text{\qquad  $\overleftrightarrow{cQ_1}$ \qquad} &\vLine 
       &\qquad  &\vLine   &\qquad  & \\
     G_{3}    &\qquad   &\chi_{3}  &   &G_{3,1}=G_{3}(\beta_1)
 &\qquad  &\chi_{3,1} & \\
  \vLine   &\qquad   &\vLine     &   &\vLine    &\qquad  &\vLine    & \\
  G_{2} &\qquad   &\chi_{2}  &   
   &G_{2,1}= P_2 \ltimes Q_{1}   &\qquad  &\chi_{2,1}=\alpha_2 \cdot \beta^e_{1}  & \\
    \vLine    &\qquad  &\vLine     & &\vLine    &\qquad  &\vLine   & \\
   G_1=Q_1  &\qquad   &\chi_1=\beta_1 &  &G_{1,1}= Q_{1} &\qquad  &\chi_{1,1}=\beta_{1} &\\
 \vLine    &\qquad  &\vLine     & &\vLine    &\qquad  &\vLine   & \\
     G_{0}=1   &\qquad   &\chi_{0}=1 & &G_{0,1} = 1 &\qquad &\chi_{0,1} = 1 &
  \end{diagram}
\end{footnotesize}

\newpage
\begin{footnotesize}
\begin{diagram}[small]
 & &\Gap{2}=G(\beta_1,\alpha_2) & & & &G_{\infty, 3}=G(\beta_1,\alpha_2,\beta_3) & & & \\
& &\vLine   &\qquad   &     &   &\vLine    &\qquad  &    &        \\ 
  & &G_{n,2}=G(\beta_1,\alpha_2)    &\qquad   &\chi_{n,2}  &  
       &G_{n,3}=G_n(\beta_1,\alpha_2,\beta_3)
 &\qquad  &\chi_{n,3} & \\
     & &\vLine   &\qquad   &\vLine     &   &\vLine    &\qquad  &\vLine    &         \\
    &  &\vdots  &\qquad   &\vdots     &   &\vdots    &\qquad  &\vdots   &          \\
& &\vLine   &\qquad   &\vLine     &   &\vLine    &\qquad  &\vLine    &                \\
& &G_{t+1,2}=G(\beta_1,\alpha_2)    &\qquad   &\chi_{t+1,2}  &  
       &G_{t+1,3}=G_{t+1}(\beta_1,\alpha_2,\beta_3)
 &\qquad  &\chi_{t+1,3} & \\
  & &\vLine   &\qquad   &\vLine     &   &\vLine    &\qquad  &\vLine    &         \\
& &G_{t,2}=G(\beta_1,\alpha_2)    &\qquad   &\chi_{t,2}  &  
       &G_{t,3}=G_{t}(\beta_1,\alpha_2,\beta_3)
 &\qquad  &\chi_{t,3} & \\
 & &\vLine   &\qquad   &\vLine     &   &\vLine    &\qquad  &\vLine    &         \\
& &G_{t-1,2}=G(\beta_1,\alpha_2)    &\qquad   &\chi_{t-1,2}  &  
       &G_{t-1,3}=G_{t-1}(\beta_1,\alpha_2,\beta_3)
 &\qquad  &\chi_{t-1,3} & \\
 & &\vLine   &\qquad   &\vLine     &   &\vLine    &\qquad  &\vLine    &         \\
    &  &\vdots  &\qquad   &\vdots     &   &\vdots    &\qquad  &\vdots   &          \\
 \text{\qquad  $\overleftrightarrow{cP_2}$ \qquad}  &\qquad &\vLine    &\qquad  &\vLine   &\text{\qquad  $\overleftrightarrow{cQ_3}$ \qquad} &\vLine &\qquad  &\vLine   &\qquad  & \\
  & &G_{3,2}=P_2 \rtimes Q_{3} &\qquad   &\chi_{3,2}=\alpha_2^e \cdot 
\beta_{3}  &   
 &G_{3,3}= P_{2,3} \times Q_{3}   &\qquad  &\chi_{3,3}=\alpha_{2,3}\times 
\beta_{3} & \\
& &\vLine   &\qquad   &\vLine     &   &\vLine    &\qquad  &\vLine    &         \\
 & &G_{2,2}=P_2 \times Q_{1,2} &\qquad   &\chi_{2,2}=\alpha_2 \times \beta_{1,2}  &   
 &G_{2,3}= P_{2,3} \times Q_{1,2}   &\qquad  &\chi_{2,3}=\alpha_{2,3}\times \beta_{1,2} & \\
    & &\vLine    &\qquad  &\vLine     & &\vLine    &\qquad  &\vLine   & \\
   & &G_{1,2}=Q_{1,2}  &\qquad   &\beta_{1,2} &  &G_{1,3}= Q_{1,2} &\qquad  &\beta_{1,2} &\\
 & &\vLine    &\qquad  &\vLine     & &\vLine    &\qquad  &\vLine   & \\
    &  &G_{0,2}=1   &\qquad   &\chi_{0,2}=1 & &G_{0,3} = 1 &\qquad &\chi_{0,3} = 1 &
  \end{diagram}
\end{footnotesize}

\newpage
\begin{footnotesize}
\begin{diagram}[small]
 &  & &\Gap{t-1}=G(\alpha_2,\dots,\alpha_{t-1},\beta_1,\dots,
\beta_{t-2}) &\qquad  & \\
& & &\vLine & & \\
& & &G_{n,t-1}=G_n(\alpha_2,\dots,\alpha_{t-1},\beta_1,\dots, 
\beta_{t-2}) & &\chi_{n,t-1}\\
& & &\vLine & &\vLine \\
& & &\vdots & &\vdots \\
& & &\vLine & &\vLine \\
& & &G_{t+1,t-1}=G_{t+1}(\alpha_2,\dots,\alpha_{t-1},\beta_1,\dots,
\beta_{t-2}) & &\chi_{t+1,t-1} \\
& & &\vLine & &\vLine \\
& & &G_{t,t-1}=P_{t-1} \rtimes Q_t & &\chi_{t,t-1}=\alpha_{t-1}^e
 \cdot \beta_t \\
& & &\vLine & &\vLine \\
& & &G_{t-1,t-1}=P_{t-1} \times Q_{t-2,t-1} & &\chi_{t-1,t-1}=\alpha_{t-1}
 \times  \beta_{t-2,t-1} \\
& & &\vLine & &\vLine \\
\overleftrightarrow{cP_4,cQ_5,cP_6, \dots, cP_{t-1}}
&\qquad & &\vdots  &\qquad &\vdots \\
& & &\vLine & &\vLine \\
& & &G_{3,t-1}=P_{2,t-2}\times Q_{3,t-1} & &\chi_{3,t-1}=\alpha_{2,t-2}
\times \beta_{3,t-1}\\
&  & &\vLine & &\vLine \\
& & &G_{2,t-1}=P_{2,t-2}\times Q_{1,t-1} & &\chi_{2,t-1}=\alpha_{2,t-2}
\times \beta_{1,t-1}\\
&  & &\vLine & &\vLine \\
& & &G_{1,t-1}=Q_{1,t-1} & &\chi_{1,t-1}= \beta_{1,t-1}\\
&  & &\vLine & &\vLine \\
& & &G_{0,t-1}=1 & &\chi_{0,t-1}= 1\\
  \end{diagram}
\end{footnotesize}

\newpage
\begin{table}[h]
\begin{footnotesize}
\begin{diagram}[small]
 &  & &\Gap{t}=G(\alpha_2,\dots,\alpha_{t-1},\beta_1,\dots,
\beta_{t-2},\beta_{t}) &\qquad  & \\
& & &\vLine & & \\
& & &G_{n,t}=G_n(\alpha_2,\dots,\alpha_{t-1},\beta_1,\dots, 
\beta_{t-2},\beta_{t}) & &\chi_{n,t}\\
& & &\vLine & &\vLine \\
& & &\vdots & &\vdots \\
& & &\vLine & &\vLine \\
& & &G_{t+1,t}=G_{t+1}(\alpha_2,\dots,\alpha_{t-1},\beta_1,\dots, 
\beta_{t-2},\beta_{t}) & &\chi_{t+1,t} \\
& & &\vLine & &\vLine \\
& & &G_{t,t}=P_{t-1,t} \times Q_t & &\chi_{t,t}=\alpha_{t-1,t}
 \times \beta_t \\
& & &\vLine & &\vLine \\
& & &G_{t-1,t}=P_{t-1,t} \times Q_{t-2,t-1} & &\chi_{t-1,t}=\alpha_{t-1,t}
 \times  \beta_{t-2,t-1} \\
& & &\vLine & &\vLine \\
\overleftrightarrow{cQ_{t}}
&\qquad & &\vdots  &\qquad &\vdots \\
& & &\vLine & &\vLine \\
& & &G_{3,t}=P_{2,t}\times Q_{3,t-1} & &\chi_{3,t}=\alpha_{2,t}
\times \beta_{3,t-1}\\
&  & &\vLine & &\vLine \\
& & &G_{2,t}=P_{2,t}\times Q_{1,t-1} & &\chi_{2,t}=\alpha_{2,t}
\times \beta_{1,t-1}\\
&  & &\vLine & &\vLine \\
& & &G_{1,t}=Q_{1,t-1} & &\chi_{1,t}= \beta_{1,t-1}\\
&  & &\vLine & &\vLine \\
& & &G_{0,t}=1 & &\chi_{0,t}= 1\\
  \end{diagram}
\caption{The $cQ_1,cQ_3,\dots,cQ_t,cP_2, cP_4, \dots, cP_{t-1}$-correspondence}
\mylabel{diagr.1-2-3...}
\end{footnotesize}
\end{table}

\newpage
We pick  $P_{t+1}$ to be any $\pi$-Hall subgroup of 
$G_{t+1,t}= G_{t+1}(\alpha_2,\dots, \alpha_{t-1},
 \beta_1,\dots,\beta_{t})$.
The factor group  $G_{t+1,t}/G_{t,t}$ is a $\pi$-group, 
while  $Q_t$ is  a $\pi'$-Hall 
subgroup of $G_{t,t}$ normalized by $G_{t+1,t}$ (see \eqref{t-t1}).
Thus $P_{t+1}$ also normalizes $Q_{t}$, which implies that 
\begin{equation}\mylabel{t-t2}
G_{t+1,t}= P_{t+1} \ltimes Q_t.
\end{equation}
Furthermore, $Q_t$ is a normal subgroup of 
$\Gap{t}$ (as $G_{t,t} \unlhd \Gap{t}$), 
  while  its irreducible character  $\beta_t$ is $\Gap{t}$-invariant.
Therefore we can apply Lemma \ref{lemma0} 
to the groups $\Gap{t}, G_{t+1,t}, Q_t$ and the character $\chi_{t+1,t}$
 in the place of the groups $G, H,N$ and the character $\theta$, respectively.
(Note that in this case $H(\theta)= H$).
Hence we conclude that $\beta_t$ has a unique canonical extension
 $\beta_t^e \in \Irr(G_{t+1,t})$. Furthermore,
since $\chi_{t+1,t}$ lies above $\beta_t$, the same lemma implies the
 existence of  a unique irreducible character
 $\alpha_{t+1} \in \Irr(P_{t+1})$ such that 
\begin{equation}\mylabel{t-t3}
\chi_{t+1,t} = \alpha_{t+1} \cdot \beta_{t}^e,
\end{equation}
while $\Gap{t}(\beta_t,\alpha_{t+1}) = N(P_{t+1} \tin G_{n,t}(\chi_{t+1,t}))$.
But  $\Gap{t}$ fixes $\beta_t$. So
\begin{equation}\mylabel{t-t4}
\Gap{t}(\alpha_{t+1}) = N(P_{t+1} \tin G_{n,t}(\chi_{t+1,t})).
\end{equation}
As $\chi_{t+1,t}$ lies above $\chi_{t,t}=\alpha_{t-1,t}\times \beta_t$, 
 equation \eqref{t-t3} obviously implies that $\alpha_{t+1}$ lies above 
$\alpha_{t-1,t}$, which is  the $Q_{t}$-Glauberman correspondent of
 $\alpha_{t-1}$.
This,  along with the fact that $P_{t+1}$ was picked as a $\pi$-Hall 
subgroup of $G_{t+1}(\alpha_2,\dots, \alpha_{t-1},
 \beta_1,\dots,\beta_{t})$, implies that the new $\pi$-group and its
 character satisfy \eqref{x3} and \eqref{x4} respectively.
As we already know that the set  $\{P_0,\dots,P_{t-1}, Q_1,\dots,Q_{t}|
\alpha_0,\dots,\alpha_{t-1},\beta_1,\dots,\beta_{t}\}$
is a triangular set,  we conclude that 
\begin{equation}\mylabel{t-t5}
\{P_0,\dots,P_{t-1},P_{t+1},Q_1,\dots,Q_{t}|
\alpha_0,\dots,\alpha_{t-1},\alpha_{t+1},\beta_1,\dots,\beta_{t}\}
\end{equation}
is a triangular set for \eqref{**}. Furthermore, it is
clear, from the way it is constructed,  that it is related to  
 the character tower \eqref{**ch} and that it respects subtowers.
Note also that the only choice for  $P_{t+1}$ was that of the Hall
$\pi$-subgroup of $G_{t+1, t}$.  Hence $P_{t+1}$ is uniquely determined 
up to conjugation by an element of 
$$
G_{t+1,t}=  N(P_0, P_2, \dots, P_{t-1}, Q_1, \dots, Q_{t}
\tin G_{t+1}(\chi_1, \dots, \chi_{t})).
$$
So $P_{t+1}$ is uniquely determined by an element of 
$G_{t+1}(\chi_1, \dots, \chi_t)$. 
This, along with the inductive hypothesis and the fact that 
$G_{t+1}(\chi_1, \dots, \chi_{t}) \geq G_t(\chi_1, \dots, \chi_{t-1})$
implies that the triangular set \eqref{t-t5} is unquely determined
 up to conjugation by an element of $G_{t+1}(\chi_1, \dots, \chi_t)$.
Hence the first part of Theorem \ref{tow--tri} is verified for the 
inductive step  in the case where  $t$ is odd.
Furthermore, as \eqref{t-t5} is a triangular set,  
Proposition \ref{prodd} implies that 
\begin{equation}\mylabel{t-t6}
\begin{aligned}
G_{t+1}(\alpha_2,\dots,\alpha_{t+1},\beta_1,\dots,\beta_{t})&= P_{t+1} \times Q_{t,t+1} 
&\quad &\text{ by \eqref{pq-15a}, }\\
G_{j}(\alpha_2,\dots,\alpha_{t+1},\beta_1,\dots,\beta_{t})&= 
P_{i,t} \times Q_{i-1,t+1} 
&\quad &\text{ if $i$ is even, by \eqref{prod2}, }\\
G_{j}(\alpha_2,\dots,\alpha_{t+1},\beta_1,\dots,\beta_{t})&= P_{i-1,t} \times Q_{i,t+1}
 &\quad &\text{ if $i$ is odd, by \eqref{prod4}},
\end{aligned}
\end{equation}
for all $j=1,\dots,t$.

To complete the proof of the theorem (at least when $t$ is odd),  it is enough 
to show that the character tower \eqref{*2} determines a character tower 
$\{1=\chi_{0,t+1}, \dots, \chi_{t+1, t+1}, \dots ,\chi_{n,t+1} \}$
for the series $1=G_{0,t+1} \unlhd G_{1,t+1} \unlhd  \dots \unlhd G_{t+1,t+1}
\unlhd  \dots \unlhd G_{n,t+1}\unlhd \Gap{t+1}$, 
where $G_{i,t+1}$ and $\chi_{i,t+1}$ satisfy 
\eqref{tow--tri3}, \eqref{tow--tri1} and \eqref{tow--tri2}.
In that direction we first observe that,  for every $i=1,\dots,n,\infty$,  we have 
\begin{equation}\mylabel{t-t7}
G_{i,t}(\chi_{t+1,t}) = G_{i,t}(\chi_{t+1}).
\end{equation}
Indeed, in view of  \eqref{t-t34} for $M= G_{i,t}$ and $n=t$, 
we get 
$G_{i,t}(\chi_1,\dots,\chi_{t})=G_{i,t}(\chi_{1,t},\dots,\chi_{t,t})$.
Hence
$G_{i,t}=G_{i,t}(\chi_1,\dots,\chi_{t})=G_{i,t}(\chi_{1,t},\dots,\chi_{t,t})$, 
as $G_{i,t} \leq G_i(\chi_1,\dots,\chi_{t})$ by \eqref{t-t1}.
So,  if we apply again \eqref{t-t34} for $M=G_{i,t}$ and
$n=t+1$,  we have
\begin{multline*}
G_{i,t}(\chi_{t+1})= G_{i,t}(\chi_1,\dots,\chi_{t},\chi_{t+1})=
G_{i,t}(\chi_{1,t},\dots,\chi_{t,t},\chi_{t+1,t})\\
=G_{i,t}(\chi_{1,t},\dots,\chi_{t,t})(\chi_{t+1,t})=G_{i,t}(\chi_{t+1,t}).
\end{multline*}
 Thus \eqref{t-t7} holds.
Hence we conclude that 
$$
N(P_{t+1} \tin G_{i,t}(\chi_{t+1,t})) = N(P_{t+1} \tin G_{i,t}(\chi_{t+1})).
$$
This, along with  \eqref{t-t4}, implies that 
\begin{equation}\mylabel{t-t8}
G_{i,t}(\alpha_{t+1}) = N(P_{t+1} \tin G_{i,t}(\chi_{t+1})),
\end{equation}
whenever $i=1,\dots,n,\infty$.
Thus if we define $G_{i,t+1} := G_{i,t}(\alpha_{t+1})$, equations \eqref{t-t1} and
 \eqref{t-t8}
imply 
\begin{multline}\mylabel{t-t9}
G_{i,t+1} = G_i(\alpha_2,\dots,\alpha_{t-1},\alpha_{t+1},\beta_1,\dots,\beta_t)\\
=G_{i, t}(\alpha_{t+1})= N(P_{t+1} \tin G_{i, t}(\chi_{t+1}))\\
=N(P_0,P_2,\dots,P_{t-1},P_{t+1},Q_1,\dots,Q_{t} \tin G_{i}(\chi_1,\dots,\chi_{t+1})).
\end{multline}
We also get (using \eqref{t-t9} and \eqref{t-t7}) 
\begin{equation}\mylabel{t-t10}
G_{i,t+1} = N(P_{t+1} \tin G_{i, t}(\chi_{t+1}))= 
N(P_{t+1} \tin G_{i,t}(\chi_{t+1,t})),
\end{equation}
for all $i=1,\dots,n,\infty$.
Note that \eqref{t-t9} proves that \eqref{tow--tri3} holds for the inductive step.
Also \eqref{t-t6}, along with \eqref{t-t2} and the inductive hypothesis
(for those $i$ with $i=1,\dots,t$), implies 
that the groups $G_{i,s}$ satisfy \eqref{tow--tri1} and \eqref{tow--tri2}
whenever $1\leq i \leq t+1$ and $i<s \leq t+1$.
In particular we have
\begin{equation}\mylabel{t-t14}
\begin{aligned}
G_{t+1,t+1} &= P_{t+1} \times Q_{t,t+1} & &\\
G_{j,t+1}&= P_{j,t} \times Q_{j-1,t+1} 
&\quad &\text{ if $j$ is even }\\
G_{j,t+1}&= P_{j-1,t} \times Q_{j,t+1}
 &\quad &\text{ if $j$ is odd},
\end{aligned}
\end{equation}
for all $j=1,\dots,t$.

To get the desired character tower for the series \eqref{*1} 
(the correspondent of 
the tower \eqref{*2}),   we first use the inductive argument  to reach the character tower 
\eqref{*3}. So it is enough to
get a  tower for \eqref{*1} that corresponds to  this  latter tower. For this we split 
 \eqref{*3} in two pieces: the tail that consists of $\chi_{i,t}$ for all $i=1,\dots,t$,
and the top  that consists of the rest, i.e., the characters $\chi_{i,t}$ where
$i=t+1,\dots,n$.

\smallskip
For the top part, we apply Lemma \ref{dade'} to the normal subgroups $G_{t+1,t},\dots,G_{n,t}$
of $\Gap{t}$ and the character $\chi_{t+1,t}=\alpha_{t+1} \cdot \beta_{t}^e$.
This way the  character tower
 $\{ \chi_{t+1,t}, \dots,\chi_{n,t}\}$ of  the normal  series 
$G_{t+1,t}\unlhd \dots \unlhd G_{n,t}$
 has a unique $cP_{t+1}$-correspondent
character tower of the   series
$$
G_{t+1,t+1} = N(P_{t+1} \tin G_{t+1,t}(\chi_{t+1,t})) \unlhd \dots \unlhd 
G_{n,t+1}= N(P_{t+1} \tin G_{n,t}(\chi_{t+1,t})).
$$
We write 
\begin{equation}\mylabel{t-t11}
\{ \chi_{t+1,t+1}, \dots , \chi_{n,t+1} \},
\end{equation}
for this $cP_{t+1}$-correspondent tower.
Note that Lemma \ref{dade'} also determines the character $\chi_{t+1,t+1}$ as 
\begin{equation}\mylabel{t-t12}
\chi_{t+1,t+1} = \alpha_{t+1} \times \beta_{t, t+1} \in 
\Irr(G_{t+1,t+1}),
\end{equation}
where $\beta_{t,t+1} \in \Irr(Q_{t,t+1})= \Irr(C(P_{t+1} \tin Q_{t}))$ is the 
$P_{t+1}$-Glauberman correspondent of $\beta_t \in \Irr^{P_{t+1}}(Q_{t})$. 
Furthermore, according to the same lemma we have that 
\begin{equation}\mylabel{t-t30}
S(\chi_{t+1,t},\dots,\chi_{n,t})= S(\chi_{t+1,t+1},\dots,\chi_{n,t+1})
\end{equation}
for any  subgroup $S$ of   $N(P_{t+1} \tin \Gap{t})$.
 
As far as the tail of \eqref{*3} is concerned,  we observe the following:
In view of \eqref{tow--tri1} and \eqref{tow--tri2} 
\begin{equation}\mylabel{t-t13}
\begin{aligned}
\chi_{t,t}&= \alpha_{t-1,t} \times \beta_{t} &\\
\chi_{j,t}&=\alpha_{j,t}\times \beta_{j-1,t-1} &\text{ if $j$ is even} \\
\chi_{j,t}&=\alpha_{j-1,t}\times \beta_{j,t-1} &\text{ if $j$ is odd },
\end{aligned}
\end{equation}
whenever $j=1,\dots,t-1$.
We define 
\begin{equation}\mylabel{t-t15}
\begin{aligned}
\chi_{t,t+1}&= \alpha_{t-1,t} \times \beta_{t,t+1} &\\
\chi_{j,t+1}&=\alpha_{j,t}\times \beta_{j-1,t+1} &\text{ if $j$ is even} \\
\chi_{j,t+1}&=\alpha_{j-1,t}\times \beta_{j,t+1} &\text{ if $j$ is odd },
\end{aligned}
\end{equation}
where $\beta_{t,t+1}, \beta_{j,t+1}$ and $\beta_{j-1,t+1}$ are the $P_{t+1}$-Glauberman 
correspondents of $\beta_t, \beta_{j,t-1}$ and $\beta_{j-1,t-1}$ respectively, for all $j=1,\dots,t-1$.
Note that all these characters are well defined characters of $Q_{t,t+1}, Q_{j,t+1}$ and 
$Q_{j-1,t+1}$, and  form a tower by Proposition \ref{pqpropo} (as 
\eqref{t-t5} is a triangular set).
Furthermore, \eqref{t-t14} implies that   
$\chi_{j,t+1}$ and $\chi_{t,t+1}$ are 
 characters of $G_{j,t+1}$ and $G_{t,t+1}$ respectively, for all $j=1,\dots,t-1$.
Thus $\{ 1=\chi_{0,t+1}, \chi_{1,t+1}, \dots \chi_{t,t+1} \}$
is a character tower of the normal series $G_{0,t+1}
\unlhd G_{1,t+1} \unlhd \dots \unlhd G_{t,t+1}$. 
Also  we pass from the $\chi_{j,t}$ to the $\chi_{j,t+1}$ through 
a $P_{t+1}$-Glauberman correspondence. Thus any subgroup of 
$G$ that normalizes the groups $G_{1,t},\dots,G_{t,t}$,  along with the
$P_{t+1}$, leaves this correspondence  invariant.
But any group $T$ with $T \leq N(P_2,\dots,P_{t-1},P_{t+1},Q_1,\dots,Q_{t}
\tin G)$ normalizes the former groups (as $G_{i,t}$ is a direct product
(see Table \ref{diagr.1-2-3...}) of groups that $T$ normalizes).
Hence for any such group $T$ and any $j=1,\dots,t$ we have 
\begin{equation}\mylabel{t-t31}
T(\chi_{j,t}) = T(\chi_{j,t+1}).
\end{equation} 
Furthermore,   $\chi_{t,t+1} = \alpha_{t-1,t} \times \beta_{t,t+1}$,  while
$\chi_{t+1,t+1} = \alpha_{t+1} \times \beta_{t,t+1}$ (by \eqref{t-t12}).
As $\alpha_{t+1}$ lies above $\alpha_{t-1,t}$ (by \eqref{x4}), 
we conlcude that $\chi_{t+1,t+1}$ lies above $\chi_{t,t+1}$.
Hence we have formed the  tower $\{ 1=\chi_{0,t+1}, \dots,\chi_{t,t+1},
 \chi_{t+1,t+1}, \dots,\chi_{n,t+1}\}$ 
of \eqref{*1}, that corresponds to the tower \eqref{*3}.
This, along with the inductive argument that provides the $P_2,\dots,P_{t-1},Q_1,\dots,
Q_t$-correspondence between \eqref{**ch} and \eqref{*3},
implies the desired correspondence  between \eqref{**ch} and the tower 
$\{ 1=\chi_{0,t+1}, \dots,\chi_{t,t+1}, \chi_{t+1,t+1},
\dots,\chi_{n,t+1}\}$.
Furthermore, \eqref{t-t3}, \eqref{t-t12} and \eqref{t-t15},  along with 
the inductive argument,  imply that 
\eqref{tow--tri1} and \eqref{tow--tri2} hold for all $i=1,\dots,t+1$ and
 $s=i+1,\dots,t+1$.

\smallskip
Also  for every $M$ with $M \leq N(P_2,\dots,P_{t-1},P_{t+1},Q_1, \dots,Q_t 
\tin G)$ we get that 
$$
M(\chi_1,\dots,\chi_t) \leq N(P_{t+1} \tin \Gap{t}),
$$ 
(see  \eqref{t-t1} for
a characterization of $\Gap{t}$).
Therefore, for all such $M$ we have
\begin{align*}
&M(\chi_1,\dots,\chi_n)& & \\
&=M(\chi_1,\dots,\chi_t)(\chi_1,\dots,\chi_n) & & \\
&=M(\chi_1,\dots,\chi_t)(\chi_{1,t},\dots,\chi_{n,t}) &  &\text{ by 
\eqref{t-t34}} \\
&=M(\chi_1,\dots,\chi_t)(\chi_{t+1,t},\dots,\chi_{n,t})(\chi_{1,t},\dots,
\chi_{t,t}) & & \\
&=M(\chi_1,\dots,\chi_t)(\chi_{t+1,t+1},\dots,\chi_{n,t+1})(\chi_{1,t},\dots,
\chi_{t,t}) & 
&\text{by \eqref{t-t30} for $S =M(\chi_1,\dots,\chi_t)$ }\\
&=M(\chi_{1,t},\dots,\chi_{t,t})(\chi_{t+1,t+1},\dots,\chi_{n,t+1})
(\chi_{1,t},\dots,\chi_{t,t}) & &\text{ by \eqref{t-t34} }\\
&=M(\chi_{1,t},\dots,\chi_{t,t})(\chi_{t+1,t+1},\dots,\chi_{n,t+1}) & & \\
&=M(\chi_{1,t+1},\dots,\chi_{t,t+1})(\chi_{t+1,t+1},\dots,\chi_{n,t+1}) & 
&\text{by \eqref{t-t31} for $T =M$ }\\
&=M(\chi_{1,t+1}, \dots,\chi_{n,t+1}).  & &
\end{align*}
This implies  that  part 3) of the theorem also holds for $m=t+1$.
Hence  the  inductive step for $m=t+1$ is verified in the case of an odd $t$.

\item[Case 2: $t$ is even]
The proof is similar to that of an odd $t$. So we will skip it.
 We only remark that we need to interchange the role of the 
$\pi$-groups with that of the $\pi'$-groups.
 So in this case for the inductive step we 
pick the $\pi'$-group $Q_{t+1}$ and its character $\beta_{t+1}$, 
as in the previous one we were picking the $\pi$-group  $P_{t+1} $ and 
its character $\alpha_{t+1}$.
We continue similarly, proving that the inductive step holds also in 
the case of an 
even $t$. 
\end{description}
This completes the inductive argument and thus proving Theorem \ref{tow--tri}.
\end{proof}

The following  remark is a  straightforward consequence of the
recursive proof of Theorem \ref{tow--tri}
\begin{remark}\mylabel{t-tremark}
Let   $\{ \chi_{i, m}\}_{i=0}^n$ 
 be  the unique 
$cP_2, \dots, cP_{2k},cQ_1, \dots, cQ_{2l-1}$-correspondent 
of the character tower \eqref{*2}. Then its subtower   $\{ \chi_{i, m}\}_{i=0}^t$
 is the unique
$cP_2, \dots, cP_{2k},cQ_1, \dots, cQ_{2l-1}$-correspondent 
of the subtower $\{1=\chi_0, \chi_1, \dots, \chi_t \}$ of \eqref{*2},
whenever  $t=0,1,\dots,n$. Also, if  $M \leq N(P_2,\dots,
P_{2k},Q_1, \dots, Q_{2l-1} \tin G)$ we have
$$M(\chi_1,\dots, \chi_{k}) = M(\chi_{1,s},\dots,\chi_{k,s}).$$
\end{remark}

\section{From triangles to towers} 
In order to complete the proof of Theorem \ref{cor:t} it suffices
to prove
\begin{theorem}\mylabel{tr--tow}
Assume that Hypothesis \ref{hyp1} holds.
Then every triangular set for \eqref{*} determines 
 a character tower of \eqref{*}, so that the tower is related to this 
triangular set via     Theorem \eqref{tow--tri}.
\end{theorem}

\begin{proof}
We will use induction on the lengh $m$ of the series \eqref{*}.
If $m=1$, then the  theorem obviously holds, as we take $\chi_1 = \beta_1$.

So assume that the theorem holds for $m=1,\dots,t$ and some integer $t \geq
0$.
We will prove it also holds for $m=t+1$.
Let 
\begin{subequations}\mylabel{tr1}
\begin{equation}\mylabel{tr1a}
1=G_0 \unlhd G_1 \unlhd \dots \unlhd G_t \unlhd G_{t+1} \unlhd 
G = G_{\infty}
\end{equation}
be a fixed normal series of $G$ that satisfies Hypothesis \ref{hyp1}.
Assume further that
\begin{equation}\mylabel{tr1b} 
\{P_{2i}, Q_{2i-1}|\alpha_{2i},\beta_{2i-1}\}
\end{equation}
\end{subequations}
is an arbitary, but fixed,  triangular set for \eqref{tr1a}.
We split the proof in two symmetric cases, according to the type of 
\eqref{tr1b}.
\begin{description}
\item[Case 1: $t$ is  odd.]
In this case the triangular set \eqref{tr1b} has the form
\begin{equation}\mylabel{tr1c}
\{P_0=1,P_2,\dots,P_{t+1},Q_1,\dots,Q_{t}| \alpha_0=1,\alpha_2,\dots,
\alpha_{t+1},\beta_1,\dots,\beta_{t}\}
\end{equation}
Note that its subset 
\begin{subequations}
\begin{equation}\mylabel{tr2a}
\{P_0=1,P_2,\dots,P_{t-1},Q_1,\dots,Q_{t}| \alpha_0=1,\alpha_2,\dots,
\alpha_{t-1},\beta_1,\dots,\beta_{t}\}
\end{equation}
is a triangular set for the series 
\begin{equation}\mylabel{tr2b}
1=G_0 \unlhd G_1 \unlhd \dots \unlhd G_t \unlhd G = G_{\infty}.
\end{equation}
Hence by the inductive hypothesis there exists a character tower 
\begin{equation}\mylabel{tr2c}
\{\chi_0,\chi_1,\dots,\chi_t\}
\end{equation}
\end{subequations}
of \eqref{tr2b}  that determines and is determined  by the set \eqref{tr2a}.
Hence,  in view of part 2) of Theorem \ref{tow--tri} (for $m=t$ and  $n=t+1$),
there is a  $cQ_1,cP_2,\dots,cP_{t-1},cQ_t$-correspondence between 
the character towers of the series \eqref{tr1a} and the character towers of
the series  
\begin{equation}\mylabel{tr3}
1=G_{0,t}, G_{1,t}, \dots,G_{t,t}\unlhd G_{t+1,t} \unlhd \Gap{t},
\end{equation}
where $G_{i,t}= G_i(\alpha_2,\dots,\alpha_{t-1},\beta_1,\dots,\beta_t)$ 
 (see \eqref{tow--tri3}),  for all 
$i=0,1,\dots,t+1,\infty$.

Let  $\Psi \in \Irr(G_{t+1}| \chi_t)$ be any  irreducible 
 character of $G_{t+1}$ lying above $\chi_t \in \Irr(G_t)$. 
Then the characters $1=\chi_0, \chi_1, \dots, \chi_t, \Psi$ form a tower 
for the series \eqref{tr1a}. Let 
$$
\chi_{0,t}, \chi_{1,t}, \dots, \chi_{t,t}, \Psi_t, 
$$
be  its  unique  $cQ_1,cP_2,\dots,cP_{t-1},cQ_t$-correspondent tower. 
So $\chi_{i, t} \in \Irr(G_{i,t})$ for all $i=0,1, \dots, t$ and 
$\Psi_{t} \in \Irr(G_{t+1, t})$. 

\smallskip
We remark that  the above is actually 
a $cQ_1,cP_2,\dots,cP_{t-1},cQ_t$-correspondence between the set
 $\Irr(G_{t+1} | \chi_{t})$ of irreducible characters $\Psi$ of $G_{t+1}$ 
lying above $\chi_{t}$ and the set  $\Irr(G_{t+1,t} |\chi_{t, t})$ of
irreducible characters $\Psi_t$ of $G_{t+1,t}$  lying above $\chi_{t,t}$.
This is clear in view of Remark \ref{t-tremark}, as the tower 
$\{\chi_{0,t}, \chi_{1,t}, \dots, \chi_{t,t}\}$ is
the unique $cQ_1,cP_2,\dots,cP_{t-1},cQ_t$-correspondent of the tower
 \eqref{tr2c}. So for any $\Psi_t \in \Irr(G_{t+1,t} |\chi_{t,t})$ the tower 
$\{\chi_{0,t}, \chi_{1,t}, \dots, \chi_{t,t}, \Psi_t,\}$ has as a 
$cQ_1,cP_2,\dots,cP_{t-1},cQ_t$-correspondent a tower of the form 
 $1=\chi_0, \chi_1, \dots, \chi_t, \Psi$ for some $\Psi  \in \Irr(G_{t+1} |\chi_t)$.

Furthermore,  according to part  4) of Theorem \ref{tow--tri} (for $i=t$ odd)
we get that $G_{t,t}=P_{t-1,t} \times Q_{t}$ while 
$\chi_{t,t} = \alpha_{t-1,t} \times \beta_t$.
Since \eqref{tr1c} is a triangular set, equation \eqref{pq15a} (for 
$r=(t+1)/2$) implies that 
$$
G_{t+1,t}=G_{t+1}(\alpha_2,\dots,\alpha_{t-1},\beta_1,\beta_{t}) =
P_{t+1} \ltimes Q_{t}. 
$$
Even more,  according to  Theorem \ref{extension}  the $P_{t+1}$-invariant  irreducible character 
$\beta_t$ of $Q_t$ has a unique canonical extension 
$\beta_t^e \in \Irr(P_{t+1} \ltimes Q_t)$. As $\alpha_{t+1} \in 
\Irr(P_{t+1})$,  the character 
$\alpha_{t+1} \cdot \beta_t^e$ is an irreducible character 
of $G_{t+1,t} =P_{t+1} \ltimes Q_t$  (see Theorem \ref{Galla}).
Also,   according to \eqref{x4} (for $r=(t+1)/2)$), 
the character $\alpha_{t+1}$ lies above the $\alpha_{t-1,t}$.
Hence the irreducible character  
$\alpha_{t+1}\cdot \beta_t^e$ of $P_{t+1} \ltimes Q_t = G_{t+1,t}$
  lies above the irreducible character $\alpha_{t-1,t}
\times \beta_t = \chi_{t, t}$ of $P_{t-1, t} \times Q_t = G_{t,t}$.
Let $\chi_{t+1} \in \Irr(G_{t+1}|\chi_t)$ be the unique 
 $cQ_1,cP_2,\dots,cP_{t-1},cQ_t$-correspondent of $\alpha_{t+1} \cdot
\beta_t^e \in \Irr(G_{t+1,t}|\chi_{t,t})$. 
So  the character tower $\{ \chi_{0,t}, \chi_{1,t}, \dots, \chi_{t,t}, 
 \alpha_{t+1} \cdot \beta_t^e \}$ of the series \eqref{tr3}
has as a  unique $cQ_1,cP_2,\dots,cP_{t-1},cQ_t$-correspondent 
the tower 
\begin{equation}\mylabel{tr4}
\{\chi_0,\chi_1,\dots,\chi_t,\chi_{t+1}\}
\end{equation}
 of the series \eqref{tr1a}.
Furthermore, the steps we followed to pick the character  $\chi_{t+1}$ 
(which are exactly the opposite  of what we used to pick $P_{t+1}$ at the inductive step
of Theorem \ref{tow--tri}) make it clear that the tower \eqref{tr4}
determines the triangular set \eqref{tr1c} in the way described in Theorem
\ref{tow--tri}.

This completes the proof of the inductive step in the case of an odd $t$.
\item[Case 2: $t$ is even.]
The proof is symmetric to that of an odd $m$, so we omit it.
\end{description}
This completes the proof of the theorem when $m=t+1$, thus proving Theorem \ref{tr--tow}.
\end{proof}

Furthermore, Theorems  \ref{tow--tri} 
and \ref{tr--tow}, along with  Corollary   \ref{daco2}, imply
\begin{remark}\mylabel{ttcong}
 Assume that the normal series  
$1=G_0 \unlhd \dots \unlhd G_{m} \unlhd G$ for $G$, along with 
the  character tower 
$\{\chi_i \in \Irr(G_i)\}_{i=0}^{m}$,  is fixed. Then 
conjugation  by any $g\in G(\chi_1, \dots, \chi_m)$ 
 leads to a new choice of the $P_{2i}$, $Q_{2i-1}$,
$\alpha_{2i}$ and $\beta_{2i-1}$ satisfying the same conditions
 (for the same $G_i$ and $\chi_i$) as the original choices.
\end{remark}

The above remark, along with Theorems  \ref{tow--tri} 
and \ref{tr--tow}, easily implies  
Theorem \ref{cor:t}.

The recursive way Theorems  \ref{tow--tri} 
and \ref{tr--tow} were  proved easily implies 
\begin{remark}\mylabel{t-tr2}
Assume that the normal series  
$1=G_0 \unlhd \dots \unlhd G_{m} \unlhd G$ for $G$, along with 
the  character tower 
$\{\chi_i \in \Irr(G_i)\}_{i=0}^{m}$,  is fixed.
Let $\{ P_{2r}, Q_{2i-1} |\alpha_{2r}, \beta_{2i-1}\}_{r=0, i=1}^{k, l}$
be a representative of the unique $G_m(\chi_1, \dots, \chi_m)$-conjugacy 
class of triangular sets that corresponds  
to the above character tower according to Theorem \ref{cor:t}.
Then $\{ \chi_i \}_{i=0}^{m-1}$ is a character tower of 
the normal series $1=G_0 \unlhd \dots \unlhd G_{m-1} \unlhd G$ of $G$.
Furthermore, the reduced  set 
$\{ P_{2r}, Q_{2i-1} |\alpha_{2r},
 \beta_{2i-1}\}_{r=0, i=1}^{[(m-1)/2], [m/2]}$,
is a  representative of the unique 
$G_{m-1}(\chi_1, \dots, \chi_{m-1})$-conjugacy class 
 of triangular sets that corresponds 
to the character tower $\{ \chi_i \}_{i=0}^{m-1}$.
\end{remark}

\begin{cor}\mylabel{t-tcor}
Assume that Hypothesis \ref{hyp1} holds. Let $\{ 1=\chi_0,\chi_1,\dots,
\chi_m\}$ be a tower of \eqref{*} and let  $\{P_{2t}, Q_{2j-1}|\alpha_{2t},
\beta_{2j-1}\}$,   for $t=0,1,\dots,k$  and $j=1,\dots,l$,  be its unique
 (up to conjugation) correspondent triangular set. 
Then 
\begin{subequations}
\begin{equation*}
Q_{2j-1} \in \Hall_{\pi'}(N(P_2, \dots, P_{2j-2},Q_1, \dots,Q_{2j-3}
\tin G_{2j-1}(\chi_1,\dots,\chi_{2j-2})), 
\end{equation*}
for all $j=1,\dots,l$,  while,  for all $t=1,\dots,k$,  we also get  that 
\begin{equation*}
Q_{2t-1} \in \Hall_{\pi'}(N(P_2, \dots, P_{2t-2},Q_1, \dots,Q_{2t-1}
\tin G_{2t}(\chi_1,\dots,\chi_{2t-1})).
\end{equation*}
Similarly,  for the $\pi$-groups we have
\begin{equation*}
P_{2t} \in \Hall_{\pi}(N(P_2, \dots, P_{2t-2},Q_1, \dots,Q_{2t-1}
\tin G_{2t}(\chi_1,\dots,\chi_{2t-1})), 
\end{equation*}
for all $t=0,1,\dots,k$,  while,  for all $j=0,1,\dots,l-1$,  we also get  that 
\begin{equation*}
P_{2j} \in \Hall_{\pi}(N(P_2, \dots, P_{2j},Q_1, \dots,Q_{2j-1}
\tin G_{2j+1}(\chi_1,\dots,\chi_{2j})).
\end{equation*}
\end{subequations}
\end{cor}

\begin{proof}                    
 Theorem \ref{tow--tri} describes completely the relations between a 
character  tower
and its corresponding triangular set.
Thus in view of \eqref{tow--tri2} (for $i=2j-1$) we  have that 
$Q_{2j-1}$ is a $\pi'$-Hall subgroup of $G_{2j-1,2j-2}$, 
 whenever $j=1,\dots,l$.
 Furthermore \eqref{tow--tri1} (for
$i=2t$) implies that $Q_{2t-1} \in \Hall_{\pi'} (G_{2t,2t-1})$,  for all
 $t=1,\dots,k$. 
But,  according to \eqref{tow--tri3}, for all such $j$ and $t$ we have 
$$G_{2j-1,2j-2}= N(P_2, \dots, P_{2j-2},Q_1, \dots,Q_{2j-3}
\tin G_{2j-1}(\chi_1,\dots,\chi_{2j-2})),
$$
while 
$$
G_{2t,2t-1} = N(P_2, \dots, P_{2t-2},Q_1, \dots,Q_{2t-1}
\tin G_{2t}(\chi_1,\dots,\chi_{2t-1})).
$$
Hence Corollary \ref{t-tcor} follows for the $\pi'$-groups.
The proof for the $\pi$-groups $P_{2t}$ is similar. 
\end{proof}

%%% Local Variables: 
%%% mode: latex
%%% TeX-master: t
%%% End: 

\section{The groups  $P^*_{2i}$: something stable in all that
mess}\mylabel{pq:sec3}

For this section we continue working under the assumptions of Hypothesis
\ref{hyp1}. So $G$ is a finite group,  while its arbitary (but fixed)
 normal series \eqref{*}
satisfies  Hypothesis \ref{hyp1}. We also fix a character tower  
\begin{subequations}
\begin{equation}\mylabel{pstar1}
1=\chi_0, \chi_1, \dots, \chi_m
\end{equation}
 of that series,  along with its uniquely determined (up to conjugation) 
triangular set
\begin{equation}\mylabel{pstar2}
\{P_0,P_2, \dots, P_{2k}, Q_1, \dots, Q_{2l-1} | \alpha_0, \alpha_2, \dots, 
\alpha_{2k}, \beta_1,\dots,\beta_{2l-1}\}, 
\end{equation} 
\end{subequations}
where $k$ and $l$ are defined as in \eqref{kl:def} 
(i.e., $2l-1$ and $2k$ are the greatest odd and even, respectively, integers in
the set $\{1, \dots, m\}$).
For the normal series \eqref{*} and its character tower \eqref{pstar1}, 
  Theorem \ref{prel:t1} can be applied. So for any $i=1,\dots, m$ we write 
 \begin{equation*}
G_i^* = G_i(\chi  _1, \chi  _2 , \dots ,\chi  _{i-1})
\end{equation*}
and 
\begin{subequations}
\begin{equation}\mylabel{pqinf}
G_{\infty}^*= G^*= G(\chi_1,\dots,\chi_{m})
\end{equation}
for the stabilizers  of $\chi  _1, \chi  _2 , \dots ,\chi  _{i-1}$ 
and  $\chi  _1, \chi  _2 , \dots ,\chi  _{m}$  in $G_i$ and $G$, respectively.
As $G_i \leq  G_j$ for all $j$ with $0\leq i \leq j \leq m$, the group $G_i$ 
fixes the characters $\chi_j$ for  all such $j$.
Hence
\begin{equation}\mylabel{pq1}
G_i^*= G_i(\chi_1,\dots,\chi_{i-1}) = G_i(\chi_1,\dots, \chi_{m}).
\end{equation}
\end{subequations} 

Then,  in view of  Theorem \ref{prel:t1},  we have that 
 $G_0^* = G_0 = 1 , G_1^* = G_1 $ and $G_j^* =G_i^* \cap G_j 
\unlhd G_i^*$,  whenever $0\leq j\leq i\leq m$.
Furthermore, there exist  unique characters $\chi _i^*$, for 
$i=1, \dots ,m$, such that
\begin{subequations} \mylabel{pq2}
 \begin{equation}\mylabel{pq2a}
\chi _i^* \in \Irr(G_i^*) \text{ lies over } \chi_1^*, \dots ,\chi_{i-1}^* 
\text{ and induces } \chi _i,  
\end{equation}
\begin{equation} \mylabel{pq2b}
G_i^* = G_i( \chi  _1^*, \chi  _2^* , \dots ,\chi  _{i-1}^*) =
 G_i(\chi  _1^*, \chi  _2^* , \dots ,\chi  _{m}^*).
\end{equation}
We note here that, according to the same theorem,   
\begin{equation}\mylabel{pq2c}
G_1^* = G_1 \text{ and } \chi_1^* = \chi_1 \in \Irr(G_1^*) =\Irr(G_1).
\end{equation}
\end{subequations}

According to \eqref{pq14'b}, for each  $i=1, \dots , k$ the group $P_{2i}$ 
normalizes all the
 previously chosen $\pi$-groups $P_2, P_4, \dots ,P_{2i-2}$. Hence the
product: 
\begin{equation}\mylabel{pqp*d}
 P_{2i}^* = P_2 \cdot  P_4 \cdots  P_{2i}
\end{equation}
is a group. We also define  $P_0^*:=1$, therefore $P_{2i}^*$ is defined for all $i=0,1,\dots,k$.
As  will become clear, these groups play the most 
important role in the construction we did in the previous section. The reason
is that they 
are the only groups that  remain unchanged when we change the $\pi'$-parts 
(groups and characters ) of the triangular set \eqref{pstar2}.   
The way we defined the groups $P_{2i}^*$ uses the individual $P_{2t}$ for
all $t$  with 
$1\leq t \leq i$. But, as the proposition that follows shows, we could
have picked the 
 groups $P^*_{2i}$ using only the groups $G^*_{2i}$.

\begin{proposition}\mylabel{pq3-}
The group $P^*_{2i}$ is a $\pi$-Hall subgroup of $G^*_{2i}$ 
 whenever  $i= 0,1, \dots ,k$. It  is also a $\pi$-Hall subgroup
of $G_{2i+1}^*$ for all $i=0,1,\dots,l-1$.
Furthermore $P^*_{2r} = P^*_{2i} \cap G^*_{2r}$ and thus  $P^*_{2r} \unlhd
P^*_{2i}$, whenever  $1 \leq r \leq i \leq k$.
\end{proposition}
To prove Proposition \ref{pq3-} we need the following lemma:
\begin{lemma}\mylabel{pq3l}
If $j=2,\dots,m,\infty$, and $s$ are  such that $2\leq 2s \leq j$, then 
\begin{multline*}
G_j^*=N(P_2, \dots ,P_{2s-2}, Q_1,  Q_3, \dots ,Q_{2s-1} \tin G^*_j)\cdot G^*_{2s-
1}  \\
 = N(P_2, \dots ,P_{2s}, Q_1, Q_3,  \dots ,Q_{2s-1} \tin G^*_j)\cdot
G^*_{2s}. 
\end{multline*}
\end{lemma}
\begin{proof}
We will use induction on $s$.
For $s= 1$ the group $N(P_2, \dots ,P_{2s-2},Q_1, \dots ,Q_{2s-1}
\tin G^*_j)$ equals $N(Q_1 \tin G^*_j) = G^*_j$, while the normalizer 
 $N(P_2, \dots ,P_{2s}, Q_1, \dots ,Q_{2s-1} \tin G^*_j)$  equals  the group 
$N(P_2,Q_1 \tin G^*_j) =
 N(P_2 \tin G^*_j)$.
According to \eqref{x2} and \eqref{x3}, we have  that 
$ P_2 \in \Hall_{\pi}(G_2(\chi_1))=\Hall_{\pi}(G_2^*)$.
Therefore, 
 for any $j \geq 2$,   the Frattini argument implies that 
$G^*_j = N(P_2 \tin G^*_j)\cdot G^*_2$,  as $G^*_2$ is a normal subgroup of
$G^*_j$.
 Thus Lemma \ref{pq3l} holds  when  $s = 1$ and $j=2,\dots,m,\infty$.

Assume  now  that Lemma \ref{pq3l}  holds for all $s = 1, 2, \dots ,t-1$, 
where $2 < 2t \leq j$. We
will prove that it also holds when $s = t$. 
By induction, for $s = t-1$, we get that $G^*_j =
 N(P_2, \dots ,P_{2t-2}, Q_1,  Q_3, \dots ,Q_{2t-3} \tin G^*_j)\cdot G^*_{2t-
2}$. But 
$$N(P_2, \dots ,P_{2t-2}, Q_1,  Q_3, \dots 
,Q_{2t-3} \tin G^*_{2t-1}) \unlhd N(P_2, \dots ,P_{2t-2}, Q_1,  Q_3, \dots ,Q_{2t-3} \tin G^*_j).
$$
Furthermore, Corollary \ref{t-tcor} implies that  
the group 
$ N(P_2, \dots ,P_{2t-2}, Q_1,  Q_3, \dots  ,Q_{2t-3} \tin G^*_{2t-1})$
has 
  $Q_{2t-1}$ as  a $\pi'$-Hall 
subgroup.
 Hence, by  the Frattini argument,  we have
\begin{multline*}
N(P_2, \dots ,P_{2t-2}, Q_1, Q_3, \dots ,Q_{2t-3} \tin G^*_j) = \\
 N(P_2, \dots ,P_{2t-2}, Q_1,  Q_3, \dots ,Q_{2t-3}, Q_{2t-1} \tin G^*_j) \cdot
 N(P_2, \dots ,P_{2t-2}, Q_1,  Q_3, \dots ,Q_{2t-3} \tin G^*_{2t-1}).
\end{multline*}
Therefore,  
\begin{multline*}
G^*_j =
 N(P_2, \dots ,P_{2t-2}, Q_1,  \dots ,Q_{2t-3} \tin G^*_j)\cdot G^*_{2t-2} \\
 =N(P_2, \dots ,P_{2t-2}, Q_1,   \dots ,Q_{2t-3}, Q_{2t-1} \tin G^*_j) \cdot
 N(P_2, \dots ,P_{2t-2}, Q_1,  \dots ,Q_{2t-3} \tin G^*_{2t-1})\cdot G^*_{2t-
2}.
\end{multline*}
By induction, $N(P_2, \dots ,P_{2t-2}, Q_1,  Q_3, \dots ,Q_{2t-3} \tin
G^*_{2t-1})\cdot G^*_{2t-2} = G^*_{2t-1}$.
Hence, 
\begin{equation}\mylabel{pq3le}
G^*_j = N(P_2, \dots ,P_{2t-2}, Q_1,  Q_3, \dots ,Q_{2t-3}, Q_{2t-1} \tin
G^*_j) \cdot G^*_{2t-1}.
\end{equation}
 This proves the first  equality in
Lemma \ref{pq3l} for $s= t$.
 
It remains to show that $G^*_j = N(P_2, \dots, P_{2t}, Q_1, Q_3, \dots Q_{2t-
1} \tin G^*_j)\cdot G^*_{2t}$.
By Corollary \ref{t-tcor}, the group $P_{2t}$ is a $\pi$-Hall
subgroup of $N(P_2, \dots ,P_{2t-2}, Q_1, 
Q_3, \dots ,Q_{2t-1} \tin G^*_{2t})$. Since $N(P_2, \dots ,P_{2t-2}, Q_1,
\dots ,Q_{2t-1} \tin
 G^*_{2t}) \unlhd  N(P_2, \dots ,P_{2t-2}, Q_1, \dots ,Q_{2t-1} \tin G^*_j)$, 
   Frattini's argument  implies that 
\begin{multline*}
N(P_2, \dots ,P_{2t-2}, Q_1, \dots ,Q_{2t-1} \tin G^*_j) =\\
 N(P_2, \dots ,P_{2t}, Q_1, \dots ,Q_{2t-1} \tin G^*_j)\cdot N(P_2, \dots
,P_{2t-2}, Q_1, \dots
 ,Q_{2t-1} \tin G^*_{2t}). 
\end{multline*}

The above equation, along with \eqref{pq3le}, implies 
\begin{multline*}
G^*_j  \\
 =N(P_2, \dots ,P_{2t}, Q_1, \dots ,Q_{2t-1} \tin G^*_j)\cdot N(P_2, \dots
,P_{2t-2}, Q_1,
 \dots ,Q_{2t-1} \tin G^*_{2t})\cdot G^*_{2t-1}  \\
= N(P_2, \dots, P_{2t}, Q_1, Q_3, \dots ,Q_{2t-1} \tin G^*_j)\cdot G^*_{2t}.
\end{multline*}
 This proves the remaining equality in  Lemma \ref{pq3l}
for $s =t$. Hence the inductive proof of that lemma is complete.
\end{proof}

\begin{proof}
[Proof of Proposition \ref{pq3-}]
 We will use induction on $i$.
For $i=0$ it is trivially true as $P_0^*=1=G_0^*$, while $G_1^*$ is a $\pi'$-group.
If $i= 1$, then $P^*_{2} = P_2$. By \eqref{x3} and \eqref{x2} the
latter group is a $\pi$-Hall subgroup of $G_2(\chi_1)=G_2^*$. Furthermore, as 
$G^*_3 / G^*_2$ is a $\pi'$-group,
$P^*_2$ is also a $\pi$-Hall subgroup  of $G^*_3$.
The rest of the proposition holds trivially for $i = 1$.

Now we assume that Proposition \ref{pq3-} holds for all 
$i = 1,2,\dots,t-1$, where  $1 < t \leq k$.  We will prove that it  holds
 for $i = t$.
We have $P^*_{2t-2} \in \Hall_{\pi}(G^*_{2t-2})$ by induction, and
$P_{2t} \in \Hall_{\pi}(N(P_2, \dots, P_{2t-2},Q_3,\dots ,Q_{2t-1} \tin
G^*_{2t}))$  by  Corollary \ref{t-tcor}.
Since $P_{2t-2}^* =P_2 \cdots P_{2t-2}$, it follows that 
 the group $P^*_{2t} = P^*_{2t-2} \cdot P_{2t}$ is a $\pi$-subgroup of
$G^*_{2t}$.
So there exists a $\pi$-Hall subgroup, $\map_{2t}$,  of $G^*_{2t}$ with 
\begin{equation}\mylabel{pq3-e2}
P^*_{2t} \leq \map_{2t}. 
\end{equation}
Since $G^*_{2t-2}$ is a normal subgroup of $G^*_{2t}$,  we conclude that 
$\map_{2t}\cap G^*_{2t-2}$ is a $\pi$-Hall subgroup of $G^*_{2t-2}$.
But $P^*_{2t-2}$ is a $\pi$-Hall subgroup of $G^*_{2t-2}$ such that 
 $P^*_{2t-2} \leq P^*_{2t} \cap G^*_{2t-2} \leq \map_{2t} \cap G^*_{2t-2}$.
Hence 
$$ 
P^*_{2t-2} = P^*_{2t} \cap G^*_{2t-2} = \map_{2t} \cap G^*_{2t-2}. 
$$
Furthermore, as  $G^*_{2t-1} /G^*_{2t-2}$ is a $\pi'$-group, $P^*_{2t-2}$ is
not only a $\pi$-Hall
 subgroup of $G^*_{2t-2}$,  but also of $G^*_{2t-1}$ (note that $G^*_{2t-1}$
exists for all $t=1,\dots,k$).
Hence
\begin{equation}\mylabel{pq3-e}
P^*_{2t-2} =  P^*_{2t} \cap G^*_{2t-1} = \map_{2t} \cap G^*_{2t-1} \in
\Hall_{\pi}(G^*_{2t-1}).
\end{equation}
Since   $G^*_{2t}/ G^*_{2t-1}$ is a $\pi$-group, 
and  $\map_{2t}$ is a $\pi$-Hall subroup of  $G^*_{2t}$ we have that 
$G^*_{2t} = \map_{2t} \cdot G^*_{2t-1}$.

\noindent
Furthermore, $G^*_{2t} = N(P_2, \dots ,P_{2t-
2}, Q_1, 
\dots ,Q_{2t-1} \tin G^*_{2t}) \cdot G^*_{2t-1}$ according to Lemma \ref{pq3l}.
 Hence,  the $\pi$-Hall subgroup $P_{2t}$ of  $ N(P_2, \dots ,P_{2t-2}, Q_1, \dots
Q_{2t-1} \tin
G^*_{2t})$ also  covers $G^*_{2t}/G^*_{2t-1}$. 
As $P^*_{2t-2} \leq G^*_{2t-2} \leq G^*_{2t-1}$, we conclude that 
$$
\map_{2t} \cdot G^*_{2t-1} = G^*_{2t} = P_{2t} \cdot G^*_{2t-1} 
= P_{2t}^*\cdot G^*_{2t-1}.
$$
This, along with \eqref{pq3-e} and \eqref{pq3-e2}, implies that $P^*_{2t}  =
\map_{2t}$. Thus $P^*_{2t}$ is a $\pi$-Hall subgroup of $G^*_{2t}$.
If $t \leq l-1$ then  the group $G^*_{2t+1}$ is defined and 
 $G^*_{2t+1}/G^*_{2t}$ is a $\pi'$-group. We conclude that  
$P^*_{2t} \in \Hall_{\pi}(G^*_{2t+1})$ for all such $t$.
The rest of Proposition \ref{pq3-} for $i = t$
 follows easily, as $G^*_{2r} \unlhd G^*_{2t}$, 
and 
$P^*_{2r} \leq P^*_{2t}$, whenever $1\leq r \leq t \leq k$.
Hence, the inductive proof of the proposition is complete.
\end{proof} 

Along with the groups $P^*_{2i}$ we have irreducible characters $\alpha_{2i}^*$
that correspond uniquely  to the irreducible characters $\alpha_{2i}$ of
$P_{2i}$,
for all $i=1, \dots  ,k$. To prove that these characters exist and show how
their  
correspondence with the $\alpha_{2i}$ works,  we will use the following lemma:
\begin{lemma}\mylabel{p*22}
If $1 \leq   i \leq  t \leq k$ then
\begin{equation}\mylabel{pq22} 
N(Q_{2i-1} \tin  P_{2i-2}\cdot P_{2i}\cdots P_{2t}) = P_{2i}\cdots P_{2t}.
\end{equation}
\end{lemma}

\begin{proof}
According to \eqref{pq14'b}, the group  $P_{2j}$ normalizes $Q_{2i-1}$ for all $j
\geq i$.
Hence the product  $ P_{2i}\cdots P_{2t}$
is contained in the normalizer $N(Q_{2i-1} \tin P_{2i-2}\cdots P_{2t})$.
 So, 
$ N(Q_{2i-1} \tin P_{2i-2} \cdots P_{2t}) = N(Q_{2i-1} \tin
 P_{2i-2}) \cdot P_{2i}\cdots P_{2n}$.
By \eqref{pq13b}, we have 
$N(Q_{2i-1} \tin P_{2i-2}) =
P_{2i-2, 2i-1} = P_{2i-2}\cap P_{2i}$. Therefore $N(Q_{2i-1} \tin
 P_{2i-2})$ is a subgroup of $P_{2i} \leq P_{2i}\cdots P_{2t}$.
This completes the proof of Lemma \ref{p*22}.
\end{proof}

\begin{proposition}
For all $i,t$ with $1 \leq i  \leq t \leq k$ we have
\begin{equation}\mylabel{pq22i}
N(Q_1, Q_3, \dots, Q_{2i-1} \tin P^*_{2t}) = P_{2i}\cdots P_{2t}.
\end{equation}
\end{proposition}
\begin{proof}
The proof is a multiple application of Lemma \ref{p*22}.
 \begin{align*}
&N(Q_1, Q_3,  Q_5,\dots ,Q_{2i-1} \tin P^*_{2t})& \\
&= N(Q_3, Q_5, \dots ,Q_{2i-1} \tin P^*_{2t})   &\text{ since $Q_1 \unlhd G$,  } \notag \\
&= N(Q_5, \dots ,Q_{2i-1} \tin N(Q_3 \tin P^*_{2t})) &\,  \notag\\
&= N(Q_5, \dots ,Q_{2i-1} \tin P_4\cdots P_{2t}) &\text{ in view of 
Lemma \ref{p*22}, } \notag\\
&= N(Q_7, \dots ,Q_{2i-1} \tin N(Q_5 \tin P_4\cdots P_{2t}))
& \notag\\
&\dots & \notag\\
&= N(Q_{2i-1} \tin P_{2i-2}\cdots P_{2t}) &\, \notag\\
&= P_{2i}\cdots  P_{2t} &\text{ in view of Lemma \ref{p*22}. }
\end{align*} 
\end{proof}

In particular, we have a way to recover the $P_{2i}$ from the products
$P^*_{2i}$ and the 
$q$-groups $Q_3, \dots ,Q_{2i-1}$, whenever $i=1,\dots,k$. 
Indeed, \eqref{pq22i} implies: 
\begin{equation}\mylabel{pq22ii}
N(Q_1, Q_3, \dots ,Q_{2i-1} \tin P^*_{2i}) = P_{2i}.
\end{equation} 

\begin{lemma}\mylabel{p*l1} 
If \, $1 \leq j < i \leq k$, then the product
$Q_{2j+1} \cdot P_{2j} \cdot P_{2j+2} \cdots P_{2i}$ is a subgroup of
$G$ having $Q_{2j+1}$ as a Hall $\pi'$-subgroup, and $P_{2j} \cdot 
P_{2j+2} \cdots P_{2i}$ as a Hall $\pi$-subgroup. Both the  $\pi$-group
$P_{2j}$ and the product  $Q_{2j+1} \cdot P_{2j}$
 are normal subgroups of $Q_{2j+1} \cdot P_{2j} \cdot P_{2j+2} \cdots P_{2i}$.
 Furthermore,
$N(Q_{2j+1} \tin P_{2j} \cdot P_{2j+2} \cdots P_{2i}) = P_{2j+2} 
\cdots P_{2i}$. Hence Theorem \ref{dade:t1}  gives us a one to one 
$Q_{2j+1}$-correspondence 
$$
 \alpha^*_{2i,2j-1} \underset {Q_{2j+1}}{\longleftrightarrow}
\alpha_{2i,2j+1}^*
 $$
between all characters 
$\alpha^*_{2i,2j+1} \in \Irr(P_{2j+2} \cdots 
P_{2i})$ and all characters   $\alpha^*_{2i,2j-1} \in \Irr(P_{2j} \cdot 
P_{2j+2} \cdots P_{2i})$  lying over some character  $\alpha^*_{2j,2j-1} 
\in \Irr^{Q_{2j+1}}(P_{2j})$. This correspondence is invariant under
conjugation by elements of any subgroup $K \le G$ normalizing both
$Q_{2j+1}$ and $P_{2j} \cdot P_{2j+2} \cdots P_{2i}$.
Furthermore, if $\alpha_{2s,2j-1}^* \in \Irr(P_{2j}\cdots P_{2s})$ is any 
character of $P_{2j} \cdots P_{2s}$ lying under $\alpha_{2i,2j-1}^*$ and above 
$\alpha_{2j,2j-1}^*$, for some $s$ with $1\leq j <s \leq i \leq k$, then 
its  $Q_{2j+1}$-correspondent $\alpha_{2s,2j+1}^*$ lies under the 
$Q_{2j+1}$-correspondent $\alpha_{2i,2j+1}^*$ of $\alpha_{2i,2j-1}^*$.
\end{lemma}

\begin{proof}
We only need to show that $Q_{2j+1}\cdot P_{2j}\cdot P_{2j+2} \cdots P_{2i}$
is a group  having  $P_{2j}$ and  $Q_{2j+1} \cdot P_{2j}$ 
as normal subgroups, while   
$P_{2j} \cdots P_{2s} \unlhd P_{2j} \cdots P_{2i}$ whenever $1\leq j < s
\leq  i \leq k$. The rest is an 
easy  application of Theorem \ref{dade:t2} and \eqref{pq22i}.

By \eqref{pq14'a} the group $Q_{2j+1}$ normalizes the group $P_{2j}$.
Furthermore, \eqref{pq14'b} implies that the groups 
$P_{2j+2}, P_{2j+4}, \dots ,P_{2i}$ normalize both $P_{2j}$ and $Q_{2j+1}$, 
while $P_{2s}$ normalizes $P_{2t}$ for all $s,t$ with $j \le t \le s \leq i$.
Hence, the products $Q_{2j+1}\cdot P_{2j}$ and  $P_{2j+2}\cdots P_{2s}$
 form groups.
Even more, the latter group normalizes the former one for every $s=j+1,\dots,i$.
Hence $Q_{2j+1}\cdot P_{2j}\cdot P_{2j+2} \cdots P_{2i}$ forms  a group
having  $P_{2j}$ and  $Q_{2j+1}\cdot P_{2j}$ as normal subgroups.
It is also clear that the product group $P_{2j} \cdots P_{2s}$ is a normal
subgroup of $P_{2j} \cdots P_{2i}$ whenever $j < s \leq i$.
 \end{proof}

\begin{theorem}\mylabel{p*t1}
For any $i=2,3,\dots,k$ we may form a chain
\begin{equation}\mylabel{p*teq}
 \alpha_{2i,1}^*  \underset {Q_{3}}{\longleftrightarrow}
\alpha_{2i,3}^*  \underset {Q_{5}}{\longleftrightarrow}
\alpha_{2i,5}^* \underset {Q_{7}}{\longleftrightarrow}
\cdots \underset {Q_{2i-3}}{\longleftrightarrow}
 \alpha^*_{2i,2i-3} \underset {Q_{2i-1}}{\longleftrightarrow}
\alpha_{2i,2i-1}^*
\end{equation}
of the $Q_{2j+1}$-correspondences in Lemma  \ref{p*l1}. The composite 
$Q_3, Q_5,\dots,Q_{2i-1}$-correspondence is  one to one between all 
characters $\alpha_{2i,2i-1}^* \in \Irr(P_{2i})$ and all characters
$\alpha_{2i,1}^* \in \Irr(P_{2i}^*)$ having the following property:
\begin{property}\mylabel{p*pp1}
There exist characters $\alpha_{2j,1}^* \in \Irr(P_{2j}^*)$, for
$j=1,2,\dots,i-1$, such that each $\alpha_{2j,1}^*$ is 
$Q_{2j+1}$-invariant and lies under $\alpha_{2j+2,1}^*$.
\end{property}
This $Q_3,Q_5,\dots,Q_{2i-1}$-correspondence is invariant under conjugation by 
elements of any subgroup $K \le G$ normalizing all the subgroups
$Q_3,Q_5,\dots,Q_{2i-1}$ and $P_{2i}^*$.
\end{theorem}

\begin{proof} 
It follows immediately from Lemma \ref{p*l1} that the composite
correspondence is one to one between all characters in $\Irr(P_{2i})$
 and some characters in $\Irr(P_{2i}^*)$. It is also clear that this
correspondence is invariant under conjugation by elements of $N(Q_3,Q_5,\dots, 
Q_{2i-1},P_{2i}^* \tin G)$. Thus it remains to show that the image of the 
composite correspondence is exactly the subset of all $\alpha_{2i,1}^* 
 \in \Irr(P_{2i}^*)$
having Property \ref{p*pp1}.

We will first show that 
any $\alpha_{2i,1}^* \in \Irr(P_{2i}^*)$ which is the
 $Q_3,Q_5,\dots,Q_{2i-1}$-correspondent of some $\alpha_{2i,2i-1}^* \in 
\Irr(P_{2i})$ must satisfy Property \ref{p*pp1}.
This can be done by induction on $i$. 
Assume that  $i=2$.
 Let   $\alpha_{4,1}^* \in \Irr(P_4^*)$ be  the 
$Q_3$-correspondent 
of some $\alpha_{4,3}^* \in \Irr(P_4)$.  
In this case, Lemma \ref{p*l1} (for $j=1$ and $i=2$), 
describes   this  $Q_3$-correspondence
as one between all characters $\alpha_{4,3}^* \in \Irr(P_4)$ 
and all characters $\alpha_{4,1}^* \in \Irr(P_4)$ lying over some character 
$\alpha_{2,1}^* \in \Irr^{Q_3}(P_2)$.
As $P_2=P_2^*$, we obviously have that $\alpha_{4,1}^*$ lies above the 
$Q_3$-invariant character $\alpha_{2,1}^*$ of $P_2^*$.
Thus $\alpha_{4,1}^*$ 
satisfies Property \ref{p*pp1}.
 This implies the $i=2$ case.

Assume that is true for all $i=2,\dots,t-1$, for some $t$ with  $2 < t \leq k$. We will prove it
also holds when $i=t$.

Let  $\alpha_{2t,1}^* \in \Irr(P_{2t}^*)$
be the $Q_3,\dots,Q_{2t-1}$-correspondent of 
  some  $\alpha_{2t,2t-1}^* \in \Irr(P_{2t})$.
Then,  according to  Lemma \ref{p*l1},  
the character $\alpha_{2t,2t-1}^*$ of $ P_{2t}$ 
has as a $Q_{2t-1}$-correspondent a character $\alpha_{2t,2t-3}^*$ 
of $P_{2t-2}\cdot P_{2t}$, that lies above some character 
$\alpha_{2t-2,2t-3}^* \in \Irr^{Q_{2t-1}}(P_{2t-2})$.
Let $\alpha_{2t-2,1}^* \in \Irr(P_{2t-2}^*)$ 
be the $Q_3, Q_5, \dots,Q_{2t-3}$-correspondent of 
$\alpha_{2t-2,2t-3}^* \in \Irr(P_{2t-2})$,
 that we get by multiple applications 
of   Lemma \ref{p*l1}. 
As $\alpha_{2t-2,2t-3}^*$ lies under $\alpha_{2t,2t-3}^*$, we get that 
the $Q_3,\dots,Q_{2t-3}$-correspondent $\alpha_{2t-2,1}^*$ 
of $\alpha_{2t-2,2t-3}^*$ lies under the $Q_3,\dots,Q_{2t-3}$-correspondent
$\alpha_{2t,1}^*$ of $\alpha_{2t,2t-3}^*$ (see Theorem \ref{dade:t2}).
Furthermore, $Q_{2t-1}$ fixes 
$\alpha_{2t-2,2t-3}^*$ and normalizes the groups 
$Q_3,\dots,Q_{2t-3}$, as well as the product group 
$P_2 \cdots P_{2t-2}$. Hence it also fixes  the 
$Q_3, \dots,Q_{2t-3}$-correspondent 
$\alpha_{2t-2,1}^*$ of $\alpha_{2t-2,2t-3}^*$. 
In conclusion,  the character $\alpha_{2t-2,1}^*$ is $Q_{2t-1}$-invariant 
and lies under $\alpha_{2t,1}^*$.
Thus $\alpha_{2t,1}^*$ satisfies Property \ref{p*pp1} for $j=t-1$. 

As $\alpha_{2t-2,1}^*$ is the $Q_3,\dots,Q_{2t-3}$-correspondent 
of $\alpha_{2t-2,2t-3}$, the inductive hypothesis applies. 
Hence  for every $j=1,\dots,t-2$, 
there exists a character
$\alpha_{2j,1}^* \in \Irr(P_{2j}^*)$ that is $Q_{2j+1}$-invariant
and lies under $\alpha_{2j+2,1}^*$.
The existence of these   characters $\alpha_{2j,1}^*$, along with $\alpha_{2t-2,1}^*$,
implies that $\alpha_{2t,1}^*$ satisfies  Property \ref{p*pp1}.
This completes the inductive proof that an $\alpha_{2i,1}^* \in \Irr(P_{2i}^*)$ which is the
 $Q_3,Q_5,\dots,Q_{2i-1}$-correspondent of some $\alpha_{2i,2i-1}^* \in 
\Irr(P_{2i})$ must satisfy Property \ref{p*pp1}.

Now assume that a  character $\alpha_{2i,1}^* \in \Irr(P_{2i}^*)$ satisfying 
Property \ref{p*pp1} is fixed, for some $i=2,3,\dots,k$. 
We want to construct some  character $\alpha_{2i,2i-1}^* \in \Irr(P_{2i})$  having
$\alpha_{2i,1}^*$ as its $Q_3,\dots,Q_{2i-1}$-correspondent.
To do this,  we will find  
characters $\alpha_{2i,2t+1}^*$, for every $t=1,\dots,i-1$,  to form the chain
\eqref{p*teq}.  So it suffices to  show that for every $t$ with 
$1\leq t \leq i-1$  we can apply
Lemma \ref{p*l1} to  get a 
$Q_{2t+1}$-correspondent $\alpha_{2i, 2t+1}^*$ of  $\alpha_{2i,2t-1}^*$.
This is done in a recursive way. 
So, for  $t=1$ we observe that Property \ref{p*pp1} implies that
$\alpha_{2i,1}^*$ lies above the $Q_3$-invariant character $\alpha_{2,1}^* \in 
\Irr(P_{2})=\Irr(P_{2}^*)$. Hence, according to  Lemma \ref{p*l1}, the character 
$\alpha_{2i,1}^* \in \Irr(P_2 \cdots P_{2i})$  has a $Q_3$-correspondent 
character $\alpha_{2i,3}^* \in \Irr(P_4 \cdots P_{2i})$.
Furthermore, as $\alpha_{2s,1}^*$  lies under $\alpha_{2s+2,1}^*$ and 
above $\alpha_{2,1}^*$,  for every $s=2,\dots,i-1$, (by Property \ref{p*pp1}),
 the same lemma implies that  the $Q_3$-correspondent
  $\alpha_{2s,3}^* \in \Irr(P_4 \cdots P_{2s})$
 of $\alpha_{2s,1}^* \in \Irr(P_2 \cdots P_{2s})$ is defined and lies
under $\alpha_{2s+2,3}^*$. Even more, the character $\alpha_{2s,1}^*$ is
$Q_{2s+1}$-invariant and $Q_{2s+1}$  normalizes
both  $Q_3$ and the product group $P_2 \cdots P_{2s}$, 
for all $s=2,\dots, i-1$. 
Hence the $Q_3$-correspondent  character $\alpha_{2s,3}^*$ of $\alpha_{2s,1}^*$
is also $Q_{2s+1}$-invariant.

We can now do the case $t=2$. Indeed, the previous comment
for $s=2$ implies that  $\alpha_{2s,3}^* \in
\Irr(P_4 \cdots P_{2i})$ lies  above the $Q_5$-invariant character
$\alpha_{4,3}^* \in \Irr(P_4)$.  Thus we can apply Lemma \ref{p*l1} again to
get a $Q_5$-correspondent character $\alpha_{2i,5}^* \in \Irr(P_6 \cdots P_{2i})$.
Note that  the $Q_5$-correspondent $\alpha_{2s,5}^* \in \Irr(P_6 \cdots
P_{2s})$ of $\alpha_{2s,3}^* \in \Irr(P_4 \cdots P_{2s})$
 is defined whenever  $s=3,\dots,i-1$. This correspondent  
lies under $\alpha_{2s+2,5}^*$, 
as $\alpha_{2s,3}^*$ lies under $\alpha_{2s+2,3}^*$. Furthermore
$\alpha_{2s,3}^*$ is $Q_{2s+1}$-invariant,  while $Q_{2s+1}$
 normalizes both $Q_5$ and the product $P_4 \cdots
P_{2s}$. Therefore the $Q_5$-correspondent  $\alpha_{2s,5}^*$ of $\alpha_{2s,3}^*$
 is also $Q_{2s+1}$-invariant.
So for $s=3$ we have that $\alpha_{2i,5} \in \Irr(P_6 \cdots P_{2i})$ lies  
over the $Q_7$-invariant character $\alpha_{6,5} \in \Irr(P_6)$. 
Hence we can apply Lemma \ref{p*l1} again  and thus get the desired
correspondence for $t=3$.
We continue similarly. At the $t$-step we have the character 
 $\alpha_{2i,2t-1}^* \in \Irr(P_{2t} \cdots P_{2i} )$ lying over the 
$Q_{2t+1}$-invariant character $\alpha_{2t,2t-1}^* \in \Irr(P_{2t})$, while
for all $s=t,\dots, i-1$ the character $\alpha_{2s,2t-1}^* \in \Irr(P_{2t}
\cdots P_{2s})$ is $Q_{2s+1}$-invariant and lies under  
$\alpha_{2s+2,2t-1}^*$. 

At the last step for $t=i-1$ we end up  with the character $\alpha_{2i,2i-3}^* 
\in \Irr(P_{2i-2}\cdot P_{2i})$ lying over the $Q_{2i-1}$-invariant character 
$\alpha_{2i-2,2i-3}^* \in \Irr(P_{2i-2})$. So the final application of Lemma 
\ref{p*l1} will provide a $Q_{2i-1}$-correspondent character
$\alpha_{2i,2i-1}^* \in \Irr(P_{2i})$ of $\alpha_{2i,2i-3}^*$, that is
actually a $Q_3,Q_5,  \dots ,Q_{2i-1}$-correspondent of $\alpha_{2i}^*$.

This completes the proof of Theorem \ref{p*t1}.
\end{proof}

\begin{remark}
The chain $Q_3, Q_5, \dots, Q_{2i-1}$  is empty when $i=1$. In that case  we
define the $Q_3,Q_5, \dots, Q_{2i-1}$-correspondence to be the 
identity correspondence between $\Irr(P_{2i})= \Irr(P_2)$ and the equal set
$\Irr(P_{2i}^*) = \Irr(P_{2}^*)$.
\end{remark}

Now  we can make the 
\begin{defn}\mylabel{p*d1}
For each $i=1,2,\dots,k$ we define $\alpha_{2i}^* \in \Irr(P_{2i}^*)$
 to be the $Q_3, Q_5, \dots Q_{2i-1}$-correspondent  of the character
$\alpha_{2i} \in \Irr(P_{2i})$.
\end{defn}
So
\begin{equation}\mylabel{p*a}
\alpha_2^* =\alpha_2
\end{equation}
by convention. Furthermore Theorem \ref{p*t1} obviously implies 
\begin{proposition}\mylabel{pqp*t}
If a subgroup $K \le G$ normalizes all the subgroups  $Q_3, Q_5, \dots, 
Q_{2i-1}$ and $P_{2i}^*$, for some $i=1,2,\dots,k$,  then $\alpha_{2i}^*$ and 
$\alpha_{2i}$ have the same stabilizer $K(\alpha_{2i}^*)= K(\alpha_{2i})$,
in $K$.
\end{proposition}

\begin{corollary}\mylabel{p*cor1}
\begin{subequations} \mylabel{pq23}
\begin{equation} \mylabel{pq23a}
Q_{2j+1}(\alpha_{2i}^*) = Q_{2j+1}(\alpha_{2i})=Q_{2j+1}
\end{equation}
and 
\begin{equation} \mylabel{pq23b}
P_{2s}(\alpha_{2i}^*) = P_{2s}(\alpha_{2i})=P_{2s}.
\end{equation}
\end{subequations}
whenever $1\leq i \le j <l$ and $1\le i \le s \leq k$.
\end{corollary} 

\begin{proof}
According to \eqref{pq14'}, the groups  $Q_{2j+1}$ and $P_{2s}$
normalize    $P_2, \dots ,P_{2i}$ and $Q_3, \dots Q_{2i-1}$, for all $j,s$ 
with $1\leq i \le j <l$ and $1\le i \le s \leq k$. 
Thus they also normalize the groups
$P_2^*, \dots ,P_{2i}^*$. Hence Proposition \ref{pqp*t}, along with the fact that 
$Q_{2j+1}$ and $P_{2s}$  fix $\alpha_{2i}$ (see (\ref{xxx}c,e)), implies 
Corollary \ref{p*cor1}.
\end{proof}

The next proposition  shows how the  characters $\alpha_{2i}^*$ are related.

\begin{lemma}\mylabel{p*l2} 
If $i=2,3,\dots,k$, then $\alpha_{2i-2}^*$ is the only character in 
$\Irr(P_{2i-2}^*)$ lying under $\alpha_{2i}^* \in \Irr(P_{2i}^*)$.
\end{lemma}

\begin{proof}
Let $i=2,3,\dots, k$ be fixed.
 We  first show that $\alpha_{2i-2}^*$
 lies under $\alpha_{2i}^*$. By \eqref{x4}, the character 
$\alpha_{2i} \in \Irr(P_{2i})$ lies over the character $\alpha_{2i-2,2i-1} \in 
\Irr(P_{2i-2,2i-1})$, where $P_{2i-2,2i-1}= C(Q_{2i-1} \tin P_{2i-2})$ 
(by \eqref{anno2})  and $\alpha_{2i-2,2i-1}$ is the $Q_{2i-1}$-Glauberman
correspondent of $\alpha_{2i-2}$ (see \eqref{pq1413def2}).
According to Theorem  \ref{dade:t2}, the $Q_{2i-1}$-correspondent 
$\alpha_{2i,2i-3}^* \in \Irr(P_{2i-2} \cdot P_{2i})$ of 
$\alpha_{2i} = \alpha_{2i,2i-1}^*$ (see  Lemma \ref{p*l1} for $j=i-1$),
 lies over the $Q_{2i-1}$-Glauberman correspondent
 $\alpha_{2i-2} \in \Irr(P_{2i-2})$ of $\alpha_{2i-2,2i-1}$.
It follows that $\alpha_{2i}^*$, which is the $Q_3,Q_5, \dots
,Q_{2i-3}$-correspondent of $\alpha_{2i,2i-3}^*$, lies over the 
$Q_3,Q_5,\dots,Q_{2i-3}$-correspondent $\alpha_{2i-2}^*$ of $\alpha_{2i-2}$
(see Theorem \ref{dade:t2}).
    
Since $P_{2i-2}^*$ is a normal subgroup of $P_{2i}^*$, Clifford's Theorem
implies that we can  prove the lemma by showing that  $\alpha_{2i-2}^*$ 
is $P_{2i}^*$-invariant. But  $P_{2i}^* = P_{2i-2}^* \cdot P_{2i}$,
and $P_{2i-2}^*$ fixes its own character, while $P_{2i}(\alpha_{2i-2}^*)
=P_{2i}$ by \eqref{pq23b}. 
We conclude that 
$\alpha_{2i-2}^*$ is $P_{2i}^*$-invariant. Thus the  lemma  holds.
\end{proof}

By induction the above lemma implies 
\begin{proposition}\mylabel{p*13}
If $1\le j\le i \le k$, then $\alpha_{2j}^*$ is the only character 
in $\Irr(P^*_{2j})$ lying under $\alpha_{2i}^* \in \Irr(P_{2i}^*)$.
\end{proposition}

\section{The groups $Q^*_{2i-1}$}\mylabel{pq:sec4}

 We can define  groups $Q^*_{2i-1}$ similar to the $P_{2i}$. 
 Indeed, 
in view of \eqref{pq14'a} the product 
\begin{equation}\mylabel{q*:e1}
Q^*_{2i-1} = Q_1\cdot Q_3 \cdot \dots \cdot Q_{2i-1},
\end{equation}
is a group whenever $1\leq i \leq l$. 

The groups $Q^*_{2i-1}$ are defined symmetrically to the $P^*_{2i}$, 
and  satisfy  results similar   to those the  $P_{2i}^*$ satisfy.
The following proposition is  analogous to Proposition \ref{pq3-} for the
 groups $Q^*_{2i-1}$.  Its proof is similar.
\begin{proposition}\mylabel{q*:p1}
The group $Q^*_{2i-1}$ is a $\pi'$-Hall subgroup of $G^*_{2i-1}$ whenever
$1\le i \le l$, while for $i=1,\dots, k$ we have, in addition, that $Q_{2i-1}^*
\in \Hall_{\pi'}(G^*_{2i})$.
Furthermore, $Q^*_{2r-1} =  Q^*_{2i-1} \cap G^*_{2r-1}$,  and thus 
$Q^*_{2r-1} \unlhd Q^*_{2i-1}$, for all $r,i$ with $1\leq r \leq i \leq l$.
\end{proposition}

\begin{proof}
The proof is done by induction, and is totally symmetric to the proof
of Proposition \ref{pq3-} for the $\pi'$-groups in the place of the $\pi$-groups. 
So we omit it.
\end{proof}

Of course  by $\pi, \pi'$-symmetry, 
 we can define irreducible characters $\beta_{2i-1}^*$ of $Q_{2i-1}^*$ 
as we did for the characters  $\alpha_{2i}^* \in \Irr(P_{2i}^*)$.
So 
\begin{defn}
For every $i=1,\dots,l$, we define $\beta_{2i-1}^* \in \Irr(Q_{2i-1}^*)$
to be the $P_2, P_4, \dots, P_{2i-2}$-correspondent of the character 
$\beta_{2i-1} \in \Irr(Q_{2i-1})$. By convention 
$$\beta_{1}^* = \beta_1.$$
\end{defn}
Results similar to Theorem \ref{p*t1}
 and Propositions \ref{pqp*t} and \ref{p*13}
hold for the $\beta^*$-characters.

According to 
Propositions \ref{pq3-} and \ref{q*:p1} we get 
\begin{corollary}\mylabel{q*:c1}
If $i=1,\dots ,k$ and $j=1,\dots,l$ then
$$G^*_{2i} = P^*_{2i}\cdot Q^*_{2i-1} \text{ and } G^*_{2j-1} = P^*_{2j-2}
\cdot Q^*_{2j-1},$$
where $P^*_0 = 1$ by convention. In particular 
$$
G^*_m = P_{2k}^*  \cdot Q_{2l-1}^*
$$
\end{corollary}

Similar to the equations \eqref{pq22}, \eqref{pq22i} and \eqref{pq22ii}
 that the groups $P^*_{2i}$ satisfy, the following equations are satisfied by the groups
$Q^*_{2i-1}$:
\begin{subequations}
\begin{equation}\mylabel{q*:e4}
N(P_{2i} \tin Q_{2i-1}\cdot Q_{2i+1} \cdot \dots \cdot Q_{2t+1}) = 
Q_{2i+1} \cdot \dots \cdot Q_{2t+1},
\end{equation}
\begin{equation}\mylabel{q*:e5}
N(P_2, P_4, \dots , P_{2i} \tin Q_{2t+1}^*) = Q_{2i+1}\cdot  \dots \cdot Q_{2t+1}
\end{equation}
and
\begin{equation}\mylabel{q*:e6}
N(P_2, \dots , P_{2i} \tin Q^*_{2i+1} ) = Q_{2i+1}, 
\end{equation}
for all $i, t $ with $1\leq i \leq t < l$.
\end{subequations}

The proof of \eqref{q*:e4} is similar to that of \eqref{pq22}, using 
\eqref{pq14b} in the place of \eqref{pq13b}.

\noindent
The equation \eqref{q*:e5} follows by repeated applications of equation \eqref{q*:e4}
 (as \eqref{pq22i} followed from \eqref{pq22}),
while  \eqref{q*:e6} is a special case (when $t=i$) of \eqref{q*:e5}.

In the proposition that follows we rewrite \eqref{pq22ii} and \eqref{q*:e6}
in a slightly different way.
\begin{proposition}\mylabel{q*:p2}
For all $i = 1, \dots ,k$ we have
\begin{equation}\mylabel{q*:e7}
N(Q^*_{2i-1} \tin P^*_{2i}) = P_{2i}.
\end{equation}
If $i=1,2,\dots,l-1$ then 
\begin{equation}\mylabel{q*:e8}
N(P^*_{2i} \tin Q^*_{2i+1}) = Q_{2i+1}.
\end{equation} 
\end{proposition} 

\begin{proof}
We use induction on $i$ to prove \eqref{q*:e7} and \eqref{q*:e8} simultaneously.
 As $Q_1^*=Q_1$ is a normal subgroup 
of $G$, it is clear that $N(Q_1^* \tin P_2^*) = P_2^* = P_2$. Hence \eqref{q*:e7} is true for $i = 1$.
Furthermore,   \eqref{q*:e6}  for $i=1$ coincides with \eqref{q*:e8} for $i = 1$.
Thus the proposition holds for $i=1$.

Assume the proposition  is true for all $i$ with $1\leq i < t $, for some
 $t=2,\dots,l-1$, (note that either $k=l$ or $k=l-1$).  We will prove it also
 holds for $i= t$, i.e., we will show that $N(Q^*_{2t-1} \tin P^*_{2t}) = P_{2t}$
and $N(P^*_{2t} \tin Q^*_{2t+1}) = Q_{2t+1}$.
  According to \eqref{pq22ii} we have that 
$$P_{2t} = N(Q_1, Q_3, \dots ,Q_{2t-1} \tin P^*_{2t}) \leq N(Q^*_{2t-1} \tin P^*_{2t}).$$
Hence the right side of \eqref{q*:e7} for $i=t$, is contained in the left side.
For the other inclusion, we observe that, by induction,
$Q_3 = N(P_2^* \tin Q^*_3), Q_5 = N(P^*_4 \tin Q^*_5 ), \dots 
,Q_{2t-1} = N(P^*_{2t-2} \tin Q^*_{2t-1})$.
In view of Propositions \ref{pq3-} and \ref{q*:p1}, 
$P^*_{2r} = P^*_{2t} \cap G^*_{2r}$ and
 $Q^*_{2r-1} = Q^*_{2t-1} \cap G^*_{2r-1}$. Hence $N(Q^*_{2t-1}
 \tin P^*_{2t})$ normalizes  both  $P^*_{2r}$ 
 (as it is a subgroup of $P^*_{2t}$) and 
$Q^*_{2r-1}$ (as it normalizes the group 
$Q^*_{2t-1}$) whenever $1\leq r\leq t$.
Hence, $N(Q^*_{2t-1} \tin P^*_{2t})$ also normalizes  the groups 
$Q_3 = N(P^*_2 \tin Q^*_3), \dots ,Q_{2t-1} = N(P^*_{2t-2} \tin Q^*_{2t-1})$.
So the left side in \eqref{q*:e7} for $i=t$  is contained in the right,  and 
the inductive step for the first equation is complete.

The proof for \eqref{q*:e8} is similar.
Indeed, according to \eqref{q*:e6} we get:
$$Q_{2t+1} = N(P_2, P_4, \dots , P_{2t} \tin Q^*_{2t+1}) \leq N(P^*_{2t} \tin Q^*_{2t+1}).$$
Thus the right side of \eqref{q*:e8} for $i=t$, is contained in the left side.

For the other inclusion,  we observe  that 
\eqref{q*:e7} holds for all $i$ with $1\leq i \leq t$.
Hence, $P_2 = N(Q_1^* \tin P^*_2), P_4 = N(Q^*_3 \tin P_4^*), \dots ,P_{2t} = 
N(Q^*_{2t-1} \tin P^*_{2t})$.
Furthermore, Propositions \ref{pq3-} and \ref{q*:p1} imply that 
$N(P^*_{2t} \tin Q^*_{2t+1})$ normalizes the groups 
$P^*_{2r}$ and $Q^*_{2r-1}$ whenever $1\leq r \leq t$.
So  it also normalizes the groups $P_2, P_4, \dots ,P_{2t}$.
Hence $N(P^*_{2t} \tin Q^*_{2t+1}) \leq N(P_2, \dots ,P_{2t} \tin Q^*_{2t+1})$.
This completes the proof of  
Proposition \ref{q*:p2} for all $i$ with $1\leq i \leq l-1$.

It remains  to show that, in the case where $k=l$, equation \eqref{q*:e7}
 holds for $i=k$. But even in this case the same argument we gave in  the inductive
proof of \eqref{q*:e7} works,  as \eqref{q*:e8} is valid for all $i=1,2,\dots,l-1=k-1$.
Hence Proposition \ref{q*:p2} holds in all cases.
\end{proof}

We close this section noticing that we have some freedom in the choice
of $P^*_{2i}$ and $Q^*_{2i-1}$, i.e.,
\begin{remark}\mylabel{pq:rcong}
As we have see in Remark \ref{ttcong}, 
conjugation  by any $g\in G^*$ leads to a new choice of the $P_{2i}$, $Q_{2i-1}$,
$\alpha_{2i}$ and $\beta_{2i-1}$ satisfying the same conditions
 (for the same $G_i$ and $\chi_i$) as the original choices. This
conjugation replaces each $P^*_{2i}$ or $Q^*_{2i-1}$
by its $g$-conjugate. In particular, we can choose $g \in G^*$ so that 
$(P^*_{2k})^g$ and $(Q^*_{2l-1})^g$ are any given $\pi$-Hall  
and $\pi'$-Hall  subgroups,
 respectively, of $G^*_{m}$. 
\end{remark}

\newpage
\section{When \textit{$\pi$-split} groups are involved}\mylabel{split}

In this section we are interested in the special case where $\pi$-split groups 
appear in the normal series \eqref{*}.
 In particular, we will  see that  the 
triangular sets, in this case,  have a very simple  form.
What we mean for a group to be  $\pi$-split is given in 
\begin{defn}\mylabel{sp.d1} 
A finite group $H$ is called $\pi$-split if it is the direct product 
$$
H = H_{\pi} \times H_{\pi'}
$$
of a $\pi$-group $H_{\pi}$,  and  a $\pi'$-group $H_{\pi'}$.
\end{defn}
Obviously,  $H_{\pi}$ and $H_{\pi'}$ are the unique $\pi$-and $\pi'$-Hall 
subgroups of $H$. Also,  if $S$ is any subgroup of $H$, then $S$ is also
 $\pi$-split. 
Furthermore, if $\chi \in \Irr(H)$ is an irreducible character of $H$, then 
$\chi$ also $\pi$-splits as 
 $$\chi= \chi_{\pi} \times \chi_{\pi'},
$$ 
where $\chi_{\pi} \in \Irr(H_{\pi}) $ and $\chi_{\pi'}\in \Irr(H_{\pi'})$ are the $\pi$-and $\pi'$-parts 
to which $\chi$ decomposes.

Assume now that the normal series \eqref{*}, in addition to its usual
 properties described in Hypothesis \ref{hyp1}, contains some 
 $\pi$-split group $G_i$,  for some $i=1, \dots,m$.  
Clearly, if $G_i$ is $\pi$-split then $G_j$ is also $\pi$-split 
for all $j=1,\dots, i$.
Let $s$ be the largest integer, with $1\leq s \leq m$, such that $G_s$ 
is $\pi$-split. Note that $s$ is necessarily bigger than
 $0$ as $G_1$ is a 
$\pi'$-group and thus a $\pi$-split group. Let 
\begin{equation}\mylabel{sp.e0}
\{\chi_0, \chi_1, \dots, \chi_m \}
\end{equation}
 be a fixed but arbitrary character 
tower for \eqref{*}. Then  
\begin{subequations}\mylabel{sp.e1}
\begin{align}
G_i = G_{i, \pi} &\times G_{i, \pi'}, \\ 
\chi_i = \chi_{i, \pi} &\times \chi_{i, \pi'},
\end{align}
\end{subequations}
whenever $0\leq i \leq s$. Furthermore, the groups $G^*_i$ 
and $G^*$ defined in 
Section \ref{pq:sec3} (see \eqref{pqinf} and \eqref{pq1}), and their
 characters $\chi_i^*$ (see \eqref{pq2a}),
satisfy
\begin{subequations}\mylabel{sp.e2}
\begin{align}
G_i^*=G_{i,\pi}^* &\times G_{i, \pi'}^*= G_{i, \pi}(\chi_1,\dots,\chi_{i-1})
\times G_{i, \pi'}(\chi_1, \dots, \chi_{i-1}), \\
\chi_{i}^*= \chi_{i, \pi}^* &\times \chi_{i, \pi'}^*,
\end{align}
\end{subequations}
whenever $0\leq i \leq s$.

Let 
\begin{equation}\mylabel{sp.e3}
\{P_{2r}, Q_{2i-1}| \alpha_{2r}, \beta_{2i-1}\}_{r=0, i=1}^{k, l},
\end{equation}
be a representative of the conjugacy class of triangular sets of \eqref{*} 
that corresponds to the tower \eqref{sp.e0}, by Theorem \ref{cor:t}.
All the groups,  their characters,  and their properties, 
that were introduced and proved in the previous 
sections with respect to a given triangular set, (like $Q_{2i-1,2r},
 P_{2r, 2i-1}, P_{2r}^*$ etc.),  
are applied to the set \eqref{sp.e3}.
Furthermore, we  write
\begin{equation}\mylabel{sp.e4}
l_s:=[(s+1)/2] \text{ and } k_s:=[s/2],
\end{equation}
for the greatest integers less than $(s+1)/2$ and $s/2$, respectively.
(This agrees with the definition that was given in \eqref{kl:def}.

In the situations where \eqref{sp.e1} occurs, the first $n$ groups 
in the  triangular set  \eqref{sp.e3} are unique and satisfy
 \begin{theorem}\mylabel{sp.t1}
Assume that \eqref{sp.e1} holds. Then 
\begin{subequations}\mylabel{sp.e5}
\begin{align}
P_{2r} &=P_{2r}^*= G_{2r, \pi}^*=G_{2r, \pi}(\chi_1,\dots,\chi_{2r-1}),\\
Q_{2i-1} &= Q_{2i-1}^*=G_{2i-1,\pi'}^*=G_{2i-1,\pi}(\chi_1,\dots,\chi_{2i-2}),
\\
G^*_{2r}&=P_{2r} \times Q_{2r-1}=G_{2r, 2r-1} = G_{2r, 2r},\\
G^*_{2i-1}&= P_{2i-2} \times Q_{2i-1} =G_{2i-1, 2i-2} =G_{2i-1, 2i-1},\\
\alpha_{2r}&= \alpha_{2r}^*=\chi_{2r, \pi}^*,\\
\beta_{2i-1} &=\chi_{2i-1,\pi'}^*, \\
\chi_{2r}^*&=\chi_{2r, 2r}= 
\alpha_{2r} \times \beta_{2r-1}, \\
\chi_{2i-1}^* &= \chi_{2i-1,2i-1}=
\alpha_{2i-2} \times \beta_{2i-1}.
\end{align}
\end{subequations}
whenever $1\leq r \leq k_s$ and $1\leq i \leq l_s$.

Furthermore, the groups  $P_2 \unlhd P_4 \unlhd \dots \unlhd P_{2k_s}$ and 
$Q_1 \unlhd Q_3 \unlhd \dots \unlhd Q_{2l_s-1}$, form a normal series for
$P_{2k_s}$ and $Q_{2l_s-1}$, respectively. In addition,  
$P_{2k_s}$ centralizes $Q_{2l_s-1}$. Thus
\begin{align}\mylabel{sp.e6} 
P_{2r, 2t+1} = P_{2r} &\text{ and } Q_{2i-1,2j} =Q_{2i-1},\\
\alpha_{2r, 2t+1}=\alpha_{2r}  &\text{ and } \beta_{2i-1,2j}= \beta_{2i-1},
\end{align}
whenever $1\leq r \leq t \leq l_s-1$ and $1\leq i \leq j \leq k_s$.
\end{theorem}

\begin{proof}
Corollary \ref{q*:c1}, along with (\ref{sp.e2}a), implies that 
$G^*_{2r}=G_{2r, \pi}^* \times G_{2r, \pi'}^*=
P_{2r}^* \cdot  Q_{2r-1}^* $ and 
$G^*_{2i-1}=G_{2i-1, \pi}^*\times G_{2i-1, \pi'}^* =
 P_{2i-2}^* \cdot  Q_{2i-1}^*$, 
whenever $1\leq r \leq k_s$ and $1\leq i \leq l_s$.
Note that the last such group is $G_s^*$,  wich satisfies
\begin{equation}\mylabel{sp.e6.5}
G_s^*= G_{s, \pi}^* \times G_{s, \pi'}^*= P_{2k_s}^* \cdot  Q_{2l_s-1}^*.
\end{equation}
Hence 
\begin{gather}\mylabel{sp.e7}
P_{2r}^*=G_{2r, \pi}(\chi_1,\dots,\chi_{2r-1})=
G_{2r, \pi}^*= G_{2r+1, \pi}^*, \\
Q_{2i-1}^*= G_{2i-1,\pi'}(\chi_1,\dots,\chi_{2i-2})=
G_{2i-1, \pi'}^*= G_{2i, \pi'}^*, \text{while}\\ 
P_{2k_n}^*= G_{2k_s, \pi}^* \text{ and } Q_{2l_s-1}^* = G_{2l_s-1, \pi'}^*,
\end{gather}
 for   all  $r, i$ with $1\leq r < l_s$ and $1\leq i \leq k_s$. 
This, along with  Proposition \ref{q*:p2}, implies 
\begin{gather*}
P_{2r}= N(Q_{2r-1}^* \tin P_{2r}^*) = 
N(G_{2r, \pi'}^* \tin G_{2r, \pi}^*)= G_{2r, \pi'}^* = P_{2r}^*,\\
Q_{2i-1}= N(P_{2i-2}^* \tin Q_{2i-1}^*)= 
N(G_{2i-1, \pi}^* \tin G_{2i-1, \pi'}^*) = G_{2i-1, \pi'}^*= Q_{2i-1}^*, 
\end{gather*}
whenever $1\leq r \leq k_s$ and $1\leq i \leq l_s$.
Hence (\ref{sp.e5}a,b) holds.

Furthermore, the fact that the groups $P_2^*, P_4^*, \dots, P_{2k_s}^*$ form a
normal series for $P_{2k_s}^*$ (see Proposition \ref{pq3-}), implies that 
$P_2 \unlhd P_4 \unlhd \dots \unlhd P_{2k_s}$ is a normal series for
 $P_{2k_s}$. Similarly,  the $\pi'$-groups $Q_{2i-1}$ form 
 a normal series $Q_1 \unlhd Q_3 \unlhd \dots \unlhd Q_{2l_s-1}$ for
 $Q_{2l_s -1}$.
According to \eqref{sp.e6.5} the $\pi$-Hall subgroup
$P_{2k_s}$ of $G_s^*$ centralizes the $\pi'$-Hall subgroup
 $Q_{2l_s-1}$ of that  same group.
Hence $P_{2r}$ centralizes $Q_{2i-1}$,  for all $r=1,\dots,k_s$ and $i=1,
\dots, l_s$. Thus (see \eqref{pq14b} and \eqref{pq13b}),
$Q_{2i-1,2j}=C(P_{2i}, \dots, P_{2j} \tin Q_{2i-1})= Q_{2i-1}$ and
$P_{2r, 2t+1}= C(Q_{2r+1}, \dots, Q_{2t+1} \tin P_{2r}) =P_{2r}$, for all
$1\leq i \leq j \leq k_s$ and $1\leq r \leq t \leq l_s-1$.
Furthermore, the $P_{2i}, \dots, P_{2j}$-Clifford correspondent 
$\beta_{2i-1, 2j}$ of $\beta_{2i-1}$ coincides with $\beta_{2i-1}$.
Similarly we get that   $\alpha_{2r, 2t+1}=\alpha_{2r}$.
Hence the last part of the theorem holds.

To prove (\ref{sp.e5}c,d) it suffices to notice that, 
according to \eqref{tow--tri1} and \eqref{tow--tri2},
\begin{align*}
G_{2r, 2r-1}&= P_{2r} \ltimes Q_{2r-1},\\
G_{2r, 2r} &= P_{2r} \times Q_{2r-1,2r},\\
G_{2i-1, 2i-2}&=P_{2i-2} \rtimes Q_{2i-1},\\
G_{2i-1,2i-1}&= P_{2i-2,2i-1} \times Q_{2i-1}, 
\end{align*}
for all $r=1,\dots, k_s$ and $i=1,\dots, l_s$.
But $Q_{2r-1,2r}= Q_{2r-1}$ and $P_{2i-2,2i-1}= P_{2i-2}$, by \eqref{sp.e6}.
Thus all the above products are direct, and we get 
$$
G_{2r, 2r-1}=G_{2r, 2r}=P_{2r} \times Q_{2r-1}= P_{2r}^* \times Q_{2r-1}^*
= G_{2r}^*.
$$
So  (\ref{sp.e5}c) holds. The proof for (\ref{sp.e5}d) is analogous.

Notice that (\ref{sp.e5}c,d) clearly imply that 
the groups $P_{2r}$ and $Q_{2i-1}$, as characteristic subgroups of $G_{2r}^*$
 and $G_{2i-1}^*$, respectively, are normal subgroups of 
$G^*=G(\chi_1,\dots, \chi_m)$, for all $r=1,\dots, k_s$ and 
$i=1,\dots, l_s$.
This  is actually the reason that the group $G_{2r, 2r}$,  defined as 
$G_{2r, 2r}=N(P_0, \dots, P_{2r}, Q_1, \dots, Q_{2r-1} \tin 
G_{2r}(\chi_1,\dots, \chi_{2r}))$ in  \eqref{tow--tri3}, coincides with
 $G_{2r}(\chi_1,\dots,\chi_{2r})=G_{2r}^*$, for all $r=1,\dots, k_s$.
(Similarly we work  for the group  $G_{2i-1, 2i-1}$).
 In addition,  this implies that the $cP_2, \dots, cP_{2r}, cQ_1,
\dots, cQ_{2r-1}$-correspondent  $\chi_{2r, 2r}$  of $\chi_{2r}$ 
is nothing else but a multiple Clifford correspondent, and  
thus coincides with $\chi_{2r}^*$, i.e., $\chi_{2r, 2r}= \chi_{2r}^*$.
Similarly, $\chi_{2i-1,2i-1}= \chi_{2i-1}^*$.
But, according to  \eqref{tow--tri1} and \eqref{tow--tri2}, we have that 
$\chi_{2r, 2r}= \alpha_{2r} \times \beta_{2r-1,2r}$ and $\chi_{2i-1,2i-1}=
\alpha_{2i-2,2i-1} \times \beta_{2i-1}$.
Hence 
\begin{align*}
\chi_{2r}^*=\chi_{2r, 2r}= \alpha_{2r} \times \beta_{2r-1,2r}=
\alpha_{2r}\times \beta_{2r-1}, \\
\chi_{2i-1}^*= \chi_{2i-1, 2i-1}=\alpha_{2i-2, 2i-1} \times \beta_{2i-1}=
\alpha_{2i-2} \times \beta_{2i-1},
\end{align*}
whenever $1\leq r \leq k_s$ and $1\leq i \leq l_s$.
This, along with (\ref{sp.e2}b), implies 
(\ref{sp.e5}f,g,h) and one equality in (\ref{sp.e5}e), namely
$\alpha_{2r} = \chi_{2r, \pi}^*$.

It remains to show that  $\alpha_{2r}^* =\alpha_{2r}$.
But $\alpha_{2r}^* \in \Irr(P_{2r}^*)$
 is the $Q_3, \dots, Q_{2r-1}$-correspondent 
of $\alpha_{2r} \in \Irr(P_{2r})$, while $P_{2r}^*=P_{2r}$ centralizes 
the $\pi'$-groups  $Q_3, \dots, Q_{2r-1}$,  for all $r=1,\dots, k_s$.
Thus this correspondence is trivial, i.e., 
$\alpha_{2r}^*= \alpha_{2r}$, for all such $r$.
This completes the proof of the theorem.
\end{proof}

The following is a straight forward application of Theorem \ref{sp.t1}.

\begin{corollary}\mylabel{sp.co1}
Assume that $G_{i}$ and $\chi_i$ satisfy \eqref{sp.e1}, for all $i=1,\dots, s$.
In addition, assume that $G$ fixes $\chi_i$ for all $i=1,\dots, s-1$. Then 
the triangular set \eqref{sp.e3} satisfies
\begin{align}\mylabel{sp.e8}
P_{2r}&=P_{2r}^*=G_{2r,\pi},\\
Q_{2i-1} &=Q_{2i-1}^*=G_{2i, \pi'},\\
\alpha_{2r}&=\chi_{2r, \pi},\\
\beta_{2i-1}&=\chi_{2i-1,\pi'},\\
\chi_{2r}&= \alpha_{2r} \times \beta_{2r-1},\\
\chi_{2i-1} &=\alpha_{2i-2} \times \beta_{2i-1}
\end{align}
whenever $1\leq r\leq k_s$ and $1\leq i \leq l_s$. 
\end{corollary}

%%% Local Variables: 
%%% mode: latex
%%% TeX-master: t
%%% End: 

%%% Local Variables: 
%%% mode: latex
%%% TeX-master: "thesis-ex"
%%% End: 

%%\include{split}

\chapter{ The  Group $G'=G(\alpha_{2k}^*)$ } \mylabel{pq:sec5}
Assume that a finite odd group $G$ is given, along with a normal
series \eqref{*} and a fixed triangular set for this series.
We have already seen in Section \ref{pq:sec3} how to define the characters 
$\alpha_{2i}^*$ of the product groups $P_{2i}^*$ whenever $1\leq i \leq k$. 
In this chapter we analyze the group $G':=G(\alpha_{2k}^*)$ with ultimate 
purpose to reach triangular sets for this group.

\section{$\pi'$-Hall subgroups of $G'$: the group $\qw$ }

The following remark is  an easy consequence of Proposition 
\ref{p*13} and the fact that $G'$ normalizes $P_{2i}^* = P_{2k}^* \cap G_{2i}$.
\begin{remark}\mylabel{b.rem1}
$G'= G(\alpha_{2k}^*)= G(\alpha_{2}^*,\alpha_{4}^*,\dots,\alpha^*_{2k})$
\end{remark}
 
\begin{proposition}\mylabel{b.pro1}
For every $i=1,\dots,k$ we have 
\begin{subequations}
\begin{align}
G'(\beta_1,\dots,\beta_{2i-1})&\leq G'(\beta_{2i-1,2k}) \text{\, and }
\mylabel{b.e1}\\
G'(\beta_1,\dots,\beta_{2i-1})&\leq G'(\chi_1,\dots,\chi_{2i})\leq
G'(\chi_1,\dots,\chi_{2i-1}) \mylabel{b.e2.5}.
\end{align}
In addition
\begin{equation}\mylabel{b.e2}
 G'(\beta_1,\dots,\beta_{2l-1})\leq G'(\chi_1,\dots,\chi_{2l-1}).
\end{equation}
\end{subequations}
\end{proposition}

\begin{proof}
Let $T_i=G'(\beta_1,\dots,\beta_{2i-1})$ for some fixed $i \in \{1,\dots,k\}$.
In view of Remark \ref{b.rem1} we have that 
$T_i$ fixes $\alpha_2^*,\dots,\alpha_{2k}^*$,   and thus normalizes the groups
$P_2^*, P_{4}^*, \dots,P_{2k}^*$. Furthermore, it normalizes the 
$\pi'$-groups $Q_1,Q_3,\dots,Q_{2i-1}$ and therefore also normalizes the 
product group
$P_{2i}\cdots P_{2k}=N(Q_1,Q_3,\dots Q_{2i-1} \tin P_{2k}^*)$
 (see \eqref{pq22i}).
As $T_i$ fixes $\beta_{2i-1} \in \Irr(Q_{2i-1})$ and normalizes 
$P_{2i}\cdots P_{2k}$,  it also fixes the $P_{2i}\cdots P_{2k}$-Glauberman 
correspondent $\beta_{2i-1,2k} \in \Irr(Q_{2i-1,2k})$ of $\beta_{2i-1}$
(see Definition \ref{pq1413def2}).
This implies \eqref{b.e1}.

We will use induction on $i$ to prove \eqref{b.e2.5}.
If $i=1$, then obviously $T_1=G'(\beta_1)= G'(\chi_1)$, as $\beta_1= \chi_1$.
Also, $T_1$ normalizes $P_{2}= P_2^*$ and fixes $\alpha_2 = \alpha_2^*$.
Hence it also normalizes $P_2 \cdot Q_1 = G_2(\beta_1)$ (see \eqref{pq3b}),
and thus fixes the canonical extension $\beta_1^e$ of $\beta_1$ to $P_2 Q_1$.
Therefore it fixes the character $\chi_{2,1}= \alpha_2 \cdot \beta_1^e$ (see \eqref{pq3c}).
According to \eqref{pqn0} we have $T_1(\chi_1,\chi_2)= T_1(\chi_{1,1},\chi_{2,1})$.
But $T_1(\chi_1)= T_1 =T_1(\chi_{2,1})$, while $\chi_{1,1} = \chi_1$. 
Hence $T_1(\chi_2)=T_1(\chi_1,\chi_2) =T_1(\chi_{1,1},\chi_{2,1})=T_1$.
Therefore, $T_1$ fixes $\chi_2$ and \eqref{b.e2.5} is proved for $i=1$.

Assume \eqref{b.e2.5} is true for $i=t-1$ and  some $t=2,3,\dots,k$.
 We will prove it  
also holds for $i=t$.
 The inductive hypothesis implies   that 
$$
T_t=G'(\beta_1,\beta_3,\dots,\beta_{2t-1}) \le 
 G'(\beta_1,\beta_3,\dots,\beta_{2t-3})\le 
G'(\chi_1,\chi_3,\dots,\chi_{2t-2}).
$$
Thus it suffices to show that $T_t$ fixes $\chi_{2t-1}$ and $\chi_{2t}$.
By \eqref{tow--tri1} and \eqref{tow--tri2} we have 
\begin{align*}
\chi_{2t-1,2t-1}&=\alpha_{2t-2,2t-1} \times \beta_{2t-1}\\
\chi_{2t,2t-1}&=\alpha_{2t} \cdot  \beta_{2t-1}^e.
\end{align*}
We have already seen that $T_t$ normalizes $P_{2t-2}^*$ and  $P^*_{2t}$, and fixes their  
characters $\alpha_{2t-2}^*$ and $\alpha_{2t}^*$. Also it  normalizes $Q_1,\dots,Q_{2t-3},Q_{2t-1}$,
  and thus normalizes $P_{2t-2}=N(Q_1,\dots,Q_{2t-3} \tin P_{2t-2}^*)$ and 
$P_{2t} = N(Q_1,\dots,Q_{2t-1} \tin P_{2t}^*)$,
(see \eqref{pq22ii}). Hence $T_t$ also fixes the $Q_1,Q_3,\dots,Q_{2t-3}$-correspondent
 $\alpha_{2t-2} \in \Irr(P_{2t-2})$ of $\alpha_{2t-2}^*$, as well as the $Q_1,\dots,
Q_{2t-1}$-correspondent $\alpha_{2t} \in \Irr(P_{2t})$ of $\alpha_{2t}^*$
(see Proposition \ref{pqp*t}). As $T_t$ also normalizes $Q_{2t-1}$, it fixes
the $Q_{2t-1}$-Glauberman correspondent $\alpha_{2t-2,2t-1} \in 
\Irr(P_{2t-2,2t-1})$ of $\alpha_{2t-2} \in \Irr(P_{2t-2})$.
So $T_t$ fixes the characters $\alpha_{2t},\alpha_{2t-2}$ and $\alpha_{2t-2,2t-1}$.

Also $T_t$ fixes $\beta_{2t-1}$ and  normalizes $P_{2t} \cdot Q_{2t-1} = G_{2t,2t-1}$ (see \eqref{tow--tri1}).
Hence it  fixes the canonical extension $\beta_{2t-1}^e \in \Irr(G_{2t,2t-1})$ of $\beta_{2t-1}$
to $G_{2t,2t-1}$. 
Therefore,  $T_t$ fixes $\alpha_{2t-2,2t-1},\beta_{2t-1},  \alpha_{2t}$ and $\beta_{2t-1}^e$, 
and thus fixes $\chi_{2t-1,2t-1}$ and $\chi_{2t,2t-1}$. This, along with 
the  inductive hypothesis on $T_t$,  implies that  
\begin{equation}\mylabel{b.e3}
T_t(\chi_1,\chi_3,\dots,\chi_{2t-2},\chi_{2t-1,2t-1},\chi_{2t,2t-1}) = T_t.
\end{equation}

We note that $T_t$  normalizes all the $\pi$-groups 
$P_2,P_4,\dots,P_{2t-4},P_{2t-2}$. This is clear,  as for every $j=1,\dots,t-1$ we have 
  $P_{2j}= N(Q_1,\dots,Q_{2j-1} \tin P_{2j}^*)$ (by \eqref{pq22ii}).
Hence we conclude that   $T_t \le N(Q_1,Q_3,\dots,Q_{2t-1},P_2,P_4,\dots,P_{2t-2} \tin G)$. 
Therefore Theorem \ref{tow--tri} (part 3 for $n=2t$ and $n=2t-2$ respectively) implies that  
$T_t(\chi_1,\dots,\chi_{2t}) = T_t(\chi_{1,2t-1},\dots,\chi_{2t,2t-1})$ and 
$T_t(\chi_1,\dots,\chi_{2t-2}) = T_t(\chi_{1,2t-1},\dots,\chi_{2t-2,2t-1})$. 
 Hence 
\begin{align*}
T_t(\chi_1,\dots,\chi_{2t})&=T_t(\chi_{1,2t-1},\dots,\chi_{2t,2t-1})& \\
&=T_t(\chi_{1,2t-1},\dots,\chi_{2t-2,2t-1})(\chi_{2t-1,2t-1},\chi_{2t,2t-1})& \\
&=T_t(\chi_1,\dots,\chi_{2t-2})(\chi_{2t-1,2t-1},\chi_{2t,2t-1})& \\
&=T_t. &\text{ by \eqref{b.e3}}
\end{align*}
So $T_t \le G'(\chi_1,\dots,\chi_{2t})$. This proves the inductive step for 
$i=t$,  and thus \eqref{b.e2.5} for every $i=1,\dots,k$.

It remains to show \eqref{b.e2}.
Observe that this additional case  is not covered by \eqref{b.e2.5}
 only when $m$ is odd, 
since for $m$ even we have $k=l$.
The arguments for this last step are similar to those we used at the inductive step.
Indeed, $T_l$ fixes $\alpha_{2l-2}^*$ and thus fixes its $Q_1,\dots,Q_{2l-3}$-correspondent 
$\alpha_{2l-2} \in \Irr(P_{2l-2})$. It also fixes the $Q_{2l-1}$-Glauberman correspondent 
$\alpha_{2l-2,2l-1}$ of $\alpha_{2l-2}$. Hence $T_l$ fixes
 $\alpha_{2l-2,2l-1}\times \beta_{2l-1} = \chi_{2l-1,2l-1}$ (see \eqref{tow--tri2}).
Therefore,
$T_l=T_l(\chi_1,\dots,\chi_{2k},\chi_{2l-1,2l-1})$. 
Furthermore, $T_l \leq N(Q_1,\dots,Q_{2l-1},P_2,\dots,P_{2l-2} \tin G )$.
 Therefore, 
Part 3 of Theorem \ref{tow--tri} implies that 
\begin{align*}
T_l(\chi_1,\dots,\chi_{2l-1})&=T_l(\chi_{1,2l-1},\dots,\chi_{2l-1,2l-1})& \\
&=T_l(\chi_{1,2l-1},\dots,\chi_{2l-2,2l-1})(\chi_{2l-1,2l-1})& \\
&=T_l(\chi_1,\dots,\chi_{2l-2})(\chi_{2l-1,2l-1})& \\
&=T_l. &
\end{align*}
So $T_l \le G'(\chi_1,\dots,\chi_{2l-1})$. This proves \eqref{b.e2}
and completes the proof of the proposition. 
\end{proof}

Our next goal is to 
show that the smallest group in \eqref{b.e2.5} has $\pi$-index in the largest 
one. The following lemma helps in this direction.

\begin{lemma}\mylabel{b.lem0}
If \,  $T \leq N(Q_1,\dots, Q_{2i-1} \tin  G')$,  for some $i$ with $1\leq i \leq l$,
then 
$$T(\chi_1,\dots,\chi_{2i-1})  \leq T(\beta_1,\dots,\beta_{2i-1}).
$$
\end{lemma}

\begin{proof}
We will use induction on $i$. 
If $i=1$,  then the lemma is obviously true,  
as $\chi_1 = \beta_1$, so that $T(\chi_1) = T(\beta_1)$ for any $T \leq G$.
Assume that the lemma holds  for all $i=1, \dots, t-1$,  and some
 $t =2,3,\dots,l$. 
We will prove it also holds for $i=t$.

Let $T$ be a subgroup of $N(Q_1,\dots, Q_{2t-1} \tin G')$. 
Then,  according to the inductive hypothesis, 
$T(\chi_1, \dots,\chi_{2t-1})\leq T(\chi_1,\dots, \chi_{2t-3})
\leq T(\beta_{1},\dots,\beta_{2t-3})$.
 Furthermore,
in view of Remark \ref{b.rem1}, the group $T$ fixes the characters 
 $\alpha_2^*,\alpha_4^*, \dots, \alpha_{2k}^*$ 
and normalizes the groups $P_{2}^*, P_4^*,\dots,P_{2k}^*$.
Hence  $T$ normalizes the groups 
$P_{2i} = N(Q_1,\dots,Q_{2i-1}  \tin P_{2i}^*)$ (see \eqref{pq22ii}), 
as well as the product groups $P_{2i}\cdots P_{2k} = N(Q_1,\dots,Q_{2i-1} 
\tin P_{2k}^*)$ (see \eqref{pq22i}),  whenever $1\leq i \leq t$.
So we get that 
$T \leq N(Q_1,\dots,Q_{2t-1}, P_2,\dots,P_{2t-2} 
\tin G)$, which, in view of Theorem \ref{tow--tri} (Part 3), implies that  
$$
T(\chi_1,\dots,\chi_{2t-1}) = T(\chi_{1,2t-1}, \dots,\chi_{2t-1,2t-1}).
$$
Hence $T(\chi_1, \dots,\chi_{2t-1})$ fixes $\chi_{2t-1,2t-1}$.
But the last  character  equals $\alpha_{2t-2,2t-1} \times \beta_{2t-1}$,
(see \eqref{tow--tri2}).
Hence $T(\chi_1, \dots,\chi_{2t-1})$ fixes $\beta_{2t-1}$.
 Therefore  $T(\chi_1, \dots,\chi_{2t-1}) \leq T(\beta_1,\dots,\beta_{2t-3},
\beta_{2t-1})$.

This completes the proof of the inductive argument,  and thus that  of 
Lemma \ref{b.lem0}.
\end{proof}

\begin{theorem}\mylabel{b.th1}
There exists a $\pi'$-Hall subgroup $\maq$ of $G'$ such that,  
for every $i=1,\dots,l$,  we have 
\begin{subequations}
\begin{align}
\maq(\chi_1,\dots,\chi_{2i-1}) &\in \Hall_{\pi'}(G'(\chi_1,\dots,\chi_{2i-1}))
 \text{\, and } \mylabel{b.te1}\\
\maq(\chi_1,\dots,\chi_{2i-1}) &\leq \maq(\beta_1,\dots,\beta_{2i-1}). 
 \mylabel{b.te2} 
\end{align}
Furthermore, whenever $1\leq i \leq l-1$ we have 
\begin{equation} \mylabel{b.te3}
\maq(\chi_{1},\dots,\chi_{2i-1}) \text{ normalizes \, }  Q_{2i+1}. 
\end{equation}
\end{subequations}
\end{theorem}

\begin{proof}
As 
$$
G'(\chi_1, \dots,\chi_{2l-1}) \leq G'(\chi_1,\dots,\chi_{2l-3})\leq 
\dots \leq G'(\chi_1) \leq G', 
$$
it is obvious that we can pick a $\pi'$-Hall subgroup $\maq$ of $G'$ 
that satisfies \eqref{b.te1} for all $i=1,\dots,l$.
We will modify $\maq$, using induction on  $i$,  so that the rest of the theorem also holds.

If $i=1$,  then we obviously have that 
 \eqref{b.te2} holds,   as $\chi_1 = \beta_1$.
Thus it suffices to show that we can modify  $\maq$ so that 
 \eqref{b.te1} holds for all $i=1,\dots, l$  while \eqref{b.te3} holds for $i=1$.
According to Remark \ref{b.rem1} and \eqref{p*a},
 the group $G'$ fixes $\alpha_{2}^* = \alpha_2$.
Hence $G'(\beta_{1}) \leq G(\alpha_2,\beta_1)$.
Therefore,  $G'(\beta_{1})$ normalizes
$G_{3}(\alpha_{2},\beta_1)$,  and thus $\maq(\beta_1)$ normalizes 
$G_3(\alpha_2,\beta_1)$.
 According to \eqref{tow--tri2}, the  semidirect product 
$Q_3 \ltimes P_2$ equals   $G_3(\alpha_2,\beta_1)$.
As  the $\pi'$-group $\maq(\beta_{1})$  normalizes $Q_3 \ltimes P_2$,
it has to normalize a $P_2$-conjugate of $Q_3$. 
Hence $\maq(\beta_{1})^{\sigma_2}$ normalizes $Q_3$, for 
some $\sigma_2 \in P_2$. 
But $P_2$ fixes $\alpha_{2k}^*$ (as $P_2 \leq P_{2k}^*$) as well as 
 $\beta_{1} = \chi_1$. It also fixes the characters
$\chi_2,\dots,\chi_{2l-1}$ as $P_2 \leq G_2 \leq \dots \leq G_{2l-1}$.
 Therefore, 
$P_2 \leq G'(\chi_1,\dots,\chi_{2i-1})$ whenever $1\leq i \leq l$.
 Hence  $\maq(\chi_1,\dots,\chi_{2i-1})^{\sigma_2} =
 \maq^{\sigma_2}(\chi_1,\dots,\chi_{2i-1} )$ and, in addition, 
  $\maq^{\sigma_2}$  and $\maq^{\sigma_2}(\chi_1,\dots,\chi_{2i-1} )$ 
 are  $\pi'$-Hall subgroups of $G'$ and $G'(\chi_1,\dots,\chi_{2i-1})$ 
respectively  (as $\maq$ and $\maq(\chi_1,\dots,\chi_{2i-1})$ are, and $P_2 
\leq G'(\chi_1,\dots,\chi_{2i-1})$). 
Furthermore $\maq^{\sigma_2}(\chi_1)$ normalizes $Q_3$. 
 Hence the group $\maq^{\sigma_2}$ satisfies \eqref{b.te1} for every 
$i=1,\dots,l$ as well as \eqref{b.te2} and \eqref{b.te3} for $i=1$.
So we can replace $\maq$ by $\maq^{\sigma_2}$ and assume that 
  \eqref{b.te1} holds  for every 
$i=1,\dots,l$ while  \eqref{b.te2} and \eqref{b.te3} hold  for $i=1$

The same type of argument as the one we gave for $i=1$ will make the 
inductive step work.
So, assume that $\maq$ has been modified so that it satisfies 
 \eqref{b.te1}  for all $i=1,\dots,l$,  
and in addition,  satisfies  the rest of the  theorem for all $i
\leq t-1$,  for some $t=2,\dots,l-1$.
We will show that there is a $G'$-conjugate of $\maq$
that satisfies \eqref{b.te1} for all $i=1,\dots,l$ and the rest of 
 Theorem \ref{b.th1} whenever $1\leq i \leq t$.

According to the inductive hypothesis
$\maq(\chi_1, \dots,\chi_{2t-1})\leq \maq(\chi_1,\dots, \chi_{2t-3})
\leq \dots \leq \maq(\chi_{1})$ while, 
$\maq(\chi_1,\dots, \chi_{2i-1})$   normalizes the group $Q_{2i+1}$
for all $i=1,\dots,t-1$.
 So  $\maq(\chi_1,\dots,\chi_{2t-1})$ normalizes the groups 
 $Q_1,\dots,Q_{2t-3}, Q_{2t-1}$. Hence Lemma \ref{b.lem0} implies that 
 $\maq(\chi_1, \dots,\chi_{2t-1}) \leq \maq(\beta_1,\dots,\beta_{2t-3},
\beta_{2t-1})$.
This, along with the inductive hypothesis, implies that 
\begin{equation}\mylabel{b.e5}
\maq(\chi_1, \dots,\chi_{2i-1}) \leq \maq(\beta_1,\dots,\beta_{2i-1}),
\end{equation}
whenever $1\leq i \leq t$.

The group $\maq(\chi_1, \dots,\chi_{2t-1})$ fixes the characters
$\alpha_{2i}^*$ and normalizes the groups $Q_{2i-1}$ and $P_{2i}$
for all $i=1, \dots,t$. Hence Proposition \ref{pqp*t} implies that 
$\maq(\chi_1, \dots,\chi_{2t-1})$ also fixes the
$Q_3,\dots,Q_{2i-1}$-correspondent 
$\alpha_{2i}\in \Irr(P_{2i})$ of $\alpha_{2i}^* \in \Irr(P_{2i}^*)$
for all such $i$.
This, along with \eqref{b.e5},  implies that 
$\maq(\chi_1, \dots,\chi_{2t-1}) \leq \maq(\beta_{1},\dots,\beta_{2t-1}, 
\alpha_{2}, \dots,\alpha_{2t})$.
Hence $\maq(\chi_1, \dots,\chi_{2t-1})$ normalizes the group 
$G_{2t+1}(\beta_1,\dots,\beta_{2t-1},\alpha_2,\dots,\alpha_{2t})$, (as
$G_{2t+1} \unlhd G$).
According to \eqref{tow--tri2} and \eqref{tow--tri3} the latter group
equals $G_{2t+1,2t} = P_{2t} \rtimes Q_{2t+1}$.
Hence the $\pi'$-group  $\maq(\chi_1,\dots,\chi_{2t-1})$ normalizes $P_{2t} 
\rtimes Q_{2t+1}$,  and thus normalizes a $P_{2t}$-conjugate of $Q_{2t+1}$.
Thus there exists an element  $\sigma \in P_{2t}$ such that 
$\maq(\chi_1,\dots,\chi_{2t-1})^{\sigma}$ normalizes $Q_{2t+1}$.
But $P_{2t}$ is a subgroup of $G_{2t+1,2t}$,  where the latter group
equals  $ N(Q_1,\dots,Q_{2t-1},P_{2},\dots,P_{2t} \tin 
G_{2t+1}(\chi_1,\dots,\chi_{2t}))$ (see \eqref{tow--tri3}).
Therefore, $P_{2t}$ fixes the characters $\chi_1, \dots, \chi_{2t}$.
Furthermore, $P_{2t} \leq G_{2t+1} \leq G_{2t+2} \leq \dots \leq G_{2l-1}$ 
which implies that $P_{2t}$ also fixes the characters
$\chi_{2t+1},\chi_{2t+2},\dots,\chi_{2l-1}$. Hence 
$\maq^{\sigma}(\chi_1,\dots,\chi_{2i-1}) =
\maq(\chi_1,\dots,\chi_{2i-1})^{\sigma}$,
 whenever $1\leq i \leq l$. As $P_{2t}$ also fixes $\alpha_{2k}^*$, 
we get that $P_{2t}\leq G'(\chi_1,\dots,\chi_{2i-1})$.
So, $\maq^{\sigma}$ and $\maq^{\sigma}(\chi_1,\dots,\chi_{2i-1})$ are 
$\pi'$-Hall subgroups of $G'$ and $G'(\chi_1,\dots, \chi_{2i-1})$
respectively, (as $\maq$ and $\maq(\chi_1,\dots,\chi_{2i-1})$ are)
for all $i=1,\dots,l$. 
Hence \eqref{b.te1} holds for the group $\maq^{\sigma}$ and all $i=1,\dots,l$.

Furthermore, $P_{2t}$ fixes the characters $\beta_1, \dots,\beta_{2t-1}$,
by  \eqref{x3}. 
So for the $\sigma$-conjugate $\maq^{\sigma}$ of $\maq$ we get that 
 $\maq^{\sigma}(\beta_1,\dots,\beta_{2i-1}) = \maq(\beta_1,\dots,
\beta_{2i-1})^{\sigma}$ whenever $i=1,\dots,t$.
Hence in view of \eqref{b.e5} we get 
$$
\maq^{\sigma}(\chi_1,\dots,\chi_{2i-1}) =
\maq(\chi_1,\dots,\chi_{2i-1})^{\sigma} \leq \maq(\beta_1,\dots,
\beta_{2i-1})^{\sigma} = \maq^{\sigma}(\beta_1,\dots,\beta_{2i-1}),
$$
whenever $1\leq i \leq t$.
Thus  \eqref{b.te2} holds for $\maq^{\sigma}$ and all $i=1,\dots,t$.
As far as \eqref{b.te3} is concerned, we note that $\sigma$ was picked so that 
$\maq^{\sigma}(\chi_1,\dots,\chi_{2t-1})$ normalizes $Q_{2t+1}$. 
Also for every $i$ with $1\leq i <t$, the inductive hypothesis, along with the
fact that $P_{2t}$ normalizes $Q_{2i+1}$, implies that 
$\maq^{\sigma}(\chi_1,\dots,\chi_{2i-1}) = 
\maq(\chi_1,\dots,\chi_{2i-1})^{\sigma}$ normalizes $Q_{2i+1}$.
Hence  $\maq^{\sigma}$ also satisfies \eqref{b.te3} for all $i=1,\dots,t$. This completes the
inductive proof. So  there exists  a $\pi'$-Hall subgroup 
$\maq$ of $G'$ that satisfies \eqref{b.te1} for $i=1,\dots,l$, as well as  
\eqref{b.te2} and  \eqref{b.te3} for $i=1,\dots,l-1$.

To complete the proof of the theorem it suffices to show that 
$\maq$ satisfies \eqref{b.te2} for $i=l$.
 The group $\maq$  
we have picked so far satisfies that  extra  condition. Indeed, 
$\maq(\chi_1,\dots,\chi_{2l-3}) \leq \maq(\beta_1,\dots,\beta_{2l-3})$, by \eqref{b.te2} for $i=l-1$. 
So 
$\maq(\chi_1,\dots,\chi_{2l-3})$ normalizes the groups 
$Q_1,\dots,Q_{2l-3}$. It also normalizes $Q_{2l-1}$ by \eqref{b.te3}
for $i=l-1$.
 Hence Lemma \ref{b.lem0} with $T=\maq(\chi_1,\dots,\chi_{2l-3})$,   implies 
 that $\maq(\chi_1,\dots,\chi_{2l-1}) \leq \maq(\beta_1,\dots, \beta_{2l-1})$.
This completes the proof of the theorem. 
\end{proof}  
The following fact was proved in the $i=1$ case of Theorem \ref{b.th1}. We state it here separately 
as we will use it again later. 
\begin{remark} \mylabel{b.rem1.1}
The group $G'(\beta_{1})$ normalizes $G_3(\alpha_2, \beta_1)$.
\end{remark}

As an easy consequence of Proposition \ref{b.pro1} and Theorem \ref{b.th1}
we get
\begin{corollary} \mylabel{b.cor1}
For every $i=1,\dots,k$  we have 
$$
\Hall_{\pi'}(G'(\beta_1,\dots,\beta_{2i-1})) \subseteq
\Hall_{\pi'}(G'(\chi_1,\dots,\chi_{2i}))\subseteq
\Hall_{\pi'}(G'(\chi_1,\dots,\chi_{2i-1})), 
$$
 while  
$$
\Hall_{\pi'}(G'(\beta_1,\dots,\beta_{2l-1}))\subseteq
\Hall_{\pi'}(G'(\chi_1,\dots,\chi_{2l-1})).
$$
\end{corollary}

\begin{proof}
We only need to note that in view of \eqref{b.e2.5} we have
$$
G'(\chi_1,\dots,\chi_{2i-1})\geq G'(\chi_1,\dots,\chi_{2i})\geq
G'(\beta_1,\dots,\beta_{2i-1}).
$$
 The rest follows  from 
\eqref{b.te1},  \eqref{b.te2} and \eqref{b.e2}.
\end{proof}

A similar statement to that of  Corollary \ref{b.cor1} 
 holds  for  $G'(\beta_{2i-1,2k})$ and 
$G'(\beta_1,\dots,\beta_{2i-1})$.  To prove it we start with the following 
general lemma.
\begin{lemma}\mylabel{pq51}
Assume that $\mathcal{P}$ is a  $\pi$-subgroup of a finite 
 group $\mathcal{G}$, and  
 that $\mathcal{S}_1, \mathcal{S}$ and $\mathcal{T}$ are $\pi'$-subgroups of
$\mathcal{G}$ such that
$\mathcal{T}$ normalizes $\mathcal{P}$, that 
 $\mathcal{S}$ normalizes the 
semidirect product $\mathcal{T} \ltimes \mathcal{P}$, and that 
 $\mathcal{S}_1$ is a subgroup
of $\mathcal{S}$  normalizing   $\mathcal{T}$.
Then there exists $t\in \mathcal{P}$ so that the following three conditions are
satisfied: 
\begin{itemize}
\item[(i)]$\mathcal{S}^t$ normalizes $\mathcal{T}$, 
\item[(ii)] $\mathcal{S}_1\leq \mathcal{S}^t$  and 
\item[(iii)] $t$ centralizes $\mathcal{S}_1$.   
\end{itemize}
\end{lemma}

\begin{proof}
As the $\pi'$-group $\mathcal{S}$ normalizes  the product $\mathcal{T} \ltimes \mathcal{P}$, it will
normalize  one of the  $\pi'$-Hall subgroups of that product.
 Therefore there exists $s \in \mathcal{P}$ such that $\mathcal{S}^s$
normalizes $\mathcal{T}$.
Thus $\mathcal{S}^s \mathcal{T}$ forms a $\pi'$-Hall subgroup of 
$\mathcal{S}^s\mathcal{T} \ltimes \mathcal{P}$.
As $\mathcal{S}_1$ normalizes $\mathcal{T}$, the product $\mathcal{S}_1 \mathcal{T}$ is  
 a $\pi'$-subgroup of $\mathcal{S}^s \mathcal{T} \ltimes \mathcal{P} = <\mathcal{S},\mathcal{T},\mathcal{P}>$.
 So there is some $x \in \mathcal{P}$ such that the $\pi'$-Hall subgroup 
$(\mathcal{S}^s\mathcal{T})^x$ of $\mathcal{S}^s \mathcal{T} \ltimes \mathcal{P}$,
 contains $\mathcal{S}_1\mathcal{T}$.
Hence
$\mathcal{T} \leq \mathcal{S}_1\mathcal{T} \leq \mathcal{S}^{sx}\mathcal{T}^x$.
As $\mathcal{T}$ is a $\pi'$-Hall subgroup of $\mathcal{T} \ltimes \mathcal{P}$, it follows that 
$$\mathcal{T} = \mathcal{S}^{sx}\mathcal{T}^x \cap (\mathcal{T} \ltimes \mathcal{P}).$$
Now, the group $\mathcal{T}^x$ clearly normalizes the intersection
$\mathcal{S}^{sx}\mathcal{T}^x \cap (\mathcal{T} \ltimes \mathcal{P}) = \mathcal{T}$. Hence  
$$\mathcal{T}^x = \mathcal{T}.$$
Even more, $\mathcal{S}_1\leq \mathcal{S}_1\mathcal{T} \leq \mathcal{S}^{sx}
\mathcal{T}^x$. So
$$\mathcal{S}_1 \leq \mathcal{S}^{sx} \mathcal{T}  \cap (\mathcal{S}\ltimes \mathcal{P}) = \mathcal{S}^{sx}.$$
Hence if we take $t= sx$, then (i) and (ii) are both satisfied.
To see that $t$ also centralizes $\mathcal{S}_1$ we observe
that $\mathcal{S}$ and its subgroup $\mathcal{S}_1$ normalize the unique 
$\pi$-Hall subgroup $\mathcal{P}$ of $\mathcal{T}\ltimes \mathcal{P}$. So 
 the commutator 
$[t, \mathcal{S}_1]$ is contained  in $\mathcal{P}$.
Furthermore, if $a \in \mathcal{S}_1$ then $a^t \in \mathcal{S}^t$, 
as $\mathcal{S}_1 \leq \mathcal{S}$.
So $a^{-1} a^t \in \mathcal{S}_1 \mathcal{S}^t = \mathcal{S}^t$. 
But $a^{-1} a^t= [a,t] \in \mathcal{P}$.
Hence $[a,t] \in \mathcal{S}^t \cap \mathcal{P} = 1$,
and the lemma follows.
\end{proof}

We can now prove 
\begin{theorem}\mylabel{b.th2}
There exists  a $\pi'$-Hall subgroup $\maq$ of $G'$ such that 
for every $i=1,\dots,k$ we have 
\begin{subequations}
\begin{align}
 T \leq \maq(\beta_{2i-1,2k}) &\in \Hall_{\pi'}(G'(\beta_{2i-1,2k}))
\mylabel{b.te4} \\
 \maq(\beta_{2i-1,2k}) &\leq G'(\beta_1,\dots,\beta_{2i-1}), \mylabel{b.te5} 
\end{align}
where $T:= Q_{2l-1,2k}=Q_{2k-1,2k}$ in case of an even $m$,  or 
$T:=Q_{2l-1}$ in case of an odd $m$.
Furthermore, for every $i=1,\dots,l-1$ we have 
\begin{equation}\mylabel{b.te6}
 \maq(\beta_{2i-1,2k}) \text{ normalizes \, }  Q_{2i+1}. 
\end{equation}
\end{subequations}
\end{theorem}

\begin{proof}
Note that the group $\maq$ in this theorem need not be the same as the group $\maq$
in Theorem \ref{b.th1}.

Assume that $m$ is even. Then 
the group $Q_{2k-1,2k}$ fixes $\alpha_{2k}$  as $Q_{2k-1,2k} = 
C(P_{2k} \tin Q_{2k-1})$. Thus   Proposition  \ref{pqp*t} implies that 
 $Q_{2k-1,2k}$ fixes $\alpha_{2k}^*$.
Hence in the case of an even $m$ we have $T=Q_{2k-1,2k} \leq G'$.
As we clearly 
 have that $Q_{2k-1,2k}$  fixes $\beta_{2k-1,2k}$, we conclude that 
  $T=Q_{2k-1,2k} \leq G'(\beta_{2k-1,2k})$.
If $m$ is odd, then  Corollary 
\ref{p*cor1} implies that $Q_{2l-1}$ fixes $\alpha_{2k}^*$
as it fixes $\alpha_{2k}$ (see \eqref{x5}).
 Hence,   in the case of an odd $m$, we have  that  $T= Q_{2l-1} \leq G'$.
Furthermore,  Proposition \ref{pqremark1'}
implies that $T=Q_{2l-1}$ fixes 
$\beta_{2k-1,2k}$. Therefore,  in the case of an odd $m$, and thus in every case, we get that 
 $T \leq G'(\beta_{2k-1,2k})$.
In view of Table \ref{pq12a} we also get that,  independent 
of the type of $T$, we have
\begin{equation}\mylabel{b.eq5.1}
Q_{1,2k}\leq Q_{3,2k} \leq \dots \leq Q_{2k-1,2k} \leq T.
\end{equation}

Furthermore,  $G'(\beta_{2i-1,2k})$ normalizes 
$Q_{2i-1,2k}$,  and thus normalizes 
$Q_{2r-1,2k} = Q_{2i-1,2k} \cap G_{2r-1}$ for all $r=1,\dots,i$ and all 
$i=1,\dots,k$ (see Remark \ref{extra}).
 Hence  $G'(\beta_{2i-1,2k})$ fixes the unique character $\beta_{2r-1,k}$
of  $Q_{2r-1,2k}$ that lies under $\beta_{2i-1,2k}$
 (see Proposition \ref{pqremark1'}).
Hence   $G'(\beta_{2i-1,2k}) \leq  G'(\beta_{2r-1,2k})$
whenever $1\leq r \leq i \leq k$. This implies that
we have the following series of subgroups
\begin{equation}\mylabel{b.e6} 
T  \leq G'(\beta_{2k-1,2k})
 \leq G'(\beta_{2k-3,2k}) \leq \dots \leq  G'(\beta_{1,2k}) \leq G',
\end{equation}
that is independent of the type of $T$.
So it is clear that there exists  a $\pi'$-Hall subgroup $\maq$ of $G'$ that
 satisfies \eqref{b.te4} for all $i=1,\dots,k$.
As in the proof of Theorem \ref{b.th1},  we will use induction on $i$ to modify $\maq$ 
so that the rest of the theorem 
also holds.
 
  For $i=1$ we note that $G'(\beta_{1,2k})$  
 normalizes the groups $P_{2k}^*$ and $Q_1$. Hence it fixes
 the $P_{2k}^*$-Glauberman correspondent $\beta_1$ of $\beta_{1,2k}$.
Thus  $G'(\beta_{1,2k}) \leq G'(\beta_{1})$.
In addition,  we have seen (see Remark \ref{b.rem1.1}), 
that $G'(\beta_{1})$ normalizes $G_3(\alpha_2, \beta_1)$.
So $\maq(\beta_{1,2k})$ normalizes $G_3(\alpha_2,\beta_1) = Q_3\ltimes P_2$.
Even more,  $T$ is a subgroup of $\maq(\beta_{1,2k})$ and 
 normalizes $Q_3$, as $Q_{2k-1}$ and $Q_{2l-1}$ do.
Hence Lemma \ref{pq51}, with $\mathcal{S} = \maq(\beta_{1,2k})$, $\mathcal{S}_1=T$, 
$\mathcal{T}=Q_3$ and $\mathcal{P}= P_2$, 
implies  that there exists some  $\sigma \in C(T \tin P_2)$
such that $\maq(\beta_{1,2k})^{\sigma}$ normalizes $Q_3 $.
As $\sigma \in C(T \tin P_2)$,  we get that $\sigma$ also centralizes 
$Q_{2i-1,2k}$ for all $i=1,\dots,k$ (see \eqref{b.eq5.1}).
Hence  $\sigma$ fixes the characters  $\beta_{2i-1,2k} \in \Irr(Q_{2i-1,2k})$ for all such $i$..
As $\sigma \in P_2 \leq G'=G(\alpha_{2k}^*)$,  we get that     $\sigma \in G'(\beta_{2i-1,2k})$.
This implies  that 
$\maq(\beta_{2i-1,2k})^{\sigma}=\maq^{\sigma}(\beta_{2i-1,2k})$, 
and  that $\maq(\beta_{2i-1,2k})^{\sigma} \in 
\Hall_{\pi'}(G'(\beta_{2i-1,2k}))$ for all  $i=1,\dots,k$. 
 Thus we can now work with  $\maq^{\sigma}$ in the place of $\maq$
and conclude that this  $\pi'$-Hall subgroup of $G'$ not only satisfies 
\eqref{b.te4}  for every $i=1,\dots,k$, but also the rest of the theorem 
for $i=1$.

We will work similarly for the inductive step.
So assume that $\maq$ has been modified so  that it satisfies  \eqref{b.te4}  
for all $i=1,\dots,k$  
and,  in addition,  satisfies  the rest of the  theorem for all $i
\leq t-1$ and some $t=1,\dots,l-1$.
We will show that there is a $G'(\beta_{2k-1,2k})$-conjugate of $\maq$
that satisfies \eqref{b.te4} for $i=1,\dots,k$  and the rest of 
Theorem \ref{b.th2} whenever $1\leq i \leq t$.
 The argument here  is very similar to 
the one we gave for the proof of Theorem \ref{b.th1}.
The only important difference is that  $P_{2t}$ is not in general a
 subgroup  of $G'(\beta_{2i-1,2k})$. So we can't just conjugate $\maq$ 
by an arbitary element of $P_{2t}$, 
 or else \eqref{b.te4} need  not hold for that conjugate.
But  Lemma \ref{pq51} will solve this difficulty as 
$\sigma$ can be picked from 
$C(T \tin P_{2t})$.

In view of \eqref{b.e6} and the inductive hypothesis (as $\maq$ satisfies 
\eqref{b.te4} for all $i=1,\dots,k$)  we have that 
$$
T \leq \maq(\beta_{2t-1,2k}) \leq \maq(\beta_{2t-3,2k}) \leq \dots \leq 
 \maq(\beta_{1,2k}).
$$ 
Therefore,  $\maq(\beta_{2t-1,2k})$ normalizes $Q_{2i+1}$ for all 
$i=1,\dots,t-1$, as  $\maq(\beta_{2i-1,2k})$ does.
In view of Remark \ref{b.rem1}, the group $\maq$ fixes the characters 
 $\alpha_2^*,\alpha_4^*, \dots, \alpha_{2k}^*$ and normalizes
the groups $P_{2}^*, P_4^*,\dots,P_{2k}^*$.
Hence  $\maq(\beta_{2t-1,2k})$ normalizes the groups 
$P_{2i}$ (see \eqref{pq22ii})
as well as the product groups $P_{2i}\cdots P_{2k}$ 
(see \eqref{pq22i}) whenever $1\leq i \leq t$.

Since   $\maq(\beta_{2t-1,2k})$ fixes the characters 
$\beta_{2i-1,2k}$ and normalizes the groups $Q_{2i-1}$ and $P_{2i} 
\cdots P_{2k}$, we conclude that it also fixes the $P_{2i} 
\cdots P_{2k}$-Glauberman correspondent $\beta_{2i-1} \in \Irr(Q_{2i-1})$
of $\beta_{2i-1,2k} \in \Irr(Q_{2i-1,2k})$, whenever $1\leq i \leq t$.
Hence $\maq(\beta_{2t-1,2k}) \leq G'(\beta_1,\dots,\beta_{2t-1})$.
Therefore,  in view of \eqref{b.e2.5},  we get that 
\begin{equation}\mylabel{b.e6.1}
\maq(\beta_{2t-1,2k}) \leq G'(\beta_1,\dots,\beta_{2t-1})\leq
 G'(\chi_1,\dots,\chi_{2t}).
\end{equation}
This, along with the inductive hypothesis and \eqref{b.e2.5}, implies  that 
\begin{equation}\mylabel{b.e7}
\maq(\beta_{2i-1,2k}) \leq G'(\beta_1,\dots,\beta_{2i-1}) \leq 
G'(\chi_1,\dots,\chi_{2i}), 
\end{equation}
for all $i=1,\dots,t$.

If we collect all the groups that 
$\maq(\beta_{2t-1,2k})$ normalizes,  and the characters it fixes, we have
$$
\maq(\beta_{2t-1,2k})\leq N(P_2,\dots,P_{2t},Q_1,\dots,Q_{2t-1} \tin 
G(\chi_1,\dots,\chi_{2t})).
$$
Hence $\maq(\beta_{2t-1,2k})$ normalizes
the group $N(P_2,\dots,P_{2t},Q_1,\dots,Q_{2t-1} \tin 
G_{2t+1}(\chi_1,\dots,\chi_{2t}))$.
In view of \eqref{tow--tri3} and \eqref{tow--tri2}, the latter group 
is $G_{2t+1,2t} = P_{2t} \rtimes Q_{2t+1}$. 
Furthermore, $T$  normalizes  $Q_{2t+1}$.
(Note here that,  at the last case where   $t=l-1$,   the last $\pi'$-group is 
$Q_{2t+1}=Q_{2l-1}$. 
 We still have that $T$ normalizes $Q_{2t+1}$ as, if $m$ is
even,
 then $Q_{2t+1}=Q_{2l-1}=Q_{2k-1}$ is normalized by  $Q_{2k-1,2k}=T$,  while 
if $m$ is odd, then clearly $T=Q_{2l-1}$ normalizes $Q_{2t+1}=Q_{2l-1}$.)
 The inductive hypothesis implies that $T \leq \maq(\beta_{2t-1,2k})$. 
Also  $\maq(\beta_{2t-1,2k})$ 
normalizes $P_{2t} \rtimes Q_{2t+1}$,  while its subgroup $T$ normalizes 
$Q_{2t+1}$. Therefore 
Lemma \ref{pq51} applies and provides an element  $s \in C(T \tin P_{2t})$ such that 
$\maq(\beta_{2t-1,2k})^s$ normalizes $Q_{2t+1}$.
As $s \in C(T \tin P_{2t})$, the inclusions  \eqref{b.eq5.1} imply  that 
$s$ centralizes $Q_{2i-1,2k}$ for all $i=1,\dots,k$.
Hence $s \in G(\alpha_{2k}^*, \beta_{2i-1,2k})$. Thus 
\begin{align*}
T \leq  \maq(\beta_{2i-1,2k})^s = \maq^s(\beta_{2i-1,2k}) \qquad \text{ and }\\
\maq(\beta_{2i-1,2k})^s &\in \Hall_{\pi'}(G'(\beta_{2i-1,2k}) ), 
\end{align*}
whenever $1\leq i \leq k$.
So $\maq^s$ satisfies \eqref{b.te4} for all such $i$.

Also $\maq^s(\beta_{2t-1,2k})$ normalizes $Q_{2t+1}$,  while 
for all $i=1,\dots,t-1$ the group $P_{2t}$ normalizes $Q_{2i+1}$. So 
$\maq(\beta_{2i-1,2k})^s $ normalizes $Q_{2i+1}$ 
as $\maq(\beta_{2i-1,2k})$ does and $s \in P_{2t}$.
Hence $\maq^s$ satisfies \eqref{b.te6} for all $i=1,\dots,t$.

Furthermore, as $P_{2t}$ fixes the characters $\beta_1,\dots,\beta_{2t-1}$, 
we get 
\begin{align*}
\maq^s(\beta_{2i-1,2k}) = \maq(\beta_{2i-1,2k})^s 
&\leq  \maq(\beta_1,\dots,\beta_{2i-1})^s =   &\qquad \text{  by \eqref{b.e7}} \\
&\maq^s(\beta_{1},\dots,\beta_{2i-1}), &\qquad \text{ as $s \in P_{2t}$,  }
\end{align*}
whenever $1\leq i \leq t$.
Thus $\maq^s$ also satisfies \eqref{b.te5} for all such $i$.
This completes the proof of the inductive step,  and thus provides a $\maq \in \Hall_{\pi'}(G')$ 
satisfying both \eqref{b.te4} and \eqref{b.te5}, 
for all $i=1,\dots,k$,     along with \eqref{b.te6},  
for all $i=1,\dots,l-1$.

To complete the proof of the theorem is enough to check that \eqref{b.te5} 
also holds for $i=k$ (as $k \leq l \leq k+1$).  
The argument is exactly the same one we used at the inductive step to
 prove \eqref{b.e6.1}.  So we omit it.
\end{proof}

An obvious 
 consequence of Proposition \ref{b.pro1} and Theorem \ref{b.th2} is 
\begin{corollary}\mylabel{b.cor2}
For all $i=1,\dots,k$ we have 
$$
 \Hall_{\pi'}(G'(\beta_1,\dots,\beta_{2i-1}))\subseteq 
\Hall_{\pi'}(G'(\beta_{2i-1,2k})).
$$
\end{corollary}

\noindent
We can now introduce the group $\qw$. 
\begin{theorem}\mylabel{hat:p1}
There exists a $\pi'$-subgroup $\qw$  of $G'= G(\alpha_{2k}^*)$
such that 
\begin{subequations}\mylabel{pq30}
\begin{align}
\mylabel{pq30a}  \qw &\in \Hall_{\pi'}(G'),  \\
 \mylabel{pq30b}\qw(\beta_{2i-1,2k}) \in \Hall_{\pi'}(G'(\beta_{2i-1,2k}))
&\cap
\Hall_{\pi'}(G'(\chi_1,\dots,\chi_{2i-1}))\cap \notag\\
&\Hall_{\pi'}(G'(\chi_1,\dots,\chi_{2i})) \cap
\Hall_{\pi'}(G'(\beta_1,\dots,\beta_{2i-1})),   \\
\mylabel{pq31a} \qw(\beta_{2i-1,2k}) =\qw(\chi_1,\dots,\chi_{2i-1}) &= 
\qw(\chi_1,\dots,\chi_{2i}) =\qw(\beta_1,\dots,\beta_{2i-1}) \text{ and }\\
\mylabel{hat:p1e} \qw(\chi_1,\dots,\chi_{2i-1}) &\leq \qw(\alpha_2, \dots,
\alpha_{2i}), 
\end{align}
\end{subequations}
for all $i=1,\dots,k$. In addition,  for all $i$ with $1\leq i \leq l-1$ we get  
\begin{equation}\mylabel{pq31b}
 \qw(\beta_{2i-1,2k}) \text{ normalizes } \, Q_{2i+1}. 
\end{equation}
\end{theorem} 

\begin{proof}
Let $\maq$ be any $\pi'$-group satisfying  the conditions in 
Theorem \ref{b.th2}. We will show that $\qw:=\maq$ is the desired group. 

Clearly $\qw$ satisfies \eqref{pq30a} and \eqref{pq31b}, for all $i=1,\dots,l-1$, 
 as $\maq$ is a $\pi'$-Hall
subgroup of $G'$ that  satisfies \eqref{b.te6}.
 Furthermore, 
\eqref{b.te5} and \eqref{b.e1} imply   that 
$$
\qw(\beta_{2i-1,2k}) \leq \qw(\beta_1,\dots,\beta_{2i-1}) \leq
\qw(\beta_{2i-1,2k}).
$$
So  $\qw(\beta_{2i-1,2k}) = \qw(\beta_1, \dots,\beta_{2i-1})$
whenever $1\leq i \leq k$.  This,  along with \eqref{b.te4} 
and Corollary \ref{b.cor2},  implies that 
$$ \qw(\beta_{2i-1,2k}) =\qw( \beta_1,\dots,\beta_{2i-1}) \in
\Hall_{\pi'}(G'(\beta_{2i-1,2k})) \cap  \Hall_{\pi'}(G'(\beta_1,\dots,\beta_{2i-1}))
$$
  for all $i=1,\dots,k$.  Which, in view of   Corollary  \ref{b.cor1}, 
implies that the group $\qw$   satisfies 
 \eqref{pq30b}. 

According to   \eqref{b.e2.5} we also get  that
$$ 
\qw(\beta_{2i-1,2k})= \qw(\beta_1,\dots,\beta_{2i-1}) \leq \qw(\chi_1,\dots,\chi_{2i})\leq
\qw(\chi_1,\dots,\chi_{2i-1}),
$$
as $\qw$ is a subgroup of $G'$. Since  $\qw(\beta_{2i-1,2k})$ is a $\pi'$-Hall subgroup of both 
$G'(\chi_1,\dots,\chi_{2i})$ and $G'(\chi_1,\dots,\chi_{2i-1})$ 
 we have equality everywhere, and \eqref{pq31a} follows.

It remains to show that \eqref{hat:p1e} also holds for $\qw$.
It follows from \eqref{pq31a} that $\qw(\beta_{2i-1,2k})$ normalizes the groups 
$Q_1,\dots,Q_{2i-1}$ and fixes the characters 
 $\alpha_2^*,
\dots,\alpha_{2k}^*$ (by Remark \ref{b.rem1}),  for all $i=1,\dots,k$. 
So it fixes the $Q_3,\dots,Q_{2j-1}$-correspondent $\alpha_{2j}$
of $\alpha_{2j}^*$,  for all $j=1,\dots,i$.   
This implies \eqref{hat:p1e}.  Thus Theorem \ref{hat:p1} holds.
\end{proof}

\begin{corollary}\mylabel{b.cor2.5}
Assume $\qw$ satisfies the conditions in  Theorem \ref{hat:p1}, and that 
\begin{subequations}
\begin{equation}
t \in G'(\beta_{2k-1,2k}) \cap G'(\chi_1,\dots, \chi_{2k}) \cap 
G'(\beta_1, \dots, \beta_{2k-1}).
\end{equation}
In addition, assume that  
\begin{equation}
t \in   N(Q_{2k+1} \tin G), 
\end{equation}
 if $m=2l-1= 2k+1$ is odd. 
\end{subequations}
Then $\qw^t$ also satisfies the conditions in  Theorem \ref{hat:p1}. 
\end{corollary}

\begin{proof}
Obviously $\qw^t$ is a $\pi'$-Hall subgroup of $G'$, as $t \in G'$. 
The fact that $\beta_{2i-1,2k}$ is the unique character of $Q_{2i-1,2k}$ that
lies under $\beta_{2k-1,2k}$, implies that 
$t \in G'(\beta_{2k-1,2k})$  fixes $\beta_{2i-1,2k}$,  for all 
$i=1,\dots, k$.
Hence $\qw^t(\beta_{2i-1,2k})= \qw(\beta_{2i-1,2k})^t$, for all such $i$. 
Furthermore, 
the definition of $t$ implies that 
$$
t \in G'(\beta_{2i-1,2k}) \cap G'(\chi_1,\dots, \chi_{2i-1})\cap 
G'(\chi_1,\dots, \chi_{2i}) \cap 
G'(\beta_1, \dots, \beta_{2i-1}), 
$$ 
for all $i =1, \dots, k$.
As $\qw(\beta_{2i-1,2k})$ satisfies \eqref{pq30b}, we conclude that 
$\qw^t(\beta_{2i-1,2k})= \qw(\beta_{2i-1,2k})^t$ also satisfies \eqref{pq30b}.

Even more, as $t$ fixes the various  characters  $\beta_{2i-1,2k}, \chi_{j}, 
\beta_{2i-1}$, for all $i=1,\dots,k$ and $j=1,\dots, 2k$, while \eqref{pq31a}
 holds for $\qw$, we get
\begin{multline*}
\qw^t(\beta_{2i-1,2k})= \qw(\beta_{2i-1,2k})^t= \qw(\chi_1,\dots,
\chi_{2i-1})^t= \qw^t(\chi_1, \dots, \chi_{2i-1}),\\
\qw^t(\beta_{2i-1,2k})= \qw(\beta_{2i-1,2k})^t= \qw(\chi_1,\dots,\chi_{2i})^t= 
 \qw^t(\chi_1, \dots, \chi_{2i}),\\
\qw^t(\beta_{2i-1,2k})= \qw(\beta_{2i-1,2k})^t= \qw(\beta_1,\dots,
\beta_{2i-1})^t=  \qw^t(\beta_1, \dots, \beta_{2i}).
\end{multline*}
Thus $\qw^t$ satisfies \eqref{pq31a}.

Also $t$ fixes $\alpha_{2i}^*$, for all $i=1,\dots, k$, 
   as  $t \in G' = G(\alpha_{2k}^*)$. Furthermore it normalizes the groups 
$Q_3, \dots, Q_{2i-1}$, as it fixes $\beta_1, \dots, \beta_{2i-1}$, for 
all such $i$.
Hence $t$ fixes the $Q_3, \dots, Q_{2i-1}$-correspondent 
$\alpha_{2i}$ of $\alpha_{2i}^*$. Hence
$$
\qw^t(\chi_1,\dots,\chi_{2i-1}) = \qw(\chi_1,\dots,\chi_{2i-1})^t \leq 
\qw(\alpha_2, \dots, \alpha_{2i})^t = \qw^t(\alpha_2, \dots, \alpha_{2i}).
$$
Thus \eqref{hat:p1e} holds for the $\qw^t$.

By hypothesis $t$ normalizes $Q_{2k+1}=Q_{2l-1}$,
 in the case of an odd $m=2l-1$.
Thus $t$ normalizes $Q_{2i+1}$, whenever $1\leq i \leq k$.
Hence \eqref{pq31b}, for $\qw$, implies that
$\qw^t(\beta_{2i-1,2k})= \qw(\beta_{2i-1,2k})^t$ normalizes
 $Q_{2i+1}^t=Q_{2i+1}$, for all $i=1,\dots, l-1$.
So $\qw^t$ satisfies \eqref{pq31b}.
This completes the proof  of the corollary.
\end{proof}

From now on, we write $\widehat{Q}$ for a $\pi'$-group that satisfies
all the conditions of Theorem  \ref{hat:p1}, for a
 fixed system character tower--triangular set.
An easy observation that follows from Theorem \ref{hat:p1} is 
\begin{corollary}\mylabel{b.cB}
Assume that the normal series  
$1=G_0 \unlhd \dots \unlhd G_{2k+1} \unlhd G$ for $G$ 
is fixed (so $m=2k+1$ is odd). Assume also that a character tower 
$\{\chi_i \in \Irr(G_i)\}_{i=0}^{2k+1}$ for that series and 
its corresponding triangular set  
\begin{subequations}\mylabel{b.cB1}
\begin{equation}
\{ P_{2i}, Q_{2i+1} |\alpha_{2i}, \beta_{2i+1}\}_{i=0}^k
\end{equation}
 are fixed.
Then the reduced set
\begin{equation} 
\{P_{2i}, Q_{2i-1}, P_0=1 |\alpha_{2i}, \beta_{2i-1}, 
\alpha_0=1 \}_{i=1}^k
\end{equation}
\end{subequations}
is a  representative of the unique 
conjugacy class  of triangular sets that corresponds 
 to the character tower $\{ \chi_i \in \Irr(G_i) \} _{i=0}^{2k}$
of the normal series $1=G_0 \unlhd \dots \unlhd G_{2k} \unlhd G$ of $G$.
Furthermore, if the group  $\qw$ is picked so as to 
satisfy the conditions in Theorem \ref{hat:p1} for the 
fixed character tower $\{\chi_i \}_{i=0}^{2k+1}$ and its associate triangular set, 
then the same group satisfies the conditions in Theorem \ref{hat:p1} for the smaller 
character tower $\{ \chi_i \}_{i=0}^{2k}$ and the reduced triangular
set (\ref{b.cB1}b).
\end{corollary}

\begin{proof}
The first part of the corollary   follows immediately 
 from Remark \ref{t-tr2}. 
As far as the group $\qw$ is concerned, first  observe that the
 triangular sets in \eqref{b.cB1} share the group $P_{2k}^*$ and 
its irreducible character $\alpha_{2k}^*$. 
Thus they also share the group $G'= G(\alpha_{2k}^*)$. 
Hence if $\maq$ was picked to satisfy the conditions in Theorem \ref{hat:p1} 
for the set (\ref{b.cB1}a), then $\maq$ satisfies both 
\eqref{pq30} and 
\eqref{pq31b} 
for all $i=1,\dots, k$.
But the conditions a $q$-group  should satisfy to be the $\qw$-group 
for the reduced set (\ref{b.cB1}b), are 
\eqref{pq30} for $i=1,\dots, k$ and \eqref{pq31b} for $i=1,\dots, k-1$
(since  in  the reduced case $l=k$). 
Clearly $\maq$ satisfies those. Hence Corollary \ref{b.cB} follows.
 \end{proof}

The following proposition describes  the relation between
$\widehat{Q}(\beta_1, \cdots, \beta_{2i-1}) = \qw(\beta_{2i-1,2k})$ and 
$Q_{2i+1, 2k}$. Note that 
 $Q_{2i+1,2k}$ fixes $\alpha_{2k}$,  normalizes
$Q_1,\dots,Q_{2i-1}$ (see \eqref{x5}) and is a subgroup of  $Q_{2j+1,2k}
 \leq Q_{2j+1}$ whenever $i \leq j \leq k-1$. Hence Proposition \ref{pqp*t}
implies
\begin{remark} \mylabel{b.rem2}
$$
Q_{1,2k} \leq Q_{2i+1,2k} \leq Q_{2i+1}(\alpha_{2k}^*) \leq G'
$$
whenever $1\leq i \leq k-1$.
\end{remark}

We can now prove 
\begin{proposition}\mylabel{pq39}
For all $i=1,\dots,k-1$ we have
\begin{multline*}
\widehat{Q}(\beta_{2i-1, 2k}) \cap G_{2i+1} = \widehat{Q}(\beta_1, \dots ,
\beta_{2i-1}) \cap 
G_{2i+1} = \\
\widehat{Q}(\chi_1, \dots , \chi_{2i-1}) \cap G_{2i+1} = \widehat{Q}(\chi_1,
\dots , 
\chi_{2i}) \cap G_{2i+1} = 
Q_{2i+1, 2k}.
\end{multline*} 
\end{proposition}

\begin{proof}
Since $G_{2i+1}$ is a normal subgroup of $G$, it follows from 
\eqref{pq31a} and \eqref{pq30b} that 
$$
\qw(\beta_{2i-1,2k}) \cap G_{2i+1} =
\widehat{Q}(\chi_1, \dots, \chi_{2i}) \cap G_{2i+1} \in \Hall_{\pi'}(G_{2i+1}
(\alpha_{2k}^*, \chi_1, \dots , \chi_{2i})), 
$$
whenever $i=1,\dots, k-1$.
In view of \eqref{pq31a} and    \eqref{pq31b} the group 
$\qw(\beta_{2i-1,2k})= \widehat{Q}(\beta_1, \dots ,\beta_{2i-1})$ normalizes 
the  groups $Q_1, Q_3, \dots
,Q_{2i+1}$. Furthermore, as  \eqref{hat:p1e} implies,  
 it also normalizes the groups $P_2, P_4, \dots, P_{2i}$. Hence
$$
\widehat{Q}(\chi_1, \dots ,\chi_{2i}) \cap G_{2i+1} \leq 
N(P_2, \dots , P_{2i}, Q_3, \dots , Q_{2i+1} \tin G_{2i+1}
(\alpha_{2k}^*, \chi_1, \dots , \chi_{2i})).
$$
Therefore
$$
\widehat{Q}(\chi_1, \dots ,\chi_{2i}) \cap G_{2i+1} \in 
 \Hall_{\pi'} (N(P_2, \dots , P_{2i}, Q_3, \dots , Q_{2i-1} \tin G_{2i+1}
(\alpha_{2k}^*, \chi_1, \dots , \chi_{2i}))).
$$
According to \eqref{tow--tri3} and \eqref{tow--tri2} we have 
\begin{multline*}
N(P_2, \dots , P_{2i}, Q_3, \dots , Q_{2i-1} \tin G_{2i+1}(\alpha_{2k}^*,
\chi_1, \dots , \chi_{2i}))=\\
N(P_2, \dots , P_{2i}, Q_3, \dots , Q_{2i-1} \tin G_{2i+1}
(\chi_1, \dots , \chi_{2i}))(\alpha_{2k}^*)\\
= G_{2i+1}(\alpha_2, \dots , \alpha_{2i}, \beta_1, \dots
,\beta_{2i-1})(\alpha_{2k}^*) =
 (P_{2i}\rtimes Q_{2i+1})(\alpha_{2k}^*).
\end{multline*}
In view of \eqref{pq14b}  and Remark \ref{b.rem2} we get  
$$
Q_{2i+1}(\alpha_{2k}^*) \leq N(P_{2k}^* \tin Q_{2i+1}) = 
C(P_{2i+2} \dots P_{2k} \tin  Q_{2i+1}) = Q_{2i+1, 2k }\leq Q_{2i+1}
(\alpha_{2k}^*).
$$
Also, $P_{2i}(\alpha_{2k}^*) = P_{2i}$,  as $P_{2i} \leq P_{2k}^*$ and 
$\alpha_{2k}^* \in \Irr(P_{2k}^*)$.
Hence
$$
\widehat{Q}(\chi_1, \dots ,\chi_{2i}) \cap G_{2i+1} \in
\Hall_{\pi'}(P_{2i}\rtimes 
Q_{2i+1, 2k}).$$
Because $\widehat{Q}(\chi_1, \dots ,\chi_{2i})$ normalizes both  $Q_{2i+1}$
 and $P_{2k}^*$,   it also
normalizes
$ N(P_{2k}^* \tin Q_{2i+1})$. The latter equals $ Q_{2i+1, 2k} $  and is a 
$\pi'$-Hall subgroup of 
 $P_{2i} \rtimes Q_{2i+1, 2k}$.
From this and the preceding statement we conclude that
$$
\widehat{Q}(\chi_1, \dots , \chi_{2i-1}) \cap G_{2i+1} = \widehat{Q}(\chi_1,
\dots , 
\chi_{2i}) \cap G_{2i+1} = Q_{2i+1, 2k}.
$$
This and \eqref{pq31a} imply the proposition.
\end{proof}

\begin{defn}\mylabel{pq40a}                            
For every $i=1, \dots , l$ we define 
$$
 \widehat{Q}_{2i-1} := \widehat{Q} \cap G_{2i-1}.
$$
\end{defn}
Let $G'_{s}$ denote the group $G'_{s} = G_s(\alpha_{2k}^*)$ for every $s=0,\dots,m$.
So $G'_1\unlhd G_2'\unlhd \dots \unlhd G'_m$ is a normal series of $G'$, 
as $G_{s}$ is a normal subgroup of $G$. Thus 
 Theorem \ref{hat:p1}  implies  that 
\begin{subequations}\mylabel{pq40}
 \begin{align} 
\mylabel{pq40b}  \widehat{Q}_{2i-1} &\in \Hall_{\pi'}(G'_{2i-1}), \\
 \mylabel{pq40c}
\widehat{Q}_{2i-1}(\beta_{2j-1,2k})&\in \Hall_{\pi'}
(G'_{2i-1}(\beta_{2j-1,2k})) \cap  \notag \\
 \Hall_{\pi'}(G'_{2i-1}( \chi_1, \dots, \chi_{2j}))\cap 
 \Hall_{\pi'}(G'_{2i-1}( \chi_1, \dots, \chi_{2j-1})) &\cap
\Hall_{\pi'}(G'_{2i-1}(\beta_1,\dots, \beta_{2j-1})),  \\
\mylabel{pq40e}
 \text{ and \, } \widehat{Q}_{2i-1}(\chi_1, \dots , \chi_{2j-1}) = \widehat{Q}_{2i-1}(\chi_1,
\dots , \chi_{2j}) =
\widehat{Q}_{2i-1}(&\beta_1, \dots, \beta_{2j-1}) = \widehat{Q}_{2i-
1}(\beta_{2j-1,2k}),
\end{align}
whenever  $1\leq i, j \leq k $.
\end{subequations}
Also, for all  $i=1,\dots, k$ and  all $j=1,\dots, l-1$ we have 
\begin{equation}\mylabel{pq40d}
\widehat{Q}_{2i-1}(\beta_{2j-1,2k}) \text{ normalizes } Q_{2j+1}.
\end{equation}

Furthermore, for all $i=1,\dots,k-1$, Proposition \ref{pq39} implies  that 
\begin{equation}\mylabel{pq40f}
\widehat{Q}_{2i+1}(\beta_{2i-1,2k}) = Q_{2i+1, 2k}.
\end{equation}
As $\qw_{2i-1}(\chi_1,\dots,\chi_{2j+1}) \leq \qw_{2i-1}(\chi_1,\dots,
\chi_{2j-1})$, 
equation  \eqref{pq40e} implies that 
\begin{equation}\mylabel{pq40s}
\qw_{2i-1}(\beta_{2j+1,2k}) \leq \qw_{2i-1}(\beta_{2j-1,2k}),
\end{equation}
whenever  $1\leq i \leq k$ and $1\leq j <k$. 
Furthermore, for all $i=1,\dots, k-1$ we have that 
\begin{align*}
\qw_{2i-1}(\beta_{2i-1,2k}) &= \qw(\beta_{2i-1,2k}) \cap G_{2i-1} &\quad \notag \\
&=\qw(\beta_{2i-1,2k})\cap G_{2i+1} \cap G_{2i-1}  &\quad \notag \\
&=Q_{2i+1, 2k} \cap G_{2i-1}  &\text{ by Proposition \ref{pq39} } \notag \\
&=Q_{2i-1,2k}. &\quad
\end{align*}
In addition, \eqref{pq40s} for $i=k$ and $j=k-1$   implies that 
$\qw_{2k-1} (\beta_{2k-1,2k}) \leq \qw_{2k-1}(\beta_{2k-3,2k})$.
The latter group equals $Q_{2k-1,2k}$, according to \eqref{pq40f} 
for $i=k-1$.
Hence $\qw_{2k-1}(\beta_{2k-1,2k})\leq  Q_{2k-1,2k}$.
We obviously have that $Q_{2k-1,2k}$ is a subgroup of $G_{2k-1}'(\beta_{2k-1})$, (as $Q_{2k-1,2k}$ fixes $\alpha_{2k}^*$).  But $\qw_{2k-1}(\beta_{2k-1,2k}) $ is a $\pi'$-Hall subgroup of 
$G_{2k-1}'(\beta_{2k-1})$ (see \eqref{pq40} for $i=j =k $). 
Therefore,  $\qw_{2k-1}(\beta_{2k-1,2k})=   Q_{2k-1,2k}$.
So we conclude that 
\begin{equation}\mylabel{pq40g}
\qw_{2i-1}(\beta_{2i-1,2k})= Q_{2i-1,2k},
\end{equation}
for all $i=1,\dots,k$.
We remark here that the the group $\qw_1$ is an old familiar,  as 
$\qw_1 \in \Hall_{\pi'}(G_1')$,  while $G_1'=G_1(\alpha_{2k}^*)=Q_1(\alpha_{2k}^*) = Q_{1,2k}$.
Hence 
\begin{equation}\mylabel{pqd40a}
\qw_1=Q_{1,2k}=G_1'.
\end{equation}

We also have
\begin{proposition}\mylabel{pqrm}
\begin{align*}
\widehat{Q}_{2i-1} \in 
 \Hall_{\pi'}(N(P_{2k}^* \tin G_{2i-1}(\alpha_{2i-2}^*)))
&=  \Hall_{\pi'}(N(P_{2k}^* \tin G_{2i}(\alpha_{2i-2}^*)))  \text{ and }\\
\qw_{2i-1}(\beta_{2i-1,2k} ) \in 
\Hall_{\pi'}(N(P_{2k}^* \tin G_{2i-1}(\alpha_{2i-2}^*, \beta_{2i-1,2k}))) &=
\Hall_{\pi'}(N(P_{2k}^* \tin G_{2i}(\alpha_{2i-2}^*, \beta_{2i-1,2k})))
\end{align*}
for all $i = 1, \cdots , k$.
\end{proposition}
\begin{proof}
Let $H$ be 
any subgroup of $N(P_{2k}^* \tin G_{2i-1}(\alpha_{2i-2}^*))$.
Then $H$ normalizes $P_{2k}^*$ and thus 
 $[P_{2k}^* , H]\leq P_{2k}^*$.
 Also,  $[P_{2k}^* , H ]\leq [P_{2k}^*, G_{2i-1} ]\leq G_{2i-1}$,  
as $P_{2k}^*$ normalizes $G_{2i-1}$ for all $i= 1, \cdots ,k$. 
This, along with the fact that  $G_{2i-1}/G_{2i-2}$ is a $\pi'$-group, implies that 
$$
[P_{2k}^* , H ] \leq  P_{2k}^*\cap G_{2i-1}   = P_{2k}^* \cap G_{2i-2} = P_{2i-
2}^*.
$$
W conclude  that $H$ 
centralizes the factor group $P_{2k}^* / P_{2i-2}^*$.
If, in addition,  $H$ is a  $\pi'$-subgroup,   then it 
fixes all the irreducible characters of $P_{2k}^*$ lying above $\alpha_{2i-
2}^*$, 
as it fixes $\alpha_{2i-2}^* $ and centralizes $P_{2k}^* / P_{2i-2}^*$ (see 
Problem 13.13 in \cite{is}).
Thus $H$ fixes $\alpha_{2k}^*$, and so is contained 
in $G_{2i-1}'=G_{2i-1}(\alpha_{2k}^*)$. 
Furthermore,
$$G'_{2i-1}=G_{2i-1}(\alpha_{2k}^*) = G_{2i-1}(\alpha_2^*, \ldots , \alpha_{2k}^*) \leq
N(P_{2k}^* \tin G_{2i-1}(\alpha_{2i-2}^*)).$$

Applying the above argument to  any $H \in \Hall_{\pi'}(N(P_{2k}^* \tin G_{2i-1}(\alpha_{2i-2}^*)))$,
 we see that   $H \leq G'_{2i-1}\leq N(P_{2k}^* \tin G_{2i-1}(\alpha_{2i-2}^*))$.
 So we get that  
$$
\Hall_{\pi'}(N(P_{2k}^* \tin G_{2i-1}(\alpha_{2i-2}^*)))= 
\Hall_{\pi'}(G'_{2i-1}).
$$
Similarly,  applying  the same arguments to any 
$H \in \Hall_{\pi'}(N(P_{2k}^* \tin G_{2i-1}(\alpha_{2i-2}^*))(\beta_{2i-1,2k}))$,
we see that $H \leq G'_{2i-1}(\beta_{2i-1,2k})
\leq N(P_{2k}^* \tin G_{2i-1}(\alpha_{2i-2}^*))(\beta_{2i-1,2k})$.
 So we have that 
$$
\Hall_{\pi'}(N(P_{2k}^* \tin G_{2i-1}(\alpha_{2i-2}^*))(\beta_{2i-1,2k}))= 
\Hall_{\pi'}(G'_{2i-1}(\beta_{2i-1,2k})).
$$
This, along with \eqref{pq40b} and \eqref{pq40c},  implies 
that $\qw_{2i-1}$ and $\qw_{2i-1}(\beta_{2i-1,2k})$ are $\pi'$-Hall 
subgroups of $N(P_{2k}^* \tin G_{2i-1}(\alpha_{2i-2}^*))$ and 
$N(P_{2k}^* \tin G_{2i-1}(\alpha_{2i-2}^*, \beta_{2i-1,2k}))$, 
respectively.                                              
The rest of the proposition is   obvious,  as $G_{2i}/G_{2i-1}$ is a $\pi$-group.
\end{proof}

\begin{lemma} \mylabel{pq41a}
Let $T$ be any subgroup of  $N(P_{2k}^* \tin G)$. Then $T$ normalizes 
$P_{2i}^*$ for all $i=0,1,\dots,k$. 
Furthermore,
\begin{subequations} 
\begin{equation}\mylabel{6.30a}
N(Q_1,Q_3, \dots , Q_{2t-1} \tin TP_{2i}^*) = 
 N(Q_1,Q_3, \dots,Q_{2t-1} \tin TP_{2t-2}^*)P_{2t}P_{2t+2}\dots P_{2i},
\end{equation}
whenever $1\leq t \leq i \leq k$, and
\begin{equation}\mylabel{6.30b}
N(Q_1,Q_3, \dots , Q_{2t-1} \tin TP_{2t-2}^*) 
=  N(Q_{2t-1} \tin 
N(Q_1, Q_3, \dots,Q_{2t-3} \tin TP_{2t-4}^*)P_{2t-2}),
\end{equation}
\end{subequations}
for all $t=2,3,\dots,k$.
\end{lemma}

\begin{proof}
We first observe that,  as $T$ normalizes both $P_{2k}^*$ and $G_{2i}$,
it also normalizes $P_{2i}^* = P_{2t}^* \cap G_{2i}^*$ for all $i=0,1,\dots,k$.
Thus $T P_{2i}^*$ is a group.

Let $i=1,\dots, k$ be fixed.
We will first prove \eqref{6.30a} using induction on $t$.
In the case that $t=1$, we clearly have that 
$N(Q_1 \tin TP_{2i}^*) = TP_{2i}^* =  N(Q_1 \tin T) P_{2i}^*$, as $Q_1$
is a normal subgroup of $G$.
Thus \eqref{6.30a}  holds  (for all  $i=1,\dots,k$) when    $t=1$.

Now  assume that \eqref{6.30a} holds for all values of $t$  with $t < s$
(for our  fixed $i$), for some  $s=2,\dots,i$. We will prove that it 
also holds for $t=s$.
By the inductive hypothesis  for $t=s-1$ we have
$
N(Q_1,\dots,Q_{2s-3} \tin TP_{2i}^*) = 
N(Q_1,\dots, Q_{2s-3} \tin TP_{2s-4}^*)P_{2s-2}P_{2s}\dots P_{2i}.
$
Furthermore, 
$
N(Q_1,\dots,Q_{2s-3},Q_{2s-1} \tin TP_{2i}^*)=
N(Q_{2s-1} \tin N(Q_1,\dots,Q_{2s-3} \tin TP_{2i}^*)).
$
This, along with the inductive hypothesis,  implies that 
\begin{multline}\mylabel{hat:41a}
N(Q_1, \dots,Q_{2s-3},Q_{2s-1} \tin TP_{2i}^*) = \\
N(Q_{2s-1} \tin N(Q_1, \dots, Q_{2s-3} \tin TP_{2s-4}^*)P_{2s-2}P_{2s}
\dots P_{2i}).
\end{multline}
According to \eqref{pq14'b} the groups $P_{2s}, P_{2s+2},\dots, P_{2i}$
normalize $Q_{2s-1}$. Furthermore,
$P_{2s-2}$ normalizes the groups $Q_1,\dots,Q_{2s-3}$. Therefore we get that 
$N(Q_1,\dots,Q_{2s-3} \tin TP_{2s-4}^*)P_{2s-2} = 
N(Q_1,\dots,Q_{2s-3} \tin TP_{2s-4}^*P_{2s-2})$.
Hence,  in view of  \eqref{hat:41a},  we get 
\begin{multline*}
N(Q_1, \dots,Q_{2s-3},Q_{2s-1} \tin TP_{2i}^*)= 
 N(Q_{2s-1} \tin N(Q_1, \dots, Q_{2s-3} \tin TP_{2s-4}^*)P_{2s-2} P_{2s}
\dots P_{2i}) = \\
N(Q_{2s-1} \tin N(Q_1,\dots,Q_{2s-3} \tin TP_{2s-4}^*)P_{2s-2})
P_{2s}\dots P_{2i} =\\
 N(Q_1,\dots,Q_{2s-3},Q_{2s-1} \tin TP_{2s-2}^*)P_{2s}\dots P_{2i}.
\end{multline*}
This  completes the proof for the inductive step,  and therefore for \eqref{6.30a}.

For the second equation of the lemma  note that, according to   \eqref{pq14'b},
the group $P_{2t-2}$ normalizes the $\pi'$-groups 
 $Q_1,Q_3,\dots,Q_{2t-3}$
whenever $t=2,\dots,k$. 
Therefore we have that 
$$N(Q_1,\dots,Q_{2t-3} \tin TP_{2t-2}^*) =
N(Q_1,\dots,Q_{2t-3} \tin TP_{2t-4}^*)P_{2t-2}.
$$
So
\begin{multline*}
N(Q_1,\dots,Q_{2t-3},Q_{2t-1} \tin TP_{2t-2}^*) \\
=N(Q_{2t-1} \tin N(Q_1,\dots,Q_{2t-3} \tin TP_{2t-2}^*)) =
N(Q_{2t-1} \tin N(Q_1,\dots,Q_{2t-3} \tin TP_{2t-4}^*)P_{2t-2}).
\end{multline*}
This completes the proof of \eqref{6.30b}. Hence Lemma \ref{pq41a} is proved.
\end{proof}
 
If the  $T$ that appears in  Lemma \ref{pq41a}  fixes 
 $\beta_{2i-1,2k}$, for some $i=1,\dots,k$, then we can prove 
\begin{lemma}\mylabel{pq41b}
If \,  $T$ is any subgroup of $N(P_{2k}^* \tin G(\beta_{2i-1,2k}))$,  for some 
$i=1,\dots,k$, then $T$ normalizes $P_{2i-2}^*$,  and 
$N(Q_1,Q_3,\dots,Q_{2t-1} \tin TP_{2i-2}^*)$ fixes $\beta_1,\beta_3,\dots,
\beta_{2t-1}$,  for all $t=1,\dots,i$.
\end{lemma}

\begin{proof}
We have already seen in  Lemma \ref{pq41a},
 that $T$ normalizes $P_{2i-2}^*$.
 So $TP_{2i-2}^*$ is a group.

To prove the rest of the lemma, i.e., that
\begin{equation}\mylabel{hat:b}
 N(Q_1,Q_3,\dots,Q_{2t-1} \tin TP_{2i-2}^*) \leq 
G(\beta_1,\dots, \beta_{2t-1}),
\end{equation}
 for all $t=1,\dots,i$, 
 we will use induction on $t$. 

For $t=1$ 
it is enough  to show that $TP_{2i-2}^* = N(Q_1 \tin TP_{2i-2}^*)$
 fixes $\beta_1$.
According to Remark \ref{pqremark1'}, the irreducible character
 $\beta_{2j-1, 2k}$ is the
only character of 
$Q_{2j-1, 2k}$ lying under $\beta_{2i-1, 2k}$, for all $j=1,\dots,i$.
Therefore, $T$ fixes $\beta_{2j-1,2k}$,  as  it fixes $\beta_{2i-1, 2k}$
 and normalizes $Q_{2j-1, 2k} = Q_{2i-
1, 2k} \cap G_{2j-1}$.
Hence, $T$ fixes $\beta_{1, 2k}$,  and normalizes $Q_1$ as well as
$P_{2k}^*$. So it fixes the
 unique $P_{2k}^*$-Glauberman correspondent $\beta_1 \in \Irr(Q_1)$ of $\beta_{1,2k} $.
 Furthermore, $P_{2i-2}^*$  fixes $\beta_1$, according to \eqref{x3} and 
 the definition  \eqref{pqp*d} of $P_{2i-2}^*$.
Hence, $T P_{2i-2}^*$ fixes $\beta_1$, and \eqref{hat:b} is proved  for $t=1$. 

We assume that \eqref{hat:b}  holds for $t=1,\dots, s-1$,
 and  some $s=2,\dots, i$. We  
 will prove it also holds for $t=s$. We need to show that 
$N(Q_1, Q_3, \dots ,Q_{2s-1} \tin TP_{2i-2}^*)$ fixes the characters
$\beta_1, \beta_3,\dots,\beta_{2s-1}$. By induction for $t=s-1$ we have 
that 
$$
N(Q_1,Q_3,\dots,Q_{2s-3} \tin TP_{2i-2}^*) \leq G(\beta_1,\beta_3, \dots,\beta_{2s-3}).
$$
As $ N(Q_1, Q_3, \dots, Q_{2s-3},Q_{2s-1} \tin TP_{2i-2}^*)\leq
 N(Q_1, Q_3, \dots ,Q_{2s-3} \tin TP_{2i-2}^*)$, we conclude that 
$N(Q_1, Q_3, \dots ,Q_{2s-1} \tin TP_{2i-2}^*)$ fixes $\beta_1, \beta_3, \dots 
,\beta_{2s-3}$. Hence it is enough to show it fixes $\beta_{2s-1}$.
By \eqref{6.30a},  we have
$$N(Q_1,Q_3, \dots ,Q_{2s-1} \tin T P_{2i-2}^*) = N(Q_1,Q_3, \dots ,Q_{2s-1} \tin TP_{2s-2}^*)P_{2s} \cdots P_{2i-2}, 
$$
where, by convention, we assume that,   in the case $s=i$, we have 
 $P_{2s}\cdots P_{2s-2} = 1$. So in that case  the  equation holds trivially. 

According to \eqref{x3}, the groups 
$P_{2s}, \dots , P_{2i-2}$ fix $\beta_{2s-1}$.
Hence it is enough to show that the group 
$N(Q_1,Q_3, \dots ,Q_{2s-1} \tin T P_{2s-2}^*) $
fixes
$\beta_{2s-1}$.

As   $N(Q_1,Q_3, \dots ,Q_{2s-1} \tin T P_{2s-2}^*)$ normalizes  both 
$P_{2k}^*$ and $Q_1,\dots, Q_{2s-1}$,
it   normalizes  the product  group $ P_{2s} \cdots P_{2k}=
N(Q_1, \dots, Q_{2s-1}  \tin P_{2k}^*) $ (see \eqref{pq22i}).
Therefore
 it also normalizes $Q_{2s-1, 2k} = N(P_{2s} \cdots P_{2k} \tin Q_{2s-1})$.
Let $\sigma \in  N(Q_1, Q_3, \dots ,Q_{2s-1} \tin T P_{2s-2}^*)$. 
Then $\sigma =
\tau \cdot p_{2s-2}$ where $\tau \in T$ and $p_{2s-2} \in P_{2s-2}^*$. 
As $\tau \in T$ and $T \leq G(\beta_{2i-1,2k}) \leq  G(\beta_{2s-1,2k})$, 
we have that $\tau $  normalizes $Q_{2s-1,2k}$.
Since   $\sigma$  also normalizes 
$Q_{2s-1, 2k}$, so does   $p_{2s-2}$.  Hence $p_{2s-2} \in N(Q_{2s-1, 2k} \tin
P_{2s-2}^*)
= C(Q_{2s-1, 2k} \tin P_{2s-2}^*)$. So $p_{2s-2}$   fixes $\beta_{2s-1, 2k}$.
As  $T$   fixes $\beta_{2s-1, 2k}$, we conclude that $\sigma$  does
as well. 
Therefore, 
  $N(Q_1, Q_3, \dots ,Q_{2s-1} \tin T P_{2s-2}^*)$ fixes  
$\beta_{2s-1,2k}$,
and normalizes both  $Q_{2s-1}$ and $P_{2s} \cdots P_{2k}$.
 Hence it also fixes the
$P_{2s} \cdots P_{2k}$-Glauberman correspondent
 $\beta_{2s-1}$  of $\beta_{2s-1,2k} \in \Irr(Q_{2s-1,2k})$.
So
$$
N(Q_1, Q_3, \dots ,Q_{2s-1} \tin T P_{2i-2}^*) \leq G_{2i}(\beta_1, \dots,
\beta_{2s-1}),$$
which completes  the inductive proof of \eqref{hat:b} for all $t=1,\dots,i$.
 Thus   Lemma \ref{pq41b} holds. 
\end{proof}

\begin{lemma}\mylabel{hat:c1}
If $\map$ is any $\pi$-subgroup of $N(P_{2k}^* \tin G(\alpha_{2i-2}^*, 
\beta_{2i-1,2k}))$,  for some $i=1,\dots,k$, then $\map$ normalizes $P_{2t}^*$, 
for all $t=0,1,\dots,k$. Furthermore, 
\begin{equation}\mylabel{pq47}
N(Q_1,Q_3, \dots , Q_{2t-1} \tin \mathcal{P}\cdot P_{2t-2}^*)\cdot P_{2t-
2} = N(Q_1,Q_3, \dots , Q_{2t-3} \tin \mathcal{P}\cdot P_{2t-4}^*)P_{2t-2},
\end{equation}
whenever  $2\leq t \leq i $. 
\end{lemma}

\begin{proof}
As $\map$ normalizes $P_{2k}^*$,  it also normalizes $P_{2t}^* = P_{2k}^* 
\cap G_{2t}$ for all $t=0,1,\dots,k$.

To prove \eqref{pq47}
we will first do the case $t = 2$ of the equation, even though it follows from
the general case, 
just to show the argument (which is nothing else but a Frattini argument) in
its easiest form.
According to Lemma \ref{pq41b}, for $\mathcal{P}$ in the place of $T$,
and for $t$ and $i$ there both equal to $1$, he character
$\beta_1$ is fixed  by $\mathcal{P}$. Obviously $\map$
 fixes $\alpha_2 = \alpha_2^*$, as $i \geq t=2$ and $\map$ fixes $\alpha_{2i-2}^*$ .  
 Therefore it  normalizes $G_3(\alpha_2, \beta_1) = Q_3 \ltimes P_2$.
Hence $\mathcal{P}(Q_3\ltimes P_2)$ is a group, with $Q_3 \ltimes P_2$
as a normal subgroup. Furthermore, all the $\pi'$-Hall subgroups of 
$Q_3 \ltimes P_2$ are the $P_2$-conjugates of 
$Q_3$. So the group $\mathcal{P}P_2$ permutes these conjugates among
themselves with $P_2 \leq \mathcal{P}P_2$ acting transitively. 
Therefore, a Frattini type argument implies that 
\begin{equation}\mylabel{hat:e7}
\mathcal{P}P_2 =  N(Q_3 \tin \mathcal{P}P_{2})P_2= N(Q_1, Q_3 \tin
 \map P_2^*)P_2.
\end{equation}
Clearly $N(Q_1 \tin \map P_0^*)P_2= \map P_2$, as $Q_1$ is normal in $G$
and $P_0^*=1$. 
Thus \eqref{pq47} is proved for $t = 2$.

In the general case, with $t\geq 3$,   we can apply  Lemma \ref{pq41b}
with the present $t-1$ as both $i$ and $t$ there.
 We get that 
 $N(Q_1, Q_3, \dots , Q_{2t-3} \tin \mathcal{P}\cdot P_{2t-4}^*)$
fixes $\beta_1,\dots ,\beta_{2t-3}$. It also fixes 
$\alpha_2^*, \dots ,\alpha_{2t-2}^*$, as $\map \leq G(\alpha_{2i-2}^*)$
and $P_{2t-2}^*$ do. In view of Proposition \ref{pqp*t}, we conclude that 
 $N(Q_1, Q_3, \dots , Q_{2t-3} \tin \mathcal{P}\cdot P_{2t-4}^*)$ fixes 
$\alpha_2,\dots,\alpha_{2t-2}$.
Therefore
$N(Q_1, Q_3, \dots , Q_{2t-3} \tin \mathcal{P} \cdot P_{2t-4}^*)$
normalizes 
$G_{2t-1}(\beta_1, \dots ,\beta_{2t-3}, \alpha_2, \dots ,\alpha_{2t-2})$,
which   equals $Q_{2t-1} \ltimes P_{2t-2}$ 
by \eqref{tow--tri3} and \eqref{tow--tri2}.                                   
Therefore, $N(Q_1, Q_3, \dots , Q_{2t-3} \tin \mathcal{P} \cdot P_{2t-4}^*) 
\cdot (Q_{2t-1} \ltimes P_{2t-2})$ is a group with 
$ Q_{2t-1} \ltimes P_{2t-2}$ as a 
normal subgroup. As all the $\pi'$-Hall subgroups of
 $Q_{2t-1} \ltimes P_{2t-2}$ are $P_{2t-2}$-conjugates of $Q_{2t-1}$, the
group 
$ N(Q_1, Q_3, \dots , Q_{2t-3} \tin \mathcal{P} \cdot P_{2t-4}^*) \cdot P_{2t-
2}$ permutes these conjugates among themselves, while its subgroup 
$P_{2t-2}$ acts transitively on them. 
So a   Frattini type argument  implies that 
\begin{multline*}
N(Q_1, Q_3, \dots , Q_{2t-3} \tin \mathcal{P} \cdot P_{2t-4}^*) \cdot P_{2t-
2} = \\
N(Q_{2t-1} \tin N(Q_1, Q_3, \dots , Q_{2t-3} \tin \mathcal{P} \cdot P_{2t-4}^*)
\cdot P_{2t-2})\cdot P_{2t-2}.
\end{multline*}
According to \eqref{6.30b} this last group equals
$N(Q_1,Q_3, \cdots , Q_{2t-1} \tin \mathcal{P} \cdot P_{2t-2}^*)\cdot P_{2t-
2}$.
This completes the proof of Lemma \ref{pq47}.
\end{proof}

\begin{proposition}\mylabel{pq43}
For all $i= 1, \dots,k$, the group  $P_{2i}^*(\beta_{2i-1, 2k})$ is the 
  unique normal $\pi$-Hall subgroup of  the normalizer 
 $N(P_{2k}^* \tin G_{2i}(\alpha_{2i-2}^*,
 \beta_{2i-1, 2k}))$.
\end{proposition}

\begin{proof}
Let $\mathcal{P}_{2i}$ be any $\pi$-Hall subgroup of $N(P_{2k}^* \tin
G_{2i}(\alpha_{2i-2}^*, \beta_{2i-1, 2k}))$,  for some  fixed $i = 1, \dots
,k$. 
$P_{2i}^*(\beta_{2i-1, 2k})$ is a $\pi$-subgroup of $N(P_{2k}^* \tin 
G_{2i}(\alpha_{2i-2}^*, \beta_{2i-1, 2k}))$,  as it is contained in $G_{2i}\cap
P_{2k}^*$ and fixes 
$\alpha_{2i-2}^*$. Furthermore, 
 $P_{2i}^*(\beta_{2i-1,2k})$ is a normal $\pi$-subgroup of
 $N(P_{2k}^* \tin G_{2i}(\alpha_{2i-2}^*,\beta_{2i-1,2k}))$, and thus 
is contained
in every $\pi$-Hall subgroup of the latter group.
  Hence  
\begin{equation}\mylabel{pq44}
 P_{2i}^*(\beta_{2i-1, 2k}) \leq \mathcal{P}_{2i}.
\end{equation}

The group $P_{2i}$ fixes $\beta_{2i-1, 2k} \in \Irr(C(P_{2i}, \dots, P_{2k}
\tin 
Q_{2i-1}))$,  and is contained in $P_{2i}^*$. So we get that 
$P_{2i} \leq P_{2i}^*(\beta_{2i-1, 2k} ) \leq \mathcal{P}_{2i}$.
According to \eqref{pq14'b} the $\pi$-group $P_{2i}$ normalizes the
 $\pi'$-groups $Q_1, \dots ,Q_{2i-1}$. 
So  
$P_{2i} \leq N(Q_1, Q_3, \dots , Q_{2i-1} \tin \mathcal{P}_{2i}) \leq 
N(Q_1, Q_3, \dots , Q_{2i-1} \tin \mathcal{P}_{2i}\cdot P_{2i-2}^*)$.
Lemma \ref{pq41b}, for $T= \map_{2i}$ and $t=i$,  implies that 
\begin{subequations}\mylabel{pq45}
\begin{equation}\mylabel{pq45a}
P_{2i} \leq N(Q_1, Q_3, \dots , Q_{2i-1} \tin \mathcal{P}_{2i}\cdot P_{2i-
2}^*)\leq G_{2i}(\beta_1, \dots ,\beta_{2i-1}).
\end{equation}
Furthermore, $\mathcal{P}_{2i}P_{2i-2}^*$ fixes $\alpha_{2i-2}^*$ and
normalizes $P_{2j}^*$ for all
$j= 1, \cdots ,k$. Hence it also fixes the unique character $\alpha_{2j}^*$
of $P_{2j}^*$
lying under $\alpha_{2i-2}^*$, whenever $1\leq j \leq i-1$. As the group 
$N(Q_1, Q_3,\dots , Q_{2i-1}
 \tin \mathcal{P}_{2i}\cdot P_{2i-2}^*)$ normalizes $Q_1, Q_3, 
\dots , Q_{2i-1}$,
it also normalizes the groups $P_{2j} = N(Q_{2j-1}^* \tin P_{2j}^*)$, 
for all  $j=1,2,\dots,i-1$,  as well as $P_{2i}$.
Hence, according to Proposition \ref{pqp*t},  we get 
$N(Q_1, Q_3, \cdots , Q_{2i-1} \tin \mathcal{P}_{2i}\cdot P_{2i-2}^*)(\alpha_{2j})
= N(Q_1, Q_3, \cdots , Q_{2i-1} \tin \mathcal{P}_{2i}\cdot P_{2i-
2}^*)(\alpha_{2j}^*)$. This
 implies that 
\begin{equation}\mylabel{pq45b}
N(Q_1, Q_3, \cdots , Q_{2i-1} \tin \mathcal{P}_{2i}\cdot P_{2i-2}^*)  
\text{ fixes } \alpha_{2j} 
\end{equation}
for all  $j= 1, \dots ,i-1$.
\end{subequations} 
Similarly we can see that for all $j$ and $t$ with $1\leq j \leq t <i$ we have 
$N(Q_1, Q_3, \dots , Q_{2t-1} \tin \mathcal{P}_{2i}
\cdot P_{2t-2}^*)(\alpha_{2j}) = 
N(Q_1, Q_3, \dots, Q_{2t-1} \tin
 \mathcal{P}_{2i}\cdot P_{2t-2}^*)(\alpha_{2j}^*)$.
Hence
\begin{equation}\mylabel{hat:e6}
N(Q_1, Q_3, \cdots , Q_{2t-1} \tin \mathcal{P}_{2i}\cdot P_{2t-2}^*)  
\text{ fixes } \alpha_{2j}, 
\end{equation}
whenever $1\leq j \leq t < i$.

The inclusions \eqref{pq45a}, along with \eqref{pq45b}, \eqref{tow--tri3}
and \eqref{tow--tri1},  imply that 
$$
P_{2i} \leq N(Q_1, Q_3, \dots , Q_{2i-1} \tin \mathcal{P}_{2i}\cdot
 P_{2i-2}^*) \leq 
G_{2i}(\alpha_2, \dots \alpha_{2i-2}, \beta_1, \dots \beta_{2i-1}) = 
 P_{2i}\ltimes Q_{2i-1}.
$$
Since $N(Q_1, Q_3, \dots , Q_{2i-1} \tin \mathcal{P}_{2i}\cdot P_{2i-2}^*)$
 is a $\pi$-group,
  and $P_{2i}$ is a $\pi$-Hall subgroup of $P_{2i}\ltimes Q_{2i-1}$, 
it follows that
\begin{equation}\mylabel{pq46}
P_{2i} = N(Q_1, Q_3, \cdots , Q_{2i-1} \tin \mathcal{P}_{2i}\cdot P_{2i-2}^*).
\end{equation}   
                                          
To finish the proof of Proposition \ref{pq43}, we only need to show, according to
\eqref{pq44}, that 
$\mathcal{P}_{2i} \leq P_{2i}^*$.
We actually have the  stronger equality
\begin{equation}\mylabel{pq48}
P_{2i}^* = \mathcal{P}_{2i}P_{2i-2}^*.
\end{equation}
To prove \eqref{pq48} we will use Lemma \ref{hat:c1} with $\map_{2i}$ 
in the place  of $\map$. Indeed,
\begin{align*}
P_{2i}^* &= P_{2i}P_{2i-2}P_{2i-4} \cdots P_2  &\quad  &\quad \\
&= N(Q_1, Q_3, \dots , Q_{2i-1} \tin \mathcal{P}_{2i}\cdot P_{2i-2}^*)  
P_{2i-2}P_{2i-4} \cdots P_2  &\quad  &( \text{by \eqref{pq46} } )\\
&= N(Q_1, Q_3, \dots , Q_{2i-3} \tin \mathcal{P}_{2i}\cdot P_{2i-4}^*) 
P_{2i-2}P_{2i-4} \cdots P_2   &\quad &( \text{by 
\eqref{pq47} for } t = i )\\ 
&= N(Q_1, Q_3, \dots , Q_{2i-3} \tin \mathcal{P}_{2i}\cdot P_{2i-4}^*)
P_{2i-4}P_{2i-6} \cdots P_2 P_{2i-2}  &\quad &\quad\\
&= N(Q_1, Q_3, \dots , Q_{2i-5} \tin \mathcal{P}_{2i}\cdot P_{2i-6}^*)
P_{2i-6}P_{2i-8} \cdots P_2 P_{2i-2}P_{2i-4}  &\quad &(\text{by 
\eqref{pq47} for  } t= i-1) \\
&= \cdots = N(Q_1, Q_3, Q_5 \tin\mathcal{P}_{2i}\cdot P_{4}^*) 
P_4 P_2 P_{2i-2}P_{2i-
4} \cdots P_8 P_6  &\quad &\quad\\
&=N(Q_1, Q_3 \tin\mathcal{P}_{2i}\cdot P_{2}^*) P_4 P_2 P_{2i-2}P_{2i-4} \cdots
P_8 P_6  &\quad &(\text{by \eqref{pq47} for  } t = 3) \\
&= N(Q_1, Q_3 \tin\mathcal{P}_{2i}\cdot P_{2}^*) P_2 P_{2i-2}P_{2i-4} \cdots
P_8 P_6 P_4  &\quad &\quad\\
&=\mathcal{P}_{2i} P_2 P_{2i-2}P_{2i-4} \cdots
 P_6 P_4  &\quad &(\text{by \eqref{pq47} for  } t = 2)\\
 &= \mathcal{P}_{2i}P_{2i-2}^*. &\quad &
\end{align*}
Hence \eqref{pq48} holds,  and Proposition \ref{pq43} is proved.
\end{proof}

We finish this section with a complete  characterization of   
$N(P_{2k}^* \tin G_{2i}(\alpha_{2i-2}^*, \beta_{2i-1,2k}))$.

\begin{theorem}\mylabel{hat:p2}
For every $i=1,\dots,k$ we have
$$N(P_{2k}^* \tin G_{2i}( \alpha_{2i-2}^*, \beta_{2i-1,2k}))=
P_{2i}^*(\beta_{2i-1,2k}) \rtimes \qw_{2i-1}(\beta_{2i-1,2k})=
P_{2i}^*(\beta_{2i-1,2k}) \rtimes Q_{2i-1,2k}.
$$
\end{theorem}

\begin{proof}
This follows easily from Proposition \ref{pqrm}, 
 Proposition \ref{pq43} and equation \eqref{pq40g}.
\end{proof}

\section{The irreducible characters $\hb_{2i-1}$ of $\qw_{2i-1}$. } 
We are now able to define irreducible characters $ \hat{\beta}_{2i-1}$ 
of $\widehat{Q}_{2i-1}$, for all $i=1, \cdots, k$,
closely related 
to the  $\beta_{2i-1} \in  \Irr(Q_{2i-1})$.
In fact, we will prove 
\begin{proposition}\mylabel{pq40x}
For every $i=1,\dots,k$ we write 
 $\hb_{2i-1}$  for the character $\beta_{2i-1,2k}^{\qw_{2i-1}}$
of $\qw_{2i-1}$ induced by  $\beta_{2i-1,2k} \in \Irr(Q_{2i-1,2k})$. 
Then $\hb_{2i-1}$ lies in $\Irr(\qw_{2i-1}|\beta_{2i-1,2k})$, while  $\hb_1=\beta_{1,2k}$. 
\end{proposition}
\begin{proof}
Let $i= 1,\dots, k$ be fixed.
For any subgroup $H$ of $G$ containing $Q_{2i-1,2k}$, we write
$\beta_{2i-1,2k}^H$ for the induced character $(\beta_{2i-1,2k})^H$
of   $H$.
We will first prove that, for all $j=1,\dots,i$, 
we have
\begin{equation}\mylabel{hat:40x}
\beta_{2i-1, 2k}^{\widehat{Q}_{2i-1}(\beta_{2j-1, 2k})} \in
\Irr(\widehat{Q}_{2i-1}(\beta_{2j-1, 2k})
| \beta_{2i-1, 2k }, \beta_{2i-3, 2k}, \dots , \beta_{1, 2k}).
\end{equation}
For the proof of \eqref{hat:40x} we will use induction on  $i-j = 0,\dots, 
i-1$.

We treat the  case where $i=1$ separately.
Let   $i=j = 1$. Then,  according to \eqref{pqd40a}, we have that   
$ \widehat{Q}_1 = Q_{1,2k}$.
Therefore
\begin{equation}\mylabel{pqd40}
\hat{\beta}_1 = \beta_{1,2k} = \beta_{1,2k}^{Q_{1,2k}},
\end{equation}
is an irreducible character of $\qw_1 = Q_{1,2k}$. So  equation \eqref{hat:40x},  as well as 
Proposition \ref{pq40x},   holds  trivially when  $i=j =1$.

For any  $i >1$ the equalities  \eqref{pq40f} and \eqref{pq40g},
imply that 
\begin{equation}\mylabel{hat.ee}
 Q_{2i-1, 2k} = \qw_{2i-1}(\beta_{2i-1,2k})= 
 \widehat{Q}_{2i-1}(\beta_{2i-3,2k}).
\end{equation}
This, along with  the inclusion  in  \eqref{pq40s},
 implies  that we can form a series 
\begin{multline}\mylabel{hat:ser}
Q_{1,2k}\leq Q_{3,2k} \leq \dots \leq Q_{2i-3, 2k}\leq Q_{2i-1, 2k} =
\qw_{2i-1}(\beta_{2i-1,2k})=\\
 \widehat{Q}_{2i-1}(\beta_{2i-3,2k})\leq \widehat{Q}_{2i-1}(\beta_{2i-5, 2k})
\leq \dots \leq \widehat{Q}_{2i-1}(\beta_{1,2k})\leq
\widehat{Q}_{2i-1}
\end{multline}
of subgroups of $\widehat{Q}_{2i-1}$. Even more,   \eqref{pq40d},
 along with the fact that $\qw_{2i-1} \leq G'$  normalizes $P_{2k}^*$,   implies  that    
\begin{equation}\mylabel{hat.eee}
\widehat{Q}_{2i-1}(\beta_{2j-3,2k}) \text{ normalizes }
Q_{2j-1, 2k} = N(P_{2k}^* \tin Q_{2j-1}),
\end{equation}
for any $j=2,\dots,i$.
 
If $i-j= 0$, i.e., $i=j$, then $\qw_{2i-1}(\beta_{2j-1,2k}) = 
\qw_{2i-1}(\beta_{2i-1,2k}) = Q_{2i-1,2k}$
 by \eqref{hat.ee}.
Thus $\beta_{2i-1, 2k}^{\qw_{2i-1}(\beta_{2j-1,2k})} = \beta_{2i-1,2k}
\in \Irr(Q_{2i-1,2k})=
\Irr(\qw_{2i-1}(\beta_{2j-1,2k}))$. Furthermore, $\beta_{2i-1,2k}$ lies above $\beta_{2i-3, 2k}, \dots,
\beta_{1,2k}$,  as we can see in Diagram  \eqref{pq12b}. Hence \eqref{hat:40x}
holds  in the case that $i-j = 0$. As 
$\qw_{2i-1}(\beta_{2i-3,2k})=Q_{2i-1,2k}$  by \eqref{hat.ee}, 
we also get that $\beta_{2i-1,2k} = \beta_{2i-1,2k}^{\qw_{2i-1}(\beta_{2i-3,2k})}$.
Hence \eqref{hat:40x} also holds for $j=i-1$.

For the inductive step it  is enough to prove that,  if \eqref{hat:40x} holds for 
some $j =2,\dots, i-1$, then it also holds for $j-1$ (as we  induct on $i-j$).
So assume that \eqref{hat:40x} holds for $j$. It suffices to  show  that 
$$ 
\beta_{2i-1, 2k}^{\widehat{Q}_{2i-1}(\beta_{2j-3, 2k})} \in
\Irr(\widehat{Q}_{2i-1}(\beta_{2j-3, 2k})
| \beta_{2i-1, 2k }, \beta_{2i-3, 2k}, \cdots , \beta_{1, 2k}).
$$
It follows from  \eqref{hat:ser}  that 
$Q_{2j-1, 2k} \leq \qw_{2i-1}(\beta_{2j-1,2k})
 \leq \qw_{2i-1}(\beta_{2j-3,2k})$,
where $Q_{2j-1,2k}$ is a normal subgroup  of  
$ \qw_{2i-1}(\beta_{2j-3,2k})$, by \eqref{hat.eee}.
We clearly have that 
$\qw_{2i-1}(\beta_{2j-1,2k}) $ equals the group 
\linebreak
$\qw_{2i-1}(\beta_{2j-3,2k})(\beta_{2j-1,2k})$.
 Furthermore,  by the inductive hypothesis
we know that $\beta_{2i-1,2k}^{\qw_{2i-1}(\beta_{2j-1,2k})}$  is an 
irreducible character 
of $\qw_{2i-1}(\beta_{2j-1,2k})$ that lies above $\beta_{2j-1,2k}$.
So  Clifford's theory can be applied  to the normal group
 $Q_{2j-1,2k}$ of  $\qw_{2i-1}(\beta_{2j-3,2k})$, the character 
$\beta_{2j-1,2k} \in \Irr(Q_{2j-1,2k})$, the stabilizer 
$\qw_{2i-1}(\beta_{2j-1,2k})$ of that character in $\qw_{2i-1}(\beta_{2j-3,2k})$, 
and the irreducible character $\beta_{2i-1,2k}^{\qw_{2i-1}(\beta_{2j-1,2k})}$  
of $\qw_{2i-1}(\beta_{2j-1,2k})$ that lies above $\beta_{2j-1,2k}$.
 Therefore we conclude that 
$\beta_{2i-1,2k}^{\qw_{2j-3,2k}}=(\beta_{2i-1,2k}^{\qw_{2i-1}(\beta_{2j-1,2k})})^{\qw_{2i-1}(\beta_{2j-3,2k})}$
is an irreducible character of $\qw_{2i-1}(\beta_{2j-3,2k})$. Furthermore, 
it lies above $\beta_{2i-1,2k}$,  and thus  also lies  above
 $\beta_{2i-3,2k},\dots,\beta_{1,2k}$. 
This completes the proof of the inductive step.  hence   \eqref{hat:40x} holds 
for all $i=1,\dots,k$ and $j=1,\dots,i$.

To complete the proof of the proposition,  we note that,  for any fixed $i=1,\dots,k$,
equation \eqref{hat:40x} for $j=1$, implies 
that $\beta_{2i-1,2k}^{\qw_{2i-1}(\beta_{1,2k})} \in
 \Irr(\qw_{2i-1}(\beta_{1,2k})| \beta_{1,2k})$.
Furthermore,  $\widehat{Q}_{2i-1}$ 
normalizes $Q_{1, 2k} =N(P_{2k}^* \tin Q_1)$.
Thus Clifford's theory, 
applied to the groups $Q_{1,2k} \unlhd \widehat{Q}_{2i-1}$,
implies that $\beta_{2i-1, 2k}^{\qw_{2i-1}(\beta_{1,2k})}$
 induces an
irreducible character of  $\widehat{Q}_{2i-1}$ , i.e., 
\begin{equation*} 
\hat{\beta}_{2i-1} = \beta_{2i-1,2k}^{\widehat{Q}_{2i-1}} = (\beta_{2i-1,2k}^{\qw_{2i-1}(\beta_{1,2k})} )^{\qw_{2i-1}} \in 
\Irr(\widehat{Q}_{2i-1} | \beta_{2i-1, 2k}).
\end{equation*}
Hence Proposition \ref{pq40x} is proved. 
\end{proof}

The way $\hb_{2i-1}$ is picked  implies
\begin{corollary}\mylabel{hat:r1}
If $i=1,\dots,k$, then 
any subgroup of $G$  that normalizes $\widehat{Q}_{2i-1}$
and fixes $\beta_{2i-1, 2k}$ also fixes $\hat{\beta}_{2i-1}$.
 Furthermore, any subgroup of $G$ that  fixes $\hat{\beta}_{2i-1}$ and 
$\beta_{2i-3,2k}$  also   fixes $\beta_{2i-1, 2k}$.
\end{corollary}

\begin{proof}
The first statement is obvious,  since $\hb_{2i-1} = \beta_{2i-1}^{\qw_{2i-1}}$.
 
The character $\hb_{2i-1}$ is obtained from $\beta_{2i-1, 2k}$ 
using a series of characters 
\begin{multline*}
\beta_{2i-1,2k}=\beta_{2i-1,2k}^{\qw_{2i-1}(\beta_{2i-1,2k})}=
\beta_{2i-1,2k}^{\qw_{2i-1}(\beta_{2i-3,2k})},
\beta_{2i-1,2k}^{\qw_{2i-1}(\beta_{2i-5,2k})},\\
\beta_{2i-1,2k}^{\qw_{2i-1}(\beta_{2i-7,2k})},\dots,
\beta_{2i-1,2k}^{\qw_{2i-1}(\beta_{1,2k})},
\beta_{2i-1,2k}^{\qw_{2i-1}}=\hb_{2i-1},
\end{multline*}
each obtained from the preceding one using Clifford theory
for the characters 
$$\beta_{2i-5,2k}, 
\beta_{2i-7,2k}, \dots,\beta_{3,2k}, \beta_{1,2k},
$$
in that order.  Since $G(\beta_{2i-3,2k})$ fixes the characters 
$\beta_{1,2k},\beta_{3,2k},\dots, \beta_{2i-3,2k}$,
Clifford theory  implies  the rest of the proof.
\end{proof}

\section{$\pi$-Hall subgroups of
 $N(P_{2k}^*, \qw_{2i-1} \tin G_{2i}(\alpha_{2i-2}^*) )$: the groups
 $\pw_{2i}$ }
The following two general lemmas,   along with Lemmas  \ref{pq51} and 
\ref{pq43},  will help us pick
``nice'' $\pi$-Hall 
subgroups $\pw_{2i}$  of 
$N( P_{2k }^* ,\widehat{Q}_{2i-1} \tin G'_{2i})$.

\begin{lemma}\mylabel{pq49}
Assume $H$ is a finite $\pi$-separable group. 
Let  $N =
N_1 \rhd N_2
\rhd \cdots \rhd N_r$ be  a  series of normal $\pi'$-subgroups of  $H$,
 for some integer $r\geq 1$.
Let $\theta_i $ be an  irreducible character of $N_i$, for each  $i= 1, \cdots , r
$, such that 
$\theta_i \in \Irr(N_i | \theta_{i+1}, \theta_{i+2}, \cdots, \theta_r)$.
Then $H(\theta_1)\geq H(\theta_1,\theta_2) \geq 
\dots \geq H(\theta_1,\dots,\theta_r)$ and the index
 $|H(\theta_1):H(\theta_1,\dots,\theta_r)|$  is the $\pi'$-number
$|N:N(\theta_1,\dots,\theta_r)|$. Hence any $\pi$-Hall subgroup 
$P$ of $H(\theta_1,\dots,\theta_r)$ is also a $\pi$-Hall subgroup of 
$H(\theta_1,\dots,\theta_i)$, for each $i=1,\dots,r$.
\end{lemma}

\begin{proof}
For every $i=1,\dots,r-1$ the group $H(\theta_1,\dots,\theta_i)$
has as   normal subgroups the groups
$N_i \rhd N_{i+1}$, and fixes the character $\theta_i \in \Irr(N_i)$.
Hence it permutes among themselves all the irreducible characters  
of $N_{i+1}$ that lie under $\theta_i$. But this set of irreducible 
characters is precisely the $N_i$-conjugacy class of $\theta_{i+1}$
by Clifford's theory. 
So 
$$
 H(\theta_1,\dots,\theta_i)=H(\theta_1,\dots,\theta_i,\theta_{i+1})N_i.
$$
Therefore, 
$$
|H(\theta_1,\dots,\theta_i):H(\theta_1,\dots,\theta_i,\theta_{i+1})|
=|N_i:N_i(\theta_{i+1})|.
$$
A similar argument with $N$ in the place of $H$ shows that 
$|N(\theta_1,\dots,\theta_i):N(\theta_1,\dots,\theta_i,\theta_{i+1})|
=|N_i:N_i(\theta_{i+1})|$.
Hence
$$
|H(\theta_1,\dots,\theta_i):H(\theta_1,\dots,\theta_i,\theta_{i+1})|
=|N(\theta_1,\dots,\theta_i):N(\theta_1,\dots,\theta_i,\theta_{i+1})|,
$$
for all $i=1,\dots,r-1$. This implies that 
\begin{multline*}
|H(\theta_1): H(\theta_1, \dots, \theta_{r})| =\prod_{i=1}^{r-1}|H(\theta_1,\dots,\theta_i):
H(\theta_1,\dots,\theta_i,\theta_{i+1})|=\\
\prod_{i=1}^{r-1}|N(\theta_1,\dots,\theta_i):
N(\theta_1,\dots,\theta_i,\theta_{i+1})|= 
|N(\theta_1):N(\theta_1,\dots,\theta_r)|.
\end{multline*}
Thus the lemma holds.
\end{proof}

The following is similar to Lemma \ref{pq51}.
\begin{lemma}\mylabel{pq50}
Assume a finite group $N$ is the semidirect product $N = P\ltimes H$ of its
$\pi$-Hall 
subgroup $P$ with its normal $\pi'$-Hall subgroup $H$.
Assume further that $\pw$ is a subgroup of $P$,  and that
$\mathcal{P}$ is any $\pi$-group of automorphisms 
of $N$ that normalizes $\pw$.
Then there exists $t \in H$ such that the following conditions are satisfied: 
\begin{itemize}
\item[(i)] $\mathcal{P}$ normalizes $P^t$ 
\item[(ii)] $\pw \leq P^t$
\item[(iii)]$t$ centralizes $\pw$.
\end{itemize}
\end{lemma}

\begin{proof}
Since $\map$ normalizes $N$, we can form the external semi--direct product product
 $N\map= N \rtimes \map$. Furthermore, as 
$\map$ normalizes $\pw$, the
 $\pi$-group $\pw \mathcal{P}$ is a subgroup 
of $N\map$, and thus  normalizes $N$.
 Hence $\pw \mathcal{P}$ normalizes a $\pi$-Hall
subgroup $P^t$ of $N$, for some  $t\in H$. Therefore $\mathcal{P}$,
 as well as
$\pw$, normalizes $P^t$.
But $\pw$ is a subgroup of $N$. so   the only way it  can normalize
 the Hall subgroup $P^t$ is to be contained 
in $P^t$.   
Or,  equivalently,  $(\pw)^{t^{-1}} \leq P$.
 Thus $s^{-1}s^{t^{-1}}\in P$ whenever $s \in \pw$.
But $\pw \leq P$ normalizes $H$. Hence   $s^{-1}s^{t^{-1}}=[s,t^{-1}]$ is also an
element of $H$. As $H\cap P=1$, we conclude that $[s,t^{-1}]=1$, 
for all $s \in \pw$. 
This implies that $t^{-1}$, and thus $t$,  centralizes $\pw$.  So
the lemma holds.  
\end{proof}

\begin{lemma}\mylabel{hat:l1}
Let $\pw_{2i}$ be a $\pi$-Hall subgroup of $N(P_{2k}^*,
\widehat{Q}_{2i-1} \tin G_{2i}(\alpha_{2i-2}^*))$, where $i= 1 , \cdots , k$.
Then 
\begin{equation}\mylabel{pq52}
N(P_{2k}^*, \widehat{Q}_{2i-1} \tin G_{2i}(\alpha_{2i-2}^*))  =
\pw_{2i}\ltimes
\widehat{Q}_{2i-1}.
\end{equation} 
Furthermore,
\begin{equation}\mylabel{pq53}
N(P_{2k}^*\tin G_{2i}(\alpha_{2i-2}^*,\hat{\beta}_1, \cdots \hat{\beta}_{2i-1}))
 =N( \widehat{Q}_{2i-1} \tin N(P_{2k}^* \tin G_{2i}(\alpha_{2i-2}^*,\beta_{2i-1,2k} ))),
\end{equation}
whenever $i=1,\dots,k$.
\end{lemma}

\begin{proof}
 According to Proposition \ref{pqrm}
the group $\widehat{Q}_{2i-1}$
is  a $\pi'$-Hall subgroup of $N(P_{2k}^* \tin G_{2i}(\alpha_{2i-2}^*))$. 
Hence
\begin{equation*}
N(P_{2k}^*, \widehat{Q}_{2i-1} \tin G_{2i}(\alpha_{2i-2}^*))  =
\pw_{2i}\ltimes
\widehat{Q}_{2i-1}.
\end{equation*}

For the second part of the lemma
 note that,  for $i=1$  we obviously have $N(P_{2k}^* \tin G_2(\alpha_0^*, \hb_1)) =
 N(\qw_1 \tin N(P_{2k}^* \tin G_2(\alpha_0^*, \beta_{1,2k})))$, as $\alpha_0^*=1$ and  
$\hb_1 = \beta_{1,2k}$,  by Proposition  \ref{pq40x}.

For any $t,i$ with $1\leq t \leq i $, 
the group $N( \widehat{Q}_{2i-1} \tin N(P_{2k}^* \tin G_{2i}(\alpha_{2i-2}^*,
\beta_{2i-1,2k} )))$ normalizes $\qw_{2t-1}=\qw_{2i-1} \cap G_{2t-1}$ 
and fixes $\beta_{2t-1,2k}$,
 as $\beta_{2t-1,2k}$ is the unique character
of $Q_{2t-1,2k}=Q_{2i-1,2k} \cap G_{2t-1}$ that lies under $\beta_{2i-1,2k}$.
So  $N( \widehat{Q}_{2i-1} \tin N(P_{2k}^*
 \tin G_{2i}(\alpha_{2i-2}^*, \beta_{2i-1,2k} )))$ fixes the characters
$\hb_1, \hb_3,\dots,\hb_{2i-1}$ by Corollary \ref{hat:r1}. Hence
$$
N( \widehat{Q}_{2i-1} \tin N(P_{2k}^* \tin G_{2i}(\alpha_{2i-2}^*,
\beta_{2i-1,2k} )))\leq 
N(P_{2k}^*\tin G_{2i}(\alpha_{2i-2}^*,\hat{\beta}_1, \dots \hat{\beta}_{2i-1})).
$$

For the other inclusion, we use  Corollary \ref{hat:r1} again, 
but in a recursive way.
We saw above 
 that the group  $N(P_{2k}^* \tin G_{2i}(\hb_1))= 
N(P_{2k}^* \tin G_{2i}(\alpha_{0}^*, \hb_1))$  fixes $\beta_{1,2k}$.
Hence, $N(P_{2k}^* \tin G_{2i}(\hb_1, \hb_3))$ normalizes the group 
$\qw_{3} (\beta_{1,2k})$. In view of \eqref{pq40f} the last group equals 
$Q_{3,2k}$. Thus Corollary \ref{hat:r1}
implies that $N(P_{2k}^* \tin G_2(\hb_1, \hb_3))$ fixes $\beta_{3,2k}$,
as it fixes $\hb_3$ and $\beta_{1,2k}$,  and normalizes $Q_{3,2k}$.
Similarly, we get that $N(P_{2k}^* \tin G_{2i}(\hb_1, \hb_3,\hb_5))$ 
normalizes $\qw_5(\beta_{3,2k}) = Q_{5,2k}$, and  fixes  both 
 $\beta_{3,2k}$ and  $\hb_5$. Thus
it also fixes $\beta_{5,2k}$, by Corollary \ref{hat:r1}.

We continue in this way. So
after $i-1$ steps we have that $N(P_{2k}^* \tin G_{2i}(\hb_1, \hb_3,\dots,\hb_{2i-3}))$ 
fixes $\beta_{2i-3,2k}$. Hence the group 
$N(P_{2k}^* \tin G_2(\hb_1, \hb_3,\dots,\hb_{2i-1}))$ 
normalizes $\qw_{2i-1}(\beta_{2i-3,2k}) = Q_{2i-3,2k}$, by \eqref{pq40f}, and fixes both 
$\beta_{2i-3,2k}$ and 
$\hb_{2i-1}$. Therefore Corollary \ref{hat:r1} implies that it fixes $\beta_{2i-1,2k}$, i.e.,
$$
N(P_{2k}^* \tin G_{2i}(\alpha_{2i-2}^*, \hb_1, \hb_3,\dots,\hb_{2i-1})) 
\leq N( \widehat{Q}_{2i-1} \tin N(P_{2k}^* \tin G_{2i}(\alpha_{2i-2}^*,
\beta_{2i-1,2k} ))).
$$
This completes the proof of the lemma.
\end{proof}
 
Let $i=1,\dots,k$.
According to Theorem \ref{hat:p2} the group
 $N(P_{2k}^* \tin G_{2i}(\alpha_{2i-2}^*, \beta_{2i-1, 2k}))$ equals
$P_{2i}^*(\beta_{2i-1,2k})\rtimes \qw_{2i-1}(\beta_{2i-1,2k})$.
This, along with \eqref{pq52} and \eqref{pq53}, implies that 
\begin{align*}
&N(P_{2k}^* \tin G_{2i}(\alpha_{2i-2}^*,\hat{\beta}_1,
\cdots  ,\hat{\beta}_{2i-1}))  &\quad \\
 &=N(\qw_{2i-1} \tin 
 N(P_{2k}^* \tin G_{2i}(\alpha_{2i-2}^*, \beta_{2i-1,2k})))
&\text{ by \eqref{pq53} } \\
&=N(\qw_{2i-1} \tin P_{2i}^*(\beta_{2i-1,2k})\rtimes \qw_{2i-1}
(\beta_{2i-1,2k})
&\text{ by Theorem \ref{hat:p2} }\\
&= N( \widehat{Q}_{2i-1} \tin P_{2i}^*(\beta_{2i-1, 2k})) \times
\widehat{Q}_{2i-1}(\beta_{2i-1,2k})  &\text{ as $\qw_{2i-1}(\beta_{2i-1,2k}) \leq \qw_{2i-1}$} \\ 
&= N( \widehat{Q}_{2i-1} \tin P_{2i}^*(\beta_{2i-1, 2k})) \times Q_{2i-1,2k}. 
&\text{ by \eqref{pq40g}. }
\end{align*}
Hence
\begin{equation}\mylabel{pq54}
 N(P_{2k}^*  \tin G_{2i}(\alpha_{2i-2}^*,\hat{\beta}_1,
\cdots  ,\hat{\beta}_{2i-1}))=  N( \widehat{Q}_{2i-1} \tin P_{2i}^*(\beta_{2i-1, 2k})) \times Q_{2i-1,2k}, 
\end{equation}
for all $i=1,\dots, k$.

The following proposition implies the existence of a ``good'' family of
groups $\pw_{2i}$,
that permit us to use Theorem \ref{cc:p1}  on the groups \
 $\widehat{Q}_{2i-1}, \qw$ and 
$\pw_{2i}$,  for  $i=1, \dots , k$.

\begin{proposition}\mylabel{pq55}
There exist $\pi$-groups $\pw_{2i}$, 
for  $i= 1, \dots , k$, such that the following 
conditions are satisfied:
\begin{subequations} \label{Eqs}
\begin{gather} \pw_{2i} \in 
\Hall_{\pi}(N(P_{2k}^*,\qw_{2i-1} \tin 
G_{2i}(\alpha_{2i-2}^*))), \mylabel{EqsA} \\
\pw_{2i}(\hat{\beta}_1, \hat{\beta}_3, \dots 
,\hat{\beta}_{2i-1}) = N(\qw_{2i-1} \tin 
P_{2i}^*(\beta_{2i-1,2k}))  \notag\\
\in \Hall_{\pi}(N(P_{2k}^*,\qw_{2i-1} \tin 
G_{2i}(\alpha_{2i-2}^*, \hat{\beta}_1, \hat{\beta}_3, \dots, 
\hat{\beta}_{2i-1}))),  \mylabel{EqsB} \\
\text{ and } \pw_{2i} \text{ normalizes } \pw_{2j},  \mylabel{EqsC}
\end{gather}
\end{subequations}
for all $i = 1,2,\dots,k$ and all $j = 1,2,\dots,i$.  Furthermore, 
any such $\pw_{2i}$ satisfy
\begin{multline} \mylabel{ExtraEq}
\pw_{2i}(\hat{\beta}_1, \hat{\beta}_3, \dots 
,\hat{\beta}_{2i-1}) = \pw_{2i}(\hat{\beta}_{2j-1}, 
\hat{\beta}_{2j+1}, \dots, \hat{\beta}_{2i-1}) =\pw_{2i}(\hat{\beta}_{2i-1})  \\ 
\in \Hall_{\pi}(N(P_{2k}^*,\qw_{2i-1} \tin G_{2i}(\alpha_{2i-2}^*, 
\hat{\beta}_{2j-1}, \hat{\beta}_{2j+1}, \dots, \hat{\beta}_{2i-1})))
\end{multline}
whenever $1 \le j \le i \le k$.
\end{proposition}

\begin{proof} 
In view of \eqref{pq54} the only Hall 
$\pi$-subgroup of $N(P_{2k}^*, \qw_{2i-1} \tin 
G_{2i}(\alpha_{2i-2}^*, \hat{\beta}_1, \hat{\beta}_3, \dots, 
\hat{\beta}_{2i-1}))$ is $N(\qw_{2i-1} \tin 
P_{2i}^*(\beta_{2i-1,2k}))$.
As $N(P_{2k}^*,  \qw_{2i-1}   \tin 
G_{2i}(\alpha_{2i-2}^*, \hat{\beta}_1, \hat{\beta}_3, \dots, 
\hat{\beta}_{2i-1}))$ is a subgroup of $N(P_{2k}^*, \qw_{2i-1} \tin G_{2i}(\alpha_{2i-2}^*))$, 
 we can certainly find 
$\pw_{2i}$ satisfying (\ref{Eqs}a,b) for $i = 1,2,\dots,k$.  We 
shall modify these $\pw_{2i}$ so as to obtain new  subgroups 
satisfying \eqref{EqsC} as well as (\ref{Eqs}a,b).

We first note that,   whenever $1 \le t \le i \le k$,  the subgroup 
$P_{2i}^*$ normalizes $P_{2t}^*$, while $P_{2i}^*(\beta_{2i-1,2k})$
 fixes    $\beta_{2t-1,2k}$ by Proposition \ref{pqremark1'}.
As $\qw_{2t-1} = \qw_{2i-1} \cap G_{2t-1}$, we get that  
 $N(\qw_{2i-1} \tin 
P_{2i}^*(\beta_{2i-1,2k}))$ normalizes $N(\qw_{2t-1} \tin 
P_{2t}^*(\beta_{2t-1,2k}))$.  
 By \eqref{EqsB} this is 
equivalent to 
\begin{equation} \mylabel{Eqs1} 
\pw_{2i}(\hat{\beta}_1, 
\hat{\beta}_3, \dots, \hat{\beta}_{2i-1}) \text{ normalizes }  
\pw_{2t}(\hat{\beta}_1, \hat{\beta}_3, \dots, \hat{\beta}_{2t-1})
\end{equation}
whenever $1 \le t \le i \le k$.

By \eqref{pq52} we have  $N(P_{2k}^*, \widehat{Q}_{2i-1} \tin G_{2i}(\alpha_{2i-2}^*))
= \pw_{2i} \ltimes
\widehat{Q}_{2i-1}$,  for each $i=1,\dots,k$. 
Also  $ N(P_{2k}^*, \qw_{2i-1} \tin G_{2i}(\alpha_{2i-2}^*))$ 
normalizes  both $P_{2t}^* =P_{2k}^* \cap G_{2t}$
 (by Proposition \ref{pq3-})
 and $\qw_{2t-1}=\qw_{2i-1} \cap G_{2t-1}$
 (by Definition \ref{pq40a}), whenever $1\leq t \leq i \leq k$. 
So  it fixes the unique 
character $\alpha_{2t-2}^* \in \Irr(P_{2t-2}^*)$ lying under 
$\alpha_{2i-2}^*$ (see Proposition \ref{p*13}). 
Hence  it  normalizes  $G_{2t}(\alpha_{2t-2}^*)$.
 So  $N(P_{2k}^*, \widehat{Q}_{2i-1} \tin G_{2i}(\alpha_{2i-2}^*))$ normalizes 
 $N(P_{2k}^*, \widehat{Q}_{2t-1} \tin G_{2t}(\alpha_{2t-2}^*))$, i.e., 
\begin{equation} \mylabel{Eqs2} 
\pw_{2i} \ltimes \qw_{2i-1} 
\text{ normalizes } \pw_{2t} \ltimes \qw_{2t-1},
\end{equation}
whenever $1\leq t \leq i \leq k$.

We are going  to modify the $\hat P_{2i}$ so 
that they satisfy
\begin{equation} \mylabel{pq58ii}
 \pw_{2i}(\hat{\beta}_1, 
\hat{\beta}_3,\dots, \hat{\beta}_{2i-1}) \text{ normalizes }
 \pw_{2t}\end{equation}
whenever $1 \le t \le i \le k$, as well as (\ref{Eqs}a,b).
For this we will use reverse induction on $t$, starting with $t=k$ and working 
down.
The group $\pw_{2k}$ requires no 
modification, since the only possible $i$ satisfying $k \le i \le 
k$ is $i = k$, and the subgroup $\pw_{2k}(\hat{\beta}_1, 
\hat{\beta}_3, \dots, \hat{\beta}_{2k-1})$ certainly normalizes 
$\pw_{2k}$.
For the inductive step 
 assume that $\pw_{2k}$, $\pw_{2k-2}$, 
\dots, $\pw_{2s+2}$, for some integer $s = 1,2, \dots,k-1$, have 
already been modified so that \eqref{pq58ii} holds whenever $s < t 
\le i \le k$, and (\ref{Eqs}a,b) hold for all $i = 1,2,\dots,k$.  
We  want to modify $\pw_{2s}$ so that (\ref{Eqs}a,b) still hold 
for $i = s$, while \eqref{pq58ii} with $t = s$ holds for all $i = 
s,s+1,\dots,k$. 

According to \eqref{Eqs1}  the product 
$$
T_{t+1} = \pw_{2k}(\hat{\beta}_1, \dots ,\hat{\beta}_{2k-1}) \cdot
\pw_{2k-2}(\hat{\beta}_1, \dots ,\hat{\beta}_{2k-3}) \cdots
\pw_{2t+2}(\hat{\beta}_1, \dots ,\hat{\beta}_{2t+1})
$$
forms a $\pi$-group, whenever $1\leq t \leq k-1$.
Each factor $\pw_{2i}(\hb_1,\dots,\hb_{2i-1})$, for $i=s+1,s+2,\dots,k$, 
in this product  
is contained in $\pw_{2i}\ltimes \qw_{2i-1}$, and hence normalizes 
  $\pw_{2s} \ltimes \qw_{2s-1}$ by \eqref{Eqs2}.
That factor also normalizes $\pw_{2s}(\hat{\beta}_1, \hat{\beta}_3, \dots 
,\hat{\beta}_{2s-1})$ by \eqref{Eqs1}.
Thus $T_{s+1}$ acts on $\pw_{2s} \ltimes \qw_{2s-1}$ and 
normalizes the subgroup $\pw_{2s}(\hat{\beta}_1, \hat{\beta}_3, \dots 
,\hat{\beta}_{2s-1})$ of $\pw_{2s}$.
Therefore  we can apply Lemma \ref{pq50} to the groups $T_{s+1},
\pw_{2s}(\hat{\beta}_1, \dots ,\hat{\beta}_{2s-1})$
and $\pw_{2s}\ltimes \widehat{Q}_{2s-1}$ to get an element $t \in
\widehat{Q}_{2s-1}$
such that 
\begin{subequations}\mylabel{pq57}
 \begin{equation}\mylabel{pq57i}
(\pw_{2s})^t \in \Hall_{\pi}(\pw_{2s} \ltimes 
\widehat{Q}_{2s-1})
\text{  is normalized by }  T_{s+1}
 \end{equation}
and
\begin{equation} \mylabel{pq57ii}
\pw_{2s}(\hat{\beta}_1, \dots ,\hat{\beta}_{2s-1}) \leq
(\pw_{2s})^{t}.
\end{equation}
\end{subequations}
Obviously the group $(\pw_{2s})^t$ satisfies \eqref{EqsA} for $i=s$,  as $\qw_{2s-1} \leq 
N(P_{2k}^*, \qw_{2s-1} \tin G_{2s}(\alpha_{2s-2}^*))$. Furthermore, 
as $\pw_{2s}(\beta_1,\dots,\beta_{2s-1})$ satisfies \eqref{EqsB} for $i=s$, 
the inclusion \eqref{pq57ii} implies that 
\begin{equation}\mylabel{pq58i} 
(\pw_{2s})^{t}(\hat{\beta}_1, \dots ,\hat{\beta}_{2s-1}) =
\pw_{2s}(\hat{\beta}_1, \dots ,\hat{\beta}_{2s-1}).
\end{equation}
Therefore, $(\pw_{2s})^{t}$ satisfies (\ref{Eqs}a,b) for $i=s$.
 In addition,    \eqref{pq57i} implies that 
$ 
\pw_{2i}(\hat{\beta}_1, \dots ,\hat{\beta}_{2i-1})$
 \text{ normalizes } $(\pw_{2s})^t $ 
whenever $i=s+1, \dots, k$.
  Hence we can work with $(\pw_{2s})^{t}$ in the place of
$\pw_{2s}$.
The new  group $\pw_{2i}$ 
satisfies \eqref{pq58ii} whenever $s\leq t \leq i \leq k$.

At  this point we have shown that we can find the groups $\pw_{2i}$ that 
satisfy (\ref{Eqs}a,b) as well as \eqref{pq58ii} whenever $1\leq t\leq i \leq k$.                    
 These groups can be modified further to satisfy, in addition, \eqref{EqsC}.
Indeed, according to \eqref{Eqs2}, we get that 
  $\pw_{2k} \leq \pw_{2k} \ltimes \widehat{Q}_{2k-3}$ normalizes
$\pw_{2k-2}\ltimes
 \widehat{Q}_{2k-3}$ while $\pw_{2k}(\hat{\beta}_1, \dots ,\hat{\beta}_{2k-1})$ 
normalizes $\pw_{2k-2}$  by \eqref{pq58ii}.
Therefore, Lemma \ref{pq51}, with $\pi'$ in the place of $\pi$, 
implies that there exists $t_{2k-3} \in \widehat{Q}_{2k-3}$ such that 
$$\pw_{2k}^{t_{2k-3}} \text{ normalizes } \pw_{2k-2}$$
and
$$\pw_{2k}(\hat{\beta}_1, \dots ,\hat{\beta}_{2k-1})\leq
\pw_{2k}^{t_{2k-3}}.$$
In view of \eqref{EqsB} (for $\pw_{2k}$ ) this  inclusion 
implies that $$\pw_{2k}(\hat{\beta}_1, \dots ,\hat{\beta}_{2k-1}) =
 \pw_{2k}^{t_{2k-3}}(\hat{\beta}_1, \dots ,\hat{\beta}_{2k-1}).$$
The above equation permits us to work with $\pw_{2k}^{t_{2k-3}}$  in
the place 
of $\pw_{2k}$. Then  $\pw_{2k}$   satisfies (\ref{Eqs}a,b) and 
\eqref{pq58ii},  and   normalizes $\pw_{2k-2}$.

Now   the product $\pw_{2k}\pw_{2k-2}$ forms a group that
normalizes  the unique $\pi'$-Hall subgroup  $\qw_{2k-5}$ 
  of $\pw_{2k-4} \ltimes \qw_{2k-5}$, 
 and has $T_{k-1}$ as a subgroup.                                                                  
Furthermore, $\pw_{2k}\pw_{2k-2} \ltimes \widehat{Q}_{2k-5}$
normalizes $\pw_{2k-4} \ltimes \widehat{Q}_{2k-5}$,  
while $T_{k-1}$ normalizes $\pw_{2k-4}$, as
$\pw_{2k}(\hat{\beta}_1,\-  \dots \- ,
\hat{\beta}_{2k-1})$ and  
$\pw_{2k-2}(\hat{\beta}_1, \dots ,\hat{\beta}_{2k-3})$ do by \eqref{pq58ii}.
Hence Lemma \ref{pq51},  with $\pi'$ in the place of $\pi$, 
implies that we can find an element $t_{2k-5} \in
\widehat{Q}_{2k-5}$ such that 
$$(\pw_{2k} \pw_{2k-2})^{t_{2k-5}} \text{ normalizes }
\pw_{2k-4}$$
and 
$$T_{k-1} \leq \pw_{2k}^{t_{2k-5}} \pw_{2k-2}^{t_{2k-5}}.$$
Also $t_{2k-5}$ centralizes $T_{k-1}$, and thus  centralizes both 
$\pw_{2k}(\hat{\beta}_1, \dots ,\hat{\beta}_{2k-1})$ 
and $\pw_{2k-2}(\hat{\beta}_1, \dots ,\hat{\beta}_{2k-3}).$
Therefore $\pw_{2k}(\hat{\beta}_1, \dots ,\hat{\beta}_{2k-1})$ and
$\pw_{2k-2}(\hat{\beta}_1, \dots ,\hat{\beta}_{2k-3})$
are subgroups of $\pw_{2k}^{t_{2k-5}}$ and $\pw_{2k-2}^{t_{2k-
5}}$ respectively.  So, in view of \eqref{EqsB},  we get that 
$$\pw_{2k}(\hat{\beta}_1, \dots,\hat{\beta}_{2k-1}) =
\pw_{2k}^{t_{2k-5}}(\hat{\beta}_1, \dots,\hat{\beta}_{2k-1})$$
and 
$$\pw_{2k-2}(\hat{\beta}_1, \dots,\hat{\beta}_{2k-3}) =
\pw_{2k-2}^{t_{2k-5}}(\hat{\beta}_1, \dots,\hat{\beta}_{2k-3}).$$
Hence we can replace $\pw_{2k}$ and $\pw_{2k-2}$ by 
$\pw_{2k}^{t_{2k-5}}$ and $\pw_{2k-2}^{t_{2k-5}}$
respectively.  Then    (\ref{Eqs}a,b) and  \eqref{pq58ii}
are satisfied by  the newly modified groups $\pw_{2k}$ and $\pw_{2k-2}$.
Furthermore,  $\pw_{2k}$ normalizes  both $\pw_{2k-2}$ and
$\pw_{2k-4}$,  
while $\pw_{2k-2}$ normalizes $\pw_{2k-4}$.

We continue similarly. At every step the product
$\pw_{2k}\pw_{2k-2}\cdots \pw_{2t}$
of the modified  groups $\pw_{2i}$, for $2\leq t\leq i \leq k$,
contains $T_{t}$. Even more,  $T_t$ is  a subgroup that 
normalizes $\pw_{2t-2}$ while $\pw_{2k}\pw_{2k-
2}\dots \pw_{2t} \ltimes \widehat{Q}_{2t-3}$
normalizes $\pw_{2t-2}\ltimes \widehat{Q}_{2t-3}$.
So Lemma \ref{pq51} implies the existence of an element  $t_{2t-3} \in
\widehat{Q}_{2t-3}$ such that, 
if we replace  $\pw_{2k}$ and $ \pw_{2t}$  by 
$\pw_{2k}^{t_{2k-3}}$ and
 $\pw_{2t}^{t_{2k-3}}$ respectively,  the new groups
satisfy (\ref{Eqs}a,b) and  \eqref{pq58ii} for $i=k,k-1,\dots,t-1$,  while 
$\pw_{2i}$ normalizes $\pw_{2j}$ whenever $k \geq i \geq j \geq t-1$.

This process stops when we reach $t=2$. This proves that we can pick 
the groups $\pw_{2i}$ to satisfy (\ref{Eqs}a,b,c). 
The additional property 
\eqref{ExtraEq} follows from \eqref{EqsB} and 
Lemma \ref{pq49}. 
This completes the proof of Proposition  \ref{pq55}.
\end{proof}

A useful property of the groups $\pw_{2i}$ is given in 
\begin{corollary}\mylabel{b.cor50}
For every $i=1,\dots,k$ we have that 
\begin{subequations}
\begin{equation}\mylabel{b.c51}
\pw_{2i}(\hb_{2i-1})=\pw_{2i}(\hb_1,\dots,\hb_{2i-1})=C(\qw_{2i-1} \tin \pw_{2i})  \text{ and }
\end{equation}
\begin{equation}\mylabel{b.c52}
\pw_{2i}(\hb_{2i-1})= C(\qw_{2i-1} \tin P_{2i}^*). 
\end{equation}
\end{subequations}
\end{corollary}

\begin{proof}
According to \eqref{ExtraEq} and \eqref{EqsB} we get that 
$$\pw_{2i}(\hb_{2i-1}) = \pw_{2i}(\hb_1,\dots,\hb_{2i-1}) = 
N(\qw_{2i-1} \tin P_{2i}^*(\beta_{2i-1,2k})).
$$
In view of   Corollary \ref{hat:r1},  the group 
$N(\qw_{2i-1} \tin P_{2i}^*(\beta_{2i-1,2k}))$ fixes $\hb_{2i-1}$, as it normalizes 
$\qw_{2i-1}$ and fixes $\beta_{2i-1,2k}$. Thus 
\begin{equation}\mylabel{n.e11}
\pw_{2i}(\hb_{2i-1}) = N(\qw_{2i-1} \tin P_{2i}^*(\beta_{2i-1,2k})) \leq 
N(\qw_{2i-1} \tin P_{2i}^*(\hb_{2i-1})).
\end{equation}
But $\qw_{2i-1} \leq \qw \leq G' =G(\alpha_{2k}^*)$.
Thus $\qw_{2i-1}$  normalizes $P_{2i}^*$,  as $\qw$ does. We conclude that 
$\qw_{2i-1}$ normalizes $N(\qw_{2i-1} \tin P_{2i}^*(\hb_{2i-1}))$.
As the latter $p$-group also normalizes the $p'$-group  $\qw_{2i-1}$, we get that 
$\qw_{2i-1}$ centralizes $N(\qw_{2i-1} \tin P_{2i}^*(\hb_{2i-1}))$.
Hence $\qw_{2i-1}$ also centralizes $\pw_{2i}(\hb_{2i-1} ) \leq N(\qw_{2i-1} 
\tin P_{2i}^*(\hb_{2i-1}))$.
Therefore, $\pw_{2i}(\hb_{2i-1}) \leq C(\qw_{2i-1} \tin \pw_{2i})$. As the other inclusion is 
obvious, we conclude that   $\pw_{2i}(\hb_{2i-1}) = C(\qw_{2i-1} \tin \pw_{2i})$.
Hence \eqref{b.c51}  holds

Furthermore, the fact that  $\qw_{2i-1}$ normalizes $P_{2i}^*(\hb_{2i-1})$
implies  that  $N(\qw_{2i-1} \tin P_{2i}^*(\hb_{2i-1})) = C(\qw_{2i-1} \tin P_{2i}^*(\hb_{2i-1}))$.
Hence we have that 
\begin{align*}
\pw_{2i}(\hb_{2i-1}) &= N(\qw_{2i-1} \tin P_{2i}^*(\beta_{2i-1,2k})) &  \\ 
&\leq N(\qw_{2i-1} \tin P_{2i}^*(\hb_{2i-1}))   &\text{ by \eqref{n.e11}  } \\
&=C(\qw_{2i-1} \tin P_{2i}^*(\hb_{2i-1})) &\\
&\leq C(\qw_{2i-1} \tin P_{2i}^*) &  \\
&= C(\qw_{2i-1} \tin P_{2i}^*(\beta_{2i-1,2k}))  &\text{ as $\beta_{2i-1,2k} \in \Irr(Q_{2i-1,2k})$ } \\
&\,   &\text{  and $Q_{2i-1,2k} \leq \qw_{2i-1}$ }\\
&\leq  N(\qw_{2i-1}  \tin P_{2i}^*(\beta_{2i-1,2k})) &\\
 &= \pw_{2i}(\hb_{2i-1}). &
\end{align*}
So $\pw_{2i}(\hb_{2i-1}) = C(\qw_{2i-1} \tin P_{2i}^*)$, 
and \eqref{b.c52} holds.
This completes the proof of the corollary.
\end{proof}

%%% Local Variables: 
%%% mode: latex
%%% TeX-master: t
%%% End: 

%\newcommand{\pns}[1]{(P_{#1}^{\nu})^*}
%\newcommand{\ans}[1]{(\alpha_{#1}^{\nu})^*}

\section{Triangular sets for $G'=G(\alpha_{2k}^*)$}
For the following two sections we are going to keep fixed the groups $P_{2i}^*$, 
along with their characters $\alpha_{2i}^*$. Even though these groups and 
characters come from a specific triangular set of \eqref{*}, we will forget
this triangular set and treat the $P_{2i}^*$ independently. 
Furthermore, we make the following assumption
%\begin{subequations}
%\begin{equation}\mylabel{bh.as1}
%\pi \text{ consists of one prime only, i.e., } \pi= \{ p \},
%\end{equation} 
\begin{equation}\mylabel{bh.as2}
m= 2k  \text{ is even. } 
\end{equation}

We have already seen that the groups $G'_{s}$, 
 defined as $G'_s=G_s(\alpha_{2k}^*)=G'\cap G_{s}$
 whenever $1\leq s \leq m=2k$, form a
 series 
\begin{equation}\mylabel{bh.*}
1 = G'_0 \unlhd G'_1 \unlhd \dots \unlhd G'_{2k-1} \unlhd G'_{2k} \unlhd G'. 
\end{equation}
This is  a series of normal subgroups of $G'$, as the series 
\eqref{*} 
is such for $G$.

\subsection{From $G$ to $G'$}

The goal of this section is to create a triangular set for the series
 \eqref{bh.*}, 
related to the triangular  set 
 $\{Q_{2i-1}, P_{2i}|\beta_{2i-1}, \alpha_{2i} \}_{i=1}^{k}$. 
We remark that the latter is a triangular set for normal series 
\begin{equation} \mylabel{bh.*1}
1 = G_0 \unlhd G_1 \unlhd \dots \unlhd G_{2k-1} \unlhd G_{2k} \unlhd G, 
\end{equation}
that is,   \eqref{*} for $m=2k$. 
As we will see, these two triangular sets are so close related 
that one determines the other uniquely, up to conjugation.

To create such a set, we first  need  to give   
groups and  characters 
that satisfy \eqref{xxx} for the series \eqref{bh.*}.
For the  $\pi'$-groups and characters  we pick  
the groups $Q_{2i-1,2k}$ for every $i=1,\dots,k$, 
 along with their irreducible  characters 
$\beta_{2i-1,2k}$.
In view of Remark \ref{b.rem2},  we have that $Q_{2i-1,2k}$ is a
subgroup of $G'$,  and thus a subgroup of $G'_{2i-1}=G' \cap G_{2i-1}$
 (as $Q_{2i-1,2k} \leq G_{2i-1}$), for all $i=1,\dots,k$.
We  define
\begin{subequations}\mylabel{bh.e2}
\begin{equation}\mylabel{bh.e2a}
Q'_{2i-1}:= Q_{2i-1,2k} \text{ and } \beta'_{2i-1} = \beta_{2i-1,2k},
\end{equation}
and 
\begin{equation}\mylabel{bh.e2b}
P'_0=1 \text{ and } P'_{2i} :=P_{2i}^*(\beta_{2i-1,2k}),
\end{equation} 
whenever  $i=1,\dots,k$.
\end{subequations}
Note that 
\begin{lemma}\mylabel{bh.l1}
For every $i=1,\dots,k$ we have 
\begin{multline*}
P'_{2i} = P^*_{2i}(\beta_{2i-1}')= P_{2i}^*(\beta_{2i-1,2k}) =\\
 N(Q_{2i-1,2k} \tin P_{2i}^*)=N(Q_{2i-1}' \tin P_{2i}^*)=C(Q_{2i-1,2k} \tin P_{2i}^*) =
C(Q_{2i-1}'\tin P_{2i}^*)  =\\
C(Q_{1,2k}, Q_{3,2k}, \dots,Q_{2i-1,2k} \tin P_{2i}^*).
\end{multline*}
Furthermore, 
$$
P_{2i}' \,  \text{is the unique $\pi$-Hall  subgroup of } N(P_{2k}^* \tin 
G_{2i}(\alpha_{2i-2}^*, \beta_{2i-1,2k})).
$$
\end{lemma}

\begin{proof}
Let $i=1,\dots,k$ be fixed.
Obviously $P_{2i}^*(\beta_{2i-1,2k}) \leq N(Q_{2i-1,2k} \tin P_{2i}^*)$.
Also $Q_{2i-1,2k}= C(P_{2i}, \dots, P_{2k} \tin Q_{2i-1})$
 normalizes $P_{2}, \dots,P_{2i-2}$,  as it is a subgroup of $Q_{2i-1}$ 
(see \eqref{x5}),  and centralizes $P_{2i}$. Hence the $\pi'$-group 
 $Q_{2i-1,2k}$ normalizes the $\pi$-group $P_{2i}^*$. Therefore
$$
N(Q_{2i-1,2k} \tin P_{2i}^*) = C(Q_{2i-1,2k} \tin P_{2i}^*) \leq 
P_{2i}^*(\beta_{2i-1,2k}). 
$$
So $$
P'_{2i}= P_{2i}^*(\beta_{2i-1,2k}) 
= N(Q_{2i-1,2k} \tin P_{2i}^*) = C(Q_{2i-1,2k} \tin P_{2i}^*). 
$$
The fact that $C(Q_{2i-1}'\tin P_{2i}^*)  =
C(Q_{1,2k}, Q_{3,2k}, \dots,Q_{2i-1,2k} \tin P_{2i}^*)$ 
 follows easily from the  fact 
(see \eqref{extra}) that $Q_{2t-1,2k}$ is a subgroup (actually normal)
of $Q_{2i-1,2k}$ whenever $1\leq t \leq i$.

The rest of the lemma follows from \eqref{bh.e2a} and  Proposition \ref{pq43}. 
\end{proof}

What about  irreducible characters $\alpha'_{2i}$ of $P_{2i}'$?
Well, there is a straightforward way to pick those characters. 
To see this,  note that  $Q_{2i-1}'=Q_{2i-1,2k}$ fixes the
 character $\alpha_{2i}^*$
of $P_{2i}^*$. Indeed, $Q_{2i-1,2k}$ fixes $\alpha_{2k}^*$ 
(as it is a  subgroup of $G'$) and normalizes $P_{2i}^*=G_{2i}\cap P_{2k}^*$. 
Hence it fixes the unique  character $\alpha_{2i}^*$ of $P_{2i}^*$ that lies
under $\alpha_{2k}^*$ (see Proposition \ref{p*13}). 
As $P_{2i}'=C(Q_{2i-1,2k} \tin P_{2i}^*)=C(Q_{2i-1}' \tin P_{2i}^*)$, we 
can make the 
\begin{defn}\mylabel{bh.d1}
For every  $i=1,\dots,k$, the character 
$\alpha'_{2i} \in \Irr(P_{2i}')$ is   the $Q_{2i-1}'$-Glauberman 
correspondent of $\alpha_{2i}^* \in \Irr(P_{2i}^*)$.
We also set  $\alpha_0':=1 \in \Irr(P_0')$.
\end{defn}

Now we can show
\begin{theorem}\mylabel{bh.th1}
The set  
\begin{equation}\mylabel{bh.set}
\{ Q'_{1}, \dots, Q'_{2k-1},P'_0,  P'_{2}, \dots, P_{2k}'|
\beta'_{1}, \dots,\beta'_{2k-1},\alpha_0', \alpha_{2}', 
\dots,\alpha_{2k}'\}, 
\end{equation}
given  by  \eqref{bh.e2} and Definition \ref{bh.d1},
is a triangular set for \eqref{bh.*}.
\end{theorem}

\begin{proof}
It is enough to check that \eqref{bh.set} satisfies \eqref{xxx}.
It obviously satisfies \eqref{x1}. 
According to \eqref{pqd40a}   we have that   $Q_1' =Q_{1,2k} = G_1'$. Hence 
 \eqref{x2} holds for the set \eqref{bh.set}.

By Lemma \ref{bh.l1}
the group  $P_{2i}'$ centralizes $Q_{2i-1}'=Q_{2i-1,2k}$,
for all $i=1,\dots,k$. Therefore the group that we would write as $Q_{2i-1,2i}'$ (see
\eqref{anno1}) is  $Q_{2i-1}'$ itself.
Furthermore, 
the $P_{2i}'$-Glauberman correspondent of $\beta'_{2i-1}$ 
is the same character  $\beta'_{2i-1}$, whenever $1\leq i \leq k$. 
Therefore the character that we would write as $\beta_{2i-1,2i}'$
(see Definition  \ref{pq1413def})
is nothing else but    the character $\beta_{2i-1}'$. 
According to \eqref{extra} and \eqref{pq14c}   the groups
$Q_{2i-1}'=Q_{2i-1,2k}$ and their characters $\beta_{2i-1}'=\beta_{2i-1,2k}$
line up. That is,  we have a series of normal subgroups 
$$
Q_{1}'\unlhd Q_3'\unlhd \dots \unlhd Q_{2k-3}'\unlhd Q_{2k-1}' 
$$
along with their characters
$$\beta_{1}',  \beta_{3}',  \dots,  \beta_{2k-3}',  \beta_{2k-1}'$$
that lie one under the other. Actually, by Proposition \ref{pqremark1'} 
  the unique character
of $Q_{2j-1}'$ that lies under $\beta_{2i-1}'$ is $\beta_{2j-1}'$,  whenever
 $1\leq j \leq i \leq k$.
Hence the characters $\beta_{2i-1}'$ satisfy  
 \eqref{x6} in the definition of a triangular set.

The group $P_{2i,2i+1}'$,  defined as $P_{2i,2i+1}'=C(Q_{2i+1}' \tin P_{2i}')$
(see \eqref{anno2}),  satisfies 
\begin{multline*}
P_{2i,2i+1}'=C(Q_{2i+1,2k} \tin P_{2i}')=
C(Q_{2i+1,2k} \tin C(Q_{2i-1,2k} \tin P_{2i}^*)) = \\
C(Q_{2i+1,2k} \tin P_{2i}^*) = C(Q_{2i+1}' \tin P_{2i}^*),
\end{multline*}
whenever $1\leq i \leq k-1$.
But 
$Q_{2i+1}' =Q_{2i+1,2k}$ normalizes $P_{2i}^*$,  and fixes its irreducible
character $\alpha_{2i}^*$,  as it fixes $\alpha_{2k}^*$.
Hence the $Q_{2i+1}'$-Glauberman correspondent of $\alpha_{2i}'$ is the 
$Q_{2i+1}'$-Glauberman correspondent of $\alpha_{2i}^*$, as $\alpha_{2i}'$ 
is the $Q_{2i-1}'$-Glauberman correspondent of $\alpha_{2i}^*$,  and 
$Q_{2i-1}' \unlhd  Q_{2i+1}'$.
This implies that the character $\alpha_{2i,2i+1}' \in \Irr(P_{2i,2i+1}')$,
 defined as the 
$Q_{2i+1}'$-Glauberman correspondent of $\alpha_{2i}'$ (see Definition 
\ref{pq1413def2}), is the $Q_{2i+1}'$-Glauberman correspondent of 
$\alpha_{2i}^* \in \Irr(P_{2i}^*)$, for all $i=1,\dots, k-1$. 
Furthermore, $\alpha_{2i+2}'$ is the $Q_{2i+1}'$-Glauberman correspondent of 
$\alpha_{2i+2}^*$. 
As  $\alpha_{2i+2}^*$ lies above $\alpha_{2i}^*$ we conclude that the same
holds for their $Q_{2i+1}'$-Glauberman correspondents. 
Thus $\alpha_{2i+2}' \in \Irr(P_{2i+2}')$ lies above 
$\alpha_{2i,2i+1}' \in \Irr(P_{2i,2i+1}')$ whenever 
$1\leq i \leq k-1$. 
We obviously have that $\alpha_{2}' \in \Irr(P_{2}')$ lies above
 $1=\alpha_{0,1}' \in \Irr(P_0')$. 
Hence the characters $\alpha_{2i}'$ satisfy \eqref{x4}.

To complete the proof of the theorem, it remains to show that the set 
\eqref{bh.set} also satisfies \eqref{x3} and \eqref{x5}.

According to \eqref{pq40c}, \eqref{pq40e}  and \eqref{pq40f} we get that 
\begin{equation}\mylabel{bh.e3}
Q_{2i-1}' = Q_{2i-1,2k} =\qw_{2i-1}(\beta_{2i-3,2k}) \in \Hall_{\pi'}(G_{2i-1}
(\alpha_{2k}^*, \beta_{2i-3,2k})),
\end{equation}
for all $i=1,\dots,k$.
Furthermore, as  $Q_{2i-1}'=Q_{2i-1,2k}$  fixes $\alpha_{2k}^*$ it also fixes 
  $\alpha_{2j}^*$ for all $j=1,\dots,k$,  by Remark \ref{b.rem1}.
In view of Proposition 
\ref{pqremark1'},  it  also fixes  $\beta_{2j-1}'=\beta_{2j-1,2k}$
for all $j=1,\dots,i-1$.  Thus $Q_{2i-1}'$ normalizes the groups $Q_{2j-1}'$ 
and fixes $\alpha_{2j}^*$. This implies that $Q_{2i-1}'$ also fixes the
$Q_{2j-1}'$-Glauberman correspondent $\alpha_{2j}'$ of $\alpha_{2j}^*$ 
whenever $1\leq j \leq i-1$.
Hence $Q_{2i-1}' \leq G_{2i-1}(\alpha_{2k}^*, \beta_1', \dots,\beta_{2i-3}', 
\alpha_2',\dots,\alpha_{2i-2}')$. 
(We actually have even more as $Q_{2i-1}'$  fixes $\alpha_{2i}'$ and
$\beta_{2i-1}'$ for all $i=1,\dots, k$, but we don't need it here.)
This, along with \eqref{bh.e3} and the fact that  
$$
G_{2i-1}(\alpha_{2k}^*,\alpha_2', 
\dots,\alpha_{2i-2}', \beta_{1}', \dots,\beta_{2i-3}')  
\leq G_{2i}(\alpha_{2k}^*, \beta_{2i-3}') = G_{2i}(\alpha_{2k}^*, \beta_{2i-3,2k}), 
$$
 implies  that 
\begin{multline}
Q_{2i-1}' \in \Hall_{\pi'}(G_{2i-1}(\alpha_{2k}^*,\alpha_2', 
\dots,\alpha_{2i-2}', \beta_{1}', \dots,\beta_{2i-3}'))=\\ 
\Hall_{\pi'}(G_{2i-1}'(\alpha_2', \dots,\alpha_{2i-2}', 
\beta_{1}', \dots,\beta_{2i-3}')),
\end{multline}
whenever $i=1,\dots,k$.
Hence  \eqref{x5} holds for the $\pi'$-groups $Q_{2i-1}'$.

As for the $\pi$-groups,  we first note that,  in view of  
Lemma \ref{bh.l1},  for every $i=1,\dots,k$ the 
 group $P_{2i}'$ centralizes 
$Q_1'=Q_{1,2k}, \dots, Q_{2i-1}'=Q_{2i-1,2k}$, and thus fixes their
characters $\beta_1', \dots,\beta_{2i-1}'$. 
 It also fixes the characters $\alpha_{2}^*, \dots, \alpha_{2i-2}^*$, 
as $P_{2i}' \leq P_{2i}^*$.  Therefore it also fixes the
$Q_{2j-1}'$-Glauberman correspondent $\alpha_{2j}'$ of $\alpha_{2j}^*$, 
for all $j=1,\dots,i$. 
Hence 
\begin{multline}\mylabel{bh.e4}
P_{2i}'= P_{2i}^*(\beta_{2i-1}') \leq G_{2i}(\alpha_{2k}^*, \alpha_{2}', 
\dots,\alpha_{2i-2}', \beta_1', \dots,\beta_{2i-1}') \leq \\
N(P_{2k}^* \tin  G_{2i}(\alpha_{2i-2}^*, \beta_{2i-1}')). 
\end{multline}
 Proposition \ref{pq43} implies that  $P_{2i}'=P_{2i}^*(\beta_{2i-1,2k})$ 
is the unique $\pi$-Hall 
 subgroup of the group  $N(P_{2k}^* \tin G_{2i}(\alpha_{2i-2}^*, 
\beta_{2i-1,2k}))$. This, along with \eqref{bh.e4},
implies that 
$P_{2i}' $  is a $\pi$-Hall  subgroup of $G_{2i}(\alpha_{2k}^*, 
\alpha_2',\dots,\alpha_{2i-2}', \beta_1', \dots,\beta_{2i-1}')$, whenever 
$1\leq i \leq k$.
As $G_{2i}(\alpha_{2k}^*) = G_{2i}'$, we conclude that \eqref{x5}
holds  for the groups $P_{2i}'$. 
Hence Theorem \ref{bh.th1} is proved.
\end{proof}

For the triangular set \eqref{bh.set} we can define, as it was described in 
Section \ref{pq:sec3},  the groups
 $(P_{2i}')^*:=P_2' \cdot P_4' \cdots P_{2i}'$,  along with their 
irreducible characters $(\alpha_{2i}')^*$ (see Definition \ref{p*d1}), 
whenever $1\leq i \leq k$.
Then it is easy to show that 
\begin{proposition}\mylabel{bh.pr1}
$$
(P_{2i}')^* = P_{2i}^*,
$$
for every $i=1,\dots,k$.
\end{proposition}

\begin{proof}
In view of  \eqref{bh.e2}, it is clear that $P_{2i}' \leq P_{2i}^*$ for 
all $i=1,\dots,k$. As $P_{2j}^* \leq P_{2i}^*$, whenever $1\leq j \leq i$,  we conclude that 
$(P_{2i}')^*$ is a subgroup of $P_{2i}^*$, whenever $1\leq i \leq k$. 

For the other inclusion,   note that, according to Lemma \ref{bh.l1},
$P_{2i}' =C(Q_{2i-1,2k} \tin P_{2i}^*)$. But $P_{2i}$ is a subgroup of
$P_{2i}^*$  and centralizes $Q_{2i-1,2k}$ (see \eqref{pq14e}).
Hence $P_{2i} \leq C(Q_{2i-1,2k} \tin P_{2i}^*) = P_{2i}'$ whenever 
$1\leq i \leq k$. 
Therefore, 
$$
P_{2i}^* =P_2 \cdot P_4 \cdots P_{2i} \leq P_2' \cdot P_4' \cdots P_{2i}'=
(P_{2i}')^*.     
$$
Hence $P_{2i}^* =( P_{2i}')^*$.
\end{proof}

\subsection{From $G'$ to $G$}
Now assume that  a triangular set for $G'$ is given. 
It would be nice if we could pass to a triangular set of $G$ in a 
``reverse '' way to that described in the previous section. 
This would  not only show  a path to pass from triangular sets of $G$ to $G'$ 
and vice versa, but also, as  Theorem \ref{cor:t} suggests, 
a path to pass from character towers  of
\eqref{bh.*1}  to character towers of \eqref{bh.*} and vice versa. 
We couldn't hope that this would   work with every triangular set of $G'$, as,
after all, the triangular set that we got in the previous section 
has  a very specific type.  
 That type we  try to reproduce in Property \ref{bh.pp1} that follows. 
In addition, we need an extra asumption for the set of primes $\pi$.
We assume that 
\begin{equation}\mylabel{bh.as1}
\pi = \{ p\}  \text{ consists of one prime only. } 
\end{equation} 
So the  various $\pi$-Hall subgroups become $p$-Sylow subgroups, 
while the $\pi'$-Hall become   $p'$-Hall. 

Now assume that  
\begin{subequations}
\begin{equation}\mylabel{bh.e5}
\{Q_{2i-1}', P_{2i}'|\beta_{2i-1}', \alpha_{2i}'\}_{i=1}^{k}
\end{equation}
is   a triangular set for \eqref{bh.*}, while 
 $Q'$ is any  $p'$-subgroup of $G'$ 
 satisfying 
\begin{equation}\mylabel{bh.e5b}
Q_{2i-1}' \unlhd  Q' \leq G'(\alpha_2', \dots, \alpha_{2k}', \beta_1', 
\dots, \beta_{2k-1}'),
\end{equation}
whenever  $1\leq i \leq k$. Note that $Q_{2k-1}'$ works for  $Q'$.
\end{subequations}
Furthermore, we assume that the  set \eqref{bh.e5} satisfies
the following property
\begin{property}\mylabel{bh.pp1}
For every $i=1,\dots,k$ we have 
\begin{subequations}
\begin{equation}\mylabel{bh.e}
P_{2i}'= N(Q_{2i-1}' \tin P_{2i}^*)=C(Q_{2i-1}' \tin P_{2i}^*)= 
P_{2i}^*(\beta_{2i-1}').
\end{equation}
In addition, 
\begin{equation}\mylabel{bh.ee}
P_{2i}' \text{ \, is the unique $p$-Sylow subgroup of  } N(P_{2k}^* \tin 
G_{2i}(\alpha_{2i-2}^*, \beta_{2i-1}')). 
\end{equation}
Furthermore,  
\begin{equation}\mylabel{bh.eee}
\alpha_{2i}' \in \Irr(P_{2i}')  \text{ is the $Q_{2i-1}'$-Glauberman
correspondent of }  \alpha_{2i}^* \in \Irr(P_{2i}^*).
\end{equation}
\end{subequations}
\end{property}
We remark that,  as  the $p'$-group
$Q_{2i}'$ normalizes the $p$-group $P_{2i}^*=P_{2k}^* \cap G_{2i}$,
we necessarily have that 
$$
N(Q_{2i-1}' \tin P_{2i}^*)=C(Q_{2i-1}' \tin P_{2i}^*)= 
P_{2i}^*(\beta_{2i-1}'),
$$
for all $i=1,\dots,k$. Thus equation \eqref{bh.e} is equivalent to 
$P_{2i}' = P_{2i}^*(\beta_{2i-1}')$. 

Lets see some of the conditions Property \ref{bh.pp1} implies for
 the triangular set \eqref{bh.e5}. Recall that, for all $i=1, \dots, k$, the groups
$Q_{2i-1,2k}'$ are defined as 
$Q_{2i-1,2k}' = C(P_{2i}', \dots, P_{2k}' \tin Q_{2i-1}' ) = 
C(P_{2i}' \cdots P_{2k}' \tin Q_{2i-1}')$ (see \eqref{pq14e}).
Furthermore, for all $i=1, \dots, k$, the character $\beta_{2i-1, 2k}'$ is the
$P_{2i}'\cdots P_{2k}'$-Glauberman correspondent of 
$\beta_{2i-1,2k}'$, by Definition \ref{pq1413def2}.

\begin{lemma}\mylabel{bh.l2}
For every $i=1,\dots,k$ we have that $Q_{2i-1,2k}' = Q_{2i-1}'$ while 
 $\beta_{2i-1,2k}' = \beta_{2i-1}'$.
 Therefore we get 
\begin{equation}\mylabel{bh.e5.1}
\begin{aligned}
 Q_1'\unlhd  Q_3' \unlhd   &\dots \unlhd  Q_{2k-1}' \unlhd  Q'  \text{while } \\
\beta_{2j-1}' &\in \Irr(Q_{2j-1}'  |\beta_{2j-3}' ) \text{\,  whenever  $2\leq j \leq k$ }.
\end{aligned}
\end{equation}

Furthermore, whenever $1\leq j \leq i \leq k$  we get 
 \begin{subequations}\mylabel{bh.e6}
\begin{equation}\mylabel{bh.e6b}
Q_{2i-1}'(\beta_{2j-1}') = Q_{2i-1}', 
\end{equation}
while 
\begin{align}
Q_1' =G_1'= G_1(\alpha_{2k}^*)&=N(P_{2k}^* \tin G_1) = C(P_{2k}^* \tin G_1), \mylabel{bh.e6a} \\
P_{2}' &= P_{2}^*=P_2, \mylabel{bh.e6c}  \\
P_{2i}' =P_{2i}^*(\beta_{2i-1}')= P_{2i}^*(\beta_1', 
\dots,\beta_{2i-1}') &=  C(Q_1',\dots,Q_{2i-1}' \tin P_{2i}^*).
 \mylabel{bh.e6d}
\end{align}
\end{subequations}
\end{lemma}

\begin{proof}
According to Property \ref{bh.pp1} the group 
 $P_{2i}'$ centralizes $Q_{2i-1}'$, for each $i=1, \dots, k$.. As 
$Q_{2i-1,2i}' :=C(P_{2i}' \tin Q_{2i-1}')$ (see \eqref{pq14e}),
we conclude that  $Q_{2i-1,2i}'= Q_{2i-1}'$. 
But,  according to \eqref{pq:norm}, the group $Q_{2i-1,2i}'$ is a normal 
subgroup of $Q_{2i+1}'$ whenever $1\leq i \leq k-1$.
 Thus $Q_{2i-1}' \unlhd Q_{2i+1}'$ for all such $i$ and the 
first part of \eqref{bh.e5.1} is proved.

As $P_{2i}'$ centralizes $Q_{2i-1}'$ and $Q_{2j-1}' \unlhd Q_{2i-1}'$,
we conclude that $P_{2i}'$ centralizes $Q_{2j-1}'$ whenever $1\leq j\leq i
\leq k$. Thus $P_{2i}'=C(Q_{2i-1}' \tin P_{2i}^*) = C(Q_1', \dots,
Q_{2i-1}'\tin P_{2i}^*)$. 
Even more, as 
 $Q_{2j-1,2i}' = C(P_{2j}', \dots, P_{2i}' \tin Q_{2j-1}')$
(see \eqref{pq14e}), we get that  $Q_{2j-1,2i}'= Q_{2j-1}'$ whenever $1\leq j
\leq i \leq k$.  Thus $\beta_{2j-1,2i}'= \beta_{2j-1}'$. 
This, along with \eqref{pq14c},  implies that  $\beta_{2j-1}'$ 
 lies above $\beta_{2j-3}'$ whenever $2\leq j\leq k$. Hence the rest of 
 \eqref{bh.e5.1}  holds.

Furthermore, as  $\beta_{2j-1,2k}'= \beta_{2j-1}'$, 
 Proposition \ref{pqremark1'},  for  the $t, j , i $ there equal to  $i , k, j $ here,  
implies that  $Q_{2i-1,2k}'(\beta_{2j-1}') = Q_{2i-1,2k}'$. 
Equation  \eqref{bh.e6b} holds, as $Q_{2i-1,2k}'= Q_{2i-1}'$.

The set \eqref{bh.e5} is a triangular set of \eqref{bh.*}. 
Hence,
(see \eqref{x2}), 
$$Q_{1}' = G_{1}' =G_1(\alpha_{2k}^*).$$
Therefore, $Q_{1}' =N(P_{2k}^* \tin G_1) =C(P_{2k}^* \tin G_1)$, as the 
$p$-group $P_{2k}^*$ normalizes the $p'$-group $G_1$. 
So \eqref{bh.e6a} holds.
Furthermore, we get that   
 $P_2 =P_2^*$ centralizes $Q_1'$.  This implies that   
$P_{2}' = C(Q_{1}' \tin P_{2}^*) = P_{2}^*=P_2$.
Thus  \eqref{bh.e6c}  holds.

It remains to show that  \eqref{bh.e6d} holds.
As $P_{2i}'$ centralizes $Q_1', \dots,Q_{2i-1}'$,  and is a subgroup of 
$P_{2i}^*$,  we obviously have 
that $P_{2i}' \leq C(Q_1' ,\dots, Q_{2i-1}' \tin P_{2i}^*)
 \leq  P_{2i}^*(\beta_1', \dots,\beta_{2i-1}') \leq 
P_{2i}^*(\beta_{2i-1}')$. 
But  $P_{2i}^*(\beta_{2i-1}') = P_{2i}'$. 
This completes the proof of the lemma. 
\end{proof}

For each $i=1,\dots,k$, let 
  $\pps{2i}$ be  the product group $\pps{2i}=P_{2}' \cdots P_{2i}'$, 
and $\aas{2i}$ its  irreducible character  that we get 
(see Definition \ref{p*d1}) from the triangular set \eqref{bh.e5}. 
Then 
 \begin{lemma} \mylabel{bh.l3}
For every $i=1,\dots,k$,
$$\pps{2i} =P_{2i}^*.$$
\end{lemma}

\begin{proof}
We will use induction on $i$. Equation  \eqref{bh.e6c} verifies the $i=1$ 
case.
Assume the lemma  is true for all $i=1,\dots,n-1$,  where $n=2, 3,\dots,k$. 
We will prove it also holds for $i=n$. 
The group $Q_{2n-1}'$ normalizes $P_{2i}^*$ for all $i=1,\dots,k$
(as it normalizes $P_{2k}^*$).
Thus  $Q_{2n-1}' \ltimes P_{2n}^*$ is a group. Furthermore, 
the semi--direct product $Q_{2n-1}'\ltimes P_{2n-2}^*$ is a normal subgroup of 
 $Q_{2n-1}' \ltimes P_{2n}^*$, as  
$Q_{2n-1}'\ltimes P_{2n-2}^* =G_{2n-1} \cap (Q_{2n-1}' \ltimes P_{2n}^*)$. 
Hence Frattini's argument implies that 
$$ 
P_{2n}^* = N(Q_{2n-1}' \tin P_{2n}^* ) \cdot P_{2n-2}^*.
$$
According to the inductive hypothesis $P_{2n-2}^* = \pps{2n-2}$. 
 Even more,  in view of \eqref{bh.e} we get that  $P_{2n}'
 = N(Q_{2n-1}' \tin P_{2n}^*)$. 
So we conclude that 
$$
P_{2n}^* = N(Q_{2n-1}' \tin P_{2n}^* ) \cdot P_{2n-2}^* = P_{2n}' \pps{2n-2}
=\pps{2n}.
$$
This completes the proof of the inductive step. Thus
Lemma \ref{bh.l3} holds.
\end{proof}

We can now prove
\begin{theorem}\mylabel{bh.th2}
Assume that the triangular set \eqref{bh.e5}  for the series 
\eqref{bh.*} satisfies Property \ref{bh.pp1}, 
while \eqref{bh.e5b} holds for a subgroup $Q'$ 
of $G'$.
Then there exists a collection of groups and characters 
\begin{equation}\mylabel{bh.e7}
\{ \pn_0 , \qn_{2i-1}, \pn_{2i} | \, \an_0,  \bn_{2i-1}, \an_{2i}  \}_{i=1}^k,
\end{equation}
with the following properties:
\begin{subequations}\mylabel{bh.al}
\begin{align}
Q' &\text{normalizes } \qn_{2i-1}, \mylabel{bh.al1} \\
Q_{2i-1}' = \qn_{2i-1}(\alpha_{2k}^*) = N(N(\qn_1,\qn_3,&\dots,\qn_{2i-1}\tin P_{2k}^*) \tin \qn_{2i-1})=\notag \\
 C(N(\qn_1,\qn_3,\dots,\qn_{2i-1}\tin P_{2k}^*) \tin \qn_{2i-1})
&=N(P_{2k}^* \tin \qn_{2i-1}),  \mylabel{bh.al2}\\
\bn_{2i-1}\in \Irr(\qn_{2i-1})  \text{ \, is the }
  N(\qn_1,\qn_3,\dots,\qn_{2i-1}\tin
P_{2k}^*)&\text{-Glauberman correspondent of } 
\beta_{2i-1}', \mylabel{bh.al3}\\
\pn_{2i} = N(\qn_1,\qn_3 , \dots,\qn_{2i-1} \tin P_{2i}')
& =N(\qn_1,\qn_3 , \dots,\qn_{2i-1} \tin P_{2i}^*), \mylabel{bh.al4}\\
\an_2 = \alpha_2^* \in \Irr(\pn_2),  &\text{ while for $i >1$ } \notag \\ 
\an_{2i}\in \Irr(\pn_{2i}) \text{ is the } 
\qn_3,\dots,\qn_{2i-1} &\text{-correspondent of } 
\alpha_{2i}^*, \mylabel{bh.al5}\\ 
\qn_1:&=G_1, \text { \, while  for $i>2$ } \notag \\
\qn_{2i-1} \in \Hall_{p'}(G_{2i-1}&(\an_2,\dots,\an_{2i-2},\bn_1,
\dots,\bn_{2i-3})), \mylabel{bh.al6}\\
\pn_0:=1 &\text{ and } \an_0:=1,  \text{ while for $i \geq 1$ } \notag \\ 
\pn_{2i} \in \Syl_{p}(G_{2i}&(\an_0, \an_2,\dots,\an_{2i-2},\bn_1,
\dots,\bn_{2i-1})), \mylabel{bh.al7}\\
\mylabel{bh.al8}P_{2j}^* =\pn_0 \cdot \pn_2\cdot \pn_4 \cdots \pn_{2i-2} \cdot
 &N(\qn_1,\qn_3 , \dots,\qn_{2i-1} \tin P_{2j}^*),  \\ 
\mylabel{bh.al9} P_{2i}^*=\pn_2 \cdot \pn_4 \cdots \pn_{2i}, 
\end{align} 
\end{subequations}
whenever $1\leq i \leq j \leq k$.
In \eqref{bh.al5}, the  $\qn_1,\dots,\qn_{2i-1}$-correspondent 
refers to the correspondence described in Lemma \ref{p*l1} and 
Theorem \ref{p*t1}. 
\end{theorem}

To prove the above theorem we will need the following easy lemma 
\begin{lemma}\mylabel{bh.l4}
Assume that $T,S,T_0$ are $p$-subgroups of a finite group $H$ such that 
$T$ normalizes $S$ while $T_0 \leq T\cap S$. Assume further that
 $T_0 =N(T \tin S)$. Then $T_0 = S$.
\end{lemma}

\begin{proof}
The group 
$TS$ is a $p$-group and thus nilpotent.
Hence, if $T$ is a proper subgroup of $TS$, then we should have that 
 $N(T \tin TS) > T$.

But $N(T \tin  TS) = T \cdot N(T \tin S) = T \cdot T_0$.
As $T_0 \leq T$ we conclude that $N(T \tin TS) =T$.
Hence $T=TS$ and thus $S \leq T$. 
Therefore  we have that $S = N(T \tin S) =T_0$. 
\end{proof}

\begin{proof}[Proof of Theorem \ref{bh.th2}]
We define $\pn_0:=1$ and $\an_0:=1$, so that the trivial part of \eqref{bh.al7} holds.  
We  will  prove theorem using induction on $i$.
Assume that $i=1$.
For  $\qn_1$ we take $\qn_1 =G_1$. Hence  the first part of \eqref{bh.al6}  holds.
 According to \eqref{bh.e6a}  we have that 
$$
Q_1' = \qn_1(\alpha_{2k}^*) =N(P_{2k}^* \tin \qn_1)=C(P_{2k}^* \tin \qn_1).
$$
As $\qn_1 =G_1 \unlhd  G$ we  get that $N(\qn_1 \tin P_{2k}^*) = P_{2k}^*$
and  that $Q'$ normalizes $\qn_{1}$.
Hence $\qn_1$ satisfies (\ref{bh.al}a,b). 
Furthermore,  $\beta_1'$ is fixed by 
$P_{2k}^*$, by \eqref{bh.e6d}. So we take  $\bn_1 \in \Irr^{P_{2k}^*}(Q_1)$ to be the
 $P_{2k}^*$-Glauberman
correspondent of $\beta_1' \in \Irr(Q_1')$. 
Thus $\bn_1$ satisfies \eqref{bh.al3}.

Let $\pn_2 :=N(\qn_1 \tin P_{2}^*)$. 
Then,  as $\qn_1=G_1 \unlhd  G$ and $P_{2}^* = P_2'$ (see \eqref{bh.e6c}),
 we have 
$$
\pn_2 = P_{2}^*=P_{2}'.
$$
So \eqref{bh.al9} holds.
 We also define $\an_2 := \alpha_2^*$. Thus $\pn_2$ and $\an_2$ satisfy 
(\ref{bh.al}d,e). 

Let $\map$ be a $p$-Sylow subgroup of $G_2(\bn_1)$. 
Then $G_2(\bn_1) = \map \ltimes G_1= \map \ltimes \qn_1$, as $G_2 /G_1$
is a $p$-group.  
Furthermore, $P_{2k}^*$ normalizes $G_2(\bn_1)$, as it fixes $\bn_1$.
Also, $P_2^* = G_2 \cap P_{2k}^*$ is a subgroup of $G_2(\bn_1)$. 
Therefore Lemma \ref{pq50} implies that we can pick $\map$ so that is
normalized by $P_{2k}^*$,   while $P_2'= P_{2}^* \leq \map$.
So $P_{2}' \leq \map \cap P_{2k}^*$. 
The group $N(P_{2k}^* \tin \map )$  fixes
the $P_{2k}^*$-Glauberman correspondent  $\bn_1$ of $\beta_1'$ 
 (as $\map$ does), and normalizes $P_{2k}^*$. Thus it also fixes 
$\beta_1'$. Hence
 $P_2'\leq N(P_{2k}^* \tin \map) \leq N(P_{2k}^* \tin G_2(\beta_1'))$.
According to \eqref{bh.ee},  the group $P_{2}'$ is a $p$-Sylow subgroup of 
$N(P_{2k}^* \tin G_2(\beta_1'))$. 
Thus Lemma \ref{bh.l4} can be applied to the groups $P_{2k}^* , \map$ and 
$ P_{2}^*=P_2'$
in the place of $T, S$ and $T_0$,  respectively. 
Therefore we get that $\map = P_{2}'=P_{2}^*$. 
As $\pn_2 :=P_2^*$, we conclude that $\pn_2 \in \Syl_{p}(G_2(\bn_1))$. 
Hence \eqref{bh.al7} holds.

We complete the proof of the $i=1$ case by observing  that $N(\qn_1 \tin P_{2j}^*)=
P_{2j}^*$ as $\qn_1 =G_1 \unlhd G$. Thus 
$$\pn_0 \cdot N(\qn_1 \tin P_{2j}^* ) =1  \cdot P_{2j}^* = P_{2j}^*,$$ 
whenever $1\leq j \leq k$. 
Hence \eqref{bh.al8} holds.

Now assume Theorem \ref{bh.th2} holds for all $i=1,\dots,t-1$,  for some 
$t =2,\dots,k$.
We will prove it also holds for $i=t$.
To simplify this proof, we give  separately the next steps 
that depend heavily on the  inductive hypothesis, 
\begin{step}\mylabel{bh.st1}
Assume that the set $\{\qn_{2i-1},\pn_{2i} | \bn_{2i-1},\an_{2i}\}_{i=1}^s$, 
for some  $ s=1, \dots,    k$, 
satisfies { \rm (\ref{bh.al}c,d,e) }  for all $i=1,\dots, s$.
 Let  $T \leq N(P_{2k}^*, \qn_1,\qn_3,\dots,\qn_{2s-1} \tin G(\alpha_{2s}^* , 
\alpha_2',\dots,\alpha_{2s}', \beta_1',\dots,\beta_{2s-1}'))$. Then 
$$T\leq  G(\an_2,\dots,\an_{2s},\bn_1,\dots,\bn_{2s-1}).$$
\end{step}
\begin{proof}
The group $T$ normalizes $P_{2r}^*$ for all $r=1,\dots,s$,  as it normalizes $P_{2k}^*$.
It also normalizes $\qn_{2r-1}$ for all such $r$. Therefore, 
it normalizes $N(\qn_1,\dots,\qn_{2r-1} \tin P_{2r}^*)$. 
But,  according to \eqref{bh.al4}, this last group equals $\pn_{2r}$ 
whenever $1\leq r \leq s$. Hence $T$ normalizes $\pn_{2r}$ and $\qn_{2r-1}$. But  
$T$ also   fixes $\alpha_{2r}^*$,  as it fixes $\alpha_{2k}^*$ 
(see Remark \ref{b.rem1}). Therefore \eqref{bh.al5}, along with Proposition 
\ref{pqp*t}, implies that  $T$ fixes $\an_{2r}$ for all $r=1,\dots,s$. 

Furthermore,  $T$ fixes $\bn_{2r-1}$,  as 
it   fixes its $N(\qn_1,\dots,
\qn_{2r-1}  \tin P_{2k}^*)$-Glauberman correspondent  $\beta_{2r-1}'$
 (see \eqref{bh.al3}) and normalizes $P_{2k}^*, \qn_1,\dots,\qn_{2r-1}$,
 whenever $1\leq r \leq s$.
Therefore
$$
T\leq    G(\an_2,\dots,\an_{2s},\bn_1,\dots,\bn_{2s-1}),
$$ and Step
\ref{bh.st1} is complete.
\end{proof}
The second step is
\begin{step}\mylabel{bh.st2}
Assume that the set $\{\qn_{2i-1},\pn_{2i} | \bn_{2i-1},\an_{2i}\}_{i=1}^s$, for some $s=1,\dots,k$,  
satisfies { \rm (\ref{bh.al}b,c,e) } for all $i=1,\dots,s$. 
Let $T \leq N(P_{2k}^*  \tin G(\an_2,\dots,\an_{2s},\bn_1,\dots,\bn_{2s-1}))$. Then 
 $$T \leq N(P_{2k}^* \tin G(\alpha_{2s}^*,\alpha_2',\dots,\alpha_{2s}',\beta_1',
\dots,\beta_{2s-1}')).
$$ 
\end{step}
\begin{proof}
The group $T$ normalizes $P_{2r}^*$ for all $r=1,\dots,s$ (even for all $r=s+1, \dots,k$,
 but this we will not need). As $T$ also  normalizes the groups 
$\qn_{2r-1}$,  and fixes the characters $\an_{2r}$, it has to fix (by
Proposition \ref{pqp*t}) the $\qn_3,\dots,\qn_{2r-1}$-correspondent 
$\alpha_{2r}^*$ of $\an_{2r}$ (see \eqref{bh.al5}) for all $r=2, \dots,s$.
It also fixes $\alpha_2^* = \an_2$.

Furthermore, $T$ normalizes $Q_{2r-1}'$, as  
\eqref{bh.al2} implies that 
 $Q_{2r-1}' =N(P_{2k}^* \tin \qn_{2r-1})$,  whenever $1\leq r
\leq s$. The group $P_{2r}'$ satisfies \eqref{bh.e} for $i=r$. 
Hence $P_{2r}'= N(Q_{2r-1}' \tin P_{2r}^*)$. Therefore, $T$ normalizes 
$P_{2r}'$,  as it normalizes both $Q_{2r-1}'$ and $P_{2r}^*$.
This, along with the fact that $T$ fixes $\alpha_{2r}^*$, 
implies that  $T$ fixes the $Q_{2r-1}'$-Glauberman correspondent 
$\alpha_{2r}'$ (see \eqref{bh.eee}) of $\alpha_{2r}^*$, 
for all $r=1,\dots,s$.

Even more, as $T$ fixes $\bn_{2r-1}$ and normalizes $Q_{2r-1}'$, 
 it must  fix the 
$N(\qn_1,\dots,\qn_{2r-1} \tin P_{2k}^*)$-Glauberman correspondent 
$\beta_{2r-1}'$ of $\bn_{2r-1}$ (see  \eqref{bh.al3}).
Hence 
$$T \leq  N(P_{2k}^* \tin G(\alpha_{2s}^*,\alpha_2',\dots,\alpha_{2s}', 
\beta_1',\dots,\beta_{2s-1}')).
$$
\end{proof}

The last step is 
\begin{step}\mylabel{bh.st3}
The group   $S:= N(\qn_1,\dots,\qn_{2t-3} \tin P_{2k}^*)$ is a subgroup of 
 $G(\an_2,\dots,\an_{2t-2},\bn_1,\dots,\bn_{2t-3})$.
\end{step}

\begin{proof}
For every $r=1,\dots,t-1$,  the group $S= N(\qn_1,\dots,\qn_{2t-3} \tin P_{2k}^*)$
normalizes $\pn_{2r}=N(\qn_1,\dots,\qn_{2r-1} \tin P_{2r}^*)$. Also $S$ fixes 
$\alpha_{2r}^*$,  as it is a subgroup of $P_{2k}^*$. 
Therefore   $S$ fixes 
the $\qn_3,\dots,\qn_{2r-1}$-correspondent
$\an_{2r} \in \Irr(\pn_{2r})$ of $\alpha_{2r}^*$, (see \eqref{bh.al5}), 
for all $r=2,\dots,t-1$,  
as well as $\an_2=\alpha_2^*$. 
 
Furthermore, \eqref{bh.al3} for $i=t-1$ implies that 
$S$ fixes $\bn_{2t-3}$. Similarly,
the inductive hypothesis for \eqref{bh.al3} implies that 
$N(\qn_1,\dots,\qn_{2r-1} \tin P_{2k}^*)$ fixes $\bn_{2r-1}$,  for all
 $r=1,\dots,t-2$. 
But $S$ is a subgroup of  
$N(\qn_1,\dots,\qn_{2r-1} \tin P_{2k}^*)$ for all such $r$.   
Hence $S$  fixes $\bn_{2r-1}$ whenever $1\leq r \leq t-1$.

Therefore $S \leq G(\an_2,\dots,\an_{2t-2},\bn_1,\dots,\bn_{2t-3})$, 
and Step \ref{bh.st3} is proved.
\end{proof}

We can now continue with the proof of the theorem.
The fact that \eqref{bh.e5} is a triangular set for \eqref{bh.*} implies that 
\begin{multline*}
Q_{2s-1}' \leq G_{2s-1}(\alpha_{2k}^*,\alpha_2',\dots,\alpha_{2t-2}',\dots
,\alpha_{2s-2}',\beta_1',\dots,\beta_{2t-3}',\dots,\beta_{2s-3}')
\leq \\
G_{2s-1}(\alpha_{2k}^*,\alpha_2',\dots,\alpha_{2t-2}',
\beta_1',\dots,\beta_{2t-3}'),
\end{multline*}
 for all
$s=t,\dots,k$.
  Also 
the inductive hypothesis,  \ref{bh.al1},  for $i\leq t-1$ implies  that 
$Q_{2s-1}'\leq Q'$ normalizes the groups $\qn_1,\dots,\qn_{2t-3}$.
Hence for all $s=t,\dots,k$ we have that 
\begin{multline*}
Q_{2s-1}'\leq N(\qn_1,\dots,\qn_{2t-3} \tin
 G_{2s-1}(\alpha_{2k}^*,\alpha_2',\dots,\alpha_{2t-2}',
\beta_1',\dots,\beta_{2t-3}')\leq \\
N(P_{2k}^*, \qn_1,\dots,\qn_{2t-3} \tin
 G_{2s-1}(\alpha_{2t-2}^*,\alpha_2',\dots,\alpha_{2t-2}',
\beta_1',\dots,\beta_{2t-3}')),
\end{multline*}
and 
\begin{multline*}
Q'\leq N(\qn_1,\dots,\qn_{2t-3} \tin
 G(\alpha_{2k}^*,\alpha_2',\dots,\alpha_{2t-2}',
\beta_1',\dots,\beta_{2t-3}')\leq \\
N(P_{2k}^*, \qn_1,\dots,\qn_{2t-3} \tin
 G(\alpha_{2t-2}^*,\alpha_2',\dots,\alpha_{2t-2}',
\beta_1',\dots,\beta_{2t-3}')).
\end{multline*}
This, along with Step \ref{bh.st1}, with the present $t-1$ in the place of $s$ there, 
 and the fact that
 $G_{2t-1} \unlhd G_{2s-1}$, implies
that both  $Q_{2s-1}'$ and $Q'$ normalize the group 
$ G_{2t-1}(\an_2,\dots,\an_{2t-2},\bn_1,\dots,\bn_{2t-3})$ and fix the characters 
$\an_2, \dots,\an_{2t-2}, \bn_1, \dots, \bn_{2t-3}$.
In particular, for $s=t$   we get  
\begin{subequations}\mylabel{bh.e8}
\begin{equation}
Q_{2t-1}' \leq  G_{2t-1}(\an_2,\dots,\an_{2t-2},\bn_1,\dots,\bn_{2t-3}),
\end{equation}
as $Q_{2t-1}' \leq G_{2t-1}$. 
Furthermore, 
\begin{equation} 
 Q' \text{ normalizes }
  G_{2t-1}(\an_2,\dots,\an_{2t-2},\bn_1,\dots,\bn_{2t-3}).
\end{equation}
\end{subequations}

Let $\maq$ be a $p'$-Hall subgroup of 
$ G_{2t-1}(\an_2,\dots,\an_{2t-2},\bn_1,\dots,\bn_{2t-3})$.
Since $\pn_{2t-2}$ satisfies \eqref{bh.al7}  for $i=t-1$, 
and since $\an_{2t-2} \in \Irr(\pn_{2t-2})$,  we have that 
$$
\pn_{2t-2} \in \Syl_{p}(G_{2t-2}(\an_2,\dots,\an_{2t-2},
\bn_1,\dots,\bn_{2t-3})).
$$
As $G_{2t-1}/G_{2t-2}$ is a $p'$-group and 
$G_{2t-2}(\an_2,\dots,\an_{2t-2}, \bn_1,\dots,\bn_{2t-3})$ normalizes 
$P_{2t-2}^{\nu}$,
 we get 
\begin{equation}\mylabel{bh.e9}
  G_{2t-1}(\an_2,\dots,\an_{2t-2},\bn_1,\dots,\bn_{2t-3})=
\maq \ltimes \pn_{2t-2}. 
\end{equation}
This, along with \eqref{bh.e8}, the 
fact that $Q_{2t-1}' \unlhd Q'$ (see \eqref{bh.e5b}),  
 and Lemma \ref{pq50}, implies that we can pick
a conjugate $\qn_{2t-1}:=\maq^s$ of  $\maq$,
so that 
\begin{equation}\mylabel{bh.e9.5}
\begin{aligned}
\qn_{2t-1} \in
 \Hall_{p'}&(G_{2t-1}(\an_2,\dots,\an_{2t-2},\bn_1,\dots,\bn_{2t-3})),\\
Q' &\text{ normalizes  } \qn_{2t-1} \text{ and }  \\
&Q_{2t-1}' \leq \qn_{2t-1}.
\end{aligned}
\end{equation}

It is obvious from the definition of $\qn_{2t-1}$ that it satisfies 
(\ref{bh.al}a,f) for $i=t$. 
Furthermore,  \eqref{bh.e9} holds for $\maq = \qn_{2t-1}$, i.e., 
\begin{equation}\mylabel{bh.e10}
G_{2t-1}(\an_2,\dots,\an_{2t-2},\bn_1,\dots,\bn_{2t-3})=
\qn_{2t-1} \ltimes \pn_{2t-2}. 
 \end{equation}
This, along with Step \ref{bh.st3}, implies that
 $N(\qn_1,\dots,\qn_{2t-3} \tin P_{2k}^*)$ normalizes 
$\qn_{2t-1} \ltimes \pn_{2t-2}$.  Hence the product
$N(\qn_1,\dots,\qn_{2t-3} \tin P_{2k}^*) \cdot \qn_{2t-1} \pn_{2t-2}$
is a group having  $\qn_{2t-1} \ltimes \pn_{2t-2}$ as a normal subgroup.
 Furthermore, \eqref{bh.al4}  for $i=t-1$ implies
that $\pn_{2t-2} \leq   
N(\qn_1,\dots,\qn_{2t-3} \tin P_{2k}^*)$.
Hence $N(\qn_1,\dots,\qn_{2t-3} \tin P_{2k}^*)$ is a $p$-Sylow 
subgroup of $N(\qn_1,\dots,\qn_{2t-3} \tin P_{2k}^*) \cdot \qn_{2t-1} \pn_{2t-2}$.
Thus Frattini's argument for the $p'$-Hall subgroup $\qn_{2t-1}$ 
of the normal subgroup $\qn_{2t-1} \ltimes \pn_{2t-2}$ implies that 
\begin{multline*}
 N(\qn_1,\dots,\qn_{2t-3} \tin P_{2k}^*) = 
\pn_{2t-2} \cdot N(\qn_{2t-1}\tin N(\qn_1,\dots,\qn_{2t-3} \tin P_{2k}^*))=\\
\pn_{2t-2} \cdot N(\qn_1,\dots,\qn_{2t-3},\qn_{2t-1} \tin P_{2k}^*).
\end{multline*}
This, along with \eqref{bh.al8} for $i=t-1$ and $j=k$, implies that
\begin{align*}
P_{2k}^* &= \pn_2 \cdots \pn_{2t-4} 
\cdot N(\qn_1,\dots,\qn_{2t-3} \tin P_{2k}^*)=\\
&\pn_2 \cdots \pn_{2t-4}\cdot \pn_{2t-2} \cdot 
N(\qn_1,\dots,\qn_{2t-3},\qn_{2t-1} \tin P_{2k}^*).
\end{align*}
Therefore,
 intersecting both sides of the above equation with $G_{2j}$, we get  
\begin{equation}\mylabel{bh.e11}
P_{2j}^* = \pn_2 \cdots \pn_{2t-2} \cdot 
N(\qn_1,\dots,\qn_{2t-3},\qn_{2t-1} \tin P_{2j}^*),
\end{equation}
whenever $t \leq j \leq k $. Hence \eqref{bh.al8} holds  for $i=t$ and $j=i , i+1, \dots,k$.

To prove \eqref{bh.al2} for $i=t$, we first note that,  according to 
 the definition of $\qn_{2t-1}$ (see \eqref{bh.e9.5}), 
 we have  $Q_{2t-1}' \leq \qn_{2t-1}$.
Hence  
$Q_{2t-1}' \leq \qn_{2t-1}(\alpha_{2k}^*)$, as
 $Q_{2t-1}'\leq G'=G(\alpha_{2k}^*)$.
Furthermore, $\qn_{2t-1}$ normalizes $\qn_1,\dots,\qn_{2t-3},\qn_{2t-1}$. 
Hence 
$$N(P_{2k}^* \tin \qn_{2t-1}) \leq 
 N(N(\qn_1,\dots,\qn_{2t-1} \tin P_{2k}^*) \tin \qn_{2t-1})  =
C(N(\qn_1,\dots,\qn_{2t-1} \tin P_{2k}^*) \tin \qn_{2t-1}),
$$
where the last equality holds as the $p$-group 
$N(\qn_1,\dots,\qn_{2t-1} \tin P_{2k}^*)$ normalizes the $p'$-group 
$\qn_{2t-1}$.
 Thus we have
\begin{equation}\mylabel{bh.e12}
\begin{aligned} 
Q_{2t-1}' \leq \qn_{2t-1}(\alpha_{2k}^*) \leq N(P_{2k}^* \tin 
 \qn_{2t-1}) \leq  \\
N(N(\qn_1,\dots,\qn_{2t-1} \tin P_{2k}^*) \tin \qn_{2t-1}) = \\
C(N(\qn_1,\dots,\qn_{2t-1} \tin P_{2k}^*) \tin \qn_{2t-1}).
\end{aligned}
\end{equation}

Let $T= C(N(\qn_1,\dots,\qn_{2t-1} \tin P_{2k}^*) \tin \qn_{2t-1})$.
Then $T$ normalizes the groups $\pn_2,\dots,\pn_{2t-2}$,
 as 
$\qn_{2t-1}$ does (it fixes their  characters  $\an_{2i}$).
Hence, $T$,   in view of \eqref{bh.e11} for $j=k$,  also normalizes $P_{2k}^*$.
Therefore, $T \leq N(P_{2k}^* \tin \qn_{2t-1})$.
This , in view of \eqref{bh.e9.5}, implies  that 
$$
T \leq N(P_{2k}^* \tin G_{2t-1}(\an_2,\dots,\an_{2t-2},\bn_1,\dots,\bn_{2t-3})).
$$
The set $\{\qn_{2i-1},\pn_{2i}|\bn_{2i-1},\an_{2i}\}_{i=1}^{t-1}$ satisfies 
(\ref{bh.al}b,c,e) (according to the inductive hypothesis).
 So  Step \ref{bh.st2} for $s=t-1$,  implies that  $T$ 
 satisfies
\begin{equation}\mylabel{bh.e14}
T  \leq 
N(P_{2k}^* \tin G_{2t-1}(\alpha_{2t-2}^*,\alpha_2',\dots,
\alpha_{2t-2}',\beta_1',\dots,\beta_{2t-3}')).
\end{equation}

Equation \eqref{bh.e11}
for $j=k$, along with \eqref{bh.al9} for $i=t-1$, implies that  
$$
\frac{P_{2k}^*}{P_{2t-2}^*} \cong 
\frac{N(\qn_1,\dots,\qn_{2t-1} \tin P_{2k}^*)}{P_{2t-2}^*\cap 
N(\qn_1,\dots,\qn_{2t-1} \tin P_{2k}^*)}.
$$
Therefore $T$ centralizes $P_{2k}^* / P_{2t-2}^*$, as it centralizes 
$N(\qn_1,\dots,\qn_{2t-3} \tin P_{2k}^*)$. 
Also $T$  fixes $\alpha_{2t-2}^*\in \Irr(P_{2t-2}^*)$,  and is a $p'$-group.
Hence (see Exercise 13.13 in \cite{is}), $T$ fixes every irreducible character 
of $P_{2k}^*$ that lies above $\alpha_{2t-2}^*$. Thus $T$ fixes
$\alpha_{2k}^*$. This, along with \eqref{bh.e14}, implies that
$$
T =C(N(\qn_1,\dots,\qn_{2t-1} \tin P_{2k}^*) \tin \qn_{2t-1})
\leq G_{2t-1}(\alpha_{2k}^*,\alpha_2',\dots,\alpha_{2t-2}',\beta_1',\dots,
\beta_{2t-1}').
$$ 

But $Q_{2t-1}'$ is a $p'$-Hall subgroup of 
$G_{2t-1}(\alpha_{2k}^*,\alpha_2',\dots,\alpha_{2t-2}',\beta_1',\dots,
\beta_{2t-1}')$,  as \eqref{bh.e5} is a triangular set for \eqref{bh.*}.
Furthermore, \eqref{bh.e12} implies that the $p'$-group $Q_{2t-1}'$ 
is contained in 
$C(N(\qn_1,\dots,\qn_{2t-1} \tin P_{2k}^*) \tin \qn_{2t-1})$. 
Thus $C(N(\qn_1,\dots,\qn_{2t-1} \tin P_{2k}^*) \tin \qn_{2t-1}) = Q_{2t-1}'$.
This, along with \eqref{bh.e12}, implies that 
\begin{equation}\mylabel{bh.e13}
\begin{aligned}
Q_{2t-1}'= \qn_{2t-1}(\alpha_{2k}^*) &= N(P_{2k}^* \tin 
\qn_{2t-1}) = \\
&N(N(\qn_1,\dots,\qn_{2t-1} \tin P_{2k}^*) \tin \qn_{2t-1}) = \\
&C(N(\qn_1,\dots,\qn_{2t-1} \tin P_{2k}^*) \tin \qn_{2t-1}).\\
\end{aligned}
\end{equation}
So \eqref{bh.al2} holds for $i=t$.
Hence we have shown that the group $\qn_{2t-1}$ satisfies 
(\ref{bh.al}a,b,f,h) for $i=t \leq j \leq k$.

As $Q_{2t-1}' = C(N(\qn_1,\dots,\qn_{2t-1} \tin P_{2k}^*) \tin \qn_{2t-1})$
we can define $\bn_{2t-1}  \in \Irr(\qn_{2t-1})$ to be  the 
$N(\qn_1,\dots,\qn_{2t-1} \tin P_{2k}^*)$-Glauberman correspondent
of $\beta_{2t-1}' \in \Irr(Q_{2t-1}')$. 
Thus $\bn_{2t-1}$ satisfies \eqref{bh.al3} for $i=t$.

To complete the inductive step  it remains to prove that we can pick 
a $p$-group $\pn_{2t}$, along with its irreducible character $\an_{2t}$,
so that (\ref{bh.al}d,e,g,i) hold for $i=t$.
In view of  \eqref{bh.e13}  we have  that $Q_{2t-1}' = 
\qn_{2t-1}(\alpha_{2k}^*)$. Hence $ N(\qn_1,\dots,\qn_{2t-1} \tin P_{2t}^* )$ 
normalizes $Q_{2t-1}'$ .
This, along with the fact that $P_{2t}' = N(Q_{2t-1}' \tin P_{2t}^*)$
by \eqref{bh.ee}, implies that 
\begin{multline*}
N(\qn_1,\dots,\qn_{2t-1} \tin P_{2t}^* ) =
N(Q_{2t-1}' \tin N(\qn_1,\dots,\qn_{2t-1} \tin P_{2t}^* ) ) = \\
N(\qn_1,\dots,\qn_{2t-1} \tin N(Q_{2t-1}' \tin P_{2t}^* )) =
 N(\qn_1,\dots,\qn_{2t-1} \tin P_{2t}' ) .
\end{multline*}

Let 
\begin{equation}\mylabel{bh.e14.5}
\begin{aligned}
M_0:&= N(\qn_1,\dots,\qn_{2t-1} \tin P_{2t}' ) =N(\qn_1,\dots,\qn_{2t-1}
 \tin P_{2t}^* )  \text{ and }\\ 
M:&= N(\qn_1,\dots,\qn_{2t-1} \tin P_{2k}^* ). 
\end{aligned}
\end{equation}
Note that $M \cap G_{2t}= M_0$,  as $P_{2k}^* \cap G_{2t} = P_{2t}^*$.

In view of Step \ref{bh.st3}, the group $M$ 
fixes the characters $\an_2,\dots,\an_{2t-2},\bn_1,\dots,\bn_{2t-3}$
 as it is a subgroup of 
$N(\qn_1,\dots,\qn_{2t-3} \tin P_{2k}^*) $.Furthermore, the definition of $\bn_{2t-1}$ (as the $M$-correspondent 
of $\beta_{2t-1}'$) implies that $M$ also fixes $\bn_{2t-1}$.
Hence $M$ normalizes
 $G_{2t}(\an_2,\dots,\an_{2t-2},\bn_1,\dots,\bn_{2t-1})$, 
while $M_0 = M\cap  G_{2t}$ 
is a subgroup of $G_{2t}(\an_2,\dots,\an_{2t-2},\bn_1,\dots,\bn_{2t-1})$. 
Let $\map$ be a $p$-Sylow subgroup of $G_{2t}(\an_2,\dots,\an_{2t-2},\bn_1,\dots,\bn_{2t-1})$, 
chosen so that $\map$ contains $M_0$. 
It is clear from the fact that $G_{2t}/G_{2t-1}$ is a $p$-group,
 and the definition of $\qn_{2t-1}$  (see \eqref{bh.e9.5}), that 
\begin{equation}\mylabel{bh.e15}
G_{2t}(\an_2,\dots,\an_{2t-2},\bn_1,\dots,\bn_{2t-1})=
\map \ltimes \qn_{2t-1}.
\end{equation}
Therefore Lemma \ref{pq51} implies that there exists a $\qn_{2t-1}$-conjugate 
of  $\map$ 
that is normalized by $M$ and contains $M_0$. So we may replace $\map$ by this conjugate and
 assume that $M_0 \leq M\cap \map$.

We can show the following 
\begin{claim}\mylabel{bh.cl1}
$N(M \tin \map) = M_0$.
\end{claim}
\begin{proof}
It is obvious that $M_0 \leq N(M \tin \map)$. 
For the other inclusion we first note that $N(M \tin \map)$ normalizes $\pn_2,\dots,
\pn_{2t-2}$  (since $\map$  does) and $M$. Hence $N(M \tin \map)$ 
normalizes $P_{2k}^*= \pn_2\cdots \pn_{2t-2}\cdot M$ (see \eqref{bh.e11}).
Hence 
\begin{align}\mylabel{bh.e16}
N(&M \tin \map) \leq  N(P_{2k}^* \tin \map)   \notag \\
&\leq N(P_{2k}^* \tin G_{2t}(\an_2,\dots,\an_{2t-2},\bn_1,
\dots,\bn_{2t-3},\bn_{2t-1}))  \\
&\leq N(\qn_1,\dots,\qn_{2t-1},P_{2k}^* \tin G_{2t}(\an_2,\dots,\an_{2t-2},\bn_1,
\dots,\bn_{2t-3}))(\bn_{2t-1} ) \notag   \\
&\leq N(\qn_1,\dots,\qn_{2t-1},P_{2k}^* \tin 
G_{2t}(\alpha_{2t-2}^*,\alpha_2',\dots,\alpha_{2t-2}',\beta_1',\dots,\beta_{2t-3}'))
 (\bn_{2t-1}), \notag 
\end{align}
where the last inclusion holds according to  Step \ref{bh.st2}  for $s=t-1$.

Also $N(M \tin \map )$ fixes the $M$-Glauberman correspondent 
$\beta_{2t-1}'$  of $\bn_{2t-1}$, as it fixes $\bn_{2t-1}$. This, along with 
\eqref{bh.e16}, implies that 
\begin{equation}\mylabel{bh.e17}
 N(M \tin \map) \leq  N(\qn_1,\dots,\qn_{2t-1}, P_{2k}^*
 \tin G_{2t}(\alpha_{2t-2}^*,
\beta_{2t-1}')).
\end{equation}
But $P_{2t}'$ satisfies \eqref{bh.ee}  for $i=t$. 
Therefore $M_0=N(\qn_1, \dots,\qn_{2t-1} \tin P_{2t}')$ is a $p$-Sylow subgroup of 
$N(\qn_1,\dots,\qn_{2t-1}, P_{2k}^* \tin G_{2t}(\alpha_{2t-2}^*,
\beta_{2t-1}'))$. As $M_0$ is contained in $N(M \tin \map)$, 
 inclusion  \eqref{bh.e17} implies that 
$M_0 = N(M \tin \map)$. Hence the claim follows.
\end{proof}

The groups $M , M_0$ and $\map$ satisfy the hypothesis of Lemma \ref{bh.l4},
in the place of $T,T_0$ and $S$ respectively. So we conclude that 
$M_0 =\map$. Therefore, $M_0$ is a 
$p$-Sylow subgroup of the group 
 $G_{2t}(\an_2,\dots,\an_{2t-2},\bn_1,\dots,\bn_{2t-1})$.
If we define  $\pn_{2t}:=M_0$, then it is clear that $\pn_{2t}$ satisfies 
\eqref{bh.al7} for $i=t$.  It also satisfies \eqref{bh.al4} for $i=t$, 
 as    \eqref{bh.e14.5}  shows.
Equation \eqref{bh.al9} for $i=t$ follows clearly from \eqref{bh.al8} for $i=j=t$, 
(that we have already proved in \eqref{bh.e11})
and \eqref{bh.al4} for $i=t$.

To complete the inductive step, it remains to show  that we can pick a
 character  $\an_{2t} \in \Irr(\pn_{2t})$ that satisfies \eqref{bh.al5} for $i=t$. That is, 
it suffices to show that the character $\alpha_{2t}^*$ has a 
$\qn_3,\dots,\qn_{2t-1}$-correspondent. 
Looking at  Lemma \ref{p*l1} and Theorem \ref{p*t1}, where this correspondence 
is described, we observe that it is enought to prove the following for every $i,j$ with 
$2\leq j < i \leq t$:
\begin{itemize}
\item[[1]]
$\qn_{2j-1} \cdot \pn_{2j-2} \cdot \pn_{2j} \cdots \pn_{2i}$ is a group 
containing $\qn_{2j-1}\cdot \pn_{2j-2}$ and $\pn_{2j-2}$ as normal subgroups.
\item[[2]]
$N(\qn_{2j-1} \tin \pn_{2j-2}\cdot \pn_{2j} \cdots \pn_{2i}) = 
\pn_{2j}\cdots \pn_{2i}.$
\item[[3]]
$\alpha_{2t}^*$ satisfies Property \ref{p*pp1}, i.e., there exist characters 
$\alpha_{2s,1}^* \in \Irr(P_{2s}^*)$, 
for  $s=1,\dots,t$,  such that $\alpha_{2t,1}^*=\alpha_{2t}^*$,  and, if $s < t$, then 
 $\alpha_{2s,1}^*$  is $\qn_{2s+1}$-invariant 
and lies under $\alpha_{2s+2,1}^*$.
\end{itemize}

Part [1] is  clear as, according to (\ref{bh.al}e,f), for every $s=j+1, \dots,i$ 
the group   $\pn_{2s}$ normalizes  $\pn_{2j-2}, \dots,\pn_{2s-2}$ and $\qn_{2j-1}$, 
 while $\qn_{2j-1}$ normalizes  $\pn_{2j-2}$.
This remark also implies that the product  $\pn_{2j}\cdots  \pn_{2i}$ normalizes 
$\qn_{2j-1}$. 
Hence 
 \begin{equation}\mylabel{bh.e18}
N(\qn_{2j-1} \tin \pn_{2j-2}\cdot \pn_{2j} \cdots \pn_{2i}) \geq 
\pn_{2j} \cdots \pn_{2i}.
\end{equation}
 But 
\begin{align*}
N(\qn_{2j-1} \tin \pn_{2j-2})&= N(\qn_{2j-1} \tin
 N(\qn_1,\dots,\qn_{2j-3} \tin P_{2j-2}^*))  &\text{ by \eqref{bh.al4} for $j-1$
 as $i$ there }\\
&=N(\qn_1,\dots,\qn_{2j-3},\qn_{2j-1} \tin P_{2j-2}^*)  &\\
&\leq N(\qn_1,\dots,\qn_{2j-3},\qn_{2j-1} \tin P_{2j}^*)   & \\
&=\pn_{2j}  &\text{ by \eqref{bh.al4} for $j$  as $i$ there }
\end{align*}
This, along with \eqref{bh.e18}, implies that 
$N(\qn_{2j-1} \tin \pn_{2j-2}\cdot \pn_{2j} \cdots \pn_{2i}) =
 \pn_{2j}\cdots \pn_{2i}$.
So [2] follows.

Part [3] holds if we take the characters $\alpha_{2s}^*$ in the place of
 $\alpha_{2s,1}^*$ for $s=1, \dots, t$.
 We only need to verify that $\qn_{2s+1}$ leaves $\alpha_{2s}^*$ 
invariant for every $s=1,\dots,t-1$.
This is clear as $\qn_{2s+1}$ fixes $\an_{2s}$  and normalizes the groups 
$\qn_{3}, \dots,\qn_{2s-1}$ (see \eqref{bh.al6} with $s-1$ as the $i$ there).
So it has to  fix the $\qn_3, \dots,\qn_{2s-1}$-correspondent 
$\alpha_{2s}^*$  (see \eqref{bh.al5} with $s$  as the $i$ there)  of $\an_{2s}$.
Hence [3] follows.
This proves that Lemma \ref{p*l1} and Theorem \ref{p*t1} can be applied.
Therefore there exists a unique character $\an_{2t} \in \Irr(\pn_{2t})$ 
that is the $\qn_3, \dots,\qn_{2t-1}$-correspondent of $\alpha_{2t}^*$.
 Thus \eqref{bh.al5} holds for $i=t$.

This completes the proof of the inductive step for $i=t$. Hence Theorem \ref{bh.th2}
holds.
\end{proof}
A useful consequence of Theorem \ref{bh.th2} is
\begin{corollary}\mylabel{bh.col1}
For every $i=1,\dots,k$ we have 
\begin{equation}\mylabel{bh.e19}
\pn_{2i} \cdot \pn_{2i+2} \cdots \pn_{2k} = N (\qn_1,\qn_3,\dots, \qn_{2i-1} 
\tin P_{2k}^*).
\end{equation}
Therefore  $\bn_{2i-1}\in \Irr(\qn_{2i-1} )$ 
is the $\pn_{2i}\cdots \pn_{2k}$-Glauberman
 correspondent of $\beta_{2i-1}' \in \Irr(Q_{2i-1}')= \Irr(C(\pn_{2i}\cdots \pn_{2k} 
\tin \qn_{2i-1}))$.
\end{corollary}

\begin{proof}
For every $j=i, \dots,k$, the group $\pn_{2j}$ normalizes
 $\qn_1,\qn_3,\dots,\qn_{2i-1} $ (see \eqref{bh.al4}).
This,  along with \eqref{bh.al9} and \eqref{bh.al4}, implies that 
\begin{multline*}
\pn_{2i} \cdots \pn_{2k} \leq N(\qn_1,\dots,\qn_{2i-1} \tin P_{2k}^*)=
 N(\qn_1,\dots,\qn_{2i-1} \tin P_{2i}^* \cdot \pn_{2i+2} \cdots \pn_{2k})\\
N(\qn_1,\dots, \qn_{2i-1} \tin P_{2i}^*)\cdot \pn_{2i+2}   \cdots \pn_{2k} =
\pn_{2i} \cdot \pn_{2i+2} \cdots \pn_{2k}.  
\end{multline*}
Thus \eqref{bh.e19}  holds.
The rest of the corollary is an obvious consequence of \eqref{bh.e19}, \eqref{bh.al2}
 and \eqref{bh.al3}.
\end{proof}

We have done all the neccesary work towards the proof of the main theorem 
of this section  which is a ``mirror'' 
of Theorem \ref{bh.th1}. That  is  
\begin{theorem}\mylabel{bh.th3}
The set \eqref{bh.e7},  constructed  in Theorem \ref{bh.th2},   forms a 
triangular set for \eqref{bh.*1}. Furthermore, 
the pair { \rm  (\ref{bh.e5}, \ref{bh.e7}) } satisfies 
\eqref{bh.e2} and Definition \ref{bh.d1}, i.e., 
\begin{subequations}\mylabel{bh.e20}
\begin{align}
Q_{2i-1}' = \qn_{2i-1,2k} &\text{ and } \beta_{2i-1}'= \bn_{2i-1,2k}, \mylabel{bh.e20a}\\ 
\pns_{2i} = P_{2i}^* =(P_{2i}')^*  &\text{ and }
 P_{2i}' = \pns_{2i}(\bn_{2i-1,2k}),  \mylabel{bh.e20b}\\
\ans_{2i} &= \alpha_{2i}^*, \mylabel{bh.e20c} \\
\alpha_{2i}' \text{ is the } Q_{2i-1}'&\text{-Glauberman correspondent 
of } \ans_{2i} \in \Irr(\pns_{2i}), \mylabel{bh.e20d}   
\end{align}
\end{subequations}
where $\pns_{2i} :=(P_{2i}^{\nu})^* = P_2^{\nu} \cdots P_{2i}^{\nu}$ and 
$\ans_{2i}:=(\alpha_{2i}^{\nu})^*$, whenever $1\leq i \leq k$.
\end{theorem}

\begin{proof}
The Properties (\ref{bh.al}j, g, f) of  that the set \eqref{bh.e7}  
 imply immediately parts (\ref{xxx}a,b,c,e) of the definition 
of a triangular set.

 Let   $\an_{2i-2,2i-1}$ denote    the 
irreducible character of $\pn_{2i-2,2i-1} :=C(\qn_{2i-1} \tin \pn_{2i-2} )$ 
that is  the $\qn_{2i-1}$-Glauberman 
 correspondent  of $\an_{2i-2} \in \Irr^{\qn_{2i-1}} (\pn_{2i-2})$, whenever $i=2,\dots,k$.  
As $\qn_{2i-1}$ normalizes  
$\pn_{2i-1}$, the $\qn_{2i-1}$-Glauberman correspondence 
coincides  with the $\qn_{2i-1}$-correspondence between
 $\Irr^{\qn_{2i-1}}(\pn_{2i-2})$ and 
$\Irr(\pn_{2i-2,2i-1}) =\Irr( C(\qn_{2i-1} \tin \pn_{2i-2}))$, by  Theorem \ref{dade:t2}.
This, along with \eqref{bh.al5}, implies that $\an_{2i-2,2i-1}$ is the 
$\qn_3, \dots,\qn_{2i-3}, \qn_{2i-1}$-correspondent of $\alpha_{2i-2}^*$.
(Note that in the case of $\an_{2, 3}$  we only have a $\qn_3$-correspondence.)
Since $\an_{2i}$ is also the $\qn_3, \dots,\qn_{2i-1}$-correspondent of 
$\alpha_{2i}^*$,   while   $\alpha_{2i}^*$ lies above $\alpha_{2i-2}^*$, 
we conclude that $\an_{2i} \in \Irr(\pn_{2i})$ lies above $\an_{2i-2,2i-1}$, 
whenever $1\leq i \leq k$. 
This proves that the set \eqref{bh.e7} satisfies  \eqref{x4}. 

We will work similarly to prove \eqref{x6}, using Corollary   \ref{bh.col1}. 
For every $i=2, \dots, k-1$, the character 
$\bn_{2i-3,2i-2} \in \Irr(\qn_{2i-3,2i-2})= \Irr(C(\pn_{2i-2} \tin \qn_{2i-3}))$ 
is defined as the $\pn_{2i-2}$-Glauberman  correspondent of
$\bn_{2i-3} \tin \Irr(\qn_{2i-3})$.
 The $\pn_{2i-2}, \pn_{2i},\dots,\pn_{2k}$-Glauberman correspondent of 
$\beta_{2i-3}' \in \Irr(Q_{2i-3}')=\Irr(C(\pn_{2i-2}, \pn_{2i},\dots,\pn_{2k}
 \tin \qn_{2i-3}))$ is the character $\bn_{2i-3}$, by  Corollary \ref{bh.col1}.
Hence $\bn_{2i-3,2i-2}$ is the 
$\pn_{2i},\dots,\pn_{2k}$-Glauberman correspondent of  $\beta_{2i-3}'$.
Furthermore, $\bn_{2i-1}$ is the $\pn_{2i},\dots,\pn_{2k}$-Glauberman
correspondent of $\beta_{2i-1}'$. As $\beta_{2i-3}'$ lies under $\beta_{2i-1}'$, 
by Lemma \ref{bh.l2},  we conclude that $\bn_{2i-3,2i-2}$ 
also lies under $\bn_{2i-1}$.  This completes the proof of \eqref{x6},  showing 
that the set \eqref{bh.e7} is a triangular set for \eqref{bh.*1}.
Hence all  the notation and the properties described in Chapter \ref{pq} can be 
 applied to this triangular set.

According to \eqref{pq14b}, the group $\qn_{2i-1,2k}$  equals  
$C(\pn_{2i}\cdot \pn_{2i+2} \cdots \pn_{2k} \tin \qn_{2i-1})$. 
This,  along with Corollary \ref{bh.col1} and \eqref{bh.al2}, implies 
$$
\qn_{2i-1,2k} = C(\pn_{2i}\cdot \pn_{2i+2} \cdots \pn_{2k} \tin \qn_{2i-1})=
C(N(\qn_1,\dots,\qn_{2i-1} \tin P_{2k}^*) \tin \qn_{2i-1}) = Q_{2i-1}'.
$$
Furthermore, according to the Definition \ref{pq1413def2}, the character
$\bn_{2i-1,2k} \in \Irr(\qn_{2i-1,2k})$
 is the $\pn_{2i}\cdot\pn_{2i+2}\cdots \pn_{2k}$-Glauberman 
correspondent of $\bn_{2i-1}$.  Therefore it coincides with $\beta_{2i-1}'$, 
 as the latter is also the $\pn_{2i} \cdots \pn_{2k}$-Glauberman correspondent 
of $\bn_{2i-1}$, by  Corollary \ref{bh.col1}. Hence \eqref{bh.e20a} holds.

If $\pns_{2i}$ denotes the product of the $\pn$-groups, i.e., 
$\pns_{2i}:= \pn_2 \cdots \pn_{2i}$, then 
in view of  \eqref{bh.al9} we get  that 
$\pns_{2i}= P_{2i}^*$, for all $i=1, \dots,k$.
This, along with  Proposition \ref{bh.pr1},  implies the first part of \eqref{bh.e20b}.  The second part  follows easily  from the first and the facts 
that $P_{2i}' = P_{2i}^* (\beta_{2i-1}')$ (see  \eqref{bh.e6d}) while 
$\beta_{2i-1}' =\bn_{2i-1,2k}$  (see  \eqref{bh.e20a}).  

The character $\ans_{2i}\in \Irr(\pns_{2i})=\Irr(P_{2i}^*)$  is constructed 
as the $\qn_3, \dots,\qn_{2i-1}$-correspondent 
of $\an_{2i}$ (see Theorem \ref{p*t1}).
 This,  along with \eqref{bh.al5}, implies \eqref{bh.e20c}.
The relation  \eqref{bh.e20d} now  follows, easily from \eqref{bh.eee}.

This completes the proof of Theorem \ref{bh.th3}.
\end{proof}

%%% Local Variables: 
%%% mode: latex
%%% TeX-master: "p-qgroups"
%%% End: 

\chapter{ The New Characters $\chi_{i}^{\nu}$ of $G_i$ }
 \mylabel{n}

Let $G$ be a finite group 
satisfying 
\begin{equation}\mylabel{n.as1}
|G|= p^a\cdot q^b, \text{ with  $p\ne q$ primes and $a,b$ non negative integers }.
\end{equation}
Assume further that 
\begin{equation}\mylabel{n.1}
1=G_0 \unlhd G_1 \unlhd \cdots \unlhd G_{2k} \unlhd G 
\end{equation}
is   a normal series 
 for $G$ satisfying Hypothesis \ref{hyp1} with $\pi = \{ p \}$, i.e., 
$G_i/G_{i-1}$  is a $p$-group if $i$ is even and a $q$-group if $i$ is odd.
Note also that \eqref{n.1} plays  the role of \eqref{*} with $m=2k$.
Let $\{ 1=\chi_0, \chi_1, \cdots, \chi_{2k} \}$ be a character tower for the series 
\eqref{n.1}.
We have seen in Chapter \ref{pq}, 
Theorem \ref{cor:t}, that there exists a unique, up to conjugation, 
triangular set 
\begin{equation}\mylabel{n.2}
\{Q_{2i-1}, P_{2r} | \beta_{2i-1}, \alpha_{2r} \}_{i=1, r=0}^{k, k} 
\end{equation}
for  \eqref{n.1} that corresponds to the above character tower. 
Of course there is no reason  for the irreducible character $\beta_{2k-1}$ to extend to its 
own stabilizer in $G$. In addition, we have seen 
how to achieve   the irreducible character $\alpha_{2k}^*$ of the product group 
$P_{2k}^* = P_2 \dot P_4 \cdots P_{2k}$, from the irreducible    character $\alpha_{2k}
\in \Irr(P_{2k})$ (see  Definition \ref{p*d1}).  We have also seen how to pick a $q$-Sylow 
subgroup $\qw$ of $G(\alpha_{2k}^*)$  satisfying  all the conditions in Theorem 
\ref{hat:p1}. (Observe that $\pi' = \{ q\}$. )

What we will prove in this chapter is  that,  under the above conditions, 
we can find a new character tower  $\{ 1=\cn_0, \cn_1,\dots,\cn_{2k} \}$ 
for the normal  series \eqref{n.1} 
of $G$ so that a  corresponding triangular set 
$
\{\qn_{2i-1},\pn_{2i}, \pn_0=1 |\bn_{2i-1}, \an_{2i},\an_0=1\}_{i=1}^{k}
$
satisfies \begin{itemize}
\item[1.] $P_{2i}^* = \pns_{2i}$ and $\alpha_{2k}^* = \ans_{2k}$, for all $i=1,\dots,k$,
\item[2.] $G(\alpha_{2k}^*) = G(\ans_{2k})$ and  $\qw= \qwn$,
\item[3.] $\bn_{2k-1,2k}$ extends to $\qw=\qwn$
\end{itemize}
where we keep the same notation as before with teh addition 
 of the superscript $\nu$ to any group that refers to the new character tower and triangular set.
So  $ \pns_{2i}= \pn_2 \cdots \pns_{2i}$, the character $\ans_{2k}$ is an
 irreducible character of $\pns$  that is achieved using the
 character $\an_{2k} \in \Irr(\pn_{2k})$ via  Definition \ref{p*d1} for the new characters.
Furthermore,   $\qwn$ is a  $q$-Sylow  subgroup of  $G(\ans_{2k})$
 that satisfies the conditions in Theorem \ref{hat:p1}.

For this we will put together all the complicated machinery we 
developed in the previous chapters.
We use the same  notation  for the groups and the characters that was 
 introduced  in those chapters, 
but applied in our specific case, i.e., where $G$ satisfies \eqref{n.as1}, 
 and the series \eqref{n.1},  its character tower and the corresponding 
triangular set \eqref{n.2} are fixed. 

Thus we can use all the information about the subgroups $\qw, \qw_{2i-1}$   and $\pw_{2i}$ 
of $G' = G(\alpha_{2k}^*)$,  given in Chapter \ref{pq:sec5}.
So we can prove 
\begin{lemma}\mylabel{lea}
The hypotheses of Theorem  \ref{cc:p1} are satisfied by the present group $G$, with the integer $n$ in the theorem equal to the present $k$, 
the $q$-subgroup $Q_{n+1}= Q$  in the theorem equal to the present $\qw$, 
the $q$-subgroup $Q_i$ in th etheorem equal to the present $\qw_{2i-1}$, 
for all $i=1, 2, \dots, n=k$, and the $p$-subgroup $P_j$ in the theorem 
equal to the present $\pw_{2j}$, for all $j=1, \dots, n=k$. 
\end{lemma}

\begin{proof}
Our present group $G$ has order $p^a q^b$, as required in Theorem \ref{cc:p1}.
By definition $\qw$ is an arbitary subgroup of $G$ satisfying all the conditions in 
Theorem \ref{hat:p1}.
In particular it is a $q$-group,  as it is a $\pi'$-Hall subgroup of 
 of $G(\alpha_{2k}^*)$ and $\pi'=q$.
For each $i=1, 2, \dots, k$,  the subgroup $\qw_{2i-1}$ is the intersection 
$G_{2i-1} \cap \qw$  by Definition \ref{pq40a}.
Since $G_1 \unlhd G_3 \unlhd \dots \unlhd G_{2k-1}$ is a series of normal 
subgroups of $G_{2k-1}$, this implies that 
$\qw_1 \unlhd \qw_3 \unlhd \dots \unlhd \qw_{2k-1}$ is a series 
of normal subgroups of $\qw_{2k-1}$, as required of
 $Q_1 \unlhd Q_2 \unlhd \dots \unlhd Q_n$ in  Theorem  \ref{cc:p1}.

For every $j=1, \dots, k$ the group $\pw_{2j}$ 
was picked  to satisfy the conditions 
in Proposition \ref{pq55}. Hence $\pw$ is a $p$-group, as it is a $\pi$-group and 
$\pi= p$.  Furthermore, according to   \eqref{EqsC} 
 the group $\pw_{2i}$  normalizes 
$\pw_{2j}$, whenever $1 \leq j \leq i \leq n=k$, as required
of $P_1, P_2, \dots, P_n$ in  Theorem  \ref{cc:p1}.
In  addition,  \eqref{EqsA} implies 
 that $\pw_{2i}$ also normalizes $\qw_{2j-1}$, whenever
 $1\leq j \leq i \leq k$.  

  According to Lemma \ref{hat:l1}  and \eqref{EqsA} we get that 
\begin{equation}\mylabel{n.e1}
N(P_{2k}^* , \qw_{2j-1} \tin G_{2j} (\alpha_{2j-2}^*))= \pw_{2j} \ltimes \qw_{2j-1}, 
\end{equation}
for $j=1, \dots, k$.
But  $\qw_{2i-1} \leq G_{2i-1}'=G_{2i-1}(\alpha_{2k}^*)$  by  \eqref{pq40b}.
 Hence  $\qw_{2i-1}$ normalizes $P_{2k}^*$ and fixes $\alpha_{2j}^*$,
 for all $j=1,\dots,i-1$.
This, along with \eqref{n.e1}, implies that $\qw_{2i-1}$ normalizes the 
semidirect product
$\pw_{2j} \ltimes \qw_{2j-1}$, whenever $1 \leq j \leq i \leq k$.
 Similarly, we use 
 $\qw \leq G(\alpha_{2k}^*)$, to see that $\qw$ also normalizes the
above  semidirect product,  $\pw_{2j} \ltimes \qw_{2j-1}$.
Therefore the groups $\qw, \qw_{2i-1}$ and
 $\pw_{2i}$ satisfy the conditions (1) and (2) 
in Theorem \ref{cc:p1} and the lemma follows.
\end{proof}

Lemma \ref{lea} implies
\begin{theorem}\mylabel{ntt1}
There exist linear  characters $\hbn_{2i-1} \in \Lin(\qw_{2i-1})$  such that the
following hold:
\begin{subequations}\mylabel{n.e}
\begin{equation}\mylabel{n.e2}
\hat{\beta}_{2i-1}^{\nu} \in \Irr (\qw_{2i-1} | \hbn_{2i-3},\dots ,\hbn_1),
\end{equation} 
\begin{equation}\mylabel{n.e3}
\qw_{2i-1}(\hbn_{2j-1}) = \qw_{2i-1},
\end{equation}
\begin{equation}\mylabel{n.e4}
\pw_{2i}(\hbn_{2i-1}) =C(\qw_{2i-1} \tin \pw_{2i})=  \pw_{2i}(\hb_{2i-1}),
\end{equation}
and
\begin{equation}\mylabel{n.e5}
\hbn_{2k-1} \text{ extends to an irreducible character of } \qw, 
\end{equation}
 whenever  $1\leq j \leq i \leq k$.  By convention $\hbn_{-1}:= 1$.
\end{subequations}
\end{theorem}

\begin{proof}
According to Lemma \ref{lea}, the groups $\{\qw,  \qw_{2i-1}, \pw_{2i}\}_{i=1}^k$.
satisfy the conditions in Theorem \ref{cc:p1}, with 
the series  $\qw_1 \unlhd \qw_3 \unlhd \dots \unlhd \qw_{2k-1}
\unlhd \qw$ here in the place of  the chain $Q_1 \unlhd Q_2
\unlhd  \dots \unlhd Q_n \unlhd Q_{n+1}=Q$ in Theorem \ref{cc:p1}, and 
 the sequence 
$\pw_2,\pw_4, \dots,  \pw_{2k}$ here, in the place of 
  the sequence  $P_1, P_2, \dots, P_n$ there.  

Note also that   Corollary  \ref{b.cor50}, and in particular \eqref{b.c51},  
provides the additional information that $C(\qw_{2i-1} \tin \pw_{2i}) = 
\pw_{2i}(\hb_{2i-1})$.
This, along with the conclusions (a) and (b) in Theorem \ref{cc:p1}, 
implies the theorem.
\end{proof}

An immediate consequence of  \eqref{n.e3} is 
\begin{remark}\mylabel{n.r1}
For every  $i=1,\dots,k$ and $j=1,\dots,i$, the character  $\hbn_{2j-1}$ is the   unique character
of $\qw_{2j-1}$ that lies under $\hbn_{2i-1} \in \Irr(\qw_{2i-1})$. 
Hence any subgroup that fixes 
$\hbn_{2i-1}$  also fixes $\hbn_{2j-1}$, as it normalizes $\qw_{2j-1} = \qw_{2i-1} \cap G_{2j-1}$.
\end{remark}
In the next lemma we collect some easy remarks that follow from the   properties \eqref{n.e}.

\begin{lemma} \mylabel{n.l1}
For every $i=1,\dots,k$ we have   
\begin{subequations}
\begin{equation}\mylabel{n.e6}
\pw_{2i}(\hbn_{2i-1})=\pw_{2i}(\hbn_1, \cdots ,\hbn_{2i-1}),
\end{equation}
and 
\begin{equation}\mylabel{n.e7}
\pw_{2i}(\hb_1, \cdots ,\hb_{2i-1}) =
C(\qw \tin P_{2i}^*)
= \pw_{2i}(\hb_{2i-1}) = C(\qw_{2i-1} \tin \pw_{2i})  =
\pw_{2i}(\hbn_{2i-1}).
\end{equation}
\end{subequations}
\end{lemma}

\begin{proof}
Equation \eqref{n.e6} follows easily from Remark \ref{n.r1}, and the fact that 
$\pw_{2i}$ normalizes $\qw_{2j-1}$  whenever $1\leq j \leq i \leq k$.
The equation \eqref{n.e4}, along with \eqref{ExtraEq} and \eqref{b.c52}, 
 implies 
\eqref{n.e7}
Therefore the lemma holds.
\end{proof}

As the next proposition shows, the group $\pw_{2i}(\hbn_{2i-1})$ has properties similar to those of
$\pw_{2i}(\hb_{2i-1})$ (see  \eqref{pq54} and \eqref{EqsB}).

\begin{proposition}\mylabel{n.p1}
For all $i=1,\dots,k$ we have 
\begin{multline}\mylabel{n.e8}
\pw_{2i}(\hbn_{2i-1})= \pw_{2i}(\hbn_1, \cdots ,\hbn_{2i-1}) \in \Syl_p(N(P_{2k}^*, \qw_{2i-1} \tin 
G_{2i}(\alpha_{2i-2}^*, \hbn_1, \cdots ,\hbn_{2i-1})))=\\
\Syl_p(N(P_{2k}^* \tin G_{2i}(\alpha_{2i-2}^*, \hbn_{2i-1}))).
\end{multline}
Even more, $\pw_{2i}(\hbn_{2i-1})$ is the unique $p$-Sylow subgroup of 
$N(P_{2k}^* \tin G_{2i}(\alpha_{2i-2}^*, \hbn_{2i-1}))$,  and 
\begin{equation}\mylabel{n.e8.5}
N(P_{2k}^* \tin G_{2i}(\alpha_{2i-2}^*, \hbn_{2i-1}))= \pw_{2i}(\hbn_{2i-1}) \times \qw_{2i-1}=
C(\qw_{2i-1} \tin \pw_{2i}) \times \qw_{2i-1},
\end{equation} 
whenever $1\leq i \leq k$.
\end{proposition}

\begin{proof}
Let $i=1,\dots, k$ be fixed. According to \eqref{pq52} we have 
$$N(P_{2k}^* , \qw_{2i-1}  \tin G_{2i}(\alpha_{2i-2}^*))= \pw_{2i} \ltimes \qw_{2i-1}.
$$
Therefore, $ \pw_{2i}(\hbn_1, \dots,\hbn_{2i-1})$ 
is a $p$-subgroup of $N(P_{2k}^* , \qw_{2i-1}  \tin G_{2i}(\alpha_{2i-2}^*, \hbn_1, \dots,
\hbn_{2i-1}))$.
 Hence there exists  an element  $s \in \qw_{2i-1}$ such that 
\begin{subequations}
\begin{align}
\pw_{2i}^s(\hbn_1, \dots, \hbn_{2i-1}) &\in \Syl_{p}(N(P_{2k}^* , \qw_{2i-1}  \tin G_{2i}(\alpha_{2i-2}^*, \hbn_1, \dots,
\hbn_{2i-1}))),  \text {and } \mylabel{n.e9}   \\
\pw_{2i} (\hbn_1, \dots, \hbn_{2i-1}) &\leq \pw_{2i}^s(\hbn_1, \dots, \hbn_{2i-1}). \mylabel{n.e9.1}
\end{align}
\end{subequations}
But $\qw_{2i-1}$ fixes $\hbn_1, \dots, \hbn_{2i-1}$, by  \eqref{n.e3}.
Thus $\pw_{2i}^s(\hbn_1, \dots, \hbn_{2i-1}) = 
(\pw_{2i}(\hbn_1, \dots, \hbn_{2i-1}))^s$. A cardinality argument, along with \eqref{n.e9.1}, 
 implies that 
$$
\pw_{2i}(\hbn_1, \dots, \hbn_{2i-1})=\pw_{2i}^s(\hbn_1, \dots, \hbn_{2i-1}) =
(\pw_{2i}(\hbn_1, \dots, \hbn_{2i-1}))^s.
$$
This,  along with \eqref{n.e9} and \eqref{n.e6}, implies that 
$$
\pw_{2i}(\hbn_{2i-1}) = \pw_{2i}(\hbn_1, \dots, \hbn_{2i-1}) \in \Syl_{p}(
N(P_{2k}^* , \qw_{2i-1}  \tin G_{2i}(\alpha_{2i-2}^*, \hbn_1, \dots,
\hbn_{2i-1}))).
$$
In view of Remark \ref{n.r1}, every subgroup of $G$ 
 that fixes $\hbn_{2i-1}$ fixes  $\hbn_{2j-1} \in 
\Irr(\qw_{2j-1})$, for all $1\leq j \leq i $. Hence
 $N(P_{2k}^* , \qw_{2i-1}  \tin G_{2i}(\alpha_{2i-2}^*, \hbn_1, \dots,
\hbn_{2i-1}))= N(P_{2k}^* \tin G_{2i}(\alpha_{2i-2}^*, \hbn_{2i-1}))$.
Therefore \eqref{n.e8}  holds.

According to Proposition \ref{pqrm}, 
we have  $\qw_{2i-1} \in \Syl_{q}(N(P_{2k}^* \tin 
G_{2i}(\alpha_{2i-2}^*)))$ . Hence 
$\qw_{2i-1}  \in \Syl_{q}(N(P_{2k}^* \tin G_{2i}(\alpha_{2i-2}^*, \hbn_{2i-1})))$.
This, along with \eqref{n.e8}, implies that 
\begin{equation}\mylabel{n.e10}
N(P_{2k}^* \tin G_{2i}(\alpha_{2i-2}^*, \hbn_{2i-1}))= \pw_{2i}(\hbn_{2i-1})\ltimes \qw_{2i-1}.
\end{equation}
But according to   \eqref{n.e4}   we have   
$\pw_{2i} (\hbn_{2i-1} ) = C(\qw_{2i-1} \tin \pw_{2i})$.
This,  along with \eqref{n.e10}, implies \eqref{n.e8.5}.
Hence Proposition \ref{n.p1}  holds.
\end{proof}

\begin{defn}\mylabel{n.d1}
For every $i=1,\dots,k$ we 
define $\ha_{2i}  \in \Irr(\pw_{2i}(\hb_{2i-1})$ to be the $\qw_{2i-1}$-Glauberman correspondent
of $\alpha_{2i}^* \in \Irr(P_{2i}^*)$.
\end{defn}
Note that the  $\ha_{2i}$ are well defined,  as $\alpha_{2i}^* \in \Irr(P_{2i}^*)$ is $\qw_{2i-1}$-invariant
(since $\qw_{2i-1} \leq G(\alpha_{2k}^*)$) and $\pw_{2i}(\hb_{2i-1}) = C(\qw_{2i-1} \tin P_{2i}^*)$
 (see \eqref{n.e7}).

We can now prove the main theorem of this section.
\begin{theorem}\mylabel{n.t1}
The set 
\begin{equation}\mylabel{n.e12}
\{ \qw_{2i-1}, \pw_{2i}(\hb_{2i-1})| \hbn_{2i-1}, \ha_{2i} \}_{i=1}^k
\end{equation}
is a triangular set for the series 
$1 = G_0' \unlhd G_1' \unlhd  \dots \unlhd G_{2k}' \unlhd G'$ in  \eqref{bh.*}.
 Furthermore,  it  satisfies  Property \ref{bh.pp1},  while 
\eqref{bh.e5b} holds for this triangular set, with $Q' = \qw$.
\end{theorem}

\begin{proof}
We will first prove that \eqref{n.e12} is a triangular set for the above series,
 i.e., we will verify the properties \eqref{xxx}
for that set and series. Assume that $i=1,\dots,k$ is fixed.

The equations in \eqref{x1}  hold trivially. As
 $\qw_{1} = Q_{1,2k}=G_1'$ (see \eqref{pqd40a})
and $\hbn_1 \in \Irr(\qw_1)$,  \eqref{x2} holds, while \eqref{x5} is trivially true
 for $i=1$.

Assume that $i \geq 2$. 
According to \eqref{pq40b} the group $\qw_{2i-1}$ is a
 $q$-Sylow subgroup of $G_{2i-1}'$. 
Furthermore, $\qw_{2i-1}$ normalizes $\qw_{2j-1}= \qw_{2i-1} \cap G_{2j-1}$, 
and fixes $\hbn_{2j-1} \in \Irr(\qw_{2j-1})$, by  \eqref{n.e3},  for all $j=1,\dots,i-1$. 
The group $\qw_{2i-1}$ also fixes $\alpha_{2j}^*$,  as it fixes $\alpha_{2k}^*$, for all such $j$. 
Therefore, it fixes the $\qw_{2j-1}$-Glauberman correspondent 
$\ha_{2j} \in \Irr(\pw_{2j}(\hbn_{2j-1}) )= \Irr(C(\qw_{2j-1} \tin P_{2j}^*)$ 
 of $\alpha_{2j}^*$in Definition
 \ref{n.d1},   for all $j=1,\dots,i-1$.
Hence 
$$\qw_{2i-1} \leq G_{2i-1}' (\ha_2, \dots,\ha_{2i-2},\hbn_1,\dots,\hbn_{2i-3}) \leq G_{2i-1}'.$$
As  $\qw_{2i-1}  \in \Syl_{q}(G_{2i-1}')$ we get that 
$$
\qw_{2i-1} \in \Syl_{q}( G_{2i-1}' (\ha_2, \dots,\ha_{2i-2},\hbn_1,\dots,\hbn_{2i-3})).
$$
Hence \eqref{x5} holds.

As $\pw_{2i-2}(\hbn_{2i-3})=C(\qw_{2i-3} \tin \pw_{2i-2})$ 
by \eqref{n.e4}, we have that 
$C(\pw_{2i-2}(\hbn_{2i-3}) \tin \qw_{2i-3})=\qw_{2i-3}$. 
Therefore, the character $\hbn_{2i-3,2i-2} \in \Irr(C(\pw_{2i-2}(\hbn_{2i-3}) 
\tin \qw_{2i-3}))$
concides with $\hbn_{2i-3}$.   
This, along with \eqref{n.e2}, makes  the condition \eqref{x6} valid.

For the $p$-groups and characters we have that $\pw_{2i}(\hbn_{2i-1})=C(\qw_{2i-1} \tin \pw_{2i})$ 
fixes $\alpha_{2j}^*$, for all $j=1,\dots, k$, 
as it is a subgroup of $G'=G(\alpha_{2k}^*)$.
 It also centralizes $\qw_{2i-1}$, and thus
 centralizes $\qw_{2j-1} \leq \qw_{2i-1}$ for all $j=1, \dots, i$.
Therefore $\pw_{2i}(\hbn_{2i-1})$ fixes 
the $\qw_{2j-1}$-Glauberman  correspondent 
 $\ha_{2j}$ of $\alpha_{2j}^*$, for all $j=1,\dots,i$. 
In view of Remark \ref{n.r1}, we also have that $\pw_{2i}(\hbn_{2i-1})$ 
 fixes $\hbn_{2j-1}$ for all such $j$.
Hence 
\begin{equation}\mylabel{n.e13}
  \pw_{2i}(\hbn_{2i-1}) \leq G_{2i}(\alpha_{2k}^*, \ha_2,\dots,\ha_{2i-2}, \hbn_1,\dots,\hbn_{2i-1})
\leq N(P_{2k}^* \tin G_{2i}(\alpha_{2i-2}^*, \hbn_{2i-1})).
\end{equation}
According to Proposition \ref{n.p1}, the unique $p$-Sylow subgroup of
$ N(P_{2k}^* \tin G_{2i}(\alpha_{2i-2}^*, \hbn_{2i-1}))$ is  the group 
 $\pw_{2i}(\hbn_{2i-1})$.
This,  along with \eqref{n.e13} and the fact that $G_{2i}' =G_{2i}(\alpha_{2k}^*)$, implies that 
$$
\pw_{2i}(\hbn_{2i-1}) \in \Syl_{p}( G_{2i}(\alpha_{2k}^*, \ha_2,\dots,\ha_{2i-2},
 \hbn_1,\dots,\hbn_{2i-1})) = \Syl_{p}( G_{2i}'( \ha_2,\dots,\ha_{2i-2},
 \hbn_1,\dots,\hbn_{2i-1})).
$$
Hence \eqref{x3} holds.

To prove \eqref{x4}, we first observe that 
$\qw_{2i-1}$ normalizes both $\qw_{2i-3} = G_{2i-3} \cap \qw_{2i-1}$ and 
$P_{2i-2}^* = P_{2k}^* \cap G_{2i-2}$, since $\qw{2i-1} \leq G'= G(\alpha_{2k}^*)$
normalizes $P_{2k}^*$. Hence the $q$-group $\qw_{2i-1}$ normalizes 
the $p$-group $C(\qw_{2i-3} \tin P_{2i-2}^*)$ . Hence 
$$
C(\qw_{2i-1} \tin P_{2i-2}^*) = C(\qw_{2i-1} \tin C(\qw_{2i-3} \tin P_{2i-2}^*)).
$$
 By convention we set $\qw_{-1}:=1$, so that  the above equation holds 
trivially for $i=1$.  Let  $\ha_{2i-2,2i-1} \in \Irr(C(\qw_{2i-1} \tin P_{2i-2}^*))$ denote the 
$\qw_{2i-1}$-Glauberman correspondent 
of the irreducible character  $\ha_{2i-2} \in \Irr(C(\qw_{2i-3} \tin P_{2i-2}^*))$.
(Note that 
since $\qw_{2i-1}$ is a subgroup of  $
G' = G(\alpha_{2k}^*)$  it fixes $\alpha_{2i-2}^* \in 
\Irr(P_{2i-2}^*)$  and normalizes $\qw_{2i-3}$, so it fixes 
the $\qw_{2i-3}$-Glauberman correspondent  $\ha_{2i-2}$ 
 of  $\alpha_{2i-2}^*$.)
Then $\ha_{2i-2, 2i-1}$ is the $\qw_{2i-1}$-Glauberman correspondent of $\alpha_{2i-2}^*$.  
Hence $\ha_{2i-2,2i-1}$ lies under the $\qw_{2i-1}$-Glauberman correspondent
$\ha_{2i}$ of $\alpha_{2i}^*$, as $\alpha_{2i-2}^*$ lies under $\alpha_{2i}^*$.
Therefore $\ha_{2i} \in \Irr(\pw_{2i}(\hbn_{2i-1}) | \ha_{2i-2,2i-1})$.  So \eqref{x4} is satisfied.
This completes the proof of \eqref{xxx}. Hence \eqref{n.e12} is a triangular set for \eqref{bh.*}.

The group $\qw_{2i-1}= \qw \cap G_{2i-1}$ is clearly  a normal subgroup of $\qw$, for all $i=1,\dots, k$. 
Furthermore, $\qw$ fixes the character $\hbn_{2k-1}$, by \eqref{n.e5}.
Hence,  Remark \ref{n.r1} implies that  $\qw$ fixes $\hbn_{2i-1}$, for all $i=1,\dots, k$.
As $\qw $ is a subgroup of $G'=G(\alpha_{2k}^*)$, it fixes 
$\alpha_{2i}^*$,  for all $i=1,\dots, k$. 
Since $\qw$ normalizes $\qw_{2i-1}$, it fixes 
the $\qw_{2i-1}$-Glauberman correspondent $\ha_{2i}$  of 
$\alpha_{2i}^*$, for all such $i$. 
Thus $\qw \leq G'(\ha_2, \dots, \ha_{2k}, \hbn_1, \dots, \hbn_{2k-1})$.
So $\qw$ satisfies \eqref{bh.e5b}. 

Definition \ref{n.d1} implies that the triangular set \eqref{n.e12} satisfies \eqref{bh.eee}.
It also satisfies \eqref{bh.ee}, according to Proposition \ref{n.p1}.
As we have already seen,  the subgroup $\qw_{2i-1}$ of $G'$   normalizes $P_{2i}^*$.
Hence 
$$
P_{2i}^*(\hbn_{2i-1}) \leq N(\qw_{2i-1} \tin P_{2i}^*) = C(\qw_{2i-1} \tin P_{2i}^*) \leq 
P_{2i}^*(\hbn_{2i-1}).
$$
So  $P_{2i}^*(\hbn_{2i-1}) =N(\qw_{2i-1} \tin P_{2i}^*) 
= C(\qw_{2i-1} \tin P_{2i}^*)$.
This,  along with \eqref{n.e4}, implies that 
$$
\pw_{2i}(\hbn_{2i-1}) =N(\qw_{2i-1} \tin P_{2i}^*) = C(\qw_{2i-1} \tin P_{2i}^*)= P_{2i}^*(\hbn_{2i-1}).
$$ 
Thus \eqref{bh.e} also holds.
Hence the set \eqref{n.e12} satisfies Property \ref{bh.pp1}.
This completes the proof of the theorem.
\end{proof}

All the work in Chapters 4-6 was done to prove the following theorem
\begin{theorem}\mylabel{n.t2}
Let $\{1=\chi_0, \chi_1, \dots, \chi_{2k} \}$ be 
a character tower for the series \eqref{n.1},
and let    \eqref{n.2} be its unique,  up to conjugation, 
 corresponding triangular set. 
Then there exists a character tower 
$\{1=\cn_0, \cn_1, \dots, \cn_{2k} \}$ for the series 
\eqref{n.1},  with corresponding  triangular set
$
\{\qn_{2i-1},\pn_{2i}, \pn_0=1 |\bn_{2i-1}, \an_{2i},\an_0=1\}_{i=1}^{k}
$, 
so that the following hold
\begin{subequations}\mylabel{n.e14}
\begin{align}
\qw_{2i-1} = \qn_{2i-1,2k} &\text{ and } \hbn_{2i-1} = \bn_{2i-1,2k},
\mylabel{n.e14a}\\
P_{2i}^*  = \pns_{2i} &\text{ and } \alpha_{2i}^*= \ans_{2i}, \mylabel{n.e15}\\
\bn_{2i-1} \text{ is the } 
 \pn_{2i} \cdot \pn_{2i+2} \cdots \pn_{2k}&\text{-Glauberman correspondent of } \hbn_{2i-1},
 \mylabel{n.e14b}  \\
\an_{2i} \text{ is the }  &\qn_3, \dots, \qn_{2i-1}\text{-correspondent of }  \alpha_{2i}^*,
 \mylabel{n.e14c} \\
\qn_{2i-1}  &\geq Q_{2i-1,2k}, \mylabel{n.e16}\\
\qw &\text{ normalizes } \qn_{2i-1}, \mylabel{n.e16b}\\ 
\bn_{2k-1,2k}  &\text{ extends  to } \qw, 
\mylabel{n.e17}  
\end{align}
\end{subequations}
whenever $1\leq i \leq k$.
\end{theorem}

\begin{proof}

For the  fixed triangular set \eqref{n.2} of the character tower 
 $\{1=\chi_0, \chi_1, \dots, \chi_{2k} \}$ 
we saw in  Chapter \ref{pq:sec5}  how to pick groups $\qw_{2i-1}$ 
and $ \pw_{2i}$, 
 along with  characters $\hb_{2i-1} \in \Irr(\qw_{2i-1})$
satisfying all the conditions in Theorem \ref{hat:p1} 
Proposition  \ref{pq55}  and Proposition \ref{pq40x}.
Furthermore, we proved at the begining of this chapter 
 that we can  replace  the characters $\hb_{2i-1}$ 
with new characters $\hbn_{2i-1} \in \Irr(\qw_{2i-1})$ 
that satisfy  \eqref{n.e}. Even more, as Theorem \ref{n.t1} shows,
 the set \eqref{n.e12}
is a triangular set for \eqref{bh.*}  that satisfies Property \ref{bh.pp1}, 
while \eqref{bh.e5b} holds for this triangular set, with 
 $\qw$ in the place of $Q'$.
According to Theorems \ref{bh.th2} and \ref{bh.th3}, 
the set \eqref{n.e12} determines a triangular set 
 \begin{equation}\mylabel{n.4}
\{\qn_{2i-1},\pn_{2i}, \pn_0=1 |\bn_{2i-1}, \an_{2i},\an_0=1\}_{i=1}^{k}
\end{equation}
for the series \eqref{n.1}, such that  \eqref{bh.al} 
and \eqref{bh.e20} hold with $\qw, \qw_{2i-1}$  and $\pw_{2i}(\hb_{2i-1})$
 in the place of $Q', Q_{2i-1}'$ and $P_{2i}'$, respectively.
In view of Theorem \ref{cor:t}, the triangular
 set \eqref{n.4} corresponds to a unique, up to conjugation, 
character tower $\{1=\cn_0, \cn_1, \dots, \cn_{2k} \}$ for the series 
\eqref{n.1}.

To complete the proof of the theorem it suffices to show that
 the set \eqref{n.4} has 
the properties \eqref{n.e14}.
As it satisfies \eqref{bh.e20}, the equations in 
 \eqref{n.e15} follow trivially from the first part of 
\eqref{bh.e20b}, and \eqref{bh.e20c}.
 It also satisfies \eqref{n.e14a} and \eqref{n.e14c} as it satisfies
\eqref{bh.e20a} and \eqref{bh.al5} respectively.
In addition  \eqref{n.e14b} holds, since the set \eqref{n.4} satisfies 
     \eqref{bh.al3} and Corollary \ref{bh.col1}.
Furthermore, according to \eqref{bh.al2}, the group $\qn_{2i-1}$ contains $\qw_{2i-1}$,  
as $\qn_{2i-1}(\alpha_{2k}^*) = \qw_{2i-1}$. 
In addition,  \eqref{pq40g} implies that  $\qw_{2i-1}  \geq  Q_{2i-1,2k}$.
We conclude that $\qn_{2i-1} \geq Q_{2i-1,2k}$ and \eqref{n.e16}  follows.

Clearly \eqref{bh.al1} implies that 
 $\qw$ normalizes $\qn_{2i-1}$, for all $i=1,\dots, k$. Thus \eqref{n.e16b}
holds.

The last part, \eqref{n.e17}, follows easily from \eqref{n.e5}, 
as $\hbn_{2k-1} = \bn_{2k-1,2k}$ by \eqref{n.e14a}.
\end{proof}

Furthermore, the new triangular set shares one more group 
 with the old one, as the next theorem shows. 

\begin{theorem} \mylabel{n.t3}
The group $\qw$ satisfies the conditions in 
 Theorem \ref{hat:p1} for the new groups.
Hence we may assume that $\qwn = \qw$. 
Then 
\begin{equation} \mylabel{n.e20}
\qwn = \qw = \qw(\bn_{2k-1,2k})= \qwn(\bn_{2k-1,2k}).
\end{equation}
\end{theorem}

\begin{proof}
For the proof 
we need to show that $\qw$ satisfies \eqref{pq30} and \eqref{pq31b}
 for the new groups.
Clearly, $\qw$ is a $\pi'$-Hall subgroup of $G(\ans_{2k})$,
 as $\ans_{2k} = \alpha_{2k}^*$, by \eqref{n.e15}, and $\qw$ is a $\pi'$-Hall 
subgroup of $G(\alpha_{2k}^*)$.
Thus \eqref{pq30a} holds.

In view of\eqref{n.e5}
 the character $\hbn_{2k-1}$ is fixed by $\qw$.
Thus,  Remark \ref{n.r1} implies that  $\qw$ fixes $\hbn_{2i-1}$, for 
all $i=1,\dots, k$. 
Hence $\qw$ fixes $\bn_{2i-1,2k}=\hbn_{2i-1}$, for all such $i$.  
Furthermore, $\qw$ normalizes the groups $\qn_{2i-1}$, by \eqref{n.e16b}, 
 and $\pns_{2i}= P_{2i}^*= P_{2k}^* \cap G_{2i}$, as it fixes 
$\alpha_{2k}^* \in \Irr(P_{2k}^*)$, whenever $1\leq i \leq k$.
But \eqref{pq22ii},  applied to the new groups,  implies that
$\pn_{2i}= N(\qn_1, \dots, \qn_{2i-1} \tin \pns_{2i})$, for all such $i$ 
We conclude that $\qw$ normalizes $\pn_{2i}$, for all  $i=1, \dots, k$.
Hence $\qw$ fixes the $\pn_{2i} \cdot \pn_{2i+2} \cdots \pn_{2k}$-Glauberman
 correspondent 
$\bn_{2i-1}  \in \Irr(\qn_{2i-1})$ of $\hbn_{2i-1}$  (see \eqref{n.e14b}), 
 as it fixes $\hbn_{2i-1}$.
So 
$$
\qw = \qw(\bn_{2i-1, 2k}) = \qw(\bn_{2i-1}),
$$
whenever $1\leq i \leq k$.

Even more, $\qw$ fixes $\alpha_{2i}^*$, as it fixes $\alpha_{2k}^*$.
We saw above that it 
and normalizes $\qn_{2i-1}$ and $\pn_{2i}$, for all $i=1,\dots, k$. 
Thus $\qw$ fixes the   $\qn_3, \dots, \qn_{2i-1}$-correspondent 
$\an_{2i}$ of $\alpha_{2i}^*$  in  \eqref{n.e14c}.
Hence
$$
\qw(\an_{2i}) = \qw,
$$
whenever $1\leq i \leq k$.

According to Theorem \ref{tow--tri}, the $c\pn_2,\dots, c\pn_{2i}, 
c\qn_1, \dots, c\qn_{2i-1}$-correspondent of $\cn_{2i} \in \Irr(G_{2i})$
is the character $\cn_{2i, 2i} = \an_{2i} \times \bn_{2i-1,2i}$, for 
all $i=1,\dots, k$. As we have already seen, 
the group $\qw$ fixes the characters $\an_{2i}$ and $\bn_{2i-1}$, and 
normalizes the groups $\qn_{2j-1}$ and $\pn_{2j}$ for $j=1, \dots, i$.
 Thus it also fixes both 
 the $\pn_{2i}$-Glauberman correspondent $\bn_{2i-1,2i}$  of 
  $\bn_{2i-1}$, and the direct product 
$\cn_{2i, 2i} = \an_{2i} \times \bn_{2i-1,2i}$.  
Therefore, $\qw$ also fixes  the  $c\pn_2,\dots, c\pn_{2i}, 
c\qn_1, \dots, c\qn_{2i-1}$-correspondent $\cn_{2i}$ of $\cn_{2i, 2i}$, 
for all $i=1,\dots, k$  (see Diagram \ref{diagr.1-2-3...}
 applied to  the new characters).
Similarly, we can see that $\qw$ fixes $\cn_{2i-1, 2i-1}= \an_{2i-2,2i-1} 
\times \bn_{2i-1}$, as well as the 
 $c\pn_2,\dots, c\pn_{2i-2}, 
c\qn_1, \dots, c\qn_{2i-1}$-correspondent $\cn_{2i-1}$ of $\cn_{2i-1,2i-1}$.

In conclusion, 
\begin{equation}\mylabel{n.e21}
\qw = \qw(\bn_{2i-1, 2k}) = \qw(\bn_{2i-1})= \qw(\an_{2i})=\qw(\cn_{2i}) =
\qw(\cn_{2i-1}),
\end{equation}
whenever $1\leq i \leq k$.
 
It is clear  that \eqref{pq31a} and \eqref{hat:p1e} hold, with the new 
characters $\bn_{2i-1}, \bn_{2i-1,2k}, \an_{2i}, \cn_{2i}, \cn_{2i-1}$ 
and $\an_{2i}$  in the place of the analogous original characters. 
(Actually, in \eqref{hat:p1e} we have equality.)
Furthermore, \eqref{n.e21} also implies that 
the group $\qw = \qw(\bn_{2i-1,2k})$ is contained in 
$G'(\bn_{2i-1,2k}) \cap G'(\cn_1, \dots, \cn_{2i-1}) \cap
G'(\cn_1,\dots, \cn_{2i}) \cap G'(\bn_1, \dots, \bn_{2i-1})$, 
Thus it is a $\pi'$-Hall  subgroup of each group in this   intersection,  
as it is a $\pi'$-Hall subgroup of $G'$. 
 Hence  $\qw$ satisfies (\ref{pq30}b, c, d)  for the new groups.

It remains to show \eqref{pq31b}. 
But, as \eqref{n.e16b} implies, 
$\qw = \qw(\bn_{2i-1,2k})$ normalizes $\qn_{2i+1}$,
for all $i=1,\dots, k-1$.
Thus \eqref{pq31b} holds. 
This completes the proof of the theorem. 
\end{proof}

An easy consequence is 
\begin{corollary} \mylabel{n.co1}
Let $\{1=\chi_0, \chi_1, \dots, \chi_{2k} \}$ be a character tower
 for the series \eqref{n.1},
and let    \eqref{n.2} be its unique,  up to conjugation,  
corresponding triangular set. Assume further that $\qw$  satisfies 
the conditions in  Theorem \ref{hat:p1} for this set and tower. 
Then there exist a character tower $\{1=\cn_0, \cn_1, \dots, \cn_{2k} \}$ 
for the series 
\eqref{n.1}, with corresponding  triangular set
$
\{\qn_{2i-1},\pn_{2i}, \pn_0=1 |\bn_{2i-1}, \an_{2i},\an_0=1\}_{i=1}^{k}
$,  and  a group  $\qwn$  that satisfy 
\begin{subequations}\mylabel{n.e22}
\begin{align}
P_{2k}^*= \pns_{2k} &\text{ and } \alpha_{2k}^*= \ans_{2k}, \\
\qw &= \qwn, \\
\qwn &\text{ fixes the characters  } \an_{2i}, \bn_{2i-1}, \cn_{j}, \\
\bn_{2i-1,2k} &\text{ extends to  } \qw=\qwn,
\end{align}
for all $i=1,\dots, k$ and $j=1,\dots, 2k$.
\end{subequations}
\end{corollary}

\begin{proof}
Follows immediately from Theorems \ref{n.t2} and \ref{n.t3}. 
\end{proof}

If instead of the series \eqref{n.1}, we consider  the bigger series
\begin{equation}\mylabel{n.3}
1=G_0 \unlhd G_1 \unlhd \dots \unlhd G_{2k} \unlhd G_{2k+1} \unlhd G,
\end{equation}
then the conclusions of Corollary \ref{n.co1} still  hold, i.e., 
\begin{corollary} \mylabel{n.co2}
Let $\{1=\chi_0, \chi_1, \dots, \chi_{2k+1} \}$ be 
a character tower for the series \eqref{n.3},
and let   $\{Q_{2i+1}, P_{2i} |\beta_{2i+1}, \alpha_{2i} \}_{i=0}^k$
 be its unique,  up to conjugation,  corresponding triangular set. 
Assume further that $\qw$ is picked to satisfy the conditions in Theorem 
\ref{hat:p1} for this set and tower.  
Then there exist a character tower
 $\{1=\cn_0, \cn_1, \dots, \cn_{2k} \}$ for the series 
\eqref{n.1},  with corresponding  triangular set
$
\{\qn_{2i-1},\pn_{2i}, \pn_0=1 |\bn_{2i-1}, \an_{2i},\an_0=1\}_{i=1}^{k}
$,  and a  group  $\qwn$  that satisfy 
\begin{subequations}\mylabel{n.e222}
\begin{align}
P_{2k}^*= \pns_{2k} &\text{ and } \alpha_{2k}^*= \ans_{2k}, \\
\qw &= \qwn, \\
\qwn &\text{ fixes the characters  } \an_{2i}, \bn_{2i-1}, \cn_{j}, \\
\bn_{2i-1,2k} &\text{ extends to  } \qw=\qwn,
\end{align}
for all $i=1,\dots, k$ and $j=1,\dots, 2k$.
\end{subequations}
\end{corollary}

\begin{proof}
Follows easily from  Corollary \ref{b.cB} 
and  Corollary  \ref{n.co1}.
\end{proof}
%%% Local Variables: 
%%% mode: latex
%%% TeX-master: "thesis-ex"
%%% End: 

\chapter{ The  $\pi, \pi'$ Symmetry and the Hall System $\{\ma, \mb\}$ }
\mylabel{sy}

\section{The group $\hap$ }
 
Let $G$ be a finite group of odd order.
As we saw in Chapter \ref{pq:sec5}, whenever we fix a normal  series 
$1=G_0 \unlhd \dots \unlhd G_m \unlhd G$ of $G$ that satisfies Hypothesis 
\ref{hyp1}, a character tower $\{ \chi_i \in \Irr(G_i) \}_{i=0}^m$ for 
 this series and its  corresponding triangular set,
then we can get a $\pi$-Hall subgroup $\qw$ of $G(\alpha_{2k}^* )$ with the
properties described in Theorem  \ref{hat:p1}. We also saw how to get the 
$\pi$-groups $\pw$. 
Furthermore, we used $\qw$ and $\pw$ in 
 Chapter \ref{n},   to replace 
the given character tower with another one having  the properties 
 described in Corollary \ref{n.co1}.

Of course the $\pi-\pi'$  symmetry in the construction of the triangular sets
 implies that  results 
similar to those for the $\pi'$-groups also hold for the $\pi$-groups. 
That is, whenever the above series, the character tower and its 
triangular set are fixed,
  we  can find a $\pi$-Hall subgroup 
$\hap$ of $G'' = G(\beta_{2l-1}^*)$ that satisfies
a modification of  Theorem \ref{hat:p1}, that is,  
\begin{subequations}\mylabel{sy.e}
\begin{align}
\hap \in \Hall_{\pi}(G(\beta_{2l-1}^*)), \\
 \hap(\alpha_{2i,2l-1}) \in \Hall_{\pi}(G''(\alpha_{2i, 2l-1}))
&\cap
\Hall_{\pi}(G''(\chi_1,\dots,\chi_{2i}))\cap \notag\\
&\Hall_{\pi}(G''(\chi_1,\dots,\chi_{2i+1})) \cap
\Hall_{\pi}(G''(\alpha_2, \dots, \alpha_{2i})),\\
 \hap(\alpha_{2i, 2l-1}) =\hap(\chi_1,\dots,\chi_{2i}) &= 
\hap(\chi_1,\dots,\chi_{2i+1}) =\hap(\alpha_2, \dots, \alpha_{2i})
 \text{ and }\\
 \hap(\chi_1,\dots,\chi_{2i}) &\leq \hap(\beta_1, \dots, \beta_{2i+1}),
\end{align}
for all $i=1,\dots, l-1$. 
Furthermore, 
\begin{equation}\mylabel{sy.ee}
\hap(\alpha_{2i, 2l-1}) \text{ normalizes } P_{2i+2},
\end{equation}
\end{subequations}
 for all $i=0, 1,\dots, k-1$. 

In the particular case of a $p^a q^b$-group  $G$ (where $p\ne q$ 
are odd primes), 
we get the analogue of Corollary \ref{n.co1}  for the $\pi$-groups, 
interchanging the roles of $p$ and $q$, that is, 
\begin{theorem}\mylabel{n.cop}
Let $\{1=\chi_0, \chi_1, \dots, \chi_{2k} \}$ be a character tower
 for the series $1=G_0\unlhd G_1 \unlhd \dots \unlhd G_{2k}$, 
and let   $\{Q_{2i-1}, P_{2i}, P_0=1 
|\beta_{2i-1}, \alpha_{2i}, \alpha_0 =1 \}_{i=1}^{k}$
be its unique,  up to conjugation,  
corresponding triangular set. 
Then there exist a character tower 
$\{1=\cn_0, \cn_1, \dots, \cn_{2k-1} \}$
 for the series  $1=G_0 \unlhd G_1 \unlhd \dots \unlhd G_{2k-1}$, 
 a  corresponding  triangular set
$
\{\qn_{2i-1},\pn_{2i-2} |\bn_{2i-1}, \an_{2i-2}\}_{i=1}^{l=k}
$,  and a  $p$-group  $\hapn$,  that satisfy 
\begin{align*}
Q_{2l-1}^*= \qns_{2l-1} &\text{ and } \beta_{2l-1}^*= \bns_{2l-1}, \\
\hap &= \hapn, \\
\hapn \text{ fixes the}  &\text{ characters  } \an_{2i}, \bn_{2i-1}, \cn_{j},
 \text{ and }  \\
\an_{2i, 2l-1} &\text{ extends to  } \hap=\hapn,
\end{align*}
for all $i=1,\dots, l-1$ and $j=1,\dots, 2k-1$.
\end{theorem}

\section{The Hall system $\{\ma, \mb \}$ of $G$ }

Let $G$ be any finite group of odd order, and $\pi$ any
 set of primes. 
If $\ma \in \Hall_{\pi}(G)$ and $\mb \in \Hall_{\pi'}(G)$, then we
 call the set
$\{\ma, \mb \}$ a \emph{ Hall $\pi, \pi'$-system } for $G$, or, more shortly, 
 a  \emph{ Hall system} for $G$. Note that $G$ has a single conjugacy class of 
such Hall systems, because it is solvable.
 Furthermore, if $H $ is a subgroup of $G$,  
we say that the Hall system $\ma, \mb$ of $G$ \emph{reduces} into $H$,
if $\ma \cap H, \mb \cap H$ form a Hall system for $H$. 

We start with a finite odd order group $G$, and we
 fix an increasing chain 
\begin{subequations}\mylabel{sy.e1}
\begin{equation}\mylabel{sy.A}
   1 = G_0 \trianglelefteq G_1 \trianglelefteq 
G_2 \trianglelefteq \dots \trianglelefteq G_n = G, 
\end{equation}
of normal subgroups $G_i$ of $G$,
that satisfy Hypothesis \ref{hyp1} with $n > 0$  in the place of $m$, i.e., 
 $G_i/G_{i-1}$ is a 
$\pi$-group if $i$ is even, and a $\pi'$-group if $i$ is odd,
 for each  $i = 1,2,\dots,n$.
 We also fix a character tower 
\begin{equation}\mylabel{sy.B}
\{ \chi_i  \in \Irr(G_i) \}_{i=0} ^n
\end{equation}
for the above series.
   
 We denote by $k'$ and $l'$ the integers
$$
 k' = [n/2] \text{ and } l' = [(n+1)/2]
$$
corresponding to $k$ and $l$ in \eqref{kl:def}, with $n$ in place of $m$.
So $2k'$ and $2l'-1$ are the greatest even and odd integers, respectively,
in the set $\{ 1,2,\dots,n\}$. As in Section 5.3, we construct a
triangular set 
\begin{equation}\mylabel{sy.e1c}
\{ P_{2r}, Q_{2i-1} | \alpha_{2r}, \beta_{2i-1}\}_{r =0 ,
i=1}^{k',l'}
\end{equation}
\end{subequations}
 corresponding to the chain \eqref{sy.A} and tower \eqref{sy.B}.
 It, in
turn, determines the groups $P^*_{2r}$ and $Q^*_{2i-1}$, for $r=
1,2,\dots,k'$ and $i = 1,2,\dots,l'$.  We know by  Corollary \ref{q*:c1}
  that $P^*_{2k}$ and
$Q^*_{2l-1}$ form a Hall system for $G_m(\chi_1,\chi_2,\dots, \chi_m)$,
whenever $m = 1,2,\dots,n$ and $k$, $l$ are related to $m$ by the usual
equations in \eqref{kl:def}. 
In particular, $P^*_{2k'}$ and 
$Q^*_{2l'-1}$
form a Hall system for 
$G(\chi_1, \chi_2, \dots, \chi_n) = G_n(\chi_1,
\chi_2,\dots,\chi_n)$. 
Furthermore, the groups $G(\chi_1,\dots,\chi_m)$, for $m=1,\dots,n$, 
 form a decreasing chain, i.e.,  
$$
G \geq G(\chi_1) \geq G(\chi_1,\chi_2) \geq \dots \geq G(\chi_1,\dots,\chi_n).
$$
So we may choose $\ma$ and $\mb$ satisfying
\begin{subequations}\mylabel{sy.e2}
\begin{gather} \mylabel{Ca} \ma \in \Hall_{\pi}(G), \mb \in \Hall_{\pi'}(G),  \\
\mylabel{Cb} \ma(\chi_1, \chi_2,\dots, \chi_h) \text{ and } 
\mb(\chi_1, \chi_2,\dots,
\chi_h) \text{ form a Hall system for } G(\chi_1, \chi_2,\dots, \chi_h),\\
\mylabel{Cc} 
\ma(\chi_1, \chi_2,\dots, \chi_n) = P^*_{2k'} \text{ and } \mb(\chi_1,
\chi_2,\dots, \chi_n) = Q^*_{2l'-1},
\end{gather}
\end{subequations}
for all $h=1,\dots, n$.
So \eqref{Cb} says that $\ma, \mb$ reduces into
 $G(\chi_1, \chi_2,\dots,\chi_h)$,
 for each $h = 1,2,\dots,n$, while \eqref{Cc} says that $\ma, \mb$ 
reduces to the Hall
system $P^*_{2k'}, Q^*_{2l'-1}$ for $G(\chi_1, \chi_2,\dots,\chi_n)$.

We fix an integer  $m =1,\dots, n$ and we consider 
the   normal series 
\begin{subequations}\mylabel{sy.e3}
\begin{equation}\mylabel{sy.e3a}
1=G_0 \unlhd G_1 \unlhd \dots \unlhd G_{m} \unlhd G.
\end{equation}
 The sub tower 
\begin{equation}\mylabel{sy.e3b}
\{\chi_i \in \Irr(G_i) \}_{i=0}^m 
\end{equation}
of \eqref{sy.A} is a character tower of the above series.
If $k$ and $l$ are defined  as in \eqref{kl:def} for  $m$, then 
 Remark \ref{t-tr2} implies that 
the subset 
\begin{equation}\mylabel{sy.e3c}
\{ P_{2r}, Q_{2i-1} | \alpha_{2r}, \beta_{2i-1} \}_{r=0, i=1}^{k, l}
\end{equation}
\end{subequations}
of \eqref{sy.e1c} is a triangular set corresponding to  the  chain 
\eqref{sy.e3a} and tower \eqref{sy.e3b}. 

As in \eqref{pqinf}, we set 
    $G^* := G(\chi_1,\chi_2,\dots,\chi_m)$. 
So we can define the intersection groups 
\begin{align}\mylabel{sy.e4}
\ma^* &:= \ma \cap G^*= \ma(\chi_1, \dots, \chi_m), \notag \\
\mb^* &:=\mb \cap G^*= \mb(\chi_1, \dots, \chi_m).
\end{align}
Note that these definitions depend heavily on $m$.
Then we can prove
\begin{proposition}\mylabel{sy.p1}
Let $m = 1, \dots, n$ be fixed,  and $k, l$ be its  associate, 
via \eqref{kl:def}, integers. Then 
\begin{subequations}\mylabel{sy.e5}
\begin{align}
P_{2k}^* \unlhd \ma^* \leq \ma &\text{ and }  Q_{2l-1}^* \unlhd \mb^* \leq \mb, \\
\mylabel{sy.G}
N(P^*_{2k} \tin \mb^* ) &\in
 \Hall_{\pi'}(G(\alpha_{2k}^*, \chi_1, \dots, \chi_m)),  \\
N(Q_{2l-1}^* \tin \ma^*) &\in \Hall_{\pi}(G(\beta_{2l-1}^*, \chi_1, \dots, \chi_m)).
\end{align}
\end{subequations}
\end{proposition}

\begin{proof}
We fix the integers $m, k$ and $l$, and the triangular set 
 \eqref{sy.e3c} corresponding to the tower 
\eqref{sy.e3b}.
Then  Corollary \ref{q*:c1} implies that 
  $P_{2k}^*$ and $Q_{2l-1}^*$ form a Hall system for 
$G^*_m$.
According to     \eqref{Cb}  the groups    
$\ma^*= \ma \cap G^*$ and $\mb^*=\mb \cap G^*$ form a Hall system for
$G^*$.  In view of \eqref{Cc} the group $\ma^*$ 
contains $P^*_{2k'}$, and hence
contains $P^*_{2k} \le P^*_{2k'}$. 
Similarly, $\mb^*$ contains $Q^*_{2l-1}$.
Since $P^*_{2k}$ and $Q^*_{2l-1}$ form a Hall system for $G^*_m$, it
follows that 
\begin{subequations}\mylabel{sy.D}
\begin{align}
 P^*_{2k} = \ma^* \cap G^*_m &\unlhd  \ma^*, \\
Q^*_{2l-1} = \mb^* \cap G^*_m &\unlhd  \mb^*.
\end{align}
\end{subequations}

     The subgroup $G^*_m = G^* \cap G_m$ is normal in $G^*$. Hence
conjugation by elements of $\mb^*$ permutes among themselves the Hall
$\pi$-subgroups of $G^*_m$.  One of those Hall $\pi$-subgroups is $P^*_{2k}$.
Since $G^*_m = P^*_{2k}Q^*_{2l-1}$, the normal subgroup $Q^*_{2l-1}$ of
$\mb^*$ acts transitively on those Hall $\pi$-subgroups. It follows that
\begin{equation}\mylabel{sy.E}
  \mb^* = N(P^*_{2k} \tin \mb^*) \cdot Q^*_{2l-1}. 
\end{equation}

This implies that $N(P^*_{2k} \tin \mb^*)\cdot G^*_m/G^*_m$ is a Hall
$\pi'$-subgroup of $G^*/G^*_m$.  Since $N(P^*_{2k} \tin \mb^*) \cap G^*_m =
N(P^*_{2k} \tin Q^*_{2l-1})$ is a Hall $\pi'$-subgroup of $N(P^*_{2k} \tin
G^*_m)$, we conclude that
\begin{equation}\mylabel{sy.F}
N(P^*_{2k} \tin \mb^*) \in \Hall_{\pi'}(N(P^*_{2k} \tin G^*)).
\end{equation}

Now the group $N(P^*_{2k} \tin \mb^*)$
normalizes $P^*_{2k}$, and thus normalizes $P^*_{2i} = P^*_{2k} \cap
G_{2i}$ for all $i = 1,2,\dots,k$. As $N(P_{2k}^* \tin \mb^*)$ 
is a subgroup of 
 $ \mb^*$,  it normalizes
$Q^*_{2l-1} = \mb^* \cap G_m$. So it normalizes $Q^*_{2j-1} = Q^*_{2l-1}
\cap G_{2j-1}$ for each $j = 1,2,\dots,l$. Since it normalizes both
$P^*_{2i}$ and $Q^*_{2i-1}$, it normalizes $P_{2i} = N(Q^*_{2i-1} \tin
P^*_{2i})$ (see \eqref{q*:e7}),  for each $i = 1,2,\dots,k$.
Similarly, it normalizes $Q_{2j-1}
= N(P^*_{2j-2} \tin Q^*_{2j-1})$ (see \eqref{q*:e8}), 
 for each $j = 2,3,\dots,l$. It also
normalizes $Q_1 = G_1$.  The definitions of $Q_{2i-1,2j}$ and $P_{2r,
2s-1}$ in \eqref{pq1413defa} and \eqref{pq1413defb} show that they, too,
are normalized by $N(P^*_{2k} \tin \mb^*)$. Thus $N(P^*_{2k} \tin \mb^*)$
normalizes every subgroup appearing in the triangles displayed as
\eqref{pq12a} and \eqref{pq11a} in Chapter 5.

     The group $N(P^*_{2k} \tin \mb^*)$ also fixes all the characters
$\chi_1, \chi_2, \dots, \chi_m$, since $\mb^*$ does. Because it also
normalizes $Q_1$, $P_2$, $Q_3$, ..., it leaves invariant the $cQ_1$- ,
$cP_2$- , $cQ_3$-, ... correspondences in Table \ref{diagr.1-2-3...}.
Hence it fixes all the characters in that table. In particular, it fixes
$\alpha_{2i}$, for $i = 1,2,\dots,k$ and $\beta_{2j-1}$, for $j =
1,2,\dots,l$. It also fixes all the characters $\alpha_{2i, 2j-1}$ and
$\beta_{2r-1,2s}$ in the displayed triangles \eqref{pq12b}
and \eqref{pq11b}. Because it 
fixes all the groups and characters entering into the definition of
$\alpha^*_{2i}$, it also fixes that character for each $i = 1,2,\dots,k$.
Similarly, it fixes $\beta^*_{2j-1}$ for $j = 1,2,\dots,l$.

     At this point we know that $N(P^*_{2k} \tin \mb^*)$ is a Hall
$\pi'$-subgroup of $N(P^*_{2k} \tin G^*)$ fixing $\alpha^*_{2k}$. Hence it is
a Hall $\pi'$-subgroup of $G^*(\alpha^*_{2k})$. Since $G^* = G(\chi_1,
\chi_2, \dots, \chi_m)$, 
 we get that \eqref{sy.G} follows immediately.

The proof of (\ref{sy.e5}c) is similar, with the roles of 
$\pi$ and $\pi'$ interchanged. So we omit it.
\end{proof}

The proof of Proposition \ref{sy.p1} implies
\begin{corollary}\mylabel{sy.co1}
Both $N(P_{2k}^* \tin \mb^*) $ and $N(Q_{2l-1}^* \tin \ma^*)$ fix  the characters 
 $\alpha_{2i}$, for $i=1, \dots, k$, 
and $\beta_{2j-1}$, for $j=1,\dots, l$. They also fix 
 $\alpha_{2l-2,2l-1}$ and $\beta_{2k-1,2k}$.
\end{corollary}
With the above notation, 
we can now  prove
\begin{theorem}\mylabel{sy.H}
Assume the series \eqref{sy.A}, the tower \eqref{sy.B} and 
the triangular set \eqref{sy.e1c} are fixed. 
Assume further,  that $m$ is any  integer with $1 \leq m \leq n $, 
 and consider the series, tower and triangular set   
appearing in \eqref{sy.e3} for that $m$.
 Then we can choose a $\pi'$-Hall subgroup $\qw$ 
of $G'= G(\alpha_{2k}^*)$, 
to satisfy the conditions \eqref{pq30} and \eqref{pq31b}
 in   Theorem \ref{hat:p1} for the set \eqref{sy.e3c} and the tower 
\eqref{sy.e3b}, along with the property
$$N(P_{2k}^* \tin \mb(\chi_1, \dots, \chi_{2k}) ) = \qw(\beta_{2k-1,2k}).$$ 
Hence  
\begin{equation}\mylabel{sy.H1}
\qw(\beta_{2k-1,2k}) \cdot Q_{2l-1}^* \leq \mb(\chi_1, \dots, \chi_{2k}) 
\leq \mb.
\end{equation}
\end{theorem}

\begin{proof}
     Suppose first that $m = 2l-1$ is odd with $l \leq l'$, so $2k = 2l-2$.
 Then \eqref{pq30b}  tells us that the
groups $G'(\beta_{2k-1,2k})$, $G'(\chi_1, \dots, \chi_{2k})$,
$G'(\chi_1,\dots, \chi_{2k-1})$ and $G'(\beta_1, \beta_3, \dots,
\beta_{2k-1})$ have a common Hall $\pi'$-subgroup, (where
 $G'= G(\alpha_{2k}^*)$). 
Proposition \ref{sy.p1}, and in particular \eqref{sy.G},
 with $2k = 2l-2$ in the place of $m$ there, implies that 
 the $\pi'$-group
$N(P^*_{2k} \tin \mb^*)=N(P^*_{2k} \tin \mb(\chi_1, \dots, \chi_{2k}))$
 is a Hall $\pi'$-subgroup of the second group on this
list. By Corollary \ref{sy.co1} the character $\beta_{2k-1,2k}$ is fixed by 
$N(P_{2k}^* \tin \mb(\chi_1, \dots, \chi_{2k}))$. Hence the latter is a subgroup of 
$G'(\beta_{2k-1,2k})$.
So it must be a Hall $\pi'$-subgroup of that group,
because of its order. Similarly, it is contained in both $G'(\chi_1,
\dots,\chi_{2k-1})$ and $G'(\beta_1, \dots, \beta_{2k-1})$. Hence it is a
Hall $\pi'$-subgroup of those groups, too. So it satisfies all the
conditions for $\qw(\beta_{2k-1,2k})$ in \eqref{pq30b}. 
Clearly it also  satisfies  the
equations \eqref{pq31a},   as these follow from \eqref{pq30b}. 
 Furthermore,
$N(P^*_{2k} \tin \mb(\chi_1,\dots,\chi_{2k}))$ fixes $\alpha_2, \alpha_4, \dots, \alpha_{2k}$, 
by Corollary \ref{sy.co1}, and thus satisfies 
equation  \eqref{hat:p1e}.

According to (\ref{sy.e5}a) (for $m=2l-1$), the group $Q_{2l-1}$ is 
a subgroup of $\mb(\chi_1,\dots, \chi_{2k}, \chi_{2l-1}) \leq 
\mb(\chi_1, \dots, \chi_{2k})$. Furthermore, $Q_{2l-1}$ normalizes 
$P_{2k}^* = P_{2l-2}^*$, by \eqref{pq14'a}. 
Thus $Q_{2l-1}$ is a subgroup of $N(P_{2k}^*  \tin \mb(\chi_1, \dots, \chi_{2k}))$.
By Corollary \ref{sy.co1}, the latter normalizes 
$G(\alpha_2, \dots, \alpha_{2k}, \beta_1, \dots, \beta_{2k-1})$. 
Hence it normalizes 
$G_{2l-1}(\alpha_2, \dots, \alpha_{2k}, \beta_1, \dots, \beta_{2k-1})= 
Q_{2l-1} \rtimes P_{2k}$, where the equality follows from \eqref{pq16a} 
as $2k= 2l-2$. Hence $N(P_{2k}^* \tin \mb(\chi_1,\dots, \chi_{2k}))$
 normalizes 
$Q_{2l-1} \rtimes P_{2k}$ and contains $Q_{2l-1}$. Therefore, 
it normalizes $Q_{2l-1}$. Thus  it satisfies   \eqref{pq31b}.
Evidently we
can choose $\qw \in \Hall_q(G')$  so that 
$N(P_{2k}^* \tin \mb(\chi_1, \dots, \chi_{2k}) ) = \qw(\beta_{2k-1,2k})$.

So $\qw(\beta_{2k-1,2k}) \leq \mb(\chi_1,\dots,\chi_{2k})$. 
But  $Q_{2l-1}^* $ is contained in $\mb(\chi_1,\dots,\chi_{2l-1})$,
 by (\ref{sy.e5}a)and \eqref{sy.e4}, as $m=2l-1$.
As $2k=2l-2 < 2l-1$ in the odd case, 
$\mb(\chi_1,\dots, \chi_{2l-1}) \leq \mb(\chi_1,\dots,\chi_{2k})$.  Hence  
Theorem \ref{sy.H} follows for any odd $m$.

If $m=2k$ is even and strictly smaller than $n$, then we can still form the 
$2k+1$ series, by adding the group $G_{2k+1}$ and its character $\chi_{2k+1}$.
Then,  according to Corollary \ref{b.cB},  the even system, (where $m=2k$), 
with the odd, (where $m=2k+1$), share the group $\qw$.
This, along  the already proved odd case of Theorem \ref{sy.H},
 implies the first   part of Theorem \ref{sy.H} when $m =2k < n$.

If $m=2k=n$, we can't form a bigger odd system, but we know exactly what
group $\qw(\beta_{2k-1,2k})$ is.
Indeed, as $G_{2k} /G_{2k-1}$ is a $\pi$-group, and 
$\qw(\beta_{2k-1,2k})$ is  a $\pi'$-Hall subgroup of 
$G'(\beta_{2k-1,2k})$, by \eqref{pq30b}, it must be a $\pi'$-Hall subgroup of 
$G_{2k-1}'(\beta_{2k-1,2k})$. This, along with \eqref{pq40c}, implies
$$
\qw(\beta_{2k-1,2k})= \qw(\beta_{2k-1,2k}) \cap G_{2k-1}= 
\qw_{2k-1}(\beta_{2k-1,2k}).
$$
By \eqref{pq40g} this gives  
$$
\qw(\beta_{2k-1,2k}) = Q_{2k-1,2k}.
$$
On the other hand, in the case $m=2k=n$ we have $2l'-1= 2k-1$. Thus
(\ref{sy.e2}c) implies $\mb (\chi_1,\dots,\chi_n)=Q_{2k-1}^*$. Hence 
$$
N(P_{2k}^* \tin \mb(\chi_1,\dots,\chi_{2k}))= N(P_{2k}^* \tin Q_{2k-1}^*)
=Q_{2k-1,2k}.
$$
So  the first part of Theorem \ref{sy.H}  holds in the  case $m=2k=n$. 

Furthermore, when $m=2k$ is even, (\ref{sy.e5}a)  and \eqref{sy.e4} imply that 
$Q_{2l-1}^* \leq \mb(\chi_1,\dots,\chi_{2k})$.
Thus \eqref{sy.H1} follows for the even case. 
This completes the proof of Theorem \ref{sy.H}
\end{proof}

     Of course, a similar result holds by $p,q$-symmetry for 
$\hap(\alpha_{2l-2, 2l-1}) \cdot P^*_{2l-2}$, whenever $l =
2,3,\dots,l'$.
\begin{theorem}\mylabel{sy.t1}
Assume the series \eqref{sy.A} the tower \eqref{sy.B} and 
the triangular set \eqref{sy.e1c} are fixed. 
Assume further that $m$ is any integer with $1 \leq m \leq n $ 
 and consider the series of subgroups, the character tower, and the 
triangular set  
appearing in \eqref{sy.e3} for that $m$.
 Then we can choose a $\pi$-Hall subgroup $\hap$ 
of $G(\beta_{2l-1}^*)$, to satisfy  \eqref{sy.e}
 for the set \eqref{sy.e3c} and the tower 
\eqref{sy.e3b}, along with the property
$$N(Q_{2l-1}^* \tin \ma(\chi_1,\dots, \chi_{2l-1}) ) = \hap(\alpha_{2l-2,2l-1}).$$ 
Hence  
\begin{equation}\mylabel{sy.H2}
\hap(\alpha_{2l-2,2l-1}) \cdot P_{2k}^* 
\leq  \ma(\chi_1,\dots,\chi_{2l-1}) \leq \ma.
\end{equation}
\end{theorem}

\section {``Shifting '' properties}
\mylabel{sh}
We assume that  the normal series 
\begin{equation}\mylabel{s1}
G_0=1 \unlhd G_1 \unlhd \dots \unlhd G_n = G
\end{equation}
 is fixed  for some $n \geq 2$,  and satisfies 
\begin{subequations}\mylabel{sh.e1}
\begin{gather}\mylabel{sh.e1a}
 \text{$G_{2}= G_{2, \pi} \times G_{2,\pi'}$
 is a $\pi$-split group (see Definition \ref{sp.d1}) \,  and  }\\
\text{ $G$ fixes the character $ \chi_1$. }
\end{gather}
Assume also that the character tower 
\begin{equation}\mylabel{sh.e1b}
\{1=\chi_0, \chi_1, \dots, \chi_n \}
\end{equation}
is fixed, while the set 
\begin{equation}\mylabel{sh.e1c}
\{Q_{2i-1}, P_{2r} |\beta_{2i-1}, \alpha_{2r} \}_{i=1, r=0}^{l', k'}
\end{equation}
 is a representative of the conjugacy class
 of triangular sets that corresponds to \eqref{sh.e1b}.
 Furthermore, we fix a Hall system 
$\{ \ma, \mb\}$ of $G$ that satisfies \eqref{sy.e2}, that is 
\begin{gather}\mylabel{sh.AB}
 \ma \in \Hall_{\pi}(G), \mb \in \Hall_{\pi'}(G),  \\
 \ma(\chi_1, \chi_2,\dots, \chi_h) \text{ and } 
\mb(\chi_1, \chi_2,\dots,
\chi_h) \text{ form a Hall system for } G(\chi_1, \chi_2,\dots, \chi_h),\\
\ma(\chi_1, \chi_2,\dots, \chi_n) = P^*_{2k'} \text{ and } \mb(\chi_1,
\chi_2,\dots, \chi_n) = Q^*_{2l'-1},
\end{gather}
\end{subequations}
for any $h= 1, 2, \dots, n$.
 According to 
Corollary \ref{sp.co1} this set satisfies
\begin{subequations}\mylabel{sh.e2}
\begin{align}
G_{2} &= P_{2} \times Q_{1},\\
\chi_{2}&= \chi_{2,\pi}\times \chi_{2,\pi'}= 
\alpha_{2} \times \beta_{1},
\end{align}
where $P_2$ and $Q_1 = G_1$ are  the $\pi$-and $\pi'$-Hall subgroups 
respectively,  
of $G_2$. 
\end{subequations}

We replace the first $\pi'$-group $G_1=Q_1$ appearing in \eqref{s1},
 by the trivial group  and consider the series  
\begin{subequations}\mylabel{sh.e3}
\begin{equation}\mylabel{sh.e3a}
1 \unlhd G_1^s:= 1 \unlhd  G_2^s:= P_2  \unlhd G_3^s:= G_{3} 
\unlhd G_4^s:=  G_{4} \unlhd  \dots \unlhd G_n^s:= G_{n}= G.
\end{equation}
We call the series \eqref{sh.e3a} a {\rm shifting  } of the series 
 \eqref{s1}.   Note that  \eqref{sh.e3a} is a  normal series 
of $G$, that satisfies
 Hypothesis \ref{hyp1} with $G_1^s = 1$.    
The characters
\begin{equation}\mylabel{sh.e3b}
 1, \chi_1^s:= 1,   \chi_2^s:= \alpha_{2},  \chi_3^s:= \chi_3, 
\chi_4^s:= \chi_4, \dots, \chi_n^s:= \chi_n,
\end{equation}
form a character tower for the series \eqref{sh.e3a}. In addition,   
the set
\begin{equation}\mylabel{sh.e3c} 
\{ Q_1^s=1= P_0^s, Q_{2i-1}^s= Q_{2i-1}, P_{2r}^s= P_{2r}|
\beta_1^s=1= \alpha_0^s,  \beta_{2i-1}^s= \beta_{2i-1},
 \alpha_{2r}^s= \alpha_{2r}\}_{i=2, r=1}^{l', k'}
\end{equation}
\end{subequations}
 is a triangular set for \eqref{sh.e3a}, 
 corresponding  to the character tower \eqref{sh.e3b} (this can be  
very easily verified using the fact that \eqref{sh.e1c} is a triangular set 
for \eqref{s1} corresponding to \eqref{sh.e1b}).
Note that the groups 
$Q_{2i-1,2j}^s, P_{2r, 2t+1}^s$ and their characters 
$\beta_{2i-1,2j}^s$ and $\alpha_{2r, 2t+1}^s$, respectively, 
 remain the same as those of \eqref{sh.e1},
 whenever $2 \leq i \leq j \leq  k'$ and $2\leq r \leq
 t \leq l'-1$, i.e., 
\begin{gather}
Q_{2i-1,2j}^s = Q_{2i-1,2j} \text{ and }  P_{2r, 2t+1}^s =P_{2r, 2t+1}, \\
\beta_{2i-1,2j}^s= \beta_{2i-1,2j} \text{ and }  
\alpha_{2r, 2t+1}^s = \alpha_{2r, 2t+1}.
\end{gather}
Also, the product groups $P_{2k}^{*, s} = P_2^s \cdots P_{2k}^s$,
 remain unchanged, for 
every $k=1,\dots,k'$, because $P_{2r}^s= P_{2r}$ whenever 
$1\leq r \leq k'$.  In addition, for any such $k$, 
 the  irreducible character
$\alpha_{2k}^* \in \Irr(P_{2k}^*)$  was chosen as 
the $Q_3 ,\dots, Q_{2k-1}$-correspondent of $\alpha_{2k}$ (see
 Definition \ref{p*d1}). As neither of the above groups
nor the character  $\alpha_{2k}$
 changes when 
passing to the shifted system \eqref{sh.e3}, we conclude that also the
character $\alpha_{2k}^{*, s} \in \Irr(P_{2k}^{*, s})$ remains unchanged,
 that is, 
\begin{equation}\mylabel{s2}
\alpha_{2k}^{*, s} = \alpha_{2k}^*,
\end{equation}
for all $k=1, \dots, k'$.

Furthermore, the group $Q_{3}$ contains  $Q_1$, as
 $Q_3 \geq Q_{1,2}= C(P_2 \tin Q_1) = Q_1$, by (\ref{sh.e2}b). Thus 
\begin{equation}\mylabel{sh.e9}
Q_{2l-1}^*= Q_3 \cdot Q_5 \cdots Q_{2l-1}, 
\end{equation}
whenever $1 < l \leq l'$.
This implies that 
\begin{equation}\mylabel{s3}
Q_{2l-1}^*= Q_3^s \cdot Q_5^s \cdots Q_{2l-1}^s= Q_{2l-1}^{*, s},
\end{equation}
for all such $l$.
Furthermore, the fact that $P_2$ centralizes both  $Q_1$ and  
 $N(P_2 \tin Q_{2l-1}^* )= Q_{2l-1}^*$  implies that the 
  $P_2, P_4, \dots, P_{2l-2}$-correspondent 
 $\beta_{2l-1}^* \in \Irr(Q_{2l-1}^*)$ 
of $\beta_{2l-1} \in \Irr(Q_{2l-1})$, 
 is actually the $P_4, \dots, P_{2l-2}$-correspondent
of $\beta_{2l-1} $, whenever $1< l \leq l'$.
Clearly,  for all such $l$, we get  
 $\beta_{2l-1}^{*,s} = \beta_{2l-1}^*$, because  $Q_{2l-1}^{*, s}
= Q_{2l-1}^*, P_2^s = P_2, \dots, P_{2l-2}^s= P_{2l-2}$ and 
$\beta_{2l-1}^s= \beta_{2l-1}$.

As $G$ fixes $\chi_1$ by (\ref{sh.e1}b), 
 while $\chi_2= \alpha_2 \times \beta_1$  by 
(\ref{sh.e2}d), we get that $G= G(\chi_1 ) = G(\beta_1)$.
We conclude that 
\begin{subequations}\mylabel{sh.e11}
\begin{equation}\mylabel{sh.e12}
G(\chi_1^s, \chi_2^s) = G(\alpha_2) = G(\chi_1, \chi_2),
\end{equation}
and thus
\begin{equation}\mylabel{sh.e10}
G(\chi_1^s, \chi_2^s, \chi_3^s, \dots, \chi_h^s) = 
G(\alpha_2, \chi_3, \dots, \chi_h) 
= G(\chi_1, \chi_2, \chi_3, \chi_4, \dots, \chi_h),
\end{equation}
for all $h=3, \dots, n$.
\end{subequations}
For any subgroup $H$ of $G$, similar equations, with $H$ in place 
of $G$, hold.
In particular, for the  Hall system $\{ \ma, \mb\}$ we get 
\begin{align*}
\ma(\chi_1^s, \chi_2^s)=\ma(\alpha_2) = \ma(\chi_1, \chi_2) &\text{ and } 
\mb(\chi_1^s, \chi_2^s) = \mb(\alpha_2) = \mb(\chi_1, \chi_2), \\
\ma(\chi_1^s, \chi_2^s, \dots, \chi_h^s) = 
&\ma(\alpha_2, \chi_3, \dots, \chi_h) = \ma(\chi_1,\chi_2,  \dots, \chi_h) 
\text{ and } \\
\mb(\chi_1^s, \chi_2^s, \dots, \chi_h^s) = 
&\mb(\alpha_2, \chi_3, \dots, \chi_h) = \mb(\chi_1,\chi_2,  \dots, \chi_h), 
\end{align*}
for all $h=3, \dots, n$.
This, along with teh conditions  (\ref{sh.e1}e,f,g) which 
 $\ma$ and $\mb$ satisfy,  implies
\begin{gather*}
\ma(\alpha_2)= \ma(\chi_1^s,\chi_2^s)  \text{ and }
 \mb(\alpha_2) = \mb(\chi_1^s,\chi_2^s)
\text{ form a Hall system for }  G(\alpha_2)= G(\chi_1^s, \chi_2^s) \\
\ma(\chi_1^s, \chi_2^s,  \dots, \chi_h^s) \text{ and } 
\mb(\chi_1^s, \chi_2^s, \dots, \chi_h^s) \text{ form a Hall system for } 
G(\chi_1^s, \chi_2^s, \dots, \chi_h^s), 
\end{gather*}
for any $h= 3, \dots, n$.
Furthermore, (\ref{sh.e1}g), along with  \eqref{sh.e9},  implies 
\begin{equation}
\ma(\chi_1^s, \chi_2^s, \dots, \chi_n^s)=  = P_{2k'}^*  \text{ and } 
\mb(\chi_1^s, \chi_2^s,\dots, \chi_{n}^s) = Q_3 \cdots Q_{2l'-1}.
\end{equation}
Therefore, the groups $\ma, \mb$ satisfy the equivalent of (\ref{sh.e1}e,f,g), 
for the shifted system \eqref{sh.e3}.

The other groups of interest that doesn't change, when we work in the 
shifted case \eqref{sh.e3},
  are  $\qw$ and $\hap$, as these are defined for 
every fixed, but arbitrary,  smaller system 
\begin{subequations}\mylabel{sh.sm}
\begin{gather}\mylabel{sh.sm1}
1= G_0  \unlhd G_1 \unlhd G_2 \unlhd \dots  \unlhd  G_m \unlhd G,\notag \\
1, \chi_1, \chi_2, \dots, \chi_m,\\
\{  Q_{2i-1}, P_{2r} | 
\beta_{2i-1}, \alpha_{2r} \}_{i=1, r=0}^{l,k} \notag
\end{gather}
where $m=3,\dots,n$, and $k, l$ are related to $m$ via \eqref{kl:def}.
Of course when we shift the original system  to get 
\eqref{sh.e3}, we also get 
the smaller shifted system  
\begin{gather}\mylabel{sh.sm2}
1 \unlhd G_1^s=1 \unlhd G_2^s=P_2 \unlhd G_3^s=G_3 \unlhd 
\dots  \unlhd  G_m^s=G_m \unlhd G, \notag \\
1, \chi_1^s, \chi_2^s, \dots, \chi_m^s,\\
\{  Q_{2i-1}^s, P_{2r}^s | \beta_{2i-1}^s, 
\alpha_{2r}^s \}_{i=1, r=0}^{l,k} \notag
\end{gather}
\end{subequations}
Indeed, it is easy   to see that the same group $\qw$, which 
 was picked among the  $\pi'$-Hall subgroups of $G(\alpha_{2k}^*)$
to  satisfy  the conditions \eqref{pq30} and \eqref{pq31b} in 
Theorem \ref{hat:p1} for the system \eqref{sh.sm1}, 
 satisfies the same conditions  for the shifted system  \eqref{sh.sm2}.
 First  note that $G'= G(\alpha_{2k}^*) = G(\alpha_{2k}^{*,s})$. 
Hence $\qw$ is also a $\pi'$-Hall subgroup of $G(\alpha_{2k}^{*,s})$.
Furthermore, 
the group $\qw$ fixes $\beta_1$, as the latter is $G$-invariant. 
This forces $\qw$ to fix 
the $P_{2k}^*$-Glauberman correspondent
 $\beta_{1,2k}$ of $\beta_1$, as $\qw$ normalizes $P_{2k}^*$.
Hence $\qw(\beta_{1,2k}) =\qw $.
For the shifted system \eqref{sh.sm2}  we have $\beta_1^s=1$. 
Thus  the $P_{2k}^{*, s}=P_{2k}^*$-Glauberman correspondent
  $\beta_{1,2k}^s$ of $\beta_1^s$
 is also trivial.
Hence $\qw(\beta_{1,2k}^s)= \qw = \qw(\beta_{1,2k})$.
In addition, $\beta_{2i-1, 2k}^s= \beta_{2i-1,2k}$, 
 for any $i=2, \dots, l$
We conclude that 
\begin{align*}
 \qw(\beta_{2i-1,2k}^s) \in \Hall_{\pi'}(G'(\beta_{2i-1,2k}^s))
&\cap
\Hall_{\pi'}(G'(\chi_1^s,  \chi_2^s, \chi_3^s, \dots,\chi_{2i-1}^s))
\cap \notag\\
&\Hall_{\pi'}(G'( \chi_1^s, \chi_2^s,\chi_{3}^s,\dots,\chi_{2i}^s)) \cap
\Hall_{\pi'}(G'(\beta_1^s, \beta_3^s,\dots,\beta_{2i-1}^s)),   \\
 \qw(\beta_{2i-1,2k}^s) =\qw(\chi_1^s, \chi_2^s, 
\chi_{3}^s,\dots,\chi_{2i-1}^s) &= \qw(\chi_1^s, \chi_2^s, 
\chi_{3}^s,\dots,\chi_{2i}^s)
=\qw(\beta_1^s, \beta_{3}^s,\beta_5^s, \dots,\beta_{2i-1}^s) \text{ and }\\
\qw(\chi_1^s, \chi_2^s, \chi_{3}^s,\dots,\chi_{2i-1}^s) &\leq
 \qw(1, \alpha_{2}^s, \dots,
\alpha_{2i}^s) 
\end{align*}
for all $i=1,2,\dots,k$. 
In addition,  for all $i$ with $1\leq i \leq l-1$ we get  
\begin{equation*}
 \qw(\beta_{2i-1,2k}^s) \text{ normalizes } \, Q_{2i+1}^s= Q_{2i+1}. 
\end{equation*}
 So  $\qw$  remains unchanged in the shifted case, as does
 $\qw(\beta_{2k-1,2k})$. Therefore
the image $I$ of $ \qw(\beta_{2k-1})$ in 
$\Aut(P_{2k}^*)$ remains unchanged.

Similarly, we can show that the group $\hap$ remains unchanged in the shifted 
case, as does  $\hap(\alpha_{2l-2,2l-1})$. So  
the image of the latter group in $\Aut(Q_{2l-1}^*)$ remains unchanged.

It is also clear that if the characters 
$\beta_{2k-1,2k} \in \Irr(Q_{2k-1,2k})$
and $\alpha_{2l-2,2l-1} \in \Irr(P_{2l-2,2l-1})$ 
extend to $\qw(\beta_{2k-1,2k})$ and $\hap(\alpha_{2l-2,2l-1})$, respectively,
 then the same property passes to the 
shifted case,  provided that $k \geq 2$, 
 as none of these groups and characters really changes. 
In conclusion we have
\begin{theorem}\mylabel{sh.t1}
Assume that the normal series $1=G_0 \unlhd G_1 \unlhd \dots \unlhd G_n 
\unlhd G$ satisfies \eqref{sh.e1}. Let \eqref{sh.e1b} be a  character
 tower for the series, \eqref{sh.e1c} the corresponding triangular set
and $\{ \ma, \mb \}$ a Hall system for $G$ that satisfies 
\eqref{sh.AB}.
Replacing  the group $G_1=Q_1$ with the trivial group   we
obtain a normal series    \eqref{sh.e3a} for $G$. Then \eqref{sh.e3b}
is a character tower for that series, and  
  \eqref{sh.e3c} its corresponding  triangular set. Furthermore, 
$\{\ma, \mb \} $ remains a Hall system for $G$ satisfying  
the equivalent of (\ref{sh.e1}e,f,g) for the series \eqref{sh.e3a}
and tower \eqref{sh.e3b}.  
 Even more,  the groups $P_{2k}^* , Q_{2l-1}^*, \qw, \hap$,
 as well as $\qw(\beta_{2k-1,2k})$ and 
$\hap(\alpha_{2l-2,2l-1})$ satisfy the same conditions for  
 the smaller system 
\eqref{sh.sm1} and the shifted one \eqref{sh.sm2}, whenever $m=3, \dots, n$.
\end{theorem}

The above theorem  makes clear that, whenever the series \eqref{s1} satisfies 
(\ref{sh.e1}a,b), we can replace the group $G_1$ with a trivial group 
without affecting any other group  or character 
involved in our constructions. 
From now on, for simplicity,   whenever   
such a shifting is performed,  we will be writing the
produced series, tower and triangular set of \eqref{sh.e3} as 
\begin{subequations}
\begin{equation}
1 \unlhd P_2 \unlhd G_3 \unlhd \dots \unlhd G_n = G
\end{equation}
\begin{equation}
\{1, \alpha_2, \chi_3, \dots, \chi_n \}
\end{equation}
\begin{equation}\mylabel{sh.e3cc} 
\{  Q_{2i-1}, P_{2r}|  \beta_{2i-1},
 \alpha_{2r}\}_{i=2, r=1}^{l', k'}
\end{equation}
\end{subequations}
Note that  the trivial groups  
$G_1^s= Q_1^s=1= P_0^s$ and their  characters, have been dropped.

%%% Local Variables: 
%%% mode: latex
%%% TeX-master: "sy"
%%% End: 

\chapter{ Normal Subgroups}\mylabel{elm}

As in Chapter \ref{sy}, we fix a normal chain
\begin{subequations}\mylabel{e0.1}
\begin{equation}\mylabel{e0.2}
   1 = G_0 \trianglelefteq G_1 \trianglelefteq 
G_2 \trianglelefteq \dots \trianglelefteq G_n = G, 
\end{equation}
for an odd order group $G$, 
such that  Hypothesis \ref{hyp1} holds with $n$ in the place of $m$, i.e., $n >0$ and  
 $G_i/G_{i-1}$ is a 
$\pi$-group if $i$ is even, and a $\pi'$-group if $i$ is odd,
 for each  $i = 1,2,\dots,n$.
 We also fix a character tower 
\begin{equation}\mylabel{e0.3}
\{ \chi_i  \in \Irr(G_i) \}_{i=0} ^n
\end{equation}
for the above series and a  corresponding triangular set 
\begin{equation}\mylabel{e0.4}
\{ P_{2r}, Q_{2i-1} | \alpha_{2r}, \beta_{2i-1}\}_{r =0 ,
i=1}^{k',l'}
\end{equation}
\end{subequations}
where $ k' = [n/2] \text{ and } l' = [(n+1)/2]$.
Along with that we fix a Hall  system $\ma, \mb$ of $G$ 
that satisfies \eqref{sy.e2}, that is,
\begin{subequations}\mylabel{e.AB}
\begin{gather}  \ma \in \Hall_{\pi}(G), \mb \in \Hall_{\pi'}(G),  \\
 \ma(\chi_1, \chi_2,\dots, \chi_h) \text{ and } 
\mb(\chi_1, \chi_2,\dots,
\chi_h) \text{ form a Hall system for } G(\chi_1, \chi_2,\dots, \chi_h),\\
\ma(\chi_1, \chi_2,\dots, \chi_n) = P^*_{2k'} \text{ and } \mb(\chi_1,
\chi_2,\dots, \chi_n) = Q^*_{2l'-1},
\end{gather}
\end{subequations}
whenever $h=1,\dots,n$.

\noindent
\section{Normal $\pi'$-subgroups inside $Q_1$}\mylabel{elm1}
We fix an integer $m= 1, \dots, n$ and we consider the  normal series 
\begin{subequations}\mylabel{el.1}
\begin{equation}\mylabel{el.01}
1=G_0\unlhd G_1 \unlhd \dots \unlhd G_{m} \unlhd  G.
\end{equation}
The sub tower
\begin{equation}\mylabel{el.e1}
\{1=\chi_0, \chi_1, \dots, \chi_m \}
\end{equation}
of \eqref{e0.3} is  a character tower for  \eqref{el.e1}, 
and the subset 
\begin{equation}\mylabel{el.e2}
\{Q_{2i-1}, P_{2r} | \beta_{2i-1}, \alpha_{2r} \}_{i=1, r=0} ^{l, \, \,  k}
\end{equation}
\end{subequations}
of \eqref{e0.4} is a representative 
 of the conjugacy class of triangular sets 
that corresponds uniquely to \eqref{el.e1}. (As usual, the integers
$k $ and $l$ are related to $m$ via \eqref{kl:def}.) 
 Thus all  the groups,  the characters  and their properties that were defined and proved 
 in Chapter \ref{pq} are valid for  \eqref{el.1}. In particular we can define
groups $G_{i,s}$  and their  characters $\chi_{i,s}$ (see Theorem \ref{tow--tri} 
for their properties). Even more,  we can pick the groups 
 $\qw$ and $\hap$ to satisfy the conditions in  Theorems \ref{sy.H} and \ref{sy.t1} respectively.
Hence if we write
\begin{subequations} \mylabel{e0.5}
\begin{align}\mylabel{e0.d}
\map &:= \hap(\alpha_{2l-2, 2l-1}) \cdot P_{2k}^*,\\ 
\maq &:= \qw(\beta_{2k-1,2k}) \cdot Q_{2l-1}^*,
\end{align}
then Theorems \ref{sy.H} and \ref{sy.t1}, and   in particular
\eqref{sy.H1} and \eqref{sy.H2},  imply 
\begin{gather}\mylabel{e0.i}
\map \leq \ma(\chi_1,\dots,\chi_{2l-1}) \leq \ma \\
 \maq \leq \mb(\chi_1,\dots,\chi_{2k}) \leq \mb.
\end{gather}
\end{subequations}
We remark  that the group  $\maq$ is well defined, as 
$\qw(\beta_{2k-1,2k}) \leq \qw(\beta_{2i-1,2k})$  normalizes 
the group  $Q_{2i+1}$, for all $i=1,\dots, l-1$, by \eqref{pq31b}.
It also normalizes $Q_1$. Thus it normalizes their product $Q_1 \cdot Q_3
\cdots Q_{2l-1}= Q_{2l-1}^*$. 
Similarly we can show that $\map$ is well defined.

For the rest of this section we assume that $S$ is a  subgroup
 of $G_1$, and $\zeta$ is a character of $S$, 
  satisfying 
\begin{subequations}\mylabel{els.1}
\begin{gather} 
S \unlhd G  \text{ and } S \leq G_1,\\
\zeta \in \Irr(S) \text{ is $G$-invariant  and lies under $\beta_{1}$. }
\end{gather}
Either $S$ or $\zeta$ may be trivial.
We also assume that $E$ is a normal subgroup of $G$ with 
\begin{equation}\mylabel{els.2}
S \leq E \leq  Q_1=G_1.
\end{equation}
\end{subequations}
Then 
\begin{lemma}\mylabel{elm:l1}
There is an irreducible character $\lambda \in \Irr(E)$ such that 
$\lambda$ is $\ma(\chi_1)$-invariant and lies under every $\chi_i$, for 
$i=1,\dots ,n$. Any such  $\lambda$ lies above $\zeta$.
\end{lemma}
\begin{proof}
Let $\lambda_1$ be an irreducible character of $E$ lying under $\chi_1$, and
thus under $\chi_i$ for any $i=1,\dots,n$.
Then Clifford's Theorem implies that $\chi_1$ lies above the $G_1$-conjugacy 
class of $\lambda_1$.
The $\pi$-group $\ma(\chi_1)$ fixes
 $\chi_1$, and normalizes $E$,  as the latter is normal in $G$. Hence
 $\ma(\chi_1)$ permutes among 
themselves the $G_1$-conjugates of $\lambda_1$. As  $(|\ma(\chi_1)|, |G_1|) = 1$,
Glauberman's Lemma (Lemma 13.8 in \cite{is}) implies that $\ma(\chi_1)$ fixes at least 
one,  $\lambda$,  of the $G_1$-conjugates of $\lambda_1$.

As $\zeta \in \Irr(S)$ is $G$-invariant and lies under $\chi_1$, 
Clifford's theorem implies that  any irreducible character of $E$ 
lying under $\chi_1$ also lies above $\zeta$.
 Thus the character $\lambda$  satisfies all  the conditions in the  lemma.
\end{proof}

Note that the  proof of   Lemma \ref{elm:l1} also shows 
\begin{remark}\mylabel{elm.r0}
Assume that  $\lambda_1 \in \Lin(E)$ is a linear character of $E$ lying 
under $\chi_1$. Then there exists 
a $G_1$-conjugate $\lambda \in \Lin(E)$ of $\lambda_1$, such that 
$\lambda$ is $A(\chi_1)$-invariant, and   lies under $\chi_1$ and above $\zeta$.
\end{remark}

\begin{remark}\mylabel{elm.r00}
The $\pi$-group $P_{2k'}^*= P_2 \cdots P_{2k'}$ fixes $\chi_1$, as every one of its 
factors $P_{2i}$  fixes $\chi_1 = \beta_1$, for all $i=1,\dots,k'$. 
Hence it is a subgroup of $ \ma(\chi_1)$. So $P_{2i}^*$
 fixes $\lambda$, for all $i=1,\dots, k'$.
\end{remark}

In the same spirit is 
\begin{remark}\mylabel{elm.r000}
The $\pi$-group $\map$ fixes $\lambda$, as it is a subgroup of 
$\ma(\chi_1)$, by  \eqref{e0.i}.
\end{remark}

From now on we fix a  character $\lambda \in \Irr(E)$ satisfying 
the conditions in Lemma \ref{elm:l1}.

The $\pi$-group  $P_{2i}^*$ acts  on the $\pi'$-group $E$,  and fixes  $\lambda$
by Remark \ref{elm.r00},  whenever $1\leq i \leq k'$. So we can define
\begin{subequations}\mylabel{el.e3.1}
\begin{align}
E_{2i} &= C(P_{2i}^* \tin E)  \text{ and } \\
\lambda_{2i} \in \Irr(E_{2i}) &\text{ is the $P_{2i}^*$-Glauberman 
correspondent of } \lambda,
\end{align}
\end{subequations}
for all $i=1,\dots,k'$.
We obviously have that 
\begin{equation}\mylabel{el.e3.2}
E_{2i} = Q_{1,2i} \cap E= C(P_{2i}^* \tin Q_1) \cap E,
\end{equation}
for each such $i$.
Furthermore,
\begin{equation}\mylabel{el.e3.3}
 \lambda_{2i} \text{ lies under  } \beta_{1,2i},
\end{equation}
as $\beta_{1,2i}$ is the $P_{2i}^*$-Glauberman correspondent of
 $\beta_1$, and $\beta_1= \chi_1$ lies 
over $\lambda$.  It follows from  the definition (\ref{el.e3.1}b) of $\lambda_{2i}$   that   
\begin{equation}\mylabel{el.e3.4}
N(\lambda)= N(\lambda_{2i}),
\end{equation}
for every group $N$ with  $N \leq N(P_{2i}^* \tin G)$.

We also  define 
\begin{subequations}\mylabel{el.e3}
\begin{align}
G_{\lambda} &:= G(\lambda), \mylabel{el.e3a}\\
\slg{i} &:= G_i(\lambda) = G_{\lambda} \cap G_i, \mylabel{el.e3b}  
\end{align}
\end{subequations}
whenever $0\leq i \leq n$. This way we can 
 form the series
 \begin{equation}\mylabel{elm:e1}
\slg{0} = G_0 = 1\unlhd \subl{G}{1} \unlhd \subl{G}{2} \unlhd \dots \unlhd \slg{n}
 = G_{\lambda}
\end{equation}
of normal subgroups of the stabilizer $G_{\lambda}$ of $\lambda$ in $G$.
Its is clear that this series satisfies Hypothesis \ref{hyp1} with $n$ in the place of $m$, as the series 
\eqref{e0.2}  does.
 Furthermore, 
 Clifford's Theorem provides  unique  irreducible 
characters $\subl{\chi}{i} \in \Irr(\subl{G}{i})$ 
 lying  above $\lambda$ and inducing 
 $\chi_i$, whenever $i= 0, 1, \dots n$, i.e.,
\begin{subequations}\mylabel{101}
\begin{equation}\mylabel{101a}
\slc{i} \in \Irr(\slg{i} | \lambda)  \text{ and } (\slc{i})^{G_i} = \chi_i.
\end{equation}
Clearly we get that 
\begin{equation}\mylabel{101b}
\slc{0} = \chi_0 = 1.
\end{equation}
\end{subequations}
 As the $\chi_i$ lie above each other, the same holds for the characters 
$\subl{\chi}{i}$, i.e.,  $\subl{\chi}{i}$  lies  above 
$\subl{\chi}{j}$ whenever $0\leq j \leq i \leq n$. 
This way we have formed a character tower 
\begin{subequations}
\begin{equation}\mylabel{el.e4a}
\{ 1=\slc{0}, \slc{1}, \dots , \slc{n} \}
\end{equation} 
for the series \eqref{elm:e1}. Hence Theorem \ref{cor:t},  applied  to the tower \eqref{el.e4a},  implies the 
existence of a unique $G_{\lambda}$-conjugacy class of triangular sets for 
\eqref{elm:e1} that correspond
to the tower \eqref{el.e4a}. Let 
\begin{equation}\mylabel{el.e4b}
\{\slq{2i-1}, \slp{2r} |\slb{2i-1}, \sla{2r} \}_{i=1, r=0}^{l', \, \, k'}
 \end{equation}
be a representative of this class. 
\end{subequations}
All the groups,  the characters and their properties  that were described 
in Chapter \ref{pq} are valid for the $\lambda$-situation. We follow the same notation as in Chapter 
\ref{pq},  with the addition of   an extra $\lambda$ in the   subscripts,   to  refer to this $\lambda$-situation.
As a small sample we give the following list:
\begin{gather*} 
\slq{2r-1,2r} =C(\slp{2r} \tin \slq{2r-1}),  \, \, \text{see \eqref{anno1},} \\ 
\slb{2r-1,2r} \in \Irr(\slq{2r-1,2r})
\text{ is the $\slp{2r}$-Glauberman correspondent of $\slb{2r-1} \in \Irr(\slq{2r-1})$ }, \\
\slp{2i,2i+1} = C(\slq{2i+1} \tin \slp{2i}),  \, \, \text{by \eqref{anno2},}  \\
\sla{2i, 2i+1}\in \Irr(\slp{2i, 2i+1})
 \text{ is the $\slq{2i+1}$-Glauberman correspondent of $\sla{2i} \in \Irr(\slp{2i})$ },\\
\slg{2i, 2i-1} = 
N( \slp{2}, \dots, \slp{2i-2}, \slq{1},\dots,\slq{2i-1} \tin \slg{2i}(\slc{1}, \dots, \slc{2i-1})), \, \, 
\text{see \eqref{tow--tri3},} \\
\slc{2i,2i-1} \text{ is the $c\slp{2}, \dots, c\slp{2i-2}, c\slq{1},\dots,
c\slq{2i-1}$-correspondent of } \slc{2i},  \, \,  \text{see } \\
\text{ Theorem \ref{tow--tri}  } \text{and  
 the Definition \ref{def.t-t} } 
\end{gather*}
for all $r=1,\dots,k'$ and $ i =1,\dots,l'-1$.

In this section we will describe the relations between the sets 
\eqref{e0.4} and \eqref{el.e4b}.  
We start with the groups $\slgs{i}$ defined as 
\begin{align}\mylabel{102c}
\slgs{0} &:= \slg{0} = 1, \notag \\
\slgs{i} &:= \slg{i}(\slc{1},\dots \slc{i-1}) = \slg{i}(\slc{1},\dots ,\slc{n}), \\
G^*_{\lambda} &:= G_{\lambda}(\slc{1}, \dots, \slc{n}), \notag
\end{align}
whenever $i=1, \dots, n$.
Note that this is  equivalent to the definitions of $G_i^*$ and $G^*$ that were given in \eqref{pq1}
and \eqref{pqinf}, respectively. Obviously we have that 
\begin{equation}\mylabel{e.e1}
\slg{i} = G^*_{\lambda} \cap G_i,
\end{equation}
for all $i=0,1,\dots,n$. Furthermore,
\begin{lemma}\mylabel{elm.l2}
\begin{subequations}\mylabel{102}
 \begin{align}
G^*_{\lambda} &= G^*(\lambda), \mylabel{102a} \\
\slgs{i} &= G^*_i(\lambda), \mylabel{102b}
\end{align}
\end{subequations}
for all $i=0,1,\dots,n$.
\end{lemma}

\begin{proof}
The fact that $\slc{i} \in \Irr(\slg{i} |\lambda)= \Irr(G_i(\lambda) | \lambda)$ is the 
$\lambda$-Clifford correspondent of $\chi_i \in \Irr(G_i)$  implies that 
$$
 G(\lambda,\subl{\chi}{1}, \dots ,\subl{\chi}{i})=  G(\chi_1, \dots ,\chi_i)(\lambda),
$$
But  $G ^* =  G (\chi_1, \dots , \chi_n)$, by \eqref{pqinf}, while 
$G(\lambda) =G_{\lambda}$ by \eqref{el.e3a}. Thus \eqref{102a} follows.

Equation \eqref{102b} follows easily from \eqref{102a} and \eqref{e.e1}.
\end{proof}

We can now prove       
\begin{proposition}\mylabel{elm:p1}
For every $r=0,1,\dots , k'$ and $i=1,\dots , l'$ we have that 
\begin{subequations}\mylabel{el.e5}
\begin{align}
\mylabel{el.e5a}  P^*_{2r} &\in \Hall_{\pi}(\slgs{2r}),     \\
\mylabel{el.e5b}
 Q^*_{2i-1}(\lambda)&\in \Hall_{\pi'}(\slgs{2i-1}). 
\end{align}
\end{subequations}
Therefore  the triangular set \eqref{el.e4b} can be chosen among the sets in its $G_{\lambda}$-conjugacy 
class so that it  satisfies 
\begin{subequations}\mylabel{104}
\begin{align}
P_{2r}^* &=\slps{2r}, \mylabel{104a}\\
Q_{2i-1}^*(\lambda) &=\slqs{2i-1}, \mylabel{104b}
\end{align}
\end{subequations}
whenever $0\leq r \leq k'$ and $1\leq i\leq l'$.
\end{proposition}

\begin{proof}
According to Proposition \ref{pq3-}, the group $P^*_{2k'}$
 is a $\pi$-Hall subgroup of $G^*_{2k'}$. 
Furthermore, $\lambda$ is  fixed by $P^*_{2k'}= \ma(\chi_1, \dots, \chi_n) \leq  \ma(\chi_1)$.
Hence, $P^*_{2k'}$ is also a $\pi$-Hall subgroup of $G^*_{2k'}(\lambda)$. 
As $P_{2r}^* = P_{2k'}^* \cap G_{2r}$,  
while $G_{2r}^*(\lambda)= G_{2k'}^*(\lambda)\cap G_{2r}
\unlhd G_{2k'}^*(\lambda)$, 
we also get that $P_{2r}^*$ is a $\pi$-Hall subgroup of $G_{2r}^*(\lambda)$ for  each  $r=0, 1, \dots, k'$ .
So \eqref{el.e5a} holds.

By Corollary \ref{q*:c1} we have that   $G^*_{2l'-1} = P^*_{2l'-2}\cdot Q^*_{2l'-1}$.
This, along with the fact that $P_{2l'-2}^* \leq  P^*_{2k'}$ fixes $\lambda$, implies that 
 $G^*_{2l'-1}(\lambda) = P^*_{2l'-2} \cdot Q^*_{2l'-1}(\lambda)$.
 Thus $Q^*_{2l'-1}(\lambda)$ 
is a $\pi'$-Hall subgroup of $G^*_{2l'-1}(\lambda) = \slgs{2l'-1}$.
Furthermore, for all $i=1,\dots, l'$, we have that $Q_{2i-1}^*(\lambda) = Q_{2l'-1}^* (\lambda)\cap G_{2i-1}$, where 
$G_{2i-1}  \unlhd G_{2l'-1}$.
Thus $Q_{2i-1}^*(\lambda)$ is a also  $\pi'$-Hall subgroup of  $G_{2i-1}^*(\lambda)$ for all $i=1,\dots, l'$.
This completes the proof of \eqref{el.e5}. 

The groups $\slps{2k'}= \slp{2} \cdots \slp{2k'}$ and $\slqs{2l'-1}=\slq{1} \cdots \slq{2l'-1}$
satisfy the conditions in  Propositions \ref{pq3-} and \ref{q*:p1} in the $\lambda$-situation, that is, 
$\slps{2k'} \in \Hall_{\pi}(\slgs{2k'})$ and $\slqs{2l'-1} \in \Hall_{\pi'}(\slgs{2l'-1})$.
Therefore, there exist  $G(\lambda)$-conjugates, 
$(\slps{2k'})^s, (\slqs{2l'-1})^s$, of $\slps{2k'}$ and $\slqs{2l'-1}$ respectively,  
such that $P_{2k'}^* = (\slps{2k'})^s$ and $Q_{2l'-1}^*(\lambda)=( \slqs{2l'-1})^s$.
Hence,  we also get that $P_{2r}^*= P_{2k'}^* \cap G_{2r}^* = (\slps{2r})^s$ and $Q_{2i-1}^*(\lambda)=
Q_{2l'-1}^* (\lambda)\cap G_{2l'-1}= ( \slqs{2i-1})^s$, whenever $0\leq r \leq k'$ and $1\leq i \leq l'$.
The set \eqref{el.e4b} was picked as any representative of a $G_{\lambda}= G(\lambda)$-conjugacy class
of triangular sets.
Thus we can pick  \eqref{el.e4b} to be the one that satisfies 
\eqref{104}.

So Proposition \ref{elm:p1} holds.
\end{proof}

The following is a straightforward but useful  lemma:
\begin{lemma}\mylabel{e.l1}
Assume that $Q_1 \leq T \leq G(\beta_1) $. Then $T = T(\lambda) \cdot Q_1$. 
Furthermore, if  $S$ satisfies 
  $Q_{1,2i} \leq  S \leq N(P_{2i}^* \tin G(\beta_{1,2i}))$,  for some $i=1,\dots,k'$, 
then $S =S(\lambda_{2i})\cdot Q_{1,2i} = S(\lambda) \cdot Q_{1,2i}$.
\end{lemma}
\begin{proof}
As the group $T$ fixes $\beta_1$, it permutes among themselves the $Q_1$-conjugacy class
of characters in $\Irr(E)$ lying under $\beta_1= \chi_1$. Since $\lambda$ is one of these characters, 
we have  $T \leq T(\lambda) \cdot Q_1$. 
The other inclusion is trivial. So $T = T(\lambda) \cdot Q_1$.

 The group $S$ normalizes $P_{2i}^*$ and thus normalizes 
$E_{2i} = C(P_{2i}^* \tin E)$. Furthermore,  it  fixes $\beta_{1,2i}$.
 Hence $S$     permutes among themselves the 
$Q_{1,2k}$-conjugacy class
of characters in $\Irr(E_{2i})$ lying under $\beta_{1,2i}$. Since $\lambda_{2i}$ lies in that class, we have 
  $S \leq S(\lambda_{2i}) \cdot Q_{1,2i}$. As the other inclusion is trivial, 
we get  $S = S(\lambda_{2i}) \cdot Q_{1,2i}$.
Since $S$ normalizes  $P_{2i}^*$,  and $\lambda_{2i}$ is the $P_{2i}^*$-Glauberman correspondent of 
$\lambda$, we obviously have that $S(\lambda_{2i})= S(\lambda)$.
Hence the lemma follows. 
\end{proof}

After these preliminary comments we are ready to state and prove
\begin{theorem}\mylabel{D}
The set \eqref{el.e4b}  chosen in Proposition \ref{elm:p1}  satisfies  
\begin{subequations}\mylabel{DD1}
\begin{align}
\slp{2r} &= P_{2r} \mylabel{D1.1} \\
\sla{2r} &= \alpha_{2r},\mylabel{D1}
\end{align}
\end{subequations}
for all $r=0,1,\dots ,k'$. And 
\begin{subequations}\mylabel{DD2}
\begin{align}
\slq{2i-1} &= Q_{2i-1}(\lambda) = Q_{2i-1}(\lambda_{2i-2}), \mylabel{D2.2} \\ 
\slb{2i-1} \in \Irr(\slq{2i-1})  &\text{ is the $\lambda_{2i-2}$-Clifford correspondent of  } 
\beta_{2i-1}\in \Irr(Q_{2i-1}),\mylabel{D2}
\end{align}
\end{subequations}
for all $i = 1, \dots ,l'$. (By convention  $\lambda_0 := \lambda$.)
Hence $\slb{2i-1}$ induces $\beta_{2i-1} \in \Irr(Q_{2i-1})$.
\end{theorem}
 \begin{proof}
According to  Lemma \ref{e.l1} for $T= Q_{2r-1}^*$ we get that $T= Q_{2r-1}^*(\lambda) Q_1$.
Since $Q_1= G_1 \unlhd G$, it follows that 
$N(Q_{2r-1}^*(\lambda)  \tin P_{2r}^*) \leq N(Q_{2r-1}^* \tin P_{2r}^* ) $.
But $P_{2r}^*$ fixes $\lambda$. So the other direction of the above inclusion also holds.
Thus
\begin{equation}\mylabel{el.e6.6}
 N(Q_{2r-1}^* \tin P_{2r}^* )  =  N(Q_{2r-1}^*(\lambda)  \tin P_{2r}^*). 
\end{equation}
The triangular set \eqref{el.e4b} satisfies the equivalent of  Proposition \ref{q*:p2} for the 
$\lambda$-situation.
Hence 
\begin{align*}
\slp{2r} &=N(\slqs{2r-1} \tin \slps{2r})  &\text{ by \eqref{q*:e7}} \\
&=N(Q_{2r-1}^*(\lambda)  \tin P_{2r}^*) &\text{ by \eqref{104}}  \\
&=N(Q_{2r-1}^* \tin P_{2r}^* ) &\text{ by \eqref{el.e6.6} } \\
&=P_{2r},  &\text{by \eqref{q*:e7}}
\end{align*}
for all $r= 1,\dots,k'$.
Thus  \eqref{D1.1} holds for all $r \geq 1$. It also holds for $r=0$, since $P_{0, \lambda}= 1= P_0$.

Similarly, for the $\pi'$-groups we have 
\begin{align*}
\slq{2i-1}&=N(\slps{2i-2} \tin \slqs{2i-1}) &\text{by \eqref{q*:e8}} \\
&=N(P_{2i-2}^* \tin Q_{2i-1}^*(\lambda) ) &\text{by \eqref{104}} \\
&=N(P_{2i-2}^* \tin Q_{2i-1}^*)(\lambda) &\\
&=Q_{2i-1}(\lambda), &\text{by \eqref{q*:e8}}
\end{align*}
for all $i=1,\dots,l'$. The group  $Q_{2i-1}$ normalizes $P_2, \dots, P_{2i-2}$ and thus normalizes 
their product  $P_{2i-2}^*$. Hence  \eqref{el.e3.4}, with $Q_{2i-1}$ in the place of $N$, 
implies that  $Q_{2i-1}(\lambda)= Q_{2i-1}(\lambda_{2i-2})$. So   \eqref{D2.2}  holds.

It remains to show \eqref{D1} and \eqref{D2} for  the $\lambda$-characters.
This will be done by induction on $i$ and $r$,   with the help of various observations that 
we write here separately  as steps.
\setcounter{step}{0}
\begin{step}\mylabel{el.st0}
For every $i=1,\dots, n$, the $c\slq{1}$-correspondent $\slc{i,1} \in \Irr(\slg{i,1})$ of 
$\slc{i} \in \Irr(\slg{i})$ induces the $cQ_1$-correspondent  $\chi_{i,1} \in \Irr(G_{i,1})$ 
of $\chi_i \in \Irr(G_i)$.
Even more,
\begin{subequations}\mylabel{el.es0}
\begin{align}
\slg{i,1}&= G_{i,1}(\lambda) \mylabel{el.es0a}  \, \text{ and } \\
 \slc{i,1} \in \Irr(\slg{i,1}) &\text{ is the $\lambda$-Clifford correspondent of
 $\chi_{i,1} \in \Irr(G_{i,1})$}. \mylabel{el.es0b}
\end{align}
\end{subequations}
\end{step}

\begin{proof}
We first remark that the $c\slq{1}$-correspondent (which is  analogous to the $cQ_1$-correspondent 
for the $\lambda$-case), is nothing else but a Clifford correspondent, as we can see in Table \ref{diagr.1}.
That is,  $\slg{i,1}= \slg{i}(\slc{1})$ and $\slc{i,1} \in \Irr(\slg{i,1})$ (see \eqref{i,1}, and \eqref{i.1}), 
 is  the $\slc{1}$-Clifford 
correspondent of $\slc{i} \in \Irr(\slg{i} |\slc{i})$, for all $i=1,\dots, n$.
Of course, $\slg{1,1}= \slg{1}$ and  $\slc{1,1} = \slc{1}$.

Furthermore, for all $i=1,\dots, n$, we have that  
\begin{align*} 
\slg{i,1} &= \slg{i}(\slc{1})  &\text{ by \eqref{i,1} for the $\lambda$-case }\\
&=G_i(\slc{1})(\lambda)  &\text{ as $\slg{1} = G_1(\lambda)$,   by \eqref{el.e3b} } \\
&=G_i(\chi_1, \lambda) &\text{ by Clifford's theory, since  $\chi_{1, \lambda}$ is the 
$\lambda$-Clifford correcpondent  of $\chi_1$  }\\
&=G_{i,1}(\lambda).  &\text{ by \eqref{i,1} }
\end{align*}
Therefore \eqref{el.es0a} holds.

Let $i=1,\dots, n$ be fixed.
The character $\slc{i,1} \in \Irr(\slg{i,1}) = \Irr(G_{i}(\chi_1, \lambda))$ induces 
$\slc{i} \in \Irr(\slg{i})=\Irr(G_i(\lambda))$, by \eqref{i.1}.
 Also 
the character $\slc{i} \in \Irr(G_i (\lambda))$  induces $\chi_i \in \Irr(G_i)$ by \eqref{101a}.
Therefore, $\slc{i,1} \in \Irr(G_{i}(\chi_1, \lambda))$ induces $\chi_i \in \Irr(G_i)$.
Hence $\slc{i,1} \in \Irr(G_{i}(\chi_1, \lambda))$ induces a character $\Psi \in \Irr(G_i(\chi_1))$.
Note that the induced character  $\Psi^{G_i}$ is $\slc{i,1}^{G_i} = \chi_i$.
Even more, the character $\Psi$ lies above $\chi_1= \beta_1$.    To see this, first note
 that, according to 
Lemma \ref{e.l1} for $T= G_i(\chi_1)$,
 we have $G_i(\chi_1) = G_i(\chi_1) (\lambda) \cdot Q_1= G_i(\chi_1, 
\lambda) \cdot G_1$, while $G_i(\chi_1,\lambda) \cap G_1= G_1(\lambda)$.
This, along with Mackey's Theorem (see Problem 5.6 in \cite{is}, or the special case in Problem 5.2 
in \cite{is}), implies  $(\slc{i,1}^{G_i(\chi_1)})|_{G_1}= (\slc{i,1}|_{G_1(\lambda)})^{G_1}$.
Hence 
$$
\langle\Psi|_{G_1}, \chi_1\rangle = \langle(\slc{i,1}^{G_i(\chi_1)})|_{G_1}, \chi_1 \rangle= 
\langle(\slc{i,1}|_{G_1(\lambda)})^{G_1}, \chi_1\rangle \ne 0, 
$$
where the last inequality holds as $\chi_1= \slc{1}^{G_1}$ (by \eqref{101a}), and
 $\slc{1} \in \Irr(G_1(\lambda))$  lies under $\slc{i,1}$.
 Thus the character $\Psi \in \Irr(G_i(\chi_1))$, induces $\chi_i \in \Irr(G_i)$ and lies 
above  $\chi_1$. Hence Clifford's theorem implies that $\Psi$ is the unique $\chi_1$-Clifford 
correspondent of $\chi_i$. This, along with \eqref{i.1}, implies that $\Psi =\chi_{i,1}$, i.e., 
$\slc{i,1}^{G_{i,1}} = \Psi = \chi_{i,1}$. 

So  $\slc{i,1} \in \Irr(\slg{i,1})=  \Irr(G_{i,1}(\lambda))$ induces $\chi_{i,1} \in \Irr(G_{i,1})$, 
and lies above $\lambda$.
 Thus $\slc{i,1}$ is the $\lambda$-Clifford correspondent of $\chi_{i,1}$.

This completes the proof of the first step. 
\end{proof}

\begin{step}\mylabel{el.st1}
For every $i=2, \dots, l'$ we have that 
\begin{equation}\mylabel{el.e7}
Q_{2i-1} = Q_{2i-1}(\lambda) \cdot Q_{1,2i-2}. 
\end{equation} 
Furthermore, for every $r=1,\dots,k'$ and every $i=1,\dots,l'-1$  we have that 
\begin{subequations}\mylabel{el.e6}
\begin{multline}\mylabel{el.e6a}
N(\slp{0}, \slp{2}, \dots, \slp{2i}, \slq{1}, \dots, \slq{2i-1} \tin \slg{2i+1}(\slc{1}, \dots, \slc{2i}))=\\
N(P_0, P_2, \dots, P_{2i}, Q_1, \dots, Q_{2i-1} \tin \slg{2i+1}(\slc{1}, \dots, \slc{2i})  ), 
\end{multline} 
and 
\begin{multline}\mylabel{el.e6b}
N(\slp{0}, \slp{2}, \dots, \slp{2r-2}, \slq{1}, \dots, \slq{2r-1} \tin \slg{2r}(\slc{1}, \dots, \slc{2r-1}))=\\
N(P_0, P_2, \dots, P_{2r-2}, Q_1, \dots, Q_{2r-1} \tin \slg{2r}(\slc{1}, \dots, \slc{2r-1})  ).
\end{multline}
 \end{subequations}
Thus, for all $i=1, \dots, l'-1$, the 
$ c\slp{2}, \dots, c\slp{2i}, c\slq{3}, \dots, c\slq{2i-1}$-correspondent 
$\slc{2i+1, 2i}$, of $\slc{2i+1, 1}$ coincides with the
   $ cP_2, \dots, cP_{2i}, cQ_3, \dots, cQ_{2i-1}$-correspondent of $\slc{2i+1, 1}$.
By convention, if  $i=1$ this is only the $cP_2= c\slp{2}$-correspondence. 
Similarly, for all $r=2,\dots, k'$, the 
$ c\slp{2}, \dots, c\slp{2r-2}, c\slq{3}, \dots, c\slq{2r-1}$-correspondent 
$\slc{2r, 2r-1}$, of $\slc{2r,1}$ coincides with the
   $ cP_2, \dots, cP_{2r-2}, cQ_3, \dots, cQ_{2r-1}$-correspondent of $\slc{2r,1}$.
\end{step}
\begin{proof}
For all $t=2,\dots,l'$ we have that 
$$
Q_{1,2t-2}  = Q_1 \cap Q_{2t-1} \leq Q_{2t-1} \leq N(P_{2t-2}^* \tin G(\beta_{1,2t-2})), 
$$
(by \eqref{pq14b} we get the   equality, while Proposition \ref{pqremark1'} shows that 
$Q_{2t-1}$ fixes $\beta_{1,2t-2}$).
Hence Lemma \ref{e.l1}, for $Q_{2t-1}$ in the place of $S$, implies that 
$$
Q_{2t-1} = Q_{2t-1}(\lambda) \cdot Q_{1,2t-2}, 
$$ 
for all $t=2,\dots, l'$. This proves \eqref{el.e7}.

We have already seen that $ \slp{2r}= P_{2r}$ while $\slq{2i-1} = Q_{2i-1}(\lambda)$, whenever
 $1\leq r \leq k'  $ and $1\leq i \leq l'$ (by \eqref{D1.1} and \eqref{D2.2}). Thus \eqref{el.e7} implies that 
$$Q_{2t-1} = \slq{2t-1} \cdot Q_{1,2t-2}, 
$$
for all $t=2,\dots, l'$.
Any subgroup of $G$ that normalizes $P_2, \dots, P_{2t-2}$ and $\slq{2t-1}$ also normalizes 
$Q_{1,2t-2}= N(P_{2t-2}^* \tin Q_1)$. Thus it normalizes 
$Q_{2t-1}=\slq{2t-1} \cdot Q_{1,2t-2} $, whenever $t=2,\dots, l'$.  
Since $Q_1 \unlhd G$,  is normalized
by any subgroup of $G$, we therefore get  

\begin{multline*}
N(\slp{0}, \slp{2}, \dots, \slp{2i}, \slq{1}, \dots, \slq{2i-1} \tin \slg{2i+1}(\slc{1}, \dots, \slc{2i}))=\\
N(P_{0}, P_{2}, \dots, P_{2i}, \slq{1}, \dots, \slq{2i-1} \tin \slg{2i+1}(\slc{1}, \dots, \slc{2i}))  \leq \\
N(P_{0}, P_{2}, \dots, P_{2i}, Q_{1}, \dots, Q_{2i-1} \tin \slg{2i+1}(\slc{1}, \dots, \slc{2i})).
\end{multline*}
 The other inclusion is trivial as  $\slg{2i+1}= G_{2i+1}(\lambda)$. So everything in $\slg{2i+1}$
 that normalizes 
$Q_{2t-1}$ also normalizes $\slq{2t-1}= Q_{2t-1}(\lambda)$, for all $t=1,\dots,l'$. Also 
$Q_1 \unlhd G$ is normalized
by any subgroup of $G$. Hence we have equality, and \eqref{el.e6a} is proved.

The proof for \eqref{el.e6b} is similar. So we omit it.

Let $i=1, \dots, l'-1$ be fixed.
For all $t=1,\dots, i$,  we have  $P_{2t}= \slp{2t}$. So to prove the next  two statements 
of  Step \ref{el.st1}, it suffices 
to show that  the $cQ_{2t-1}$-correspondence 
coincides with the $c\slq{2t-1}$-correspondence,  for all $t=2\dots, i$.
We remark that  the $c\slq{2t-1}$-correspondence was used,
in Theorem \ref{tow--tri} 
 and Table \ref{diagr.1-2-3...} for the $\lambda$-situation, 
 to get the irreducible character 
$\slc{2i+1, 2t-1} $ of  $ N(\slq{2t-1} \tin \slg{2i+1, 2t-2}(\slc{2t-1}))
= \slg{2i+1, 2t-1} $ 
from $\slc{2i+1, 2t-2} \in \Irr(\slg{2i+1,2t-2})$,  whenever $t=2,\dots, i$.
 That is, we applied Lemma \ref{dade'} to the groups 
$$
\slg{2t-1, 2t-2}= \slq{2t-1} \ltimes \slp{2t-2} \unlhd \slg{2t, 2t-2} \unlhd \dots \unlhd 
\slg{2i-1,2t-2}
$$
and their irreducible  characters
$$
\slc{2t-1,2t-2}=\slb{2t-1} \cdot \sla{2t-2}^e, \slc{2t,2t-2},  \dots ,  \slc{2i-1,2t-2},
$$
 in the place of $N \unlhd K_1 \unlhd \dots \unlhd K_r$ and 
$\{\chi=\alpha \cdot \beta^e, \chi_1, \dots, \chi_r\}$ respectively.
On the other hand, for all  $t=2, \dots, i$,   the character $\slc{2i+1,2t-2}$
has a $cQ_{2t-1}$-correspondent. Indeed, the group $Q_{2t-1}$ acts on $P_{2t-2} = 
\slp{2t-1}$, while the  semidirect product 
$Q_{2t-1} \rtimes \slp{2t-1}$ is normalized by 
all the groups in the normal series 
 $$
P_{2t-1, \lambda} \unlhd \slg{2t,2t-1} \unlhd  \dots \unlhd \slg{2i-1, 2t-1}.
$$
Notice that $\slq{2t-1}$ and $Q_{2t-1}$ have the same image in $\Aut(\slp{2t-2})$,  as 
$Q_{2t-1} = \slq{2t-1} \cdot Q_{1,2t-2}$ (by \eqref{el.e7}), and  $Q_{1,2t-2}$ centralizes
 $\slp{2t-2}=P_{2t-2}$, whenever $2\leq t \leq i$. 
 Furthermore,  $\slg{2i+1, 2t-2}$ normalizes 
$P_2, \dots, P_{2t-2}$ (see \eqref{tow--tri3}),   and thus normalizes $Q_{1,2t-2}$. Hence 
 $N(\slq{2t-1} \tin \slg{2i+1, 2t-2})  \leq N(Q_{2t-1} \tin \slg{2i+1, 2t-2}) $. The other inclusion holds 
trivially,  as $\slg{2i+1, 2t-2} \leq \slg{2i+1} = G_{2i+1}(\lambda)$ fixes $\lambda$ and 
$Q_{2t-1, \lambda} =  Q_{2t-1}(\lambda)$.
Hence $N(\slq{2t-1} \tin \slg{2i+1, 2t-2})  = N(Q_{2t-1} \tin \slg{2i+1, 2t-2}) $, for all $t=2, \dots, i$.
Therefore Proposition \ref{dade:p1.5} 
 implies that the $c\slq{2t-1}$-correspondent $\slc{2i+1,2t-1}$ of
$\slc{2i+1,2t-2}$ coincides with its  $cQ_{2t-1}$-correspondent, whenever $2\leq t \leq i$.

We conclude  that the 
$ c\slp{2}, \dots, c\slp{2i}, c\slq{3}, \dots, c\slq{2i-1}$-correspondent 
$\slc{2i+1, 2i}$, of $\slc{2i+1,1}$ coincides  with the
   $ cP_2, \dots, cP_{2i},cQ_3, \dots, cQ_{2i-1}$-correspondent of $\slc{2i+1,1}$, for all
 $i=1, \dots, l'-1$.

Similarly we can work with the character $\slc{2r,2r-1}$,  the group $\slg{2r, 2t-1}$ and its 
normal subgroup $\slg{2t, 2t-1}= \slp{2t} \ltimes \slq{2t-1}$, for all $t=2, \dots, r$,
 for some fixed $r=2,\dots, k'$. Thus we get that   that 
the $c\slq{2t-1}$-correspondent $\slc{2r, 2t-1}$ of $\slc{2r, 2t-2}$ coincides with the 
$cQ_{2t-1}$-correspondent.
We conclude similarly  that the $ c\slp{2}, \dots, c\slp{2r-2},
 c\slq{3}, \dots, c\slq{2r-1}$-correspondent 
$\slc{2r, 2r-1}$, of $\slc{2r,1}$ coincides with the
   $ cP_2, \dots, cP_{2r-2}, cQ_3, \dots, cQ_{2r-1}$-correspondent of $\slc{2r,1}$, for all
 $r=2,\dots, k'$.

This completes the proof of Step \ref{el.st1}.
\end{proof}

\begin{step}\mylabel{el.st2}
For all $i=1,\dots, l'-1$, we have that $\slg{2i+1, 2i} =G_{2i+1,2i}(\lambda)$, while 
  the character $\slc{2i+1,2i} \in \Irr(\slg{2i+1,2i})$  induces  $\chi_{2i+1, 2i} \in 
\Irr(G_{2i+1, 2i})$. Similarly,  for all $r=1,\dots, k'$, we get that
 $\slg{2r,2r-1} =G_{2r, 2r-1}(\lambda)$, while 
the character $\slc{2r,2r-1} \in \Irr(\slg{2r, 2r-1})$
induces $\chi_{2r,2r-2} \in \Irr(G_{2r,2r-1})$. 
\end{step}

\begin{proof}
According to \eqref{tow--tri3} for the $\lambda$-case we have that 
$$
\slg{2i+1,2i} =
 N(\slp{0}, \slp{2}, \dots, \slp{2i}, \slq{1}, \dots, \slq{2i-1} \tin \slg{2i+1}(\slc{1}, \dots, \slc{2i})). 
$$
This,  along with \eqref{el.e6a}, \eqref{102c}  and \eqref{102b},  implies that 
\begin{multline*}
\slg{2i+1, 2i} = N(P_0, P_2, \dots, P_{2i}, Q_1, \dots, Q_{2i-1} \tin \slg{2i+1}(\slc{1}, \dots, \slc{2i})  )=\\
 N(P_0, P_2, \dots, P_{2i}, Q_1, \dots, Q_{2i-1} \tin \slgs{2i+1}   ) = 
 N(P_0, P_2, \dots, P_{2i}, Q_1, \dots, Q_{2i-1} \tin G_{2i+1}^*(\lambda)   )=\\
  N(P_0, P_2, \dots, P_{2i}, Q_1, \dots, Q_{2i-1} \tin G_{2i+1}^*)(\lambda).
\end{multline*}
As $G_{2i+1,2i} =  N(P_0, P_2, \dots, P_{2i}, Q_1, \dots, Q_{2i-1} \tin G_{2i+1}^*)$ (by \eqref{tow--tri3}), 
we have that  $\slg{2i+1, 2i} = G_{2i+1,2i}(\lambda)$ , for all $i=1,\dots, l'-1$.
Furthermore, 
according to Step \ref{el.st1}, the character 
$\slc{2i+1,2i}$ is the $cP_2, \dots, cP_{2i}, cQ_3, \dots, cQ_{2i-1}$-correspondent of $\slc{2i+1,1}$, 
whenever $i=1,\dots, l'-1$.
But  $\slc{2i+1,1}\in \Irr(G_{2i+1,1}(\lambda))$ induces $\chi_{2i+1,1}$  in $G_{2i+1,1}$ according 
to Step \ref{el.st0},  for all such $i$. Furthermore,   the $cA$-correspondence 
(for arbitrary $A$) respects induction (see  Theorem \ref{dade:t2}).
Hence $\slc{2i+1,2i} $ induces the 
$ cP_2, \dots, cP_{2i}, cQ_3, \dots, cQ_{2i-1}$-correspondent character of $\chi_{2i+1,1}$
in the normalizer $N(P_0, P_2, \dots, P_{2i}, Q_1, \dots,
 Q_{2i-1} \tin G_{2i+1}(\chi_1, \dots, \chi_{2i}))= G_{2i+1,2i}$. 
As this correspondent character of $\chi_{2i+1,1}$ is  $\chi_{2i+1,2i}$ (by Theorem \ref{tow--tri}), we conclude that $\slc{2i+1, 2i}$ induces $\chi_{2i+1,2i}$, for all $i=1,\dots,l'-1$.
This completes the first part of Step \ref{el.st2}. 

For the second part, we first remark that the case $r=1$ has been done in Step \ref{el.st0}. 
Indeed, by \eqref{el.es0a}, 
we have that $\slg{2,1}= G_{2,1}(\lambda)$, while by \eqref{el.es0b} the character 
$\slc{2,1}$ induces $\chi_{2,1}$. 
The rest of  proof for $r=2,\dots,k'$ is analogous to the proof of the first part, with the use of 
 \eqref{el.e6b} in the place of \eqref{el.e6a}.
 So we omit it. 
\end{proof}

We can now continue with the  proof of \eqref{D1} and \eqref{D2}  for the $\lambda$-characters.
If $r=0$ then $\sla{0}= 1= \alpha_0$. Hence \eqref{D1} holds trivially for $r=0$.
Furthermore, $\slb{1} = \slc{1}$, by \eqref{x2}. But $\slc{1}$ is the $\lambda$-Clifford correspondent 
of $\chi_1 \in \Irr(Q_1 | \lambda)$. Thus \eqref{D2} holds for $i=1$.

Using an inductive argument we will prove that,  if \eqref{D2} holds when $i$ is some integer
 $t=1, \dots, k'-1$,  then \eqref{D1}
holds for $r=t$.
 Symmetrically if \eqref{D1} holds when $r$ is   some integer $s=1, \dots , l'-1$, 
then  \eqref{D2} holds for $i=s+1$.
This is enough to prove that \eqref{D1} and \eqref{D2} hold for all $r=0, \dots, k'$ and all $i=1,\dots , l'$, 
respectively.

Assume that \eqref{D2} holds for $i=t$. That is, 
$\slb{2t-1}$ is the $\lambda_{2t-2}$-Clifford correspondent of $\beta_{2t-1}$. Therefore
\begin{equation}\mylabel{el.e8}
\slb{2t-1}^{Q_{2t-1}} = \beta_{2t-1}.
\end{equation}

Furthermore, Theorem \ref{tow--tri},  and in particular \eqref{tow--tri1}, implies that  
the character $\sla{2t}$ was picked as the unique character of
 $\slp{2t} $ 
that satisfies
\begin{subequations}\mylabel{el.e9}
\begin{equation}\mylabel{el.e9a}
\slc{2t, 2t-1}= \sla{2t}  \cdot \slb{2t-1}^e, 
\end{equation}
where $\slb{2t-1}^e$ is the canonical extension of $\slb{2t-1} \in \Irr(\slq{2t-1})$ to 
$\slg{2t, 2t-1}$. Similarly, $\alpha_{2t}$ is the unique character of $P_{2t}  $ such that 
\begin{equation}\mylabel{el.e9b}
\chi_{2t,2t-1}= \alpha_{2t} \cdot \beta_{2t-1}^e.
\end{equation}
\end{subequations}

According to Step \ref{el.st2} the character $\slc{2t, 2t-1}$ induces  $\chi_{2t, 2t-1}  \in
 \Irr(G_{2t, 2t-1})$. Also $G_{2t, 2t-1} = P_{2t} Q_{2t-1}$,  see \eqref{tow--tri1}. 
This, along with \eqref{el.e9}, implies 
\begin{equation} \mylabel{el.e10}
 \chi_{2t, 2t-1} = (\slc{2t, 2t-1})^{G_{2t,2t-1} }
 = (\sla{2t}  \cdot \slb{2t-1}^e) ^{G_{2t,2t-1}}.
\end{equation}
    Lemma \ref{100d}  can be applied to the group $G_{2t, 2t-1} = P_{2t} \ltimes Q_{2t-1}$ and the 
 characters $\sla{2t} \in \Irr(\slp{2t})= \Irr(P_{2t})$ and $\slb{2t-1}^e \in \Irr(P_{2t} \ltimes \slq{2t-1})$. Thus    
\begin{equation}\mylabel{el.e11}
 (\sla{2t}  \cdot \slb{2t-1}^e) ^{G_{2t,2t-1}}= \sla{2t} \cdot (\slb{2t-1}^e)^{G_{2t,2t-1}}.
\end{equation}
The groups $G_{2t, 2t-1}= P_{2t} \ltimes Q_{2t-1} = (P_{2t} \slq{2t-1}) \cdot Q_{2t-1}, 
Q_{2t-1} , P_{2t} \slq{2t-1}=\slg{2t,2t-1}$ and $\slq{2t-1}$, along with the character $\slb{2t-1} \in 
\Irr(\slq{2t-1})$, satisfy   the hypothesis of Proposition \ref{prel:100B} in the place of the groups
$G , N , K $ and $H$ respectively. So we conclude that 
$$
(\slb{2t-1}^e)^{G_{2t,2t-1}} = (\slb{2t-1}^{Q_{2t-1}}) ^e, 
$$
where $ (\slb{2t-1}^{Q_{2t-1}}) ^e$ is the canonical 
extension of $\slb{2t-1}^{Q_{2t-1}} \in \Irr(Q_{2t-1})$ to 
$G_{2t, 2t-1} $. But, according to the inductive hypothesis,  the character $\slb{2t-1}$ satisfies
 \eqref{el.e8}.
Therefore
$$
(\slb{2t-1}^e)^{G_{2t,2t-1}} = (\slb{2t-1}^{Q_{2t-1}}) ^e= \beta_{2t-1}^e, 
$$
where $\beta_{2t-1}^e$ is the canonical extension of $\beta_{2t-1}$ to $G_{2t, 2t-1}$.
This,  along with \eqref{el.e11} and \eqref{el.e10}, implies that 
$$
\chi_{2t,2t-1}=
(\sla{2t}  \cdot \slb{2t-1}^e) ^{G_{2t,2t-1}}= \sla{2t} \cdot \beta_{2t-1}^e.
$$
Therefore \eqref{el.e9b} holds with $\sla{2t}$ in the place of $\alpha_{2t}$.
As $\alpha_{2t}$ is the unique character of $P_{2t}$  satisfying \eqref{el.e9b}, 
we must have $\alpha_{2t}=\sla{2t}$. Therefore, \eqref{D1} holds for $r=t$.

Now assume that \eqref{D1} holds when $r$ is  some integer $s=1, \dots, l'-1$. 
We will prove, in an argument  similar  to the one we just gave, that 
\eqref{D2} also holds for $i=s+1$.  

According to \eqref{D1} for $r=s$ we get 
\begin{equation}\mylabel{el.e12}
\alpha_{2s} = \sla{2s}.
\end{equation}
Furthermore, equations \eqref{tow--tri2} imply that 
$\slb{2s+1}$ and $\beta_{2s+1}$ are the unique characters of $\slq{2s+1}$ and $Q_{2s+1}$, respectively, 
that satisfy 
\begin{subequations}\mylabel{el.e13}
\begin{align}
 \slc{2s+1,2s}&= \sla{2s}^e \cdot \slb{2s+1} \mylabel{el.e13a} \\
\chi_{2s+1,2s} &=\alpha_{2s}^e \cdot \beta_{2s+1} \mylabel{el.e13b},
\end{align}
\end{subequations}
where $\sla{2s}^e$ and $\alpha_{2s}^e$ are the canonical  extensions  of $\sla{2s}\in \Irr(\slp{2s})$ and
 $\alpha_{2s} \in \Irr(P_{2s})$ to 
$\slg{2s+1,2s}= \slp{2s} \rtimes \slq{2s+1}$ and $G_{2s+1,2s}= P_{2s} \rtimes Q_{2s+1}$, respectively.
But $\slp{2s}=P_{2s}$ and $\sla{2s}=\alpha_{2s}$ by \eqref{el.e12}.
Furthermore, $\slg{2s+1,2s} \leq G_{2s+1,2s}$ (see Step \ref{el.st2}).
Therefore 
\begin{equation}\mylabel{el.e14}
\sla{2s}^e= \alpha_{2s}^e|_{\slg{2s+1,2s}}.
\end{equation}
The fact that $\slc{2s+1, 2s}$ induces $\chi_{2s+1, 2s}$ (see Step \ref{el.st2}), along with \eqref{el.e13}
 and \eqref{el.e14}, implies 
\begin{equation}\mylabel{el.e15}
\chi_{2s+1,2s} = (\alpha_{2s}^e|_{\slg{2s+1,2s}}\cdot \slb{2s+1})^{G_{2s+1,2s}}.
\end{equation}
Using the isomorphism  $\slq{2s+1}\cong \slq{2s+1} \ltimes P_{2s} / P_{2s}$,  we 
denote the inflation of $\slb{2s+1}$ to $\slq{2s+1} \ltimes \slp{2s} = 
\slq{2s+1} \ltimes \slp{2s}= \slg{2s+1, 2s}$ as $\slb{2s+1}^i$.
So \eqref{el.e15} becomes
$$
\chi_{2s+1,2s} = (\alpha_{2s}^e|_{\slg{2s+1,2s}}\cdot \slb{2s+1}^i)^{G_{2s+1,2s}}.
$$
Hence  we can apply  Lemma \ref{100c} to  get  
$$
 (\alpha_{2s}^e|_{\slg{2s+1,2s}}\cdot \slb{2s+1}^i)^{G_{2s+1,2s}}= \alpha_{2s}^e \cdot 
(\slb{2s+1}^i)^{G_{2s+1,2s}}.
$$
This, along with \eqref{el.e15}, implies that 
$$
\chi_{2s+1,2s} = \alpha_{2s}^e \cdot (\slb{2s+1}^i)^{G_{2s+1,2s}}= 
 \alpha_{2s}^e \cdot (\slb{2s+1}^i)^{P_{2s} \rtimes Q_{2s+1}}= 
\alpha_{2s}^e \cdot (\slb{2s+1}^{Q_{2s+1}})^i,
$$
where  $(\slb{2s+1}^{Q_{2s+1}})^i$ denotes the inflation of  $\slb{2s+1}^{Q_{2s+1}} \in \Irr(Q_{2s+1})
= \Irr(P_{2s}Q_{2s+1}/P_{2s})$ to $P_{2s} \rtimes Q_{2s+1}$.
Note that  $\alpha_{2s}^e \cdot (\slb{2s+1}^{Q_{2s+1}})^i$  is equal to 
  $\alpha_{2s}^e \cdot (\slb{2s+1}^{Q_{2s+1}})$  by the definition of the ltter product.
But  $\beta_{2s+1}$ is the unique character of $Q_{2s+1}$ that satisfies \eqref{el.e13b}.
 Hence 
\begin{equation}\mylabel{el.e16}
(\slb{2s+1}^{Q_{2s+1}})^i= \beta_{2s+1}^i,   \text{  and thus \, } 
\slb{2s+1}^{Q_{2s+1}}= \beta_{2s+1}.
\end{equation}
Furthermore, the character $\chi_{2s+1, 2s}$ lies above $\chi_{1, 2s} \in \Irr(G_{1,2s})= \Irr(Q_{1,2s})$.
Also  $\chi_{1,2s}= \beta_{1,2s}$ (as $\chi_1 = \beta_1$), lies above $\lambda_{2s}$ by \eqref{el.e3.3}.
This, along with the fact that $\alpha_{2s, \lambda}^e$ is trivial on $Q_{1, 2s} = G_{1, 2s}$, 
implies that  $\slb{2s+1}\in \Irr(\slq{2s+1})$  lies above $\lambda_{2s}$ and induces $\beta_{2s+1}
\in \Irr(Q_{2s+1})$.
As $\slq{2s+1}= Q_{2s+1}(\lambda)= Q_{2s+1}(\lambda_{2s})$, by \eqref{el.e3.4} with $N= Q_{2s+1}$, 
 we conclude that $\slb{2s+1}$ is the $\lambda_{2s}$-Clifford 
correspondent of $\beta_{2s+1}$. Hence \eqref{D2} holds for $i=s+1$.

This completes the proof of Theorem \ref{D}. 
\end{proof}

Assume that $i, j$ satisfy $1\leq i \leq j \leq k'$.
According to Lemma 2.5 in \cite{wo1}, Glauberman correspondence 
is compatible with 
Clifford theory. This,  along with 
 \eqref{D2} and \eqref{D1.1},  implies that 
the $P_{2i}\cdots P_{2j}= \slp{2i}\cdots
 \slp{2j}$-Glauberman correspondent 
$\beta_{2i-1,2j}$ of $\beta_{2i-1}$  is induced by  
 the $\slp{2i}\cdots \slp{2j}$-Glauberman correspondent $\slb{2i-1,2j}$
 of $\slb{2i-1}$.
 Furthermore, $\slb{2i-1,2j}$ lies above the $P_{2i}\cdots 
P_{2j}$-Glauberman 
correspondent $\lambda_{2j}$ of $\lambda_{2i-2}$ (see \eqref{el.e3.1}), as $\slb{2i-1}$ lies above
 $\lambda_{2i-2}$.
Since $\slq{2i-1}= Q_{2i-1}(\lambda)$, we also have that 
$$
\slq{2i-1,2j} = C(P_{2i} \cdots P_{2j} \tin \slq{2i-1})= 
C(P_{2i}\cdots P_{2j} \tin Q_{2i-1})(\lambda) = Q_{2i-1,2j}(\lambda),
$$
whenever $1\leq i \leq j \leq k'$.
But $Q_{2i-1,2j}(\lambda) = Q_{2i-1,2j}(\lambda_{2j})$, as $Q_{2i-1,2j}$ normalizes $E$ and $P_{2j}^*$.
Therefore, 
\begin{remark}\mylabel{el.rem1}
For every $i, j$ with $1\leq i \leq j \leq k'$, we have 
$$
\slq{2i-1,2j}= Q_{2i-1,2j}(\lambda) = Q_{2i-1,2j}(\lambda_{2j}),
$$
while 
 the character $\slb{2i-1,2j}$ is the $\lambda_{2j}$-Clifford 
correspondent of $\beta_{2i-1,2j}$.
\end{remark}

Furthermore, $Q_{1,2i-2}$ centralizes $P_2, \dots, P_{2i-2}$  by \eqref{pq14e},
and  $Q_{2i-1}= Q_{2i-1}(\lambda) \cdot Q_{1,2i-2}$ by \eqref{el.e7}.
We conclude that 
\begin{equation}\mylabel{elm.ee1}
C(Q_{2i-1}(\lambda)  \tin P_{2r}) = C(Q_{2i-1} \tin P_{2r}), 
\end{equation}
whenever  $1\leq r < i \leq l'$.
Hence, 
\begin{align*}
P_{2r, 2i-1, \lambda}&= C(\slq{2r+1}, \dots, \slq{2i-1} \tin \slp{2r}) &\\
&=C(Q_{2r+1}(\lambda), \dots, Q_{2i-1}(\lambda) \tin P_{2r}) &\text{
 by \eqref{D1.1} and \eqref{D2.2} }\\
&=C(Q_{2r+1}, \dots, Q_{2i-1} \tin P_{2r} ) &\text{ by \eqref{elm.ee1} }\\
&=P_{2r,2i-1},
\end{align*}
whenever $1\leq r < i \leq l'$.

Even more, the $Q_{2r+1} \cdots Q_{2i-1}$-Glauberman correspondent
 $\alpha_{2r, 2i-1}
\in \Irr(P_{2r, 2i-1})$ of  the irreducible character $\alpha_{2r}$ of $ P_{2r}$,
 (see \eqref{pq13c}),
coincides  with the $Q_{2r+1}(\lambda)  \cdots   Q_{2i-1}(\lambda)$-Glauberman 
correspondent  of $\alpha_{2r}$, by  Corollary \ref{coo.d1}.
In conclusion,
\begin{remark}\mylabel{el.rem11}
For every $r, i$ with $1\leq r < i  \leq l'$, we have 
$$
\slp{2r, 2i-1}= P_{2r, 2i-1} \text{ and } \sla{2r,2i-1}=\alpha_{2r,2i-1}.
$$
\end{remark}

The relation between  $\alpha^*_{2i}$ and $\sla{2i}^*$ is an easy corollary of 
Theorem \ref{D}.
\begin{corollary}\mylabel{elm:co1}
For all $r=0,\dots ,k'$  we have
$$\sla{2r}^* = \alpha_{2r}^* \in \Irr(P^*_{2r}) = \Irr(\slps{2r}).$$
\end{corollary}
\begin{proof}
By \eqref{104a} and \eqref{D1}, we have $P^*_{2r} = \slps{2r}$ and 
$\sla{2r}  = \alpha_{2r}$, for all $r=0, 1,\dots,k'$.
Furthermore,  the character $\alpha^*_{2r}$ is uniquely determined
by $\alpha_{2r}$ and the one  to one 
 $Q_{2j+1}$-correspondence 
$$
\alpha_{2r, 2j-1}^* \overleftrightarrow{Q_{2j+1}} \alpha_{2r, 2j+1}^*.
$$
The latter is a correspondence 
between all characters $\alpha_{2r, 2j+1}^* 
\in \Irr(P_{2j+2} \cdots P_{2r})$ and all characters 
$\alpha_{2r,2j-1}^* \in \Irr(P_{2j} \cdot P_{2j+2} \cdots P_{2r})$  
lying over some $Q_{2j+1}$-invariant character of $P_{2j}$,
for all 
$j=1,\dots,r-1$, as 
Lemma \ref{p*l1} and  Theorem \ref{p*t1} imply. (Note that 
$\alpha^*_{2r}\in \Irr(P^*_{2r})$ is the 
$Q_3,  Q_5, \dots , Q_{2r-1}$-correspondent
of $\alpha_{2r}\in \Irr(P_{2i})$.)

Because $P_{2j+2} \cdots P_{2r}$ normalizes $Q_{2j+1}$ and fixes $\lambda$, 
(by Remark \ref{elm.r00}), it normalizes $Q_{2j+1}(\lambda)=Q_{2j+1,\lambda}$.
The subgroup $Q_{1,2j}=C(P_2 \cdots P_{2j} \tin Q_1)$ centralizes $P_{2j}$.
But $ Q_{2j+1} =Q_{2j+1}(\lambda) \cdot Q_{1,2j} =  \slq{2j+1} \cdot Q_{1,2j}$,
according to \eqref{el.e7},
for all $j=1, \dots, r-1$. 
Hence
\begin{multline*}
P_{2j} \cap P_{2j+2} \cdots P_{2r} = N(Q_{2j+1} \tin P_{2j}) \\
= C(Q_{2j+1} \tin P_{2j}) = C(\slq{2j+1} \tin P_{2j}) = 
N(\slq{2j+1} \tin P_{2j}).
\end{multline*}
Therefore
 \begin{multline*}
N(\slq{2j+1} \tin P_{2j} \cdot P_{2j+2} 
\cdots P_{2r})= N(\slq{2j+1} \tin P_{2j})\cdot P_{2j+2} \cdot P_{2r} \\
=N(Q_{2j+1} \tin P_{2j}) \cdot P_{2j+2} \cdots P_{2r}=  N(Q_{2j+1} \tin 
 P_{2j} \cdot P_{2j+2}  \cdots P_{2r}).
\end{multline*}
This, along with Proposition \ref{dade:p1.5}, ,
 
implies that the above $Q_{2j+1}$-correspondence coincides with the
$\slq{2j+1}$-correspondence, for all $j=1,\dots, r-1$.
Hence the 
$Q_3, Q_5, \dots, Q_{2r-1}$-correspondent, $\alpha^*_{2r}\in \Irr(P^*_{2r})$,
 of $\alpha_{2r}=\sla{2r} \in \Irr(P_{2r})$
coincides with the $\slq{3}, \slq{5},\dots, \slq{2r-1}$-correspondent $\sla{2r}^* \in \Irr(\slps{2r})$
of $\alpha_{2r} = \sla{2r} \in\Irr(\slp{2r})$, i.e., 
 $\alpha^*_{2r} = \sla{2r}^*$. So the  corollary follows.
\end{proof}

What about the groups $\ma$ and $\mb$? How is a Hall system for $G(\lambda)$ 
that satisfies the analogue  of  \eqref{e.AB}  for the $\lambda$-case 
related to $\ma,\mb$? The answer is given in 
\begin{theorem}\mylabel{el.tAB}
We can find  $\ma_{\lambda} \in \Hall_{\pi}(G_{\lambda})$ and $\mb_{\lambda} \in \Hall_{\pi'}(G_{\lambda})$ 
satisfying  the 
equivalent  of \eqref{e.AB} for the $\lambda$-groups, along with 
\begin{subequations}\mylabel{el.eAB}
\begin{gather}
\ma_{\lambda} (\slc{1}, \dots, \slc{h}) = \ma(\chi_1, \dots, \chi_h), \\
\mb_{\lambda}(\slc{1} , \dots, \slc{h}) = \mb(\chi_1, \dots, \chi_h, \lambda),
\end{gather}
\end{subequations}
for all $h=1,\dots,n$.
\end{theorem}

\begin{proof}
It suffices to show that $\ma(\chi_1,\dots, \chi_h)$ and 
$\mb(\chi_1, \dots, \chi_h, \lambda)$ satisfy (\ref{e.AB}b,c) for the
 $\lambda$-groups.
We already know, by (\ref{e.AB}c), 
 that $\ma(\chi_1,\dots, \chi_n)= P_{2k'}^*$, while 
$\mb(\chi_1, \dots, \chi_n)(\lambda)= Q_{2l'-1}^*(\lambda)$.  
But $P_{2k'}^*= \slps{2k'}$ and $Q_{2l'-1}^*(\lambda)= \slqs{2l'-1}$, 
by \eqref{el.e5}. Thus 
$$
\ma(\chi_1,\dots, \chi_n)= \slps{2k'}, \text{ and }  
\mb(\chi_1, \dots, \chi_n)(\lambda)= \slqs{2l'-1}.
$$
Thus they satisfy (\ref{e.AB}c) for the $\lambda$-groups.

The fact that $\slc{i}$ is the $\lambda$-Clifford correspondent of 
$\chi_i$, whenever $i=1, \dots, h$,  implies 
$$
G(\chi_1, \dots, \chi_h, \lambda)=G_{\lambda}(\slc{1}, \dots, \slc{h}) 
\leq G(\chi_1, \dots, \chi_h),
$$
for all $h=1, \dots, n$.
As $\ma(\chi_1)$ fixes $\lambda$, the group $\ma(\chi_1, , \dots, \chi_h)$ is 
a subgroup of the first group in the above list.
 It is also a $\pi$-Hall subgroup of 
$G(\chi_1,\dots, \chi_h)$, by (\ref{e.AB}b). 
Hence it is a $\pi$-Hall subgroup of $G_{\lambda}(\slc{1}, \dots, \slc{h})$.
Furthermore, 
$$G(\chi_1, \dots, \chi_h) = \ma(\chi_1,\dots, \chi_h) \cdot
\mb(\chi_1,\dots, \chi_h), 
$$
by (\ref{e.AB}b). 
Hence $G_{\lambda}(\slc{1}, \dots , \slc{h})= 
G(\chi_1, \dots, \chi_h , \lambda)= \ma(\chi_1,\dots, \chi_h) \cdot
\mb(\chi_1, \dots, \chi_h, \lambda)$. 
So $\ma(\chi_1, \dots, \chi_h) $  and $\mb(\chi_1, \dots, \chi_h, \lambda)$ 
form a Hall system for $G_{\lambda}(\slc{1}, \dots, \slc{h})$, for all $h=1,\dots,n$. 
This completes the proof of Theorem \ref{el.tAB}.
\end{proof}

From now until  the end of the section we 
 restrict our attention to the smaller system  \eqref{el.1}. 
It is clear that the subset 
$$\{  \slq{2i-1}, \slp{2r} | \slb{2i-1} , \sla{2r} \}_{i = 1, r=0}^{l, k},
$$
of \eqref{el.e4b}, is a triangular set for 
the normal series $G_0 \unlhd \slg{1} \unlhd \dots \unlhd 
 \slg{m} \unlhd G_{\lambda}$, 
  corresponding  to the tower $\{\slc{i} \}_{i=0}^m$.
Hence Theorem \ref{D} implies that the above triangular set 
 satisfies  \eqref{DD1} 
and \eqref{DD2} for all $r=0, \dots, k$ and all $i=1,\dots, l$, respectively, since 
\eqref{el.e4b} satisfies them.
The groups in question now are $\qw(\beta_{2k-1, 2k})$ and $\hap(\alpha_{2l-2,2l-1})$
along with their corresponding groups $\slqw(\slb{2k-1,2k})$ and 
$\hap_{\lambda}(\sla{2l-2,2l-2})$,  in the $\lambda$-situation.
Their relation is described in the next two theorems.
\begin{theorem}\mylabel{el.t2}
Assume that $\{ \ma_{\lambda}, \mb_{\lambda}\}$ is a  Hall system for
$G_{\lambda}$ that is derived from $\{\ma, \mb \}$ and satisfies  the conditions in Theorem 
\ref{el.tAB}. Assume further   that,  for every $m=1, \dots, n$, the group 
 $\qw$ is picked to satisfythe conditions in  Theorem \ref{sy.H} for the smaller 
system \eqref{el.1}, while  the group $\slqw$ is picked to satisfy 
 similar conditions  for the $\lambda$-groups. Then  
\begin{equation}\mylabel{el.t2e}
\slqw(\slb{2k-1,2k}) = \qw(\beta_{2k-1,2k}, \lambda).
\end{equation}
So $\maq_{\lambda} = \slqw(\slb{2k-1, 2k}) \cdot \slqs{2l-1} 
\leq  \maq(\lambda)$, where  $\maq= \qw(\beta_{2k-1,2k}) \cdot Q_{2l-1}^*$.
\end{theorem}

\begin{proof}
Assume that  $\slqw$ satisfies the conditions in
  Theorem \ref{sy.H}  for the $\lambda$-situation. 
Of course it satisfies  the equivalent   of  the conditions in Theorem \ref{hat:p1}
 for the $\lambda$-groups.
Furthermore, 
\begin{align*}
\slqw(\slb{2k-1,2k}) &=
N(\slps{2k} \tin \mb_{\lambda}(\slc{1}, \dots, \slc{2k})) 
&\text{ by Theorem \ref{sy.H} }\\
&=N(P_{2k}^* \tin \mb(\chi_{1}, \dots, \chi_{2k}, \lambda)) 
&\text{ by \eqref{104a} and \eqref{el.eAB} } \\
&=N(P_{2k}^* \tin \mb(\chi_1, \dots, \chi_{2k})) (\lambda)
& \\
&=\qw(\beta_{2k-1,2k})(\lambda).
&\text{ by Theorem \ref{sy.H} }
\end{align*}
This proves the first part of the theorem. 
The last part follows from the first and 
\eqref{104b}.
\end{proof}

Similarly to the above Theorem \ref{el.t2} we have 
\begin{theorem}\mylabel{el.t22}
Assume that $\{ \ma_{\lambda}, \mb_{\lambda}\}$ is a  Hall system for
$G_{\lambda}$ that is derived from $\{\ma, \mb \}$ and satisfies the conditions in 
Theorem  \ref{el.tAB}. Assume further  that, for every $m=1, \dots, n$, the group 
 $\hap$ is picked to satisfy the conditions in  Theorem \ref{sy.t1} for the smaller 
system \eqref{el.1}, while  the group $\hap_{\lambda}$ is picked to satisfy 
the similar conditions  for the $\lambda$-groups. Then  
 Then  
\begin{equation}\mylabel{el.t22e}
\hap_{\lambda}(\sla{2l-2,2l-1})= \hap(\alpha_{2l-2,2l-1}).
\end{equation}
So $\map_{\lambda} = \hap_{\lambda}(\sla{2l-2,2l-1}) \cdot \slps{2k}=  \map
=\hap(\alpha_{2l-1,2l-1}) \cdot P_{2k}^*$.
\end{theorem}

\begin{proof}
Assume that  $\slqw$ satisfies  the conditions in 
Theorem \ref{sy.t1}  for the $\lambda$-situation. 
Of course it satisfies  the equivalent of the conditions in Theorem \ref{hat:p1}
 for the $\lambda$-groups.
Furthermore, 
\begin{align*}
\hap_{\lambda}(\sla{2l-2,2l-1})
&=N(\slqs{2l-1} \tin \ma_{\lambda}(\slc{1}, \dots, \slc{2l-1})) 
&\text{ by Theorem \ref{sy.t1} }\\
&=N(Q_{2l-1}^*(\lambda) \tin \ma(\chi_{1}, \dots, \chi_{2l-1})) 
&\text{ by \eqref{104b} and \eqref{el.eAB} } 
\end{align*}
Clearly $Q_1  \leq Q_{2l-1}^* \leq G(\beta_1)$. Thus Lemma \ref{e.l1} implies 
$$
Q_{2l-1}^* = Q_{2l-1}^* (\lambda) \cdot Q_1.
$$
As $\ma(\chi_1, \dots, \chi_{2l-1}) $ is contained in $ \ma(\chi_1)$, it fixes $\lambda$
since $\ma(\chi_1)$ fixes $\lambda$.
It also normalizes $Q_1 \unlhd G$.
Therefore
$$
N(Q_{2l-1}^*(\lambda) \tin \ma(\chi_{1}, \dots, \chi_{2l-1})) =
N(Q_{2l-1}^* \tin \ma(\chi_{1}, \dots, \chi_{2l-1})). 
$$
According to Theorem \ref{sy.t1}, the latter group  is 
$\hap(\alpha_{2l-2,2l-1})$.
 Hence
$$
\hap_{\lambda}(\sla{2l-2,2l-1})= N(Q_{2l-1}^*(\lambda) \tin \ma(\chi_{1}, 
\dots, \chi_{2l-1})) = \hap(\alpha_{2l-2,2l-1}). 
$$
So the first part of  Theorem \ref{el.t22} follows. This,
 along with \eqref{104a} implies the rest of the theorem.
\end{proof}

The fact that $\qw \leq G'=G(\alpha_{2k}^*)$ normalizes $P_{2k}^*$, along with
 \eqref{el.e3.4}, obviously implies  
\begin{equation}\mylabel{el.e30}
\qw(\beta_{2k-1,2k}, \lambda)= \qw(\beta_{2k-1,2k}, \lambda_{2k}).
\end{equation}

We define 
\begin{equation}\mylabel{el.e31}
I:=\text{the image of $\qw(\beta_{2k-1,2k})$ in } \Aut(P_{2k}^*).
\end{equation}
Obviously, the group $I$ is well defined, as $\qw \leq G(\alpha_{2k}^*)$
normalizes $P_{2k}^*$.
Then as  an easy corollary of Theorem \ref{el.t2}  we get 
\begin{corollary} \mylabel{el.co2}
The groups $\slqw(\slb{2k-1,2k})$ and $\qw(\beta_{2k-1,2k})$,
chosen in Theorem \ref{el.t2}, 
 have the 
same image in $\Aut(P_{2k}^*)$. 
In particular,
\begin{equation}\mylabel{el.e32}
I = I_{\lambda}, 
\end{equation}
where $I_{\lambda}$ is  the image of $\slqw(\slb{2k-1,2k})$
 in $\Aut(\slps{2k})$.    

\end{corollary}

\begin{proof}
The group $\qw(\beta_{2k-1,2k})$ is a subgroup of $G'(\beta_{2k-1,2k}) = 
G(\alpha_{2k}^*,\beta_{2k-1,2k})$. By Proposition \ref{pqpropo}, the
 character $\beta_{1,2k}$ is the unique character of $Q_{1, 2k}$ 
lying under $\beta_{2k-1,2k}$. 
Thus $G'(\beta_{2k-1,2k})$ fixes $\beta_{1,2k}$. 
Hence 
$$
\qw(\beta_{2k-1,2k}) \leq G'(\beta_{2k-1,2k}) 
\leq N(P_{2k}^* \tin G(\beta_{1,2k})).
$$ 
Furthermore, according to \eqref{pqd40a} and the Definition \ref{pq40a}, 
the group $Q_{1,2k}$ is a subgroup of $\qw$. 
It also fixes $\beta_{2k-1,2k} \in \Irr(Q_{2k-1, 2k})$,
 as $Q_{1,2k} \leq Q_{2k-1,2k}$. Thus $Q_{1,2k}$ is a subgroup 
of $\qw(\beta_{2k-1,2k})$. 
In conclusion,
$$Q_{1,2k} \leq \qw(\beta_{2k-1,2k}) \leq N(P_{2k}^* \tin G(\beta_{1,2k})).
$$
Therefore Lemma \ref{e.l1} can be applied with  $\qw(\beta_{2k-1,2k})$ in the 
place of $S$.  So
$$
\qw(\beta_{2k-1,2k}) = \qw(\beta_{2k-1,2k}, \lambda) \cdot Q_{1,2k}.
$$
This, along with  Theorem \ref{el.t2} implies 
$$
\qw(\beta_{2k-1,2k}) = \slqw(\slb{2k-1,2k}) \cdot Q_{1,2k}.
$$
The fact that $Q_{1,2k}= C(P_{2k}^* \tin Q_1)$ centralizes $P_{2k}^*$ 
implies the first part of the  corollary immediately.
This,  along with \eqref{104a} and the definitions of
 $I$ and $I_{\lambda}$, 
implies   equation \eqref{el.e32}.  
\end{proof}

Theorem \ref{el.t22} easily implies
\begin{corollary}\mylabel{el.co22}
The groups $\hap(\alpha_{2l-2,2l-1})$ and $\hap_{\lambda}(\sla{2l-2,2l-1})$,
chosen in Theorem \ref{el.t22},
have the same images in $\Aut(Q_{2l-1}^*)$ and the same images in 
 $\Aut(\slqs{2l-1})$. 
\end{corollary}

We finish this section by  checking a special case of extendibility
in  the $\lambda$-situation.
The following is a well known result.
\begin{theorem}\mylabel{prel.t5}
Assume that $G$ is a finite group and that $S \leq H $ are subgroups of
 $G$ with  $S$ normal in $G$. Assume further that 
$\theta \in \Irr(H)$ lies above $\lambda \in \Irr(S)$. Let
 $\theta_{\lambda} \in \Irr(H(\lambda))$ 
denote the unique $\lambda$-Clifford correspondent 
of $\theta$.

If $\theta$ extends to its stabilizer $G(\theta)$ in $G$, then 
$\theta_{\lambda}$ also  extends to 
$ G(\theta, \lambda)$.
\end{theorem}

\begin{proof}
A straight forward application of Clifford Theory implies that 
$G(\theta, \lambda) \leq  G(\theta_{\lambda})$.
 Furthermore, as $G(\theta)$ fixes $\theta$ it permutes among themselves the 
members of the  $H$-conjugacy 
class of characters in $\Irr(S)$ lying under $\theta$. Since $\lambda \in \Irr(S)$
lies under $\theta$ we get 
\begin{subequations}\mylabel{prel.t5.1}
\begin{equation}
G(\theta) = H \cdot G(\theta, \lambda) \leq  H \cdot G(\theta_{\lambda}).
\end{equation}
In addition, 
\begin{equation}
G(\theta, \lambda) \cap H = H( \lambda).
\end{equation}
\end{subequations}

Let $\theta^i \in \Irr(G(\theta))$ be an extension of $\theta$ to $G(\theta)$. 
Then  $\theta^i$ lies above $\lambda $. Let $\Psi \in \Irr(G(\theta, \lambda))$ 
denote the unique $\lambda$-Clifford correspondent of $\theta^i \in \Irr(G(\theta))$.
So $\Psi $ lies above $\lambda$ and induces $\theta^i$.
 Therefore, 
$$(\Psi^{G(\theta)})|_H = \theta^i |_H = \theta.$$
Mackey's Theorem, along with \eqref{prel.t5.1},  implies that 
$$(\Psi^{G(\theta)})|_H = (\Psi|_{H(\lambda)})^H.$$
Hence $ (\Psi|_{H(\lambda)})^H= \theta$ is an irreducible character of
$H$. So the restriction $\Psi|_{H(\lambda)}$  
is an irreducible character of $H(\lambda)$ that induces $\theta$ 
and lies above $\lambda$ (as $\Psi $ lies above $\lambda$).
We conclude that $\Psi|_{H(\lambda)} $ is the $\lambda$-Clifford correspondent of 
$\theta$.  Hence $\Psi|_{H(\lambda)} = \theta_{\lambda}$. 
Thus $\Psi$ is an extension of $\theta_{\lambda}$  to $G(\theta, \lambda)$, 
and Theorem \ref{prel.t5} follows.
\end{proof}

As a consequence of Theorem \ref{prel.t5} we get
\begin{theorem}\mylabel{el.t3}
Assume that $\beta_{2k-1,2k} \in \Irr(Q_{2k-1,2k})$ extends to 
$\qw(\beta_{2k-1,2k})$. Let $\slqw$ be a  $\pi'$-Hall subgroup of
 $G_{\lambda}'=G(\lambda, \sla{2k}^*)$ that satisfies the conditions in  Theorem \ref{el.t2}.
Then the character $\slb{2k-1,2k}$ extends to $\slqw(\slb{2k-1,2k})$.
\end{theorem}

\begin{proof}
According to \eqref{pq40g} we get that 
$$
Q_{2k-1,2k}= \qw_{2k-1,2k}(\beta_{2k-1,2k})=\qw(\beta_{2k-1,2k}) \cap G_{2k-1}
\unlhd \qw(\beta_{2k-1,2k}).
$$
In view of Remark \ref{el.rem1}, this  implies  
$$E_{2k} \unlhd \slq{2k-1,2k}=Q_{2k-1,2k}(\lambda_{2k})
 \leq Q_{2k-1,2k} \unlhd \qw(\beta_{2k-1,2k}).
$$
As $\slb{2k-1,2k}\in \Irr(\slq{2k-1,2k})$ is the $\lambda_{2k}$-Clifford 
correspondent of $\beta_{2k-1,2k} \in \Irr(Q_{2k-1,2k})$, Theorem 
\ref{prel.t5} implies that $\slb{2k-1,2k}$ extends to $\qw(\beta_{2k-1,2k}, 
\lambda_{2k})$, when $\beta_{2k-1,2k}$ extends to $\qw(\beta_{2k-1,2k})$.
But
$$
\qw(\beta_{2k-1,2k},\lambda_{2k})= \qw(\beta_{2k-1,2k}, \lambda)=
\slqw(\slb{2k-1,2k}), 
$$
by \eqref{el.e30} and Theorem \ref{el.t2}.

This completes the proof of the theorem.
\end{proof} 
Furthermore, Remark \ref{el.rem11} easily implies 
\begin{theorem}\mylabel{el.t33}
Assume that $\alpha_{2l-2,2l-1}  \in \Irr(P_{2l-2,2l-1})$ extends to 
$\hap(\alpha_{2l-2,2l-1})$. Let $\hap_{\lambda}$ be the 
$\pi$-Hall subgroup of
 $G_{\lambda}(  \beta_{2l-1}^*)$ that satisfies the conditions in Theorem \ref{el.t22}.
Then the character $\sla{2l-2, 2l-1}$ extends to 
$\hap_{\lambda}(\sla{2l-2,2l-1})$.
\end{theorem}

%%% Local Variables: 
%%% mode: latex
%%% TeX-master: t
%%% End: 

\section{Normal $\pi$-subgroups inside  $P_2$}
\mylabel{elp}
Assume now that we are in a situation where the 
fixed  normal series \eqref{e0.2}, i.e.,  
$1=G_0 \unlhd G_1\unlhd \dots \unlhd G_{n} \unlhd G$, satisfies \eqref{sh.e1}. 
So teh following two onditions hold
\begin{gather}\mylabel{lp.1}
G_2= G_{2, \pi} \times G_{2,\pi'}, \\
G \text{ fixes $\chi_1$. } 
\end{gather}
 We assume fixed the character tower \eqref{e0.3} and 
its corresponding  triangular set \eqref{e0.4}. 
We also fix the Hall system $\{\ma, \mb \}$ that satisfies \eqref{e.AB}.

The additional hypothesis on the group $G_2$  give  
 more specific information on $\chi_2$ (see \eqref{sh.e2}). Thus we have 
\begin{subequations}\mylabel{117}
\begin{equation}\mylabel{117a}
 G_2 = P_2 \times Q_1 = P_2\times G_1,
\end{equation}
\begin{equation}\mylabel{117b}
\chi_2 = \alpha_2 \times \beta_1.
\end{equation}
\end{subequations}
Furthermore
\begin{equation}\mylabel{el.ch}
G(\chi_2) = G(\alpha_2). 
\end{equation}

We also fix an integer $m= 2, \dots, n$ and we consider the normal series 
\begin{subequations}\mylabel{elp.e1}
\begin{equation}
1=G_0\unlhd G_1 \unlhd \dots \unlhd G_{m} \unlhd  G.
\end{equation}
We also fix the subtower
\begin{equation}\mylabel{elp:e1a}
\{1=\chi_0, \chi_1, \dots, \chi_m \}
\end{equation}
of \eqref{e0.3} 
and the subset 
\begin{equation}\mylabel{elp.e1b}
\{Q_{2i-1}, P_{2r} | \beta_{2i-1}, \alpha_{2r} \}_{i=1, r=0} ^{l, \, \,  k}
\end{equation}
\end{subequations}
of \eqref{e0.4} . The subtower \eqref{elp:e1a}  is a character tower of 
 (\ref{elp.e1}a), and the subset  \eqref{elp.e1b} is 
 is a representative 
 of the conjugacy class of triangular sets 
that corresponds uniquely to \eqref{elp:e1a}. 
We also assume known the groups $\qw$ and $\hap$ 
satsfying \eqref{pq30} and  \eqref{sy.e}  respectively
for the system \eqref{elp.e1}. Also known are assume the groups 
 $\map$ and $\maq$, as these were defined in 
  \eqref{e0.d}, i.e., 
\begin{equation}\mylabel{elp.e2}
\map = \hap(\alpha_{2l-2, 2l-1}) \cdot P_{2k}^* \text{ and } 
\maq = \qw(\beta_{2k-1,2k}) \cdot Q_{2l-1}^*.
\end{equation}
So $\map \leq \ma(\chi_1,\dots,\chi_{2l-1})$ and $\maq \leq
\mb(\chi_1,\dots,\chi_{2k})$, by \eqref{e0.1}.

As the series \eqref{e0.2} satisfies \eqref{lp.1}, we can apply 
all the results of Section \ref{sh}.
So  according to Theorem \ref{sh.t1}, we can drop the group $Q_1= G_1$
 (i.e. shift the above  series  by one), without any loss.
Thus the original system \eqref{e0.1} reduces to 
\begin{subequations}\mylabel{llp.e0}
\begin{gather}\mylabel{llp.e1}
1 \unlhd P_2 \unlhd G_3 \unlhd \dots  \unlhd  G_n =  G,\notag \\
1, \alpha_2, \chi_3, \dots, \chi_n,\\
\{ Q_1^s=1, Q_{2i-1}, P_{2r} |\beta_1^s=1, \beta_{2i-1}, \alpha_{2r} \}_{i=2, r=0}^{l',k'}\notag
\end{gather}
Similarly the  subsystem \eqref{elp.e1} reduces to 
\begin{gather}\mylabel{llp.e2}
1 \unlhd P_2 \unlhd G_3 \unlhd \dots  \unlhd  G_m \unlhd G,\notag \\
1, \alpha_2, \chi_3, \dots, \chi_m,\\
\{ Q_1^s, Q_{2i-1}, P_{2r} |\beta_1^s, \beta_{2i-1}, \alpha_{2r} \}_{i=2, r=0}^{l,k}\notag
\end{gather}
\end{subequations}
They both satisfy the conditions in  Theorem \ref{sh.t1}. Therefore the Hall system
 $\{ \ma, \mb \}$ remains the same as well as the groups $P_{2k}^*, 
Q_{2l-1}^*, \qw, \qw(\beta_{2k-1, 2k}),  \hap$ and $\hap(\alpha_{2l-1,2l-1})$.
Note also that the shifted systems satisfy 
\begin{gather}\mylabel{llp.e3}
G_1^s= Q_1^s = 1 \text{ and }  \chi_1^s=  \beta_1^s =1  \\
G_2^s = P_2 = P_2^* \text{ and }  \chi_2^s= \alpha_2= \alpha_2^*. 
\end{gather}

As with the previous section, the goal  is to understand the behavior 
of the  above systems when we apply Clifford's theory to a 
normal subgroup of $G$. This time the normal subgroup is contained 
in $P_2$, and not $Q_1$ as was the case with Section \ref{elm1}.
The results we obtain are similar to those of the previous section, and their
 proofs are identical, (the only modification being, whenever necessary, 
the interchanged role of $\pi$ and $\pi'$).
As with the group $S$ and the character $\zeta$  earlier, we fix 
 (until the end of this  section)  a subgroup $R$ of $G_2$
and a character $\eta$ of $R$  
 that satisfy 
\begin{subequations}\mylabel{eps.1}
\begin{gather} 
R \unlhd G  \text{ and } R\leq P_2,\\
\eta \in \Irr(R) \text{ is $G$-invariant  and lies under $\alpha_2$. }
\end{gather}
We also assume that $M$ is a normal subgroup of $G$ with 
\begin{equation}\mylabel{eps.2}
R \leq M \leq  P_2.
\end{equation}
\end{subequations}
The role of $\ma(\beta_1)$ is played  here by the group
 $\mb(\alpha_2)$. So, as 
 in the case of $\lambda$ and $E$ in the previous 
section, we have
\begin{lemma}\mylabel{elmp:l1}
There is an irreducible character $\mu \in \Irr(M)$  such that $\mu$ is
 $\mb(\alpha_2)$-invariant and lies under $\alpha_2$ 
and under  $\chi_i$, for all
$i=3,\dots,n$.  Any such  $\mu$ lies above $\eta \in \Irr(R)$.
\end{lemma}
\begin{proof}
The proof is similar to that of Lemma \ref{elm:l1}.
If $\mu_1$ is any irreducible character of $M$ lying under $\alpha_2$,
then $\mb(\alpha_2)$ permutes 
among themselves the $P_2$-conjugates of $\mu_1$,
as it  normalizes  $P_2 \unlhd G$ and fixes $\alpha_2$.
So Glauberman's Lemma  implies that $\mb(\alpha_2)$ fixes at least 
one $P_2$-conjugate of $\mu_1$.
The lemma follows if we observe that any character of $M$ that lies under 
$\alpha_2$ also lies 
 under $\chi_i$, for all $i=3, \dots,n$, while any character of 
$M$ lying under $\alpha_2$ lies necessarily above 
$\eta \in \Irr(R)$, 
as $\eta$ is $G$-invariant.
\end{proof}

Note also that 
\begin{remark}\mylabel{elmp.r0}
Assume that  $\mu_1 \in \Lin(M)$ is a linear character of $M$ lying 
under $\alpha_2$. Then there exists 
a $P_2$-conjugate $\mu \in \Lin(M)$ of $\mu_1$, such that 
$\mu$ is $\mb(\alpha_2)$-invariant, and   
lies under $\alpha_2$ and above $\eta$.
\end{remark}

\begin{remark}\mylabel{elp.r00}
  The $\pi'$-group $Q_{2l'-1}^*=1 \cdot  Q_3 \cdot Q_5   \cdots Q_{2l'-1}$ 
  fixes $\alpha_2$, as every one of its factors does,  by  \eqref{x5}.
  Thus $Q_{2l'-1}^* \leq \mb(\alpha_2)$.  So  $Q_{2l-1}^*$  fixes $\mu$.
\end{remark}

\begin{remark}\mylabel{elp.r000}
For any $m \geq 2$, the $\pi'$-group $\maq$ is a subgroup of
$\mb(\chi_1,\chi_2, \dots, \chi_{2k})$, by (\ref{e0.5}d). Thus
$\maq \leq \mb(\chi_1,\chi_2) = \mb(\chi_2) = \mb(\alpha_2)$.
 So $\maq$  fixes $\mu$.
\end{remark}

From now on we fix a character $\mu \in \Irr(M)$ having
 all the properties in Lemma \ref{elmp:l1}.

Since the $\pi'$-group $Q_{2i-1}^*$, as a subgroup
 of $Q_{2l'-1}^* \leq \mb(\alpha_2)$, fixes the character $\mu$,
we can define $M_{2i-1}$ and $\mu_{2i-1}$ by 
\begin{subequations}\mylabel{elp.e3}
\begin{align}
M_{2i-1} &= C(Q_{2i-1}^* \tin M )  \text{ and } \\
\mu_{2i-1} \in \Irr(M_{2i-1}) &\text{ is  the $Q_{2i-1}^*$-Glauberman correspondent of } \mu\in \Irr(M),
\end{align}
for all $i=2,\dots,l'$. 
We also define $M_1$ and $\mu_1$ as 
\begin{equation}\mylabel{elp.e3.2}
M_1= M    \text{ and }  \mu_1=\mu.
\end{equation}
\end{subequations}
Furthermore, 
\begin{subequations}\mylabel{elp.e4}
\begin{equation}\mylabel{elp.e4.1}
M_{2i-1} = P_{2,2i-1} \cap M= C(Q_{2i-1}^* \tin P_2) \cap M,
\end{equation}
for all $i=2,\dots, l'$, while 
\begin{gather}\mylabel{elp.e4.2}
\mu_1 \text{ lies under } \alpha_2 \in \Irr(P_2), \\ 
\mu_{2i-1} \text{ lies under  } \alpha_{2, 2i-1}\in \Irr(P_{2, 2i-1}),
\end{gather}
\end{subequations}
 whenever $2\leq i \leq l'$,
as $\alpha_{2,2i-1}$ is the $Q_{2i-1}^*$-Glauberman correspondent of $\alpha_2$, and $\alpha_2$ lies 
over $\mu$. In view of \eqref{elp.e3}  we have 
\begin{equation}\mylabel{elp.e4.3}
N(\mu)= N(\mu_{2i-1}),
\end{equation}
for all groups $N$ with  $N \leq N(Q_{2i-1}^* \tin G)$, and all 
 $i=2, \dots, l'$. 

As in the previous section, we define 
\begin{subequations}\mylabel{elp.e5}
\begin{align}
G_{\mu} &:= G(\mu), \mylabel{elp.e5a}\\
\smg{i} &:= G_i(\mu) = G_{\mu} \cap G_i, \mylabel{elp.e5b}  \\
 P_{2, \mu} &:= P_2(\mu),
\end{align}
whenever $3\leq i \leq n$. This way we can 
 form the series
\begin{equation}\mylabel{elmp:e1}
1 \unlhd \smp{2} \unlhd \smg{3} \unlhd \dots \unlhd \smg{n}= G_{\mu},
\end{equation}
of normal subgroups of  $G_{\mu} = G(\mu)$.
\end{subequations}

We also write   
\begin{subequations}\mylabel{elmp.e3}
\begin{equation}\mylabel{elmp:e3}
G_{1, \mu}:= Q_1^s = 1,
\end{equation}
and
\begin{equation}\mylabel{elmp:e4}
G_{2, \mu} :=  P_{2}(\mu)= P_{2, \mu}.
\end{equation}
\end{subequations}

Furthermore, we  can apply Clifford's Theorem to the
groups $G_i$
and the characters $\chi_i$,  for all $i=3, \dots,n$.
Thus there exist unique irreducible characters $\smc{i} \in \Irr(\smg{i})$
lying above $\mu$ and inducing $\chi_i$
\begin{subequations}\mylabel{118}
\begin{equation}\mylabel{118a}
\smc{i} \in \Irr(\smg{i} | \mu) \text{ and } (\smc{i})^{G_i}= \chi_i,
\end{equation}
for all $i=3,\dots,n$.
 We complete this list by setting
\begin{equation}\mylabel{118b}
 \smc{1}=\beta_1^s=  1
\end{equation}
We also write
\begin{equation}\mylabel{118c}
\smc{2} = \alpha_{2, \mu} \in \Irr(P_2(\mu)) 
\end{equation}
for the  $\mu$-Clifford correspondent 
of $\alpha_2$, 
 i.e.,  $\smc{2}$ lies above $\mu$ and induces $\alpha_2 \in \Irr(P_2)$. 
\end{subequations}
Clearly $\smc{2}$ lies above $\smc{1}=1$. Furthermore, 
as the $\chi_i$ lie  above each other,  the same holds for the characters 
$\smc{i}$, by Clifford's theory, for all $i=2,\dots,n$.
Therefore, we have formed a character tower 
\begin{subequations}\mylabel{elp.e6}
\begin{equation}\mylabel{elp.e6a}
\{1, \smc{2}=\sma{2}, \smc{3}, \dots, \smc{n}\}
\end{equation}
for the series \eqref{elmp:e1}. Hence Theorem \ref{cor:t}, applied  to the tower \eqref{elp.e6a}, implies the 
existence of a unique $G_{\mu}$-conjugacy class of triangular sets of \eqref{elmp:e1} that corresponds 
to the tower \eqref{elp.e6a}. Let 
\begin{equation}\mylabel{elp.e6b}
\{\smq{2i-1}, \smp{2r} |\smb{2i-1}, \sma{2r} \}_{i=1, r=0}^{l', \, \, k'}
 \end{equation}
be a representative of this class. 
\end{subequations}

 All the groups,  the characters and their properties  that were described 
in Chapter \ref{pq} are valid for the $\mu$-situation. We follow the same notation as in
the previous section, with 
   an extra $\mu$ in the place of  the $\lambda$ there. The goal 
is the same as in the previous section, that is, to compare the triangular
 set in  \eqref{llp.e1} 
with that in \eqref{elp.e6b}. 
As many of the results here have proofs analogous  to those
 in Section \ref{elm1}, 
we give them briefly or skip them. 

The first steps in that  direction  follow from 
\eqref{elmp.e3} and \eqref{118}. In particular, these relations  clearly imply 
\begin{subequations}\mylabel{elp.e7}
\begin{align}
\smq{1}= \smg{1} &= G_1^s = Q_1^s= 1, \mylabel{elp.e7a}\\
\smb{1} &= \beta_1^s= 1, \mylabel{elp.e7b}\\
\smp{2}&= P_2(\mu),  \mylabel{elp.e7c}\\
\sma{2}  \in \Irr(\smp{2})
&\text{ is the $\mu$-Clifford correspondent of $\alpha_2 \in \Irr(P_2)$}. \mylabel{elp.e7d}
\end{align}
\end{subequations}

As in \eqref{pq1} and \eqref{pqinf}, we get the groups  $\smgs{i}$  and $G^*_{\mu}$, 
defined as 
\begin{align}\mylabel{119}
\smgs{0} &:= \smg{0} = 1, \notag \\
\smgs{i} &:= \smg{i}(\smc{1},\dots \smc{i-1}) = \smg{i}(\smc{1},\dots ,\smc{n}), \\
G^*_{\mu} &:= G_{\mu}(\smc{1}, \dots, \smc{n}), \notag
\end{align}
for all $i=1,\dots, n$. Clearly we have 
\begin{subequations}\mylabel{119'}
\begin{align}
\smg{i} &= G^*_{\mu} \cap G_i , \\
 \smgs{1}&=\smg{1}= Q_{1, \mu}= 1, \\
\smgs{2}&=\smg{2}(\chi_1)= \smp{2}(1) = \smp{2},
\end{align}
\end{subequations}
whenever $3\leq i \leq n$.

As with the $\lambda$-groups and  Lemma \ref{elm.l2}, the following holds:
\begin{lemma}\mylabel{elmp:l2}
For any $i=3, \dots,n $ we have
\begin{subequations}\mylabel{elp.e8}
\begin{align}
G^*_{\mu}&= G^*(\mu) \mylabel{elp.e8a} \\
\smgs{i} &= G^*_i(\mu). \mylabel{elp.e8b}
\end{align}
\end{subequations}
\end{lemma}
 
\begin{proof} 
The proof is similar to that of  Lemma \ref{elm.l2}.
 So we  omit it.
\end{proof}

We can now prove
\begin{proposition}\mylabel{elmp:p1}
For every $r=0,1,\dots , k'$ and $i=2,\dots , l'$ we have that 
\begin{subequations}\mylabel{elp.e9}
\begin{align}
\mylabel{elp.e9a}  P^*_{2r}(\mu) &\in \Hall_{\pi}(\smgs{2r}),     \\
\mylabel{elp.e9b}
 Q^*_{2i-1}&\in \Hall_{\pi'}(\smgs{2i-1}). 
\end{align}
\end{subequations}
Therefore  the triangular set \eqref{elp.e6b} can be chosen among the sets in its $G_{\mu}$-conjugacy 
class so that it  satisfies 
\begin{subequations}\mylabel{121}
\begin{align}
P_{2r}^*(\mu) &=\smps{2r}, \mylabel{121a}\\
Q_{2i-1}^* &=\smqs{2i-1}, \mylabel{121b}
\end{align}
\end{subequations}
whenever $0\leq r \leq k'$ and $2\leq i\leq l'$.
\end{proposition}
\begin{proof}
The cases $P_0^*$ and $ P_2^*$  are easy (see \eqref{elp.e7c} ).
 The rest of the proof is similar to that of  
 Proposition \ref{elm:p1}, with the roles of
$\smps{2i}$ and $\smqs{2i-1}$ interchanged.
\end{proof}

The analogue of Lemma \ref{e.l1} is 
\begin{lemma} \mylabel{elp.l1}
Assume that $P_2 \leq T \leq G(\alpha_2) $. Then $T = T(\mu) \cdot P_2$. 
Furthermore, if  $S$ satisfies 
  $P_{2,2i-1} \leq  S \leq N(Q_{2i-1}^* \tin G(\alpha_{2,2i-1}))$, for some $i=2,\dots,l'$, 
Then $S =S(\mu_{2i-1})\cdot P_{2, 2i-1} = S(\mu) \cdot P_{2, 2i-1}$.
\end{lemma} 
\begin{proof}
Similar to  the proof of Lemma \ref{elp.l1}.
We only  remark that the role of $P_2$ here 
is the same as that of $Q_1$ there,  as 
the group $P_2$ is a normal subgroup of $G$.
\end{proof}

To compare  the groups $\smp{2i}$ and $\smq{2i-1}$ 
with $P_{2i}$ and $Q_{2i-1}$ we have, as we would have guessed
\begin{theorem}\mylabel{P}
The set \eqref{elp.e6b}  chosen in Proposition \ref{elmp:p1}  satisfies 
\begin{subequations}\mylabel{P1}
\begin{align}
\smq{2i-1}&= Q_{2i-1}, \mylabel{P1a} \\
\smb{2i-1} &= \beta_{2i-1},\mylabel{P1b}
\end{align}
\end{subequations}
for all $i=2,\dots ,l'$, and 
\begin{subequations}\mylabel{P2}
\begin{align}
\smp{2r} &= P_{2r}(\mu)= P_{2r}(\mu_{2r-1}), \mylabel{P2a} \\ 
\sma{2r}\in \Irr(\smp{2r}) &\text{ is the $\mu_{2r-1}$-Clifford correspondent of  } 
\alpha_{2r}\in \Irr(P_{2r}),\mylabel{P2b}
\end{align}
\end{subequations}
for all $r = 1, \dots ,k'$.
Hence $\sma{2r}$ induces $\alpha_{2r} \in \Irr(P_{2r})$.
\end{theorem}  
\begin{proof}
The proof is similar to that of   Theorem \ref{D}. So we omit it.
We only remark that we need to interchange the role of the $\pi$-and
 $\pi'$-groups.
Note also that  the case $i=1$, that is omitted here,  was  already done, along with the case  $r=1$,  
in \eqref{elp.e7}.
\end{proof}

The analogue of \eqref{el.e7}, that would also appear as a step if we
 had  given the proof of the above 
theorem in detail,  is
\begin{remark}\mylabel{elp.r2}
For all $r=2,\dots, k'$, we get 
$$
P_{2r} = P_{2r}(\mu)\cdot P_{2,2r-1}= P_{2r}(\mu_{2r-1}) \cdot P_{2,2r-1}.
$$
\end{remark}

\begin{proof}
For all $r=2, \dots, k'$ we have that 
$$
P_{2, 2r-1} = P_2 \cap P_{2r} \leq P_{2r} \leq N(Q_{2r-1}^* \tin G(\alpha_{2, 2r-1})), 
$$
by \eqref{pq13b} and Proposition \ref{pqremark1'}.
Hence we can apply Lemma \ref{elp.l1}, with $P_{2r}$ in the place of  $S$.
So the remark follows.
\end{proof}

Furthermore, the analogue of Remark \ref{el.rem1}  for the $\pi$-groups is
\begin{remark}\mylabel{elp.rem1}
For all  $r, i$ with $1\leq r < i \leq l'$ we have 
$$
\smp{2r,2i-1}= P_{2r,2i-1}(\mu)= P_{2r,2i-1}(\mu_{2i-1}),
$$
while the character $\sma{2r,2i-1}$ is the $\mu_{2i-1}$-Clifford correspondent of 
$\alpha_{2r,2i-1}$.
\end{remark}

\begin{proof}
The proof is the same as that of Remark \ref{el.rem1}, so we only sketch it.
\begin{multline}
\smp{2r,2i-1} = C(\smq{2r+1}, \dots, \smq{2i-1} \tin \smp{2r})= 
C(Q_{2r+1} \cdots Q_{2i-1} \tin P_{2r}(\mu)) \\
= C(Q_{2r+1} \cdots Q_{2i-1} \tin P_{2r})(\mu) = P_{2r, 2i-1}(\mu)= P_{2r, 2i-1}(\mu_{2i-1}), 
\end{multline}
where the last equation follows from \eqref{elp.e4.3}, along with  the fact
 that $P_{2r,2i-1}$ normalizes $Q_{2i-1}^*$.

Furthermore, the fact that Clifford theory is compatible with Glauberman 
correspondence (Lemma 2.5 in \cite{wo1}), along with \eqref{P2b}, 
implies that the $Q_{2r+1}\cdots 
Q_{2i-1}= \smq{2r+1} \cdots \smq{2i-1}$-Glauberman correspondent 
$\alpha_{2r,2i-1}$ of $\alpha_{2r}$ is induced by the 
$\smq{2r+1}\cdots \smq{2i-1}$-Glauberman correspondent 
$\sma{2r, 2i-1}$ of $\sma{2r}$, and, in addition, 
 $\sma{2r,2i-1}$ lies above the $Q_{2r+1} \cdots Q_{2i-1}$-Glauberman 
correspondent $\mu_{2i-1}$ of $\mu_{2r-1}$.
Hence Remark \ref{elp.rem1} follows.
\end{proof}

Now we can prove two corollaries that follow from  Theorem \ref{P},
\begin{corollary}\mylabel{elp.co1}
For all $i, j$ with $2\leq i \leq j \leq k'$ we have that 
$$
\smq{2i-1,2j} = Q_{2i-1,2j}  \text{ and } \smb{2i-1,2j}= \beta_{2i-1,2j}.
$$
In addition, $\smq{1, 2j} 
= Q_{1,2j}^s =1$ and $\smb{1,2j}= \beta_{1,2j}^s=1$, 
whenever $1\leq j \leq k'$. 
\end{corollary}
 
\begin{proof}
For all $t, i$ with  $2\leq  i \leq t  \leq k'$, the group  $Q_{2i-1}$ centralizes 
 $P_{2, 2t-1}= C(Q_3, \dots, Q_{2t-1} \tin P_2)$.
  As $P_{2t}= P_{2t}(\mu)  \cdot P_{2, 2t-1} $, we conclude that 
\begin{equation}\mylabel{elp.e10}
C(P_{2t} \tin Q_{2i-1}) = C(P_{2t}(\mu) \tin Q_{2i-1}),
\end{equation} 
whenever $2\leq i \leq t \leq k'$.
So we get 
\begin{align*}
\smq{2i-1,2j} &=C(\smp{2i} \cdots \smp{2j} \tin \smq{2i-1}) &\text{ by \eqref{pq14b}, for the 
  $\mu$-case }\\
&=C(P_{2i}(\mu)  \cdots P_{2j}(\mu) \tin Q_{2i-1}) &\text{ by \eqref{P1a} and \eqref{P2a} }\\
&=C(P_{2i} \cdots P_{2j} \tin Q_{2i-1}) &\text{ by \eqref{elp.e10} }\\
&=Q_{2i-1,2j} &\text{ by \eqref{pq14b} }
\end{align*}
whenever  $2\leq i \leq j \leq k'$.

The character $\smb{2i-1,2j}$  is the
 $\smp{2i} \cdots \smp{2j}$-Glauberman correspondent of 
$\smb{2i-1}$, by Definition \ref{pq1413def2} for the $\mu$-case.  
Hence it is also the 
 $P_{2i}(\mu) \cdots P_{2j}(\mu)$-Glauberman correspondent of 
$\beta_{2i-1}$, according to   Theorem \ref{P}.  This, along with \eqref{elp.e10}, implies that 
$\smb{2i-1,2j}$  is the
 $P_{2i} \cdots P_{2j}$-Glauberman correspondent of 
$\beta_{2i-1}$. Thus   $\smb{2i-1,2j}$ equals $  \beta_{2i-1,2j}$, as the latter was also  defined as the 
  $P_{2i} \cdots  P_{2j}$-Glauberman correspondent of 
$\beta_{2i-1}$. 
This completes the proof of the first part of the corollary.
The rest holds trivially.  So Corollary \ref{elp.co1} holds.
\end{proof}

\begin{corollary}\mylabel{elp.co2}
For all $r=1,\dots, k'$, the character  $\sma{2r}^* \in \Irr(\smps{2r})$ is the
 $\mu$-Clifford correspondent 
of $\alpha_{2r}^* \in \Irr(P_{2r}^*)$.
\end{corollary}

\begin{proof}
The character $\sma{2r}^*$ is defined (see Definition \ref{p*d1}), as the $\smq{3}, \dots, 
\smq{2r-1}$-correspondent of $\sma{2r}$, whenever $1\leq r \leq k'$. 
Hence,  \eqref{P1a} implies that $\sma{2r}^*$ is the $Q_3, \dots, Q_{2i-1}$-correspondent of 
$\sma{2r}$. But  $\sma{2r}$ is the $\mu$-Clifford correspondent of
 $\alpha_{2r}$, and the groups 
$Q_3, \dots, Q_{2r-1}$ fix $\mu$. 
According to  Proposition \ref{dade:p1.6}  the $A$-correspondence is compatible with the Clifford correspondence.
Thus, taking  $A$ as $  Q_3,  Q_5, \dots, Q_{2r-1}$ in turn, 
  we conclude  that the $Q_3, \dots,
 Q_{2r-1}$-correspondent $\sma{2r}^*$
  of $\sma{2r} \in \Irr(P_{2r}(\mu))$ is the $\mu$-Clifford correspondent of 
the $Q_3, \dots, Q_{2r-1}$-correspondent $\alpha_{2r}^*$ of $\alpha_{2r} \in \Irr(P_{2r})$.
Hence Corollary \ref{elp.co2} follows.
\end{proof}

As far as the Hall system $\{ \ma, \mb \}$ is concerned, 
we have, similarly to Theorem \ref{el.tAB}, the following 
\begin{theorem}\mylabel{elp.tAB}
We can find new $\ma_{\mu} \in \Hall_{\pi}(G_{\mu})$ 
and $\mb_{\mu} \in \Hall_{\pi'}(G_{\mu})$ 
satisfying  the 
equivalent of \eqref{e.AB} for the $\mu$-groups, along with 
\begin{subequations}\mylabel{elp.eAB}
\begin{gather}
\ma_{\mu} (\smc{1}, \dots, \smc{h}) = \ma(\chi_1, \dots, \chi_h, \mu), \\
\mb_{\mu}(\smc{1} , \dots, \smc{h}) = \mb(\chi_1, \dots, \chi_h),
\end{gather}
\end{subequations}
for all $h=2,\dots,n$. Hence 
 \begin{gather*}
\ma_{\mu} (\smc{2}, \dots, \smc{h}) = \ma(\chi_2, \dots, \chi_h, \mu), \\
\mb_{\mu}(\smc{2} , \dots, \smc{h}) = \mb(\chi_2, \dots, \chi_h),
\end{gather*}
for all such $h$.
\end{theorem}

\begin{proof}
The   proof is similar to that of Theorem \ref{el.tAB}, 
with the roles of $\ma$ and $\mb$ interchanged.
Just observe, for the last part, that $\smc{1}=1$, while $\chi_1$
is $ G$-invariant (and thus $\ma$-and $\mb$-invariant).
\end{proof}

We restrict our attention to the smaller system \eqref{elp.e1}. The subset 
$$
\{\smq{2i-1}, \smp{2r} |\smb{2i-1}, \sma{2r} \}_{i=1, r=0}^{l, \, \, k}
$$
of \eqref{elp.e6b}, is clearly a triangular set of the normal series
$1=G_0 \unlhd \smg{1} \unlhd  \smg{2} \unlhd \dots \unlhd \smg{m} \unlhd G_{\mu}$,  and  the tower $\{ \smc{i} \}_{i=0}^m$. 
Of course, \eqref{elmp.e3} and  \eqref{elp.e7} imply
that $\smg{1}=1= \smq{1}$ and  $\smc{1}=1 = \smb{1}$, while 
$\smg{2}= \smp{2}$ and $\smc{2}=\sma{2}$.  
In view of Theorem \ref{P},  the above 
system can be chosen to satisfy \eqref{P1} and \eqref{P2}, for all $r=0, \dots, k$ and 
all  $i=1, \dots, l$.  As in the previous section, we can
chose the groups $\qw$ and $\hap$ along with their corresponding in the $\mu$-case
 groups $\smqw$ and  $\hap_{\mu}$ to satisfy theorems analogous to 
Theorems \ref{el.t2} and \ref{el.t22}, that is, 

\begin{theorem}\mylabel{elp.t2}
Assume that $\{ \ma_{\mu}, \mb_{\mu}\}$ is a  Hall system for
$G_{\mu}$ that is derived from $\{\ma, \mb \}$ and satisfies the conditions in 
  Theorem 
\ref{elp.tAB}. Assume further  that for every $m=1, \dots, n$, the group 
 $\qw$ is picked to satisfy the conditions in Theorem \ref{sy.H} 
 for the smaller 
system \eqref{elp.e1}, while  the group $\slqw$ is picked to satisfy 
the equivalence of teh conditions in  Theorem \ref{sy.H} 
for the $\mu$-groups. Then  
$$
\smqw(\smb{2k-1,2k}) = \qw(\beta_{2k-1,2k}).
$$
So $\maq_{\mu} = \smqw(\smb{2k-1, 2k}) \cdot \smqs{2l-1}  = \maq$.
\end{theorem}

\begin{proof}
Same as the proof of Theorem \ref{el.t22}, with the roles of 
$\hap$ and $\qw$ interchanged. Note that, when passing to the 
shifted system, 
the groups $\hap$ and $\qw$ remain the same by Theorem \ref{sh.t1}.
\end{proof}

We also have 
\begin{theorem}\mylabel{elp.t22}
Let $\{\ma_{\mu}, \mb_{\mu}\}$ be
 as above, and let $m=1, \dots, n$ be fixed.
If   $\hap$ is picked to satisfy the conditions in Theorem \ref{sy.t1}, 
for the smaller system 
\eqref{elp.e1}, while $\hap_{\mu}$ is picked to satisfy   the 
equivalent  of the conditions in   Theorem \ref{sy.t1} 
for the $\mu$-groups, then  
$$
\hap_{\mu}(\sma{2l-2,2l-1})= \hap(\alpha_{2l-2,2l-1}, \mu).
$$
So $\map_{\mu} = \hap_{\mu}(\sma{2l-1,2l-1}) \cdot \smps{2k}
\leq  \map(\mu)$.
\end{theorem}

\begin{proof}
See the proof of Theorem \ref{el.t2}.
\end{proof}

As a corollary of 
Theorem \ref{elp.t2}, we get  
\begin{corollary} \mylabel{elp.co3}
The groups $\smqw(\smb{2k-1,2k})$ and $\qw(\beta_{2k-1,2k})$
 have the  same image in $\Aut(P_{2k}^*)$ . They also have the same image in 
$\Aut(P_{2k, \mu}^*)$.
Thus  
\begin{equation}\mylabel{elp.e17}
I = \text{ the image of  
$\smqw(\smb{2k-1,2k})$ in 
$\Aut(P_{2k}^*)$.} 
\end{equation}
\end{corollary}
 
\begin{proof}
This  is  trivially true, as 
according to Theorem \ref{elp.t2} the two groups 
$\smqw(\smb{2k-1, 2k}) $ and $\qw(\beta_{2k-1,2k})$ coincide.
Hence \eqref{elp.e17} follows.
\end{proof}

We define 
\begin{equation}\mylabel{elp.e31}
J:= \text{ the  image of $\hap(\alpha_{2l-2,2l-1})$ in } \Aut(Q_{2l-1}^*).
\end{equation}
Note that this is the analogue to the definition  of $I$ in \eqref{el.e31}.
Of course $J$ is well defined as $\hap$ is a subgroup of 
 $ G(\beta_{2l-1}^*)$, 
and thus normalizes $Q_{2l-1}^*$. 
Furthermore, Theorem \ref{elp.t22} easily  implies
\begin{corollary}\mylabel{elp.co33}
The groups $\hap_{\mu}(\sma{2l-2,2l-1})$ and $\hap(\alpha_{2l-2,2l-1})$
have the same image inside 
\linebreak
$\Aut(Q_{2l-1}^*)= \Aut(Q_{2l-1,\mu}^*)$.
Hence 
$$J = J_{\mu},$$
where $J_{\mu}$ is the image of $\hap_{\mu}(\sma{2l-1,2l-1})$ 
in $\Aut(Q_{2l-1, \mu}^*)$. 
\end{corollary}
\begin{proof}
Same as the proof of Corollary \ref{el.co2}.
So, in view of  \eqref{121b}, it is clear that 
$\Aut(Q_{2l-1}^*) = \Aut(Q_{2l-1,\mu}^*)$, 
as $Q_{2l-1}^* = Q_{2l-1, \mu}^*$.
While 
$$
P_{2, 2l-1} \leq \hap(\alpha_{2l-2,2l-1})  \leq N(Q_{2l-1}^* \tin G(\alpha_{2,2l-1})).
$$
Thus 
$$
\hap(\alpha_{2l-2,2l-1}) = \hap_{\mu}(\sma{2l-2,2l-1}) \cdot P_{2,2l-1},
$$
where $ P_{2,2l-1} = C(Q_3\cdot Q_5 \cdots Q_{2l-1} \tin P_2)$.
As $Q_1$ centralizes $P_2$, we conclude that $P_{2,2l-1} = C(Q_{2l-1}^*  \tin P_2)$.
Therefore, Corollary \ref{elp.co33} follows.
\end{proof}

We conclude this section with the analogue to Theorems  \ref{el.t3} and 
\ref{el.t33}. 
\begin{theorem}\mylabel{elp.t3}
Assume that  the character $\beta_{2k-1,2k} \in \Irr(Q_{2k-1,2k})$ extends  
to $\qw(\beta_{2k-1,2k})$. Let  $\smqw$ be the group picked 
in Theorem \ref{elp.t2}. Then 
 the character 
$\smb{2k-1, 2k} \in \Irr(\smq{2k-1})$ extends to the group 
$\smqw(\smb{2k-1,2k})$.
\end{theorem}

\begin{proof}
The proof here is  as trivial was that of Theorem \ref{el.t33}.
So, $\smb{2k-1,2k}= \beta_{2k-1,2k}$, by Corollary 
\ref{elp.co1}, and $\smqw(\smb{2k-1,2k})= \qw(\beta_{2k-1,2k})$, by Theorem 
\ref{elp.t2}. Thus Theorem \ref{elp.t3} holds.
\end{proof}                                       

Theorem \ref{prel.t5} and Remark \ref{elp.rem1} imply
\begin{theorem}\mylabel{elp.t33}
Assume that $\alpha_{2l-2,2l-1}  \in \Irr(P_{2l-2,2l-1})$ extends to 
$\hap(\alpha_{2l-1,2l-1})$. Let $\hap_{\mu}$ be  a
$\pi$-Hall subgroup of
 $G_{\mu, \beta_{2l-1}^*}$,  chosen so that the conditions in  
 Theorem \ref{elp.t22} hold.
Then the character $\sma{2l-2l-1}$ extends to 
$\hap_{\mu}(\sma{2l-2,2l-1})$.
\end{theorem}

\begin{proof}
Same as that of Theorem \ref{el.t3}, so we omit it.
\end{proof}

%%% Local Variables: 
%%% mode: latex
%%% TeX-master: "kernel"
%%% End: 

\section{Kernels } \mylabel{ker}

\subsection{Inside  $Q_1$ } \mylabel{ker1}
In our common situation, we have a finite group $G$, and 
the fixed system described in \eqref{e0.1} and \eqref{e.AB}, that is 
\begin{subequations}\mylabel{ker.1}
\begin{gather}
\text{ the normal series: } 1=G_0\unlhd G_1 \unlhd \dots \unlhd G_{n} =  G,\\
\text{ the character tower:} \{1=\chi_0, \chi_1, \dots, \chi_n \} 
\mylabel{ker.e1} \\
\text{ the triangular set: }
\{Q_{2i-1}, P_{2r} | \beta_{2i-1}, \alpha_{2r}
 \}_{i=1, r=0} ^{l', \, \,  k'} \mylabel{ker.e2}
\end{gather}
and the Hall system 
\begin{equation}
\{ \ma, \mb \},
\end{equation}
that satisfies \eqref{e.AB}.
\end{subequations}
As usual, we assume known all the groups and their characters that 
accompany the above setting. In particular, we assume known the groups 
$P_{2k'}^*, Q_{2l'-1}^*, P_{2i, 2j+1},
 Q_{2i-1,2j}$,  as well as  $\qw, \hap$, and the characters 
$\alpha_{2k'}^*, \beta_{2l'-1}^*, 
 \alpha_{2i,2j+1}$ and  $\beta_{2i-1,2j}$.

We also assume that $S$ and $\zeta \in \Irr(S)$  satisfy \eqref{els.1}, 
i.e., $S$ is a normal subgroup of $G$ contained in $G_1$, 
and $\zeta$ is a $G$-invariant character of $S$ lying under $\beta_1=\chi_1$.
If $K$ denotes the kernel of $\zeta$, then $K$ is a normal subgroup of $G$, 
as $\zeta$ is $G$-invariant, and  $K$ is contained in $G_1 \leq G_i$, for all 
$i=1,\dots,n$.
Thus we can define the factor groups 
\begin{subequations}\mylabel{ker.e3}
\begin{equation}\mylabel{ker.e3a}
G_{K}=G/K \text{ and }  \skg{i}= G_i/K,
\end{equation}
for all $i=1,\dots,n$. Then  $\skg{i}$ is the image of $G_i$ in the factor group
$G_{K}$. This way we have created a normal series
\begin{equation}\mylabel{ker.e3b}
\skg{0}=1 \unlhd \skg{1} \unlhd \skg{2} \unlhd \dots \unlhd \skg{n} = G_{K}
\end{equation}
of $G_{K}$, that clearly satisfies Hypothesis \ref{hyp1}.
Along with that series we can associate a character tower that arises from 
\eqref{ker.e1}. 
Indeed, $K= \Ker(\zeta)$ is contained in the kernel of 
$\chi_i$, for all $i=1,\dots, n$,
as $\chi_i$ lies above the $G$-invariant character  $\zeta$.
Thus there exists a unique character $\skc{i}$ of the 
factor group  $\skg{i}=G_i/K$, that inflates to $\chi_i \in \Irr(G_i)$, whenever
$1\leq i \leq n$. 
Hence the set 
\begin{equation}\mylabel{ker.e3c}
\{1=\skc{0}, \skc{1}, \dots,  \skc{n} \},
\end{equation}
forms a character tower for  the series \eqref{ker.e3b}.
\end{subequations} 

As in the earlier sections, we  fix an integer $m=1,\dots,n$, and consider 
the smaller system
\begin{subequations} \mylabel{ker.sm1}
\begin{gather} 
 1=G_0\unlhd G_1 \unlhd \dots \unlhd G_{m} \unlhd G, \\
 \{1=\chi_0, \chi_1, \dots, \chi_m \}, \\
\{Q_{2i-1}, P_{2r} | \beta_{2i-1}, \alpha_{2r} \}_{i=1, r=0} ^{l, \, \,  k}. 
\end{gather}
\end{subequations}
Clearly the series \eqref{ker.e3b} and the tower \eqref{ker.e3c}
provide the smaller reduced system 
\begin{subequations}\mylabel{ker.smk}
\begin{gather}
\skg{0}=1 \unlhd \skg{1} \unlhd \skg{2} \unlhd \dots \unlhd \skg{m} \unlhd G_{K}\\
\{1=\skc{0}, \skc{1}, \dots,  \skc{m} \}.
\end{gather}
\end{subequations}

The aim of this section is to give a ``nice'' triangular set for 
\eqref{ker.e3b}, that corresponds to the tower \eqref{ker.e3c},
so that we can control the groups  $P_{2k}^*, Q_{2l-1}^*,I $ and $J$
of the system \eqref{ker.sm1}.  
This is done using the natural group epimorphism 
\begin{subequations}\mylabel{ker.e3.1}
\begin{equation}
\rho : G \to G/K.
\end{equation}
Assume that  $H$ is a subgroup of $G$ and 
 $\theta $ an irreducible character of $H$ 
having $K\cap H $ in its kernel. Let $\Irr(H | H \cap K )$ 
be the set of all such characters of $H$.
Then, as it  is well known,  any  $\theta \in \Irr(H | H\cap K)$
 determines, in a natural way, 
a unique  irreducible  character  $\rho(\theta) \in \Irr((HK)/K)$ 
such that 
\begin{equation}\mylabel{ker.e3.2}
[\rho(\theta)](\sigma K)= \theta(\sigma),
\end{equation}
for all $\sigma \in H$.
\end{subequations}
We remark that an arbitrary element $\tau \in G$
 fixes $\theta$ if and only if it normalizes $H$ and
its image $\rho(\tau)$ in $G/K$ fixes $\rho(\theta)$.
(Note that we could have $\tau \in G$ such that  $\rho(\tau)$ fixes
 $\rho(\theta)$, and thus normalizes $HK$,
but  $\tau$   moves $H$ around inside $HK$.)
Therefore
\begin{remark}\mylabel{ker.r1}
The homomorphism 
$$
\rho: G(\theta) \to (G/K)(\rho(\theta))=G_K(\rho(\theta))
$$
is onto if and only if for every $x \in G$ with   $\rho(x) \in
 G_K(\rho(\theta))$, there exists an element $\tau$  such that 
 $\tau \in xK \cap N(H \tin G)$.
\end{remark}
Some sufficient conditions that make the above homomorphism onto 
are given in the next lemmas.
\begin{lemma}\mylabel{ker.l2}
If $\theta \in \Irr(H | H \cap K )$ and $H$ is a $\hat{\pi}$-Hall 
subgroup of $HK$, for some set of primes $\hat{\pi}$, then the map
$$
\rho: G(\theta) \to G_K(\rho(\theta))
$$
is an epimorphism.
\end{lemma}
\begin{proof}
Let $x\in G$ be such that  $\rho(x)$ fixes $\rho(\theta) \in \Irr(\rho(H))$. 
Then  $\rho(x)$ normalizes $\rho(H)$. So $x$ normalizes the inverse image 
$HK$ of $\rho(H)$ in $G$.
Therefore 
$$
H^x \leq H^xK =(HK)^x = HK.
$$
As $H$ is a $\hat{\pi}$-Hall subgroup of $HK$, we
conclude that $H^x =H^k$, for some $k \in K$.
Hence the element $\tau= xk^{-1}$ normalizes $H$ and lies in $xK$.
 In view of Remark \ref{ker.r1} the proof is complete.  
\end{proof}
As a generalization of Lemma \ref{ker.l2} we have 
\begin{lemma}\mylabel{ker.l3}
Assume that $H_i$ is a subgroup  of $G$ and $\theta_i \in \Irr(H_i | H_i \cap K)$,
for all $i=1,\dots,s$, and some  integer $s \geq 1$.
Assume further that the product $H_1 \cdot H_2 \cdots H_s$ is a 
$\hat{\pi}$-Hall subgroup   of  $H_1 \cdot H_2 \cdots H_s \cdot K$, 
for some set of primes $\hat{\pi}$,  and that  
$H_1 \cdots H_s \cap H_iK = H_i$ for all $i=1,\dots, s$.
Let $x $ be an element in $G$ that normalizes $H_iK$, for all $i=1,\dots,s$. Then 
there exists some  $\tau \in xK$ such that $\tau$ normalizes $H_i$, 
for all $i=1,\dots, s$. 
Hence the map
$$
 \rho:G(\theta_1, \dots, \theta_s) \to 
G_K(\rho(\theta_1), \dots, \rho(\theta_s))
$$
 is an epimorphism.
\end{lemma}
\begin{proof}
Since  $x$ normalizes each  $H_iK$, it also  normalizes the product group 
$H_1 \cdot H_2 \cdots H_s \cdot K$.
As $H_1\cdots H_s$ is a $\hat{\pi}$-Hall subgroup of $H_ 1 \cdots H_s \cdot K$,
 we get that $(H_1 \cdots H_s)^x= (H_1 \cdots H_s)^k$, for 
some $k \in K$. Thus the element $\tau = x k^{-1}$ normalizes the product 
$H_1 \cdots H_s$. But $\tau$ also  normalizes $H_iK$, as $x$ does, 
for all $i=1,\dots,s$. Hence $\tau$ normalizes the intersection
$H_i=H_1 \cdots H_s \cap H_iK$, for each $i=1,\dots,s$.  This completes the proof 
of the first part of the lemma.

To show that the desired map is an epimorphism
it suffices to see, (according to Remark \ref{ker.r1}), 
that for any $x \in G$ with 
$\rho(x) \in G_K(\rho(\theta_1), \dots, \rho(\theta_s))$, there exists some  
$\tau   \in xK$ that normalizes the
groups $H_i$, for all $i=1,\dots,s$. 
But any such element $x$ normalizes $H_iK$, as $\rho(x)$ fixes
 $\rho(\theta_i) \in \Irr(H_iK/K)$, for all $i=1,\dots,s$.
Therefore, the first part of the lemma applies, and guarantees the 
existence of such a $\tau$.
\end{proof}

As we did with the system \eqref{ker.e3}, we 
will follow the same, standard by now, notation as in  Chapters \ref{pq} and
 \ref{pq:sec5}, with the addition of an extra $K$ in the subscripts.
We first observe 
\begin{lemma}\mylabel{ker.l1}
For all $j=1,\dots, n$, all $i=1,\dots, l'$, and  
all $r=1, \dots, k'$, we have that 
\begin{gather*}  
\skgs{j}:=\skg{j}(\skc{1}, \dots, \skc{j-1}) = G_{j}^*/K, \\
(Q_{2i-1}^*K )/K \in \Hall_{\pi'}(\skgs{2i-1}), \\ 
P_{2r}^* \cong (P_{2r}^*K) /K  \in \Hall_{\pi}(\skgs{2r}).
\end{gather*}
\end{lemma}

\begin{proof}
Since each $\chi_i$, for $i=1,\dots,n$, is a character of a normal subgroup $G_i$ of $G$ 
containing $K$, its stabilizer $G(\chi_i)$ is the unique image of $G_K(\chi_{i,K})=
G(\chi_i)/K$ in $G$.
The first part follows  from this and the definition of $G_j^*$ in \eqref{pq1}.
The other two parts are implied by the first, and the fact that $Q_{2i-1}^*$ 
and  $P_{2r}^*$ are $\pi'$-and $\pi$-Hall subgroups of $G_{2i-1}^*$ and 
$G_{2r}^*$, respectively. 
\end{proof}

 We also note
\begin{lemma}\mylabel{ker.r2}
The intersection    $K\cap Q_{2i-1}$ is a subgroup of  the kernel
 $\Ker(\beta_{2i-1})$ of $\beta_{2i-1} \in \Irr(Q_{2i-1})$, whenever
 $1\leq i \leq l'$.
\end{lemma}
\begin{proof}
Indeed, assume $i=1,\dots, l'$ is fixed. Then  $K \cap Q_{2i-1} \leq Q_1 \cap Q_{2i-1}= Q_{1,2i}$, (see \eqref{pq14b}
for the last equality). So $K \cap Q_{2i-1} \leq \Ker(\beta_{1}|_{Q_{1,2i}})$, 
as $K \leq \Ker(\beta_1) =\Ker(\chi_1)$. Since  $\beta_{1,2i}$ is the 
$P_{2i}^*$-Glauberman correspondent of $\beta_1$, it is a constituent of the 
restriction $\beta_1|_{Q_{1,2i}}$ of $\beta_1$ to $Q_{1,2i}$. 
Thus  $K\cap Q_{2i-1}$ is also a subgroup of the kernel, 
$\Ker(\beta_{1,2i})$, of this constituent $\beta_{1,2i}$. 
The character $\beta_{1,2i}$ is the unique character of $Q_{1,2i}$ lying 
under $\beta_{2i-1}$, according to Proposition \ref{pqremark1'}.
Thus $\Ker(\beta_{1,2i}) \leq \Ker(\beta_{2i-1})$. We conclude that 
\begin{equation}\mylabel{ker.e4}
K \cap Q_{2i-1} \leq \Ker(\beta_{1,2i}) \leq \Ker(\beta_{2i-1}),
\end{equation}
whenever $1\leq i \leq l'$. Thus Lemma \ref{ker.r2} follows.
\end{proof}
In view of above lemma, we see that the character   $\beta_{2i-1}$ determines a unique character 
\begin{subequations}\mylabel{ker.e4.1}
\begin{equation}
\beta_{2i-1,K}:= \rho(\beta_{2i-1}) \in \Irr(\rho(Q_{2i-1})),
\end{equation}
 by \eqref{ker.e3.2}.

Since   $P_{2r}$ is a $\pi$-group, for all $r=0, 1, \dots, k'$, 
we have that $P_{2r} \cong (P_{2r}K) /K$, as $K$ is a $\pi'$-group. 
Thus the character $\alpha_{2r} \in \Irr(P_{2r})$    
determines, under the above isomorphism, a unique character 
\begin{equation}
\ska{2r}:=\rho(\alpha_{2k}) \in \Irr((P_{2r}K )/K). 
\end{equation}
\end{subequations}

Now we can define the desired $K$-triangular set for \eqref{ker.e3b}.
\begin{theorem}\mylabel{ker.t1}
For every $r=0, 1, \dots, k'$ and $i=1, \dots, l'$ we define
\begin{align}\mylabel{ker.5}
\skp{2r} =\rho(P_{2r})= (P_{2r}K) /K \cong P_{2r} , \notag \\
 \skq{2i-1} =\rho(Q_{2i-1})=(Q_{2i-1}K )/K.
\end{align}
Then the set 
\begin{equation}\mylabel{ker.e5}
\{\skq{2i-1}, \skp{2r}| \skb{2i-1}, \ska{2r} \}_{i=1, r=0}^{l', \, k'}
\end{equation}
is  representative of the unique $G_{K}$-conjugacy class of triangular
 sets that corresponds  to \eqref{ker.e3c}.
\end{theorem}

\begin{proof}
It suffices to verify all the relations
in \eqref{xxx} and in Theorem \ref{tow--tri}, for the $K$ case.
We will do that using the map $\rho$ and the fact that the same relations 
hold for the set \eqref{ker.e2}.

The first two relations (\ref{xxx}a, b) of \eqref{xxx}, hold trivially.
To see that \eqref{x3} and \eqref{x5} hold, it is enough to show that the 
maps
\begin{equation}\mylabel{ker.e5.5}
\begin{aligned}
\rho:G_{2r}(\alpha_2, \dots, \alpha_{2r-2}, \beta_1, \dots, \beta_{2r-1})
&\to G_{2r,K}(\rho(\alpha_{2}), \dots, 
\rho(\alpha_{2r-2}), \rho(\beta_1), \dots, \rho(\beta_{2r-1}))\\
\rho:G_{2i-1}(\alpha_2, \dots, \alpha_{2i-2}, \beta_1, \dots, \beta_{2i-3})
&\to G_{2i-1,K}(\rho(\alpha_{2}), \dots, 
\rho(\alpha_{2i-2}), \rho(\beta_1), \dots, \rho(\beta_{2i-3}))\\
\end{aligned}
\end{equation}
are onto, whenever $1\leq r \leq k'$ and $2\leq i \leq l'$.
Indeed, that would be enough to guarantee that if we apply  
$\rho$ to \eqref{x3} and \eqref{x5} we get 
\begin{align*}
\rho(P_{2r})=P_{2r,K} &\in \Hall_{\pi}(G_{2r,K}(\rho(\alpha_{2}), \dots, 
\rho(\alpha_{2r-2}), \rho(\beta_1), \dots, \rho(\beta_{2r-1}))),\\
\rho(Q_{2i-1}) =Q_{2i-1,K}  &\in \Hall_{\pi'}(G_{2i-1,K}(\rho(\alpha_{2}),
 \dots, \rho(\alpha_{2i-2}), \rho(\beta_1), \dots, \rho(\beta_{2i-3}))).
\end{align*}

We first fix some $i=2, \dots, l'$ and  consider the map 
\begin{equation}\mylabel{ker.e8}
\rho:G_{2i-1}(\alpha_2, \dots, \alpha_{2i-2}, \beta_1, \dots, \beta_{2i-3})
\to G_{2i-1,K}(\rho(\alpha_{2}), \dots, 
\rho(\alpha_{2i-2}), \rho(\beta_1), \dots, \rho(\beta_{2i-3})).
\end{equation}
The groups $P_{2j}$ and their characters $\alpha_{2j}$, with $1\leq j
 \leq i-1$, satisfy the hypotheses of Lemma \ref{ker.l3}. That is,
the product   $P_{2} \cdots P_{2i-2}=P_{2i-2}^*$ 
forms a group, that is actually
 a $\pi$-Hall subgroup of $P_{2i-2}^* K$. Furthermore, 
$P_{2i-2}^*\cap P_{2j}K = P_{2j}$, for all $j=1, \dots, i-1$.
We conclude that the map 
\begin{equation}\mylabel{ker.e9}
\rho:G(\alpha_2, \dots, \alpha_{2i-2})
\to G_K(\rho(\alpha_{2}), \dots, \rho(\alpha_{2i-2})),
\end{equation}
is an epimorphism.

Let $x \in G$ with $\rho(x) \in  G_{2i-1,K}(\rho(\alpha_{2}),
 \dots, \rho(\alpha_{2i-2}), \rho(\beta_1), \dots, \rho(\beta_{2i-3}))$.
Then the epimorphism of \eqref{ker.e9} allows to assume that 
$x$ is an element of $G(\alpha_2, \dots \alpha_{2i-2})$.
Thus to prove that the map in \eqref{ker.e8} is onto, it suffices to show 
that $x$ normalizes $Q_{2j-1}$, for all  $j=1, \dots, i-1$. 
Clearly $x$ normalizes $Q_1 \unlhd G$. For the rest we will induct on $j$.

Assume that $x$ normalizes $Q_{1}, \dots, Q_{2r-1}$,
 for some $r$ with $1 \leq r < i-1$. Then $x$ fixes 
$\beta_1, \dots, \beta_{2r-1}$.  Hence $x$  normalizes the group 
$G_{2r+1}(\alpha_{2} , \dots, \alpha_{2r}, \beta_1, \dots, \beta_{2r-1})$.
 But $Q_{2r+1}$ is a $\pi'$-Hall subgroup of this latter group, by \eqref{x5},
 and $x$
 normalizes the $\pi'$-group  $Q_{2r+1} K $, as $\rho(x)$ fixes $\rho(\beta_{2r+1})
\in \Irr((Q_{2r+1}K)/K)$.
Therefore $x$ normalizes the intersection 
$$
Q_{2r+1}K \cap 
G_{2r+1}(\alpha_{2} , \dots, \alpha_{2r}, \beta_1, \dots, \beta_{2r-1})= 
Q_{2r+1}.
$$
This completes the proof of the inductive argument, and thus 
shows  that the map in \eqref{ker.e8} is onto. 

A similar proof shows that the other map in \eqref{ker.e5.5} is also onto.

To prove \eqref{x6}, 
we have to show that $\skb{2i-1}$ lies above 
 the $\skp{2i-2}$-Glauberman correspondent 
$\skb{2i-3,2i-2} \in \Irr(\skq{2i-3, 2i-2})$
 of $\skb{2i-3}$, for all $i=2, \dots, l'$.
 First notice that 
\begin{equation*}
\skq{2i-3,2i-2}=C(\skp{2i-2} \tin \skq{2i-3})= 
C((P_{2i-2}K)/K \tin (Q_{2i-3}K)/K).
\end{equation*}
The fact that $(|Q_{2i-3}|, |P_{2i-2}|)=1$ and Glauberman's lemma 
 imply  that 
$$
C((P_{2i-2}K)/K \tin (Q_{2i-3}K)/K)=(C(P_{2i-2} \tin Q_{2i-3})K)/K=
(Q_{2i-3,2i-2} K)/K=\rho(Q_{2i-3,2i-2}).
$$
Hence $\skq{2i-3,2i-2}= \rho(Q_{2i-3,2i-2})$, whenever $2\leq i \leq l'$.
Furthermore, $\beta_{2i-3,2i-2}$ is the $P_{2i-2}$-Glauberman correspondent 
of $\beta_{2i-3}$. Thus $\rho(\beta_{2i-3,2i-2})$ is the 
$\rho(P_{2i-2})=P_{2i-2,K}$-Glauberman correspondent of 
$\rho(\beta_{2i-3})$, by Proposition  \ref{daco}. In conclusion, 
\begin{equation}\mylabel{ker.e10}
\rho(\beta_{2i-3,2i-2})=\skb{2i-3,2i-2},
\end{equation}
for all $i=2, \dots, l'$. This, along with the fact that $\beta_{2i-1}$ lies 
above $\beta_{2i-3,2i-2}$, implies that $\skb{2i-1}=\rho(\beta_{2i-1})$
lies above $\beta_{2i-3,2i-2,K}=\rho(\beta_{2i-3,2i-2})$. 
Thus \eqref{x6} holds.

Similarly we can show
\begin{equation}\mylabel{ker.e5.1}
\rho(\alpha_{2r-2,2r-1})=\ska{2r-2,2r-1},
\end{equation}
for all $r=1,\dots, k'$. From this \eqref{x4} follows.

As far as Theorem \ref{tow--tri} is concerned, 
the relations in  \eqref{tow--tri1}, and \eqref{tow--tri2},
are easily translated to  the $K$-case using $\rho$ (for groups and characters). 
For example, we have  
$$
\skg{2i, 2i-1} :=\rho(G_{2i, 2i-1}) = \rho(P_{2i} \rtimes Q_{2i-1})
=\skp{2i} \rtimes \skq{2i-1},
$$
whenever $1\leq  i\leq k'$.
So the first part of 
\eqref{tow--tri1} holds for the $K$-case.
The proof for the rest is analogous, and we leave it to the reader.

We need to work more to show that \eqref{tow--tri3} holds for the $K$-case, 
i.e., to show that 
\begin{subequations}\mylabel{ker.e11}
\begin{multline}
G_{i, 2j-1, K}:= G_{i,K}(\alpha_{2,K}, \dots, \ska{2j-2}, \skb{1}, \dots, 
\skb{2j-1})\\=N(P_{2,K}, \dots, P_{2j-2,K}, Q_{1,K}, \dots, \skq{2j-1}
\tin G_{i,K}(\skc{1}, \dots, \skc{2j-1})), 
\end{multline}
and 
\begin{multline}
G_{i,2r,K}:= G_{i,K}(\alpha_{2,K}, \dots, \ska{2r}, \skb{1}, \dots, 
\skb{2r-1})\\
=N(P_{2,K}, \dots, P_{2r,K}, Q_{1,K}, \dots, \skq{2r-1}
\tin G_{i,K}(\skc{1}, \dots, \skc{2r})),
\end{multline}
\end{subequations}
whenever $j=1,\dots, l', \,  r=1,\dots, k'$, and $i=1, \dots, n$.
(Note that we have separated the odd from the even case in \eqref{tow--tri3}.)

We have already seen that the maps in \eqref{ker.e5.5} are onto.
So if we prove that the maps
\begin{subequations}\mylabel{ker.e12}
\begin{multline}\mylabel{ker.e12a}
\rho:N(P_2, \dots, P_{2j-2}, Q_1, \dots, Q_{2j-1} 
\tin G_i(\chi_1, \dots, \chi_{ 2j-1})) \\
\to
 N(P_{2,K}, \dots, P_{2j-2,K}, Q_{1,K}, \dots, \skq{2j-1}
\tin G_{i,K}(\skc{1}, \dots, \skc{2j-1})),
\end{multline}
and 
\begin{multline} \mylabel{ker.e12b}
\rho:N(P_2, \dots, P_{2r}, Q_1, \dots, Q_{2r-1} 
\tin G_i(\chi_1, \dots, \chi_{ 2r})) \\
\to
 N(P_{2,K}, \dots, P_{2r,K}, Q_{1,K}, \dots, \skq{2r-1}
\tin G_{i,K}(\skc{1}, \dots, \skc{2r})),
\end{multline}
\end{subequations}
are onto, then the equations in \eqref{ker.e11} hold, as 
we can apply  $\rho$ to the equation \eqref{tow--tri3}.

We will prove that the map in \eqref{ker.e12a} is onto 
and leave the proof of \eqref{ker.e12b} to the reader.
The idea for the proof is the same as that used to prove the maps in 
 \eqref{ker.e5.5} were onto.

So assume that  $x \in G$ is such that 
$\rho(x)$ lies $G_i(\chi_{1,K}, \dots, \skc{2j-1})$ and normalizes the groups 
 $P_{2,K}, \dots, P_{2j-2,K}$ and $Q_{1,K}, \dots, \skq{2j-1}$, for some fixed
 $i=1, \dots, n$ and $j=1, \dots, l'$. Then it is easy to see that 
$x$ lies in $G_i$ and  fixes $\chi_1, \dots, \chi_{2j-1}$.
Furthermore, $x$ normalizes $P_{2}K , \dots, P_{2j-2}K$, while the groups 
$P_{2}, \dots, P_{2j-2}$ and the characters $\alpha_2, \dots, \alpha_{2i-2}$, 
  satisfy the hypotheses of Lemma \ref{ker.l3}.
Hence we can assume that  $x$ normalizes the groups $P_2, \dots, P_{2j-2}$.

It remains to show that $x$ normalizes the groups $Q_1, \dots, Q_{2j-1}$.
We use induction on $j$. Clearly $x$ normalizes  $Q_1 \unlhd G$.
Now  assume that $x$ normalizes $Q_1, \dots, Q_{2r-1}$ for some $r$ with 
$1\leq r <j$. Then $x$ normalizes 
$N(P_2, \dots, P_{2r}, Q_1, \dots, Q_{2r-1} \tin G_{2r+1}(\chi_1, \dots, 
\chi_{2r-1}) )$, as $G_{2r+1} \unlhd G$. 
But this normalizer equals 
$G_{2r+1, 2r}=P_{2r} \rtimes Q_{2r+1}$, having $Q_{2r+1}$ as a $\pi'$-Hall
subgroup. As $x$ also normalizes the $\pi'$-group  $Q_{2r+1}K$, (since $\rho(x)$  normalizes
 $Q_{2r+1,K}$), we conclude that $x$ normalizes the intersection 
$(P_{2r} \rtimes Q_{2r+1}) \cap Q_{2r+1}K =Q_{2r+1}$. 
This completes the proof of  the inductive step,  thus 
proving that the map \eqref{ker.e12a} is an epimorphism.  

So, with analogous proofs left to the reader, 
 Theorem \ref{ker.t1} follows.
\end{proof}

Clearly we have 
\begin{remark}\mylabel{ker.rema1}
For any fixed $m=1,\dots,n$, the smaller set 
$$
\{\skq{2i-1}, \skp{2r}| \skb{2i-1}, \ska{2r} \}_{i=1, r=0}^{l, \, k}
$$
is a triangular set for (\ref{ker.smk}a) that corresponds to 
the character tower (\ref{ker.smk}b). So now we have a complete smaller  $K$-system 
\begin{subequations}\mylabel{smk}
\begin{gather}
\skg{0}=1 \unlhd \skg{1} \unlhd \skg{2} \unlhd \dots \unlhd \skg{m} \unlhd G_{K},\\
\{1=\skc{0}, \skc{1}, \dots,  \skc{m} \},\\
\{\skq{2i-1}, \skp{2r}| \skb{2i-1}, \ska{2r} \}_{i=1, r=0}^{l, \, k}. 
\end{gather}
\end{subequations}
\end{remark}

For any $k=1,\dots,k'$, we write 
 $\alpha_{2k,K}^*$    for the 
$\skq{3}, \dots, \skq{2k-1}$-correspondent of 
$\ska{2k}$ (see Definition \ref{p*d1}). 
Then the above theorem implies
\begin{corollary}\mylabel{ker.co1}
\begin{gather*}
P_{2k,K}^* =\rho(P_{2k}^*) \cong P_{2k}^*, \\
\alpha_{2k, K}^* =\rho(\alpha_{2k}^*) \in \Irr((P_{2k}^* K)/K),
\end{gather*}
for all $k=1,\dots,k'$.
\end{corollary}

\begin{proof}
In view of  Lemma \ref{ker.l1} and \eqref{ker.5} we have that 
  $$
P_{2k, K}^* =P_{2, K} \cdots P_{2k,K}  =\rho(P_{2}) \cdots \rho( P_{2k})= \rho(P_{2k}^*)
\cong  P_{2k}^*,
$$ 
for all $k=1,\dots,k'$.
The character  $\alpha_{2k}^*$ was defined as  the
 $Q_3, \dots, Q_{2k-1}$-correspondent of $\alpha_{2k}$, for all such $k$. 
Hence  Proposition  \ref{daco} implies that  $\rho(\alpha_{2k}^*)$ is the 
$\rho(Q_3), \dots, \rho(Q_{2k-1})$-correspondent of 
$\rho(\alpha_{2k})$, whenever $k=1,\dots, k'$.
But $\rho(\alpha_{2k}) = \ska{2k}$, by (\ref{ker.e4.1}b), 
and $\rho(Q_{2i-1})= \skq{2i-1}$, by \eqref{ker.5}, for all 
$i=1,\dots, k' \leq l'$.
Thus, the $\rho(Q_3), \dots, \rho(Q_{2k-1})$-correspondent of 
$\rho(\alpha_{2k})$ is nothing else but the $\skq{3}, \dots, 
\skq{2k-1}$-correspondent $\alpha_{2k, K}^*$ of $\alpha_{2k,K}$. 
We conclude that $\alpha_{2k,K}^* = \rho(\alpha_{2k}^*)$,
for all $k=1,\dots, k'$, 
and the corollary follows.
\end{proof}

Furthermore, 
\begin{corollary}\mylabel{ker.co11}
For all $l=1,\dots, l'$, the character $\beta_{2l-1, K}^* \in \Irr(Q_{2l-1,K}^*)$ 
is the unique character of $Q_{2l-1,K}^*= (Q_{2l-1}^*K)/K$ that inflates to 
$\beta_{2l-1}^* \in \Irr(Q_{2l-1}^*)$. Hence $ \beta_{2l-1, K}^* = 
\rho(\beta_{2l-1}^*)$.
\end{corollary}

\begin{proof}
In view of \eqref{ker.5},
 the product group $Q_{2l-1, K}^*= \skq{1} \cdots \skq{2l-1}$
is the image under $\rho$ of $Q_{2l-1}^*= Q_1 \cdots Q_{2l-1} $, that is, 
$$
Q_{2l-1, K}^* = \rho(Q_{2l-1}^*) = (Q_{2l-1}^*K)/K,
$$
for all $l=1,\dots,l'$.
Furthermore,  $K $ is a subgroup of $Q_1 = Q_1^*$, and is 
contained in $\Ker(\beta_1)$. 
As $\beta_1 = \beta_1^*$ is the unique character of $Q_1^*$ lying 
under $\beta_{2l-1}^*$,  we conclude that $K $ is a subgroup
of $\Ker(\beta_{2l-1}^*)$, for all $l=1,\dots,l'$.
Hence there is a unique character, $\rho(\beta_{2l-1}^*)$ 
of $(Q_{2l-1}^*K)/K$ that inflates to $\beta_{2l-1}^*$.
It suffices to show that $\rho(\beta_{2l-1}^*) = \beta_{2l-1, K}^*$.

Indeed, 
since $\beta_{2l-1}^*$ is the $P_2, \dots, P_{2l-2}$-correspondent  of
 $\beta_{2l-1}$, we get that 
 $\rho( \beta_{2l-1}^*)$  is the $\rho(P_2), \dots, \rho(P_{2l-2})$-correspondent
of $\rho(\beta_{2l-1})$, for all $l=1,\dots,l'$, (see Proposition  \ref{daco}).
But $\rho(P_{2i})= P_{2i, K}$, for all $i=1,\dots,k'$,  by \eqref{ker.5}, 
and $\rho(\beta_{2l-1}) = \skb{2l-1}$, by \eqref{ker.e4.1}.
Hence, $\rho( \beta_{2l-1}^*)$  is the $\skp{2}, \dots, \skp{2l-2}$-correspondent
of $\skb{2l-1}$, for all $l=1,\dots,l'$.
This completes the proof of the corollary. 
\end{proof}

Even more,
the Hall system $\{ \ma, \mb\}$ is  nicely transfered via $\rho$ to a Hall 
system of \eqref{ker.e3}, as the next theorem shows.
\begin{theorem}\mylabel{ker.t4}
Let  
$$
\ma_K :=\rho(\ma) = (\ma K)/K \cong \ma
  \text{ and } \mb_K= \rho(\mb)= (\mb K)/K = \mb/K.
$$
Then $\{ \ma_K, \mb_K \}$ is  a Hall system for $G_K$ that satisfies 
the equivalent of \eqref{e.AB} 
for the $K$-case.
\end{theorem}

\begin{proof}
The maps 
\begin{gather*}
\rho: G(\chi_1, \dots, \chi_h) \to G_K(\skc{1}, \dots, \skc{h}),\\
\rho: \ma(\chi_1, \dots, \chi_h) \to  \ma_K(\skc{1}, \dots, \skc{h}),\\
\rho: \mb(\chi_1, \dots, \chi_h) \to  \mb_K(\skc{1}, \dots, \skc{h}),
\end{gather*}
are clearly onto, as  $\chi_i \in \Irr(G_{i})$ with $G_i \unlhd G$, for all 
$i=1,\dots, h$, and all $h=1,\dots,n$.
This,  along with (\ref{e.AB}a, b), implies 
\begin{gather*}
\ma_K \in \Hall_{\pi}(G_K), \quad  \mb_K \in \Hall_{\pi'}(G_K),  \\
 \ma_K(\skc{1}, \skc{2},\dots, \skc{h}) \text{ and } 
\mb(\skc{1}, \skc{2},\dots,
\skc{h}) \text{ form a Hall system for }\\
 G_K(\skc{1}, \skc{2},\dots, \skc{h}),
\end{gather*}
for all $h=1, \dots, n$.

In addition,  (\ref{e.AB}c) and Corollary \ref{ker.co1}  imply 
$$\ma_K(\skc{1}, \dots, \skc{n})= \rho(\ma(\chi_1, \dots, \chi_n))= 
\rho(P_{2k'}^*)= P_{2k', K}^*.
$$ 
Similarly we get that $\mb_K(\skc{1}, \dots, \skc{n})= \rho(Q_{2l'-1}^*)
= Q_{2l'-1,K}^*$.
Hence the groups $\ma_K, \mb_K$ form a Hall system for $G_K$,
 and  satisfy \eqref{e.AB} for the $K$-case.

Since  $(|\ma|, |K|)=1$, we clearly have $(\ma K)/K \cong \ma$.
This completes the proof of the theorem.
\end{proof}

As a corollary of the above theorem we have 
\begin{corollary}\mylabel{ker.co5}
For any $k=1,\dots,k'$ and any $l=1,\dots,l'$ we have 
\begin{gather*}
\rho(N(P_{2k}^* \tin \mb(\chi_1,\dots,\chi_{2k}))) = 
N(P_{2k, K}^*  \tin \mb_K(\skc{1}, \dots, \skc{2k})), \\
\rho(N(Q_{2l-1}^* \tin \ma(\chi_1,\dots, \chi_{2l-1})))=
N(Q_{2l-1,K}^* \tin \ma_K(\skc{1},\dots,\skc{2l-1})).
\end{gather*}
\end{corollary}

\begin{proof}
As $(|P_{2k}^*|, |K |)=1$, we have 
$$
(N(P_{2k}^* \tin \mb(\chi_1,\dots,\chi_{2k} ))K)/K =
N( (P_{2k}^*K)/K  \tin (\mb(\chi_1, \dots, \chi_{2k})K)/K).
$$
Also,   $N( (P_{2k}^*K)/K  \tin (\mb(\chi_1, \dots, \chi_{2k})K)/K)$ equals 
$N(P_{2k, K}^* \tin \mb_K(\skc{1},\dots, \skc{2k}))$, 
as 
$P_{2k, K}^* = (P_{2k}^* K) / K$, by Corollary \ref{ker.co1}, 
 and $(\mb(\chi_1,\dots,\chi_{2k})K)/K = \mb_K(\skc{1},\dots, \skc{2k})$,
by Theorem \ref{ker.t4}.
Furthermore, as 
$(N(P_{2k}^* \tin \mb(\chi_1,\dots,\chi_{2k} ))K)/K =
\rho(N(P_{2k}^* \tin \mb(\chi_1,\dots,\chi_{2k})))$, the first part of the
 corollary follows.
 
The proof for the  second equation is similar.
\end{proof}

Working on the smaller system \eqref{ker.sm1}, we can now prove
\begin{theorem}\mylabel{ker.t2}
Let $\{ \ma_K, \mb_K \}$ be a  Hall system of $G_K$ that 
arises from $\{ \ma, \mb \}$ via Theorem \ref{ker.t4}. For 
 any fixed $m=1,\dots,n$, 
we choose the groups $\qw$ and $\qw_K$ to satisfy the conditions in Theorem \ref{sy.H}
for the systems \eqref{ker.sm1} and \eqref{smk}, respectively. Then 
\begin{equation}\mylabel{ker.e6}
 \qw_K(\beta_{2k-1,2k,K})=\rho(\qw(\beta_{2k-1,2k})) =(\qw(\beta_{2k-1,2k})K )/K.
\end{equation}
\end{theorem}

\begin{proof}
We choose $\qw$ and $\qw_K$ to satisfy the conditions  Theorem \ref{sy.H} for
the systems \eqref{ker.sm1} and \eqref{smk}, respectively.
Therefore  $N(P_{2k}^* \tin \mb(\chi_1,\dots,\chi_{2k}))= \qw(\beta_{2k-1,2k})$.
In addition, we have 
\linebreak
   $N(P_{2k,K}^* \tin \mb_K(\skc{1},\dots,\skc{2k}))=\qw_K(\skb{2k-1,2k})$.
This, along with Corollary \ref{ker.co5}, implies Theorem \ref{ker.t2} 
\end{proof}

Similarly we can show
\begin{theorem}\mylabel{ker.t22}
Assume that $\{\ma_K, \mb_K\}$ are   as above. Assume further that,  
for any fixed $m=1,\dots,n$, 
we choose the groups $\hap$ and $\hap_K$ to satisfy the conditions in  Theorem \ref{sy.t1}
for the systems \eqref{ker.sm1} and \eqref{smk}, respectively. 
Then  
$$
\hap_K(\ska{2l-2,2l-1})= \rho(\hap(\alpha_{2l-2,2l-1})) = 
(\hap(\alpha_{2l-2,2l-1})K)/K 
$$
Hence 
$$
\hap_K(\ska{2l-2,2l-1})\cong \hap(\alpha_{2l-2,2l-1}).
 $$
\end{theorem}

\begin{proof}
Choose $\hap$ and $\hap_K$ to satisfy  the conditions in Theorem \ref{sy.t1}, for the systems
\eqref{ker.sm1} and \eqref{smk}, respectively. Then the first part of 
Theorem \ref{ker.t22} follows from Corollary \ref{ker.co5}.
Note that $(\hap(\alpha_{2l-2,2l-1})K)/K \cong \hap(\alpha_{2l-2,2l-1})$,
as $(|\hap(\alpha_{2l-2,2l-1})|, |K|)=1$. Hence the theorem follows.
\end{proof}

As 
$\qw_K(\skb{2k-1,2k}) = \rho(\qw(\beta_{2k-1,2k}))$,
while $P_{2k}^* \cong P_{2k, K}^* =\rho(P_{2k}^*)$, 
 the action  of $\qw(\beta_{2k-1,2k})$ on $P_{2k}^*$ 
is carried onto the action of $\qw_K(\skb{2k-1,2k})$ on 
$P_{2k, K}^*\cong P_{2k}^*$, 
in the sense that
 \begin{equation}\mylabel{ker.e7.1}
\rho(\sigma^{\tau}) = \rho(\sigma)^{\rho(\tau)} \in P_{2k,K}^*
\end{equation}
for any $\sigma \in P^*_{2k}$ and any $\tau \in \qw(\beta_{2k-1,2k})$.
Let    $I_K$ be   the image of $\qw_K(\skb{2k-1,2k})$ in the
automorphism group 
 $\Aut(P_{2k, K}^*)$. As the  isomorphism $\rho$ of $P_{2k}^*$ 
onto $P_{2k, K}^*$ induces an isomorphism 
of $\Aut(P_{2k}^*)$ onto $\Aut(P_{2k, K}^*)$, we conclude that 
this isomorphism carries 
the image of  $\qw(\beta_{2k-1,2k})$  in the former automorphism group onto 
the image of $\qw_K(\beta_{2k-1,2k, K})$ in the latter such group.
So we have an isomorphism 
\begin{multline}\mylabel{ker.e7}
\rho_K: I=\text{ the image of $\qw(\beta_{2k-1,2k})$  in $\Aut(P_{2k}^*)$ } \\
\to
I_K= \text{  the image of $\qw_K(\skb{2k-1,2k})$ in } \Aut(P_{2k, K}^*). 
\end{multline}
Furthermore, identifying $P_{2k}^*$ with $P_{2k, K}^*$ and 
$\Aut(P_{2k}^*)$ with $\Aut(P_{2k, K}^*)$ we conclude
\begin{corollary}\mylabel{ker.co2}
For any fixed $m=1,\dots,n$, 
the groups $\qw(\beta_{2k-1,2k})$ and $\qw_K(\skb{2k-1,2k})$ have the same image in 
$\Aut(P_{2k, K}^*) = \Aut(P_{2k}^*)$.
\end{corollary}

Similarly, Theorem \ref{ker.t22} implies
\begin{corollary}\mylabel{ker.co22}
For any fixed $m=1,\dots,n$, 
the groups $\hap(\alpha_{2l-2,2l-1})$ and $\hap_K(\ska{2l-2,2l-1})$
have the same image in $\Aut(Q_{2l-1,K}^*)$.
\end{corollary}

\begin{proof}
As we have seen in Corollary \ref{ker.co11} and 
 Theorem \ref{ker.t22}
 \begin{gather*}
Q_{2l-1,K}^* = \rho(Q_{2l-1}^*) = (Q_{2l-1}^*K)/K, \\
\hap_K(\ska{2l-2,2l-1}) = \rho(\hap(\alpha_{2l-2,2l-1})) 
\cong \hap(\alpha_{2l-2,2l-1}).
\end{gather*}
So the action of $\hap_K(\alpha_{2l-2,2l-1, K})$ on $Q_{2l-1,K}^*$ is given 
(similarly to  \eqref{ker.e7.1}) as 
\begin{equation}\mylabel{ker.e7.2}
\rho(x)^{\rho(y)} = \rho(x^{y}) \in Q_{2l-1,K}^*
\end{equation}
for any $x \in Q_{2l-1}^*$ and $y \in \hap(\alpha_{2l-2,2l-1})$.
Furthermore, $\hap(\alpha_{2l-2,2l-1})$ acts also on $Q_{2l-1,K}^*$ 
via
$$
\rho(x)^y = \rho( x^y),
$$
for any $x$ and $y$ as above. As the map 
$$
\rho: \hap(\alpha_{2l-2,2l-1}) \to \hap_K(\ska{2l-2,2l-1}),
$$
sending $y \in \hap(\alpha_{2l-2,2l-1})$ to $\rho(y) \in \hap_K(\ska{2l-2,2l-1})$,
is an isomorphism, Corollary \ref{ker.co22} follows.
\end{proof}

We conclude this section with 
\begin{theorem}\mylabel{ker.t3}
Assume that  the character $\beta_{2k-1,2k} \in \Irr(Q_{2k-1,2k})$ extends  
to $\qw(\beta_{2k-1,2k})$. Then 
 the character 
$\skb{2k-1, 2k} \in \Irr(\skq{2k-1,2k})$ extends to the group 
$\skqw(\skb{2k-1,2k})$.
\end{theorem}

\begin{proof}
Obvious, since $\ beta_{2k-1,2k, K} = \rho(\beta_{2k-1,2k})$  by \eqref{ker.e10}.
\end{proof}

Similarly we have. 
\begin{theorem}\mylabel{ker.t33}
Assume that  the character $\alpha_{2l-2,2l-1} \in \Irr(P_{2l-2,2l-1})$
extends to 
$\hap(\alpha_{2l-2,2l-1})$.  Then 
 the character $\ska{2l-2,2l-1} \in \Irr(\skp{2l-2,2l-1})$
 extends to the group 
$\hap_K(\ska{2l-2,2l-1})$.
\end{theorem}

\subsection{Inside  $P_2$} \mylabel{ker2}
Assume now that the system \eqref{ker.1} for the normal series 
$1=G_0 \unlhd G_1\unlhd \dots \unlhd G_{n} = G$ is fixed, but in addition 
$G_2$ satisfies  \eqref{lp.1}, i.e., 
$G_2$ is the direct product  of a $\pi$-and a $\pi'$-group, while 
$\chi_1 = \beta_1$ is $G$-invariant.
Then the triangular set \eqref{ker.e2} satisfies \eqref{117}, i.e., 
\begin{subequations}\mylabel{ke.e1}
\begin{align}
G_2=G_2(\beta_1)= P_2 \times Q_1 = P_2 \times G_1,\\
\chi_2= \alpha_2 \times \beta_1.
\end{align}
\end{subequations}
Assume further that $R$  is a   normal subgroup 
of $G$ and $\eta$  is a character in $\Irr(R)$  satisfying \eqref{eps.1}, i.e., 
$R$ is a subgroup of $P_2$  while the irreducible character $\eta$  of 
 $R$ is $G$-invariant and lies under $\alpha_2$.
We write $K$ for the kernel of $\eta$. (Note this is not the same $K$ as in Section 
\ref{ker1}.)
Then $K$ is a normal subgroup of $G$, as $\eta$ is $G$-invariant.
As in the previous section, we are interested at the factor group
 $G_K=G/K$. (Note that this time $K$ is a $\pi$-group.)
So we define the factor groups
 $\skg{i}=G_i/K$,  for all $i=2,\dots,n$.
We also write $\skg{1}= (G_1 K)/K$. Then (\ref{ke.e1}a) implies 
\begin{equation}\mylabel{ke.e2}
\skg{1}=(G_1 \times K ) /K \cong G_1.
\end{equation}
The series 
\begin{subequations}\mylabel{ke.e3}
\begin{equation}\mylabel{ke.e3a}
1=\skg{0} \unlhd \skg{1} \unlhd \skg{2} \unlhd \dots \unlhd \skg{n} = G_K
\end{equation}
is a normal series of $G$ that satisfies Hypothesis \ref{hyp1}.
Furthermore, for every $i=2,\dots, n$ the character $\chi_i$ lies 
above the $G$-invariant character  $\zeta$. Hence $K$ is a subgroup of 
$\Ker(\chi_i)$, for all such $i$.
So we can define again the character $\skc{i} \in \Irr(\skg{i})$, 
to be the unique character of $G_i/K$ that inflates to 
$\chi_i$, whenever $2\leq i \leq n$.
As $\skg{1} \cong G_1$, we denote by  $\skc{1}$ 
the unique irreducible character of $\skg{1}$ that corresponds to 
$\chi_1$,  via that isomorphism. So we get a character tower for 
\eqref{ke.e3}, 
\begin{equation}\mylabel{ke.e3b}
\{1=\skc{0}, \skc{1}, \skc{2}, \dots, \skc{n} \},
\end{equation}
\end{subequations}
that arises from the original tower \eqref{ker.e1}.

In conclusion, we have created a similar system to that of Section
 \ref{ker1}, with the only important difference being that $K$ is  a $\pi$-group
 instead of a $\pi'$-group.
The natural map  $\rho$, on groups and characters, that was defined at 
\eqref{ker.e3.1}, is carried unchanged in this situation.
Of course Remark \ref{ker.r1} and Lemmas \ref{ker.l2} and \ref{ker.l3} are still valid.
With the help of the same map $\rho$  we will define 
a triangular set for \eqref{ke.e3a} that corresponds to  \eqref{ke.e3b}. As we would expect, this set is going to be the 
mirror of the set \eqref{ker.e5}, with the roles of 
the $\pi$-and the  $\pi'$-groups interchanged.

We start with 
\begin{lemma}\mylabel{ke.l1}
For all $j=2,\dots,n$, all $i=1,\dots,l'$  and all $r=1,\dots k'$, 
we have 
\begin{gather*}
\skg{j}^*:=G_{j,K}(\skc{1},\dots,\skc{j-1})= G_j^*/K,\\
Q_{2i-1}^* \cong (Q_{2i-1}^*K)/K  \in \Hall_{\pi'}(G^*_{2i-1,K}),\\
P_{2r}^*/K = (P_{2r}^* K ) / K  \in \Hall_{\pi}(G^*_{2r,K}).
\end{gather*}
\end{lemma}
\begin{proof}
See the proof of  Lemma \ref{ker.l1}.
\end{proof}

Also in the same way we worked to prove  Lemma \ref{ker.r2} and \eqref{ker.e4.1}, 
but this time using the characters $\alpha_{2r}, \alpha_2$ and  
$\alpha_{2,2r-1}$ in the place of $\beta_{2i-1}, \beta_1$ and 
$\beta_{1,2i}$ respectively, we can see that 
\begin{remark}\mylabel{ke.r1}
The intersection $K\cap P_{2r}$ is a subgroup of the kernel
$\Ker(\alpha_{2r})$ of $\alpha_{2r} \in \Irr(P_{2r})$, whenever 
$1\leq r \leq k'$.
Thus,  for all such $r$, there exists a unique character  
\begin{equation}\mylabel{ke.e4}
 \alpha_{2r,K}:=\rho(\alpha_{2r}) \in \Irr((P_{2r}K)/K),
\end{equation}
that inflates to $\alpha_{2r} \in \Irr(P_{2r})$. 
\end{remark}

Furthermore, $Q_{2i-1}$ is a $\pi'$-group, and thus has order  coprime to $|K|$, for all 
$i=1,\dots, l'$. Hence $Q_{2i-1} \cong (Q_{2i-1}K)/K$, for all such $i$.
So the character $\beta_{2i-1} \in \Irr(Q_{2i-1})$
 determines, under the above isomorphism, a unique character 
\begin{equation}
\skb{2i-1}:=\rho(\beta_{2i-1}) \in \Irr((Q_{2i-1}K) /K),
\end{equation}
for all $i=1,\dots, l'$.

  We can now prove the main theorem of this section, the analogue 
of Theorem \ref{ker.t1}.
\begin{theorem}\mylabel{ke.t1}
For every $r=0, 1, \dots, k'$ and $i=1, \dots, l'$ we define 
\begin{gather}\mylabel{ke.e5}
\skp{2r}=\rho(P_{2r})=(P_{2r}K)/K, \\
\skq{2i-1}= \rho(Q_{2i-1})= (Q_{2i-1}K)/K \cong Q_{2i-1}.
\end{gather}
Then the  set 
\begin{equation}\mylabel{ke.e6}
\{\skq{2i-1}, \skp{2r}|\skb{2i-1}, \ska{2r}\}_{i=1, r=0}^{l',k'} 
\end{equation}
is a representative of the unique $G_K$-conjugacy class of triangular 
sets that corresponds to \eqref{ke.e3b}.
\end{theorem}

\begin{proof}
The proof is  the same as that of Theorem \ref{ker.t1}, 
if we interchange the roles of $P_{2r}$ and $ \alpha_{2r}$
with those of $Q_{2i-1}$ and $\beta_{2i-1}$, respectively. 
\end{proof}

We also get 
\begin{corollary}\mylabel{ke.co1}
For all $k=1,\dots,k'$,
we have $P_{2k, K}^* = \rho(P_{2k}^*)$. Furthermore, 
the character  $\alpha_{2k, K}^* $ is the unique character $\rho(\alpha_{2k}^*)$ 
of $P_{2k, K}^*$ that inflates to $\alpha_{2k}^* \in \Irr(P_{2k}^*)$.
\end{corollary}
\begin{proof}
Same as that of Corollary \ref{ker.co11}, 
with the roles of $\pi$ and $\pi'$ -interchanged.
\end{proof}

and 

\begin{corollary}\mylabel{ke.co11}
\begin{gather*}
Q_{2l-1, K}^* = \rho(Q_{2l-1}^*) \cong Q_{2l-1}^*, \\
\beta_{2l-1,K}^* = \rho(\beta_{2l-1}^*),
\end{gather*}
for all $l=1,\dots,l'$.
\end{corollary}

\begin{proof}
See Corollary \ref{ker.co1}.
\end{proof}

The same argument as that of Theorem \ref{ker.t4} implies 
\begin{theorem}\mylabel{ke.t4}
Let  
$$
\ma_K :=\rho(\ma) = (\ma K)/K 
  \text{ and } \mb_K= \rho(\mb)= (\mb K)/K \cong \mb.
$$
Then $\{ \ma_K, \mb_K \}$ forms a Hall system for $G_K$ that satisfies 
the equivalent of \eqref{e.AB} 
for the $K$-case.
\end{theorem}

Until the end  of the section,
  we fix an integer $m=2,\dots,n$ and consider the  smaller system \eqref{ker.sm1}.
Of course, as before, we get a smaller $K$-system \eqref{smk}, where now 
the triangular set is picked to be a subset of  \eqref{ke.e6}.
So as in  Theorems \ref{ker.t2} and \ref{ker.t22}, we have 
\begin{theorem}\mylabel{ke.t2}
Let $\{ \ma_K, \mb_K \}$ be tha Hall system of $G_K$ that 
arises from $\{ \ma, \mb \}$ via Theorem \ref{ke.t4}. 
For any fixed $m=1,\dots,n$, 
we choose the groups $\qw$ and $\qw_K$ to satisfy  the conditions in Theorem \ref{sy.H}
for the systems \eqref{ker.sm1} and \eqref{smk}, respectively. Then
\begin{equation}
 \qw_K(\beta_{2k-1,2k,K})=\rho(\qw(\beta_{2k-1,2k})) =(\qw(\beta_{2k-1,2k})K )/K.
\end{equation}
Hence $$
\qw_K(\beta_{2k-2,2k, K}) \cong \qw(\beta_{2k-1,2k}).$$ 
\end{theorem}
\begin{proof}
See Theorem \ref{ker.t22}
\end{proof}

and 

\begin{theorem}\mylabel{ke.t22}
Assume that $\{\ma_K, \mb_K\}$ are   as above. Assume further that,  
for any fixed $m=2,\dots,n$, 
we choose the groups $\hap$ and $\hap_K$ to satisfy the conditions in Theorem \ref{sy.t1}
for the systems \eqref{ker.sm1} and \eqref{smk}, respectively. 
Then 
$$
\hap_K(\ska{2l-2,2l-1})= \rho(\hap(\alpha_{2l-2,2l-1})) = 
(\hap(\alpha_{2l-2,2l-1})K)/K. 
$$
\end{theorem}

\begin{proof}
Same as that of Theorem \ref{ker.t2}.
\end{proof}

So we get the next two corollaries
\begin{corollary}\mylabel{ke.co2}
The groups  $\qw(\beta_{2k-1,2k})$ and 
$\qw_K(\skb{2k-1,2k})$, have the same image in  
the group of automorphisms  $\Aut(P_{2k, K}^*)$.
\end{corollary}
\begin{proof}
See the proof of Corollary \ref{ker.co22}.
\end{proof}

and 
\begin{corollary}\mylabel{ke.co22}
For any fixed $m=2,\dots,n$, 
the groups $\hap(\alpha_{2l-2,2l-1})$ and $\hap_K(\ska{2l-2,2l-1})$
have the same image in $\Aut(Q_{2l-1,K}^*)= \Aut(Q_{2l-1}^*)$.
\end{corollary}

\begin{proof}
See Corollary \ref{ker.co2}.
\end{proof}

We conclude the section and the chapter with 
\begin{theorem}\mylabel{ke.t3}
Assume that $\beta_{2k-1,2k} \in \Irr(Q_{2k-1,2k})$
 extends to $\qw(\beta_{2k-1,2k})$. Then the character
$\skb{2k-1,2k} \in \Irr(\skq{2k-1,2k})$ extends to 
$\qw_K(\skb{2k-1,2k})$.
\end{theorem}

\begin{proof}
Obvious.
\end{proof}

and 
\begin{theorem}\mylabel{ke.t33}
Assume that  the character $\alpha_{2l-2,2l-1} \in \Irr(P_{2l-2,2l-1})$
extends to 
$\hap(\alpha_{2l-2,2l-1})$.  Then 
 the character $\ska{2l-2,2l-1} \in \Irr(\skp{2l-2,2l-1})$
 extends to the group 
$\hap_K(\ska{2l-2,2l-1})$.
\end{theorem}

%%% Local Variables: 
%%% mode: latex
%%% TeX-master: t
%%% End: 

\chapter{ Linear  Limits }\mylabel{lim}

\section{Basic properties}\mylabel{lim1}
We say that $(G, A , \phi, N , \psi)$ is a {\em linear     quintuple}
if $A \leq N$ are normal subgroups of a finite group $G$, 
  $\phi\in \Lin(A)$  is a $G$-invariant linear character of $A$   and   
$\psi \in \Irr(N |\phi)$. Note that as 
$\phi$ is $G$-invariant,  $\Ker(\phi)$ is a normal subgroup of
$G$. Furthermore, $\Ker(\phi)=\Ker(\psi|_{A}) \leq \Ker(\psi)$. Hence
$\Ker(\phi)$ is contained in the largest normal subgroup $M$ of $G$
contained in $\Ker(\psi)$. 
(Note that  $M = \bigcap_{x \in G}(\Ker(\psi))^x$.) 
This, along with the fact that  $\phi$ is  linear, 
implies  that   $A$ is abelian modulo $M$, i.e., 
 $(AM) /M \cong A/(A\cap M)$ is abelian.
 
Let   $A' \unlhd G$  with $A \leq A' \leq N$. Let  $\phi' \in \Irr(A')$ be a 
linear character of $A'$ extending $\phi$ and lying under $\psi$.
Then $(A'M) /M \cong A'/(A'\cap M ) $ is also abelian.
Indeed,  $[A',A']$ is contained in $\Ker((\phi')^n)$, for every  $n \in
N$,  since  $\phi'$ is linear. 
Thus   $[A', A'] \leq \bigcap_{n \in N } \Ker((\phi')^n) $.
As the   restriction of $\psi $ to $A'$ is a sum of $N$-conjugates of
$\phi'$, 
we conclude that $[A',A']$ is contained in $\Ker(\psi|_{A'})= \bigcap_{n
\in N} \Ker((\phi')^n) $. So $[A', A'] \leq \Ker(\psi)$.
Furthermore, $[A', A']$ is a normal subgroup of $G$, as $A' \unlhd G $.
Hence $[A', A'] $ is a subgroup of $M$, which implies that 
$A' / (A'\cap M) $ is abelian.

Furthermore, we can use  Clifford theory to form a new linear quintuple
$(G', A', \phi', N', \psi')$, where $G'= G(\phi')$ and $  N'= N(\phi')$ 
 are  the stabilizers  of $\phi'$ in $G$ and $N$ respectively, 
and $\psi'$ is the $\phi'$-Clifford correspondent 
of $\psi \in \Irr(N | \phi')$. 
We say that $(G', A', \phi', N', \psi')$ is a { \em linear  reduction } 
of $(G, A , \phi, N , \psi)$. We call this reduction {\em  proper } if 
the reduced linear  quintuple is different from the original one, i.e., 
if and only if $A < A'$.
We can repeat this process and consider a linear  reduction of
 the linear  reduction  $(G', A', \phi', N', \psi')$.
Any linear  quintuple that  we reach after a series of such linear 
reductions, is called  a   { \em multiple linear reduction }
of   $(G, A , \phi, N , \psi)$.
A ``minimal'' multiple linear reduction 
  is a linear  quintuple
 that has no proper linear  reductions
We call such a minimal linear  quintuple a  {\em linear  limit } 
of $(G, A , \phi, N , \psi)$. We denote by $LL(G, A , \phi, N , \psi)$
the set of all the linear  limits of $(G, A , \phi, N , \psi)$, 
and by  
$$
(l(G),  l(A),  l(\phi),  l(N), l(\psi) )
$$
  any element  of that set. 
Assume that  $H$ is a subgroup of $G$ with $N \leq H \leq G$.
Then the quintuple $(H, A , \phi, N , \psi)$ is clearly a linear one. 
The definition of linear limits clearly implies 
\begin{remark}\mylabel{lim.r2}
If $(l(G), l(A), l(\phi), l(N),  l(\psi))$ is a linear limit of 
$(G, A, \phi, N , \psi)$, then 
$(l(G) \cap H , l(A), l(\phi), l(N), l(\psi))$  is a linear limit of 
$(H, A , \phi, N , \psi)$.
\end{remark}
The following is also straight forward:
\begin{remark}\mylabel{lim.r1}
If  we reach the linear quintuple $(G', A', \phi', N', \psi')$ after a series of  linear reductions 
starting with  the  quintuple   $(G, A, \phi, N, \psi)$, then 
$LL(G', A' , \phi', N' , \psi') \subseteq  LL(G, A , \phi, N , \psi)$.
\end{remark}

If  $(l(G), l(A) , l(\phi), l(N) , l(\psi))$ is a linear  limit of 
$(G, A , \phi, N , \psi)$, 
then we can form the quintuple $(l(G)/K, l(A)/K , l(\phi)/K, l(N)/K , 
l(\psi)/K)$
where $K=\Ker(l(\phi)))$  is the kernel 
of $l(\phi)$ (and thus a  normal subgroup of $l(G)$) 
 while  $l(\phi)/K$ and $l(\psi)/K$ 
are the unique characters of the factor groups $l(A)/K$ and $l(N)/K$  
that  inflate to  $l(\phi)$ and $l(\psi)$, respectively.
It is clear that the quintuple  $(l(G)/K, l(A)/K , l(\phi)/K, l(N)/K , 
l(\psi)/K)$ is 
linear, and that $l(\phi)/K$ is faithful. We 
  call  this triple a {\em faithful linear  limit} of  $(G, A, \phi, N, \psi)$
We denote by $FLL(G, A , \phi, N , \psi)$ 
the set of all faithful linear  limits of $(G, A , \phi, N , \psi)$, 
and by 
\begin{equation}\mylabel{lmf}
(fl(G),  fl(A),  fl(\phi),  fl(N), fl(\psi) ) 
\end{equation}
 any element  of that set.
  
Now assume that $N \unlhd H \unlhd G$ and  $\chi \in \Irr(H|\psi)$.
Then any  linear  reduction  $(G', A' , \phi', N' , \psi')$ of 
$(G, A , \phi, N , \psi)$ provides an irreducible character
 $\chi' \in \Irr(H')$, where $H'=G'\cap H = H(\phi')$ 
and $\chi'$ is the $\phi'$-Clifford correspondent of 
$\chi$, i.e.,  $\chi'$  lies above  $\phi'$ and induces  $\chi$.
We can repeat this process and consider the  Clifford correspondent
for $\chi'$ in the next linear  reduction of $(G', A' , \phi', N' ,
\psi')$
that we perform. When we reach a linear  limit 
$(l(G), l(A) , l(\phi), l(N), l(\psi))$, 
of  $(G, A , \phi, N , \psi)$ 
we have also reached  a character $\theta \in \Irr(l(G) \cap H)$ that induces 
$\chi$. Any  such character $\theta$, that arises by repeated  
Clifford correspondences on linear  reductions,
 we call a  {\em linear  limit}
of $\chi$. We write it as $\theta= l_{\phi, \psi}(\chi)$,
or more simply  as $l(\chi)$ if  the starting linear quintuple is clear.
 We also write as $l_{\phi, \psi}(H)=l(H)$ the domain of 
$l(\chi)$,  i.e., $l(H)=l(G)\cap H$.  Clearly 
$l(\chi)$ lies above $l(\phi)$ and $l(\psi)$, and induces $\chi$.
The collection of all linear limits of $\chi$ we write as $LL(\chi)$.
Note that $LL(\chi)$ is a subset of $CCC_N(\chi)$ as this was defined in  
\cite{isa2}. 
Furthermore, let   $(fl(G),  fl(A),  fl(\phi),  fl(N), fl(\psi) )$
 be a faithful linear limit, i.e., 
 \begin{equation}\mylabel{10.35}
(fl(G),  fl(A),  fl(\phi),  fl(N), fl(\psi) )
=(l(G)/K, l(A)/K , l(\phi)/K, l(N)/K , 
l(\psi)/K)
\end{equation}
where  $(l(G), l(A), l(\phi), l(N), l(\psi))$ is  a linear limit of 
$(G, A, \phi, N, \psi)$ and $K=\Ker(l(\phi))$. 
Then $K$ is a subgroup of $\Ker(l(\chi))$ the kernel of the linear limit
$l(\chi) \in \Irr(l(H))$ of $\chi$, as $\chi$ lies above the $l(G)$-invariant
 character $l(\phi) \in \Irr(l(A))$.
Thus $l(\chi)$ is  inflated   from 
a unique character $l(\chi)/K$ of the factor group
$l(H)/K$   
that we call {\em faithful linear limit } of $\chi$ and 
 write as $fl(\chi)$. We also write $fl(H)$ for the domain of $fl(\chi)$, 
i.e., $fl(H)= l(H)/K$.
The set of all faithful linear limits of $\chi$ we
 denote by $FLL(\chi)$.

We conclude these preliminary definitions of linear 
limits  with  the following straight forward observations.
\begin{remark}\mylabel{lm.rem1}
Let $(l(G), l(A), l(\phi), l(N), l(\psi))$ be a linear  limit  
of $(G, A, \phi, N, \psi)$. Let   $K$ be the kernel  $\Ker(l(\phi))$, 
and let 
 $(fl(G), fl(A), fl(\phi), fl(N), fl(\psi))$ be   the faithful 
linear limit defined in \eqref{10.35}. 
Then $l(G)$ is a subgroup of $G$ while $fl(G)$ is the section $l(G)/K$ 
of $G$.
Furthermore, any subgroup $H$ of $G$ with  $N \unlhd H$, has a 
limit and a faithful limit  group $l(H)= l(G) \cap H $ and 
$fl(H)= l(H)/K $, respectively, that
 satisfy  $l(N) \unlhd l(H) < l(G)$ and $fl(N) \unlhd fl(H) < fl(G)$. 
In addition,  $(l(H), l(A), l(\phi), l(N), l(\psi))$ and 
 $(fl(H), fl(A), fl(\phi), fl(N), fl(\psi))$ are  a linear and a faithful
 linear limit  respectively, of  $(H, A, \phi, N, \psi)$.
If  $H$ is normal in $G$, then $l(H)$ and $fl(H)$ are normal in 
$l(G)$ and $fl(G)$, respectively.
\end{remark}
\begin{defn}\mylabel{lm.d1}
By convention, whenever $N \leq H \leq G$ and $(\mbg, \mba, \VP, \mbn, \VPS)$
is a faithful linear limit of 
$(G, A, \phi, N, \psi)$,
we write $(\mbg \cap H, \mba, \VP, \mbn, \VPS)$ for the faithful linear limit 
of $(H, A, \phi, N, \psi)$ (described in Remark \ref{lm.rem1})
  that  $(\mbg, \mba, \VP, \mbn, \VPS)$ induces.
\end{defn}

We can say a little more for a special type of subgroup  $B$ of $G$. 
Assume that $B \leq C_G(N)$, i.e., $B$ centralizes $N$.
Then $B$ centralizes any $A'$ with $A \leq A' \leq N$. Hence it fixes any 
character $\phi' \in \Irr(A')$, and, in particular, those that extend $\phi$.
Hence $B$ is also a subgroup of $G(\phi')$. Furthermore, it centralizes 
$N(\phi') \leq N$.  Repeating the same argument
 at every linear reduction we perform (note that all such are inside $N$), 
we see  that $B$ is a subgroup of $l(G)$,  and  centralizes $l(N)$.
Furthermore, as $K \leq l(N)$, we get that 
$fl(B):=(B K) /K$ is a subgroup of $fl(G)$ that centralizes 
$fl(N) = l(N) /K$. Thus we have shown 
\begin{remark}\mylabel{lm.rem2}
If  $B$ is a subgroup of $G$ that centralizes $N$, i.e., 
$B \leq C_G(N)$, then $ B \leq l(G)$ centralizes $l(N)$ 
while $fl(B) = (BK)/K \leq fl(G)$  centralizes $fl(N)$. 
If in addition   $(|B|, |K|)=1$ then   $B \cong fl(B) \leq fl(G)$. 
\end{remark}

The next two lemmas are straight forward applications of the above
 definitions.   
\begin{lemma}\mylabel{lim.l2}
Any faithful linear limit of $(G, A, \phi, N, \psi)$ is a  minimal
linear quintuple, that is, no proper linear reductions
can be made to a faithful linear limit of $(G, A, \phi, N, \psi)$.
\end{lemma}

\begin{proof}
Let 
$$
(fl(G),  fl(A),  fl(\phi),  fl(N), fl(\psi) ) =
(\mfg/K, \mfa/K, \Phi/K, \mfn/K, \Psi/K),
$$
be a faithful linear limit of 
$(G, A, \phi, N, \psi)$, where $(\mfg, \mfa, \Phi, \mfn, \Psi)$
is a linear limit of the latter, and
$K =\Ker(\Phi)$. 
It is not hard to see that any linear reduction of 
$(\mfg/K, \mfa/K, \Phi/K, \mfn/K, \Psi/K)$ provides a linear reduction of 
$(\mfg, \mfa, \Phi, \mfn, \Psi)$.
Indeed, if $\hat{\gamma}$ is a linear extension  of 
$\Phi/K$ to a normal subgroup $\hat{\Gamma}$ of $\mfg/K$,  and lies under 
$\Psi/K$, then 
$\hat{\Gamma} = \Gamma/ K$, where  $\Gamma $ is a normal subgroup of $\mfg$. 
Furthermore, $\hat{\gamma} $ inflates to a unique character 
$\gamma \in \Irr(\Gamma)$, i.e., 
$\hat{\gamma} = \gamma/K$. 
Also, $\gamma$ is linear, as $\hat{\gamma} $ is, and 
lies under $\Psi \in \Irr(\mfn)$, as $\hat{\gamma}$ lies under
 $\Psi/K \in \Irr(\mfn/K)$. Hence we can form a linear reduction of 
 $(\mfg, \mfa, \Phi, \mfn, \Psi)$, using  the extension $\gamma$ of
 $\Phi$ to $\Gamma$.
As no proper linear reductions can be made to the quintuple
 $(\mfg, \mfa, \Phi, \mfn, \Psi)$, 
the lemma follows.
\end{proof}

\begin{corollary}\mylabel{lim.co2}
Let $(\mbg, \mba, \VP, \mbn, \VPS)$ be a faithful linear limit of $(G, A, \phi, N, \psi)$. 
Then  $\mba$ is a cyclic central subgroup of $\mbg$ (it could be trivial), and
is  maximal among the abelian $\mbg$-invariant subgroups 
of $\mbn$. Hence $\mba = Z(\mbn)$.
\end{corollary}

\begin{proof}
Clearly $\mba$ is a normal subgroup of $\mbg$, 
while $\VP$  is a linear faithful 
$\mbg$-invariant  character of $\mba$.
Then   $\mba$ is  cyclic, as it affords  a faithful linear character. 
The additional fact that $\VP$ is $\mbg$-invariant  implies that 
$\mba$ is a subgroup of the center $Z(\mbg)$ of $\mbg$. 

Now assume that  $B$ is an abelian $\mbg$-invariant subgroup of 
$\mbn$ that contains $\mba$.
Let $\beta \in \Irr(B)$ be any character of $B$ that lies
 above $\VP$ and under $\VPS$. 
(Clearly such a character exists if the group $B$ exists.)
Then $\beta$ is an extension of $\VP$ to $B$.
 Furthermore, if $\VPS_{\beta}$ is the $\beta$-Clifford 
correspondent  of $\VPS$, then   the quintuple 
$(\mbg(\beta), B , \beta, \mbn(\beta), \VPS_{\beta})$ 
is a linear reduction of $(\mbg, \mba, \VP, \mbn, \VPS)$.
According to  Lemma \ref{lim.l2}  the latter quintuple can't have any proper reductions. 
Therefore, $B = \mba$, and $\mba$ is maximal among the
 abelian  $\mbg$-invariant 
subgroups of $\mbn$. 

As the center $Z(\mbn)$ is an abelian  characteristic subgroup of $\mbn$, it is clearly
$\mbg$-invariant. Furthermore, $\mba$ is a subgroup of $Z(\mbn)$, as $\mba$ is a subgroup of 
$Z(\mbg)$.  So $\mba= Z(\mbn)$, and   Corollary \ref{lim.co2} follows.
\end{proof}

\begin{corollary}\mylabel{lim.co22}
Let $(\mbg, \mba, \VP, \mbn, \VPS)$ be a faithful
 linear limit of $(G, A, \phi, N, \psi)$.
 Assume further that $N$ is a $p$-group, for some odd prime $p$.
Then either $\mbn= Z(\mbn)$ is cyclic, or $\mbn$  is the central product 
$\mbn =E \odot \mba$, of a non-trivial 
 extra special $p$-group $E$  of exponent $p$, and $\mba = Z(\mbn)$.
In both  cases  the irreducible character $\VPS \in \Irr(\mbn)$ 
is zero on $\mbn - \mba$ and  a multiple of $\VP$  on $\mba$. Hence
 $\VPS$ is $\mbg$-invariant.
\end{corollary}

\begin{proof}
According to Corollary \ref{lim.co2}, the group $\mba= Z(\mbn)$ 
is cyclic, central in $\mbg$ and maximal among the abelian 
$\mbg$-invariant subgroups of $\mbn$.
Furthermore, the fact that  $(\mbg, \mba, \VP, \mbn, \VPS)$ is a linear
 quintuple, implies that $\VP$ is $\mbg$-invariant.

If $\mbn=\mba=Z(\mbn)$,
 then  $\mbn$ is cyclic and 
$\VP = \VPS$. In this case  the corollary  holds trivially.
 
If $\mbn > \mba$, then  $\mba \ne 1$, or else $\mba=\mbn=1$, since  $\mbn$ is a 
$p$-group. Furthermore, $\VP$ is a faithful linear
character of the cyclic group $\mba$. 
The fact that the cyclic group 
 $\mba= Z(\mbn)$ is maximal among the abelian  subgroups of $\mbn$, 
normal in $\mbg$ implies that every 
characteristic abelian  subgroup of the $p$-group $\mbn$
is cyclic and central. Hence  P.Hall's theorem (see Theorem 4.22
pp.75 in \cite{su}) implies that $\mbn$ is the central product
 $$
\mbn = E  \odot Z(\mbn)= E \odot \mba,
$$
where $E$ is a non-trivial (or else $\mbn$ is abelian), 
 extra special $p$-group  of exponent $p$. 
Hence  $\VPS$ must be the  unique character of $\mbn$
 that lies above the faithful character $\VP$ of the center,
(see Theorem 7.5 in \cite{hu2})
  Even more, this unique character is 
zero outside  $\mba$, while its restriction to $\mba$ is a multiple 
of $\VP$.
 In addition, the fact that $\VP$ is 
$\mbg$-invariant, while $\mbn$ is a normal subgroup of $\mbg$, 
makes $\VPS$  also $\mbg$-invariant. 
So Corollary \ref{lim.co22}  follows.
\end{proof}

\begin{lemma}\mylabel{lim.l3}
Let $(\mbg, \mba, \VP, \mbn, \VPS)$ be a faithful linear limit of 
the linear  quintuple $(G, A, \phi, N, \psi)$. If $\T \in \Irr(\mbg| \VP)$,
then there exists a unique  $ \chi \in \Irr(G| \phi)$, so that $\T$ is a faithful
 linear limit of $\chi$, that is, 
$\T= fl(\chi)$.
If, in addition, $\T$ lies above $\VPS$, then $\chi$ lies above $\psi$. 
\end{lemma}

\begin{proof}
Let 
$$
(\mbg,  \mba,  \VP,  \mbn, \VPS ) =
(\mfg/K, \mfa/K, \Phi/K, \mfn/K, \Psi/K),
$$
 where $(\mfg, \mfa, \Phi, \mfn, \Psi)$
is a linear limit of $(G, A, \phi, N, \psi)$, and
$K =\Ker(\Phi)$. 
Then $\T$ inflates to a unique character $\theta \in \Irr(\mfg)$. 
Clearly $\theta$ lies above $\Phi$. 
If 
$(\mfg, \mfa, \Phi, \mfn, \Psi)= (G, A, \phi, N, \psi)$, i.e., the starting 
quintuple was already minimal, then  the lemma obviously holds with 
  $\chi = \theta$.
 
If $(\mfg, \mfa, \Phi, \mfn, \Psi)$ is a linear reduction of 
$ (G, A, \phi, N, \psi)$, that is we reach the limit quintuple after only one 
proper reduction, then $\mfg = G(\Phi)$. 
Hence Clifford's theorem implies that 
$\theta \in \Irr(\mfg| \Phi)$   induces irreducibly to $G$. 
Therefore the character $\theta^G = \chi$ is 
the only character in $\Irr(G | \phi)$ having $\theta$ as its 
$\Phi$-Clifford correspondent.
Hence $\chi= \theta^G $ is an irreducible character of  $G$ that lies above
$\phi$, since $\Phi$ is an extension of $\phi$.
It is also obvious that $l(\chi) = \theta$, while $fl(\chi) = \T$. 

If we need a series of linear reductions to reach the limit quintuple 
 $(\mfg, \mfa, \Phi, \mfn, \Psi)$, then we repeat Clifford's theorem 
 as many times  as the number of proper linear reductions we perform. 
In conclusion, the character $\theta^G = \chi$ is an irreducible character of 
$G$ that lies above $\phi$, and satisfies the conditions in  Lemma \ref{lim.l3}.

If, in addition, $\T$ lies above $\VPS$, then the inflation $\theta$ of $\T$ 
to $\mfg$, lies above $\Psi$.  Since $l(\psi) =\Psi$, the character 
$\Psi$ induces $\psi$ in $N$, i.e.,  $\Psi^N = \psi$. 
We conclude that $\theta^G= \chi$  lies above $\psi$ whenever 
$\theta$ lies above $\Psi$. This completes the proof of Lemma \ref{lim.l3}.
\end{proof}

The following is  straight forward.
\begin{lemma}\mylabel{da.ll0}
Assume that $(G, A, \phi, N, \psi)$  and $(H, B, \beta, M, \mu)$  are 
two linear quintuples. Assume further that there exists an epimorphism 
of linear quintuples 
$$
\rho: (G, A, \phi, N, \psi) \to (H, B, \beta, M, \mu).
$$
By this we mean that $\rho$ is an epimorphism 
of the group $G$ onto $H$  sending $A$ onto $B$ and $N$ onto $M$.
Furthermore, 
$\phi = \beta \circ \rho_A$ and $\psi = \mu \circ  \rho_N$.
The restriction of $\rho$ to any linear reduction (or multiple 
linear reduction) $(G', A', \phi', N', \psi')$ of  $(G, A, \phi, N, \psi)$ 
is an epimorphism  onto a linear reduction  (or multiple linear reduction 
respectively)
 $(H', B', \beta', M', \mu')$ of $ (H, B, \beta, M, \mu)$.
In this way $\rho$  induces  a one to one correspondence between 
linear reductions ( or multiple linear reductions) of  
 $(G, A, \phi, N, \psi)$  and linear reductions (respectively 
multiple linear reductions) 
of $(H, B, \beta, M, \mu)$.
Hence $\rho$ induces a one to one correspondence 
between linear limits of  $(G, A, \phi, N, \psi)$ and linear limits of 
$(H, B, \beta, M, \mu)$.
\end{lemma}

\begin{proof}
If   $(G', A', \phi', N', \psi')$ is a linear reduction of 
 $(G, A, \phi, N, \psi)$, it is easy to check  that its image
$\rho((G', A', \phi', N', \psi'))$
 under  $\rho$ is a linear reduction of $(H, B, \beta, M, \mu)$.

Let  $(H', B', \beta', M', \mu')$ be a linear reduction of
 $(H, B, \beta, M, \mu)$. If $G', A'$ and $N'$ are  the inverse images,
under $\rho$, of $H', B'$ and $M'$, respectively, then the quintuple
  $(G', A', \beta' \circ  \rho_A', N',  \mu' \circ   \rho_N')$ is a linear 
reduction of  $(G, A, \phi, N, \psi)$,  and its  image under  $\rho$ 
equals  $(H', B', \beta', M', \mu')$. 
We conclude that there exists a one to one correspondence between 
linear reductions of $(G, A, \phi, N, \psi)$  
and linear reductions of  $(H, B, \beta, M, \mu)$.

Because a multiple linear reduction is reached after a series 
of  linear reductions, repeated applications
 of the one  to one correspondence on  linear reductions  implies 
the existence of a one to one correspondence between 
multiple linear reductions of $(G, A, \phi, N, \psi)$  and 
$(H, B, \beta, M, \mu)$. Furthermore, since  
any  linear limit of   $(G, A, \phi, N, \psi)$  is a
minimal multiple linear reduction of the latter quintuple,
we also get a one to one correspondence between the 
linear limits of $(G, A, \phi, N, \psi)$  and those of 
$(H, B, \beta, M, \mu)$.
Hence  the lemma holds.
\end{proof}

\begin{corollary}\mylabel{lico1}
Assume that the linear quintuples 
  $(G, A, \phi, N, \psi)$  and $(H, B, \beta, M, \mu)$ satisfy the
 hypothesis in Lemma \ref{da.ll0}. 
Then any faithful linear limit of  $(H, B, \beta, M, \mu)$
is isomorphic to a faithful linear limit of $(G, A, \phi, N, \psi)$.
\end{corollary} 

\begin{proof}
Let  $(H', B', \beta', M', \mu')$ be a linear limit of $(H, B, \beta, M, \mu)$.
Then according to Lemma \ref{da.ll0} it corresponds to a linear limit 
$(G', A', \phi', N', \psi')$  of $(G, A, \phi, N, \psi)$ . Because $\rho$ maps 
the latter linear limit onto the former,  we get that 
$\rho$ maps $ A'$ onto  $ B'$,   while 
$\phi' = \beta' \circ  \rho_A'$.
We conclude that the kernel $\Ker(\phi')$ 
of $\phi'$ is mapped, under $\rho$,  onto the kernel $\Ker(\beta')$
 of $\beta'$, i.e., $\rho( \Ker(\phi')) = \Ker(\beta')$.
Furthermore, 
if $S$ equals the kernel of $\rho_{G'}$  then $S$ is a normal
 subgroup of $G'$ that  is contained in $\Ker(\phi')$, 
(since for all $s \in S$  we get $\phi' (s) = 
\beta'(\rho(s)) = \beta'(1) = 1$).
Hence the following holds
\begin{gather*}
G'/ \Ker(\phi') \cong H'/ \Ker(\beta'), \\
A'/ \Ker(\phi') \cong B'/ \Ker(\beta'), \\
N'/ \Ker(\phi') \cong M'/ \Ker(\beta').
\end{gather*}
In addition, the unique  characters $\phi'/ \Ker(\phi')$ and 
$\psi'/ \Ker(\phi')$ of the factor groups 
 $ A'/ \Ker(\phi')$  and  $N'/ \Ker(\phi')$  that inflate  to 
$\phi'$  and  $\psi'$, respectively, correspond under the above 
isomorphisms,  to the unique characters $\beta'/ \Ker(\beta')$ and 
$\mu'/ \Ker(\beta')$ of the factor groups 
 $ B'/ \Ker(\beta')$  and  $M'/ \Ker(\beta')$  that inflate  to 
$\beta'$  and  $\mu'$, respectively.
This completes the proof of  the corollary.
\end{proof}

\begin{proposition}\mylabel{lim.p2}
Let $(G, A, \phi, N, \psi)$ be a linear quintuple, and $T= \Ker(\phi)$. 
Then the factor quintuple $(G/T, A/T, \phi/T, N/T, \psi/T)$ is well defined. 
Furthermore,  any faithful linear limit of the factor quintuple is isomorphic 
to  a faithful linear limit
of the original one. 
\end{proposition}

\begin{proof}
Observe that the natural epimorphism  from $G $ to $ G/T$ 
provides an epimorphism 
$$
\rho: (G, A, \phi, N, \psi) \to (G/T, A/T, \phi/T, N/T, \psi/T)
$$
of linear quintuples, so that $(G, A, \phi, N, \psi)$ and 
 $  (G/T, A/T, \phi/T, N/T, \psi/T)$  satisfy the hypothesis in Lemma \ref{da.ll0}.
The rest of  proof is a simple application of Lemma \ref{da.ll0}
 and Corollary \ref{lico1}. 
\end{proof}

From now on, whenever needed, we will identify  any  two  
faithful linear limits  of  $(G, A, \phi, N, \psi)$ and 
$(G/T, A/T, \phi/T, N/T, \psi/T)$ which are isomorphic under the preceding
 proposition.

\begin{corollary}\mylabel{lim.co1}
Assume that the linear quintuple  $(G', A', \phi', N', \psi')$ is a multiple 
linear reduction of $(G, A, \phi, N, \psi)$, and let  $T= \Ker(\phi')$.
Then  the factor quintuple 
 $(G'/T, A'/T, \phi'/T, N'/T, \psi'/T)$ is well defined. 
Furthermore,  any faithful linear limit of the factor quintuple is, under 
some identification,  
also a faithful linear limit
of $(G, A, \phi, N , \psi)$,  i.e., 
$$
FLL(G'/T, A'/T, \phi'/T, N'/T, \psi'/T) \leq 
FLL(G, A, \phi, N, \psi).
$$
\end{corollary}
\begin{proof}
Follows immediately from Remark \ref{lim.r1} and Proposition \ref{lim.p2}.
\end{proof}

For the next proposition we  will need a  nice observation of I.M.Isaacs, that is actually 
the exercise (6.11) in \cite{is}.   
\begin{lemma}\mylabel{lim.l1}
Let $B$  be a  normal subgroup of a finite group $G$, 
$\gamma \in \Lin(B)$ a linear character of $B$ and $\chi \in 
\Irr(G|\gamma)$ an irreducible  character of $G$ lying above
$\gamma$.
If $\chi_{\gamma} \in \Irr(G(\gamma))$ is the $\gamma$-Clifford correspondent 
of $\chi$ in the stabilizer $G(\gamma)$ of $\gamma$ in $G$, then 
$\chi$ is monomial if and only if 
$\chi_{\gamma}$ is  monomial.
\end{lemma}

\begin{proof}
It is clear that if $\chi_{\gamma}$ is monomial then $\chi$ is monomial, as
$\chi_{\gamma}$ induces $\chi$ in $G$. 

So we assume that $\chi$ is a monomial character, and we will show that 
$\chi_{\gamma}$  is also monomial. 
Let $K= \Ker(\chi)$.  Of course  $K \unlhd G$.
It is clear that $\chi_{\gamma}$ is monomial if and only if the irreducible
character  $\chi_{\gamma}/ K$ of the factor group $G(\gamma)/K$ that
inflates to $\chi_{\gamma}$, is monomial.
Hence it  suffices to prove the lemma 
in the case of a faithful irreducible character $\chi$, as we can pass to 
the factor groups $G/K$ and $(BK)/ K$. 
So in the rest of the proof we assume that $K=1$.

Clifford's Theorem implies that the restriction  $\chi|_B$ of $\chi$ to
$B$  is a sum of $G$-conjugates of
$\gamma$. Thus  $1=\Ker(\chi|_B)= \bigcap_{s \in
G/G(\gamma)}(\Ker(\gamma^s))$.
But the derived group $[B, B]$  of $B$ is contained in the kernel of
$\gamma^s$ for every $s\in G$, as $\gamma $ is linear. 
Thus $[B, B] \leq \Ker(\chi|_B)=1$. So $B$ is abelian.

We can now  follow the hint of problem 6.11 in \cite{is}. 
As $\chi$ is monomial,  there exists $H \leq G$ and $\lambda \in \Lin(H)$ 
with $\chi= \lambda^G$. Thus the irreducible character $\lambda^{HB}$ of
$HB$ lies above a $G$-conjugate  $\gamma^s$ of $\gamma$, where $s \in
G $. As the $G$-conjugate  $\lambda^{s^{-1}} \in
\Lin(H^{s^{-1}})$  of $\lambda$ also induces $\chi$, 
we can replace $H$ by $H^{s^{-1}}$ and $\lambda$ by $\lambda^{s^{-1}}$.
This way $\lambda^{HB} $ is replaced 
by $(\lambda^{s^{-1}})^{H^{s^{-1}} B}=(\lambda^{HB})^{s^{-1}}$, which lies above $\gamma$.

According to 
  Mackey's Theorem   
\begin{equation}\mylabel{lim.e1}
\lambda^{HB}|_B = (\lambda|_{H\cap B}) ^B. 
\end{equation}
As $B$ is abelian, the right hand side of
 \eqref{lim.e1} equals the sum of $|B:H \cap B|$ 
distinct character extensions of $\lambda|_{H \cap B}$ to $B$, 
each one appearing with multiplicity one.
Thus every irreducible  constituent of   
$\lambda^{HB}|_B $ appears with multiplicity one.
This, along with Clifford's  theorem, (as 
 $\lambda^{HB}$ lies above $\gamma$),
implies that 
$$
\lambda^{HB}|_B= e\cdot\sum_{s \in  S}
 \gamma^s=\sum_{s \in S} \gamma^s ,
$$
where $S$ is a family of representatives for the cosets $H(\gamma)Bs$ of 
$H(\gamma) B = (HB)(\gamma)$ in $HB$, and 
$e$ is a positive integer.
Furthermore, Clifford's theorem implies the existence of
 an irreducible character
$\theta \in \Irr(HB(\gamma))$ lying above $\gamma$ 
and inducing $\lambda^{HB}$. The fact that $e=1$ 
implies that $\theta|_{B}= \gamma$, i.e., 
$\theta \in \Irr(HB(\gamma))$ is an extension 
of $\gamma \in \Irr(B)$ to $HB(\gamma)$.
Thus $\theta \in \Lin(HB(\gamma)| \gamma)$
  induces $\lambda^{HB}$. Hence $\theta^G= \chi$, as 
$\lambda$ induces $\chi$. Therefore, $ \theta^{G(\gamma)}$ 
is an irreducible character of $G(\gamma)$ 
lying above $\gamma$ and inducing $\chi$.  
As the $\gamma$-Clifford correspondent $\chi_{\gamma}$  of $\chi$
is unique, we conclude that  $\theta^{G(\gamma)} = \chi_{\gamma}$.
Hence $\chi_{\gamma}$ is  induced from the linear character  $\theta$,
 and thus is monomial.  

This completes the proof of the lemma in the case of an abelian $B$.
So the lemma follows.
\end{proof}

\begin{proposition}\mylabel{lim.p1}
Let $(G, A , \phi, N , \psi)$ be a linear  quintuple.
Let    $\chi \in \Irr(G |\phi)$ an irreducible character of $G$ lying above
$\phi$,  $l(\chi) \in \Irr(l(G) |l(\phi))$be
a linear limit of $\chi$,  and $fl(\chi ) \in \Irr(fl(G) | fl(\phi))$
be the corresponding  faithful linear limit of $\chi$. 
Then the following are equivalent
\begin{itemize}
\item[{1)}]$\chi$  is monomial
\item[{2)}]$l(\chi)$ is monomial
\item[{3)}]$fl(\chi)$ is monomial
\end{itemize}
\end{proposition} 

\begin{proof}
Let  $(G', A' , \phi', N' , \psi')$ be  a linear  reduction of 
 $(G, A , \phi, N , \psi)$. 
According to Lemma \ref{lim.l1},  the character $ \chi \in \Irr(G|\phi)$
is monomial
if and only its $\phi'$-Clifford correspondent $\chi'$  is  monomial. 
This is true for every linear  reduction, so at the end we get that 
$\chi$ is monomial if and only if any   linear  limit $l(\chi)$ of $\chi$
is monomial. 

Let $fl(\chi) \in \Irr(fl(G))= \Irr(l(G)/K) $ be  the faithful linear limit of
$\chi$ corresponding to $l(\chi)$. 
It is obvious that $l(\chi)$ is monomial if and only if $fl(\chi)$ is
monomial.
This, along with the already proved first  equivalence, implies
that 
$fl(\chi)$ is monomial if and only if $\chi$ is monomial.
As this is true for any faithful linear limit $fl(\chi)$ of $\chi$, 
the proof of Proposition \ref{lim.p1} is complete.
\end{proof}

\section{Linear limits of characters of $p$-groups }
\mylabel{lp}

Assume  that $(G , A, \phi, N , \psi)$ is a linear quintuple.
For the rest of this section we suppose that $N$ is a $p$-group, for some
odd prime $p$.   
The main result of this section is 
\begin{theorem}\mylabel{lp.t1}
Suppose that $(G, A, \phi, N, \psi)$  is a linear quintuple with $N$ a
$p$-group, for some odd prime $p$. Assume further that 
$(\mbg, \mba, \VP, \mbn, \VPS)$ and $(\mbg', \mba', \VP', \mbn', \VPS')$ 
are two faithful linear limits of $(G, A, \phi, N, \psi)$. 
Then both $\mbn/\mba$ and $\mbn'/\mba'$ are naturaly  symplectic
 $\ZZ_p(G(\psi)/N)$-modules.
Furthermore, $\mbn/\mba$ is isomorphic to $\mbn'/\mba'$ as a  symplectic
$\ZZ_p(G(\psi)/N)$-module.
\end{theorem}
To prove it we will use strongly Theorem 8.4 in \cite{da3}.
We remark here that the definitions of a ``stabilizer  limit'' and an 
``elementary stabilizer limit'' that were given in Sections 2 and 3, 
respectively, in \cite{da3}, are related  to
 but  not the same as our definition of linear limits.

We start with an elementary construction of  symplectic modules,
  and the associate notation.
Assume that a $p$-group   $R$ is a normal subgroup of some finite group $X$, 
where $p$ is an odd prime. Assume further that 
the center $Z(R)$  of $R$ is central in $X$, while  $R$ is a  central product
%\begin{subequations*}\mylabel{m.e10}
\begin{equation*}
R = E \odot Z(R), 
\end{equation*}
where  either $E$ is  $1$,  or else $E$ is  an extra
 special group of exponent $p$ and 
\begin{equation*}
E \cap Z(R) = Z(E).
\end{equation*}
%\end{subequations*}
Then $R/ Z(R)$ is an elementary abelian $p$-group, which may be trivial.
So $R/ Z(R)$, when written additively,  can  be considered  as a vector space 
over  the field $\ZZ_p$ of $p$ elements. This  way $R / Z(R)$ becomes a $\ZZ_p(X)$-module. 
Moreover, $[R, R] = [E, E] =   Z(E)$ is either trivial or a cyclic group of order $p$. 
If $Z(E)= 1$, i.e., $R / Z(R) = 1$, then $R/ Z(R)$ becomes trivially a symplectic 
$\ZZ_p(X)$-module.  If $|Z(E)|=p$,  then we can still make $R / Z(R)$ a symplectic $\ZZ_p(X)$-module. 
Indeed, if $\lambda \in \Irr(Z(E))$ is any faithful linear character of $Z(E)$, then we can    define
a  bilinear form $<, >$ from 
$(R/Z(R)) \times (R/Z(R))$ to the multiplicative group
$\CC_p$ of complex  $p$-roots of  unity
  as
\begin{equation}\mylabel{m.d1}
<\bar{x}, \bar{y}> = \lambda([x, y]) \in \CC_p, 
\end{equation} 
for all $  x,y \in R$, 
where $\bar{x}$ denotes the image of $x \in R$   in the factor group 
$R / Z(R)$, and 
$[x,y]$  is the commutator of $x$ and $ y$ in $R$.
Note that,  as the multiplicative group 
$\CC_p$ of $p$-roots  of unity is isomorphic 
to the additive  group  $\ZZ_p^+$  of $\ZZ_p$, we can  
identify these   two isomorphic groups, and consider 
the bilinear form  $< , >$ as a 
symplectic  form in $\ZZ_p$. 
As $Z(R)$ is central in $X$,   this   form is $X$-invariant. 
So  $R/Z(R)$ is a symplectic $\ZZ_p(X)$-module.

We assume that $R$ and $ X$ are as above, 
with  $Z(R)$  central in $X$. 
 Let $U$ be a subgroup of $R$ containing $Z(R)$, 
and normal  in $X$. 
Then,  (see the notation in \cite{da}), 
  we call the symplectic $\ZZ_p(X)$-submodule 
$U /Z(R)$ of $R /Z(R)$  {\bf isotropic } 
if $U\leq R$ is an  abelian subgroup of $R$.
 We call $U /Z(R)$ { \bf anisotropic  } 
if $0$ is its only   isotropic
$\ZZ_p(X)$-submodule, i.e., every abelian  subgroup of
 $U$ which is normal in  $X$ is contained in $Z(R)$.
Observe that $0$ is an anisotropic symplectic $\ZZ_p(X)$-module.

Now we go back to the linear quintuple 
$(G, A, \phi, N, \psi)$. 
Let $(\mfg, \mfa, \Phi, \mfn, \Psi)$ be a linear limit of $(G, A, \phi,
 N, \psi)$, and 
\begin{equation}\mylabel{lp.0}
 (\mbg, \mba, \VP, \mbn, \VPS)=(\mfg/K, \mfa/K, \Phi/K, \mfn/K, \Psi/K),
\end{equation}
be the corresponding faithful linear limit of $(G, A, \phi, N, \Psi)$, 
where $K= \Ker(\Phi)$.  Assume further that 
 $(G_i , A_i , \phi_i, N_i , \psi_i)$,
is  a chain of linear quintuples, for all $i=0, \dots, n$, 
such that
\begin{subequations}\mylabel{lp.e1}
\begin{gather}
(G_0, A_0, \phi_0, N_0 , \psi_0)= (G, A, \phi, N, \psi), \\
(G_n, A_n, \phi_n, N_n , \psi_n)= (\mfg, \mfa, \Phi, \mfn, \Psi),
\text{ and }  \\
(G_{i}, A_{i}, \phi_{i}, N_{i} , \psi_{i}) \text{ is a
proper  linear reduction of }
(G_{i-1}, A_{i-1}, \phi_{i-1}, N_{i-1} , \psi_{i-1}),
\end{gather}
whenever $i=1,\dots,n$.
\end{subequations}
 These objects stay fixed until the end of the section.
We also keep fixed an arbitrary  $\CC G(\psi)$-elementary stabilizer limit
$\Lambda \in ESL(\psi| \CC G(\psi))$ of $\psi$,  in the sense of \cite{da3},
and in particular (3.7) of that paper. 
(Note that the ordered triple $(G(\psi), N, \psi)$ is a member of the family 
defined in (2.1) of \cite{da3}. Thus we can define a $\CC G(\psi)$-elementary 
stabilizer  limit of $\psi$.)

We start with some results following \eqref{lp.e1}.
\begin{lemma}\mylabel{lp.l1}
Let $M$ be a subgroup of $G$ with $A_i \leq M $, for some $i=0,1, \dots, n$.
Assume further that an  
irreducible character $\chi \in \Irr(M)$, when restricted to 
$A_i$, is a multiple of $\phi_i$. 
Then   $G(\chi) = G(\chi, \phi_1, \dots, \phi_i) = G_i(\chi)$.
In particular,   $G(\phi_i)= G_i$, for all $i=0, 1, \dots, n$. 
\end{lemma}
\begin{proof}
Clearly (\ref{lp.e1}c) implies that $G_i= G(\phi_0, \phi_1, \dots, \phi_i)$ for every $i=0, 1,\dots, n$.
For all $i=1, \dots, n$, the linear character $\phi_i \in \Irr(A_i)$ is an extension of $\phi_{i-1} \in \Irr(A_{i-1})$. 
Even more,  for all such $i$ the group $A_i$ is a normal subgroup of $G(\phi_1, \dots, \phi_{i-1})= G_{i-1}$.
 
Assume that $i=0, 1,\dots,n$ is fixed. Let $M \geq A_i$  and $\chi \in \Irr(M)$ 
with $\chi |_{A_i} = m \phi_i$  for some integer $m \geq 1$. 
For any $j=0, 1,\dots, i-1$, the linear character 
 $\phi_i$ is an extension of  $\phi_j \in \Irr(A_j)$ to $A_i$. 
Hence $\chi|_{A_j} = m \phi_j$ for all $j=0,1, \dots, i$.
Therefore  $G(\chi)$ fixes any such $\phi_j$  if and only if it normalizes $A_j$. 
Clearly $G(\chi)$ fixes the $G$-invariant character $\phi_0$ of $A_0=A$. 
It also   normalizes $A_1$, as the latter is normal in $G$. Thus 
$G(\chi)$  fixes $\phi_1$. 
Suppose  now that  $G(\chi)$ normalizes $A_1, \dots, A_{j-1}$, where $j=2, \dots, i$.
Then it fixes $\phi_1, \dots, \phi_{j-1}$, i.e., 
$G(\chi)  \leq G(\chi)(\phi_0, \phi_1, \dots, \phi_{j-1})$.
But $A_j$ is a normal subgroup of $G(\phi_0, \phi_1, \dots, \phi_{j-1}) = G_{j-1}$.
Hence $G(\chi)$ normalizes $A_j$ as well, and therefore 
 also fixes $\phi_j$ . 
As this holds for all $j=1, \dots, i$, the 
first statement of the  lemma follows.

The second part of the lemma follows from the first, if we take  $\chi= \phi_i $ and $M= A_i$.
\end{proof}

\begin{proposition}\mylabel{lp.p0}
For every $i=0, 1, \dots, n$, we have 
\begin{equation}\mylabel{lp.e0}
G(\psi_i)= G_i(\psi_i) = G_i(\psi) \text{ and }
G(\psi_i) N = G(\psi). 
\end{equation}
Hence $G(\Psi) = \mfg(\Psi)= \mfg(\psi)$ and $G(\Psi)N = G(\psi)$.
\end{proposition}

\begin{proof}
As $\mfg=G_n$ and  $\Psi= \psi_n$, it suffices to prove \eqref{lp.e0}.
For this proof  we will use induction on $i$.  
As $G_0=G$ and $\psi_0= \psi$, the equations in \eqref{lp.e0} hold 
trivially  for $i=0$.
Suppose \eqref{lp.e0} is true for all
 $i=0, \dots, t-1$,  where  $t=1, \dots, n$.
 We will show it holds for $i=t$. 

By (\ref{lp.e1}c), both groups  $A_t$ and $N_{t-1}$ are normal subgroups 
of $G_{t-1}$.  Furthermore, $\phi_t \in \Irr(A_t)$ is a linear extension 
of $\phi_{t-1} \in \Irr(A_{t-1})$ and lies  under $\psi_{t-1}$. In addition, 
$\psi_t \in \Irr(N_t)$ is the $\phi_t$-Clifford  correspondent of 
 $\psi_{t-1} \in \Irr(N_{t-1})$.  Hence 
\begin{equation}\mylabel{lp.e2}
G_{t-1}(\psi_t)= G_{t-1}(\psi_{t-1}, \phi_t).
\end{equation}
As $G_t$ is the subgroup $G_{t-1}(\phi_t)$ of $G_{t-1}$,
 both sides of this equation are equal to 
\begin{equation}\mylabel{lp.e3}
G_t(\psi_{t-1}) = G_t(\psi_t).
\end{equation}
Furthermore, any element of $G_{t-1}$ that fixes $\psi_{t-1}$ permutes among 
themselves the $N_{t-1}$-conjugates of $\phi_t$, as $A_t$ is normal 
in $G_{t-1}$. We conclude that $G_{t-1}(\psi_{t-1}) = 
G_{t-1}(\psi_{t-1}, \phi_t) N_{t-1}$. 
This, along with \eqref{lp.e2}, implies
\begin{equation}\mylabel{lp.e4}
G_{t-1}(\psi_{t-1})= G_{t-1}(\psi_{t-1}, \phi_t)N_{t-1}=
G_{t-1}(\psi_t) N_{t-1}. 
\end{equation}
Hence 
\begin{align*}
G_t(\psi_t) &= G_t(\psi_{t-1})  &\text{ by \eqref{lp.e3} }\\
&=G_{t-1}(\phi_t, \psi_{t-1} )  &\text{ as $G_t = G_{t-1}(\phi_t)$ } \\
&=G_{t-1}(\phi_t, \psi) &\text{ by induction for $i=t-1$ }\\
&=G_t(\psi). &
\end{align*}
Lemma \ref{lp.l1} implies that $G(\psi_t)= G_t(\psi_t)$,  as  $\psi_t|_{A_t}$ is a multiple 
of $\phi_t$. We conclude that 
\begin{equation}\mylabel{lp.e5}
G_t(\psi_t)=G_t(\psi)=G(\psi_t).
\end{equation}
Hence the first part of \eqref{lp.e0} follows for  the inductive step.
For the second part we get 
\begin{align*}
G(\psi) &= G (\psi_{t-1})N  &\text{ by induction for $i=t-1$ }\\
&=G_{t-1}(\psi_{t-1})N   &\\
&=G_{t-1}(\psi_{t-1},\phi_t)N_{t-1}N  &\text{ by \eqref{lp.e4} }\\
&=G_{t-1}(\psi_{t-1}, \phi_t)N &\\
&=G_t(\psi_{t-1}) N  &\text{ as $G_t= G_{t-1}(\phi_t)$ } \\
&=G_t(\psi_t)N    &\text{ by \eqref{lp.e3}} \\
&=G(\psi_t)N  &\text{  by \eqref{lp.e5}} 
\end{align*}
This completes the inductive  proof of \eqref{lp.e0} for $i=t$. 
Hence Proposition \ref{lp.p0} follows.
\end{proof}

\begin{corollary}\mylabel{lp.c1}
The inclusion $\mfg(\Psi) \to G(\psi)$ induces an isomorphism 
  of $\mfg(\Psi)/\mfn $ onto  $G(\psi)/N$, where  $\bar{s}$ maps to 
$ \bar{s}N$, for any  $\bar{s} \in \mfg(\Psi)/\mfn $.
\end{corollary}

\begin{proof}
Obvious, as $\mfg(\Psi)N = G(\Psi)N = G(\psi)$, while
$\mfg(\Psi) \cap N = \mfg  \cap N = \mfn$.
\end{proof}

\begin{corollary}\mylabel{lp.c0}
The character $\Psi \in \Irr(\mfn)$ is $\mfg$-invariant.
Furthermore,   $\Psi$ is zero on $\mfn - \mfa$,
 and it is   a multiple of $\Phi$ on $\mfa$.
  Hence
$$
\mfg= \mfg(\Psi)=  \mfg(\psi)= G(\Psi).
$$
\end{corollary}

\begin{proof}
  According to \eqref{lp.0}, we have $\mbg= \mfg/K, \mbn = \mfn/K$ 
and $\mba = \mfa/K$,  
where $K = \Ker(\Phi)$. Furthermore, 
$\VPS$ is the unique character of the factor group  $\mfn/K$  that inflates to 
$\Psi \in \Irr(\mfa)$. But as $N$ is a $p$-group, 
  Corollary \ref{lim.co22} implies that  the character $\VPS \in \Irr(\mbn)$ 
is $\mbg$-invariant. Hence $\Psi$ is $\mfg$-invariant.
The same corollary implies that    $\VPS \in \Irr(\mbn)$ 
vanishes outside $\mba$ and is a multiple of $\VP$ on $\mba$. 
Thus a similar  property holds for its unique inflation $\Psi \in \Irr(\mfn)$.
The rest of the corollary follows easily from Proposition \ref{lp.p0}.
\end{proof}

The next result follows immediately from the above two corollaries. 
\begin{corollary}\mylabel{lp.c4}
The isosmorphism of $G(\psi)/N$ onto $\mfg(\Psi)/\mfn$ in Corollary 
\ref{lp.c1}, composed with the natural isomorphism of 
$\mfg / \mfn= \mfg(\Psi)/ \mfn$ onto  $\mbg / \mbn $, provides  
an isomorphism $j$ from $G(\psi)/ N$  onto the factor group 
$\mbg/ \mbn$. So  any coset $t \in G(\psi)/ N$ gets mapped under 
$j$, to the image of the coset $(t \cap \mfg(\Psi))/ \mfn$
under the natural isomorphism of $\mfg/ \mfn$ onto 
$\mbg/ \mbn$. 
\end{corollary}

\begin{proposition}\mylabel{lp.p1}
The factor group $\mbn/\mba$
is an anisotropic symplectic $\ZZ_p(\mbg/\mbn)$-module, which may be  $0$, 
 with respect to the bilinear $\mbg$-invariant 
  form defined, as in  \eqref{m.d1}, by
\begin{equation}\mylabel{lp.form1}
< x\mba, y\mba > = \VP([x, y]) \in \CC_p \cong \ZZ_p^+, 
\end{equation}
for any $x , y \in \mbn$. Here the $\mbg/ \mbn$-action is induced 
by conjugation in $\mbg$. 
\end{proposition}

\begin{proof}
As $N$ is a $p$-group, Lemmas \ref{lim.co2} and \ref{lim.co22}
imply that either $\mbn= \mba$ or  $\mbn$ is the central product 
$\mbn = E \odot Z(\mbn) = E \odot \mba$, of a non-trivial extra special 
$p$-group $E$ of exponent $p$, and $\mba = Z(\mbn)$ which is central in $\mbg$.
Furthermore, $\VP \in \Irr(\mba)$ is a faithful linear character of 
$\mba$.         In both cases, 
 the factor group $\mbn/\mba$, 
 becomes a  symplectic $\ZZ_p(\mbg)$-module, where 
$\mbg$ acts on $\mbn/\mba$ by conjugation,
 and the symplectic form is defined 
via commutation in $\mbn$, 
 (see \eqref{m.d1} and the paragraph that follows it). 
In addition,  Lemma \ref{lim.co2} implies that $\mba$   is maximal 
among the abelian subgroups of $\mbn$ which are normal in 
$\mbg$. 
Hence $\mbn/\mba$ is an anisotropic symplectic $\ZZ_p(\mbg)$-module, which may be $0$.

Clearly $\mbn$ centralizes the factor group $\mbn/\mba$, as $[\mbn, \mbn]
\leq \mba$. 
Hence the action of $\mbg$ on $\mbn/\mba$ induces one of $\mbg/\mbn$ 
on that symplectic group.
So Proposition \ref{lp.p1} follows.
\end{proof}

Proposition \ref{lp.p1}, along with the isomorphism $j$ defined in Corollary 
\ref{lp.c4}, implies  
\begin{corollary}\mylabel{lp.c3}
The factor group $\mbn/ \mba$ is an anisotropic  symplectic 
$\ZZ_p(G(\psi)/N)$-module, with respect to the bilinear form defined in 
\eqref{lp.form1}.  Here the $G(\psi)/ N$-action is defined as through the 
isomorphism $j$, defined  in Corollary \ref{lp.c4} as 
$$
(x\mba)^{\bar{s}}= x^ {j(\bar{s})} \mba \in \mbn/ \mba  
$$
for all $x \in \mbn$ and $\bar{s}  \in G(\psi)/ N$.
\end{corollary}

\begin{corollary}\mylabel{lp.c2}
The factor group  $\mfn/\mfa$ 
is an anisotropic symplectic $\ZZ_p(\mfg/\mfn)$-module, that may be $0$,  with respect
 to the bilinear form defined by 
\begin{equation}\mylabel{lp.form2}
< s\mfa, t\mfa > = \Phi([s, t]) \in \CC_p \cong  \ZZ_p^+, 
\end{equation}
for any $s , t \in \mfn$. Here the $\mfg/\mfn$-action is induced by 
conjugation in $\mfg$. 
Hence with respect to the same form, $\mfn/\mfa$  is  an anisotropic symplectic
$\ZZ_p(G(\psi)/N)$-module, where 
  the action of $G(\psi)/N \cong \mfg(\Psi)/\mfn = \mfg/\mfn $ is defined by 
$$
(s\mfa)^{rN}= (s^r)\mfa \in \mfn/\mfa,
$$
for any $s, t \in \mfn$ and $r \in \mfg$.
\end{corollary}

\begin{proof}
Let $K =\Ker(\Phi)$. Then, as we have already seen, 
\begin{equation}\mylabel{lp.e6}
\mfn/\mfa \cong (\mfn/K)/ (\mfa/K) = \mbn/\mba,
\end{equation}
where the isomorphism is $\mfg$-invariant.
Furthermore, $\VP$ is the unique character of the factor group 
$\mfa/K = \mba$ that inflates to $\Phi \in \Irr(\mfa)$.
Hence, under the  isomorphism in \eqref{lp.e6},
  Proposition \ref{lp.p1} implies that 
 $\mfn/\mfa$ is  an anisotropic 
symplectic $\ZZ_p(\mfg/\mfn)$-module, with respect to the bilinear form 
that \eqref{lp.form1} determines. (Note that this bilinear form  
translates to \eqref{lp.form2}.)
Furthermore, $\mfg$ acts on $\mfn/\mfa$ by conjugation, while $\mfn$ 
centralizes it.
But $\mfg = \mfg(\Psi) = G(\Psi)$, by  Corollary \ref{lp.c0}, while 
 Corollary \ref{lp.c1} implies that 
$G(\psi)/N$  is naturally isomorphic to $\mfg(\Psi)/\mfn  = \mfg/\mfn$.
Hence $\mfn/\mfa$ becomes  an anisotropic $G(\psi)/N$-module. 
This completes the proof of the corollary.
\end{proof}

\begin{proposition}\mylabel{lp.p2}
If $N$ is a $p$-group, then $\Psi$  is a $\CC G(\psi)$-stabilizer limit 
of $\psi$ as this is  defined in \cite{da3}, that is, 
$\Psi \in \text{ SL}(\psi | \CC G(\psi))$. 
\end{proposition}

\begin{proof}
According to (\ref{lp.e1}c), for every $i=1, \dots, n$, we have a normal 
subgroup $A_{i}$ of $G_{i-1}$ contained in $N_{i-1}$, and a linear 
character $\phi_i \in \Irr(A_i)$ lying under $\psi_{i-1} \in \Irr(N_{i-1})$.
Furthermore,  $\psi_i \in \Irr(N_i)$ is the $\phi_i$-Clifford correspondent 
of $\psi_{i-1}$. As the ordered triple $(G_{i-1}(\psi_{i-1}), N_{i-1},
 \psi_{i-1})$ is a member of the family defined in (2.1) of \cite{da3},
while  $A_{i}$ is a normal subgroup of $G_{i-1}(\psi_{i-1})$, we conclude that 
$\psi_{i}$  is an element of the set $\DCC(\psi_{i-1} |\CC G_{i-1}(\psi_{i-1}))$
defined in (2.2) of \cite{da3}. 
According to Proposition  \ref{lp.p0}
 we have $G_{i-1}(\psi_{i-1})= G(\psi_{i-1})$.
Hence we get a sequence of characters $\psi= \psi_0, \psi_1, \dots, 
\psi_n= \psi$, such that 
$$
\psi_{i} \in \DCC(\psi_{i-1} | \CC G(\psi_{i-1}) ),
$$
for all $i=1, \dots, n$.
Hence $\Psi$ lies in the set $\text{CC}(\psi| \CC G(\psi))$,
 (see (2.3) in \cite{da3}).

According to the definition of stabilizer limits, in 
 (2.16)  of  \cite{da3}, we can complete the proof of the proposition by 
showing that $\Psi \in \Irr(\mfn)$
is the only member of  $\DCC (\Psi |\CC G(\Psi))$.
By (2.14) and (2.15) in \cite{da3}, it  
suffices to show that whenever $M$ is a normal subgroup of $G(\Psi)$
contained in $\mfn$, 
 the restriction $\Psi|_{M}$ is a multiple of 
of a single irreducible character.
 Suppose such an $M \leq \mfn$ is fixed. Let 
$\theta \in \Irr(M)$  be an irreducible character of $M$ that lies under 
$\Psi$. It is enough to show that $\theta$ is $G(\Psi)$-invariant. 
We know from Corollary \ref{lp.c0} that $G(\Psi)=\mfg(\Psi) = \mfg$. 
So it suffices to show that $\mfg(\theta)= \mfg$.

$\Psi$ lies above the $\mfg$-invariant linear character 
$\Phi  \in \Irr(\mfa)$. Hence  we can replace $M$ with $M\cdot \mfa$ and
 $\theta \in \Irr(M)$ with $\theta \cdot \Phi \in \Irr(M\mfa)$ 
(where $(\theta \cdot  \Phi)(ma) = \theta(m) \Phi(a)$, for all $m \in M$ and $a \in \mfa$ ). This way $\mfg(\theta)=
 \mfg(\theta \cdot   \Phi)$
remains  the same. So we may assume that $\mfa \leq M \leq \mfn$, and that 
$\theta$ lies  above $\Phi$. Then $M / \mfa$ is a $\ZZ_p(\mfg/\mfn)$-submodule
of $\mfn /\mfa$. But the latter is an anisotropic $\ZZ_p(\mfg/\mfn)$ module by 
Corollary \ref{lp.c2}. 
Hence its symplectic form $<, >$ (see \eqref{lp.form2}), restricts to a 
non-singular  bilinear alternating  form on $(M/\mfa ) \times (M/\mfa)$. 
It follows that $\theta$ is zero on $M - \mfa$ and a multiple of 
$\Phi$ on $\mfa$. 
Therefore $\mfg(\theta)= \mfg$, and the proposition follows. 
\end{proof}

According to (2.12) in \cite{da3}, we may define another triple, 
denoted by $(G(\psi)\{\Psi\}^*, N\{\Psi\}^*, \Psi^*)$,   using  
 the $\CC G(\psi)$-stabilizer limit $\Psi$  of $\psi$. The star groups 
are  defined   in (2.12) of  \cite{da3}, 
 as the  factor  groups we get when we divide the triple
$(G(\psi)\{\Psi\}, N\{\Psi\}, \Psi)$ by $\Ker(\Psi)$. 
So $\Psi^*$ in \cite{da3} denotes 
 the unique character   $\Psi/ \Ker(\Psi)$ from which 
$\Psi$ is inflated.
  Note also that $X\{\theta \}$ denotes in \cite{da3}
the stabilizer $X(\theta)$ of $\theta$ in $X$, for any group $X$ and any irreducible character 
$\theta$ of any  subgroup of $X$.
In our case, where $N$ is a $p$-group, the kernel $\Ker(\Psi)$ of $\Psi$
 coincides with  $K = \Ker(\Phi)$, by Corollary \ref{lp.c0}. 
Furthermore, the same corollary implies that 
$G(\psi)\{\Psi\}= G(\psi)(\Psi)= \mfg(\Psi)= \mfg$. 
Of course $N\{\Psi\}= N(\Psi)= \mfg(\Psi) \cap N = \mfg \cap N = \mfn$.
Hence the star triple  $(G(\psi)\{\Psi\}^*, N\{\Psi\}^*, \Psi^*)$ in \cite{da3}, 
is what we write as $(\mbg, \mbn , \VPS)$ (see \eqref{lp.0}). 

Even more, according to (2.13a) in \cite{da3} the stabilizer limit 
$\Psi$ of $\psi$ defines a natural isomorphism denoted by $\cdot / \Psi$ 
from $G(\psi)/ N$ to $G(\psi, \Psi)^*/ N(\Psi)^* = \mbg/ \mbn$.
Observe that this is exactly the isomorphism $j$ defined in Corollary 
\ref{lp.c4}.
Having explained this, we can now prove
\begin{theorem}\mylabel{lp.t2} 
Let  $\Lambda \in \text{ ESL}(\psi |\CC G(\psi))$ 
be a $\CC G(\psi)$-elementary 
stabilizer limit of $\psi$, with   $K_0= \Ker(\Lambda)$
and $N(\Lambda)^* = N(\Lambda)/K_0$.
Then $N(\Lambda)^*  / Z(N(\Lambda)^*)$
 is isomorphic to $\mbn/\mba$ as symplectic $\ZZ_p(G(\psi)/N)$-modules.
\end{theorem}
\begin{proof}
We are going to apply Theorem 8.4 in \cite{da3}, for the triple 
$(G(\psi), N, \psi)$ here in the place of $(G, N, \psi)$ there,  the 
$\CC G(\psi)$-elementary stabilizer limit $\Lambda$ of  $\psi$ here,
 in the place of the $\CC G$-elementary stabilizer limit $\phi$ of 
$\psi$ there,  and the $\CC G(\psi)$-stabilizer limit $\Psi$ of $\psi$ here, 
in the place of $\theta$ there.
(Note that the hypotheses (7.1) and (7.2a) in \cite{da3}
 are satisfied.) Observe also that $\Lambda$ is an irreducible 
character of $N(\Lambda)$, by (2.4c) in \cite{da3}, as $\Lambda \in
 \text{ CC}(\psi | \CC G(\psi))$. 
Hence Theorem  8.4 gives us a monomorphism $\mu$ of the group 
$G(\psi)\{\Lambda \}^*=G(\psi, \Lambda)/K_0$ into 
the group $G(\psi)\{\Psi\}^*= \mbg$ that satisfies the equivalent 
of (6.1) in \cite{da3}.
 Furthermore, the relations (8.5)in \cite{da3} tell us that
\begin{equation}\mylabel{lp.e7}
\mbg = \mba \;   \mu(G(\psi, \Lambda)^*)=
 \mba \;  \mu(G(\psi, \Lambda)/K_0) \text{ and } 
\mbn = \mba \; \mu(N(\Lambda)^*)= \mba \;  \mu(N(\Lambda)/K_0).
\end{equation}
(Note that in our case  $Z(N\{\Psi\}^*)= Z(\mbn) = \mba$.)
Furthermore, $\mu$ satisfies (6.1),  and,  in particular,  (6.1a), 
 of \cite{da3}. Hence 
the triple $(\mu(G(\psi, \Lambda)^*), \mu(N(\Lambda)^*), \mu(\Lambda^*))$ 
is a restrictor of  $(\mbg, \mbn , \VPS)$, in the sense of (5.1) in \cite{da3}.  (Where 
the irreducible  character $\Lambda^*$ is the unique character
 of the factor group $N(\Lambda)^*= N(\Lambda)/ K_0$ 
from which $\Lambda \in \Irr(N(\Lambda))$ is inflated, 
and $\mu(\Lambda^*) \in \Irr(\mu(N(\Lambda)^*))$
is the unique character of $\mu(N(\Lambda)^*)$
whose composition with $\mu$ is $\Lambda^*$.)
 Therefore (5.1) of \cite{da3} implies 
\begin{equation}\mylabel{lp.e8}
\mu(N(\Lambda)^*) = \mbn \cap \mu(G(\psi, \Lambda)^*)
\text{ and } 
\mu(\Lambda^*) = \VPS |_{\mu(N(\Lambda)^*)}.
\end{equation}
Hence $\mu$ restricts to an isomorphism  
\begin{equation}\mylabel{lp.e9}
N(\Lambda)/K_0 = N(\Lambda)^* \cong 
\mu(N(\Lambda)^*) = \mbn \cap \mu(G(\psi, \Lambda)^*)
\end{equation}
 that sends the irreducible character 
$\Lambda^* \in \Irr(N(\Lambda)^*)$ to the restriction of 
$\VPS$ to $\mu(N(\Lambda)^*)$. 
Even more, in view of  \eqref{lp.e7} we have 
 $\mba \cap \mu(N(\Lambda)^*)= Z(\mu(N(\Lambda)^*))$, and thus 
\begin{equation*}
%\mylabel{lp.e11}
\mbn / \mba \cong  \mu(N(\Lambda)^*)/ (\mba \cap \mu(N(\Lambda)^*))= 
\mu(N(\Lambda)^*)/ Z(\mu(N(\lambda)^*)). 
\end{equation*}
 According to 
 \eqref{lp.e9}, the group $N(\Lambda)^* = N(\Lambda)/K_0$ is isomorphic to 
$\mu(N(\Lambda)^*)$.  Hence the inverse image
 under $\mu$  of $Z(\mu(N(\Lambda)^*))= 
\mba \cap \mu(N(\Lambda)^*)$ in $N(\Lambda)^*$ is the center 
$Z(N(\Lambda)^*)$. 
Furthermore,  
\begin{equation}\mylabel{lp.e10}
N(\Lambda)^*/ Z(N(\Lambda)^*) \cong 
\mu(N(\Lambda)^*)/ Z(\mu(N(\Lambda)^*))\cong  
\mbn/\mba.
\end{equation}
Let $i$ be the above isomorphism that sends the factor group 
$N(\Lambda)^*/ Z(N(\Lambda)^*)$ onto $\mbn/ \mba$. 
(Of course $i$ is induced by  the restriction of $\mu$  to $N(\Lambda)^*$.)
As we have seen (at \eqref{lp.e8}), the character $\Lambda^*$ maps, 
 under
 $\mu$,  to the  restriction of $\VPS$ to $\mu(N(\Lambda)^*)$. 
Hence $\Lambda^*$  has a structure  similar  to that 
of $\VPS$, ie., $\Lambda^* \in \Irr(N(\Lambda)^*)$ 
lies above the unique linear character $\lambda^*$ of 
$Z(N(\Lambda)^*)$ that is carried, under $\mu$, to the restriction of 
$\VP$ to $\mba \cap \mu(N(\Lambda)^*)= 
Z(\mu(N(\Lambda)^*))$.
There is a natural alternating  bilinear form
 on  $N(\Lambda)^*/ Z(N(\Lambda)^*) \times N(\Lambda)^*/ Z(N(\Lambda)^*)$ 
defined by  
\begin{equation}\mylabel{lp.form3}
<x\, Z(N(\Lambda)^*),\;  y\, Z(N(\Lambda)^*)> = 
\lambda^*([x, y])= \VP([\mu(x), \mu(y)]) \in \ZZ_p,
\end{equation} 
for all $x, y  \in N(\Lambda)^*$.
The isomorphism $i$ carries this bilinear form onto the 
form $<, >$ of $\mbn/\mba \times \mbn/\mba$,
  defined in \eqref{lp.form1}.
Hence $N(\Lambda)^*/ Z(N(\Lambda)^*)$ is a symplectic group isomorphic 
to the symplectic group $\mbn/\mba$.

In view of   \eqref{lp.e7} and \eqref{lp.e8}, we get 
a natural isomorphism  between the groups $\mbg/ \mbn$ and 
$\mu(G(\psi, \Lambda)^*)/ \mu(N(\Lambda)^*)$.
This,  composed with $\mu$,  provides an isomorphism $\mu^*$ 
of $G(\psi, \Lambda)^*/ N(\Lambda)^*$ onto 
$\mbg/ \mbn$.
The group $G(\psi, \Lambda)^*/ N(\Lambda)^*$ acts on 
$N(\Lambda)^*/ Z(N(\Lambda)^*)$ via conjugation in $G(\psi, \Lambda)^*$, 
and leaves the form \eqref{lp.form3} invariant. As $\mu$ preserves 
conjugation, and  induces the isomorphism $i$, it follows 
that $\mu^*$ and $i$  send the action of 
$G(\psi, \Lambda)^*/ N(\Lambda)^*$ on $N(\Lambda)^*/ Z(N(\Lambda)^*)$
to the action of $\mbg/ \mbn $ on $\mbn / \mba$ in the sense that 
\begin{equation} \mylabel{lp.e11}
i(\bar{x}^{\bar{s}}) = i(\bar{x})^{\mu^*(\bar{s})} \in \mbn/ \mba,
\end{equation}
for all $\bar{x} \in N(\Lambda)^*/ Z(N(\Lambda)^*)$ and $\bar{s}
\in G(\psi, \Lambda)^*/ N(\Lambda)^*$.

The group $G(\psi)/ N$ is naturally isomorphic to 
the factor group $G(\psi, \Lambda)^*/ N(\Lambda)^*$,
 via the isomorphism  $\cdot / \Lambda$ in (2.13a) of  \cite{da3}.
Any coset $\gamma \in G(\psi)/ N$ gets mapped under $\cdot  /\Lambda$, to 
the image $\gamma / \Lambda$ of the coset $\gamma \cap
 G(\psi, \Lambda) in  
G(\psi, \Lambda)/ N(\Lambda)$ under the natural epimorphism of 
$G(\psi, \Lambda)$ onto $G(\psi, \Lambda)^*$.
We use this isomorphism to make the 
symplectic $\ZZ_p(G(\psi, \Lambda)^*/ 
N(\Lambda)^*)$-module $N(\Lambda)^*/ Z(N(\Lambda)^*)$
 into a symplectic 
$\ZZ_p(G(\psi)/ N)$-module. 
As we have already seen in Corollary \ref{lp.c3}, 
we may  turn the $\ZZ_p(\mbg / \mbn)$-module $\mbn/ \mba$ 
into a $\ZZ_p(G(\psi)/ N)$-module, using the isomorphism $j$ 
 of Corollary \ref{lp.c4}. But $j$ is the natural isomorphism  $\cdot / \Psi$,
as this is defined in (2.13a) of\cite{da3}.
 According to (6.1b) in \cite{da3}, the isomorphism $\cdot / \Lambda$ 
is the composition of   $ \cdot / \Psi = j $ with  $\mu^*$.
We conclude that $i$ is an isomorphism of $N(\Lambda)^*/ Z(N(\Lambda)^*)$ 
onto $\mbn/ \mba$ as symplectic $\ZZ_p(G(\psi)/ N)$-modules.
So Theorem \ref{lp.t2} follows.
\end{proof}

Theorem \ref{lp.t1} is now an easy corollary of Theorem \ref{lp.t2}, 
as $$
\mbn/ \mba \cong N(\Lambda)^*/ Z(N(\Lambda)^*) \cong 
\mbn' / \mba',
$$
as symplectic $G(\psi)/ N$-modules.

We conclude this section with a characterization of  
any   faithful linear limit $(\mbg, \mba, \VP, \mbn, \VPS)$ of 
$(G, A, \phi, N, \psi)$ when $N$ is nearly extra  special.
\begin{proposition}\mylabel{lp.p3} 
Assume that $(G, A, \phi, N, \psi)$ is a linear quintuple with $N$ a $p$-group, 
such that   $A= Z(N)$ is cyclic and 
 central in $G$.  Assume further that $A$ is maximal
 among the abelian characteristic subgroups of $N$, while $\phi$ is a 
faithful linear character of $A$. 
Then $V= N/ A$ is a symplectic $\ZZ_p(G/ N)$-space
with the symplectic form $<wA, yA> = \phi([w,y])$, 
 for any $w, y \in N$.
If $(\mbg, \mba, \VP, \mbn, \VPS)$ is a faithful linear limit of 
$(G, A, \phi, N, \psi)$, then $\mbn/ \mba$ is isomorphic as a symplectic 
 $\ZZ_p(G/ N)$-module  to 
$W^{\perp}/ W$,
 where $W$ is a maximal $G/ N$-invariant totally isotropic 
subspace of $V$, and $W^{\perp}$ is the perpendicular subspace to $W$ with 
respect to the above bilinear form.
\end{proposition}

\begin{proof}
As $A= Z(N)$ is maximal characteristic abelian subgroup of $N$, 
we conclude that $N$ is the central product of $A$ with an extra special 
$p$-group of exponent $p$. Hence the factor group $V= N/ A$ is 
a $\ZZ_p(G)$-module  and thus a $\ZZ_p(G/ N)$ module (see the discussion 
after Theorem \ref{lp.t1}).
Note also that $\psi$ is the unique character of 
$N$ that lies above $\phi$, and thus is  $G$-invariant as $\phi$ is.

Let $W$ be  maximal among the 
 $G/ N$-invariant totally isotropic subspaces  of $V$.  
If $X$ is the inverse image of $W$ in $N$, then $X$ is an  abelian normal 
subgroup of $N$ that contains $A= Z(N)$. (Note that $X$ could be $A$.)
Then  $\phi \in \Lin(A)$ extends  to a linear character $\lambda$ of $X$. 
 In addition, the stabilizer $X'$ 
of $\lambda$ in $N$ is the inverse image in $N$ of 
$W^{\perp}$, while  $X'/ X$ is naturally isomorphic to   the 
factor symplectic space $W^{\perp}/W$. 
Furthermore, if $\psi_{\lambda} \in \Irr(X')$ 
is the $\lambda$-Clifford correspondent of $\psi$, then   
the quintuple $(G(\lambda), X, \lambda, X', \psi_{\lambda})$ 
is a linear reduction of $(G, A, \phi, N, \psi)$. 
Now, $G= G(\lambda) \cdot N$ as $G$ fixes the unique character $\psi$ 
of $N$ that lies above $\lambda$. 
As $W$  is a maximal $G/ N$-invariant totally isotropic 
subspace of $V$, we conclude that
 $(G(\lambda), X, \lambda, X', \psi_{\lambda})$ is a linear limit of 
$(G, A, \phi, N, \psi)$. (Or else, $\lambda$ would be extended to an 
abelian  normal subgroup $B$ of $G(\lambda)$ contained in  $N(\lambda)$. 
Thus the image of $B$ in $V$ would be a $G(\lambda)$-invariant, 
and thus $G$-invariant, totally isotropic subspace of $V$, contradicting 
the maximality of $W$.) Hence, if $K= \Ker(\lambda)$, then 
 $(G(\lambda)/K, X/K, \lambda/K, X'/K, \psi_{\lambda}/K)$ is 
a faithful linear limit of $(G, A, \phi, N, \psi)$.
But  $(X'/K) / (X/K)$ is isomorphic to $X'/X$ (see \eqref{lp.e6}), and this 
isomorphism is $G(\lambda)$-, and thus $G$-,invariant. 
We  conclude that, for the faithful linear limit  
$(G(\lambda)/K, X/K, \lambda/K, X'/K, \psi_{\lambda}/K)$ of 
 $(G, A, \phi, N, \psi)$, the proposition holds, that is, 
$(X'/K) / (X/K)$ is  isomorphic to $W^{\perp}/ W$ 
for some maximal  $G/N$-invariant totally isotropic subspace of $N/A$.

According to Theorem \ref{lp.t1}, if  $(\mbg, \mba, \VP, \mbn, \VPS)$ is
 another faithful linear limit of $(G, A, \phi, N, \psi)$, then 
$\mbn/ \mba$ is isomorphic to $(X'/K) / (X/K)$, and this isomorphism is invariant 
under $G(\psi)/N = G/N$.
This completes the proof of Proposition \ref{lp.p3}.
\end{proof}

%%% Local Variables: 
%%% mode: latex
%%% TeX-master: t
%%% End: 

\section{Linear limits,  character towers and triangular sets} 
\mylabel{lmt}
Assume that we have the same situation as in Chapter \ref{elm}.
That is, we have a fixed  normal series 
\begin{subequations}\mylabel{lmt.e1}
\begin{equation}\mylabel{lmt.e1a}
1=G_0 \unlhd G_1 \unlhd \dots \unlhd G_n = G,  
\end{equation}
of $G$ that satisfies Hypothesis \ref{hyp1}.
We also  fix a character tower 
\begin{equation}\mylabel{lmt.e1b}
\{1=\chi_0, \chi_1, \dots, \chi_n \}
\end{equation}
for that  series, along with a representative 
of its corresponding   conjugacy class of triangular sets
\begin{equation}\mylabel{lmt.e1c}
\{Q_{2i-1}, P_{2r} | \beta_{2i-1}, \alpha_{2r} \}_{i=1, r=0} ^{l', \, \,  k'}.
\end{equation}
Along with the above system we fix a Hall system $\{ \ma, \mb\}$ 
of $G$  that satisfies \eqref{e.AB}, that is, 
\begin{gather}\mylabel{lmt.AB}  
\ma \in \Hall_{\pi}(G), \mb \in \Hall_{\pi'}(G),  \\
 \ma(\chi_1, \chi_2,\dots, \chi_h) \text{ and } 
\mb(\chi_1, \chi_2,\dots,
\chi_h) \text{ form a Hall system for } G(\chi_1, \chi_2,\dots, \chi_h),\\
\ma(\chi_1, \chi_2,\dots, \chi_n) = P^*_{2k'} \text{ and } \mb(\chi_1,
\chi_2,\dots, \chi_n) = Q^*_{2l'-1},
\end{gather}
\end{subequations} 
for all $h=1, \dots, n$.
The  way the above character tower, its triangular set, and the Hall system
change, if we take a linear limit with respect to a subgroup $G_i$ of $G$,
is in general arbitrary. 
In some special cases we can control these changes,
 as we will see  in the next two subsections. The basic results  were 
already proved  in Chapter \ref{elm}. Here we will apply them multiple times 
and  translate them into  the language of ``linear limits''.

For the rest of the chapter,
 we fix an integer $m=1,\dots,n$. Whenever necessary we consider the
smaller system
\begin{align}\mylabel{lmt.sm}
1=G_0 & \unlhd G_1 \unlhd \dots \unlhd G_m  \unlhd G,  \notag\\ 
\{1&=\chi_0, \chi_1, \dots, \chi_m \},\\
\{Q_{2i-1},& P_{2r} | \beta_{2i-1}, \alpha_{2r} \}_{i=1, r=0} ^{l, \, \,  k},  \notag
\end{align}
where the integers $k, l$ are related to 
$m$ via \eqref{kl:def}.
Of  course, as always,  along with the above system the groups
$\qw(\beta_{2k-1,2k})$ and $\hap(\alpha_{2l-2,2l-1})$ are uniquely
defined, up to conjugation, 
via Theorem \ref{sy.H} and Theorem \ref{sy.t1}, respectively.

We first work, as in Chapter \ref{elm},  inside a $\pi'$-group.
\subsection{ $``\ma(\beta_1)"$-invariant linear reductions }  
Assume that the normal series \eqref{lmt.e1a}, its character
 tower \eqref{lmt.e1b}, the triangular set \eqref{lmt.e1c} and the Hall 
system \eqref{lmt.AB}  are fixed.
In addition, we assume that $S$ is a subgroup of $G$ satisfying
\begin{subequations}\mylabel{lmt.e3} 
\begin{align}\mylabel{lmt.e3a}
S\unlhd G &\text{ with } S\leq Q_1, \text{ and } \\
\zeta \in \Lin(S) &\text{ is  $G$-invariant and lies under $\beta_1$. } 
\end{align}
\end{subequations}
Note that (\ref{lmt.e3}a,b) are  the conditions in \eqref{els.1}. 
Furthermore,   the quintuple $(G, S, \zeta, Q_1, \beta_1)$ is 
a linear one.  

Let $E$ be a normal subgroup of $G$ with $S \leq E \leq Q_1$, 
and $\lambda_1 \in \Lin(E)$ be a  linear character of $E$
lying above $\zeta$  and under $\beta_1$. So $ \lambda_1$ is  an 
extension of $\zeta $ to $E$.
Then we can use all the results of the first section of Chapter 
\ref{elm}. We also use the same notation as that introduced in 
Section \ref{elm1}.
 In particular, Remark \ref{elm.r0} implies that 
some  $G_1$-conjugate  
$\lambda \in \Irr(E)$ of $\lambda_1$  is $\ma(\chi_1)= \ma(\beta_1)$-invariant
 and  extends $\zeta$. 
So  the quintuple
$(G_{\lambda}, E, \lambda, Q_{1,\lambda}, \beta_{1,\lambda})$
is a linear reduction of $(G, S, \zeta, Q_{1}, \beta_{1})$.
We  call it an {\em $``\ma(\beta_1)''$-invariant  linear reduction, }
as  $\lambda$ was picked, among its $G_{1}$-conjugates, 
 to be $\ma(\beta_1)$-invariant.  
We saw  in \eqref{elm:e1}  that  
the series $1=\slg{0} \unlhd \slg{1} \unlhd \dots \unlhd \slg{n} 
= G_{\lambda}$, formed by the stabilizers  of $\lambda$ in the various 
subgroups $G_i$ of $G$, is a normal series of $G_{\lambda}$.
Along with that series of groups, we get the tower of characters 
$\{\slc{i} \in \slg{i}\}_{i=0}^n$, where $\slc{i}$ is the $\lambda$-Clifford 
correspondent of $\chi_i$  (see \eqref{el.e4a}).
As in Section \ref{elm1}, we add a subscript $\lambda$  to any object
 such as  $P_{2r}, Q_{2i-1},  \alpha_{2r}, \beta_{2i-1}$  etc,   to indicate the
 corresponding object for the $\lambda$-situation. 
We pick the  groups $\{ P_{2r, \lambda}, Q_{2i-1, \lambda} \}$, 
for all $r=1, \dots , k$ and all $i=1, \dots, l$, to satisfy the conditions in 
Proposition \ref{elm:p1}.  In particular,  we get $P_{2r}^* 
= \slps{2r}$ by \eqref{104a},  while
$Q_{2i-1}^*= Q_{2i-1,\lambda}^*$ by \eqref{104b},
 whenever $r=1,\dots, k'$ and 
$i=1,\dots, l'$, respectively.  Then  the triangular set 
$$
\{Q_{2i-1, \lambda}, P_{2r, \lambda} | \beta_{2i-1, \lambda}, \alpha_{2r, 
\lambda} \}_{i=1, r=0} ^{l', \, \,  k'}
$$
satisfies the conditions  in  Theorem \ref{D}. 
  In addition, the $\lambda$-Hall system  
$\{\ma_{\lambda}, \mb_{\lambda} \}$ for $G_{\lambda}$  can be chosen
 to satisfy the conditions in Theorem \ref{el.tAB}. In particular \eqref{el.eAB} implies 
\begin{equation}\mylabel{lm.eeA}
\ma_{\lambda}(\chi_{1,\lambda}) = \ma(\chi_1).
\end{equation}

Also for the fixed smaller system \ref{lmt.sm}, all the 
conclusions of  Theorems \ref{el.t2} and \ref{el.t22}
hold. Hence the groups $\qw(\beta_{2k-1,2k})$
 and $\hap(\alpha_{2l-1,2l-1})$ and their $\lambda$-correspondents
 can be chosen to satisfy the conditions  in 
Theorems \ref{sy.H} and \ref{sy.t1}, respectively,
along with \eqref{el.t2e} and \eqref{el.t22e}.
Thus
\begin{subequations}\mylabel{lm.eA}
\begin{gather}
 \qw(\beta_{2k-1,2k})(\lambda)=
\qw_{\lambda}(\beta_{2k-1,2k, \lambda}),\\
 \hap(\alpha_{2l-2,2l-1}) = \hap_{\lambda}(\alpha_{2l-2,2l-1, \lambda}),\\
\maq (\lambda) \geq \maq_{\lambda} \\
\map = \map_{\lambda}, 
\end{gather}
where the groups $\maq$ and $\map$ are defined as in \eqref{e0.5}.
Also  \eqref{e0.i}, along with \eqref{lm.eeA}, implies 
\begin{equation}\mylabel{lm.eA5}
\hap(\alpha_{2l-2,2l-1}) \leq \map \leq \ma(\chi_1)= \ma_{\lambda}(\slc{1}).
\end{equation}
\end{subequations}
Furthermore,
Corollary \ref{el.co2} implies that 
 the image  $I$ of $\qw(\beta_{2k-1,2k})$ in $\Aut(P_{2k}^*)$ 
equals the image $I_{\lambda}$ of $\qw_{\lambda}(\beta_{2k-1,2k, \lambda})$
in $\Aut(P_{2k,\lambda}^*) =\Aut(P_{2k}^*)$.  Similarly, 
the images
 of $\hap(\alpha_{2l-2,2l-1})$ and $\hap_{\lambda}(\sla{2l-2,2l-1})$ in 
$\Aut(Q_{2l-1,\lambda}^*)$  coincide, by Corollary \ref{el.co22}.

The following observation turns out to be very important for the ultimate proof of
 Main Theorem 1.
We define the group  $U$ as 
\begin{equation}\mylabel{U}
U:= Q_{2l-1}^* \rtimes J, 
\end{equation}
where $J$ is   the image of $\hap(\alpha_{2l-2, 2l-1})$ in $\Aut(Q_{2l-1}^*)$, as this 
was defined in \eqref{elp.e31}.
(Clearly the group $U$ depends on the smaller system \eqref{lmt.sm} and thus 
on $m$).
We observe that the quintuple $(U, S, \zeta, Q_{2l-1}^*, \beta_{2l-1}^*)$
is a linear one. The $``\ma(\beta_1)"$-invariant linear reduction 
$(G_{\lambda}, E, \lambda, Q_{1, \lambda}, \beta_{1,\lambda})$ of 
$(G, S, \zeta, Q_1, \beta_1)$  determines naturally the linear reduction 
$(U(\lambda), E, \lambda, Q_{2l-1,\lambda}^*, \beta_{2l-1,\lambda}^*)$
of the quintuple  $(U, S, \zeta, Q_{2l-1}^*, \beta_{2l-1}^*)$. Note that 
 $U(\lambda)= Q_{2l-1,\lambda}^* \rtimes J$, as $Q_{2l-1, \lambda}^*= 
Q_{2l-1}^*(\lambda)$ and $\lambda$ is 
$\ma(\beta_1) \geq \hap(\alpha_{2l-2,2l-1})$-invariant.
We call such a reduction a { \em $G$-associate linear reduction } of
$(U, S, \zeta, Q_{2l-1}^*, \beta_{2l-1}^*)$, as the $E$-group we are choosing for this 
reduction is not only normal in $U$, as it is the common case in linear reductions,
but it is normal in $G$. We remark here that  the group $U$ is isomorphic 
to a section of $G$.

Assume now that   there exists another pair 
$(E', \lambda'_1)$, where  $E '$ is a subgroup of $Q_{1, \lambda}$, 
normal in $ G_{\lambda}$,  such that 
$S\leq E \leq E'\leq \slq{1}$, and 
$\lambda'_1$ is a  
linear character of $E'$ that extends $\lambda$ and lies under $\slb{1}$.
Then, by Remark \ref{elm.r0}, we can replace $\lambda'_1$ with 
one of its  $\slq{1}$-conjugates $\lambda'$,
that is $\ma_{\lambda}(\slc{1})$-invariant and also lies under $\slb{1}$.
By \eqref{lm.eeA} the above character is  
$\ma(\beta_1)= \ma(\chi_1)=\ma_{\lambda}(\slc{1})$-invariant.
So we can repeat the same process and consider  an 
$``\ma(\beta_1)"$-invariant linear reduction 
$(G_{\lambda'}, E', \lambda', Q_{1, \lambda, \lambda'},
 \beta_{1, \lambda, \lambda'})$ 
 of the $``\ma(\beta_1)"$-linear reduction  
$ (G_{\lambda}, E, \lambda, Q_{1,\lambda}, \beta_{1,\lambda}).$
That is, we apply again the methods of Section \ref{elm1}, 
but this time for the normal series 
$1=\slg{0} \unlhd \slg{1} \unlhd \dots \unlhd \slg{n} 
= G_{\lambda}$, the normal subgroup $E$ of $G_{\lambda}$
 in the place of $S$, and the 
$G_{\lambda}$-invariant character $\lambda$ in the place of $\zeta$.
Clearly $E'$ satisfies \eqref{els.2}.
So Proposition \ref{elm:p1}, Theorems \ref{D}, \ref{el.t2} and \ref{el.t22}
along with their 
Corollaries \ref{el.co2} and \ref{el.co22} can be applied.
We  conclude that 
\begin{subequations}\mylabel{lm.0} 
\begin{gather}
P_{2k, \lambda, \lambda'}^*=P_{2k, \lambda}^*= P_{2k}^*, \mylabel{lm.0a}\\
Q_{2l-1, \lambda, \lambda'}^*= Q_{2l-1,\lambda}^*(\lambda')= 
Q_{2l-1}^*(\lambda, \lambda'), \mylabel{lm.0b} 
\end{gather}
 We also have 
a Hall system $\{ \ma_{\lambda, \lambda'}, 
\mb_{\lambda, \lambda'} \}$ for $G_{\lambda, \lambda'}$ 
that satisfies the conditions  Theorem \ref{el.tAB} and 
is derived from $\{ \ma_{\lambda}, \mb_{\lambda} \}$.  
For any  fixed 
$m=1, \dots, n$, the groups $\qw_{\lambda, \lambda'}$  and 
$\hap_{\lambda, \lambda'}$  can be chosen with respect to the above Hall 
system.  Hence  they  satisfy
   \begin{gather}
\qw_{\lambda,\lambda'} (\slb{2k-1,2k, \lambda, \lambda'})= 
\qw(\beta_{2k-1,2k})(\lambda, \lambda'),\mylabel{lm.01} \\
\hap_{\lambda, \lambda'}(\sla{2l-2,2l-1, \lambda, \lambda'}) 
= \hap_{\lambda}(\sla{2l-2,2l-1, \lambda})=
\hap(\alpha_{2l-2,2l-1}),\mylabel{lm.02}\\
\maq_{\lambda, \lambda'} \leq 
\maq(\lambda, \lambda'),\mylabel{lm.03} \\
\map_{\lambda, \lambda'}= \map. \mylabel{lm.04}
\end{gather}
Hence 
 \eqref{lm.02} and repeated applications of  \eqref{e0.i} and (\ref{el.eAB}a)
imply
\begin{equation}
\hap(\alpha_{2l-2,2l-1}) \text{ is a subgroup of } 
\ma(\beta_1)=\ma(\chi_1)=\ma_{\lambda}(\slc{1})=
\ma_{\lambda, \lambda'}(\chi_{1,\lambda'}).
\end{equation}
Furthermore, \eqref{lm.0a} implies
 that $\Aut(P_{2k, \lambda, \lambda'}^*)= \Aut(P_{2k}^*)$.
Thus Corollary \ref{el.co2} implies 
\begin{equation}\mylabel{lm.05}
I= I_{\lambda}=  \text{ the image of $\qw_{\lambda, \lambda'}(\beta_{2k-1,2k, \lambda, \lambda'})$
in $ \Aut(P_{2k, \lambda, \lambda'}^*)$},
\end{equation}
where $I$ is  the image  of 
$\qw(\beta_{2k-1,2k})$ in $\Aut(P_{2k}^*)$ and $I_{\lambda}$ that of 
$\qw_{\lambda}(\slb{2k-1,2k})$ in $\Aut(P_{2k, \lambda}^*)$.
 Also 
\begin{multline}\mylabel{lm.06}
\text{ The subgroups   
$ \hap_{\lambda, \lambda'}(\alpha_{2l-2,2l-1,\lambda, \lambda'})$ and  
$\hap(\alpha_{2l-2,2l-1})$  have the same images   in } \\
\Aut(Q_{2l-1,\lambda, \lambda'}^*),  
\Aut(Q_{2l-1,\lambda}^*) \text{ and } \Aut(Q_{2l-1}^*).
\end{multline}
\end{subequations}

As far as the quintuple $(U(\lambda), E, \lambda, Q_{2l-1, \lambda}^*, \beta_{2l-1, \lambda}^*)$ is
 concerned, we clearly have that one of its  $G_{\lambda}$-associate linear reductions is the quintuple
$(U(\lambda, \lambda'), E', \lambda', 
Q_{2l-1, \lambda, \lambda'}^*, \beta_{2l-1,\lambda, \lambda'}^*  )$.
Furthermore, the group $\hap(\alpha_{2l-2,2l-1})$ fixes $\lambda$, as it is a subgroup 
of $\ma(\chi_1)$, by (\ref{lm.0}h).
It also fixes $\lambda'$, as $\hap(\alpha_{2l-2,2l-1})$ 
is a subgroup of $\ma_{\lambda}(\chi_{1,\lambda})$, 
by (\ref{lm.0}h).
Hence the image $J$ of $\hap(\alpha_{2l-2,2l-1})$ in $\Aut(Q_{2l-1}^*)$,
fixes $\lambda, \lambda'$. Thus $J \leq U(\lambda, \lambda')$.
Even more, $Q_{2l-1, \lambda, \lambda'}^* = Q_{2l-1}^*(\lambda, \lambda')$, 
by \eqref{lm.0b}.
We conclude that 
$$
U(\lambda, \lambda') = Q_{2l-1,\lambda, \lambda'}^* \rtimes J.
$$

We can continue this process until we reach a linear limit 
\begin{equation}\mylabel{lmt.e4}
(l(G), l(S), l(\zeta), l(Q_{1}), l(\beta_{1})) 
\in LL(G, S, \zeta, Q_{1}, \beta_{1})
\end{equation}
of 
$(G, S, \zeta, Q_{1}, \beta_{1})$.
As this was done in a very specific way, at every linear reduction 
we were using a character that is invariant, 
under $\ma(\beta_1)$, we call  any such  limit
an {\em $``\ma(\beta_1)"$-invariant linear limit } of
 $(G, S, \zeta, Q_{1}, \beta_{1})$.
The fact that 
 we only consider $``\ma(\beta_1)"$-invariant linear characters in every 
linear reduction  does not restrict our options in the possible linear 
reduction we can perform, as,  according to Remark  \ref{elm.r0},
we can always replace any given linear character with one of its conjugates 
that is  $``\ma(\beta_1)"$-invariant. 

Of course, along with the limit group $l(G)$,
we reach a limit normal series of $l(G)$   
\begin{equation}\mylabel{lmt.e5}
1=l(G_0) \unlhd l(G_1) \unlhd l(G_2)
\unlhd \dots \unlhd l(G_n) = l(G), 
\end{equation}
where $l(G_i) = G_i \cap l(G)$. 

Along with the series  \eqref{lmt.e5} we get a limit character tower 
\begin{subequations}\mylabel{lmt.e6}
\begin{equation}\mylabel{lmt.e6a}
\{l(\chi_i)  \in \Irr(l(G_i))\}_{i=0}^n,
\end{equation}
where $l(\chi_0)=1$ and 
\begin{equation}\mylabel{lmt.e6c}  
l(\chi_i) \in LL(\chi_i) \text{ is a linear limit of } \chi_i, 
\end{equation}
for all $i=1,\dots,n$.
\end{subequations}
We also write 
\begin{equation}\mylabel{lmt.e7} 
\{l(Q_{2i-1}), l(P_{2r}) |l(\beta_{2i-1}), l(\alpha_{2r})
\}_{i=1, r=0}^{l', k'}
\end{equation}
for a representative of the unique $l(G)$-conjugacy class of triangular sets
of \eqref{lmt.e5} that corresponds to the character tower \eqref{lmt.e6a}, 
and is derived  from the original triangular set  \eqref{lmt.e1c} following the
rules in  Theorem \ref{D}.
In addition we write $\{l(\ma), l(\mb)\}$ for a Hall system of $l(G)$ 
that satisfies \eqref{sy.e2} for the limit case.

Of course, the above system restricts to the smaller
\begin{align}\mylabel{lmt.sml}
1=l(G_0) \unlhd &l(G_1) \unlhd l(G_2)
\unlhd \dots \unlhd l(G_m) \unlhd l(G), \notag \\ 
\{l(\chi_i)  &\in \Irr(l(G_i))\}_{i=0}^m, \\
\{l(Q_{2i-1}),& l(P_{2r}), l(P_0)|l(\beta_{2i-1}), l(\alpha_{2r})
\}_{i=1, r=0}^{l, k}. \notag 
\end{align}
We also  write $l(P_{2k}^*)$ for the product group 
$l(P_0)\cdot l(P_2) \cdots l(P_{2k})$,  and
$l(Q_{2l-1}^*)$ for the product $l(Q_1) \cdots l(Q_{2l-1})$.
Similarly, working for a fixed $m$ 
and looking at the smaller system \eqref{lmt.sml}, we denote by 
$l(\qw)$ the analogue of $\qw$ in this limit case, and by 
$l(\hap)$ the analogue of $\hap$,
i.e., 
$l(\qw)$ and $l(\hap)$, satisfy Theorems \ref{hat:p1} and \ref{sy.e}, respectively, 
  for the limit case.

Using this notation we can easily see that results 
 similar to \eqref{lm.0} hold. In particular,
\begin{theorem}\mylabel{lmt.t1}
Assume that the normal series,  the character tower and the triangular set 
in \eqref{lmt.e1} satisfy the conditions \eqref{lmt.e3}.  
Assume further that \eqref{lmt.e4} is an $``\ma(\beta_1)"$-invariant linear limit of 
$(G, S, \zeta, Q_{1}, \beta_{1})$, and 
\eqref{lmt.e6a} a  character tower that arises as a linear limit of 
 \eqref{lmt.e1b} (see \eqref{lmt.e6c}). 
Then the  triangular set \eqref{lmt.e7}, that corresponds  to 
the tower \eqref{lmt.e6a},  can be chosen to satisfy 
\begin{subequations}\mylabel{lmt.e8}
\begin{equation}
l(P_{2k}^*)= P_{2k}^*.\mylabel{lmt.e8a} 
\end{equation}
 Furthermore, a Hall system 
$\{ l(\ma), l(\mb) \}$ of $l(G)$  can be  derived from  $\{ \ma, \mb\}$, so that at
 every linear reduction the conditions in  Theorem \ref{el.tAB}   hold. 
Then, for any $m=1,\dots, n$, the groups $l(\qw)(l(\beta_{2k-1,2k}))$ 
and $l(\hap)(l(\alpha_{2l-2,2l-1}))$,
for the smaller system \eqref{lmt.sml}, can be chosen in association 
to $l(\mb)$ and $l(\ma)$ respectively. Therefore we have 
\begin{align}
l(\hap)(l(\alpha_{2l-2,2l-1}))&=\hap(\alpha_{2l-2,2l-1}), \\ 
l(\map) &= \map,\\
\mylabel{lmt.e8b}l(I) &= I, 
\end{align}
where $l(I)$ is  the image of $l(\qw)(l(\beta_{2k-1,2k}))$ in 
$\Aut(l(P_{2k}^*))$.
\end{subequations}
\end{theorem}

\begin{proof}
Follows immediately by repeated applications of Proposition \ref{elm:p1},
Theorems \ref{D}, \ref{el.tAB},  \ref{el.t2}, \ref{el.t22} and Corollaries \ref{el.co2}
and \ref{el.co22},
 at every $``\ma(\beta_1)"$-invariant linear reduction that we perform. 
\end{proof}

\begin{remark}\mylabel{lmt.r1}
The fact   $l(P_{2k}^*)=P_{2k}^*$ implies that the groups 
$\qw(\beta_{2k-1,2k})$ and $l(\qw)(l(\beta_{2k-1,2k}))$ have the same image, 
that is $I$,  in the automorphism group $\Aut(P_{2k}^*) = \Aut(l(P_{2k}^*))$.
\end{remark}

Similarly, 
\begin{remark}\mylabel{lmt.rr1}
The equation 
 $l(\hap)(l(\alpha_{2l-2,2l-1}))= \hap(\alpha_{2l-2,2l-1})$ implies that 
both  these groups have the same images  in the automorphism groups 
$\Aut(Q_{2l-1}^*)$ and $\Aut(l(Q_{2l-1}^*))$.
\end{remark}

Along with the limit \eqref{lmt.e4}, we reach the quintuple
\begin{equation}\mylabel{lm.eC1}
(l(U), l(S),l(\zeta), l(Q_{2l-1}^*), l(\beta_{2l-1}^*)), 
\end{equation}
This  is a multiple linear reduction, but not necessarily a linear limit, of 
 $(U, S, \zeta, Q_{2l-1}^*, \beta_{2l-1}^*)$,
 as we could possibly  reduce it further using a normal 
subgroup of  $l(T)$, that  is not normal in $l(G)$, and a linear 
extension of $l(\zeta)$ to that normal subgroup. We call 
$(l(U), l(S),l(\zeta), l(Q_{2l-1}^*), l(\beta_{2l-1}^*))$ a
{\em  $G$-associate  limit } of
$(U, S, \zeta, Q_{2l-1}^*, \beta_{2l-1}^*)$.  Note that $l(U)$ is isomorphic
 to a section of $l(G)$ .
We clearly have 
\begin{remark}\mylabel{lm.rr}
The  $G$-associate linear  limit  
  $(l(U), l(S),l(\zeta), l(Q_{2l-1}^*), l(\beta_{2l-1}^*))$  of  the quintuple
\linebreak
  $(U, S, \zeta, Q_{2l-1}^*, \beta_{2l-1}^*)$ satisfies
\begin{subequations}\mylabel{lm.eC}
\begin{align}
l(Q_{2l-1}^*)&= l(G) \cap Q_{2l-1}^*, \\
 l(U) &= l(Q_{2l-1}^*) \rtimes J.
\end{align}
\end{subequations}
\end{remark}

Repeated applications of Theorems \ref{el.t3} and \ref{el.t33}, at every 
$``\ma(\beta_1)"$-linear reduction that we perform, imply
\begin{theorem}\mylabel{lmt.t2}
If  $\beta_{2k-1, 2k}$ extends to $\qw(\beta_{2k-1,2k})$,
then the character  $l(\beta_{2k-1,2k})$ also extends to
 $l(\qw)(l(\beta_{2k-1,2k}))$.
Similarly, if $\alpha_{2l-2,2l-1}$ extends to $\hap(\alpha_{2l-2,2l-1})$, then 
$l(\alpha_{2l-2,2l-1})$ also  extends to $l(\hap)(l(\alpha_{2l-2,2l-1}))$. 
\end{theorem}

Now let $K$ be the kernel of the limit character $l(\zeta)$.
As we have seen in the previous section (see \eqref{lmf}), 
 we can form a faithful linear limit 
$$
(fl(G), fl(S), fl(\zeta), fl(Q_1), fl(\beta_1))=
(l(G)/K, l(S)/K, l(\zeta)/K, l(Q_1)/K,l(\beta_1)/K)
$$
of the linear quintuple $(G, S, \zeta, Q_1, \beta_1)$.
We call such a limit an  { \em $``\ma(\beta_1)"$-invariant faithful linear limit },
as it is obtained from an $``\ma(\beta_1)"$-invariant  linear limit.
Along with that we have a normal series of $fl(G)$
\begin{subequations}\mylabel{lmt.e10}
\begin{equation}\mylabel{lmt.e10a}
0=fl(G_0) \unlhd fl(G_1) \unlhd fl(G_2)
\unlhd \dots \unlhd fl(G_n)  = fl(G), 
\end{equation}
where $fl(G_i) = l(G_i)/K$, for all $i=1, \dots, n$. 
Along with the series  \eqref{lmt.e10a} we get a  character tower 
\begin{equation}\mylabel{lmt.e10b}
\{fl(\chi_i)  \in \Irr(fl(G_i))\}_{i=0}^n
\end{equation}
where  $fl(\chi_i)$ is the unique character of 
$fl(G_i)=l(G_i)/K$ that inflates to $l(\chi_i) \in \Irr(l(G_i))$,
for each  $i=1,\dots,n$. Let $fl(\chi_0)=1$, then  
\begin{equation}\mylabel{lmt.e10c}
fl(\chi_i) \in FLL(\chi_i) \text{ is a faithful linear limit of } \chi_i 
\end{equation}
for all $i=0, 1, \dots, n$.
Let 
\begin{equation}\mylabel{lmt.e10d}
\{ fl(Q_{2i-1}), fl(P_{2r}) |fl(\beta_{2i-1}), fl(\alpha_{2r}) \}_{i=1, r=0}^{l',k'}
\end{equation}
be the  representative of the unique $fl(G)$-conjugate class 
of triangular sets that corresponds  to \eqref{lmt.e10b}, that is derived from 
the set \eqref{lmt.e7}.
\end{subequations}
 
The fact that the  quintuple $(l(G), l(S), l(\zeta), l(Q_1), l(\beta_1))$ 
is a linear one  implies that  the group $l(S)$ 
and its irreducible character $l(\zeta)$ satisfy 
\eqref{els.1}. Thus we can apply the results of Section \ref{ker}.
In particular,  Theorem \ref{ker.t1} 
implies that the set \eqref{lmt.e10d}  satisfies 
\begin{subequations}\mylabel{lm.eD}
\begin{gather}
fl(Q_{2i-1})= (l(Q_{2i-1})K)/K, \\
fl(P_{2r})= (l(P_{2r})K)/K \cong l(P_{2r}).
\end{gather}
whenever $1\leq i \leq l'$ and $1\leq r \leq k'$.
Hence
\begin{gather}\mylabel{lmt.e11}
fl(Q_{2l-1}^*) = (l(Q_{2l-1}^*)K)/K, \\
fl(P_{2k}^*) = (l(P_{2k}^*) K) /K  \cong  l(P_{2k}^*).
\end{gather}
\end{subequations}
Even more, we can pick a Hall system $\{ fl(\ma), fl(\mb) \}$ 
of $fl(G)$ to satisfy  the conditions in Theorem \ref{ker.t4}, i.e., 
\begin{equation}\mylabel{lm.eha}
fl(\ma)= (l(\ma) K)/K  \cong l(\ma) \text{ and }  fl(\mb) = (l(\mb)K )/K.
\end{equation}

 For every fixed $m=1, \dots, n$, 
the smaller  limit system \eqref{lmt.sml}  provides the faithful limit system
\begin{subequations}\mylabel{lm.eE}
\begin{gather}
0=fl(G_0) \unlhd fl(G_1) \unlhd fl(G_2) \unlhd \dots \unlhd fl(G_m) \unlhd fl(G), \\
\{fl(\chi_i)  \in \Irr(fl(G_i))\}_{i=0}^m \\
\{ fl(Q_{2i-1}), fl(P_{2r}) |fl(\beta_{2i-1}), fl(\alpha_{2r}) \}_{i=1, r=0}^{l,k}
\end{gather}
\end{subequations}
where its triangular set (\ref{lm.eE}c)  satisfies \eqref{lm.eD}.
Even more, having fixed the Hall system 
\linebreak
$\{ fl(\ma), fl(\mb)\}$, 
Theorem \ref{ker.t2} implies that the group  $fl(\qw)$ can be chosen, (in 
relation  to $fl(\mb)$), to satisfy 
the conditions of Theorem \ref{hat:p1} for the faithful linear 
situation \eqref{lm.eE},
along with 
\begin{equation}\mylabel{lmt.e12}
fl(\qw)(fl(\beta_{2k-1,2k}))= (l(\qw)(l(\beta_{2k-1,2k}))K)/K.
\end{equation}
 Even more, if 
\begin{equation}\mylabel{lmt.e12.5}
fl(I) := \text { the image of $fl(\qw)(fl(\beta_{2k-1,2k}))$ 
in the automorphism group $\Aut(fl(P_{2k}^*))$, }
\end{equation}
 then identifying  $l(P_{2k}^*)$ with $fl(P_{2k}^*)$,   Corollary \ref{ker.co2} implies that 
 $$
fl(I) \cong  \text{ Image of $l(\qw)(l(\beta_{2k-1,2k}))$ in }
\Aut(l(P_{2k}^*)).
$$
But this last  group equals the image $I$ of 
 $\qw(\beta_{2k-1,2k})$ in $\Aut(P_{2k}^*)$, by Remark \ref{lmt.r1}. Hence 
\begin{equation}\mylabel{lmt.e13}
fl(I)\cong  l(I)= I.
\end{equation}

In addition, Theorem \ref{ker.t22} implies that we may choose  
the group $fl(\hap)$, (in relation  to $fl(\ma)$), for the system 
\eqref{lm.eE} to satisfy \eqref{sy.e}, along with 
\begin{equation}\mylabel{lmt.e122}
fl(\hap)(fl(\alpha_{2l-2,2l-1})) \cong l(\hap)(l(\alpha_{2l-2,2l-1})).
\end{equation}
Thus (see Corollary \ref{ker.co22}), the above two isomorphic 
groups  have the same image in  the automorphism group $\Aut(fl(Q_{2l-1}^*))$.
This, along with Remark \ref{lmt.rr1}, implies that the groups 
$fl(\hap)(fl(\alpha_{2l-2,2l-1}))$ and $\hap(\alpha_{2l-2,2l-1})$ have the same image in 
$\Aut(fl(Q_{2l-1}^*))$.

In conclusion we get 
\begin{theorem}\mylabel{lmt.t3}
Assume that the normal series, the character tower, the triangular set 
and the Hall system in \eqref{lmt.e1} are fixed. Along with them we fix 
 $S$ and $\zeta \in \Irr(S)$ to  satisfy \eqref{lmt.e3}. 
Let  $(fl(G), fl(S), fl(\zeta), fl(Q_1), fl(\beta_1))$
be  an $``\ma(\beta_1)"$-invariant faithful linear limit of
 $(G, S, \zeta, Q_1, \beta_1)$ and \eqref{lmt.e10b} be the  
 character tower for the normal series \eqref{lmt.e10a}, that  arises as 
the  faithful linear limit of  the tower \eqref{lmt.e1b}.
Then we can pick the  triangular set \eqref{lmt.e10d}
to satisfy \eqref{lm.eD}. In particular 
$P_{2r}^*$  is naturally isomorphic to $fl(P_{2r}^*)$, for all $r=1, \dots, k'$.
 We also  
derive  a  Hall system $\{ fl(\ma), fl(\mb)\}$ of $fl(G)$ from 
the original $\{\ma, \mb\}$, via \eqref{lm.eha} and Theorem \ref{lmt.t1}. 
Then for any $m=1,\dots, n$, 
the groups $fl(\qw)$ and $fl(\hap)$, for the smaller 
faithful  system \eqref{lm.eE}, can be chosen, (in association 
to $fl(\mb)$ and $fl(\ma)$),   to satisfy 
\begin{itemize}
\item[{1.}]
the associated isomorphism of $\Aut(P_{2k}^*)$ onto 
$\Aut(fl(P_{2k }^*))$ sends the image of
$\qw(\beta_{2k-1,2k})$ in $\Aut(P_{2k}^*)$ onto that of 
$fl(\qw)(fl(\beta_{2k-1,2k}))$ in $\Aut(fl(P_{2k}^*))$, i.e., 
$$fl(I)\cong I.$$
\item[{2.}]
$\hap(\alpha_{2l-2,2l-1}) \cong fl(\hap)(fl(\alpha_{2l-2,2l-1}))$,
and they both have the same image in $\Aut(fl(Q_{2l-1}^*))$.
\end{itemize}
\end{theorem}
\begin{proof}
Follows from  Theorem \ref{lmt.t1}, and equations (\ref{lm.eE}d),
 \eqref{lmt.e13} and \eqref{lmt.e122}.
\end{proof}

Furthermore, Theorems \ref{ker.t3} and \ref{ker.t33},  along with  Theorem \ref{lmt.t2}, easily imply
\begin{theorem}\mylabel{lmt.t4}
If the character $\beta_{2k-1,2k} \in \Irr(Q_{2k-1,2k})$ extends to 
$\qw(\beta_{2k-1,2k})$, then the character $fl(\beta_{2k-1,2k}) \in 
\Irr(fl(Q_{2k-1,2k}))$ extends to $fl(\qw)(fl(\beta_{2k-1,2k}))$.
Similarly,  if the character $\alpha_{2l-2,2l-1} \in \Irr(P_{2l-2,2l-1})$ extends to 
$\hap(\alpha_{2l-2,2l-1})$, then the irreducible character 
$fl(\alpha_{2l-2,2l-1})$ of  $fl(P_{2l-2,2l-1})$
extends to $fl(\hap)fl(\alpha_{2l-2,2l-1})$.
\end{theorem}

Finally,  the group $K= \Ker(l(\zeta))$
 is a normal subgroup of $l(G)$, as  $l(\zeta)$ is $l(G)$-invariant.
Furthermore, $K \leq l(Q_1) \leq l(Q_{2l-1}^*)$, 
thus $K \unlhd l(Q_{2l-1}^*)$. Since  the group $l(U)$ is isomorphic to 
a section of $l(G)$, and $K$ is a subgroup of $l(Q_{2l-1}^*) \leq l(U)$, 
 we conclude that    $K$  is also  
 a normal subgroup of $l(U)$.
 Hence we can form the faithful linear quintuple 
\begin{subequations}\mylabel{lm.eG}
\begin{equation}
(l(U)/K, l(S)/K, l(\zeta)/K, l(Q_{2l-1}^*)/K, l(\beta_{2l-1}^*)/K)= 
(fl(U), fl(S), fl(\zeta), fl(Q_{2l-1}^*), fl(\beta_{2l-1}^*). 
\end{equation}
We call the above quintuple a { \em $G$-associate faithful linear limit }
of $(U, S, \zeta, Q_{2l-1}^*, \beta_{2l-1}^*)$.
The fact that $K$ is a $\pi'$-normal subgroup of 
$l(Q_{2l-1}^*)$, while $J$ is the image of the $\pi$-group 
 $\hap(\alpha_{2l-2,2l-1})=l(\hap)(l(\alpha_{2l-2,2l-1}))$ in $\Aut(Q_{2l-1}^*)$, 
along with (\ref{lm.eC}b), implies 
\begin{equation}
fl(U)=l(U)/K= (l(Q_{2l-1}^*)/K )\rtimes J= fl(Q_{2l-1}^*) \rtimes J.
\end{equation}
\end{subequations}
Corollary \ref{lim.co1} clearly implies
\begin{proposition}\mylabel{lm.pr33}
Let $(fl(Q_{2l-1}^*) \rtimes J,   fl(S), fl(\zeta), fl(Q_{2l-1}^*), fl(\beta_{2l-1}^*)$
be a $G$-associate faithful linear limit of $(U, S, \zeta, Q_{2l-1}^*, \beta_{2l-1}^*)$. 
Then any faithful linear limit of the former quintuple 
is also a faithful linear limit of  
$(U, S, \zeta, Q_{2l-1}^*, \beta_{2l-1}^*)$. 
\end{proposition}

\subsection{ $``\mb(\alpha_2)"$-invariant linear reductions }
Assume that the normal series \eqref{lmt.e1a}, its character
 tower \eqref{lmt.e1b}  and the triangular set \eqref{lmt.e1c}
 are fixed.
In addition, we assume that $G_2$ is a direct product 
\begin{subequations}\mylabel{lmt.e13.5}
\begin{equation}
 G_2 = G_{2,\pi} \times G_{2, \pi'}
\end{equation}
  while 
\begin{equation}
\chi_1 \text{ is $G$-invariant, } 
\end{equation}
that is, \eqref{lp.1} holds.
Hence \eqref{117} holds for the triangular set 
\eqref{lmt.e1c}. In particular we have 
\begin{gather}
G_2= P_2 \times G_1 = P_2 \times Q_1, \\
\chi_2 = \alpha_2 \times \beta_1, \\
G(\chi_2) = G(\alpha_2).
\end{gather}
Furthermore, we assume that the normal subgroup $R$ of $G$
and its irreducible character $\eta \in \Irr(R)$, satisfy 
\eqref{eps.1}, that is,
\begin{gather}
 R \unlhd G \text{ with } R \leq P_2, \\
\eta \in \Lin(R) \text{ is  $G$-invariant and lies under $\alpha_2$. }
 \end{gather}
\end{subequations}
The quintuple $(G, R, \eta, P_2, \alpha_2)$ is clearly a linear one.
As with the $``\ma(\beta_1)"$-invariant linear reductions, we will 
get a linear limit of the above quintuple with respect to the 
group $\mb(\alpha_2) = \mb(\chi_2)$. 

To get a linear reduction of $(G, R, \eta, P_2, \alpha_2)$ 
we start with a normal subgroup $M$ of $G$ contained in 
$P_2$ and a linear character 
$\mu_1$ of $M$ that  extends  $\eta$ and lies under $\alpha_2$.
Note that all the hypothesis of Section \ref{elp} are satisfied,
 and therefore all the results of that section hold.
Thus,  according to   Remark \ref{elmp.r0}, there exists a $P_2$-conjugate
 $\mu \in \Lin(M)$  of $\mu_1$, such that  $\mu$
 is $\mb(\alpha_2)$-invariant,
 extends $\eta$, 
and lies under $\alpha_2$.  We proceed using the same notation as that of 
Section \ref{elp}. As in \eqref{elmp:e1}  we form the series 
\begin{subequations}\mylabel{lmt.13}
\begin{equation}
1=\smg{0} \unlhd \smg{1} \unlhd \dots \unlhd \smg{n} = G_{\mu},
\end{equation}
 consisting of the stabilizer of $\mu$ in the 
groups $G_i$ and $G$  for $i=2, 3,  \dots, n$.
In addition (see  \eqref{elmp.e3}),  we write 
\begin{align}
G_{1,\mu} &=1,\\ 
G_{2, \mu} & = P_2(\mu).
\end{align}
 Along with that we get, as in 
  \eqref{118} and  \eqref{elp.e6a}, the 
$\mu$-character tower $\{\smc{i} \}_{i=0}^n$, where 
\begin{align}
\smc{0}&=1,\\   
\smc{1}&=1,\\
\smc{2} &= \alpha_{2, \mu}.
\end{align}
\end{subequations}
Furthermore, $\alpha_{2, \mu}$ and   $\smc{i}$ are the
$\mu$-Clifford correspondents of $\alpha_2$ and $\chi_i$, respectively, 
for all $i=3,\dots,n$.  Proposition \ref{elmp:p1}
and Theorem \ref{P} 
 show  that we can choose a triangular set 
$\{\smq{1}=1, \smq{2i-1}, \smp{2r}|\smb{1}=1, \smb{2i-1}, \sma{2r}\}_{i=2, r=0}^{l',k'}$,
  that corresponds to 
 the above $\mu$-character tower,  so that 
$P_{2r, \mu}^*= P_{2r}^*(\mu)$, while 
$\smq{2i-1}^*=Q^*_{2i-1}$, whenever $1\leq r \leq k'$ and $1\leq i \leq l'$.
In addition, Theorem \ref{elp.tAB} implies that  the $\mu$-Hall  
system $\{ \ma_{\mu}, \mb_{\mu}\}$ for $G_{\mu}$ can be chosen to satisfy 
\eqref{elp.eAB}. Then 
$\mb_{\mu}(\smc{1},\smc{2})= \mb(\chi_1,\chi_2)$.
As  $\chi_2 = \alpha_2 \times \chi_1$, where $\chi_1$ is $G$-invariant and 
 $\smc{2}= \sma{2}$, 
 we conclude that 
\begin{equation}\mylabel{lmp.eeA}
\mb_{\mu}(\sma{2})= \mb_{\mu}(\smc{2})= \mb(\chi_1,\chi_2)= \mb(\alpha_2).
\end{equation}

We assume fixed the smaller system  \eqref{lmt.sm}. In addition, we 
assume that $m$ is any integer  so that 
$$
m \geq 2.
$$
Then  Theorems 
\ref{elp.t2} and \ref{elp.t22} hold for this smaller system.  
Hence the groups $\qw(\beta_{2k-1,2k})$ and $\hap(\alpha_{2l-2,2l-1})$,
 along with 
their $\mu$-correspondents  $\qw_{\mu}(\smb{2k-1,2k})$ and 
$\hap_{\mu}(\sma{2l-2,2l-1})$, can be chosen to satisfy 
\begin{subequations}\mylabel{lmp.eA}
\begin{gather}
 \qw(\beta_{2k-1,2k}) = \qw_{\mu}(\smb{2k-1,2k}), \\
\hap(\alpha_{2l-2,2l-1}, \mu) = \hap_{\mu}(\sma{2l-2,2l-1}),\\
\maq= \maq_{\mu}, \mylabel{lmt.e14} \\
\map(\mu) \geq  \map_{\mu}.
\end{gather}
Equation \eqref{lmp.eeA}, along with Remark \ref{elp.r000}, implies 
\begin{equation}\mylabel{lmp.eA5}
\qw(\beta_{2k-1,2k}) \leq \maq \leq \mb(\chi_2) = \mb_{\mu}(\sma{2}).
\end{equation}
\end{subequations}

Furthermore,
Corollary  \ref{elp.co3} implies that 
$\qw_{\mu}(\beta_{2k-1,2k, \mu})$ and $ \qw(\beta_{2k-1,2k})$
 have the same images in both  $\Aut(P_{2k}^*)$ and 
$\Aut(P_{2k, \mu}^*)$. Similarly Corollary \ref{elp.co33} implies that the image
 $J$ of $\hap(\alpha_{2l-2,2l-1})$ in $\Aut(Q_{2l-1}^*)$ equals
 the image $J_{\mu}$ of $\hap_{\mu}(\sma{2l-2,2l-1})$  in $\Aut(Q_{2l-1, \mu}^* )=
\Aut(Q_{2l-1}^*)$.  
The quintuple $(G_{\mu},M, \mu, P_{2, \mu}, \sma{2})$ is clearly a linear 
reduction of  $(G, R, \eta, P_2, \alpha_2)$. We call it a 
{ \em $``\mb(\alpha_2)"$-invariant linear reduction,} 
as $\mu$ is $\mb(\alpha_2)$-invariant.

Similarly to the group $U$, we write $T$ for group 
\begin{equation}\mylabel{T}
T= P_{2k}^* \rtimes I
\end{equation}
where (as always) 
$$
I = \text{ Image of $\qw(\beta_{2k-1,2k})$  in } \Aut(P_{2k}^*).
$$
It is clear that  the quintuple 
$(T, R, \eta, P_{2k}^*, \alpha_{2k}^*)$ 
is a linear one. Furthermore,  the $``\mb(\alpha_2)"$-invariant 
linear reduction 
 $(G_{\mu},M, \mu, P_{2, \mu}, \sma{2})$ of  
 $(G, R, \eta, P_2, \alpha_2)$  determines naturally a  linear reduction 
$(T(\mu), M, \mu, P_{2k,\mu}^*, \sma{2k}^*)$ of $(T, R, \eta, P_{2k}^*,
 \alpha_{2k}^*)$. Note that, as $\qw(\beta_{2k-1,2k})$ fixes $\mu$  
by \eqref{lmp.eA5}, its image $I$ in $\Aut(P_{2k}^*)$ also fixes $\mu$.
As $P_{2k, \mu}^*= P_{2k}^*(\mu)$, we conclude that 
the stabilizer $T(\mu)$ of $\mu $ in $T$
 satisfies  $T(\mu)= P_{2k, \mu}^* \rtimes I$.
 We call such a reduction a
{\em  $G$-associate linear reduction } of $(T, R, \eta, P_{2k}^*, \alpha_{2k}^*)$,  
as the 
 group $M$ we are choosing for this reduction is  normal in  $G$
 We also remark that $T$ and $T_{\mu}$  are 
 isomorphic to a section of $G$ and $G_{\mu}$, respectively.

Now we can repeat the procedure. 
So assume that  there exists another pair $(M', \mu'_1)$, such that 
 $M'$ is  a  normal subgroup of  $G_{\mu}$  satisfying 
$R \leq M \leq M' \leq P_{2, \mu}$, and 
$\mu'_1 \in \Lin(M')$ is an extension of $\mu$, and thus an extension 
of $\eta$, that lies under $\alpha_{2,\mu}$.
Again, using Remark \ref{elmp.r0}, we can replace 
$\mu'_1$ with a $\smp{2}$-conjugate $\mu'$ of  $\mu'_1$ that is 
$\mb_{\mu}(\sma{2})$-invariant, 
extends $\mu$,  and lies under $\alpha_{2,\mu}$. (So $\mu'$ is
 $\mb(\alpha_2) = \mb_{\mu}(\sma{2})$-invariant, by \eqref{lmp.eeA}).
We  apply   the results of Section \ref{elp} to  the series 
$\smg{0}= 1 \unlhd  \smg{1}=1 \unlhd \smg{2}=  \smp{2} \unlhd \smg{3} 
\unlhd   \dots \unlhd \smg{n} = G_{\mu}$,
  its character tower  $\{\smc{i} \}_{i=0}^n$,
and the triangular set 
$\{\smq{2i-1}, \smp{2r}|\smb{2i-1}, \sma{2r}\}_{i=1, r=0}^{l',k'}$, 
already picked at the previous reduction.
We also use the normal subgroup $M$ of $G_{\mu}$ in the 
place of $R$, the $G_{\mu}$-invariant character $\mu$ in the place 
of $\eta$, and the normal subgroup $M'$ of $G_{\mu}$ in the place of $M$. 
Notice that  \eqref{eps.1} holds, with $M'$ here 
 in the place of $M$ there, and  $M$ here in the  place of $R$ there.
 Furthermore, the group $G_{2, \mu}$ splits trivialy 
 as the product $G_{2, \mu} \times 1$  of a $\pi$-and 
a $\pi'$-group.  Thus the conditions 
\eqref{lp.1} and \eqref{eps.1} are satisfied. 
Hence all the results of Section \ref{elp} hold.
In particular, we have a normal series 
$$
1=G_{0,\mu,\mu'} \unlhd G_{1,\mu,\mu'} \unlhd 
G_{2, \mu, \mu'} \unlhd \dots \unlhd G_{n,\mu,\mu'} =
G_{\mu,\mu'},
$$
 of the stabilizer $G_{\mu, \mu'}= G(\mu, \mu')$
of $\mu'$ in $G_{\mu}=G(\mu)$. In addition,    
\begin{subequations}\mylabel{lmt.e144}
\begin{align}
G_{1, \mu, \mu'} &= G_{1, \mu} = 1, \\
G_{2, \mu, \mu'} &= P_{2}( \mu , \mu').
\end{align}
and  
$$
G_{i, \mu, \mu'}= G_{i, \mu} \cap
G_{ \mu, \mu'}= G_i(\mu, \mu'),
$$
 for all $i=2,3,\dots,n$.
We also get a character tower $\{\chi_{i,\mu, \mu'}\}_{i=0}^n$ for that series, 
where $\chi_{i,\mu,\mu'}$ is the $\mu'$-Clifford correspondent of 
$\chi_{i,\mu}$, for each $i=2, \dots, n$. 
Furthermore,  as in \eqref{lmt.13}, we have  
\begin{align}
\chi_{1, \mu , \mu'}&=1, \\
\chi_{2, \mu , \mu'} &=\alpha_{2, \mu , \mu'},
\end{align}
\end{subequations}
where $\alpha_{2, \mu, \mu'}$ is the $\mu'$-Clifford correspondent of 
$\alpha_{2, \mu}$.

Hence  
 Proposition \ref{elmp:p1}, Theorems \ref{P}, \ref{elp.t2} and \ref{elp.t22}, along 
with their Corollaries \ref{elp.co3} and \ref{elp.co33}, 
imply that we can pick a  triangular set 
$\{Q_{2i-1,\mu, \mu'}, P_{2r,\mu, \mu'}|\beta_{2i-1,\mu,\mu'},
 \alpha_{2r,\mu,\mu'}\}_{i=1, r=0}^{l',k'}$
that corresponds to the character tower 
$\{\chi_{\mu, \mu'}\}_{i=0}^n$, so that 
\begin{subequations}\mylabel{lmt.e14.5} 
\begin{gather}
P_{2k, \mu, \mu'}^*=P_{2k, \mu}^*(\mu')=P_{2k}^*( \mu , \mu')=P_{2k}^* \cap G_{\mu, \mu'},
\mylabel{lmp.0a}\\ 
Q_{1,\mu, \mu'}^*= Q_{1, \mu}^* = Q_{1, \mu}= 1, \\
Q_{2l-1,\mu, \mu'}^* =Q_{2l-1,\mu}^*=Q_{2l-1}^*.\mylabel{lmp.0b}
\end{gather}
We also pick a Hall system $\{ \ma_{\mu,\mu'}, \mb_{\mu, \mu'} \}$ of
 $G_{\mu, \mu'}$, that satisfies the conditions in 
 Theorem \ref{elp.tAB} and is derived from 
$\{ \ma_{\mu}, \mb_{\mu}\}$. So it is  
derived from the original $\{\ma, \mb \}$. 
Therefore, for any fixed $m=1,\dots,n$,  the groups    
$\qw(\beta_{2k-1,2k})$ and $\hap(\alpha_{2l-2,2l-1})$ 
for the smaller system \eqref{lmt.sm}, 
can be chosen via $\mb_{\mu, \mu'}$ and $\ma_{\mu, \mu'}$, respectively,
 (see Theorems \ref{sy.H} and \ref{sy.t1}).
Hence, as in Theorem \ref{elp.t2} and \ref{elp.t22}, they satisfy
\begin{gather}
  \qw_{\mu, \mu'}(\beta_{2k-1,2k, \mu,\mu'}) =\qw_{\mu}(\beta_{2k-1,2k,\mu})
=\qw(\beta_{2k-1,2k}),\mylabel{lmp.0c}\\
\hap_{\mu,\mu'}(\alpha_{2l-2,2l-1,\mu,\mu'})= 
\hap(\alpha_{2l-2,2l-1})(\mu,\mu'), \mylabel{lmp.0d}\\
\maq_{\mu, \mu'} =\maq_{\mu}= \maq, \mylabel{lmp.0e}\\
\map_{\mu,\mu'}  \leq  \map(\mu,\mu'). \mylabel{lmp.0f}
\end{gather}
Furthermore,
\begin{gather}
\qw(\beta_{2k-1,2k})  \text{ is a subgroup of }
 \mb(\alpha_2) = \mb_{\mu}(\sma{2})= 
\mb_{\mu,\mu'}(\alpha_{2,\mu, \mu'}),\mylabel{lmp.e.}\\
G_{2, \mu , \mu'} = P_{2, \mu, \mu'}.
\end{gather}
\end{subequations}

Equation \eqref{lmp.0c} implies that the groups $\qw(\beta_{2k-1,2k})$ and 
$\qw_{\mu,\mu'}(\beta_{2k-1,2k,\mu,\mu'})$ have the same image in 
$\Aut(P_{2k, \mu,\mu'}^*)$. Also, \eqref{lmp.0b}, along with 
 Corollary \ref{elp.co33}, implies that
$J_{\mu} = J_{\mu, \mu'}$, where
$J_{\mu, \mu'}$ denotes the image of 
$\hap_{\mu, \mu'}(\alpha_{2l-2,2l-1, \mu, \mu'})$
 in $\Aut(Q_{2l-1, \mu, \mu'}^*)$. So 
\begin{equation}\mylabel{lmp.05}
J = J_{\mu} = J_{\mu, \mu'}.
\end{equation}

As far as the linear reductions are concerned,  we have that  
$(G_{\mu,\mu'} M', \mu', P_{2, \mu, \mu'}, \alpha_{2, 
\mu, \mu'})$ is a $``\mb(\alpha_2)"$-invariant 
 linear reduction  of  
$(G_{\mu}, M, \mu, P_{2,\mu},
 \alpha_{2,\mu})$.
Furthermore,  the reduced quintuple 
 $(T( \mu, \mu'), M', \mu',
 P_{2k, \mu ,\mu'}^*, \alpha_{2k, \mu, \mu'}^*)$ is a $G_{\mu}$-associate 
linear reduction of $(T( \mu), M, \mu , P_{2k, \mu}^*, \alpha_{2k,\mu}^*)$.
Note that
\begin{equation}\mylabel{lmp.e00}
T(\mu, \mu') =( P_{2k}^* \rtimes I) (\mu, \mu')=
P_{2k,\mu, \mu'}^* \rtimes I,
\end{equation}
as both $\mu$ and $\mu'$, are 
$\mb(\alpha_2) \geq \qw(\beta_{2k-1,2k})$-invariant,
(by \eqref{lmp.e.}), 
and thus  $I$-invariant.

We continue this process until we reach a linear limit
\begin{equation}\mylabel{lmt.e15}
(l(G), l(R), l(\eta), l(P_{2}), l(\alpha_{2}))
\in LL(G, R, \eta, P_{2}, \alpha_{2}),
\end{equation}
 that we call  a  {\em $``\mb(\alpha_2)"$-invariant 
linear limit } of the linear quintuple 
$(G, R, \eta, P_{2}, \alpha_{2})$.
We also reach  a limit normal series for the
group $l(G)$
\begin{subequations}\mylabel{lmt.e16}
\begin{equation}\mylabel{lmt.e16a}
1 =l(G_0) \unlhd l(G_1) \unlhd l(G_2)
\unlhd \dots \unlhd l(G_n) =l(G), 
\end{equation}
where $l(G_i) = G_i \cap l(G)$, for all $i=2, \dots, n$, and 
\begin{align}\mylabel{lmt.e16aa}
l(G_1)&= 1, \\
l(G_2) &= l(P_2), 
\end{align}
as the same holds  at every linear reduction.
Observe also  that the above normal series has the same notation as the one in 
 \eqref{lmt.e5},
 but of course  is produced in a different way.

Along with the series  \eqref{lmt.e16a} we get a  character tower 
\begin{equation}\mylabel{lmt.e16b}
\{l(\chi_i)  \in \Irr(l(G_i))\}_{i=0}^n
\end{equation}
where  
\begin{align}\mylabel{lmt.e16c}
l(\chi_1)&=1 = l(\chi_0),\notag \\
l(\chi_2) &= l(\alpha_2), \\
l(\chi_i) &\in LL(\chi_i) \text{ is a  linear limit of } \chi_i \notag
\end{align}
for all $i=3,\dots,n$. 
Let 
\begin{equation}\mylabel{lmt.e16d}
\{ l(Q_{2i-1}), l(P_{2r}) |l(\beta_{2i-1}), l(\alpha_{2r}) \}_{i=1, r=0}^{l',k'}
\end{equation}
be  the representative of the unique $l(G)$-conjugate class 
that corresponds  to \eqref{lmt.e16b}, and is derived from  the 
original triangular set \eqref{lmt.e1c} following the rules in Theorem \ref{P}. .
We also denote by 
\begin{equation}
\{ l(\ma), l(\mb) \}
\end{equation}
a Hall system for $l(G)$ that satisfies \eqref{sy.e},
 for the above limit case. 
\end{subequations}

Of course the above system restricts to the smaller 
\begin{subequations}\mylabel{lmp.sml}
\begin{gather}
1=l(G_0) \unlhd l(G_1) \unlhd l(G_2)
\unlhd \dots \unlhd l(G_m) \unlhd l(G), \notag \\ 
\{l(\chi_i)  \in \Irr(l(G_i))\}_{i=0}^m \\
\{l(Q_{2i-1}), l(P_{2r}) |l(\beta_{2i-1}), l(\alpha_{2r})
\}_{i=1, r=0}^{l, k} \notag 
\end{gather}
\end{subequations}
Note that we have the same  notation as that in \eqref{lmt.sml}.
Similar to the notation there, we write $l(P_{2k}^*)$ and $l(Q_{2l-1}^*)$
 for the product groups 
$l(P_0)\cdot l(P_2) \cdots l(P_{2k})$ and 
$l(Q_1) \cdots l(Q_{2l-1})$, respectively.
 Also for any  fixed $m$, we denote by 
$l(\qw)$ the analogue of $\qw$ in this limit case, and by 
$l(\hap)$ the analogue of $\hap$, for the smaller system \eqref{lmp.sml}
i.e., 
$l(\qw)$ and $l(\hap)$, satisfy the conditions in  Theorems \ref{hat:p1} and \ref{sy.e}, respectively, 
  for the limit case.

 Then
\begin{theorem}\mylabel{lmt.t5}
Assume that the normal series,  the character tower, the triangular set  and the Hall system 
in \eqref{lmt.e1} satisfy the conditions \eqref{lmt.e13.5}.  
Assume further that \eqref{lmt.e15} is a $``\mb(\alpha_2)"$-invariant linear limit of 
$(G, R, \eta, P_2, \alpha_2)$ and 
\eqref{lmt.e16a} a  character tower that arises as a linear limit of 
 \eqref{lmt.e1b} (see \eqref{lmt.e16c}). 
Then the  triangular set \eqref{lmt.e16d}, that corresponds  to 
the tower \eqref{lmt.e16a},  can be chosen to satisfy 
\begin{subequations}\mylabel{lmt.e17}
\begin{equation}\mylabel{lmt.e17e}
l(Q_{2i-1}^*) = Q_{2i-1}^*,
\end{equation}
for all $i=1,\dots,l'$. Furthermore, a Hall
 system  $\{ l(\ma), l(\mb) \}$  for $l(G)$, 
can be derived  from $\{ \ma, \mb \}$  so that at every linear reduction 
theconditions in Theorem \ref{elp.tAB} hold.
Then, for every $m=1, \dots, n$, the groups
 $l(\qw)(l(\beta_{2k-1,2k}))$ and  
$l(\hap)(l(\alpha_{2l-2,2l-1}))$ for the smaller system  \eqref{lmp.sml}, 
 can be chosen, using  the groups $l(\mb)$ and $l(\ma)$, respectively, to satisfy
 \begin{align}
\qw(\beta_{2k-1,2k}) &=l(\qw)(l(\beta_{2k-1,2k})), \mylabel{lmt.e17c} \\
\maq &= l(\maq), \mylabel{lmt.e17d}\\
J &=l(J), \mylabel{lmt.e17f}
\end{align}
\end{subequations}
where $l(J)$ is the image of $l(\hap)(l(\alpha_{2l-2,2l-1}))$ in $\Aut(l(Q_{2l-1}^*))$.
\end{theorem}

\begin{proof}
We reach the linear limit \eqref{lmt.e15} doing, at every step,  
$``\mb(\alpha_2)"$-invariant 
linear reductions.  Therefore at every step we are picking a 
triangular set  that satisfies the conditions in   Proposition \ref{elmp:p1} and
 Theorem \ref{P}. We also pick, at every linear redution,
  a Hall system that satisfies theconditions in 
Theorem \ref{elp.tAB}. 
 Furthermore, for any fixed $m=1, \dots, n$, the groups $\qw$ 
and $\hap$  satisfy the conditions in Theorems \ref{elp.t2} and \ref{elp.t22}. 
Hence  at every step  equations \eqref{lmt.e14.5} hold.
In particular, repeated applications of \eqref{lmp.0b}, \eqref{lmp.0c} and 
\eqref{lmp.0e} imply (\ref{lmt.e17}a,b) and  \eqref{lmt.e17d}, respectively.
Similarly, repeated applications of \eqref{lmp.05} imply \eqref{lmt.e17f}.
Hence Theorem \ref{lmt.t5} follows. 
\end{proof}

As an easy consequence of \eqref{lmt.e17e}  we have 
\begin{remark}\mylabel{lmp.rm}
The groups $\qw(\beta_{2k-1,2k})$ and $l(\qw)(l(\beta_{2k-1,2k}))$ have
 the same images in 
both automorphism groups  $\Aut(P_{2k}^*)$  and $\Aut(l(P_{2k}^*))$. 
\end{remark}

Also repeated applications of 
Theorems \ref{elp.t3} and \ref{elp.t33}
 at every $``\mb(\alpha_2)"$-invariant linear reduction 
imply
\begin{theorem}\mylabel{lmt.t8}
If $\beta_{2k-1,2k}$ extends to $\qw(\beta_{2k-1,2k})$, then 
the character $l(\beta_{2k-1,2k})$ also extends to the  limit group 
$l(\qw)(l(\beta_{2k-1,2k}))$. Similarly, if $\alpha_{2l-2,2l-1}$
extends to $\hap(\alpha_{2l-2,2l-1})$ then 
$l(\alpha_{2l-2,2l-1})$ also extends to $l(\hap)(l(\alpha_{2l-2,2l-1}))$.
\end{theorem}

Notice that, along with the limit in \eqref{lmt.e15}, 
we reach the quintuple
\begin{equation}\mylabel{lmt.e15.5}
(l(T), l(R), l(\eta), l(P_{2k}^*), l(\alpha_{2k}^*) ),
\end{equation}
that we call a
{\em $G$-associate   limit } of $(T, R, \eta, P_{2k}^*, \alpha_{2k}^*)$.
Note that, as with \eqref{lm.eC1}, the $G$-associate limit
 is a multiple linear reduction,   but  
not a linear limit, of 
$(T, R, \eta, P_{2k}^*, \alpha_{2k}^*)$.
Because \eqref{lmp.e00} holds  for 
every $``\mb(\alpha_2)"$-invariant linear reduction,  we have  
\begin{proposition}\mylabel{lmp.opr}
The $G$-associate linear limit 
$(l(T), l(R), l(\eta), l(P_{2k}^*), l(\alpha_{2k}^*) )$
of the quintuple
\linebreak
 $(T, R, \eta, P_{2k}^*, \alpha_{2k}^*)$,
satisfies 
\begin{subequations}
\begin{align}
l(P_{2k}^*)&=l(G) \cap P_{2k}^*, \mylabel{lmt.e15.6b} \\
l(T) &= l(P_{2k}^*) \rtimes I. \mylabel{lmt.e15.6c} 
\end{align}
\end{subequations}
\end{proposition}

We want to pass to a faithful linear limit of $(G, R, \eta, P_2, \alpha_2)$,
as we did with the $``\ma(\beta_1)"$-invariant case in  \eqref{lmt.e10}. 
So we  first note that  $l(G_2)= l(P_2)$. So $l(G_2)$  is the  product of 
a $\pi$-and a trivial $\pi'$-group.  Furthermore, 
 $l(R)$ is a normal subgroup of $l(G)$,  while 
$l(\eta) \in \Lin(l(R))$ is an $l(G)$-invariant linear character that lies under 
$l(\alpha_2) \in \Irr(l(P_2))$, as $(l(G),l(R), l(\eta), l(P_2), l(\alpha_2))$
is a linear quintuple. Hence all the conditions of Section \ref{ker2} are satisfied. 
Thus if  $K=\Ker(l(\eta))$ is the kernel of $l(\eta)$,
 then we can form the  faithful linear limit
$$
(fl(G), fl(R), fl(\eta), fl(P_2), fl(\alpha_2))=
(l(G)/K, l(R)/K, l(\eta)/K, l(P_2)/K,l(\alpha_2)/K)
$$
of the linear quintuple $(G, R, \eta, P_2, \alpha_2)$.
We call this a { \em $``\mb(\alpha_2)"$-invariant faithful linear limit },
as it is obtained from a $``\mb(\alpha_2)"$-invariant linear limit.
Along with it we have a normal series of $fl(G)$, as in  \eqref{ke.e3},
\begin{subequations}\mylabel{lmt.e18}
\begin{equation}\mylabel{lmt.e18a}
1=fl(G_0) \unlhd fl(G_1) \unlhd fl(G_2)
\unlhd \dots \unlhd fl(G_n) = fl(G), 
\end{equation}
where  
\begin{align}\mylabel{lmt.e18aa} 
fl(G_1) = l(G_1) K / K &= 1, \\
fl(G_2)  = l(G_2)K / K &=l(P_2) / K \\
fl(G_i) &= l(G_i)/K,
\end{align}
 for all $i=3, \dots, n$. 
Along with the series  \eqref{lmt.e18a} we get a  character tower, see \eqref{ke.e3b}, 
\begin{equation}\mylabel{lmt.e18b}
\{fl(\chi_i)  \in \Irr(fl(G_i))\}_{i=0}^n
\end{equation}
where 
\begin{align}\mylabel{lmt.e18c}
fl(\chi_1) &= 1, \\
fl(\chi_2) &= fl(\alpha_2), \\
fl(\chi_i) &\in FLL(\chi_i) \text{ is a faithful linear limit of } \chi_i 
\end{align}
for all $i=3,\dots,n$. That is, $fl(\chi_i)$ is the unique character of 
$fl(G_i)=G_i/K$ that inflates to $\chi_i \in \Irr(G_i)$.
Let 
\begin{equation}\mylabel{lmt.e18d}
\{ fl(Q_{2i-1}), fl(P_{2r}) |fl(\beta_{2i-1}), fl(\alpha_{2r}) \}_{i=1, r=0}^{l',k'}
\end{equation}
be a representative of the unique $fl(G)$-conjugate class 
that corresponds  to \eqref{lmt.e18b}.
\end{subequations}
 Then Theorem \ref{ke.t1} implies that we can 
pick the set \eqref{lmt.e18d} so that 
\begin{subequations}\mylabel{lmt.e119}
\begin{gather}
fl(Q_{2i-1})= (l(Q_{2i-1})K)/K \cong l(Q_{2i-1}), \notag  \\
fl(P_{2r})= (l(P_{2r})K)/K.
\end{gather}
whenever $1\leq i \leq l'$ and $1\leq r \leq k'$.
Hence
\begin{align}\mylabel{lmt.e19}
fl(Q_{2l-1}^*) &= (l(Q_{2l-1}^*)K)/K  \cong l(Q_{2l-1}^*),\\
fl(P_{2k}^*)&=   (l(P_{2k}^*)K) /K= l(P_{2k}^*)/K.
\end{align}
\end{subequations}
Furthermore, we can pick a Hall system $\{ fl(\ma), fl(\mb) \}$ for $fl(G)$ 
to satisfy  Theorem \ref{ke.t4}, that is 
\begin{equation}\mylabel{lmp.eha}
fl(\ma) = (l(\ma)K)/K  \text{ and }  fl(\mb) = (l(\mb)K)/K \cong l(\mb).
\end{equation}

For every fixed $m=1,\dots,n$ we have the smaller faithful 
limit system 
\begin{subequations}\mylabel{lmp.eE}
\begin{gather}
0=fl(G_0) \unlhd fl(G_1) \unlhd fl(G_2) \unlhd \dots \unlhd fl(G_m) \unlhd fl(G), \\
\{fl(\chi_i)  \in \Irr(fl(G_i))\}_{i=0}^m \\
\{ fl(Q_{2i-1}), fl(P_{2r}) |fl(\beta_{2i-1}), fl(\alpha_{2r}) \}_{i=1, r=0}^{l,k}
\end{gather}
\end{subequations}
withn the triangular set picked so that \eqref{lmt.e19} holds.
Furthermore, if $fl(\qw)$ denotes the corresponding to $\qw$ for the system 
 \eqref{lmp.eE}, then  
Theorem \ref{ke.t2} implies that, having fixed the Hall system \ref{lmp.eha}, 
 we can choose $fl(\qw)$ so that 
\begin{equation*}
fl(\qw)= (l(\qw)K)/K \cong l(\qw).
\end{equation*}
Even more, Corollary \ref{ke.co2} implies that 
\begin{multline}\mylabel{lmt.e20}
fl(\qw)(fl(\beta_{2k-1,2k})) \text{ and } l(\qw)(l(\beta_{2k-1,2k})) \\
 \text { have the same image 
in the automorphism group $\Aut(fl(P_{2k}^*))$.}
\end{multline}
In addition, Theorem \ref{ke.t22} implies that we may choose the group $fl(\hap)$ for the smaller system 
\eqref{lmp.eE} so that 
\begin{subequations}\mylabel{lmp.e20}
\begin{equation}
fl(\hap)(fl(\alpha_{2l-2,2l-1}))=(l(\hap)(l(\alpha_{2l-2,2l-1}))K)/K. 
\end{equation}
Thus, as in  Corollary \ref{ke.co22}, identifying  $l(Q_{2l-1}^*)$  with 
the isomorphic group $fl(Q_{2l-1}^*)$, 
we conclude that  the image  $fl(J)$ of $fl(\hap)(fl(\alpha_{2l-2,2l-1}))$ in $\Aut(fl(Q_{2l-1}^*))$ 
is isomorphic to the image $l(J)$ of $l(\hap)(l(\alpha_{2l-2,2l-1}))$ in $\Aut(l(Q_{2l-1}^*))$.
But the latter equals $J$, by \eqref{lmt.e17f}. Hence 
\begin{equation}fl(J) \cong l(J) = J.
\end{equation}
\end{subequations}

In conclusion we get 
\begin{theorem}\mylabel{lmt.t7}
Assume that the normal series, the character tower, the triangular set  and the Hall system 
in \eqref{lmt.e1} are fixed and satisfy (\ref{lmt.e13.5}a). Along with them we fix 
 $R$ and $\eta \in \Irr(R)$ to  satisfy (\ref{lmt.e13.5}b,c). 
Let  $(fl(G), fl(R), fl(\eta), fl(P_2), fl(\alpha_2))$
be  a $``\mb(\alpha_2)"$-invariant faithful linear limit of
 $(G, R, \eta, P_2, \alpha_2)$,  and \eqref{lmt.e18b} be 
a character tower for the normal series \eqref{lmt.e18a}, arising  as 
a faithful linear limit of  the tower \eqref{lmt.e1b}.
Then we can pick the  triangular set \eqref{lmt.e18d} to satisfy \eqref{lmt.e119}. In particular, 
$Q_{2i-1}^*$ is naturally isomorphic to $fl(Q_{2i-1}^*)$, for all
 $i=1,\dots,l'$. 
We also reach  a Hall system $\{ fl(\ma), fl(\mb)\}$ for $fl(G)$,   from $(\ma, \mb)$   via 
\eqref{lmp.eha} and Theorem \ref{lmt.t5}. Then for any $m=1,\dots,n$, the groups 
$fl(\qw)$ and $fl(\hap)$ for the smaller system \eqref{lmp.eE},
 can be chosen, 
with respect to $fl(\mb)$ and $fl(\ma)$ respectively, to satisfy
\begin{itemize}
\item[{1.}]
$\qw(\beta_{2k-1,2k}) $ is isomorphic to 
$fl(\qw)(fl(\beta_{2k-1,2k}))$,  and they both have the same image in 
$\Aut(fl(P_{2k}^*))$.
\item[{2.}]
the associated isomorphism of   $\Aut(Q_{2l-1}^*)$ onto the
 group of automorphisms 
$\Aut(fl(Q_{2l-1}^*))$  sends the image of  $\hap(\alpha_{2l-2,2l-1})$  
inside  $\Aut(Q_{2l-1}^*)$
onto the image of $fl(\hap)(fl(\alpha_{2l-2,2l-1})) $  inside 
 $\Aut(fl(Q_{2l-1}^*))$, i.e., 
$$fl(J)\cong J.$$
\end{itemize}
 \end{theorem}
\begin{proof}
We have already seen that we can pick the set \eqref{lmt.e18d} so that \eqref{lmt.e119}  holds.
The rest follows from \eqref{lmt.e20} and \eqref{lmp.e20}.
\end{proof}

Furthermore, Theorems \ref{ker.t3} and \ref{ke.t33}, 
 along with  Theorem \ref{lmt.t8}, easily imply
\begin{theorem}\mylabel{lmt.t44}
If the character $\beta_{2k-1,2k} \in \Irr(Q_{2k-1,2k})$ extends to 
$\qw(\beta_{2k-1,2k})$, then the character $fl(\beta_{2k-1,2k}) \in 
\Irr(fl(Q_{2k-1,2k}))$ extends to $fl(\qw)(fl(\beta_{2k-1,2k}))$.
Similarly,  
if the character $\alpha_{2l-2,2l-1} \in \Irr(P_{2l-2,2l-1})$ extends to 
$\hap(\alpha_{2l-2,2l-1})$, then the irreducible character 
$fl(\alpha_{2l-2,2l-1})$ of  $fl(P_{2l-2,2l-1})$
extends to $fl(\hap)fl(\alpha_{2l-2,2l-1})$.
\end{theorem}

The character $l(\eta)$ is $l(T)$-invariant and thus $I$-invariant. 
 Furthermore, $K = \Ker(l(\eta))$ is a subgroup of $l(P_2) \leq P_2$.
Hence  equation \eqref{lmt.e15.6c} implies 
$$
fl(T)=l(T)/ K =( l(P_{2k}^* )/K ) \rtimes I= 
fl(P_{2k}^*) \rtimes I.
$$
So we can form the quintuple
\begin{equation}\mylabel{lmt.ef}
(fl(T)=fl(P_{2k}^*)\rtimes I, fl(R), fl(\eta),
fl(P_{2k}^*), fl(\alpha_{2k}^*)), 
\end{equation}
that we call a $G$-associate faithful linear limit of 
$(T, R, \eta, P_{2k}^*, \alpha_{2k}^*)$.
Corollary \ref{lim.co1} clearly implies 
\begin{proposition}\mylabel{lmt.pr3}
Any faithful linear limit of $(fl(T), fl(R), fl(\eta),
fl(P_{2k}^*), fl(\alpha_{2k}^*))$ is also a faithful linear limit of 
 $(T, R, \eta, P_{2k}^*, \alpha_{2k}^*)$.
\end{proposition}

%%% Local Variables: 
%%% mode: latex
%%% TeX-master: t
%%% End: 

\chapter{ Main Theorem}
\mylabel{m}
\section{An outline of the proof }
\mylabel{ma.se1}
We start with a monomial group $G$ of order $p^a q^b$, for distinct 
odd primes $p, q$ and 
integers $a, b \ge 0$. Of course $G$ is solvable. 
Hence there exists some chain 
\begin{subequations}\mylabel{ma.e1}
\begin{equation}\mylabel{ma.A}
1=G_0 \unlhd G_1 \unlhd G_2 \unlhd  \dots \unlhd G_n = G, 
\end{equation}
of normal subgroups $G_i$ of $G$ that satisfy Hypothesis \ref{hyp1}. 
So $G_1$ is a $q$-group, while $G_i / G_{i-1}$ is a $p$-group if $i$ is even, 
and a $q$-group if $i$ is odd, for each $i=2, \dots, n$.
Let $\chi_0, \chi_1,  \dots, \chi_n$ satisfy 
\begin{equation}\mylabel{ma.B}
 \chi_i \in \Irr(G_i)  \text{ lies under  } \chi_j \in \Irr(G_j)
\end{equation}
for any $i,j = 0, 1,2,\dots,n$ with $i \le j$, i.e., the 
$\chi_i$ form a character tower for the 
series \eqref{ma.A}. Assume further that the integers $k'$ and $l'$ are 
related to $n$ via \eqref{kl:def}, while the set 
\begin{equation}\mylabel{ma.BB}
\{ Q_{2i-1}, P_{2r} | \beta_{2i-1}, \alpha_{2r} \}_{i=1, r=0}^{l', k'}
\end{equation}
is a representative of the unique conjugacy class of triangular sets 
 that correspond  to \eqref{ma.B}.

 We fix a Sylow 
 system $\{ \ma, \mb\}$ of $G$ satisfying \eqref{sy.e2} with
 $\pi = \{ p \}$,   that is,
\begin{align}\mylabel{ma.C}
\ma \in \Syl_p(G) &\text{ and } \mb \in \Syl_q(G), \\
 \ma(\chi_1, \chi_2, \dots, \chi_i) \in
\Syl_p(G(\chi_1,\chi_2, \dots, \chi_i)) &\text{ and } \mb(\chi_1,
\chi_2, \dots, \chi_i) \in \Syl_q(G(\chi_1,\chi_2, \dots,\chi_i)), \\
\ma(\chi_1,\dots,\chi_n) = P_{2k'}^*= P_2 \cdot P_4 \cdots P_{2k'}
 &\text{ and } 
 \mb(\chi_1,\dots,\chi_n)= Q_{2l'-1}^* =Q_1 \cdot Q_{3} \cdots Q_{2l'-1}, 
\end{align}
for each $i = 1,2,\dots,n$.  Therefore 
\begin{equation}\mylabel{ma.C'}
 G(\chi_1, \chi_2, \dots, \chi_i) = \ma(\chi_1,\chi_2, \dots,
\chi_i)\mb(\chi_1, \chi_2, \dots, \chi_i) 
\end{equation}
for $i = 1,2,\dots,n$.      
\end{subequations}
 
We are going to perform a series of linear reductions of a very special form.
 We set $S = 1$, and let $\zeta$ be the
trivial linear character of $S$. We also set $R = 1$, and let
$\eta$ be the trivial character of $R$. Then $(G, S, \zeta, G_1,
\chi_1)$ is a quintuple satisfying \eqref{lmt.e3}.
 So we may pass to an arbitrary
$``\ma(\beta_1)"$-invariant faithful linear limit
 $(\mbgi{1}, \mbsi{1}, \bzi{1}, 
\mbgi{1}_1, \Ti{1}_1)$
of the quintuple  $(G,S, \zeta, G_1, \chi_1)$.
Of course along with that we obtain (see \eqref{lmt.e10})
a normal series
\begin{equation}\mylabel{ma.e3}
1 \unlhd \mbqi{1}_1= \mbgi{1}_1 \unlhd \dots \unlhd \mbgi{1}_n = \mbgi{1},
\end{equation}
of $\mbgi{1}$, from the series \eqref{ma.A}. In addition 
we  reach  a tower of characters 
$\Ti{1}_i  \in \Irr(\mbgi{1}_i)$, for $i=0, 1, \dots, n$, such that 
$\Ti{1}_0 = 1$ and  $\Ti{1}_i$ lies under $\Ti{1}_j$ and above $\bzi{1}$, 
whenever  $1 \leq i \leq j \leq n$.
We also   get a triangular set
 \begin{equation}\mylabel{ma.e4}
\{ \mbqi{1}_{2i-1}, \mbpi{1}_{2r} |\bbi{1}_{2i-1}, \bai{1}_{2r}\}_{i=1, 
r=0}^{k', l'}
\end{equation}
for \eqref{ma.e3} that corresponds uniquely to the tower 
$\{ \Ti{1}_i \} _{i=0}^n$ and satisfies Theorem \ref{lmt.t3}.
 Furthermore, the original  Sylow system $\{\ma, \mb\}$ for $G$, 
(see \eqref{ma.C}), provides a Sylow system $\{ \mai{1}, \mbi{1} \}$
for $\mbgi{1}$, that also 
 satisfies Theorem \ref{lmt.t3}.

\noindent 
 Obviously the trivial group $R$ is still a subgroup of $\mbpi{1}_2$, and 
its trivial character $\eta$ lies under $\bai{1}_2$. We denote by 
$\mbri{1}= 1= R$ the trivial group, seen inside $\mbpi{1}_2$, and 
by $\bei{1} \in \Irr(\mbri{1})$ its trivial character. Hence 
\begin{subequations}\mylabel{m.rs1}
\begin{gather}
1 \times \mbsi{1}=\mbri{1} \times \mbsi{1} 
\text{ is a central subgroup of } \mbgi{1}\\
1\times \bzi{1}= \bei{1} \times \bzi{1} 
\in \Irr(\mbri{1} \times \mbsi{1}).
\end{gather}
If  $n \geq 2$ then  in addition we have 
\begin{gather}
\mbri{1} \times \mbsi{1} \unlhd \mbpi{1}_2 \cdot \mbqi{1}_1= \mbgi{1}_2  \\
\bei{1} \text{ and }   \bzi{1}
 \text{ lie under 
$\bai{1}_2$ and  $\bbi{1}_1$, respectively } .
\end{gather}
\end{subequations}
Notice that, as  $(\mbgi{1}, \mbsi{1}, \bzi{1}, 
\mbgi{1}_1, \Ti{1}_1)$ is a faithful linear limit of 
 $(G,S, \zeta, G_1, \chi_1)$, Corollaries \ref{lim.co2} and \ref{lim.co22} 
imply
\begin{remark}\mylabel{ma.r1}
$\mbsi{1}= Z(\mbgi{1}_1)$ is a  cyclic central subgroup of $\mbgi{1}$,
maximal among the abelian $\mbgi{1}$-invariant 
subgroups of $\mbgi{1}_1$. Furthermore, 
the character $\bbi{1}_1= \Ti{1}_1$ is $\mbgi{1}$-invariant.
\end{remark}  
Note also that Corollary \ref{lim.co1} easily implies
\begin{remark}\mylabel{ma.rrr1}
Any faithful linear limit of $(\mbgi{1}, \mbsi{1}, \bzi{1}, 
\mbgi{1}_i, \Ti{1}_i)$ is also a faithful linear limit of 
$(G, 1, 1, G_i , \chi_i)$, for all $i=1,2, \dots, n$.
\end{remark}

As far as the monomial characters 
of $G$ are concerned we have 
\begin{proposition}\mylabel{ma.p1}
Any character $\Ti{1} \in \Irr(\mbgi{1})$  that lies above 
$1\times \bzi{1}= \bei{1} \times  \bzi{1} \in \Irr(\mbri{1} \times \mbsi{1})$ is monomial.  
\end{proposition}
\begin{proof}
Let $\Ti{1}$ be an irreducible character of $\mbgi{1}$ that lies 
above $\bei{1} \times \bzi{1}$, and thus above $\bzi{1}$. 
According to Lemma \ref{lim.l3}, 
there exists an irreducible character $\chi \in \Irr(G)$,  
that lies above $\zeta=1$ and satisfies $\Ti{1}= fl(\chi)$, that is, $\Ti{1}$ 
is the faithful linear limit of $\chi$. But $\chi$ is monomial, as 
 $G$ is a monomial group.
Therefore,  Proposition \ref{lim.p1} implies that 
 $\Ti{1}$  is also monomial. 
\end{proof}

The first critical result, which we will prove in Section \ref{ma.se2}, is
\begin{theorem}\mylabel{ma.t1}
After the above reduction, the group $\mbgi{1}_2$ is
nilpotent, if it exists, i.e., if $n \ge 2$.
\end{theorem}
This theorem implies  that $\mbgi{1}_2$ is the 
direct product
\begin{subequations}\mylabel{ma.e2}
\begin{equation}
\mbgi{1}_2=\mbpi{1}_2 \times \mbqi{1}_1 = \mbpi{1}_2 \times \mbgi{1}_1
\end{equation}
of its $p$-Sylow subgroup $\mbpi{1}_2$  and its $q$-Sylow subgroup
 $\mbqi{1}_1 = \mbgi{1}_1$.
 It also implies that $\Ti{1}_2$  is the direct product
\begin{equation}
\Ti{1}_2= \bai{1}_2 \times \bbi{1}_1
\end{equation}
of $\bai{1}_2 \in \Irr(\mbpi{1}_2)$ and
 $\bbi{1}_1 = \Ti{1}_1 \in \Irr(\mbqi{1}_1)= 
\Irr(\mbgi{1}_1)$.
\end{subequations}
Therefore, in the case of $n \geq 2$, the relations in 
(\ref{m.rs1}c,d) imply
\begin{gather}\mylabel{m.rs11}
\mbri{1} \times \mbsi{1} \unlhd \mbpi{1}_2 \times  \mbqi{1}_1 \notag  \\
\bei{1} \times  \bzi{1}
 \text{ lies under $\bai{1}_2 \times \bbi{1}_1 = \Ti{1}_2$ } .
\end{gather}

Furthermore, the character  $\bbi{1}_1$ is $\mbgi{1}$-invariant.
Hence the  normal series \eqref{ma.e3}, along with its character tower
$\{\Ti{1}_i\}_{i=0}^n$, satisfies \eqref{sh.e1} of Section \ref{sh}.
So we shift the series by one, and consider the series
\begin{subequations}\mylabel{ma.sh1}
\begin{equation}\mylabel{ma.sh1a}
1\unlhd \mbpi{1}_2 \unlhd \mbgi{1}_3 
\unlhd \dots \unlhd \mbgi{1}_n = \mbgi{1}, 
\end{equation}
of normal subgroups of $\mbgi{1}$, that satisfies Hypothesis \ref{hyp1}.
 Then, as we have seen in
 Section \ref{sh}, and in particular Theorem \ref{sh.t1},  the characters 
\begin{equation}\mylabel{ma.sh1b}
1, \bai{1}_2, \Ti{1}_3, \dots, \Ti{1}_n,
\end{equation}
form a tower for the above series, with corresponding triangular set 
\begin{equation}
\{ \mbqi{1}_{2i-1}, \mbpi{1}_{2r}| \bbi{1}_{2i-1}, 
\bai{1}_{2r}\}_{i=2, r=1}^{l', k'}.
\end{equation}
\end{subequations}
(Note that we have dropped the first $q$ and $p$ groups,  
 $\mathbb{Q}^{(1),s}_1$ and  $ \mathbb{P}^{(1),s}_0$ (see \eqref{sh.e3c}
 and \eqref{sh.e3cc}) respectively, along with their characters,   
as these are assume trivial for the shifted system.) 
Furthermore, as Theorem \ref{sh.t1} implies,  
the above shifted system and the one for \eqref{ma.e3} have in common 
 the Sylow system $\{ \mai{1}, \mbi{1}\}$, i.e., this  Sylow 
system satisfies
\eqref{sy.e2} for  \eqref{ma.sh1}.

The quintuple  $( \mbgi{1}, \mbri{1}, \bei{1}, \mbpi{1}_2,
 \bai{1}_2)$ is clearly a linear one.
Therefore we may pass to a $\mbi{1}(\bai{1}_2)$-invariant 
faithful linear limit 
$(\mbgi{2}, \mbri{2}, \bei{2}, \mbpi{2}_2, \bai{2}_2 )$
of the former  quintuple. So, (see \eqref{lmt.e18}), 
 the chain \eqref{ma.sh1a} reduces 
to a chain 
\begin{subequations}\mylabel{ma.eE}
\begin{equation}\mylabel{ma.eE1}
1\unlhd  \mbpi{2}_2\unlhd \mbgi{2}_3 \unlhd 
\mbgi{2}_4 \unlhd  \dots  \unlhd \mbgi{2}_n = \mbgi{2}, 
\end{equation}
of normal subgroups $\mbgi{2}_i$  of $\mbgi{2}$. The character tower
\eqref{ma.sh1b}  reduces  (see \eqref{lmt.e18b})  to the tower 
\begin{equation}\mylabel{ma.eE2}
\{1, \bai{2}_2,  \Ti{2}_i \}_{i=3}^n 
\end{equation}
where  $\bai{2}_2 \in \Irr(\mbpi{2}_2)$ and  
$\Ti{2}_i  \in \Irr(\mbgi{2}_i )$ , for all $i=3,\dots, n$.
According to the Theorem  \ref{lmt.t7},
the Sylow system $\{ \mai{1} , \mbi{1}\}$ for $\mbgi{1}$ reduces
 to a Sylow system 
$\{ \mai{2}, \mbi{2} \}$  for $\mbgi{2}$. Furthermore, the same theorem 
provides a unique, up to conjugation,   triangular set 
\begin{equation}\mylabel{ma.eE3}
\{ \mbqi{2}_{2i-1}, \mbpi{2}_{2r} |\bbi{2}_{2i-1},
 \bai{2}_{2r}\}_{i=2,r=1}^{l', k'},
\end{equation}
that corresponds to  \eqref{ma.eE2}.
(Note that we have dropped  two trivial groups 
and their characters.)  
\end{subequations}

Clearly,   the characters  $\bai{2}_2$ and $\Ti{2}_i$ lie above 
the limit character $\bei{2}$, for all $i=3,\dots,n$.   
As $(\mbgi{2}, \mbri{2}, \bei{2}, \mbpi{2}_2, \bai{2}_2 )$
is a faithful linear limit of  $( \mbgi{1}, \mbri{1}, \bei{1}, \mbpi{1}_2,
 \bai{1}_2)$, Corollaries \ref{lim.co2} and \ref{lim.co22} imply 
\begin{remark}\mylabel{ma.r2}
$\mbri{2}= Z(\mbpi{2}_2)$ is a cyclic central subgroup of $\mbgi{2}$, 
maximal among the abelian $\mbgi{2}$-invariant subgroups of $\mbpi{2}_2$, 
while the character $\bai{2}_2$ is $\mbgi{2}$-invariant.
\end{remark}
In addition, the way we perform the linear reductions, along with 
 Remark \ref{ma.rrr1} and  Corollary \ref{lim.co1},  implies 
\begin{remark}\mylabel{ma.rrr2}
Any faithful linear limit of  $(\mbgi{2}, \mbri{2}, \bei{2}, 
\mbgi{2}_i, \Ti{2}_i)$ is also a faithful linear limit of 
$(G, 1, 1, G_i, \chi_i)$, for all $i=3,4, \dots, n$.
\end{remark}

Furthermore, the fact that $\mbsi{1}$ is a normal subgroup of $\mbgi{1}$
 that centralizes $\mbpi{1}_2$  implies, by  Remarks \ref{lm.rem1}
 and \ref{lm.rem2}, that there exists a subgroup $\mbsi{2}$ of $\mbgi{2}$ 
such that 
\begin{equation}\mylabel{ma.rs1b} 
\mbsi{1} \cong \mbsi{2} \unlhd \mbgi{2}
\end{equation}
and $\mbsi{2}$ centralizes $\mbpi{2}_2$. Thus, it also centralizes
 $\mbri{2} \leq  \mbpi{2}_2$. In addition, $\mbsi{2}$ is a central subgroup
 of $\mbgi{2}$, as $\mbsi{1}$ is a central subgroup 
of $\mbgi{1}$, and $\mbgi{2}$ is a section of $\mbgi{1}$.
 Under the isomorphism in \eqref{ma.rs1b}, the 
irreducible character $\bzi{1} \in \Irr(\mbsi{1})$
maps to an irreducible character $\bzi{2} \in \Irr(\mbsi{2})$.  
If $n \geq 3$ we can say more about $\mbsi{2}$ and its character
 $\bzi{2}$.
Indeed,  the fact we used a $\mbi{1}(\bai{1}_2)$-invariant 
faithful limit to get \eqref{ma.eE1}, implies 
(see Theorem \ref{lmt.t7}) that $\mbqsi{1}_{2l-1} \cong \mbqsi{2}_{2l-1}$
for all $l=2, \dots, l'$. In particular we have 
$$
\mbqi{1}_3 \cong \mbqi{2}_3.
$$
The group $\mbsi{1}$ is a subgroup of $\mbqi{1}_1$. But the 
latter is a subgroup of $\mbqi{1}_3$, as $\mbgi{1}_2 = \mbpi{1}_2 
\times \mbqi{1}_1 \leq \mbgi{1}_3$.  Hence $\mbsi{1}$ is a subgroup
 of $\mbqi{1}_3$. 
We conclude that its isomorphic image $\mbsi{2}$ is a 
subgroup of $\mbqi{2}_3$.
Furthermore, its irreducible character $\bzi{1} \in \Irr(\mbsi{1})$ 
lies under $\bbi{1}_1$. But $\bbi{1}_1$ lies under $\bbi{1}_3$,
 as $\mbqi{1}_1 \leq \mbqi{1}_3$. Hence $\bzi{1}$ lies under $\bbi{1}_3$. 
As  $\bzi{1}$ maps to $\bzi{2} \in \Irr(\mbsi{2})$, while  
$\bbi{1}_3$ maps to $\bbi{2}_3 \in \Irr(\mbqi{2}_3)$,   we conclude that 
$\bzi{2}$ lies under $\bbi{2}_3$.
Hence 
\begin{subequations}\mylabel{m.rs2}
\begin{gather}
\mbri{2} \times \mbsi{2} \text{ is a central subgroup of } 
 \mbgi{2} \\
\bei{2} \times \bzi{2} \in \Irr(\mbri{2} \times \mbsi{2}).
\end{gather}
If in addition $ n \geq 3$ then  
 \begin{gather}
\mbri{2} \times \mbsi{2} \unlhd \mbpi{2}_2 \cdot \mbqi{2}_3= \mbgi{2}_3, \\
  \bei{2}  \text{ and }   \bzi{2} \text{ lie under $ \bai{2}_2$ and 
 $\bbi{2}_3$,  respectively. } 
\end{gather}
\end{subequations}
Now we can extend  Proposition \ref{ma.p1}  to
\begin{proposition}\mylabel{ma.p2}
Every irreducible character $\Ti{2} \in \Irr(\mbgi{2})$ that lies above
$\bei{2} \times \bzi{2}  \in \Irr(\mbri{2} \times \mbsi{2})$ is monomial.
\end{proposition}
\begin{proof}
Obviously any  $\Ti{2} \in \Irr(\mbgi{2} |\bei{2} \times \bzi{2})$
 lies above  $\bei{2} \in \Irr(\mbri{2})$. As the quintuple
$(\mbgi{2}, \mbri{2}, \bei{2}, \mbpi{2}_2, \bai{2}_2 )$
is a faithful linear limit of 
$(\mbgi{1}, \mbri{1}, \bei{1}, \mbpi{1}_2, \bai{1}_2 )$,
Lemma \ref{lim.l3} implies  
the existence of an irreducible character $\Ti{1} \in \Irr(\mbgi{1})$, lying  
above $\bei{1}$, so that $\Ti{2}$ is a faithful linear limit of $\Ti{1}$,
under the $\mbi{1}(\bai{1}_2)$-invariant reductions we performed.
As we have already seen, under those reductions
 the $q$-subgroup $\mbsi{1}$ of $\mbgi{1}$ 
maps isomorphically to the subgroup    $\mbsi{2}$ of $\mbgi{2}$. 
Hence the only way the faithful limit character $\Ti{2} \in \Irr(\mbgi{2})$ 
can lie above $\bzi{2} \in \Irr(\mbsi{2})$, is if $\Ti{1}$ lies above
$\bzi{1} \in \Irr(\mbsi{1})$.
In conclusion $\Ti{1} \in \Irr(\mbgi{1})$ lies above 
$\bei{1} \times \bzi{1} \in \Irr(\mbri{1} \times \mbsi{1})$. 
Now we can apply Proposition \ref{ma.p1} to conclude that $\Ti{1}$ is
 monomial. But $\Ti{2}$ is a faithful linear limit of $\Ti{1}$. 
Hence Proposition \ref{ma.p1} implies that $\Ti{2}$ is also monomial.
This completes the proof of the proposition.
\end{proof}

The next   important theorem, that is proved in Section  \ref{ma.se3} is 
 \begin{theorem}\mylabel{ma.t2}
After the above reductions, the group $\mbgi{2}_3$ is nilpotent 
if it exists, i.e., if $n \ge 3$.
\end{theorem}
Hence $\mbgi{2}_3$ is the direct product
\begin{subequations}\mylabel{ma.e5}
\begin{equation}\mylabel{ma.e5a}
\mbgi{2}_3 = \mbpi{2}_2 \times \mbqi{2}_3, 
\end{equation}
of its $p$-Sylow subgroup $\mbpi{2}_2$ and its $q$-Sylow subgroup
 $\mbqi{2}_3$. Furthermore, its irreducible character $\Ti{2}_3$ is the 
direct product
\begin{equation}\mylabel{ma.e5b}
\Ti{2}_3 = \bai{2}_2 \times \bbi{2}_3, 
\end{equation}
of $\bai{2}_2 \in \Irr(\mbpi{2}_2)$ and $\bbi{2}_3 \in \Irr(\mbqi{2}_3)$.
\end{subequations}
Hence \eqref{m.rs2} in the case of $n \geq 3$ becomes
\begin{gather}\mylabel{m.rs22}
\mbri{2} \times \mbsi{2} \unlhd \mbpi{2}_2 \times \mbqi{2}_3= \mbgi{2}_3, \\
  \bei{2}\times  \bzi{2} \text{ lies under $ \bai{2}_2 \times \bbi{2}_3$. } 
\end{gather}
The fact that $\mbgi{2}_3$ is a nilpotent group permits us to shift 
the series \eqref{ma.eE1} by one, and apply all the results of Section
 \ref{sh}. (Note that the roles of $p$ and $q$ are interchanged.)
 Thus we get the series 
\begin{subequations}\mylabel{ma.sh2}
\begin{equation}
1 \unlhd \mbqi{2}_3 \unlhd \mbgi{2}_4 \unlhd \dots \unlhd \mbgi{2}_n 
= \mbgi{2}, 
\end{equation}
of normal subgroups of $\mbgi{2}$.
Then, according to Section  \ref{sh} and, in particular, Theorem \ref{sh.t1},
the characters 
\begin{equation}
1, \bbi{2}_3, \Ti{2}_4, \dots, \Ti{2}_n,
\end{equation}
form a tower for the above series,  with 
corresponding triangular set 
\begin{equation}\mylabel{ma.sh2c}
\{  \mbpi{2}_{2r}, \mbqi{2}_{2i-1} |
 \bai{2}_{2r}, \bbi{2}_{2i-1} \}_{r=2,i=2}^{k',l'}.  
\end{equation}
(As expected, (see \eqref{sh.e3c} and \eqref{sh.e3cc}), 
  there are 3  trivial groups (the $\mathbb{P}^{(2),s}_{0},
\mathbb{P}^{(2),s}_2$ and  
$\mathbb{Q}^{(2),s}_1$)   in \eqref{ma.sh2c} 
that have been  dropped.)
Furthermore, the Hall system $\{ \mai{2}, \mbi{2}\}$ for $\mbgi{2}$, that
 was obtained via the second faithful linear limit, satisfies the
 equivalent of (\ref{ma.C}-f) for \eqref{ma.sh2} (see Theorem \ref{sh.t1}).
\end{subequations}
Of course the groups $\mbsi{2}$ and $\mbri{2}$ remain central subgroups 
of $\mbgi{2}$ that satisfy \eqref{m.rs2}.

At this point 
 we can repeat the process from the beginning, with $\mbri{2}$ and 
$\mbsi{2}$ in the place of $R$ and $S$ respectively.
(So the next  step would be to take  an
 $\mai{2}(\bbi{2}_3)$-invariant
 faithful linear limit $(\mbgi{3}, \mbsi{3}, \bzi{3}, \mbqi{3}_3, \bbi{3}_3)$
  of the quintuple  $(\mbgi{2},\mbsi{2}, \bzi{2}, \mbqi{2}_3, \bbi{2}_3)$.)

Suppose,  for the sake of our inductive hypothesis, that 
 we have repeated this process $t-1$ times, 
for some integer $t$ with $2<  t \leq  n$, i.e., we have taken $t-1$
invariant  faithful limits and, after each such, we have shifted 
our series by one.
So we arrive at a group $\mbgi{t-1}=\mbgi{t-1}_n$, that is a section of
the original group $G$. Note that according to 
 the inductive hypothesis, the last
group  shown to be nilpotent is the group $\mbgi{t-1}_{t}$. 
Therefore, depending on the parity of $t$, the group $\mbgi{t-1}_t$
equals 
\begin{subequations}\mylabel{ma.pt}
\begin{align}
\mbgi{t-1}_t &= \mbpi{t-1}_t \times \mbqi{t-1}_{t-1},
 &\text{ when $t$ is even } \\
\mbgi{t-1}_{t} &= \mbpi{t-1}_{t-1} \times \mbqi{t-1}_t, 
&\text{ when $t$ is odd. }
\end{align}
\end{subequations}
According to the  inductive hypothesis, we can also generalize 
\eqref{m.rs11} and Remark \ref{ma.r2}. Thus we can assume 
that 
\begin{subequations}\mylabel{ma.e7}
\begin{gather}
\mbri{t-1} \times \mbsi{t-1} \unlhd \mbpi{t-1}_t \times \mbqi{t-1}_{t-1} = 
\mbgi{t-1}_t  \mylabel{ma.e7a}\\
\bei{t-1} \times \bzi{t-1} \text{ lies under } 
\bai{t-1}_t \times \bbi{t-1}_{t-1}= \Ti{t-1}_t, \mylabel{ma.e7b}\\
\mbsi{t-1} = Z(\mbqi{t-1}_{t-1}) \text{ and } 
\bbi{t-1}_{t-1} \text{ is $\mbgi{t-1}$-invariant, }\mylabel{ma.e7c}
\end{gather}
in the case of an even $t$. 
If $t$ is odd then 
\begin{gather}
\mbri{t-1} \times \mbsi{t-1} \unlhd \mbpi{t-1}_{t-1} \times \mbqi{t-1}_{t}=
\mbgi{t-1}_t \mylabel{ma.e7aa}\\
\bei{t-1} \times \bzi{t-1} \text{ lies under } 
\bai{t-1}_{t-1} \times \bbi{t-1}_{t}= \Ti{t-1}_t, \mylabel{ma.e7bb}\\
\mbri{t-1} = Z(\mbpi{t-1}_{t-1}) \text{ and }
\bai{t-1}_{t-1} \text{ is $\mbgi{t-1}$-invariant. }\mylabel{ma.e7cc}
\end{gather}
\end{subequations}
Furthermore, the last group   dropped is 
the $q$-group $\mbqi{t-1}_{t-1}$ in the case of an even $t$, 
or  the $p$-group $\mbpi{t-1}_{t-1}$ in the case of an odd $t$.

{\bf  Case 1: $t$ is even } 

Assume first that $t$ is even. So we reach the series (after the $q$-group 
$\mbqi{t-1}_{t-1}$ is dropped)
\begin{subequations}\mylabel{ma.e6}
\begin{equation}\mylabel{ma.e6a}
1 \unlhd \mbpi{t-1}_t  \unlhd \mbgi{t-1}_{t+1} \unlhd 
\dots \unlhd  \mbgi{t-1}_n = \mbgi{t-1},
\end{equation}
of normal subgroups of $\mbgi{t-1}$. Then  
$\mbgi{t-1}_i/ \mbgi{t-1}_{i-1}$ is a $p$-group if $i$ is even, and a 
$q$-group if $i$ is odd, for each $i=t+2, \dots, n$, while  for $i=t+1$
we get $\mbgi{t-1}_{t+1}/ \mbpi{t-1}_t$ is a $q$-group with
 $\mbpi{t-1}_t$ a $p$-group.  
Along with the above series, we reach the characters  
\begin{equation}\mylabel{ma.e6b}
1, \bai{t-1}_t \in \Irr(\mbpi{t-1}_{t}),   \Ti{t-1}_i \in \Irr(\mbgi{t-1}_i)
\end{equation}
for each  $i=t+1, \dots, n$, that form a tower for \eqref{ma.e6a}. 
Furthermore, we have the triangular set
\begin{equation}
\mylabel{ma.e6c}
\{ \mbqi{t-1}_{2i-1}, \mbpi{t-1}_{2r} | \bbi{t-1}_{2i-1},
 \bai{t-1}_{2r} \}_{i=(t/2) +1, r=t/2}^{l', k'}
\end{equation}
that corresponds uniquely, up to conjugation,   to \eqref{ma.e6b}.
(Note the first $t$ groups in \eqref{ma.e6c}, have been dropped 
 as they are   trivial.)
Also the Sylow system $\{ \ma, \mb \}$ has been transfered to a Sylow system 
$\{ \mai{t-1}, \mbi{t-1} \}$  of $\mbgi{t-1}$ that satisfies the  properties 
in \eqref{sy.e2}, for \eqref{ma.e6}.
\end{subequations}
 In addition, we reach  two central 
 subgroups of $\mbgi{t-1}$, the $p$-group $\mbri{t-1}$ 
and the $q$-group $\mbsi{t-1}$, along with
 their characters $\bei{t-1} \in \Irr(\mbri{t-1})$ and $\bzi{t-1}
 \in \Irr(\mbsi{t-1})$. We assume that $n \geq t$, so that our inductive step
will be the $t$-th  step. 
So  we get the linear quintuple (see \eqref{ma.e7})
\begin{equation}
(\mbgi{t-1}, \mbri{t-1}, \bei{t-1}, \mbpi{t-1}_t, \bai{t-1}_t). 
\end{equation}
Furthermore, our inductive hypothesis implies that  every irreducible 
character $ \Ti{t-1} \in \Irr(\mbgi{t-1})$ that lies above 
$\bei{t-1} \times \bzi{t-1} \in \Irr(\mbri{t-1} \times \mbsi{t-1})$
 is monomial.
Also  any faithful linear limit of the quintuple 
$(\mbgi{t-1}, \mbri{t-1}, \bei{t-1}, \mbgi{t-1}_i, \Ti{t-1}_i)$ or
 the
$(\mbgi{t-1}, \mbsi{t-1}, \bzi{t-1}, \mbgi{t-1}_i, \Ti{t-1}_i)$ 
is also a faithful linear limit of 
$(G, 1, 1, G_i, \chi_i)$ for all $i=t-1, t, \dots, n$.

For the inductive step, 
 we take a $\mbi{t-1}(\bai{t-1}_t)$-invariant faithful
 linear limit 
$$
(\mbgi{t}, \mbri{t}, \bei{t}, \mbpi{t}_t, \bai{t}_t)
$$ 
 of $(\mbgi{t-1}, \mbri{t-1}, \bei{t-1}, \mbpi{t-1}_t, \bai{t-1}_t) $.
The series  
 \eqref{ma.e6a} reduces to  
\begin{subequations}\mylabel{ma.e8}
\begin{equation}\mylabel{ma.e8a}
1 \unlhd \mbpi{t}_t \unlhd \mbgi{t}_{t+1} \unlhd \dots \unlhd
 \mbgi{t}_n = \mbgi{t}. 
\end{equation}
 Furthermore, the character tower \eqref{ma.e6b}
 reduces to the tower 
\begin{equation}
\mylabel{ma.e8b}
1, \bai{t}_t,   \Ti{t}_i \in \Irr(\mbgi{t}_i)
\end{equation}
Along with the above character tower we get a triangular set 
\begin{equation}
\mylabel{ma.e8c}
\{ \mbqi{t}_{2i-1}, \mbpi{t}_{2r} | \bbi{t}_{2i-1},
 \bai{t}_{2r} \}_{i=(t/2) +1, r=t/2}^{l', k'}
\end{equation}
that satisfies Theorem \ref{lmt.t7}.  Hence 
 the Sylow system $\{\mai{t-1}, \mbi{t-1}\}$ 
for $\mbgi{t-1}$, reduces to a Sylow system 
$\{ \mai{t}, \mbi{t} \}$ for $\mbgi{t}$, that satisfies \eqref{sy.e2} 
for the  above reduced $t$-system.
\end{subequations}
As with Remarks \ref{ma.r1} and \ref{ma.r2}, the fact that 
we have taken faithful linear limits, 
  along with Corollaries \ref{lim.co2} and 
\ref{lim.co22}, implies 
\begin{remark}\mylabel{ma.r3}
The group  $\mbri{t}= Z(\mbpi{t}_t)$ is a cyclic 
central  subgroup of $\mbgi{t}$, maximal among the abelian 
$\mbgi{t}$-invariant subgroups of $\mbgi{t}_t$, 
 while the character $\bai{t}_t$ is 
$\mbgi{t}$-invariant.
\end{remark}
In addition, the way the reductions are done,  along with 
Corollary \ref{lim.co1}, implies 
\begin{remark}\mylabel{ma.rrr3} 
Any faithful linear limit of the quintuple 
 $(\mbgi{t}, \mbri{t}, \bei{t}, \mbgi{t}_i, \Ti{t}_i)$ or
 the quintuple 
$(\mbgi{t}, \mbsi{t}, \bzi{t}, \mbgi{t}_i, \Ti{t}_i)$ 
is also a faithful linear limit of 
$(G, 1, 1, G_i, \chi_i)$ for all $i=t,t+1,  \dots, n$.
\end{remark}

Even more,  Theorem \ref{lmt.t7} implies that 
\begin{subequations}\mylabel{ma.e20}
\begin{equation}\mylabel{ma.e200}
\mbqsi{t-1}_{2l-1} \cong \mbqsi{t}_{2l-1}, 
\end{equation}
 for all $l=t+1,\dots, l'$, 
where $\mbqsi{t-1}_{2l-1}$ and $\mbqsi{t}_{2l-1}$ denote the product  groups  
$\mbqsi{t-1}_{2l-1}= \mbqi{t-1}_{t+1} \cdot \mbqi{t-1}_{t+3} \cdots 
\mbqi{t-1}_{2l-1}$  and 
$\mbqsi{t}_{2l-1}= \mbqi{t}_{t+1} \cdot \mbqi{t}_{t+3} \cdots 
\mbqi{t}_{2l-1}$, respectively.
In particular,
\begin{equation}\mylabel{ma.e20a}
\mbqi{t}_{t+1} \cong \mbqi{t-1}_{t+1}.
\end{equation}
 (Observe that $\mbqsi{t-1}_{t+1}= \mbqi{t-1}_{t+1}$ and 
$\mbqsi{t}_{t+1}= \mbqi{t}_{t+1}$.)
The group $\mbsi{t-1}$ is a central  subgroup of $\mbgi{t-1}$ that 
centralizes $\mbpi{t-1}_t$ (see \eqref{ma.e7a}). Hence Remarks
\ref{lm.rem1} and  \ref{lm.rem2}
imply  that $\mbsi{t-1}$ maps isomorphically to a normal  subgroup
\begin{equation}\mylabel{ma.e20b}
\mbsi{t} \cong \mbsi{t-1}
\end{equation} 
of $\mbgi{t}$, that centralizes $\mbpi{t}_t$. 
In addition, $\mbsi{t}$ is a central subgroup of 
$\mbgi{t}$, as $\mbsi{t-1}$ is a central subgroup 
of $\mbgi{t-1}$, and $\mbgi{t}$ is a section of $\mbgi{t-1}$. 
Furthermore, under  the group  isomorphism in \eqref{ma.e20b}, the 
irreducible character $\bzi{t-1} \in \Irr(\mbsi{t-1})$ maps 
to 
\begin{equation}\mylabel{ma.e20c}
\bzi{t} \in \Irr(\mbsi{t}).
\end{equation} 
\end{subequations}
As  $n \geq t+1$,  the group  $\mbqi{t-1}_{t+1}$ exists. 
Furthermore, as $\mbgi{t-1}_t = \mbpi{t-1}_t \times \mbqi{t-1}_{t-1}$ 
is nilpotent we get that 
$\mbqi{t-1}_{t-1}$ is a subgroup
 of $\mbqi{t-1}_{t+1}$, while its irreducible character
$\bbi{t-1}_{t-1} \in \Irr(\mbqi{t-1}_{t-1}) $ lies under
$\bbi{t-1}_{t+1} \in \Irr(\mbqi{t-1}_{t+1})$. 
  Hence $\mbsi{t-1} \leq \mbqi{t-1}_{t-1}$ is a
 subgroup  of  $\mbqi{t-1}_{t+1}$. This along with \eqref{ma.e20} implies 
that $\mbsi{t}$ is a subgroup of $\mbqi{t}_{t+1}$. 
Even more, its irreducible character $\bzi{t} \in \Irr(\mbsi{t})$  
 lies under $\bbi{t}_{t+1}$, as $\bzi{t-1}$ lies under 
$\bbi{t-1}_{t-1}$ (see \eqref{ma.e7b}) and that under $\bbi{t-1}_{t+1}$.
Hence 
\begin{proposition}\mylabel{ma.p4}
We have two central subgroups of $\mbgi{t}$, the $p$-group 
$\mbri{t}$ and the $q$-group $\mbsi{t}$. 
Along with them we get their irreducible characters 
$\bei{t} \in \Irr(\mbri{t})$ and $\bzi{t} \in \Irr(\mbsi{t})$. 
These  groups and characters satisfy
\begin{align*}
\mbri{t} \times  \mbsi{t} \unlhd  \mbpi{t}_{t}\cdot
 \mbqi{t}_{t+1}= \mbgi{t}_{t+1}, \\
\bei{t} \text{ and }  \bzi{t} \text{ lie under } \bai{t}_t
 \text{ and $\bbi{t}_{t+1}$, 
respectively.}
\end{align*}
\end{proposition}

We can also show 
\begin{proposition}\mylabel{ma.p5}
Every irreducible character $\Ti{t} \in \Irr(\mbgi{t})$ that lies above 
$\bei{t} \times \bzi{t} \in \Irr(\mbri{t} \times \mbsi{t})$, is  
monomial.
\end{proposition}
\begin{proof}
The quintuple 
$(\mbgi{t}, \mbri{t}, \bei{t}, \mbpi{t}_t, \bai{t}_t )$
is a faithful linear limit of the linear quintuple 
$(\mbgi{t-1}, \mbri{t-1}, \bei{t-1}, \mbpi{t-1}_t, \bai{t-1}_t )$.
If $\Ti{t} \in \Irr(\mbgi{t})$ lies above $\bei{t} \times \bzi{t}$, then 
it lies above $\bei{t}$, 
 hence 
Lemma \ref{lim.l3} implies  
the existence of an irreducible character $\Ti{t-1} \in \Irr(\mbgi{t-1})$ 
above $\bei{t-1}$, so that $\Ti{t}$ is a faithful linear limit of $\Ti{t-1}$,
under the $\mbi{t-1}(\bai{t-1}_t)$-invariant reductions we performed.
As we have already seen, under those reductions,
 the $q$-subgroup $\mbsi{t-1}$ of $\mbgi{t-1}$ 
maps isomorphically to the subgroup    $\mbsi{t}$ of $\mbgi{t}$, 
see \eqref{ma.e20b}. 
Hence the only way the faithful limit character $\Ti{t} \in \Irr(\mbgi{t})$ 
can lie above $\bzi{t} \in \Irr(\mbsi{t})$, is if $\Ti{t-1}$ lies above
$\bzi{t-1} \in \Irr(\mbsi{t-1})$.
In conclusion $\Ti{t-1} \in \Irr(\mbgi{t-1})$ lies above 
$\bei{t-1} \times \bzi{t-1} \in \Irr(\mbri{t-1} \times \mbsi{t-1})$. 
According to our inductive hypothesis, $\Ti{t-1}$ is monomial.
As $\Ti{t}$ is  a faithful linear limit of $\Ti{t-1}$,
 Proposition \ref{ma.p1} implies that $\Ti{t}$ is also monomial.
This completes the proof of the proposition.
\end{proof}

The main theorem for the inductive step, in the case of an even $t$,  will be
\begin{theorem}\mylabel{ma.main}
After the above reductions, the group $\mbgi{t}_{t+1}$ 
is nilpotent if it exists, i.e., if $n \ge t+1$.
\end{theorem}

 Observe that,
as before, the fact $\mbgi{t}_{t+1}$ is nilpotent  
 allows us to shift the series \eqref{ma.e8a} 
 by one. So, provided that  $t+2 \leq n$,  we get the series 
\begin{subequations}\mylabel{ma.e9}
\begin{equation}\mylabel{ma.e9a}
1 \unlhd \mbqi{t}_{t+1} \unlhd \mbgi{t}_{t+2} \unlhd \dots \unlhd
 \mbgi{t}_n = \mbgi{t}.
\end{equation}
The character tower and its  corresponding triangular set are carried
 over to the shifted series, as in  the Sylow system.
Hence the characters 
\begin{equation}
1, \bbi{t}_{t+1}, \Ti{t}_{t+2}, \dots , \Ti{t}_n 
\end{equation}
form a tower for the series \eqref{ma.e9a}.
 Its  corresponding triangular set is  
\begin{equation}
\mylabel{ma.e9c}
\{ \mbqi{t}_{2i-1}, \mbpi{t}_{2r} | \bbi{t}_{2i-1},
 \bai{t}_{2r} \}_{i=(t/2) +1, r=(t/2)+1}^{l', k'}
\end{equation}
\end{subequations}

{ \bf Case 2: $t$ is odd }

We work similarly if $t$ is odd. 
Note that, 
 after we drop the  $p$-group $\mbpi{t-1}_{t-1}$ (see \eqref{ma.pt}),
we reach the system    
\begin{subequations}\mylabel{ma.e6'}
\begin{equation}\mylabel{ma.e6aa}
1 \unlhd \mbqi{t-1}_t  \unlhd \mbgi{t-1}_{t+1} \unlhd 
\dots \unlhd  \mbgi{t-1}_n = \mbgi{t-1},
\end{equation}
\begin{equation}\mylabel{ma.e6bb}
1, \bbi{t-1}_t \in \Irr(\mbqi{t-1}_{t}),   \Ti{t-1}_i \in \Irr(\mbgi{t-1}_i)
\end{equation}
\begin{equation}
\mylabel{ma.e6cc}
\{ \mbqi{t-1}_{2i-1}, \mbpi{t-1}_{2r} | \bbi{t-1}_{2i-1},
 \bai{t-1}_{2r} \}_{i=(t+1)/2, r=(t+1)/2}^{l', k'}
\end{equation}
\end{subequations}
that is equaivalent to \eqref{ma.e6} for the even case.

Observe that we can adjoint the trivial group at the bottom of 
\eqref{ma.e6aa} and consider the series 
\begin{equation}\mylabel{ad.1}
1 \unlhd 1 \unlhd  \mbqi{t-1}_t  \unlhd \mbgi{t-1}_{t+1} \unlhd 
\dots \unlhd  \mbgi{t-1}_n = \mbgi{t-1},
\end{equation}
This way $t$ has become even. We can now interchange the roles of $p$ and $q$,
and apply the already proved results of Case 1. 
For clarity  we remark  that, to prove the inductive step in this case 
 of an odd $t$, 
 we take 
 an $\mai{t-1}(\bbi{t-1}_t)$-invariant faithful linear limit
$$
(\mbgi{t}, \mbsi{t}, \bzi{t}, \mbqi{t}_t, \bbi{t}_t)
$$ 
 of $(\mbgi{t-1}, \mbsi{t-1}, \bzi{t-1}, \mbqi{t-1}_t, \bbi{t-1}_t)$.
So the system \eqref{ma.e6'} reduces to 
\begin{subequations}\mylabel{ma.e8'}
\begin{equation}\mylabel{ma.e8aa}
1 \unlhd \mbqi{t}_t \unlhd \mbgi{t}_{t+1} \unlhd \dots \unlhd
 \mbgi{t}_n = \mbgi{t}. 
\end{equation}
\begin{equation}
\mylabel{ma.e8bb}
1, \bbi{t}_t,   \Ti{t}_i \in \Irr(\mbgi{t}_i)
\end{equation}
\begin{equation}
\mylabel{ma.e8cc}
\{ \mbqi{t}_{2i-1}, \mbpi{t}_{2r} | \bbi{t}_{2i-1},
 \bai{t}_{2r} \}_{i=(t+1)/2, r=(t+1)/2}^{l', k'}
\end{equation}
\end{subequations}
All the conclusions of the  even case are  transfered to the odd case. 
In particular we get,  
\begin{remark}\mylabel{ma.r33}
The group  $\mbsi{t}= Z(\mbqi{t}_t)$ is a cyclic 
central  subgroup of $\mbgi{t}$, maximal among the abelian 
$\mbgi{t}$-invariant subgroups of $\mbgi{t}_t$,
while the character $\bbi{t}_t$ is 
$\mbgi{t}$-invariant.
\end{remark}
 
\begin{remark}\mylabel{ma.rrr33} 
Any faithful linear limit of the quintuple 
 $(\mbgi{t}, \mbri{t}, \bei{t}, \mbgi{t}_i, \Ti{t}_i)$ or
 the quintuple 
$(\mbgi{t}, \mbsi{t}, \bzi{t}, \mbgi{t}_i, \Ti{t}_i)$ 
is also a faithful linear limit of 
$(G, 1, 1, G_i, \chi_i)$ for all $i=t,t+1,  \dots, n$.
\end{remark}

\begin{proposition}\mylabel{ma.p44}
We have two central subgroups of $\mbgi{t}$, the $p$-group 
$\mbri{t}$ and the $q$-group $\mbsi{t}$. 
Along with them we get their irreducible characters 
$\bei{t} \in \Irr(\mbri{t})$ and $\bzi{t} \in \Irr(\mbsi{t})$. 
These  groups and characters satisfy
\begin{align*}
\mbri{t} \times \mbsi{t} \unlhd \mbpi{t}_{t+1} \cdot
 \mbqi{t}_t = \mbgi{t}_{t+1}\\
\bei{t} \text{ and } 
 \bzi{t} \text{ lie under $\bai{t}_{t+1}$ and $\bbi{t}_t$, 
respectively.  }
\end{align*}
\end{proposition}
 
\begin{proposition}\mylabel{ma.p55}
Every irreducible character $\Ti{t} \in \Irr(\mbgi{t})$ that lies above 
$\bei{t} \times \bzi{t} \in \Irr(\mbri{t} \times \mbsi{t})$, is  
monomial.
\end{proposition}

The main theorem for the inductive step in the case of an odd $t$,  will be
\begin{theorem}\mylabel{ma.main'}
After the above reductions, the group $\mbgi{t}_{t+1}$ 
is nilpotent if it exists, i.e., if $n \ge t+1$.
\end{theorem}

 Observe that,
as before, the fact $\mbgi{t}_{t+1}$ is nilpotent  
 allows us to shift the series \eqref{ma.e8aa} 
 by one. So, provided that  $t+2 \leq n$,  we get the series 
\begin{subequations}\mylabel{ma.e9'}
\begin{equation}\mylabel{ma.e9aa}
1 \unlhd \mbpi{t}_{t+1} \unlhd \mbgi{t}_{t+2} \unlhd \dots \unlhd
 \mbgi{t}_n = \mbgi{t}.
\end{equation}
The character tower and its  corresponding triangular set are carried
 over to the shifted series, as in  the Sylow system.
Hence the characters 
\begin{equation}
1, \bai{t}_{t+1}, \Ti{t}_{t+2}, \dots , \Ti{t}_n 
\end{equation}
form a tower for the series \eqref{ma.e9aa}.
 Its  corresponding triangular set is  
\begin{equation}
\mylabel{ma.e9cc}
\{ \mbqi{t}_{2i-1}, \mbpi{t}_{2r} | \bbi{t}_{2i-1},
 \bai{t}_{2r} \}_{i=(t+3)/2, r=(t+1)/2}^{l', k'}
\end{equation}
\end{subequations}

Now we can repeat the process  at the $t+1$-st step.

\section{Conclusions for the  smaller systems } \mylabel{ma.seI}

Before giving   the proofs of the three theorems stated in the previous 
section, we will analyze the behavior, under the described 
reductions,  of the smaller system 
\begin{subequations}\mylabel{ma.sm}
\begin{gather}
1\unlhd G_1 \unlhd G_2 \unlhd \dots \unlhd  G_m \unlhd G,\\
\{ \chi_i \in \Irr(G_i) \}_{i=0}^m\\
\{ Q_{2i-1}, P_{2r}, |\beta_{2i-1}, \alpha_{2r}\}_{i=1, r=0}^{l, k}
\end{gather}
\end{subequations}
for any fixed, but arbitrary, $m=1, \dots, n$.
The integers  $k$ and $l$ are, as usual, related to $m$ via \eqref{kl:def}.
Associated to this system are the groups 
 $ T=P_{2k}^* \rtimes I$ and $ U=Q_{2l-1}^* \rtimes J$,
where $I$ is the image of $\qw(\beta_{2k-1,2k})$ in $\Aut(P_{2k}^*)$ 
and $J$ that of $\hap(\alpha_{2l-2,2l-1})$ in $\Aut(Q_{2l-1}^*)$,
(for their definitions see \eqref{T} and  \eqref{U}).

After the first $\ma(\beta_1)$-invariant reductions the  above smaller system
 reduces, along with the original \eqref{ma.e1}, to the system
\begin{subequations}\mylabel{ma.sm1}
\begin{gather}
1\unlhd \mbgi{1}_{1}
 \unlhd \mbgi{1}_{2} \unlhd \dots \unlhd  \mbgi{1}_{m} \unlhd \mbgi{1},\\
\{ \Ti{1}_i \in \Irr(\mbgi{1}_i) \}_{i=0}^m\\
\{ \mbqi{1}_{2i-1}, \mbpi{1}_{2r}, |\bbi{1}_{2i-1},
 \bai{1}_{2r}\}_{i=1, r=0}^{l, k}
\end{gather}
\end{subequations}
(Note that this last triangular set is a subset of \eqref{ma.e4}.)
Of course the groups $T$ and $U$ reduce 
 to $\mbti{1}$ and $\mbui{1}$, respectively. 
So $\mbti{1}= \mbpsi{1}_{2k} \rtimes \mbii{1}$, where 
$\mbii{1}$ is the image of $\mbqwi{1}(\bbi{1}_{2k-1,2k})$ 
in $\Aut(\mbpsi{1}_{2k})$. Similarly, 
$\mbui{1}= \mbqsi{1}_{2l-1} \rtimes \mbji{1}$, where $\mbji{1}$
is the image of $\mbhapi{1}(\bai{1}_{2l-2,2l-1})$ in 
$\Aut(\mbqsi{1}_{2l-1})$. 

Then, according to  Theorem \ref{lmt.t3},
\begin{theorem}\mylabel{f.t1}
The  groups $\mbqwi{1}(\bbi{1}_{2k-1,2k})$ and 
$\mbhapi{1}(\bai{1}_{2l-1,2l-1})$   can be chosen to satisfy  
\begin{equation*}
P_{2k}^* \text{ is naturally isomorphic to } \mbpsi{1}_{2k}.
\end{equation*}
This  isosmorphism sends the image of $\qw(\beta_{2k-1,2k})$ 
in $\Aut(P_{2k}^*)$
onto that of $\mbqwi{1}(\bbi{1}_{2k-1,2k})$ in $\Aut(\mbpsi{1}_{2k})$. 
So $I \cong \mbii{1}$. In addition, 
\begin{gather*}
T \cong \mbti{1}, \\
\text{ Image of $\hap(\alpha_{2l-2,2l-1})$  
in $\Aut(\mbqsi{1}_{2l-1})$ } = \mbji{1}= 
\text{ Image of $\mbhapi{1}(\bai{1}_{2l-2,2l-1})$  
in $\Aut(\mbqsi{1}_{2l-1})$ }.
\end{gather*}
\end{theorem}
Note that as $P_{2k}^*$ is naturally isomorphic to $ \mbpsi{1}_{2k}$,  we
get that
\begin{equation}\mylabel{f.eq0}
\text{ Image of $\qw(\beta_{2k-1,2k}) $  in $\Aut(\mbpsi{1}_{2k})$ } =
  \mbii{1}=
\text{ Image of $\mbqwi{1}(\bbi{1}_{2k-1,2k}) $  in $\Aut(\mbpsi{1}_{2k})$ }
\end{equation}

As we have already seen, the group 
$R=1$ seen inside $\mbpi{1}_2 \leq \mbpsi{1}_{2k}$ is denoted by $\mbri{1}$, 
and its irreducible character $\eta=1$, by $\bei{1}$. So we can form the
 linear quintuple $(\mbti{1}, \mbri{1}, \bei{1},
 \mbpsi{1}_{2k}, \basi{1}_{2k})$. Furthermore,  $\mbti{1} \cong T$, while   
 $\mbpsi{1}_{2k}$ is naturally isomorphic to 
 $P_{2k}^*$, (so that the character $\basi{1}_{2k}$ maps to 
$\alpha_{2k}^*$).  
This, along with the equations in \eqref{lm.eG},  implies
\begin{corollary}\mylabel{f.co1}
After the first $\ma(\beta_1)$-invariant  reductions we have 
\begin{equation}\mylabel{f.q1}
(\mbti{1}, \mbri{1}, \bei{1}, \mbpsi{1}_{2k}, \basi{1}_{2k}) \cong 
(T, R , \eta, P_{2k}^*, \alpha_{2k}^*).
\end{equation}
Furthermore, the  quintuple 
 $(\mbui{1}, \mbsi{1}, \bzi{1}, \mbqsi{1}_{2l-1}, \bbsi{1}_{2l-1})$ 
is a $G$-associate faithful linear limit  of  
 $(U, S, \zeta, Q_{2l-1}^*,
 \beta_{2l-1}^*)$. Hence any faithful linear limit of the former quintuple 
is also a faithful linear limit of the latter one. 
\end{corollary}
From now on, we will identify the quintuples in \eqref{f.q1}.

Note that when we shift the series \eqref{ma.e3} by one in \eqref{ma.sh1},
 the same happens to 
 the series (\ref{ma.sm}a). Thus we get the  system 
\begin{subequations}\mylabel{ma.sm2}
\begin{gather}
1\unlhd \mbpi{1}_{2}
 \unlhd \mbgi{1}_{3} \unlhd \dots \unlhd  \mbgi{1}_{m} \unlhd \mbgi{1},\\
1, \bai{1}_2,  \Ti{1}_3, \dots,   \Ti{1}_m, \\
\{ \mbqi{1}_{2i-1}, \mbpi{1}_{2r}, |\bbi{1}_{2i-1},
 \bai{1}_{2r}\}_{i=2, r=1}^{l, k}
\end{gather}
\end{subequations}
that is a smaller system for \eqref{ma.sh1}. 
Observe also that, according to Theorem \ref{sh.t1}, 
 none of the groups $\mbii{1}, \mbji{1}, \mbti{1}$ and $\mbui{1}$ changes
passing to the shifted case. 

When the next series of $\mbi{1}(\bai{1}_2)$-invariant reductions is 
performed 
the system \eqref{ma.sm2} reduces to 
\begin{subequations}\mylabel{ma.sm3}
\begin{gather}
1\unlhd \mbpi{2}_{2}
 \unlhd \mbgi{2}_{3} \unlhd \dots \unlhd  \mbgi{2}_{m} \unlhd \mbgi{2},\\
1, \bai{2}_2,  \Ti{2}_3, \dots,   \Ti{2}_m, \\
\{ \mbqi{2}_{2i-1}, \mbpi{2}_{2r}, |\bbi{2}_{2i-1},
 \bai{2}_{2r}\}_{i=2, r=1}^{l, k}
\end{gather}
\end{subequations}
The groups $\mbii{1}, \mbji{1}$ and $\mbti{1}, \mbui{1}$ reduce to the groups 
$\mbii{2}, \mbji{2}$ and $\mbti{2}, \mbui{2}$, respectively, so that Theorem 
\ref{lmt.t7} holds.
(The notation is as expected, that is 
 $\mbti{2}= \mbpsi{2}_{2k} \rtimes \mbii{2}$, where 
$\mbii{2}$ is the image of $\mbqwi{2}(\bbi{2}_{2k-1,2k})$ 
in $\Aut(\mbpsi{2}_{2k})$, and  
$\mbui{2}= \mbqsi{2}_{2l-1} \rtimes \mbji{2}$, where $\mbji{2}$
is the image of $\mbhapi{2}(\bai{2}_{2l-2,2l-1})$ in 
$\Aut(\mbqsi{2}_{2l-1})$.)  In view of   
Theorem \ref{lmt.t7}, we get 
\begin{theorem}\mylabel{f.t2}
The reduced  system 
\eqref{ma.sm3} satisfies
\begin{equation*}
 \mbqsi{1}_{2l-1} \text{ is naturally isomorphic to } \mbqsi{2}_{2l-1},
\end{equation*}
This  isosmorphism sends the image of
 $\mbhapi{1}(\bai{1}_{2l-2,2l-1})$
in $\Aut(\mbqsi{1}_{2l-1})$ 
onto that of $\mbhapi{2}(\bai{2}_{2l-2,2l-1})$ in 
$\Aut(\mbqsi{2}_{2l-1})$. So $\mbji{1} \cong \mbji{2}$. In addition, 
\begin{gather*}
\mbui{1}\cong \mbui{2}, \\
\text{ Image of $\mbqwi{1}(\bbi{1}_{2k-1,2k})$  
in $\Aut(\mbpsi{2}_{2k})$ } =\mbii{2}=
 \text{ Image of $\mbqwi{2}(\bbi{2}_{2k-1,2k})$  
in $\Aut(\mbpsi{2}_{2k})$ }.
\end{gather*}
\end{theorem}
Note that similar to \eqref{f.eq0}, the fact that $\mbqsi{1}_{2l-1}$ is
 naturally  isomorphic to $\mbqsi{2}_{2l-1}$, along with Theorem \ref{f.t2}, 
 implies 
  \begin{equation}\mylabel{f.eq00}
\text{ Image of $\mbhapi{1}(\bai{1}_{2l-2,2l-1}) $ 
 in $\Aut(\mbqsi{2}_{2l-1})$ } =
  \mbji{2}=
\text{ Image of $\mbhapi{2}(\bai{1}_{2l-2,2l-1}) $  
in $\Aut(\mbqsi{2}_{2l-1})$. }
\end{equation}
Furthermore, the group $\mbsi{1}$ is isomorphic to 
 $\mbsi{2}$, (see \eqref{ma.rs1b}), 
while its irreducible character $\bzi{1}\in \Irr(\mbsi{1})$ 
maps under the above isomorphism  to $\bzi{2} \in \Irr(\mbsi{2})$.
This, along with the isomorphism between $\mbui{1}$ and $\mbui{2}$, 
and the remarks following Theorem \ref{lmt.t7} (see \eqref{lmt.ef}), implies 
\begin{corollary}\mylabel{f.co2}
After the second  $\mbi{1}(\bai{1}_{2})$-invariant  reductions we have 
\begin{equation}\mylabel{f.q2}
(\mbui{1}, \mbsi{1}, \bzi{1}, \mbqsi{1}_{2l-1}, \bbsi{1}_{2l-1}) \cong 
(\mbui{2}, \mbsi{2}, \bzi{2}, \mbqsi{2}_{2l-1}, \bbsi{2}_{2l-1}).
\end{equation}
Furthermore, the quintuple
 $(\mbti{2}, \mbri{2}, \bei{2}, \mbpsi{2}_{2k}, \basi{2}_{2k})$ 
is a $\mbgi{1}$-associate faithful linear limit  of
 $(\mbti{1}, \mbri{1}, \bei{1}, \mbpsi{1}_{2k}, \basi{1}_{2k})$. 
Hence any faithful linear limit of the former quintuple is also a faithful 
linear limit of the latter one. 
\end{corollary}
From now on we will identify the quintuples in \eqref{f.q2}.

Corollary \ref{f.co1} and \ref{f.co2} combined 
(with the appropriate identifications) imply
\begin{theorem}\mylabel{f.t3}
Any faithful linear limit of the quintuple $(\mbti{2}, \mbri{2}, 
\bei{2}, \mbpsi{2}_{2k},  \basi{2}_{2k})$, or the  $(\mbui{2}, \mbsi{2}, 
\bzi{2}, \mbqsi{2}_{2l-1},  
\bbsi{2}_{2l-1})$ is also a faithful linear limit of the quintuple
$(T, R, \eta, P_{2k}^*, \alpha_{2k}^*)$, or 
$(U, S, \zeta, Q_{2l-1}^*, \beta_{2l-1}^*)$, respectively.
\end{theorem}

Furthermore, for the image $\mbji{2}$ of $\mbhapi{2}(\bai{2}_{2l-2,2l-1})$ 
in $\Aut(\mbqsi{2}_{2l-1})$ we have 
\begin{align*}
\mbji{2}
&=\text{ Image of $\mbhapi{1}(\bai{1}_{2l-2,2l-1})$ in 
$\Aut(\mbqsi{2}_{2l-1})$ }
 &\text{by \eqref{f.eq00}} \\
&= \text{ Image of $\hap(\alpha_{2l-2,2l-1})$ in $\Aut(\mbqsi{2}_{2l-1})$ }
&\text{ by Theorem \ref{f.t1}, since  $\mbqsi{1}_{2l-1} \cong \mbqsi{2}_{2l-1}$}
\end{align*}
Even more, equation \eqref{f.eq0}, along with the fact 
that    $\mbpsi{2}_{2k}$ is a factor of $\mbpsi{1}_{2k}$, implies   
$$
\text{ Image of $\qw(\beta_{2k-1,2k}) $  in $\Aut(\mbpsi{2}_{2k})$ } =
\text{ Image of $\mbqwi{1}(\bbi{1}_{2k-1,2k})$  in $\Aut(\mbpsi{2}_{2k})$. }
$$
Therefore, for the image $\mbii{2}$ of $\mbqwi{2}(\bbi{2}_{2k-1,2k})$
in $\Aut(\mbpsi{2}_{2k})$ we have 
\begin{align*}
\mbii{2}&= \text{ Image of $\mbqwi{1}(\bbi{1}_{2k-1,2k})$  in 
$\Aut(\mbpsi{2}_{2k})$ } &\text{by  Theorem \ref{f.t2} }  \\
&\text{= Image of $\qw(\beta_{2k-1,2k})$ in $\Aut(\mbpsi{2}_{2k})$ .}
\end{align*}
In conclusion, we get 
\begin{theorem}\mylabel{f.t4}
$$
\mbji{2} =  \text{ Image of $\hap(\alpha_{2l-2,2l-1})$
 in $\Aut(\mbqsi{2}_{2l-1})$ },
$$
$$ 
\mbii{2} =  \text{ Image of $\qw(\beta_{2k-1,1k})$ in
 $\Aut(\mbpsi{2}_{2k})$ }.
$$ 
\end{theorem}

When we shift the series \eqref{ma.eE1} by one  to  get (\ref{ma.sh2}a), 
the system \eqref{ma.sm3} is also shifted  by one. Thus we reach 
\begin{subequations}\mylabel{ma.sm4}
\begin{gather}
1\unlhd \mbqi{2}_{3}
 \unlhd \mbgi{2}_{4} \unlhd \dots \unlhd  \mbgi{2}_{m} \unlhd \mbgi{2},\\
1, \bbi{2}_3,  \Ti{2}_4, \dots, \Ti{2}_m, \\
\{ \mbqi{2}_{2i-1}, \mbpi{2}_{2r}, |\bbi{2}_{2i-1},
 \bai{2}_{2r}\}_{i=2, r=2}^{l, k}
\end{gather}
\end{subequations}
According to  Theorem \ref{sh.t1},
 the groups $\mbii{2}, \mbji{2}, \mbti{2}$ and
$\mbui{2}$ remain unchanged passing to the shifted  system \eqref{ma.sm4}.

According to our inductive hypothesis, after $t-1$ steps 
(where $2 < t \leq m$),
we reach the system 
\begin{subequations}\mylabel{ma.smt}
\begin{gather}
1\unlhd \mbpi{t-1}_{t} \unlhd \mbgi{t-1}_{t+1}
 \unlhd \dots \unlhd  \mbgi{t-1}_{m} 
\unlhd \mbgi{t-1},\\
1, \bai{t-1}_{t},  \Ti{t-1}_{t+1}, \dots,  \Ti{t-1}_m, \\
\{ \mbqi{t-1}_{2i-1}, \mbpi{t-1}_{2r}, |\bbi{t-1}_{2i-1},
 \bai{t-1}_{2r}\}_{i=(t/2) +1, r=t/2}^{l, k}
\end{gather}
\end{subequations}
when $t$ is even. In the case of an odd $t$  we get   
\begin{subequations}\mylabel{ma.smtt}
\begin{gather}
1\unlhd \mbqi{t-1}_{t} \unlhd \mbgi{t-1}_{t+1} 
\unlhd \dots \unlhd  \mbgi{t-1}_{m} 
\unlhd \mbgi{t-1},\\
1, \bbi{t-1}_{t},  \Ti{t-1}_{t+1}, \dots,  \Ti{t-1}_m, \\
\{ \mbqi{t-1}_{2i-1}, \mbpi{t-1}_{2r}, |\bbi{t-1}_{2i-1},
 \bai{t-1}_{2r}\}_{i=(t+1)/2 +1, r=(t+1)/2}^{l, k}
\end{gather}
\end{subequations}
Furthermore, the groups $I , J, T$ and $U$ are reduced to 
the groups  $\mbii{t-1}, \mbji{t-1}, \mbti{t-1}$ and $\mbui{t-1}$, where 
 $\mbti{t-1}= \mbpsi{t-1}_{2k} \rtimes \mbii{t-1}$,  and
$\mbii{t-1}$ is the image of $\mbqwi{t-1}(\bbi{t-1}_{2k-1,2k})$ 
in $\Aut(\mbpsi{t-1}_{2k})$. Similarly,
$\mbui{t-1}= \mbqsi{t-1}_{2l-1} \rtimes \mbji{t-1}$ where $\mbji{t-1}$
is the image of $\mbhapi{2t-1}(\bai{t-1}_{2l-2,2l-1})$ in 
$\Aut(\mbqsi{t-1}_{2l-1})$.
These groups satisfy
\begin{theorem}\mylabel{f.t5}
Any faithful linear limit of  $(\mbti{t-1}, \mbri{t-1}, 
\bei{t-1}, \mbpsi{t-1}_{2k},  \basi{t-1}_{2k})$
 or  the quintuple $(\mbui{t-1}, \mbsi{t-1}, \bzi{t-1}, \mbqsi{t-1}_{2l-1},  
\bbsi{t-1}_{2l-1})$ is also a faithful linear limit of 
$(T, R, \eta, P_{2k}^*, \alpha_{2k}^*)$ or 
$(U, S, \zeta, Q_{2l-1}^*, \beta_{2l-1}^*)$, respectively.
\end{theorem}
Furthermore for the groups $\mbii{t-1}$ and $\mbji{t-1}$ we have 
\begin{theorem}\mylabel{f.t6}
$$ 
\mbji{t-1} =  \text{ Image of $\hap(\alpha_{2l-2,2l-1})$
 in $\Aut(\mbqsi{t-1}_{2l-1})$ },
$$
$$ 
\mbii{t-1} =  \text{ Image of $\qw(\beta_{2k-1,1k})$ in
 $\Aut(\mbpsi{t-1}_{2k})$ }.
$$
\end{theorem}

At the inductive step (the $t$-th step), either we perform a series of 
$\mbi{t-1}(\bai{t-1}_t)$-invariant reductions,if $t$ is even, or 
 a series of $\mai{t-1}(\bbi{t-1}_t)$-invariant reductions, if $t$ is odd.
The systems \eqref{ma.smt} and \eqref{ma.smtt} reduce to 
\begin{subequations}\mylabel{f.smt}
\begin{gather}
1\unlhd \mbpi{t}_{t} \unlhd \mbgi{t}_{t+1}
 \unlhd \dots \unlhd  \mbgi{t}_{m} 
\unlhd \mbgi{t},\\
1, \bai{t}_{t},  \Ti{t}_{t+1}, \dots,  \Ti{t}_m, \\
\{ \mbqi{t}_{2i-1}, \mbpi{t}_{2r}, |\bbi{t}_{2i-1},
 \bai{t}_{2r}\}_{i=(t/2) +1, r=t/2}^{l, k}
\end{gather}
\end{subequations}
when $t$ is even, and   
\begin{subequations}\mylabel{f.smtt}
\begin{gather}
1\unlhd \mbqi{t}_{t} \unlhd \mbgi{t}_{t+1} 
\unlhd \dots \unlhd  \mbgi{t}_{m} 
\unlhd \mbgi{t},\\
1, \bbi{t}_{t},  \Ti{t}_{t+1}, \dots,  \Ti{t}_m, \\
\{ \mbqi{t}_{2i-1}, \mbpi{t}_{2r}, |\bbi{t}_{2i-1},
 \bai{t}_{2r}\}_{i=(t+1)/2 +1, r=(t+1)/2}^{l, k}
\end{gather}
\end{subequations}
when $t$ is odd.

{\bf Case 1: $t$ is even } 

Suppose first that $t$ is even.  Then the groups $\mbii{t-1}, 
\mbji{t-1}, \mbti{t-1}$ and $\mbui{t-1}$ (for the system 
\eqref{ma.smt}), reduce to the groups
$\mbii{t}, \mbji{t}, \mbti{t}$ and $\mbui{t}$, 
(for the system \eqref{f.smt}), 
so that Theorem \ref{lmt.t7} holds. (Observe we can use Theorem \ref{lmt.t7}
as the reductions we used to get \eqref{f.smt} were 
  $\mbi{t-1}(\bai{t-1}_t)$-invariant.)  In particular we get  
\begin{subequations}\mylabel{f.e6}
\begin{gather}
 \mbqsi{t-1}_{2l-1} \text{ is naturally isomorphic to } \mbqsi{t}_{2l-1},\\
\text{ this  isosmorphism sends the image of
 $\mbhapi{t-1}(\bai{t-1}_{2l-2,2l-1})$ } \notag  \\
\text{ in $\Aut(\mbqsi{t-1}_{2l-1})$
onto that of $\mbhapi{t}(\bai{t}_{2l-2,2l-1})$ in 
$\Aut(\mbqsi{t}_{2l-1})$,i.e.,  }
\mbji{t-1} \cong \mbji{t}, \\
\mbui{t-1} \cong \mbui{t},  \\
\text{ Image of $\mbqwi{t-1}(\bbi{t-1}_{2k-1,2k})$ in
 $\Aut(\mbpsi{t}_{2k})$} =
\text{ Image of $\mbqwi{t}(\bbi{t}_{2k-1,2k})$ in
 $\Aut(\mbpsi{t}_{2k})$}. 
\end{gather}
\end{subequations}
Similar to \eqref{f.eq00}, the above implies
  \begin{equation}\mylabel{f.eq000}
\text{ Image of $\mbhapi{t-1}(\bai{t-1}_{2l-2,2l-1}) $ 
 in $\Aut(\mbqsi{t}_{2l-1})$ } =
  \mbji{t}=
\text{ Image of $\mbhapi{t}(\bai{t}_{2l-2,2l-1}) $  
in $\Aut(\mbqsi{t}_{2l-1})$. }
\end{equation}
This, along with the fact that $\mbsi{t} \cong \mbsi{t-1}$
 (see \eqref{ma.e20b}),
implies that  
\begin{equation}\mylabel{f.e2}
(\mbui{t-1}, \mbsi{t-1}, \bzi{t-1},\mbqsi{t-1}_{2l-1}, \bbsi{t-1}_{2l-1})
\cong (\mbui{t}, \mbsi{t}, \bzi{t},\mbqsi{t}_{2l-1}, 
\bbsi{t}_{2l-1}).
\end{equation}
From now on we identify these two quintuples.
Observe  that 
 $(\mbti{t}, \mbri{t}, \bei{t}, \mbpsi{t}_{2k}, \basi{t}_{2k})$
is a $\mbgi{t-1}$-associate faithful linear limit of 
 $(\mbti{t-1}, \mbri{t-1}, \bei{t-1}, \mbpsi{t-1}_{2k}, \basi{t-1}_{2k})$, 
(see \eqref{lmt.ef}).
Hence, Proposition \ref{lmt.pr3} implies that any faithful linear limit of
 $(\mbti{t}, \mbri{t}, \bei{t}, \mbpsi{t}_{2k}, \basi{t}_{2k})$ is also a 
faithful linear limit of
  $(\mbti{t-1}, \mbri{t-1}, \bei{t-1}, \mbpsi{t-1}_{2k}, \basi{t-1}_{2k})$.
This, along with \eqref{f.e2} and Theorem \ref{f.t3} implies  
\begin{theorem}\mylabel{f.t7}
Any faithful linear limit of the quintuple  $(\mbti{t}, \mbri{t}, 
\bei{t}, \mbpsi{t}_{2k},  \basi{t}_{2k})$
 or the quintuple  $(\mbui{t}, \mbsi{t}, \bzi{t}, \mbqsi{t}_{2l-1},  
\bbsi{t}_{2l-1})$ is also a faithful linear limit of the quintuple 
$(T, R, \eta, P_{2k}^*, \alpha_{2k}^*)$ or 
$(U, S, \zeta, Q_{2l-1}^*, \beta_{2l-1}^*)$, respectively.
\end{theorem}
So Theorem \ref{f.t5} holds for the inductive $t$-th step, in the
case of an even $t$. In addition, 
 Theorem \ref{f.t6} is still valid  for an even $t$, i.e., 
 \begin{theorem}\mylabel{f.t8}
\begin{subequations}\mylabel{f.e3}
\begin{equation}\mylabel{f.e3a}
\mbji{t} =  \text{ Image of $\hap(\alpha_{2l-2,2l-1})$
 in $\Aut(\mbqsi{t}_{2l-1})$ },
\end{equation}
\begin{equation}\mylabel{f.e3b}
\mbii{t} = \text{ Image of $\qw(\beta_{2k-1,2k})$ in
 $\Aut(\mbpsi{t}_{2k})$ }. 
\end{equation}
\end{subequations}
\end{theorem}

\begin{proof}
For the proof of \eqref{f.e3a},  note
that Theorem \ref{f.t6}, along with the fact that the section $\mbqsi{t}_{2l-1}$ of  $ \mbqsi{t-1}_{2l-1}$ is isomorphic to 
  $\mbqsi{t-1}_{2l-1}$(see \eqref{ma.e200}), implies that 
$$
\text{ Image  of $\mbhapi{t-1}(\bai{t-1}_{2l-2,2l-1})$ in
 $\Aut(\mbqsi{t}_{2l-1})$ }= \text{ Image of $\hap(\alpha_{2l-2, 2l-1})$ in
 $\Aut(\mbqsi{t}_{2l-1})$ }.
$$
This, along with \eqref{f.eq000}, implies 
\begin{align*}
\mbji{t} &=   \text{ Image  of $\mbhapi{t-1}(\bai{t-1}_{2l-2,2l-1})$ in
 $\Aut(\mbqsi{t}_{2l-1})$ }  \\
&= \text{ Image  of $\hap(\alpha_{2l-2,2l-1})$ in
 $\Aut(\mbqsi{t}_{2l-1})$ }.
\end{align*}
Hence \eqref{f.e3a} follows.

As far as \eqref{f.e3b} is concerned, first observe that 
\begin{align}\mylabel{f.e4}
\mbii{t} &=  \text{ Image of $\mbqwi{t}(\bbi{t}_{2k-1, 2k})$
 in $\Aut(\mbpsi{t}_{2k})$ } & \notag \\
&= \text{ Image of $\mbqwi{t-1}(\bbi{t-1}_{2k-1,2k})$
 in $\Aut(\mbpsi{t}_{2k})$ } &\text{ by (\ref{f.e6}d).} 
\end{align}
On the other hand, Theorem \ref{f.t6} implies that 
\begin{equation}\mylabel{f.e5}
\mbii{t-1}=  \text{ Image of $\mbqwi{t-1}(\bbi{t-1}_{2k-1,2k})$
 in $\Aut(\mbpsi{t-1}_{2k})$ }=\text{ Image of $\qw(\beta_{2k-1,2k})$
 in $\Aut(\mbpsi{t-1}_{2k})$ }.
\end{equation}
But $\mbpsi{t}$ is a section of $\mbpsi{t-1}$. Hence \eqref{f.e5} implies that 
$$
 \text{ Image of $\mbqwi{t-1}(\bbi{t-1}_{2k-1,2k})$
 in $\Aut(\mbpsi{t}_{2k})$ }=\text{ Image of $\qw(\beta_{2k-1,2k})$
 in $\Aut(\mbpsi{t}_{2k})$ }.
$$
This, along with \eqref{f.e4}, implies the desired 
equation in \eqref{f.e3a}.
\end{proof}

{ \bf Case 2: $t$ is odd }
We can work similarly in the case of an odd $t$, to prove that Theorems
 \ref{f.t7} and \ref{f.t8} are still valid.
This time we use Theorem \ref{lmt.t3} on the system \ref{f.smtt}, as 
the reductions we used to get  \ref{f.smtt} were
 $\mai{t-1}(\bbi{t-1}_t)$-invariant. Note that the analogue
 to \eqref{f.e6}
in this case is 
\begin{subequations}\mylabel{f.e1}
\begin{gather}
 \mbpsi{t-1}_{2k} \text{ is naturally isomorphic to } \mbpsi{t}_{2k},\\
\text{ this  isosmorphism sends the image of
 $\mbqwi{t-1}(\bbi{t-1}_{2k-1,2k})$ } \notag  \\
\text{ in $\Aut(\mbpsi{t-1}_{2k})$
onto that of $\mbqwi{t}(\bbi{t}_{2k-1,1k})$ in 
$\Aut(\mbpsi{t}_{2k})$,i.e., }
\mbii{t-1} \cong \mbii{t} \\ 
\mbti{t-1} \cong \mbti{t}, \\
\text{ Image of $\mbhapi{t-1}(\bai{t-1}_{2l-2,2l-1})$ in
 $\Aut(\mbqsi{t}_{2l-1})$} =
\text{ Image of $\mbhapi{t}(\bai{t}_{2l-2,2l-1})$ in
 $\Aut(\mbqsi{t}_{2l-1})$}. 
\end{gather}
\end{subequations}
As with the even  case (note that this time  $\mbri{t} \cong \mbri{t-1}$), 
 we have 
\begin{equation}
(\mbti{t-1}, \mbri{t-1}, \bei{t-1}, \mbpsi{t-1}_{2k}, \basi{t-1}_{2k})\cong
(\mbti{t}, \mbri{t}, \bei{t}, \mbpsi{t}_{2k}, \basi{t}_{2k}).
\end{equation}
Furthermore, Proposition \ref{lm.pr33} implies that 
any faithful linear limit of
 $(\mbui{t}, \mbsi{t}, \bzi{t}, \mbqsi{t}_{2k}, \bbsi{t}_{2k})$ is also a 
faithful linear limit of
  $(\mbui{t-1}, \mbsi{t-1}, \bzi{t-1}, \mbqsi{t-1}_{2k}, \bbsi{t-1}_{2k})$.
These are enough for the proof of Theorem \ref{f.t7}, when $t$ is odd.

The proof of Theorem \ref{f.t8} in the case of an odd $t$,
 is in the same spirit as the one we have already given for even $t$, 
so we omit it.

We can now give the proofs of 
Theorems \ref{ma.t1}, \ref{ma.t2} and \ref{ma.main}.

%%% Local Variables: 
%%% mode: latex
%%% TeX-master: t
%%% End: 

\section{The $t=1$ case}
\mylabel{ma.se2}
In this section we handle the $t=1$ case, i.e., we  give the proof of 
Theorem \ref{ma.t1}. 
For that  we need a result
stated as Theorem \ref{m.t2} below, which is interesting on its own.

We use the same notation about symplectic modules as the one introduced in 
Section \ref{lp}.
The ``magic'' theorem of  this chapter is Theorem 3.2  in \cite{da}.
An equivalent form of this theorem we give here as Theorem \ref{m.t1}. 
\begin{theorem}\mylabel{m.t1}
Suppose that $F$ is a finite field  of odd characteristic $p$,
 that $G$ is a finite $p$-solvable group, 
that $H$ is a subgroup of $p$-power index in $G$, 
 that $B$ is an anisotropic  symplectic  $FG$-module
 and that $C$ is an $FG$-submodule of $B$. Then the
$G$-invariant symplectic form on $B$ restricts to a $G$-invariant
symplectic form on $C$. If $C$, with this form, restricts to a hyperbolic
symplectic $FH$-module $C|_H$, then $C = 0$.
\end{theorem}
\begin{proof}
  Since $B$ is symplectic and $FG$-anisotropic, so is its
$FG$-submodule $C$. Theorem 3.2 of \cite{da}, applied to $C$, 
tells us that $C$ is $FG$-hyperbolic 
if $C|_H$ is $FH$-hyperbolic. In that case $C$ is both
$FG$-anisotropic and $FG$-hyperbolic. So it must be $0$.
\end{proof}

After these preliminaries we can prove the first important theorem of 
this chapter.

\begin{theorem}\mylabel{m.t2}
Assume that $G$ is a $p, q$-group, where $p$ and $q$ are distinct odd primes, 
 and that $N, M$ are  normal subgroups of $G$.
Let $M=P \times S$ and  $N = P \rtimes Q$, where $P$ is a $p$-group, 
and $S,  Q $ are $q$-groups with $S \leq Q$.
  Assume that the center  $Z(P)$ of $P$ is  maximal among the  abelian 
$G$-invariant subgroups of $P$.
 Let $\chi, \alpha,\beta $ and $  \zeta$  be irreducible characters of
 $G, P, S$ and 
$Z(P)$ respectively that satisfy 
\begin{subequations}
\begin{align}
\chi \in \Irr(G | \alpha \times \beta)
 &\text{ and } \alpha \in \Irr(P| \zeta),  \mylabel{m.e1}\\ 
 \zeta \text{ is a faithful } &G\text{-invariant character of  }Z(P),  \mylabel{m.e2} \\
G(\beta ) &= G \mylabel{m.e3}, \\
\chi \text{ is a monomial } &\text{character of $G$   with } \chi(1)_{q} = \beta(1),  \mylabel{m.e4}
\end{align}
\end{subequations}
where $\chi(1)_{q}$ denotes the $q$-part of the integer $\chi(1)$.
Then $Q$  centralizes $P$.
\end{theorem}

\begin{proof}
First observe that 
$S, P $ are   normal subgroups of $G$, as they are  characteristic subgroups
 of $M \unlhd G$.
The existence of the  $G$-invariant
  faithful  character $\zeta \in \Irr(Z(P))$ 
implies   that  $Z(P)$ is a cyclic central subgroup of $G$, i.e., 
$Z(P)  \leq Z(G)$.   
Furthermore, the fact that $Z(P)$ is maximal among the abelian
subgroups of $P$ that are  normal in $G$ implies that 
every characteristic abelian subgroup of $P$ is contained in $Z(P)$ and thus
 is cyclic.  Hence P. Hall's  theorem  
implies that either  $P$ is an abelian group or it 
is the central product 
\begin{subequations}\mylabel{m.e11}
\begin{equation}
P = T \odot Z(P), 
\end{equation}
where $T$ is an extra special $p$-group of exponent $p$ and 
 \begin{equation}
T \cap Z(P) = Z(T).
\end{equation}
Note that  the group $T$ is unique, as $T = \Omega(P)$.
\end{subequations}

In the case that $P = Z(P)$ is an abelian group, Theorem \ref{m.t2} holds trivially, 
as $P = Z(P) \leq Z(G)$ is centralized by $G$. 
Thus we may assume that $P > Z(P)$ and \eqref{m.e11} holds.

Since $\chi$ is monomial, there exists a subgroup $H$ of $G$, 
  and a linear character $\lambda \in \Lin(H)$ that 
induces $\chi = \lambda^G$. Then the product $HS$ forms  a subgroup of $G$. 
Furthermore,
\begin{claim}\mylabel{m.st1}
$|G: HS| $ is a $p$-number, and $(\lambda^{HS}) |_S = \beta$. 
\end{claim}
\begin{proof}
Indeed, 
as $S \unlhd G$, Clifford's theory,  along with \eqref{m.e3}, implies that
 $$
\chi |_S = m \cdot \beta, 
$$
for some integer $m \geq 0$.  Hence $\deg (\chi)= m \deg (\beta)$.
Since $\chi(1)_{q}= \beta(1)$, we have  that $m$ is a power of $p$.
As $ H \leq HS  \leq G$, the induced character $\lambda^{HS}$
lies in $\Irr(HS)$ and induces $(\lambda^{HS})^G = \lambda^G = \chi$.
 So 
$$
\deg(\lambda^{HS}) \cdot |G:HS| = \deg(\chi)= m \deg(\beta).
$$
Clifford's theorem also implies that $\lambda^{HS}|_S = s \beta$,
 for some integer $s$.
As $ \deg(\lambda^{HS}) = |HS :H|= |S : H \cap S|$ we get that
both  $\deg(\lambda^{HS}|_S)$ and  $s$ are  $q$-numbers.   
But 
$$
s \deg(\beta) |G :HS| = \deg(\lambda^{HS}) \cdot |G:HS| = m \deg(\beta),
$$
with $m$ a $p$-number.
Hence $s=1$, while  $|G:HS|$ is a power of  $p$. This completes the 
proof of the claim.
\end{proof}

The  fact that $\lambda \in \Lin(H)$ induces irreducibly to  $G$ implies that 
the center, $Z(G)$,  
of $G$ is a subgroup  of $H$. This, along with the fact that $Z(P ) \leq Z(G)$, implies
\begin{equation}\mylabel{m.e6}
Z(P) \leq Z(G) \leq H.
\end{equation}
Let $E := [P, Q]$. Then $E$ is a characteristic subgroup of
 $N$ and thus a normal subgroup of $G$. 
Even more, we have 
\begin{claim}\mylabel{n.st2}
$E = [P, Q]$ is an abelian group.
\end{claim}
\begin{proof}
 Suppose not. Then $E$ is a non-abelian normal subgroup of $G$ 
contained in  $P = T \cdot Z(P)$. 
 As $Z(P) \leq Z(G)$ (by \eqref{m.e6}),we have 
$E = [P, Q] = [T, Q] \leq T$. 
 Furthermore, $Z(E)$ is an abelian  normal subgroup of $G$, 
contained in $T \leq P$. Hence $Z(E)$ is contained in 
$T \cap Z(P) = Z(T)$. As $E$ is non-abelian and $Z(T)$ has order $p$, 
we conclude that  $Z(E)= Z(T) \leq Z(P) \leq Z(G)$. 
Therefore $E= [T, Q]$ is an extra special subgroup of $T$ of exponent $p$, 
whose center is central in $G$.
Hence the group $E$ satisfies condition (4.3a)
 in \cite{da}. In addition, $Q$ is a $p'$-subgroup of $G$ such that 
$QE$ is normal in $G$ (as  $P =[P,Q]C(Q \tin P)$ and thus
$G = EN(Q \tin G)$). 
Clearly the commutator subgroup $[E, Q]=[[P,Q], Q]$ coincides with $E=[P,Q]$. 
Hence  (4.3b) in \cite{da} holds with $Q$ here, in the place of $K$ there.

As  the  index of $HS$ in $G$ is a power of $p$, and $PQ=N$ is a normal
 subgroup of $G$, we conclude that $HS$ contains a $q$-Sylow subgroup of 
$PQ$. Hence $HS$ contains a $P$-conjugate of $Q$. Therefore, we may 
replace $H$ and $\lambda$  by some $P$-conjugates, and  
assume that $HS$ contains $Q$.

The subgroup  $H \cap (E \times S)$ of $E \times S$ is equal to
$(H \cap E ) \times  (H \cap S)$, since $|E|$
and $|S|$ are relatively prime. Note that $S$ centralizes $E$, as the 
latter is a subgroup of $P$.
 This implies that 
$$
HS \cap (E \times S)= 
(H \cap (E \times S)) \cdot S= (H \cap E) \times S.
$$
Hence 
$H \cap E = HS  \cap E$.
 Thus  $H \cap E$ is a normal subgroup of $HS$.
Furthermore, the restriction $\lambda|_{H\cap E}$ of $\lambda$ to $H\cap E$, 
is a linear character of $H \cap E$ that is clearly $H$-invariant. 
It is also $S$-invariant, as $S$ centralizes $E  \geq H \cap E$. 
Hence $\lambda|_{H \cap E}$ is $HS$-invariant.  We conclude that the 
restriction  of the irreducible character  $ \lambda^{HS} $
of $HS$ to $H \cap E$ is a multiple of the linear character 
$\lambda|_{H \cap E}$. 
Of course the irreducible character $\lambda^{HS}$  of $HS$ induces
 irreducibly to $\chi \in \Irr(G)$, and lies above a non-trivial character 
of $Z(E)$ (as $Z(E) \leq Z(P)$ and $\zeta \in \Irr(Z(P))$ is non-trivial). 
Hence we can apply Lemma (4.4) and its Corollary   (4.8) of \cite{da}, 
using $HS$ here in the place of $H$ there,  and $\lambda^{HS}$ here
 in the place of $\phi$ there. We conclude that 
$HS \cap E = H \cap E$ is a maximal abelian subgroup of $E$.

Let $\bar{P}:=P / Z(P)$. Then $\bar{P}$ is a symplectic $\ZZ_p(G)$-module.
According to the hypotheses of the theorem, 
 $Z(P)$ is the maximal abelian  $G$-invariant  subgroup of $P$.
 Hence $\bar{P}$ is an anisotropic 
$\ZZ_p(G)$-module. If $\bar{E}$ is the image of $E$ in $\bar{P}$, i.e., 
$\bar{E} \cong E /Z(E)$, then $\bar{E}$ is a symplectic $\ZZ_p(G)$-submodule 
of $\bar{P}$, as $E$ is normal in $G$. Furthermore, $\bar{E}$
 is $\ZZ_p(HS)$-hyperbolic as $HS \cap E$ is a maximal abelian $HS$-invariant 
subgroup of $E$. Since  the index $[G: HS]$ is a power of $p$, Theorem 
\ref{m.t1} forces $\bar{E}$ to be trivial. 
Hence $E = Z(E)$ is abelian, and the claim follows. 
\end{proof}

Now  $E = [P,Q]$ is an abelian subgroup of $P$  normal in $G$.
According to the hypotheses of the theorem,  $Z(P)$   is the  maximal 
such subgroup.
Therefore  $1 \leq [P,Q] \leq Z(P) \leq Z(G)$.
So  $Q$ centralizes $[P,Q]$, which implies that 
$[P,Q,Q]=1$ and thus $[P,Q]=[P,Q,Q]=1$.

This completes the proof of the theorem.
\end{proof}

\begin{lemma}\mylabel{m.l1}
Assume that $G$ is a  $p, q$-group, 
where $p$ and $q$ are distinct odd primes.
Let  $M=P \times S$  be 
  a normal subgroup of $G$, where $P$ is a $p$-group and $S$ is a $q$-group.
Assume further that $ \beta$ is a $G$-invariant irreducible character of $S$ that can be extended 
to a $q$-Sylow subgroup $Q$ of $G$. 
Let $\alpha \in \Irr(P)$ be a $Q$-invariant character of $P$.
Then there exists an irreducible character $\chi \in \Irr(G | \alpha \times \beta)$ with 
$\chi(1)_q = \beta(1)$.
\end{lemma}

\begin{proof}
Let $A/S$ be a $p$-Sylow subgroup of $G/S$. Then the 
$q$-special character $\beta \in \Irr(S)$ can be extended to $A$, as 
$(|A:S|, \beta(1)o(\beta))=1$, (see Corollary 8.16 in \cite{is}). 
As $\beta$ is also extendible to a $q$-Sylow subgroup of $G$, we conclude that 
it is extendible to $G$ (see Corollary  11.31 in \cite{is}). 
Let $\beta^e$ be an extension of $\beta$ to $G$.

The irreducible character $\alpha$ of $P$ is $Q$-invariant.
Hence Proposition 21.5 in \cite{wo} implies that the canonical  extension $\alpha^e$ of 
$\alpha$ to $P\rtimes Q$ is the unique $p$-special character of $P\rtimes Q$ lying above $\alpha$.
Furthermore, the index of the group $P \rtimes Q$ in $G$ is a $p$-number, as
 $Q$ is a $q$-Sylow subgroup of $G$. Therefore, the same proposition 
implies that any character of $G$ above $\alpha^e$ is $p$-special.
Let $\Psi \in \Irr(G | \alpha^e)$ be such. 

As $\Psi$ is a $p$-special character,  while $\beta^e$ is a $q$-special
 character of $G$, Theorem 21.6 in \cite{wo} 
implies that the product $\chi:= \Psi \cdot \beta^e$ is an irreducible 
character of $G$.
Obviously $\chi$ lies above $\alpha \times \beta$, while 
$\chi(1)_q = \beta(1)$, as $\Psi(1)$ is a $p$-number.
So the lemma follows.
\end{proof}

We can now complete the proof of the $t=1$ case.
\begin{proof}[Proof of Theorem \ref{ma.t1}]
According to Remark  \ref{ma.r1}, the character $\bbi{1}_1$ is 
$\mbgi{1}$-invariant, while $\mbsi{1} = Z(\mbgi{1}_1)$ is maximal 
among the abelian  $\mbgi{1}$-invariant subgroups of $\mbgi{1}_1$. 
Furthermore, $\mbgi{1}_1$ is a $q$-group (as a section of $G_1=Q_1$), while 
$\mbgi{1}_2/ \mbgi{1}_1$ is a $p$-group. 
Hence $\mbgi{1}_2= \mbpi{1}_2 \ltimes \mbgi{1}_1$, where $\mbpi{1}_2$ is a
$p$-Sylow subgroup of $\mbgi{1}_2= \mbgi{1}_2(\bbi{1})$.
  
We can now apply Theorem \ref{m.t2}, with 
the roles of $p$ and $q$ interchanged. 
So we work with 
$\mbgi{1}, \mbgi{1}_2, \mbgi{1}_1, \mbsi{1}$ and $\mbpi{1}_2$ 
in the place of 
$G, N, M=P,  Z(P)$ and $Q$, respectively. 
As $S$ there, we take the trivial group here.
 We also use $\bbi{1}_1 =\Ti{1}_1$ and $\bzi{1}$ here, 
 in the place of $\alpha$ and $\zeta$ there. 
 Of course $\beta \in \Irr(S)$ there, is the trivial character here.
Clearly,  $\bzi{1}$ is a faithful $\mbgi{1}$-invariant character of
 $\mbsi{1}= Z(\mbgi{1}_1)$.  Thus \eqref{m.e2} holds.
Also \eqref{m.e3} holds trivially  as $1 \in \Irr(1)$ is $\mbgi{1}$-invariant.

We are  missing  a monomial  character of $\mbgi{1}$  
that will play the role of $\chi \in \Irr(G)$ there. 
For that we will use some baby $\pi$-special theory.
Indeed, the character $\bbi{1}_1$
 is $q$-special, as a character of the $q$-group $\mbgi{1}_1$.
Furthermore, $\mbgi{1}_1$ is a normal subgroup of $\mbgi{1}$, 
while $\bbi{1}_1$ is $\mbgi{1}$-invariant. 
Hence Corollary 4.8 in \cite{ga}, implies  that there exists
a $q$-special character $\Ti{1} \in \Irr(\mbg)$  
that lies above $\bbi{1}_1$, and thus above $\bzi{1}$.
According to Proposition \ref{ma.p1}, the character $\Ti{1}$ is monomial.
Now all the hypothesis of Theorem \ref{m.t2} are satisfied. 
We conclude that $\mbpi{1}_2$ centralizes $\mbgi{1}_1$. 
Hence $\mbgi{1}_2= \mbpi{1}_2 \ltimes \mbgi{1}_1= \mbpi{1}_2 \times
\mbgi{1}_1 $.
This completes the proof of  Theorem \ref{ma.t1}.
\end{proof}

\section{The $t=2$ case } \mylabel{ma.se3}

In this section we handle the $t=2$ case. So we ultimate prove that
 the group   $\mbgi{2}_3$ is nilpotent.

For the rest of the section we assume that $n \geq 3$.
We also fix  the subsystem 
\begin{subequations}\mylabel{f.e10}
\begin{gather}
1= G_0 \unlhd G_1 \unlhd G_2 \unlhd G_3  \unlhd G,\\
\{ \chi_i \in \Irr(G_i) \}_{i=0}^3\\
\{ Q_{2i-1}, P_{2r} |\beta_{2i-1}, \alpha_{2r}\}_{i=1, r=0}^{2,1}
\end{gather}
\end{subequations}
of the system \eqref{ma.e1}.
(Note that this is the $m=3$ case for \eqref{ma.sm}. So $k=1$, 
 $P_{2k}^*= P_2$ and $\beta_{2k-1,2k}= \beta_{1,2}$.)
We also  pick  and fix the groups $\qw$  and $\hap$, for the above system,
so that Theorems \ref{sy.H} and \ref{sy.t1} hold.
This way the groups  $T$ and $U$ are also fixed. 

According to  \eqref{ma.sm1}, after the first set of reductions  the above
system reduces to 
\begin{subequations}\mylabel{fi.e1}
\begin{gather}
1= \mbgi{1}_0 \unlhd \mbgi{1}_1 \unlhd \mbgi{1}_2 \unlhd  \mbgi{1}_3 \unlhd \mbgi{1} \\
\{ \Ti{1}_i \in \Irr(\mbgi{1}_i) \}_{i=0}^3\\
\{ \mbqi{1}_{2i-1}, \mbpi{1}_{2r} |\bbi{1}_{2i-1}, \bai{1}_{2r}\}_{i=1, 
r=0}^{2,1}
\end{gather}
\end{subequations}
According to Remark \ref{ma.r1}, the character 
$\bbi{1}_1$ is  $\mbgi{1}$-invariant. Furthermore, 
Theorem \ref{ma.t1} implies that  the group $\mbgi{1}_2$ is nilpotent.
Hence $\mbqi{1}_{1,2} = C(\mbpi{1}_2 \tin \mbqi{1}_1) = \mbqi{1}_1$.
So the character $\bbi{1}_{1,2}$ coincides with $\bbi{1}_1$
(its $\mbpi{1}_2$-Glauberman correspondent). Thus
$\mbgi{1}(\bbi{1}_{1,2}) = \mbgi{1}(\bbi{1}_1) = \mbgi{1}$.
In particular, for the $q$-Sylow subgroup of $\mbgi{1}(\bai{2})$ we have 
$\mbqwi{1}(\bbi{1}_{1,2}) = \mbqwi{1}$.
This, along with \eqref{f.eq0},   implies (note that $k=1$)
\begin{equation}\mylabel{fi.e2}
\text{ Image of $\qw(\beta_{1,2})$ in $\Aut(\mbpi{1}_2)$ }=
\text{ Image of $\mbqwi{1}$  in $\Aut(\mbpi{1}_2)$.}
\end{equation}
Furthermore,  \eqref{f.q1} holds.

For the next set of reductions, we first have to shift the 
system \eqref{fi.e1} (see \eqref{ma.sm2}), and then we perform the
 $\mbi{1}(\bai{1}_2)$-invariant reductions.
We end up with the system \eqref{ma.sm3}, which for $m=3$ (i.e., the case here) 
gives
\begin{subequations}\mylabel{fi.e3}
\begin{gather}
1\unlhd \mbpi{2}_2 \unlhd \mbgi{2}_3  \unlhd \mbgi{2},\\
  1, \bai{2}_2, \Ti{2}_3\\
\{ \mbqi{2}_{3}, \mbpi{2}_{2} |\bbi{2}_{3}, \bai{2}_{2} \}
\end{gather}
\end{subequations}
Note that,  as $\mbqi{1}_1$ is a normal subgroup of $\mbgi{1}$ that 
centralizes $\mbpi{1}_2$, Remark \ref{lm.rem2} 
implies that the group $\mbqi{1}_1$ maps isomorphically to a normal subgroup
 $\mbqi{2}_1$ of $\mbgi{2}$ that centralizes $\mbpi{2}_2$. 
In addition, its irreducible character $\bbi{1}_1 \in \Irr(\mbqi{1}_1)$ maps to
$\bbi{2}_1 \in \Irr(\mbqi{2}_1)$.
The fact that $\bbi{1}_1$ is $\mbgi{1}$-invariant, while $\mbgi{2}$ is a
 section of $\mbgi{1}$,  implies that 
$\bbi{2}_1$ is $\mbgi{2}$-invariant, and thus $\mbqwi{2}$-invariant.
 This, along with Theorem \ref{f.t4},  implies
\begin{subequations}\mylabel{fi.e4}
\begin{equation}\mylabel{fi.e4a}
\text{ Image of $\mbqwi{2}$ in $\Aut(\mbpi{2}_2)$ } = 
\text{ Image of $\qw(\beta_{1,2})$ in $\Aut(\mbpi{2}_2)$. }
\end{equation}
Even more, Theorem \ref{f.t3} implies 
\begin{gather}\mylabel{fi.e4b}
\text{ Any faithful linear limit of $(\mbti{2}, \mbri{2}, 
\bei{2}, \mbpi{2}_2, \bai{2}_2)$ is also a } \notag\\
\text{ faithful linear limit of  $(T, R,  \eta, P_{2}, \alpha_2)$}.
\end{gather}
\end{subequations}
According to Remark \ref{ma.r2},  the character $\bai{2}_2$ is 
$\mbgi{2}$-invariant.  Hence the $q$-Sylow subgroup $\mbqwi{2}$ 
of $\mbgi{2}(\bai{2}_2)$ satisfies
\begin{equation}\mylabel{fi.e5}
\mbqwi{2} \in \Syl_q(\mbgi{2}).
\end{equation}

Having fixed the system \eqref{f.e10}, and the groups $\qw$ and $ T$,
  we can get a new   system,  via
 Corollary \ref{n.co2}. That is, we get a new character tower  
\begin{subequations}\mylabel{f.e11}
\begin{equation}
1=\cn_0, \cn_{1}, \cn_{2}
\end{equation}
for the series 
\begin{equation}
1=G_0 \unlhd G_1 \unlhd G_2  \unlhd  G,
\end{equation}
along with a triangular set 
\begin{equation}
\{1=\pn_0, \pn_{2}, \qn_1 | 1=\an_0, \an_2, \bn_{1} \}
\end{equation}
\end{subequations}
and a $q$-Sylow subgroup $\qwn$ of $G(\an_2)$, 
so that  \eqref{n.e222} holds. In particular, we have 
\begin{subequations}\mylabel{f.e12}
\begin{gather}
P_{2}= \pn_2= \pns_2  \text{ and } \alpha_2= \an_2= \ans_2, \\ 
\bn_{1,2} \text{ extends to } \qw= \qwn.
\end{gather}
\end{subequations}
Observe that   the character $\bn_{1,2}$ is $q$-special, as it is an 
irreducible character of the $q$-group $\qn_{1,2}$. Hence the fact that it 
 extends to   a $q$-Sylow subgroup of $G(\an_2)$ implies (see Corollary
11.31 in \cite{is}) that it extends to $G(\an_2)$.
As expected, 
 we write $I^{\nu}$ for the image of $\qw^{\nu}(\bn_{1,2})= \qw^{\nu} $ 
 $\Aut(P_{2}^{*,\nu})$. Therefore  \eqref{f.e12} implies that $I^{\nu}$ is the image of $\qw $ 
in $\Aut(P_{2})$.

Assume  we perform the reduction procedure described in Section \ref{ma.se1}
for this new system. We assume that $m=n=2$ here.
We keep the same notation as earlier, with the addition 
of the  superscript $\nu$ to any group that
 refers to the new system.  So we first start with 
$\ma^{\nu}(\bn_{1})$-invariant linear reductions of the quintuple
$(G, S^{\nu}, \zeta^{\nu}, G_1^{\nu}, \cn_1)$, where $S^{\nu}=1=S$ and
 $\zeta^{\nu}=1= \zeta$. Furthermore, $G_1^{\nu}=G_1$ and  $\cn_1 = \bn_1$. 
Let
\begin{equation}\mylabel{f.e12.5}
 (\mbgni{1}, \mbsni{1}, \bzni{1}, \mbgni{1},  \bbni{1}_1)
\end{equation}
be an $\ma^{\nu}(\bn_1)$-invariant faithful limit of  $(G, S^{\nu}, \zeta^{\nu}, G_1^{\nu}, \cn_1)$.
Note that,  even though  we start with  $S^{\nu}=S$ and $G_1^{\nu}=G_1$, 
 the
limit groups $\mbsni{1}, \mbgni{1}_1$ and $\mbsi{1}, \mbgi{1}_1$ 
are not the same, 
as the reductions we perform may be different. 
At this point the system \eqref{f.e11}  has been reduced to 
\begin{subequations}\mylabel{f.e13}
\begin{gather}
1= \mbgni{1}_0  \unlhd \mbgni{1}_1 \unlhd \mbgni{1}_2 \unlhd    \mbgni{1},\\
1, \bbni{1}_1, \Tni{1}_2  \\
\{ 1= \mbpni{1}_0, \mbpni{1}_2, \mbqni{1}_1 |
 1=\bani{1}_0, \bani{1}_2, \bbni{1}_1\}
\end{gather}
\end{subequations}
where,
(according to Theorem \ref{ma.t1}), the limit group 
$\mbgni{1}_2$ is nilpotent, i.e., 
$\mbgni{1}_2= \mbpni{1}_{2} \times \mbqni{1}_1$. 
So 
$$
\mbqni{1}_1 = \mbqni{1}_{1,2},
$$
and the limit character $\bbni{1}_{1}$ coincides with
 its $\mbpni{1}_2$-Glauberman correspondent, i.e.
$$
\bbni{1}_1 = \bbni{1}_{1,2}.
$$
According to Theorem \ref{lmt.t4}, 
the character $\bbni{1}_{1,2}=\bbni{1}_1$ extends 
to $\mbqwni{1}(\bbni{1}_{1,2})$ , 
as the character $\bn_{1,2}$ extends to $\qw= \qwn$. 
But 
 $\bbni{1}_1$ is $\mbgni{1}$-invariant by   Remark \ref{ma.r1}. 
Hence $\bbni{1}_1$ extends to $\mbqwni{1}$.
Even more, in this case we have 
 $\mbpsni{1}_{2k}= \mbpsni{1}_{2}= \mbpni{1}_2$. Therefore, 
 if  $\mbini{1}$ denotes the image of $\mbqwni{1}(\bbni{1}_{1,2})$ in
 $\Aut(\mbpsni{1}_2)$, 
the above equations, along with 
\eqref{f.eq0} for  the new system, imply 
\begin{equation}\mylabel{f.e14}
\text{ Image of ($\qw=\qwn$) in $ \Aut(\mbpni{1}_2)$ } =  \mbini{1}
=\text{ Image of $\mbqwni{1}$ in $ \Aut(\mbpni{1}_2)$. }
\end{equation}
Furthermore, \eqref{m.rs1} implies that $\mbrni{1} \unlhd \mbpni{1}_2$, while 
its irreducible character $\beni{1} $ lies under $\bani{1}_2$.

To perform the next set of reductions on \eqref{f.e13},  we have to shift 
it first (see \eqref{ma.sm2} with  $m=2$). So we get 
\begin{subequations}\mylabel{f.e15}
\begin{gather}
1 \unlhd \mbpni{1}_2 \unlhd   \mbgni{1},\\
1, \bani{1}_2  \\
\{ \mbpni{1}_2 | \bani{1}_2\}
\end{gather}
\end{subequations} 
Now we can take 
a  $\mbni{1}(\bani{1}_2)$-invariant faithful linear  limit
\begin{equation}\mylabel{f.e15.5}
(\mbgni{2}, \mbrni{2}, \beni{2}, \mbpni{2}_2, \bani{2}_2),   
\end{equation}
 of the quintuple
$(\mbgni{1}, \mbrni{1}=1, \beni{1}=1, \mbpni{1}_2, \bani{1}_2)$. This way 
  \eqref{f.e15}  is reduced to 
\begin{subequations}\mylabel{f.e16}
\begin{gather}
1 \unlhd \mbpni{2}_2 \unlhd   \mbgni{2},\\
1, \bani{2}_2  \\
\{ \mbpni{2}_2 | \bani{2}_2\}
\end{gather}
\end{subequations} 
Then, according to Remark \ref{ma.r2}, the center 
$\mbrni{2}$ of $\mbpni{2}_2$ is maximal among the 
abelian $\mbgni{2}$-invariant subgroups of $\mbpni{2}_2$, while 
$\bani{2}_2$ is $\mbgni{2}$-invariant. Hence 
\begin{equation}\mylabel{f.e17}
\mbgni{2}(\bani{2}_2) = \mbgni{2}(\basni{2}_2)= \mbgni{2}.
\end{equation}
Thus the $q$-Sylow subgroup   $\mbqwni{2}$ of  $\mbgni{2}(\bani{2}_2)$ is actually 
a $q$-Sylow subgroup of $\mbgni{2}$.
Furthermore, Remark \ref{lm.rem2} implies that the normal subgroup
 $\mbqni{1}_1$ of $\mbgni{1}$,  that centralizes $\mbpni{1}_2$, maps 
isomorphically to a normal subgroup $\mbqni{2}_1$ of the limit group 
$\mbgni{2}$. So 
\begin{equation}\mylabel{f.e18}
\mbqni{1}_1 \cong \mbqni{2}_1 \unlhd  \mbgni{2}.
\end{equation}
Therefore the group $M:= \mbpni{2}_2 \times \mbqni{2}_1$ is a normal subgroup 
of $\mbgni{2}$. Under the isomorphism in \eqref{f.e18}, 
 the character $\bbni{1}_1$ of $\mbqni{1}_1$ maps 
 to the  character $\bbni{2}_1 \in \Irr(\mbqni{2}_1)$. 
Note that $\mbqni{2}_1$ is the faithful linear limit of 
$\mbqni{1}_1$, under  the $\mbni{1}(\bani{1}_2)$-invariant linear 
reductions we perform. In addition,  $\bbni{2}_1$ is the faithful 
linear limit of $\bbni{1}_1$ and $M= \mbpni{2}_2 \times \mbqni{2}_1$ that  of 
$\mbpni{1}_2 \times \mbqni{1}_1 = \mbgni{1}_2$.
The character $\bbni{2}_1$ is $\mbgni{2}$-invariant as $\bbni{1}_1$ is
 $\mbgni{1}$-invariant. Hence $\bbni{2}_1$ is $\mbqwni{2}$-invariant. 
This, along with Theorem \ref{f.t4}, implies
\begin{subequations}\mylabel{f.e19}
\begin{equation}\mylabel{f.e19a}
\text{ Image of $\mbqwni{2}$ in $\Aut(\mbpni{2}_2)$ } = 
\text{ Image of $\qw^{\nu}$ in $\Aut(\mbpni{2}_2)$. }
\end{equation}
Furthermore, Theorem \ref{f.t3} implies for the new system
\begin{gather}\mylabel{f.e19b}
\text{ Any faithful linear limit
 of $(\mbtni{2}, \mbrni{2}, 
\beni{2}, \mbpni{2}_2, \bani{2}_2)$ is also a } \notag\\
\text{ faithful linear limit of  $(T^{\nu}, R^{\nu}, 
\eta^{\nu}, P_{2}^{\nu}, \an_2)$}.
\end{gather}
(Observe that \eqref{f.e19} is the analogue  of \eqref{fi.e4} for the 
new system \eqref{f.e11}.)
Note that $R^{\nu} = 1 = R$ while $\eta^{\nu}= 1= \eta$. 
Furthermore, $P_{2}^{\nu}= P_2$ and $\an_2 = \an$, by (\ref{f.e12}a).
In addition, the image $I^{\nu}$ of $\qwn(\bn_{1,2})$ in $\Aut(\pn_2)$
equals the image of $\qw$ in $\Aut(\pn_2)$, as $\qw= \qwn = \qwn(\bn_{1,2})$
 by (\ref{f.e12}b). Hence the image $I$ of $\qw(\beta_{1,2})$ 
in $\Aut(P_{2})$,  is a subgroup of $I^{\nu}$. 
So $T^{\nu}= P_2^{\nu} \rtimes I^{\nu} = P_2 \rtimes I^{\nu}$
contains $T = P_2 \rtimes I$, while $T \geq \pn_2 = P_2$.
This, along with Remark \ref{lm.rem1}, Definition \ref{lm.d1} 
and \eqref{f.e19b}, implies that
\begin{gather}\mylabel{f.e19c}
\text{If  $(\mbtn, \mbrn, \ben, \mbpn_{2}, \ban_2)$ is a  faithful
 linear limit 
 of $(\mbtni{2}, \mbrni{2}, 
\beni{2}, \mbpni{2}_2, \bani{2}_2)$  } \notag\\
\text{ then $(\mbtn \cap T, \mbrn, \ben, \mbpn_{2}, \ban_2)$ is a 
 faithful linear limit of  $(T, R, \eta, P_{2}, \alpha_2 )$}.
\end{gather}
\end{subequations}

We give the rest of the proof  in a series  of steps
\setcounter{step}{0}
\begin{step} \mylabel{f.s1}
The character $\bbni{2}_1$  is $\mbgni{2}$-invariant, and extends to a $q$-Sylow subgroup of 
$\mbgni{2}$.
\end{step}
\begin{proof}
The group  $\mbqwni{2}$ is both a faithful invariant limit of 
$\mbqwni{1}$, and isomorphic to the latter group. Similarly, 
 $\mbqni{2}_1$ is both  a  faithful invariant limit
of, and isomorphic to,    $\mbqni{1}_1$. 
In addition, $\bbni{2}_1$ is the faithful linear limit of $\bbni{1}_1$.
Hence the fact that  $\bbni{1}_1$ is
 $\mbgni{1}$-invariant  implies that  $\bbni{2}_1$  is 
$\mbgni{2}$-invariant. 
Furthermore, as $\bbni{1}_1$ extends to $\mbqwni{1}$, we conclude that 
$\bbni{2}_1$ extends to $\mbqwni{2}$, that is a $q$-Sylow subgroup of 
$\mbgni{2}$.
So the first step follows.
\end{proof}

\begin{step}\mylabel{f.s2}
There exists a monomial character $\Tni{2} \in \Irr(\mbgni{2})$  lying above
$\bani{2}_2 \times \bbni{2}_1$  and satisfying  
$\Tni{2}(1)_q= \bbni{2}_1(1)$.
\end{step}
\begin{proof}
As  $\bbni{2}_1$ extends to a $q$-Sylow subgroup of $\mbgni{2}$, while 
$\bani{2}_2$ is $\mbgni{2}$-invariant, and $\mbpni{2}_2 \times \mbqni{2}_1$ 
is a normal subgroup of $\mbgni{2}$,  Lemma \ref{m.l1}
implies the existence of an irreducible character
 $\Tni{2} \in \Irr(\mbgni{2})$ that lies above  $\bani{2}_2 \times \bbni{2}_1$ 
and satisfies $\Tni{2}(1)_q = \bbni{2}_1(1)$.
It suffices to show that $\Tni{2}$ is monomial.
The character $\bbni{1}_1$ lies above $\bzni{1}$ (see \eqref{f.e12.5}), 
hence its faithful linear limit
$\bbni{2}_1$ lies above the faithful linear limit $\bzni{2}$ of $\bzni{1}$, 
(see \eqref{ma.rs1b} and the following remarks for 
the definition of $\bzni{2}$). In addition, 
 $\bani{2}_2$ lies above $\beni{2}$ (see  \eqref{f.e15.5}). 
Hence $\Tni{2} \in \Irr(\mbgni{2})$ lies above $\beni{2} \times \bzni{2}
\in \Irr(\mbrni{2} \times \mbsni{2})$. Therefore, 
 Proposition \ref{ma.p2} implies that $\Tni{2}$ is monomial.
\end{proof}

\begin{step}\mylabel{f.s3}
$$
\text{ Image of $\qw \cap G_3 $  in $\Aut(\mbpni{2}_2) $ } = 1.
$$
\end{step}
\begin{proof}
Let $N$ be any normal subgroup of $\mbgni{2}$ with 
$M= \mbpni{2}_2 \times \mbqni{2}_1 \leq N \unlhd \mbgni{2}$, and  $N / M$ a $q$-group. 
Then  $N = \mbpni{2}_2 \rtimes Q$ for some $q$-group $Q$, with $\mbqni{2}_1 \unlhd Q$. 
Furthermore,  as we have seen, $Z(\mbpni{2}_2) = \mbrni{2}$ is maximal among 
the abelian 
 $\mbgni{2}$-invariant subgroups of $\mbpni{2}_2$, while  $\bani{2}_2 
\in \Irr(\mbpni{2}_2)$
 lies above the $\mbgni{2}$-invariant faithful irreducible character $\beni{2} \in \Irr(\mbrni{2})$.
This, along with Step \ref{f.s1} and \ref{f.s2},  implies
 that we can apply Theorem \ref{m.t2} 
with the groups $\mbgni{2}, \mbpni{2}_2, \mbqni{2}_1,  Q$ here,  in the place of the groups 
$G, P, S, Q$ there, and the characters $\Ti{2}, \bani{2}_2, \bbni{2}_1$ and $\beni{2}$ here, 
in the place of  $\chi, \alpha, \beta$ and $\zeta$ there.
We conclude that any such  normal subgroup $N$ of $\mbgni{2}$ is nilpotent, i.e., 
\begin{gather}\mylabel{f.e20}
\text{ $Q \in \Syl_q(N)$ centralizes $\mbpni{2}_2 \in \Syl_p(N)$, whenever }\\
\text{ $M \unlhd N \unlhd \mbgni{2}$ with $N/ M$ a $q$-group. } \notag
\end{gather}

The group $G_3$  is a normal subgroup of $G$ that contains $G_1= G_1^{\nu} $ and $G_2 = G_{2}^{\nu}$. 
Hence, see Remark \ref{lm.rem1},  when the normal series (\ref{f.e11}b) reduces to 
(\ref{f.e16}a), the group $G_3$ reduces to a normal subgroup $\mbgni{2}_3$ 
of $\mbgni{2}$. Furthermore,  $\mbgni{2}_3 /M$ is a $q$-group 
as $G_3 / G_2$ is a $q$-group, and $M $ is the limit of $\mbgni{1}_2$ and thus of $G_2$.
As $\mbqwni{2}$ is a $q$-Sylow subgroup of $\mbgni{2}$, we get that 
$\mbqwni{2} \cap \mbgni{2}_3$ is a $q$-Sylow subgroup of $\mbgni{2}_3$. 
Hence \eqref{f.e20} implies 
\begin{equation*}
\text{  $\mbqwni{2} \cap \mbgni{2}_3$ centralizes $\mbpni{2}_2 \in \Syl_p(\mbgni{2}_3)$.}   
\end{equation*}
This, along with \eqref{f.e19a} implies 
$$
\text{ Image of  $\qw^{\nu} \cap G_3 $ in $\Aut(\mbpni{2}_2) $ } =
\text{ Image of $\mbqwni{2}  \cap \mbgni{2}_3 $ in $\Aut(\mbpni{2}_2) $ } =1.
$$
As  $\qw= \qw^{\nu}$,  Step \ref{f.s3} follows.
\end{proof}

Assume that $(\mbt, \mbr, \be, \mbp_2, \ba_2)$  and 
$(\mbtn, \mbrn, \ben, \mbpn_2, \ban_2)$ are faithful linear limits of 
the quintuples $(\mbti{2}, \mbri{2}, \bei{2}, \mbpi{2}_2, \bai{2}_2)$ and 
$(\mbtni{2}, \mbrni{2}, \beni{2}, \mbpni{2}_2, \bani{2}_2)$, respectively. 
Then according to \eqref{fi.e4b} and \eqref{f.e19c} the quintuples
$(\mbt, \mbr, \be, \mbp_2, \ba_2)$ and 
$(\mbtn \cap T, \mbrn, \ben, \mbpn_2, \ban_2)$ are  both  faithful 
linear limits of  $(T, R, \eta, P_2, \alpha_2)$.
Hence Theorem \ref{lp.t1} and Corollary \ref{lp.c3}
 imply that $\mbp_2/ \mbr$ and $\mbpn_2/ \mbrn$ are
 isomorphic anisotropic symplectic $\ZZ_p(T(\alpha_2)/ P_2)$-modules.
The group $T= P_2 \rtimes I$ fixes $\alpha_2 \in \Irr(P_2)$, as $\qw$ 
fixes $\alpha_2 = \alpha_2^*$ and $I$ is the image of $\qw(\beta_{1,2})$ 
in $\Aut(P_2)$. Hence $T(\alpha_2)/ P_2 $ is naturally isomorphic to 
the $q$-group $I$. 
So 
\begin{equation}\mylabel{f.e22}
\mbp_2 / \mbr \cong \mbpn_2 / \mbrn \text{ as anisotropic symplectic 
$\ZZ_p(I)$-modules. }
\end{equation}
In view of Step \ref{f.s3},  the group 
$\qw \cap G_3$ centralizes $\mbpni{2}_2$. As $\mbpn_2$ is a section of 
$\mbpni{2}_2$, we conclude that $\qw(\beta_{1,2})  \cap G_3$ centralizes
$\mbpn_2$. Hence it also centralizes the factor group $\mbpn_2/ \mbrn$. 
As the latter factor group is isomorphic to $\mbp_2/ \mbr$ as $I$-modules,
and $I$ is the image of $\qw(\beta_{1,2})$  in $\Aut(P_2)= \Aut(\pn_2)$, 
  we conclude that  $\qw(\beta_{1,2}) \cap G_3$ also  centralizes
$\mbp_2/ \mbr$. (The action of $I$ on the above two isomorphic  factor groups 
is defined in Corollary \ref{lp.c3}.)
If  $\mbgi{2}_3$ denotes the limit group to which $G_3$ reduces, 
as $G$ reduces to $\mbgi{2}$, then \eqref{fi.e4a} implies 
\begin{equation}\mylabel{f.e23}
I_3:= \text{ Image of $\qw(\beta_{1,2}) \cap G_3$ in $\Aut(\mbpi{2}_2)$}=
 \text{ Image of $\mbqwi{2} \cap \mbgi{2}_3 $ in $\Aut(\mbpi{2}_2)$.} 
\end{equation}
Since $\mbp_2$ is a section of $\mbpi{2}_2$,
the fact that $\qw(\beta_{1,2}) \cap G_3$  centralizes $\mbp_2/ \mbr$, 
 implies  
\begin{equation}\mylabel{f.e24}
I_3 \text{ centralizes } \mbp_2/ \mbr.
\end{equation}

Let $V := \mbpi{2}_2 / \mbri{2}$. Then $V$ is anisotropic
 $\ZZ_p(\mbgi{2})$-module (by Remark \ref{ma.r2}). 
Thus $V$, when written additively,  is the direct sum 
\begin{subequations}\mylabel{f.e25}
\begin{equation}\mylabel{f.e25i}
V = V_1 \dotplus V_2,
\end{equation}
of the perpendicular $\ZZ_p(\mbgi{2})$-modules, 
 $V_1 = C( \mbgi{2}_3 \tin V)$ and $V_2= [V, \mbgi{2}_3]$.
Note that, as $\mbqwi{2}$ is a $q$-Sylow subgroup 
of $\mbgi{2}$, see \eqref{fi.e5}, we get that 
$\mbqwi{2}_3 = \mbqwi{2} \cap \mbgi{2}_3$   is a $q$-Sylow subgroup 
of $\mbgi{2}_3$. 
Thus $\mbgi{2}_3= \mbqwi{2}_3 \ltimes \mbpi{2}_2$. We conclude that the 
direct summands in  \eqref{f.e25i} are 
\begin{equation}\mylabel{f.e25ii}
V_1 = C(\mbqwi{2}_3 \tin V) \text{ and } V_2= [V, \mbqwi{2}_3].
\end{equation} 
\end{subequations}
Both $V_1$ and $V_2$ are anisotropic $\ZZ_p(\mbgi{2})$-modules.

As  $(\mbt, \mbr, \be, \mbp_2, \ba_2)$ is a faithful linear limit of 
$(\mbti{2}, \mbri{2}, \bei{2}, \mbpi{2}_2, \bai{2}_2)$,  we have that 
 $U := \mbp_2/ \mbr$ is isomorphic to a factor subgroup of 
$V$. Furthermore, $U$ is isomorphic, as a symplectic $\ZZ_p(I)$-module, 
to the orthogonal direct sum $U= U_1 \dotplus U_2$, 
where $U_i$ is a limit module  for $V_i$, for each $i=1, 2$.
In view of \eqref{f.e25ii} we have 
\begin{equation}\mylabel{f.e26}
U_1= C(I_3 \tin U) \text { and } U_2 = [U, I_3],
\end{equation}
where $I_3$ is the image of $\mbqwi{2}_3$ in $\Aut(\mbpi{2}_2)$ 
(see \eqref{f.e23}).
In view of \eqref{f.e24} we get $U_2= 0$.

According to Remark \ref{ma.r2}, the center 
$\mbri{2}$ of  $\mbpi{2}_2$ is a cyclic central subgroup of 
$\mbgi{2}$,  maximal among the abelian $\mbgi{2}$-invariant 
subgroups of $\mbpi{2}_2$.
As $\mbti{2} = \mbpi{2}_2  \rtimes \mbii{2}$, where $\mbii{2}$ is the image of 
$\mbqwi{2}$ in $\Aut(\mbpi{2}_2)$, we get that
 $\mbri{2}$ is a central subgroup of $\mbti{2}$. Even more, it is maximal 
among the characteristic abelian subgroups of $\mbpi{2}_2$.
Thus Proposition \ref{lp.p3} applies to the faithful linear limit
$(\mbt, \mbr, \be, \mbp_2, \ba_2)$. So 
 $U$ is isomorphic to 
$W^{\perp}/ W$ for some maximal $I$-invariant totally isotropic subspace $W$
of $V$. Then   $W = W_1 \dotplus W_2$,
where $W_i$ is a maximal totally isotropic $I$-invariant subspace of
$V_i$, for $i = 1,2$. But  $U_2 = 0$. Hence   $W_2^{\perp}$ must equal 
 $ W_2$. Thus
$V_2$ contains a self perpendicular $I$-invariant subspace. We conclude that 
$V_2$ is  hyperbolic as a $\ZZ_p(I)$-module, and so as a 
  $\ZZ_p(\mbqwi{2})$-module.
As  $\mbqwi{2}$ is a $q$-Sylow subgroup  of $\mbgi{2}$, it  has 
$p$-power index in $\mbgi{2}$. Since  $V_2$ is an anisotropic  
symplectic  $\ZZ_p(\mbgi{2})$-module, 
Theorem 3.2 in \cite{da}, implies that $V_2$ is both  hyperbolic and 
anisotropic. Therefore $V_2$ is $0$. 
So $ V = V_1 = C( \mbqwi{2}_3 \tin V)$. 
Thus $\mbqwi{2}_3$ centralizes  $\mbpi{2}_2/ \mbri{2}$. 
We conclude that the $q$-Sylow subgroup $\mbqwi{2}_3$ 
 of $\mbgi{2}_3$ centralizes 
the $p$-Sylow subgroup  $\mbpi{2}_2$ of the same group. Hence $\mbgi{2}_3$
 is nilpotent, and Theorem \ref{ma.t2} follows.

%%% Local Variables: 
%%% mode: latex
%%% TeX-master: "thesis-ex"
%%% End:

 \section{The general case } 
\mylabel{ma.se4}

The aim in this section is to prove Theorem \ref{ma.main}, i.e. to show
 that the group  $\mbgi{t}_{t+1}$ is nilpotent.  
The ideas for the proof are already given in  
 the $t=2$ case, whose proof is a demonstration of the general argument. 
That's the reason we leave hidden some of the details of the general 
proof, already discussed in the previous section.

Assume the system \eqref{ma.e1} is  fixed, with $n \geq t+1$.
Assume further, using an inductive argument, that the groups $\mbgi{i}_{i+1}$
are nilpotent for all $i=1, \dots, t-1$. 
Thus  we can   perform  all the reductions 
described in Section  \ref{ma.se1}, until we reach  the group in question, 
i.e., the group  $\mbgi{t}_{t+1}$.  
As the last step in the reductions  depends on the parity of $t$, we first 
assume that $t=2k$ is even. (As expected,  we will see that it is enough to 
prove Theorem \ref{ma.main} in the even case.)
In addition,  we assume fixed the subsystem 
\begin{subequations}\mylabel{fg.e1}
\begin{gather}
1\unlhd G_1 \unlhd G_2 \unlhd \dots \unlhd  G_{t+1}  \unlhd G,\\
\{ \chi_i \in \Irr(G_i) \}_{i=0}^{t+1},\\
\{ Q_{2i-1}, P_{2r} |\beta_{2i-1}, \alpha_{2r}\}_{i=1, r=0}^{l,k},
\end{gather}
\end{subequations}
where, in our  case (that of an even $t$), $k=t/2$ and $l=t/2 +1$.
(Note that $k$ and $l$ are related to $t+1$ via \eqref{kl:def}.) 
 Along with that system we pick and fix the groups $\qw$ and $\hap$
so  as to satisfy the conditions in   Theorems \ref{sy.H} and \ref{sy.t1}.
This way the groups $T$ and $U$ are also fixed. 
After $t$ steps of reductions  the above system reduces to
 (see \eqref{f.smt} with  $m= t+1$)  
\begin{subequations}\mylabel{fg.e2}
\begin{gather}
1\unlhd \mbpi{t}_t \unlhd \mbgi{t}_{t+1}  \unlhd  \mbgi{t},\\
  1, \bai{t}_t, \Ti{t}_{t+1}\\
\{ \mbqi{t}_{t+1},\mbpi{t}_{t} |\bbi{t}_{t+1}, \bai{t}_{t} \}
\end{gather}
\end{subequations}
In addition note that 
the  group $\mbqi{t}_{2k-1}$ and all the groups with indices smaller 
than $t=2k$,  have become trivial. (To be precise,  all the
$q$-groups $\mbqi{t}_{2i-1}$,  with indices $2i-1$ smaller than $2k$ 
have been dropped by repeated 
shifts of the original series, as they are normal subgroups of $\mbgi{t}$
 that are  contained in  $\mbqi{t}_{2k+1}$ and  are  centralized by 
$\mbpi{t}_{2k}$,  while  their characters $\bbi{t}_{2i-1}$
 are $\mbgi{t}$-invariant. 
The same holds for the $p$-groups $\mbpi{t}_{2j}$ 
 with indices $2j$  smaller than $2k=t$.
They are normal subgroups of $\mbgi{t}$ that are  contained in 
$\mbpi{t}_{2k}$ and centralized by $\mbqi{t}_{2k+1}$, while their 
characters are  $\mbgi{t}$-invariant. This is the reason
 we drop them on the way.)
  So 
$\mbpsi{t}_{2k} = \mbpi{t}_{2k}= \mbpi{t}_{t}$.
Furthermore the last group that is been dropped is the $q$-group 
$\mbqi{t-1}_{2k-1}$. Note that $\mbqi{t-1}_{2k-1}$ is centralized by 
$\mbpi{t-1}_{2k}$. 
In addition, the irreducible character 
$\bbi{t-1}_{2k-1}$ is $\mbgi{t-1}$-invariant.  After the last
 set of reductions is been performed, the group $\mbqi{t-1}_{2k-1}$ maps
 isomorphically to a normal subgroup $\mbqi{t}_{2k-1}$ of the limit group 
 $\mbgi{t}$  that centralizes $\mbpi{t}_{2k}$, by Remark \ref{lm.rem2}.
In addition,  the irreducible character  $\bbi{t-1}_{2k-1}$ 
of $\mbqi{t-1}_{2k-1}$ maps, under the above isomorphism, to 
an irreducible character  $\bbi{t}_{2k-1}$ of $\mbqi{t}_{2k-1}$. 
Note that $\bbi{t}_{2k-1}$ is $\mbgi{t}$-invariant, as $\bbi{t-1}_{2k-1}$ 
is $\mbgi{t-1}$-invariant and $\mbgi{t}$ is a section of $\mbgi{t-1}$.
Even more, the centralizer  $\mbqi{t}_{2k-1,2k}$ of $\mbpi{t}_{2k}$ 
in $\mbqi{t}_{2k-1}$ equals $\mbqi{t}_{2k-1}$. Thus the character 
$\bbi{t}_{2k-1,2k}$ coincides with $\bbi{t}_{2k-1}$. 
So $\bbi{t}_{2k-1,2k}$ is $\mbgi{t}$-invariant.
 This implies that $\mbqwi{t}(\bbi{t}_{2k-1,2k})= 
\mbqwi{t}$.  Hence the image 
 $\mbii{t}$ of $\mbqwi{t}(\bbi{t}_{2k-1,2k})$ 
in $\Aut(\mbpsi{t}_{2k})$ is reduced to 
\begin{equation}\mylabel{fg.e3}
\mbii{t}= \text{ Image of $\mbqwi{t}$ in $\Aut(\mbpi{t}_{2k})$. } 
\end{equation}
This,  along with Theorem \ref{f.t8} and in particular \eqref{f.e3b}, implies
\begin{subequations}\mylabel{fg.se4}
\begin{equation}\mylabel{fg.e4}
\text{ Image of $\mbqwi{t}$ in $\Aut(\mbpi{t}_{2k})$ }= 
 \text{ Image of $\qw(\beta_{2k-1,2k})$ in $\Aut(\mbpi{t}_{2k})$. } 
\end{equation}
If $\mbti{t}= \mbpsi{t}_{2k} \rtimes \mbii{t}= \mbpi{t}_{2k} \rtimes 
\mbii{t}$, then  Theorem \ref{f.t7} implies 
\begin{gather}\mylabel{fg.e5}
\text{ Any faithful linear limit of $(\mbti{t}, \mbri{t}, 
\bei{t}, \mbpi{t}_{2k}, \bai{2}_{2k})$ is also a } \notag\\
\text{ faithful linear limit of  $(T, R,  \eta, P_{2k}^*, \alpha_{2k}^*)$}.
\end{gather}
Observe also that $\basi{t}_{2k}= \bai{t}_{2k}$ is $\mbgi{t}$-invariant, 
by Remark \ref{ma.r3}, as $2k=t$. Hence the fact that $\mbqwi{t}$ is 
a $q$-Sylow  subgroup of $\mbgi{t}(\basi{t}_{2k})$ implies
\begin{equation}\mylabel{fg.e6}
\mbqwi{t} \in \Syl_{q}(\mbgi{t}).
\end{equation}
\end{subequations}
(Note that \eqref{fg.se4} is the analogue of \eqref{fi.e4} and \eqref{fi.e5}
for the general case.)

According to  
 Corollary \ref{n.co2} we get a new character tower  
\begin{subequations}\mylabel{fg.e7}
\begin{equation}
1=\cn_0, \cn_{1}, \dots ,  \cn_{t}
\end{equation}
for the series 
\begin{equation}
1=G_0 \unlhd G_1 \unlhd \dots  \unlhd G_{t}  \unlhd G_n= G,
\end{equation}
along with a triangular set 
\begin{equation}
\{\qn_{2i-1}, \pn_{2r}| \bn_{2i-1},  \an_{2r} \}_{i=1, r=0}^{k, k}
\end{equation}
\end{subequations}
and a $q$-Sylow subgroup $\qwn$ of $G(\ans_{2k})$, 
so that  \eqref{n.e222} holds. In particular, we have 
\begin{subequations}\mylabel{fg.e8}
\begin{gather}
P_{2k}^*= \pns_{2k}  \text{ and } \alpha_{2k}^*= \ans_{2k}, \\ 
\bn_{2k-1,2k} \text{ extends to } \qw= \qwn.
\end{gather}
\end{subequations}
We proceed with the reductions described in Section \ref{ma.se1}
for the new system \eqref{fg.e7}. We follow the same notation 
as that of Sections \ref{ma.se1} and \ref{ma.seI}
 with the addition of a supescript 
$\nu$ on anything that refers to the new system \eqref{fg.e7}. 
So, we reduce the above new system, using  first reductions that are   
$\ma^{\nu}(\bn_1)$-invariant,  then  
$\mb^{(1), \nu}(\bani{1}_2)$-invariant that are followed by 
$\ma^{(2),\nu}(\bbni{2}_3)$-invariant and so on. 
Hence all the conclusions of Section \ref{ma.seI} are valid for the system 
\eqref{fg.e7}. 
In particular, after $t$ steps of reductions,  what is left 
of the system \eqref{fg.e7}  is (see \eqref{f.smt} with $m= t=2k$)
 \begin{subequations}\mylabel{fg.e9}
\begin{gather}
1\unlhd \mbpni{t}_{2k}   \unlhd \mbgni{t},\\
1, \bani{t}_{2k} \\
\{ \mbpni{t}_{2k} |\bani{t}_{2k}\}
\end{gather}
\end{subequations}  
According to Remark \ref{ma.r3} the character $\basni{t}_{2k}= \bani{t}_{2k}$
is $\mbgni{t}$-invariant. 
Thus 
\begin{equation}\mylabel{fg.e09} 
\mbgni{t}(\basni{t}_{2k})= \mbgni{t}(\bani{t}_{2k})= \mbgni{t}.
\end{equation}
Note that the last group dropped after $t-1$ steps of reductions is
the $q$-group $\mbqni{t-1}_{t-1}$, as according to the inductive hypothesis the
last group proved  to be nilpotent is  $\mbgni{t-1}_{t}= \mbpni{t-1}_t 
\times \mbqni{t-1}_{t-1}$ (see \eqref{ma.e7}).
It is clear that the centralizer $\mbqni{t-1}_{2k-1,2k}$ of 
$\mbpni{t-1}_{2k}= \mbpni{t-1}_{t}$ in  $\mbqni{t-1}_{2k-1}=
 \mbqni{t-1}_{t-1}$ equals $\mbqni{t-1}_{2k-1}$.
Furthermore, the irreducible character $\bbni{t-1}_{2k-1,2k}$ of 
$\mbqni{t-1}_{2k-1,2k}$ coincides with $\bbni{t-1}_{2k-1}$. 
In addition, repeated applications (at every step of reductions)
of Theorems \ref{lmt.t4}  and \ref{lmt.t8} imply that the character 
$\bbni{t-1}_{2k-1}$ extends to $\mbqwni{t-1}(\bbni{t-1}_{2k-1})$, 
(where $\mbqwni{t-1}$ is a $q$-Sylow subgroup of 
$\mbgni{t-1}(\bani{t-1}_{2k})$), 
as the character $\bn_{2k-1,2k}$ extends to $\qwn(\bn_{2k-1,2k})
= \qwn = \qw$ (see \eqref{fg.e8}). 
But the character $\bbni{t-1}_{2k-1}= \bbni{t-1}_{t-1}$ is 
$\mbgni{t-1}$-invariant, according to the inductive hypothesis (see 
\eqref{ma.e7c}). Hence  
\begin{equation}\mylabel{fg.e10}
\bbni{t-1}_{2k-1} \in \Irr(\mbqni{t-1}_{2k-1})
\text{ extends to } \mbqwni{t-1}. 
\end{equation}
Each reduction in the  $t$-th set  of reductions is
$\mbni{t-1}(\bani{t-1}_{t})$-invariant   (see the comments
following \eqref{ma.e7}). After this set of reductions is performed, 
the normal subgroup   $\mbqni{t-1}_{2k-1}$ of $\mbgni{t-1}$ that centralizes 
$\mbpni{t-1}_{2k}$,
maps isomorphically to a normal subgroup $\mbqni{t}_{2k-1}$ 
of the limit group  $\mbgni{t}$, by Remark \ref{lm.rem2}.
 Under this  isomorphism the character $\bbni{t-1}_{2k-1}$
of $\mbqni{t-1}_{2k-1}$ maps to the character $\bbni{t}_{2k-1} 
\in \Irr(\mbqni{t}_{2k-1})$. So 
\begin{gather}\mylabel{fg.e11}
\mbqni{t-1}_{2k-1} \cong \mbqni{t}_{2k-1} \unlhd \mbgni{t},
 \text{ and } \notag\\
\bbni{t-1}_{2k-1} \to \bbni{t}_{2k-1}.
\end{gather}
 (Observe this is the analogue of \eqref{f.e18}.)
Clearly $\bbni{t}_{2k-1}$ is $\mbgni{t}$-invariant as $\bbni{t-1}_{2k-1}$ is 
$\mbgni{t-1}$-invariant,  and $\mbgni{t}$ is a section of $\mbgni{t-1}$.
According to Theorem \ref{lmt.t7} (see its first part),
under this last set of $\mbni{t-1}(\bani{t-1}_{t})$-invariant reductions,
the $q$-Sylow subgroup $\mbqwni{t}$ of $\mbgni{t}(\basni{t}_{2k})$ 
satisfies 
$$\mbqwni{t}(\bbni{t}_{2k-1,2k}) \cong \mbqwni{t-1}(\bbni{t-1}_{2k-1,2k}).
$$
We conclude that 
\begin{equation}\mylabel{fg.e13}
\mbqwni{t} \cong \mbqwni{t-1}, 
\end{equation}
as $\bbni{t-1}_{2k-1,2k}= \bbni{t-1}_{2k-1}$ and 
$\bbni{t}_{2k-1,2k}= \bbni{t}_{2k-1}$ are  $\mbgni{t-1}$-
and $\mbgni{t}$-invariant respectively.
Note also that \eqref{fg.e09} implies 
\begin{equation}\mylabel{fg.e12}
\mbqwni{t}(\bbni{t}_{2k-1,2k}) = \mbqwni{t} \in \Syl_q(\mbgni{t}).
\end{equation}  
We can know easily prove
\setcounter{step}{0}
\begin{step}\mylabel{fg.s1}
The character $\bbni{t}_{2k-1} \in \Irr(\mbqni{t}_{2k-1})$ is 
$\mbgni{t}$-invariant and extends to  $\mbqwni{t} \in \Syl_q(\mbgni{t})$.
\end{step}
\begin{proof}
Follows immediately from \eqref{fg.e10}-\eqref{fg.e12}. 
\end{proof}
Of course the last set of reductions send $\mbpni{t-1}_{2k}$ to 
$\mbpni{t}_{2k}$. Thus the group 
 $M:=\mbpni{t}_{2k} \times \mbqni{t}_{2k-1}$ is a normal subgroup of 
$\mbgni{t}$. (The group $M$ is the image of 
$\mbgni{t-1}_{t} \unlhd \mbgni{t-1}$ in $\mbgni{t}$ 
   under the last set of reductions.)

Completing the list of the general properties about the new system 
\eqref{fg.e9},  we remark that if $\mbini{t}$ denotes the image of 
$\mbqwni{t}(\bbni{t}_{2k-1,2k})= \mbqwni{t}$ in 
$\Aut(\mbpsni{t}_{2k})= \Aut(\mbpni{t}_{2k})$, then 
Theorem \ref{f.t8} implies 
\begin{equation}\mylabel{fg.e14}
\mbini{t}= \text{ Image of $\mbqwni{t}$ in $\Aut(\mbpni{t}_{2k})$ } = 
\text{ Image of $\qw^{\nu}$ in $\Aut(\mbpni{t}_{2k})$. }
\end{equation}
Furthermore,  Theorem \ref{f.t7} implies 
\begin{gather}\mylabel{fg.e15}
\text{ Any faithful linear limit
 of $(\mbtni{t}, \mbrni{t}, 
\beni{t}, \mbpni{t}_{2k}, \bani{t}_{2k})$ is also a } \notag\\
\text{ faithful linear limit of  $(T^{\nu}, R^{\nu}, 
\eta^{\nu}, \pns_{2k}, \ans_{2k})$}.
\end{gather}
(We remind the reader  the definition of  $\mbtni{t}$  as the product 
$\mbtni{t}= \mbpni{t}_{2k} \rtimes \mbini{t}$.)
Of course $R^{\nu}= 1 = R$ and $\eta^{\nu}=1 = \eta$. In addition, 
$\pns_{2k} = P_{2k}^*$ and $\ans_{2k}= \alpha_{2k}^*$, by \eqref{fg.e8}.
Note also that the image $I$ of $\qw(\beta_{2k-1,2k})$in $\Aut(P_{2k}^*)$ 
is a subgroup of the image $I^{\nu}$ of $\qwn(\bn_{2k-1,2k})$ 
in $\Aut(\pns_{2k})= \Aut(P_{2k}^*)$, as $\qwn(\bn_{2k-1,2k})= \qwn =\qw$. 
Hence $T = P_{2k}^* \rtimes I$ is a subgroup of $T^{\nu}= \pns_{2k} 
\rtimes I^{\nu}= P_{2k}^* \rtimes I^{\nu}$. This, along with 
Remark \ref{lm.rem1}, Definition \ref{lm.d1} and 
\eqref{fg.e15}, implies 
\begin{gather}\mylabel{fg.e16}
\text{If  $(\mbtn, \mbrn, \ben, \mbpn_{2k}, \ban_{2k})$ is a  faithful
 linear limit 
 of $(\mbtni{t}, \mbrni{t}, 
\beni{t}, \mbpsni{t}_{2k}, \basni{t}_{2k})$  } \notag\\
\text{ then $(\mbtn \cap T, \mbrn, \ben, \mbpn_{2k}, \ban_{2k})$ is a 
 faithful linear limit of  $(T, R, \eta, P_{2k}^*, \alpha_{2k}^* )$}.
\end{gather}

The next two  steps will complete the proof of the theorem, and are identical
 to those proved in the case $t=2$.
\begin{step}\mylabel{fg.s2}
There exists a monomial character $\Tni{t} \in \Irr(\mbgni{t})$  lying above
$\bani{t}_{2k} \times \bbni{t}_{2k-1}$  and satisfying  
$\Tni{t}(1)_q= \bbni{t}_{2k-1}(1)$.
\end{step}
\begin{proof}
As  $\bbni{t}_{2k-1}$ extends to a $q$-Sylow subgroup of $\mbgni{t}$, while 
$\bani{t}_{2k}$ is $\mbgni{t}$-invariant, and $M= \mbpni{t}_{2k} 
\times \mbqni{t}_{2k-1}$ 
is a normal subgroup of $\mbgni{t}$,  Lemma \ref{m.l1}
implies the existence of an irreducible character
 $\Tni{t} \in \Irr(\mbgni{t})$ that lies above  $\bani{t}_{2k} 
\times \bbni{t}_{2k-1}$ 
and satisfies $\Tni{t}(1)_q = \bbni{t}_{2k-1}(1)$.
It suffices to show that $\Tni{t}$ is monomial.
The character $\bbni{t-1}_{2k-1}$ lies above $\bzni{t-1}$
 by \eqref{ma.e7b}.
Hence its faithful linear limit
$\bbni{t}_{2k-1}$ lies above the faithful linear limit 
$\bzni{t}$ of $\bzni{t-1}$, 
(for the definition of $\bzni{t}$ see \eqref{ma.e20}). 
In addition,  $\bani{t}_{2k}$ lies above $\beni{t}$ (see  
Proposition \ref{ma.p4}, with $t=2k$ even). 
Hence $\Tni{t} \in \Irr(\mbgni{t})$ lies above $\beni{t} \times \bzni{t}
\in \Irr(\mbrni{t} \times \mbsni{t})$. Therefore, 
 Proposition \ref{ma.p5} implies that $\Tni{t}$ is monomial.
\end{proof}

\begin{step}\mylabel{fg.s3}
$$
\text{ Image of $\qw \cap G_{t+1} $  in $\Aut(\mbpni{t}_{2k}) $ } = 1.
$$
\end{step}
\begin{proof}
Let $N$ be any normal subgroup of $\mbgni{t}$ with 
$M= \mbpni{t}_{2k} \times \mbqni{t}_{2k-1} \unlhd N \unlhd \mbgni{t}$, 
and  $N / M$ a $q$-group. 
Then  $N = \mbpni{t}_{2k} \rtimes Q$ for some $q$-group $Q$, 
with $\mbqni{t}_{2k-1} \unlhd Q$. 
Furthermore,  as we have seen in Remark \ref{ma.r3}, 
$Z(\mbpni{t}_{2k}) = \mbrni{t}$ is maximal among 
the abelian 
 $\mbgni{t}$-invariant subgroups of $\mbpni{t}_{2k}$. Furthermore, 
  $\bani{t}_{2k} \in \Irr(\mbpni{t}_{2k})$
 lies above the $\mbgni{t}$-invariant faithful irreducible 
character $\beni{t} \in \Irr(\mbrni{t})$, by Proposition \ref{ma.p4}.
This, along with Steps \ref{fg.s1} and \ref{fg.s2},  implies
 that we can apply Theorem \ref{m.t2} 
with the groups $\mbgni{t}, \mbpni{t}_{2k}, \mbqni{t}_{2k-1},  Q$ here,  
in the place of the groups 
$G, P, S, Q$ there, and the characters $\Ti{t}, \bani{t}_{2k}, 
\bbni{t}_{2k-1}$ and $\beni{t}$ here, 
in the place of  $\chi, \alpha, \beta$ and $\zeta$ there.
We conclude that any such  normal subgroup 
$N$ of $\mbgni{2}$ is nilpotent, i.e., 
\begin{gather}\mylabel{fg.e20}
\text{ $Q \in \Syl_q(N)$ centralizes $\mbpni{t}_{2k} \in \Syl_p(N)$,
 whenever }\\
\text{ $M \unlhd N \unlhd \mbgni{t}$ with $N/ M$ a $q$-group. } \notag
\end{gather}

The group $G_{t+1}$  is a normal subgroup of $G$ that contains
 $G_i= G_i^{\nu}$ for all $i=1, \dots, t$.
Hence, by  Remark \ref{lm.rem1},  the normal series (\ref{fg.e7}b) 
reduces to (\ref{fg.e9}a), the group $G_{t+1}$ 
reduces to a normal subgroup $\mbgni{t}_{t+1}$ 
of $\mbgni{t}$. Furthermore,  $\mbgni{t}_{t+1} /M$ is a $q$-group 
as $G_{t+1} / G_t$ is a $q$-group, and $M $ is the limit of
 $\mbgni{t-1}_{t}$ and thus of $G_t$.
As $\mbqwni{t}$ is a $q$-Sylow subgroup of $\mbgni{t}$, we get that 
$\mbqwni{t} \cap \mbgni{t}_{t+1}$ is a $q$-Sylow subgroup 
of $\mbgni{t}_{t+1}$. 
Hence \eqref{fg.e20} implies 
\begin{equation*}
\text{  $\mbqwni{t} \cap \mbgni{t}_{t+1}$
 centralizes $\mbpni{t}_{2k} \in \Syl_p(\mbgni{t}_{t+1})$.}   
\end{equation*}
This, along with \eqref{fg.e14} implies 
$$
\text{ Image of  $\qw^{\nu} \cap G_{t+1} $ in $\Aut(\mbpni{t}_{2k}) $ } =
\text{ Image of $\mbqwni{t}  \cap \mbgni{t}_{t+1} $ in $\Aut(\mbpni{t}_{2k}) $ } =1.
$$
As  $\qw= \qw^{\nu}$,  Step \ref{fg.s3} follows.
\end{proof}

We continue as in the case $t=2$. So assume that 
  $(\mbt, \mbr, \be, \mbp_{2k}, \ba_{2k})$  and 
$(\mbtn, \mbrn, \ben, \mbpn_{2k}, \ban_{2k})$ are faithful linear limits of 
the quintuples $(\mbti{t}, \mbri{t}, \bei{t}, \mbpi{t}_{2k}, \bai{t}_{2k})$ 
and 
$(\mbtni{t}, \mbrni{t}, \beni{t}, \mbpni{t}_{2k}, \bani{t}_{2k})$, 
respectively. 
 Then 
 $(\mbt, \mbr, \be, \mbp_{2k}, \ba_{2k})$  and 
$(\mbtn \cap T, \mbrn, \ben, \mbpn_{2k}, \ban_{2k})$
are both faithful linear limits of 
$(T, R, \eta, P_{2k}^*, \alpha_{2k}^*)$. Hence, by
Theorem \ref{lp.t1}, 
$$
\mbp_{2k}/ \mbr  \cong  \mbpn_{2k}/ \mbrn, 
$$
as  anisotropic symplectic $\ZZ_p(T(\alpha_{2k}^*)/ P_{2k}^*)$-modules.
But  $T= P_{2k}^* \rtimes I$ fixes $\alpha_{2k}^*$,  as $I$ is the image of 
$\qw(\beta_{2k-1,2k})$  in $\Aut(P_{2k}^*)$ and the group 
$\qw \in G(\alpha_{2k}^*)$ fixes $\alpha_{2k}^*$. So 
$T(\alpha_{2k}^*)/ P_{2k}^*$ is naturally isomorphic to $I$.
 We conclude that  
\begin{equation}\mylabel{fg.e22}
\mbp_{2k} / \mbr \cong \mbpn_{2k} / \mbrn \text{ as anisotropic symplectic 
$\ZZ_p(I)$-modules. }
\end{equation}
In view of Step \ref{f.s3},  the group 
$\qw \cap G_{t+1}$  centralizes 
$\mbpni{t}_{2k}$ and thus  also centralizes its  section $\mbpn_{2k}$.
This, along with \eqref{fg.e22} and the definition of $I$, implies that 
$\qw \cap G_{t+1}$ centralizes $\mbp_{2k} / \mbr$.
Let $\mbgi{t}_{t+1}$ denote the limit group to which $G_{t+1}$ reduces, 
as $G$ reduces to $\mbgi{t}$. Then \eqref{fg.e4} implies 
\begin{equation}\mylabel{fg.e23}
I_{t+1}:= \text{ Image of $\qw(\beta_{2k-1,2k}) \cap G_{t+1}$ in 
$\Aut(\mbpi{t}_{2k})$ }=
 \text{ Image of $\mbqwi{t} \cap \mbgi{t}_{t+1} $ in $\Aut(\mbpi{t}_{2k})$.} 
\end{equation}
The fact that $\qw(\beta_{2k-1,2k}) \cap G_{t+1}$  centralizes
 $\mbp_{2k}/ \mbr$, while  $\mbp_{2k}$ is a section of $\mbpi{t}_{2k}$,
 implies 
\begin{equation}\mylabel{fg.e24}
I_{t+1} \text{ centralizes } \mbp_{2k}/ \mbr.
\end{equation}

This time we define $V := \mbpi{t}_{2k} / \mbri{t}$. So $V$ is an  anisotropic
 $\ZZ_p(\mbgi{t})$-module. 
Thus $V$, when written additively,  is the direct sum 
\begin{subequations}\mylabel{fg.e25}
\begin{equation}\mylabel{fg.e25i}
V = V_1 \dotplus V_2,
\end{equation}
of the perpendicular $\ZZ_p(\mbgi{t})$-modules, 
\begin{gather}\mylabel{fg.e25ii}
V_1 =C(\mbgi{t}_{t+1} \tin V) =  C(\mbqwi{t}_{2k+1} \tin V) \\
 V_2= [V, \mbgi{t}_{t+1}]= [V, \mbqwi{t}_{2k+1}],
\end{gather} 
\end{subequations}
where  $\mbqwi{t}_{2k+1}= \mbqwi{t} \cap \mbgi{t}_{2k+1}$ is a 
$q$-Sylow subgroup 
of $\mbgi{t}_{2k+1}= \mbgi{t}_{t+1}$, as $\mbqwi{t}$ is a 
$q$-Sylow subgroup of $\mbgi{t}$.

As  $(\mbt, \mbr, \be, \mbp_{2k}, \ba_{2k})$ is a faithful linear limit of 
$(\mbti{t}, \mbri{t}, \bei{t}, \mbpi{t}_{2k}, \bai{t}_{2k})$,  we have that 
 $U := \mbp_{2k}/ \mbr$ is isomorphic to a section  of 
$V$. Furthermore, $U$ is isomorphic as a symplectic $\ZZ_p(I)$-module, 
to the orthogonal direct sum $U= U_1 \dotplus U_2$, 
where $U_i$ is a limit module  for $V_i$, for each $i=1, 2$.
In view of \eqref{fg.e25ii} we have 
\begin{equation}\mylabel{fg.e26}
U_1= C(I_{t+1} \tin U) \text { and } U_2 = [U, I_{t+1}],
\end{equation}
where $I_{t+1}$ is the image of $\mbqwi{t}_{2k+1}$ in $\Aut(\mbpi{t}_{2k})$ 
(see \eqref{fg.e23}).
In view of \eqref{fg.e24} we get $U_2= 0$.

According to Remark \ref{ma.r3}, the center 
$\mbri{t}$ of  $\mbpi{t}_{2k}$ is a cyclic central subgroup of 
$\mbgi{t}$,  maximal among the characteristic abelian subgroups of 
$\mbpi{t}_{2k}$. So $\mbri{t}$ is also a central subgroup of 
$\mbti{t} = \mbpi{t}_{2k} \rtimes \mbii{t}$,
as  $\mbii{t}$  is the image of 
$\mbqwi{t}\leq \mbgi{t}$ in $\Aut(\mbpi{t}_{2k})$. 
Thus Proposition \ref{lp.p3} applies to the faithful linear limit
$(\mbt, \mbr, \be, \mbp_{2k}, \ba_{2k})$. So 
 $U$ is isomorphic to 
$W^{\perp}/ W$ for some maximal $I$-invariant totally isotropic subspace $W$
of $V$. Then   $W = W_1 \dotplus W_2$,
where $W_i$ is a maximal totally isotropic $I$-invariant subspace of
$V_i$, for $i = 1,2$. But  $U_2 = 0$. Hence   $W_2^{\perp}$ must equal 
 $ W_2$. Thus
$V_2$ contains a self perpendicular $I$-invariant subspace. We conclude that 
$V_2$ is  hyperbolic as both  a $\ZZ_p(I)$- and  a
  $\ZZ_p(\mbqwi{t})$-module. The  group 
 $\mbqwi{t}$ has  $p$-power index in $\mbgi{t}$, since it is 
 a  $q$-Sylow subgroup  of the latter group. 
 Since  $V_2$ is an anisotropic  
symplectic  $\ZZ_p(\mbgi{t})$-module, 
Theorem (3.2) in \cite{da}, implies that $V_2$ is both  hyperbolic and 
anisotropic. Therefore $V_2$ is $0$. 
So $ V = V_1 = C( \mbqwi{t}_{2k+1} \tin V)$. 
Thus $\mbqwi{t}_{2k+1}$ centralizes  $\mbpi{t}_{2k+1}/ \mbri{t}$. 
We conclude that the $q$-Sylow subgroup $\mbqwi{t}_{2k+1}$ 
 of $\mbgi{t}_{2k+1} = \mbgi{t}_{t+1}$ centralizes 
the $p$-Sylow subgroup  $\mbpi{t}_{2k}$ of the same group. Hence 
$\mbgi{t}_{t+1}$ is nilpotent.
So  Theorem \ref{ma.main} follows
in the case  of an even $t$.
 
The proof for an odd $t$ follows immediately 
from the already proved case of an even $t$.
Indeed, assume that  $t$ is odd. Then  
 we can adjoin a trivial group and  character at the bottom 
of \eqref{fg.e1} so that $t$ becomes even. 
We now  interchange $p$ and $q$, and apply the 
already proved result. 
(Note that the normal series  (\ref{fg.e1}a) becomes
$1 \unlhd 1=H_1  \unlhd G_1= H_2 \unlhd G_2=H_3 \unlhd \dots \unlhd 
G_{t+1}=H_{t+2} \unlhd G$, 
so that  $1=H_1$ is assumed to be the first $p$-group of order $p^0=1$, 
while $G_1/1= H_2/H_1$  is a $q$-group and $H_{i+1}/H_i$  is either 
a $q$-group if $i$ is odd, or a $p$-group if $i$ is even, for all 
$i=1, \dots, t+1$.)

Hence Theorem \ref{ma.main} follows.

\section{Corollaries }

Below we list a series of corollaries following  Theorem \ref{ma.main}. 

\begin{corollary}\mylabel{fg.1}
Let $G$ be a finite $p^a q^b$-monomial  group, for some odd primes $p$ and $q$.
Assume that $N$ is a normal subgroup of $G$ and that 
 $\psi$ is an irreducible character of $N$. Consider the linear  quintuple
$(G, 1, 1, N, \psi)$.
Then there exists a faithful linear limit $(\mbg, \mba, \VP, \mbn, \VPS)$
of $(G, 1, 1, N, \psi)$ such that $\mbn$ is a nilpotent group.
\end{corollary}

Observe that  this is Theorem \ref{main} of the introduction.  
\begin{proof}
As $G$ is a solvable group and $N$ is normal subgroup of $G$, we can 
form a series 
\begin{equation}\mylabel{fg.e30}
1=G_0 \unlhd G_1 \unlhd G_2 \unlhd \dots 
\unlhd G_t \unlhd G_{t+1}= N \unlhd G_{t+2} \unlhd \dots \unlhd  G_n=G, 
\end{equation}
such that $G_i$ is a normal subgroup of $G$,  while the order of 
$G_{i+1}/ G_i$ is a power of a prime, for all $i=0, 1, \dots, t$.
Furthermore, without any loss of generality we can assume 
that $G_{i+1}/ G_i$ is 
$p$-group if $i$ is odd and a $q$-group if 
$i$ is even, for all $i=0,1, \dots, t$.
We also form recursively a character tower
\begin{equation}\mylabel{fg.e31}
\{ \chi_i \in \Irr(G_i) \}_{i=0}^{t+1}, 
\end{equation} 
for \eqref{fg.e30},  so that $\chi_{t+1} = \psi$ and 
$\chi_{i}$ is  any irreducible character of $G_i$ that lies under 
$\chi_{i+1} \in \Irr(G_{i+1})$, for all $i=0, 1, \dots, t$. 
In addition, we fix a representative of the unique conjugacy class of 
triangular sets that corresponds to the above character tower, along with 
a Sylow system for $G$ so that \eqref{sy.e2} holds.
 
We proceed with the series of reductions described in Section \ref{ma.se1}.
So after $t$ steps, we reach the limit groups 
$\mbgi{t}= \mbgi{t}_n,  \mbgi{t}_{t+1}$ and 
$\mbri{t}, \mbsi{t}$ along with their  limit characters 
$\Ti{t}_n, \Ti{t}_{t+1}$ and $\bei{t}$ and $\bzi{t}$.
 According to Theorem \ref{ma.main}
the group $\mbgi{t}_{t+1}$ is nilpotent.
Furthermore,  Remark \ref{ma.rrr3} implies that 
any faithful linear limit $(\mbg, \mba, \VP, \mbn, \VPS)$
 of $(\mbgi{t}, \mbri{t}, \bei{t}, \mbgi{t}_{t+1}, 
\Ti{t}_{t+1})$ is also a faithful linear limit of 
$(G, 1, 1, G_{t+1}, \chi_{t+1})= (G, 1, 1, N, \psi)$.
Hence if we take  $(\mbg, \mba, \VP, \mbn, \VPS)$ to be any faithful 
linear limit of  $(\mbgi{t}, \mbri{t}, \bei{t}, \mbgi{t}_{t+1}, 
\Ti{t}_{t+1})$, then $\mbn$ is nilpotent as a section  of $\mbgi{t}_{t+1}$.
So Corollary \ref{fg.1} follows.
\end{proof}

As any nilpotent  group is monomial, Proposition \ref{lim.p1}, along 
Corollary \ref{fg.1} applied to any irreducible character 
$\psi $ of $N$,   easily implies
\begin{corollary}\mylabel{fg.2}
Let $G$ be a finite $p^a q^b$-monomial  group, for some odd primes $p$ and $q$.
Assume that $N$ is a normal subgroup of $G$. Then $N$ is a monomial group.
\end{corollary}
Observe that this is Theorem \ref{main2} of the Introduction.

If, in addition, we take $N=G$  and $\chi \in \Irr(G)$,
 then Corollary \ref{fg.1} implies 
the existence of a faithful linear limit
$(\mbg, \mba, \VP, \mbg, \VPS)$ of $(G, 1, 1, G, \chi)$ so that 
$\mbg$ is nilpotent.
In view of Corollary \ref{lim.co2} the group $\mba= Z(\mbg)$ is maximal 
among the abelian $\mbg$-invariant subgroups of $\mbg$. As $\mbg$ is 
nilpotent this forces $Z(\mbg)=\mba= \mbg$. Hence  $\VPS = \VP \in \Irr(\mba)$
 is  linear. We conclude
\begin{corollary}
Let $G$ be an odd order monomial $p^a q^b$-group and let $\chi \in \Irr(G)$.
Then there exists a faithful linear limit  $\VPS$ of $\chi$ such that 
$\VPS(1)=1$, i.e., $\VPS$  is a linear character. 
\end{corollary}
Hence Theorem \ref{main3} follows.

%%% Local Variables: 
%%% mode: latex
%%% TeX-master: "485"
%%% End: 

\backmatter

\bibliographystyle{amsplain}

%%% Local Variables: 
%%% mode: latex
%%% TeX-master: t
%%% End: 

\appendix

%\include{A-data.tex}

%\vita
%Maria Loukaki was born in Iraklion Crete, Greece. She got her 
%Diploma degree
%B.A.  and  M.S.  degrees  
%from the Department of Mathematics at 
%the  University of Crete.

\end{document}